\documentclass[10pt]{amsbook}
\usepackage[letterpaper,body={16cm,22cm}, mag=1000]{geometry}

\usepackage{amsthm, amscd, imakeidx, tikz, tikz-cd, latexsym,hyperref, amssymb,enumerate, mathtools, wasysym, dcpic, pictexwd, graphicx, comment, adjustbox, yfonts, upgreek, float, tcolorbox, array, epsfig, import, dsfont, bbm, bbold, mathrsfs, chngcntr}
\usepackage[all,cmtip]{xy}
\makeindex[title=Index, columns=2]

\makeatletter
\@namedef{subjclassname@2020}{\textup{2020} Mathematics Subject Classification}
\makeatother
\newcommand{\para}{\par\vspace{.25cm}}

\newcommand\sbullet[1][.5]{\mathbin{\vcenter{\hbox{\scalebox{#1}{$\bullet$}}}}}

\counterwithout{figure}{chapter}
\counterwithout{table}{chapter}
\numberwithin{equation}{chapter}
\newtheorem{corollary}[equation]{Corollary}
\newtheorem{lemma}[equation]{Lemma}
\newtheorem{problem}[equation]{Problem}
\newtheorem{proposition}[equation]{Proposition}
\newtheorem{question}[equation]{Question}
\newtheorem{exercise}[equation]{Exercise}
\newtheorem{conjecture}[equation]{Conjecture}
\newtheorem{theorem}[equation]{Theorem}
\newtheorem{definition}[equation]{Definition}
\newtheorem{example}[equation]{Example}

\newtheorem{remark}[equation]{Remark}

\newcommand{\dlabel}[1]{\ifmmode \text{\ttfamily \upshape [#1] } \else
{\ttfamily \upshape [#1] }\fi \label{#1}}

\newcommand{\B}{\operatorname{B} }
\newcommand{\WB}{\operatorname{WB} }
\newcommand{\C}{\operatorname{C} }

\newcommand{\Stab}{\operatorname{Stab} }
\newcommand{\Fix}{\operatorname{Fix} }

\newcommand{\GL}{\operatorname{GL} }
\newcommand{\Sl}{\operatorname{SL} }
\newcommand{\BQ}{\operatorname{BQ} }

\newcommand{\Ho}{\operatorname{H} }
\newcommand{\T}{\operatorname{T} }
\newcommand{\IA}{\operatorname{IA} }

\newcommand{\J}{\operatorname{J} }
\newcommand{\Oo}{\operatorname{O} }
\newcommand{\Uu}{\operatorname{U} }
\newcommand{\Tt}{\operatorname{T} }
\newcommand{\UT}{\operatorname{UT} }
\newcommand{\rank}{\operatorname{rank} }
\newcommand{\charac}{\operatorname{char} }
\newcommand{\M}{\operatorname{M} }
\newcommand{\R}{\operatorname{R} }
\newcommand{\Ret}{\operatorname{Ret} }
\newcommand{\cl}{\operatorname{cl}}
\newcommand{\Autcent}{\operatorname{Autcent}}
\newcommand{\Z}{\operatorname{Z} }
\newcommand{\id}{\operatorname{id}}

\newcommand{\cw}{\operatorname{cw}}
\renewcommand{\det}{\operatorname{det} }
\newcommand{\Inn}{\operatorname{Inn} }
\newcommand{\Ann}{\operatorname{Ann} }
\newcommand{\QInn}{\operatorname{QInn} }
\newcommand{\Aut}{\operatorname{Aut} }
\newcommand{\diag}{\operatorname{diag} }
\newcommand{\Trans}{\operatorname{Trans} }
\newcommand{\Dis}{\operatorname{Dis} }
\newcommand{\Out}{\operatorname{Out} }
\newcommand{\Hol}{\operatorname{Hol} }

\newcommand{\End}{\operatorname{End} }
\newcommand{\Hom}{\operatorname{Hom} }
\newcommand{\Fib}{\operatorname{Fib} }
\newcommand{\lcm}{\operatorname{lcm} }
\newcommand{\Tot}{\operatorname{Tot} }
\newcommand{\Soc}{\operatorname{Soc} }
\newcommand{\Ext}{\operatorname{Ext} }
\newcommand{\Tor}{\operatorname{Tor} }
\newcommand{\im}{\operatorname{im} }

\newcommand{\Conj}{\operatorname{Conj} }
\newcommand{\Map}{\operatorname{Map} }
\newcommand{\map}{\operatorname{map} }
\newcommand{\Core}{\operatorname{Core} }
\newcommand{\Alex}{\operatorname{Alex} }
\newcommand{\Adj}{\operatorname{Adj} }

\begin{document}
\bigskip
\bigskip
\bigskip
\title{\Huge{Yang--Baxter Equation\\ and\\ Related Algebraic Structures}
\bigskip
\bigskip
\bigskip
\bigskip
\bigskip
\bigskip
\begin{figure}[!ht]
 \begin{center}
\includegraphics[height=12cm]{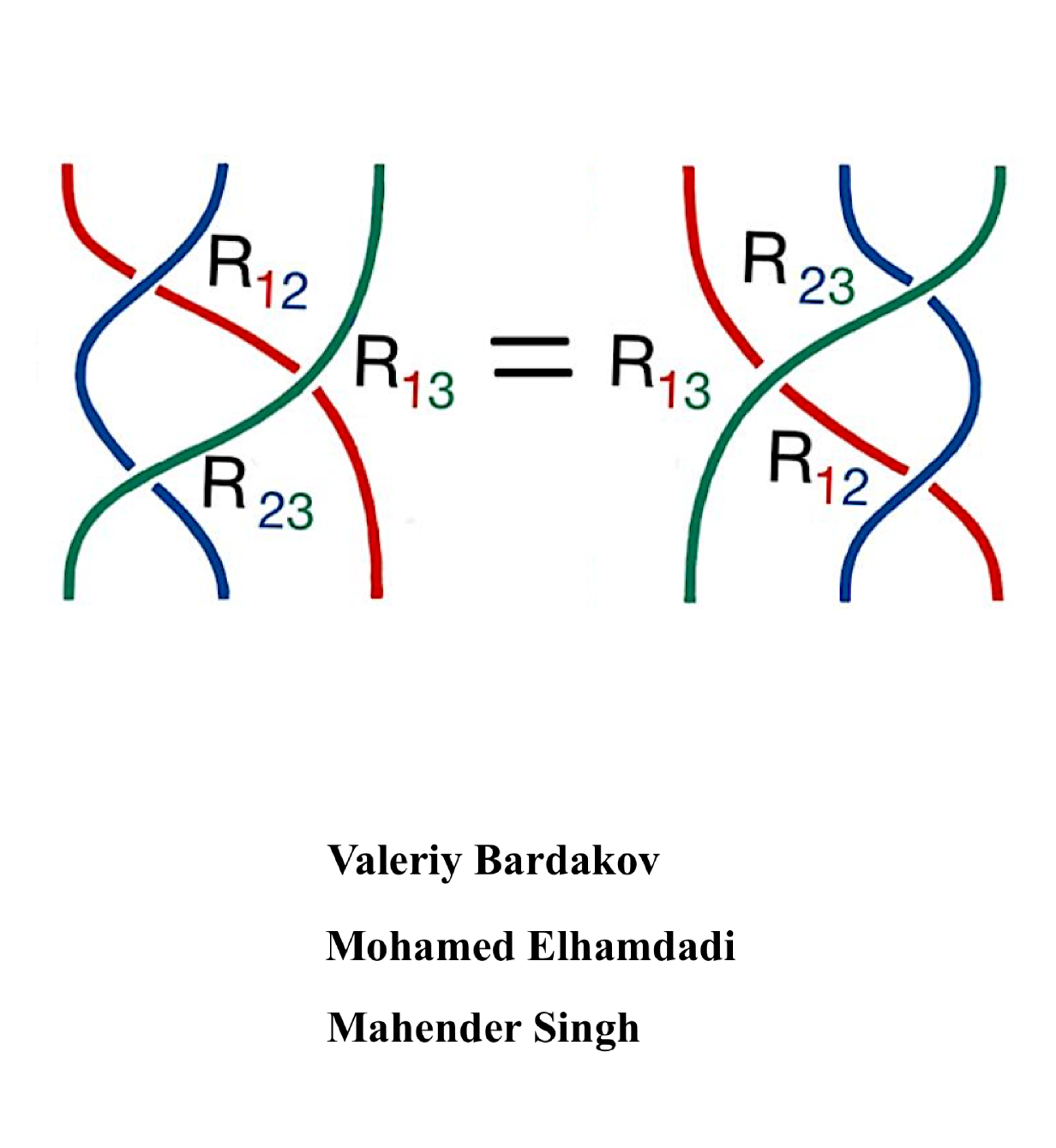}
\end{center}
\end{figure}
}
\maketitle

\noindent \author{Valeriy G. Bardakov\\
Sobolev Institute of Mathematics, Siberian Branch of the Russian Academy of Sciences, 4 Acad. Koptyug avenue, 630090 Novosibirsk Russia.\\
bardakov@math.nsc.ru}
\bigskip
\bigskip
\bigskip

\noindent\author{Mohamed Elhamdadi\\
Department of Mathematics and Statistics, University of South Florida, Tampa, FL 33620, USA.\\
emohamed@usf.edu}
\bigskip
\bigskip
\bigskip

\noindent\author{Mahender Singh\\
Department of Mathematical Sciences, Indian Institute of Science Education and Research Mohali, Sector 81, Knowledge City, SAS Nagar, Punjab 140306, India.\\
mahender@iisermohali.ac.in}
\bigskip
\bigskip
\bigskip

\bigskip
\bigskip

\noindent Keywords: Affine solution, Alexander quandle, annihilator of skew brace, adjoint group of quandle, braiding operator, commutator width, conjugation quandle, cycle set, Dehn quandle, derived structure group, dihedral quandle, dynamical cocycle, exact factorisation, free product of quandles, flat quandle, free quandle, free rack, group automorphism, group cohomology, holomorph of group, homotopy group, idempotent, involutive solution, indecomposable solution, irreducible 3-manifold, isomorphism problem, lambda map, latin quandle, Lie ring, linear solution,  link group, link quandle, left brace, morphism of solutions, non-degenerate solution, orderable quandle, permutation rack, quandle, quandle algebra, quandle automorphism, quandle cocycle, quandle cohomology, quandle covering, quandle coloring, quandle extension, quandle module, quandle ring, quandle space, orderable quandle, permutation rack, pre-Lie algebra, rack, rack space, radical ring, regular subgroup, residual finiteness, Rota--Baxter group, Rota--Baxter operator, skew left brace, socle of skew brace, $\square$-set, structure group of solution, symmetric skew brace, Takasaki quandle, torus link, trunk, Wells-type sequence, Yang--Baxter cohomology, Yang--Baxter equation, zero-divisor.
\maketitle

\tableofcontents

\section*{Notation}\par\vspace{1cm}
\begin{tabular}{l@{\quad\dots\quad}p{.70\textwidth}}
$\mathbb{N}$ & The set of natural numbers\\
$\mathbb{Z}$ & The ring of integers\\
$\mathbb{Q}$ & The field of rational numbers\\
$\mathbb{R}$ & The field of real numbers\\
$\mathbb{C}$ & The field of complex numbers\\
$\mathbb{C}^\times$ & The multiplicative group of non-zero complex numbers\\
$\mathbb{S}^1$ & The unit circle in the Euclidean 2-space\\
$\mathbb{S}^3$ & The unit sphere in the Euclidean 4-space\\
$|X|$ & The cardinality of the set $X$\\
$X \times Y$ & The cartesian product of $X$ and $Y$\\ 
$X^n$ & The $n$-fold cartesian product of $X$\\ 
$\prod_i X_i$ & The cartesian product of the family $\{X_i \}$ of sets\\
$\sup(X)$ & The supremum of the subset $X$ of real numbers\\
$\lceil x \rceil$ & The smallest integer that is greater than or equal to $x$.\\
$\lfloor x \rfloor$ & The largest integer that is less than or equal to $x$.\\
$f|_A$ & The restriction of the map $f$ on the subset $A$ of its domain\\
$\Map(X, Y)$ & The set of all maps from the set $X$ to the set $Y$\\
$\Fix(f)$ & The fixed-point set of the map $f:X \to X$\\
$gf(x)$ & The element $g \big(f(x)\big)$, where $f: X \to Y$ and $g:Y \to Z$ are maps and $x \in X$\\
$\det (A)$ & The determinant of the square matrix $A$\\
$1$ & The trivial group\\
$\mathbb{Z}_n$ & The cyclic group of order $n$\\
$\Sigma_n$ & The symmetric group on $n$ symbols\\
$\Sigma_X$ & The symmetric group on the set $X$\\
$\B_n$ & The braid group on $n$ strands\\
$\WB_n$ & The welded braid group on $n$ strands\\
$\Hol(G)$ & The holomorph of the group $G$\\
$\ker(\phi)$ & The kernel of the group homomorphism $\phi:G \to H$\\
$\im(\phi)$ & The image of the map $\phi:X \to Y$\\
$[x,y]$ &  The commutator $x y x^{-1} y^{-1}$ of $x$ and $y$\\
$[x_1, \dots, x_n]$ & The left-normed commutator $[\dots [[x_1, x_2], x_3], \dots, x_n]$ of $x_1, \dots, x_n$\\
$H \le G$ & The subgroup $H$ of the group $G$\\
$H \trianglelefteq G$ & The normal subgroup $H$ of the group $G$\\
$[A, B]$ & The subgroup generated by commutators $aba^{-1}b^{-1}$, where $a \in A$ and $b \in B$\\
$[G, G]$ & The commutator subgroup of the group $G$\\
$\gamma_i(G)$ & The $i$-th term of the lower central series of $G$, where $\gamma_2(G)=[G,G]$\\
$\Z(G)$ & The center of the group $G$\\
$[G:H]$ & The index of the subgroup $H$ in the group $G$\\
$G_{\rm ab}$ & The abelianisation of the group $G$\\
$\C_G(H)$ & The centralizer of the subset $H$ in the group $G$\\
$\C_G(x)$ & The centralizer of the element $x$ in the group $G$\\
$G_x$ & The stabilizer of the element $x$ under an action of the group $G$\\
$\BQ(m, n)$ & The Burnside quandle\\
$\Fib(n)$ & The Fibonacci  quandle generated by $n$ elements\\
$\mathcal{D}(A^G)$ & The Dehn quandle of $G$ with respect to the subset $A$\\
$\R_n$ & The dihedral quandle of order $n$\\
$\T_n$ & The trivial quandle of order $n$\\
$\T_\infty$ & The trivial quandle with countably infinite number of elements\\
$\Conj_n(G)$ & The $n$-th conjugation quandle of the group $G$\\
$\Conj(G)$ & The conjugation quandle of the group $G$\\
\end{tabular}

\section*[Notation]{}\par\vspace{1cm}
\begin{tabular}{l@{\quad\dots\quad}p{.70\textwidth}}
$\Core(G)$ & The core quandle of the group $G$\\
$\Alex(G, \phi)$ & The generalised Alexander quandle of the group $G$ with respect to the automorphism $\phi$\\
$\Adj(X)$ & The adjoint group of the quandle $X$\\
$G(X, r)$ & The structure group of the solution $(X, r)$ to the Yang--Baxter equation\\
$\textswab{a}_x$ & The generator of $\Adj(X)$ corresponding to the quandle element $x \in X$\\
$\mathbb{Z}^{X}$ & The free abelian group on the set $X$\\
$\mathbb{Z}^\infty$ & The free abelian group of countably infinite rank\\
$\mathbb{Z}^n$ & The free abelian group of rank $n$\\
$F(S)$ & The free group on the set $S$\\
$F_n$ & The free group of rank $n$\\
$F_{\infty}$ & The free group of countably infinite rank\\
$FQ(S)$ & The free quandle on the set $S$\\
$FQ_n$  & The free quandle of rank $n$\\
$FQ_\infty$  & The free quandle of countably infinite rank\\
$FR(S)$ & The free rack on the set $S$\\
$FR_n$ & The free rack of rank $n$\\
$Q(L)$ & The link quandle of the link $L$ in $\mathbb{S}^3$\\
$G(L)$ & The link group $\pi_1(\mathbb{S}^3 \setminus L)$ of the link $L$ in $\mathbb{S}^3$\\
$X \star Y$ & The free product of quandles $X$ and $Y$\\ 
$G \hexstar H$ & The free product of groups $G$ and $H$\\
$G \rtimes H$ & The semi-direct product of groups $G$ and $H$, where $H$ acts on $G$ by automorphisms \\
$G \ltimes H$ & The semi-direct product of groups $G$ and $H$, where $G$ acts on $H$ by automorphisms \\
$\iota_x$ & The inner automorphism of the group $G$ induced by the element $x \in G$\\
$S_x$ & The inner automorphism of the quandle (or rack) $X$ induced by the element $x \in X$\\
$L_x$ & The left multiplication of the quandle (or rack) $X$ induced by the element $x \in X$\\
$V \otimes W$ & The tensor product of vector spaces $V$ and $W$ over a field\\
$G \otimes H$ & The tensor product of abelian groups $G$ and $H$ viewed as $\mathbb{Z}$-modules\\
$\Tor(G, H)$ & The Tor group of abelian groups $G$ and $H$ viewed as $\mathbb{Z}$-modules\\
$\Ext(G, H)$ & The Ext group of abelian groups $G$ and $H$ viewed as $\mathbb{Z}$-modules\\
$\M(n, \,\mathbb{k})$ & The ring of $n \times n$ matrices over the ring $\mathbb{k}$\\
$\GL(n, \,\mathbb{k})$ & The group of $n \times n$ invertible matrices over the field $\mathbb{k}$\\
$\Oo(n, \,\mathbb{k})$ & The group of $n \times n$ orthogonal matrices over the field $\mathbb{k}$\\
$\Sl(n, \,\mathbb{k})$ & The group of $n \times n$ special linear matrices over the field $\mathbb{k}$\\
$\Uu(n, \,\mathbb{C})$ & The group of $n \times n$ unitary matrices over $\mathbb{C}$\\
$\Tt(n, \,\mathbb{C})$ & The group of $n \times n$ upper triangular matrices over $\mathbb{C}$\\
$\UT(n, \mathbb{Z})$ & The group of  $n \times n$ upper unitriangular matrices over $\mathbb{Z}$\\
$\Ho^n_{\rm Grp}(G;\,A)$ & The $n$-th  group cohomology of the group $G$ with coefficients in the $G$-module $A$\\
$\Ho_n^{\rm Grp}(G;\,A)$ & The $n$-th  group homology of the group $G$ with coefficients in the $G$-module $A$\\
$\Ho^n_{\rm R}(X;\,A)$ & The $n$-th  rack cohomology of the rack $X$ with coefficients in the group $A$\\
$\Ho_n^{\rm R}(X;\,A)$ & The $n$-th  rack homology of the rack $X$ with coefficients in the group $A$\\
$\Ho^n_{\rm Q}(X;\,A)$ & The $n$-th  quandle cohomology of the quandle $X$ with coefficients in the group $A$\\
$\Ho_n^{\rm Q}(X;\,A)$ & The $n$-th  quandle homology of the quandle $X$ with coefficients in the group $A$\\
\end{tabular}

\section*[Notation]{}\par\vspace{1cm}
\begin{tabular}{l@{\quad\dots\quad}p{.70\textwidth}}
$\Ho^n_{\rm Cell}(X;\,A)$ & The $n$-th  cellular cohomology of the CW-complex $X$ with coefficients in the group $A$\\
$\Ho_n^{\rm Cell}(X;\,A)$ & The $n$-th  cellular homology of the CW-complex $X$ with coefficients in the group $A$\\
$\Autcent(G)$ & The group of central automorphisms of the group $G$\\
$\End(G)$ & The endomorphism ring of the abelian group $G$\\
$\Hom_{\mathcal{G}}(G,\, H)$ & The set of group homomorphisms from the group $G$ to the group $H$\\
$\Hom_{\mathcal{Q}}(X,\, Y)$ & The set of quandle homomorphisms from the quandle $X$ to the quandle $Y$\\
$\mathcal{S}_g$ & The closed orientable surface of genus $g \ge 1$\\  
$\mathcal{M}_g$ & The mapping class group of the surface $\mathcal{S}_g$\\
$T_\alpha$ & The Dehn twist along the simple closed curve $\alpha$ on a surface\\
$\pi_1(X, x)$ & The fundamental group of the space $X$ based at the point $x \in X$\\
$K(\pi, 1)$ & The Eilenberg-Maclane space with fundamental group $\pi$ and trivial higher homotopy groups\\
$\mathbb{k}^\times$ & The group of units of the associative ring $\mathbb{k}$ with unity\\
$\charac(\mathbb{k})$ & The characteristic of the ring $\mathbb{k}$\\
$\mathbb{k}\langle\langle X_1, \ldots, X_n \rangle\rangle$  & The ring of formal power series over $\mathbb{k}$ in non-commutative variables $X_1, \ldots, X_n$\\
$\mathbb{k}[G]$ & The group ring of the group $G$ with coefficients in the ring $\mathbb k$ \\
$\mathbb{k}[X]$ & The quandle ring of the quandle $X$ with coefficients in the ring $\mathbb k$ \\
$\mathbb{k}[X]'$ & The commutator subalgebra  of the quandle ring $\mathbb{k}[X]$\\
$e_x$ & The basis element of $\mathbb{k}[X]$ corresponding to the element $x \in X$\\
$\hat{\varphi}$ & The homomorphism $\mathbb{k}[X] \to \mathbb{k}[Y]$ of quandle rings induced by the quandle homomorphism $\varphi:X \to Y$\\
$\cw \big(\mathbb{k}[X] \big)$ & The commutator width of the quandle ring  $\mathbb{k}[X]$\\
$\Aut(X)$ & The group of automorphisms of $X$, where $X$ is a group, quandle, rack or  skew brace\\
$\Inn(X)$ & The group of inner automorphisms of $X$, where $X$ is a group, quandle or  rack\\
$\QInn(X)$ & The group of quasi-inner automorphisms of $X$, where $X$ is a group, quandle or rack\\
$\mathcal{O}(X)$ & The set of orbits of action of $\Inn(X)$ on the quandle $X$\\
$\Ret(X, r)$ & The retraction of the solution $(X, r)$\\
$\mathcal{S}$ & The category of sets\\
$\mathcal{G}$ & The category of groups\\
$\mathcal{A}$ & The category of abelian groups\\
$\mathcal{Q}$ & The category of quandles\\
$\mathcal{SLB}$ & The category of skew left braces\\
$\mathcal{RBG}$ & The category of Rota--Baxter groups\\

\end{tabular}



\chapter*{}
\vspace*{5cm}
\begin{center}
\begin{tcolorbox}[width=15cm]
\it 
This monograph is dedicated \\
by Valeriy to his beloved children, Oksana and Sasha; \\
by Mohamed to his cherished children, Lina and Hamza; and \\
by Mahender to his loving daughters, Aashna and Avni.
\end{tcolorbox}
\end{center}

\setcounter{chapter}{-1}

\chapter{Introduction}

The quantum Yang--Baxter equation is a fundamental equation in mathematical physics that arises in the study of integrable systems, quantum groups, and quantum information theory. The equation first appeared in the works of Yang \cite{MR0261870} and Baxter \cite{MR0290733}. A solution to the quantum Yang--Baxter equation is a linear map 
$$R:V \otimes V \to V \otimes V,$$
where $V$ is a vector space over a field, such that the equality 
$$R_{12} \,R_{13}  \,R_{23} = R_{23} \,R_{13}  \,R_{12}$$
of operators hold in the space $\End(V \otimes V \otimes V)$. Here, $R_{ij}$ denotes the linear map $V\otimes V \otimes V \to V\otimes V \otimes V$ that acts as $R$ on the $(i,j)$-th tensor factors and as the identity on the remaining third factor.  The equation $R_{12} \,R_{13}  \,R_{23} = R_{23} \,R_{13}  \,R_{12}$ is referred to as the quantum Yang--Baxter equation. The significance of this equation lies in its connection to important concepts in physics and mathematics. For instance, it has applications in the study of exactly solvable models in statistical mechanics and quantum field theory, and has connections to knot theory,  braid theory and Hopf algebras. A central problem in the subject is to construct and classify solutions of this equation.
\para

\begin{figure}[!ht]
 \begin{center}
\includegraphics[width=5.5cm]{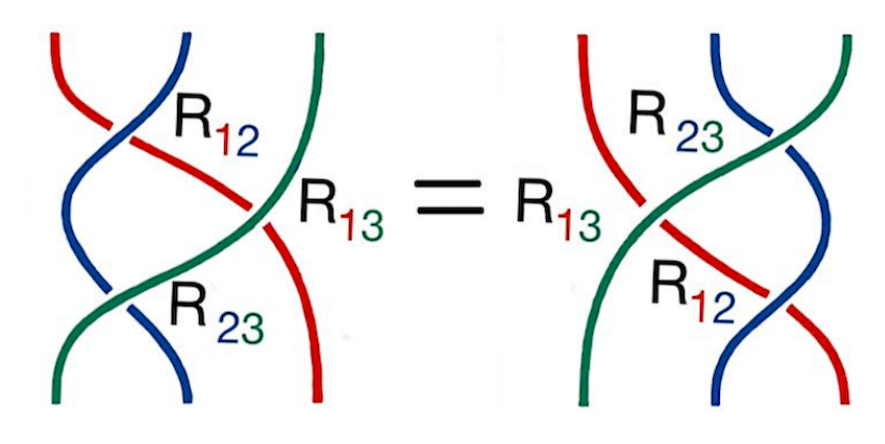}
\end{center}
\caption{The third Reidemeister move.}\label{fig1}
\end{figure}

The linear map $\mathtt{t}: V \otimes V  \to V \otimes V$ given by $\mathtt{t}(v \otimes w)= w \otimes v$ is clearly a solution to the quantum Yang--Baxter equation. A direct check shows that a linear map $R: V \otimes V  \to V \otimes V$ is a solution to the quantum Yang--Baxter equation if and only if the linear map $S = \mathtt{t}~ R$ satisfies the braid equation
$$S_{12}  \,S_{23}  \,S_{12} = S_{23} \, S_{12} \,S_{23}.$$
In this case, $S$ is said to be a solution to the Yang--Baxter equation or the braid equation. Topologically, the braid equation is simply the third \index{Reidemeister move}{Reidemeister move} of planar diagrams of links in the Euclidean 3-space as shown in Figure \ref{fig1}.
\para

Let $X$ be a basis of a vector space $V$ over some field. Then a map $r : X \times X \to X \times X$ satisfying the braid equation
\begin{equation}\label{set-theoretic solution of quantum YBE intro}
r_{12} \, r_{23}  \, r_{12} = r_{23}  \, r_{12}  \, r_{23},
\end{equation}
where $r_{ij}: X \times X \times X \to X \times X \times X$ acts as $r$ on the $(i,j)$-th factors and as the identity on the remaining third factor, induces a solution to the Yang--Baxter equation on the vector space spanned by $X$. The pair $(X, r)$ is referred to as a set-theoretic solution to the Yang--Baxter equation, and the map $r$ is called a braiding. In \cite[Section 9]{MR1183474}, Drinfel'd proposed to investigate set-theoretic solutions to the quantum Yang--Baxter equation. Clearly, equation \eqref{set-theoretic solution of quantum YBE intro} can be adapted to various contexts. For instance, when $X$ is a manifold, then we desire that  $r : X \times X \to X \times X$ is a smooth map. For $n \ge 2$, let 
$$\B_n = \big\langle \sigma_1,\ldots,\sigma_{n-1}  \, \mid \,  \sigma_i\sigma_{i+1}\sigma_i = \sigma_{i+1} \sigma_i \sigma_{i+1}  ~\textrm{for}~1 \le i \le n-2~\textrm{and}~\sigma_i\sigma_j = \sigma_j\sigma_i ~\textrm{for}~|i-j| \ge 2 \big\rangle$$
be the Artin braid group on $n$-strands.  If $(X, r)$ is a bijective set-theoretic solution to the Yang--Baxter equation (that is, the braiding $r$ is a bijection) and $r_i : X^n \to X^n$ is given by
$$
r_i = \id_{X^{i-1}} \times r \times \id_{X^{n-i-1}},
$$
then the assignment $\sigma_i \mapsto r_i$ extends to an action of $\B_n$ on $X^n$ for each $n \geq 2$. Passing to the vector space $V$ spanned by $X$ gives a linear representation $\B_n \to \Aut(V ^{\otimes n})$. In \cite{MR0908150, MR0766964}, Jones introduced a new polynomial invariant of links using certain linear representations of the braid group $\B_n$ and a Markov trace. In general, using trigonometric solutions to the Yang--Baxter equation, a series of linear representations of $\B_n$ of the aforementioned form can be obtained. In \cite{MR0939474}, Turaev systematically investigated these representations to obtain a state functional (or the Markov trace), which was employed to establish an invariant for links. It was discovered that this invariant corresponds to a specialization of both the two-variable Jones polynomial and the Kauffman polynomial.
\para 

Shortly after the publication of \cite{MR1183474}, Weinstein and Xu \cite{MR1178147} constructed a solution in the setting of symplectic geometry. Let $\mathcal{R}$ be a Lagrangian submanifold in the Cartesian square $\Gamma \times \Gamma$ of a symplectic groupoid $\Gamma$ associated to a Poisson Lie group $G$. Then, for any symplectic leaf $H$ in $G$, the submanifold $\mathcal{R}$ induces a symplectic automorphism $\sigma: H \times H \to H \times H$, which is a set-theoretic solution to the Yang--Baxter equation. 
\para 

If $(X, r)$ is a set-theoretic solution to the Yang--Baxter equation, then we can write
$$r(x, y) = \big(\sigma_x(y), \, \tau_y(x)\big)$$
for $x,y \in X$. The solution $(X, r)$ is said to be right non-degenerate if the map $\sigma_x$ is bijective for all $x \in X$. Similarly, $(X, r)$ is said to be left non-degenerate if the map $\tau_x$ is bijective for all $x \in X$. And, the solution is called non-degenerate if it is both left as well as right non-degenerate. It turns out that if $(X,r)$ is a set-theoretic solution to the Yang--Baxter equation with $X$ finite, then the solution is non-degenerate if and only if it is either left or right non-degenerate. The class of non-degenerate solutions is studied for applications of the Yang--Baxter equation in physics as well as for its connection with various algebraic structures in mathematics. A set-theoretic solution $(X, r)$ is said to be involutive if $r^2= \id_{X \times X}$. Note that, in this case, the action of the braid group $\B_n$ descends to an action of the symmetric group $\Sigma_n$ on $X^n$. This action is twisted, and in general, is different from the usual permutation action of $\Sigma_n$, but the two actions are conjugate \cite{MR1722951}. 
\para

Let $S$ be the semigroup generated by $\{x_1, \ldots, x_n\}$ and $\mathcal{U}$ the multiplicative free commutative semigroup generated by $\{u_1, \ldots, u_n\}$. Then $S$ is said to be of $I$-type if there is a bijection $\phi:\mathcal{U}\to S$ such that $ \big(\phi(u_1a), \ldots, \phi(u_na) \big)= \big(x_1 \phi(a), \ldots, x_n \phi(a) \big)$ for all $a\in \mathcal{U}$. In \cite[Theorem 1.3]{MR1637256},
Gateva-Ivanova and Van den Bergh proved  that if $S=\langle X \,\mid\, R \rangle$ is a semigroup of $I$-type, then there is an involutive set-theoretic solution $(X, r)$ to the Yang--Baxter equation with some conditions. In fact, they proved that the converse also holds, and that these semigroups are related to Bieberbach groups.
\para

In \cite{MR1722951}, Etingof, Schedler, and Soloviev conducted a study of set-theoretic solutions to the quantum Yang--Baxter equation, focusing on solutions that satisfy additional conditions of bijectivity, involutivity, and non-degeneracy. They discussed the geometric and algebraic interpretations of such solutions, introduced several new constructions of them, and gave their classification in terms of groups. Most notably, they introduced the structure group 
$$G(X, r)= \big\langle X ~\mid ~ xy = \sigma_x(y) \tau_y(x) ~\textrm{for all}~x, y \in X \big\rangle$$
of a set-theoretic solution $(X, r)$ to the Yang--Baxter equation, which turned out to be the most prominent invariant of a solution. The structure group $G(X, r)$ has two actions on the underlying set $X$ which are conjugate to each other. They proved that the structure group of a finite non-degenerate involutive solution is solvable \cite[Theorem 2.15]{MR1722951}.  A set-theoretic solution is called indecomposable if it cannot be written as a union of two non-empty disjoint invariant subsets that are themselves solutions. It was proved in \cite[Theorem 2.13]{MR1722951} that for $X=\mathbb{Z}_p$, where $p$ is a prime, there is a unique (up to isomorphism) indecomposable  non-degenerate involutive solution. In fact, the solution is $(\mathbb{Z}_p , r)$, where  $r: \mathbb{Z}_p \times \mathbb{Z}_p  \to \mathbb{Z}_p \times \mathbb{Z}_p $ is given by $$r(x, y) = (y-1, ~x+1).$$ Further, a classification of non-degenerate involutive solutions was given in terms of groups with bijective 1-cocycles. In a subsequent work  \cite{MR1848966}, it was proved that all indecomposable non-degenerate set-theoretic solutions to the Yang--Baxter equation on a set of prime order are affine.
\para

The presentation of the structure group of a set-theoretic solution to the Yang--Baxter equation, in fact, gives the presentation of a monoid. In a rather curious connection, Chouraqui proved in \cite{MR2764830} that the structure group (or the monoid) associated to a non-degenerate involutive set-theoretic solution to the  Yang--Baxter equation is Garside. And conversely, any Garside monoid $M$ that has a monoid presentation $M=\langle X\mid R \rangle$, under certain conditions on the set of relations $R$, can be realized in this way. Furthermore, it was established that the (co)homological dimension of the structure group $G(X, r)$ equals the cardinality $|X|$ of the set $X$. In \cite{MR2994025}, the connection between subgroups of the structure group and corresponding solutions was explored. It was proved that, for standard parabolic subgroups of the structure group, there is a one-to-one correspondence with non-degenerate invariant subsets of the solution set.
\para

The constructions presented in \cite{MR1722951, MR1178147} were generalised to a broader setting by Lu, Yan, and Zhu in  \cite[Theorem 1]{MR1769723}.  Let $G$ be a group equipped with a left-action $\xi$ and a right-action $\zeta$ of $G$ on itself, denoted by $(u,v) \mapsto \xi(u) \cdot v$ and $(u,v) \mapsto u \cdot \zeta(v)$, respectively. If the two actions satisfy the compatibility condition
$$uv = \big(\xi(u) \cdot v \big) \big(u \cdot \zeta(v)\big),$$
then the map $r: G \times G \to G \times G$, given by $$r(u, v) = \big(\xi(u) \cdot v,~ u\cdot \zeta(v)\big),$$
yields an invertible set-theoretic solution $(G, r)$ to the Yang--Baxter equation. They also proved that the structure group $G(X, r)$ of a non-degenerate bijective set-theoretic solution $(X,r)$ to the Yang--Baxter equation satisfies some universal property \cite[Theorem 4]{MR1769723}. More precisely, if $i_X : X \to G(X,r)$ is the natural map, then there is a unique braiding $r^G$ on $G(X, r)$ such that $$r^G ~(i_X \times i_X) = (i_X \times i_X)~r.$$ Further, if $s$ is a braiding on another group $H$ and $f: X \to H$ is a braiding-preserving map, then there is a unique braiding-preserving group homomorphism $\tilde{f}:G(X, r) \to H$ such that $f=\tilde{f}~ i_X$. In case of involutive solutions, the natural map $i_X : X \to G(X,r)$ is injective. Further, it was noticed in \cite{MR1769723} that in the absence of involutivity condition on the solution, the natural map $i_X : X \to G(X,r)$ need not be injective.
\para

The ideas developed in \cite{MR1722951, MR1769723} were extended by Soloviev in \cite{MR1809284} to encompass non-degenerate set-theoretic solutions to the Yang--Baxter equation. In particular, the involutivity condition assumed in \cite{MR1722951} was relaxed in order to obtain a group-theoretical characterization of such solutions \cite[Theorem 2.7]{MR1809284}. The same work also investigated injective solutions $(X, r)$, that is, ones for which the natural map $i_X : X \to G(X,r)$ is injective. A combinatorial criterion was provided to describe the class of injective non-degenerate solutions \cite[Theorem 2.9]{MR1809284}. Since injectivity is a generalization of involutivity,  the classification given in \cite{MR1722951} was generalised to injective solutions.
\para

In \cite{MR2132760}, Rump proved that left non-degenerate involutive solutions are equivalent to cycle sets, where a cycle set is a set $X$ equipped with a left invertible binary operation $\cdot$ satisfying the condition
\begin{equation}\label{cycle set condition}
(x \cdot y) \cdot (x \cdot z) = (y \cdot x) \cdot (y \cdot z)
\end{equation}
for all $x, y, x \in X$. In fact, a finite cycle set $X$ gives rise to a (left and right) non-degenerate solution and can be naturally extended to the free abelian group $\mathbb{Z}^{X}$ generated by $X$. The condition \eqref{cycle set condition} is then replaced by the condition
\begin{equation}\label{linear cycle set condition}
(a + b) \cdot c = (a \cdot b) \cdot (a \cdot c)
\end{equation}
for all $a, b, c \in \mathbb{Z}^{X}$. An abelian group $A$ with a left distributive multiplication which makes $A$ into a cycle set satisfying condition \eqref{linear cycle set condition} is called a linear cycle set. Rump showed in \cite{MR2278047} that linear cycle sets are closely related to radical rings, which led him to introduce braces, which are at the center-stage of research on the Yang--Baxter equation at the time of writing of this monograph. More precisely, Rump defined a brace as an abelian group $A$ with a right distributive multiplication such that $A$ is a group with respect to the operation
$$a \circ b := ab + a + b.$$
It turns out that a brace is left distributive if and only if it is a radical ring. A Galois correspondence between ideals of braces and quotient cycle sets was also established in \cite{MR2278047}. It is worth noting that Kurosh had already examined braces during his lectures \cite[Section 10]{MR0392756}.
\para

In \cite{MR3177933}, Ced\'{o}, Jespers, and Okni\'{n}ski reformulated Rump's definition of a brace into a more accessible form. A \index{left brace} left brace $(B, +, \circ)$ is a set $B$ with two operations $+$ and $\circ$ such that $(B,+)$ is an abelian group, $(B, \circ)$ is a group and the condition
\begin{equation}\label{left brace defining axiom}
  a \circ (b + c)  = (a\circ b) - a + (a \circ c)
\end{equation}
holds for all $a,b,c \in B$. The groups $(B,+)$ and $(B, \circ)$ are called  the additive group and the multiplicative group of the left brace $(B, +, \circ)$, respectively. A right brace is defined analogously, but instead of \eqref{left brace defining axiom}, we require the condition
\begin{equation}\label{right brace defining axiom}
  (b + c) \circ a = (b\circ a) - a + (c \circ a)
\end{equation}
for all $a,b,c \in B$.  If $(B,+,\circ)$ is both a left and a right brace, then we say that it is a two-sided brace. It turns out that two-sided braces are in a bijective correspondence with Jacobson radical rings. Also, there is also a bijective correspondence between left and right braces, which allows us to only consider left braces. Each left brace gives a non-degenerate involutive set-theoretic solution to the Yang--Baxter equation. More precisely, if $(B,+,\circ)$ is a left brace and $\lambda: (B, \circ) \to \Aut(B, +)$ is the homomorphism given by $\lambda_a(b)= -a+ (a\circ b)$, then the map $r:B \times B \to B \times B$ defined by
$$r(a,b)=\big(\lambda_a(b),~\lambda_{\lambda_a(b)}^{-1}(a)\big)$$
gives a non-degenerate involutive set-theoretic solution $(B,r)$ to the Yang--Baxter equation. Conversely, given a non-degenerate involutive set-theoretic solution $(X, r)$ to the Yang--Baxter equation, one can define a left brace structure on its structure group $G(X, r)$. It turns out that the structure of a group arising as the multiplicative group of a left brace cannot be arbitrary. It follows from \cite[Theorem 2.25]{MR1722951} and \cite[Theorem 2.1]{MR2584610} that  if $(B,+,\circ)$ is a finite left brace, then $(B,\circ)$ is solvable. In \cite{MR3574204}, Ced\'{o}, Gateva-Ivanova, and Smoktunowicz introduced left and right nilpotent left braces. Later, in \cite{MR3814340}, Smoktunowicz proved  that a finite left brace $(B,+,\circ)$ is nilpotent if and only if its multiplicative group $(B,\circ)$ is nilpotent. 
\para

The papers \cite{MR1769723} and \cite{MR1809284} examined non-degenerate set-theoretic solutions to the Yang--Baxter equation that are not necessarily involutive. This naturally led to the search for an algebraic framework capable of encompassing such solutions, extending beyond the class of left braces, which yield involutive solutions. In this direction, Guarnieri and Vendramin introduced the concept of a skew left brace in \cite{MR3647970}. A \index{skew left brace}{skew left brace} $(B, \cdot, \circ)$ is a set $B$ with two operations $\cdot $ and $\circ$ such that both $(B,\cdot )$ and $(B, \circ)$ are groups and the condition
$$a \circ (b \cdot c)  = (a\circ b) \cdot a^{-1} \cdot (a \circ c)$$
holds for all $a,b,c \in B$, where $a^{-1}$ is the inverse of $a$ in $(B,\cdot)$. Note that, if $(B,\cdot)$ is an abelian group, then $(B, \cdot, \circ)$ is a left brace. It was shown in \cite{MR3647970} that every skew left brace gives rise to a non-degenerate bijective set-theoretic solution to the Yang--Baxter equation, and the solution is involutive if and only if $(B, \cdot, \circ)$ is a left brace. Further, an algorithm implemented both in GAP and Magma was provided to construct all non-isomorphic skew  left braces of small size. Skew braces have now become central to the classification of set-theoretic solutions to the Yang--Baxter equation, with a significant body of results consistently emerging. In \cite{MR4023387}, the structure of skew left braces was studied via their strong left ideals. A notable result was an analogue of It\^{o}'s theorem in the context of skew left braces. Specifically, if a skew left brace $(B, \cdot, \circ)$ factorizes through two strong left ideals that are trivial as skew left braces, then $(B, \cdot, \circ)$ is right nilpotent of class at most three. 
\para 

Another generalization of (skew) left braces, called \index{left cancellative left semi-brace}{left cancellative left semi-braces}, was introduced by Catino, Colazzo, and Stefanelli in \cite{MR3649817}. A left cancellative left semi-brace is a set $B$ with two operations $\cdot$ and $\circ$ such that $(B,\cdot)$ is a left cancellative semigroup and $(B,\circ)$ is a group such that
$$a \circ (b \cdot c) = (a \circ b) \cdot  \big(a \circ (\Bar{a} \cdot c) \big)$$
holds for all $a,b,c \in B$, where $\Bar{a}$ denotes the inverse of $a$ in $(B,\circ)$. Later, Jespers and Van Antwerpen introduced \index{left semi-brace}{left semi-braces} in \cite{MR3898225}, which generalize left cancellative left semi-braces by assuming that $(B, \circ)$ is just a semigroup. They also showed that, under some condition, these algebraic structures give left non-degenerate set-theoretic solutions to the Yang--Baxter equation. It is worth looking into the survey article \cite{MR3824447} by Ced\'{o} on braces and their connection with other algebraic structures, such as the holomorph of a group, a monoid of $I$-type, bijective 1-cocycles and Garside groups. The reader is also referred to the survey articles  \cite{MR3974481, MR4807456} by Vendramin, which focus on a range of open problems and conjectures in the field, with particular emphasis on skew left braces.
\para

Rota--Baxter operators on Lie algebras have been studied for a long time. In \cite{MR4271483}, Guo, Lang, and Sheng introduced Rota--Baxter operators on Lie groups, and showed that the inverse of the projection to the second factor in a direct product of  Lie groups gives rise to such operators. This construction leads to examples ranging from the Iwasawa decomposition to the Langlands decomposition. It was also shown in \cite[Theorem 2.9]{MR4271483} that the differentiation of a Rota--Baxter operator on a Lie group yields a Rota--Baxter operator on its Lie algebra.  In \cite{MR4370524}, Bardakov and Gubarev investigated the connection between skew left braces and Rota--Baxter groups. They proved that every Rota--Baxter group gives rise to a skew left brace, and conversely, every skew left brace can be embedded into a Rota--Baxter group. Moreover, when the additive group of a skew left brace is complete, then the brace structure is induced by a Rota--Baxter group. The findings establish a connection between Rota--Baxter groups and set-theoretic solutions to the Yang--Baxter equation.
\para

Another prominent class of solutions to the Yang--Baxter equation emerged in low-dimensional topology. A \index{quandle}{quandle} is a set equipped with a binary operation satisfying three axioms which can be seen as algebraic formulation of the three Reidemeister moves of planar diagrams of links in the Euclidean 3-space. More precisely, a quandle is a set $X$ with a right-invertible binary operation $(x,y) \mapsto x *y$ such that
$$x*x=x \quad \textrm{and} \quad (x*y)*z=(x*z) * (y*z)$$
for all $x, y, z \in X$. These structures first appeared in the fundamental work of Joyce \cite{MR2628474, MR0638121}, who introduced the term quandle, and independently in the work of Matveev \cite{MR0672410} under the name  distributive groupoid. 
They simultaneously proved that to each tame link $L$, one can associate a quandle $Q(L)$ that serves as a link invariant.  Furthermore, if $L_1$ and $L_2$ are two non-split tame links with isomorphic quandles $Q(L_1) \cong Q(L_2)$, then there is a homeomorphism of the 3-space mapping $L_1$ onto $L_2$, not necessarily preserving the orientation of the ambient space and the links.  By omitting the idempotency condition $x*x=x$, one obtains a weaker structure known as a \index{rack}{rack}.  If only the right-distributivity condition is required, the resulting structure is called a \index{shelf}{shelf}. These algebraic structures not only provide set-theoretic solutions to the Yang--Baxter equation, but also capture one of the essential properties of group conjugation. Over the years, racks, quandles and their generalizations have been extensively studied as tools for constructing invariants of links in classical, virtual, flat and other generalised knot theories.
\para

Racks constitute a special class of non-degenerate bijective set-theoretic solutions to the Yang--Baxter equation. Specifically, if $(X, *)$ is a rack, then the map $r: X \times X \to X \times X$ given by $$r(x, y)=(y,~ x*y)$$ gives a non-degenerate bijective set-theoretic solution $(X, r)$ to the Yang--Baxter equation.
Conversely,  as shown in \cite[Theorem 2.3]{MR1809284}, if $(X, r)$ is a non-degenerate set-theoretic solution to the Yang--Baxter equation with $r(x, y)=\big(\sigma_x(y), \,\tau_y(x)\big)$, then the binary operation $$x*y = \sigma_x \tau_{\sigma_y^{-1}(x)}(y)$$ defines a rack structure on $X$. 
\para

Historically, a special class of quandles, now known as  involutory quandles, first appeared in the work of Takasaki on finite geometry \cite{MR0021002}. Later, in an unpublished 1959 correspondence \cite{MR2885229}, Wraith and Conway discussed the idea of groups acting on themselves by conjugation, an idea that eventually gave rise to conjugation quandles. In \cite{MR0239005}, Loos introduced a differentiable reflection space as a differentiable manifold equipped with a differentiable binary operation satisfying the quandle axioms. For instance, if $G$ is a connected Lie group with an involutive automorphism $\varphi$ and $H$ a subgroup of the fixed-point subgroup containing the identity component, then the symmetric space $G/H$ acquires the structure of a differentiable reflection space with the binary operation $$Hx *Hy=H\varphi(x y^{-1}) y.$$ Moreover, any differentiable reflexion space can be obtained from a symmetric space \cite{MR0239005}. Results of this kind have led to transfer of ideas, for example the notion of flatness, from the theory of Riemannian symmetric spaces to the theory of racks and quandles \cite{MR3544543}.
\para

Interestingly, quandles arise naturally in a variety of mathematical contexts. Joyce \cite{MR0638121} related quandles to the theory of loops by observing that the core of a Moufang loop forms an involutory quandle.  Racks, a related structure, were connected to pointed Hopf algebras by Andruskiewitsch and Gra\~{n}a \cite{MR1994219}. A key step in the classification of finite-dimensional complex pointed Hopf algebras involves identifying all finite-dimensional Nichols algebras of braided vector spaces arising from groups, where a braided vector space is a vector space equipped with a braiding. In \cite{MR1994219}, Nichols algebras arising from braided vector spaces of the form $(\mathbb{C}Y,c)$ were studied, where $Y$ is a finite rack and $c$ is a 2-cocycle. Further, racks, quandles and  their cohomology theories have been related to  categorical groups and broader notions of categorification  \cite{MR2395367, MR2368886}. Another intriguing appearance of quandles is in the study of simple closed curves on closed orientable surfaces. Specifically,  the set of isotopy classes of such curves carries a natural quandle structure. From the perspective of the mapping class group, this quandle structure corresponds to the conjugation of Dehn twists. The idea was explored in detail  by Zablow in a series of papers \cite{MR2699808, MR1967241, MR2377276}. In a related direction, the papers \cite{MR1865704, Yetter2002, MR1985908} examined a similar quandle structure on the set of isotopy classes of simple closed arcs on an orientable surface with at least two punctures. For a broader overview of the historical development and connections to other areas of mathematics, we refer the reader to the survey articles \cite{MR2885229, MR2002606, MR2896084}.  For an accessible introduction focused on knot theory, see the text  \cite{MR3379534} by Elhamdadi and Nelson. For more advanced topics, especially those related to (co)homology, see the monograph \cite{MR3729413} by Nosaka.
\para

The past three decades saw many interesting works related to (co)homological aspects of quandles and racks, with a record of diverse applications to many areas of contemporary mathematics.  Fenn, Rourke, and Sanderson introduced a homology theory for racks in \cite{MR1194995, MR1257904, MR2063665, FennRourkeSanderson, MR2255194, MR1364012} as homology of the rack space of a rack, where a rack space is an appropriate analogue of the classifying space of a group. This motivated development of a (co)homology theory for racks and quandles in \cite{MR1990571} by Carter, Jelsovsky, Kamada, Langford, and Saito. They defined state-sum invariants using quandle cocycles as weights and used low-dimensional cocycles to construct invariants of knots and knotted surfaces. Later, in \cite{MR2128041}, Carter, Elhamdadi, and  Saito  defined a homology theory for set-theoretic solutions to the Yang--Baxter equation as a generalization of the homology theory developed in \cite{MR1990571}. These homology theories were further developed and generalised by Lebed and Vendramin in \cite{MR3558231}, which encompassed all the previously known homology theories. More precisely, the paper dealt with left non-degenerate set-theoretic solutions to the Yang--Baxter equation with a focus on racks and cycle sets. Two homology theories were considered which extended previous works known for biracks. An interpretation of the second cohomology group of left non-degenerate set-theoretic solutions was given in terms of group cohomology, which generalized a result of Etingof and Gra\~{n}a for racks \cite{MR1948837}.  There are still very few examples of racks whose homology is fully understood. In \cite{arXiv:2011.04524}, Lebed and Szymik computed the entire integral homology of all permutation racks using methods from homotopical algebra. The homotopical study of racks and quandles was carried further by Lawson and Szymik in \cite{arXiv:2106.01299}. They proved an analogue of Milnor's theorem on free groups by identifying homotopy types of free racks and free quandles on spaces of generators. As an application, they showed that the stable homotopy of the link quandle of a link is, in general, more complicated than what any Wirtinger presentation coming from a diagram predicts.
\para 

It is natural to consider extensions in order to produce more examples of racks and quandles. This was carried out by Andruskiewitsch and Gra\~{n}a in \cite{MR1994219}, who introduced modules for racks and 2-cocycles for racks with coefficients in a module. The 2-cocycles considered there fits into a general cohomology theory, which includes all other cohomology theories for racks found in the literature as special cases. Szymik showed in \cite{MR3937311} that quandle cohomology is a Quillen cohomology, which is the cohomology group of a functor from the category of models (or algebras) to that of complexes. In most algebraic theories, including that of groups and Lie algebras, 2-cocycles are known to be closely related to extensions of corresponding algebraic structures \cite{MR0672956}. Analogous theories were developed for quandle cocycles by Carter,  Elhamdadi, and Saito in \cite{MR1885217}. Besides being of independent algebraic interest, these proved useful in computing cocycle knot invariants. Abelian extensions of quandles were investigated using a diagrammatic approach by Carter, Kamada, and Saito in \cite{MR1973510}. Further generalizations of quandle extensions by dynamical cocycles were presented in \cite{MR1994219}, and a more general abelian extension theory for racks and quandles, which also encompassed other variants in the literature, was given by Jackson in  \cite{MR2155522}. Using the generalised quandle cohomology theory of \cite{MR1994219} some new knot invariants were introduced in \cite{MR2166720}. Quandle cocycle invariants were computed from quandle coloring numbers in \cite{MR3488311}, and several properties of invariants based on 2-cocycles were investigated in \cite{MR3582881}. For example, it was proved that if the quandle extension is a conjugation quandle, then the state-sum invariant given by the 2-cocycle is constant. An explicit description of the second quandle cohomology of a finite indecomposable quandle was given by Iglesias and Vendramin in \cite{MR3671570}. This was achieved by further developing ideas from \cite{MR1948837}, which related the second cohomology of a quandle to the first cohomology of its adjoint group.
\para

An explicit description of abelian group objects in the category of quandles and racks was given by Jackson in \cite{MR2155522}, leading to the construction of an abelian category of modules over these objects. Automorphisms of quandles have been investigated in much detail. In \cite{MR2175299}, Ho and Nelson determined automorphism groups of quandles of order less than 6. The investigation was carried forward by Elhamdadi, Macquarrie, and Restrepo in \cite{MR2900878}, where the automorphism group of the dihedral quandle $\R_n$ was shown to be isomorphic to the group of invertible affine transformations of $\mathbb{Z}_n$. In \cite{MR2831947}, Hou gave a precise description of automorphism groups of Alexander quandles. In \cite{MR3718201}, some structural results for the group of automorphisms and inner automorphisms of generalised Alexander quandles of finite abelian groups were obtained. The work was extended by Bardakov, Nasybullov, and Singh in \cite{MR3948284}, where several interesting subgroups of automorphism groups of conjugation quandles of groups were determined. Further, necessary and sufficient conditions were determined for these subgroups to coincide with the group of all automorphisms.
\para

While link quandles are recognized as robust invariants of links,  determining whether two quandles are isomorphic remains a challenging problem. 
This difficulty has motivated the study of additional structural properties of quandles, with particular emphasis on link quandles. It is well-known that residual finiteness and related residual properties play a crucial role in combinatorial group theory and low dimensional topology. 
The residual finiteness of link groups, and more generally, of 3-manifold groups, has been a longstanding topic of interest. Neuwirth \cite{MR0176462} showed that knot groups of fibered knots are residually finite. This result was extended by Mayland \cite{MR0295329} to the case of twist knots, and further generalized by Stebe  \cite{MR0237621} to include certain classes of non-fibered knots. Following Perelman's proof of the geometrization conjecture  \cite{Perelman1, Perelman2, Perelman3}, Hempel proved in \cite{MR0895623} that all finitely generated 3-manifold groups, in particular link groups, are residually finite. Motivated by these developments, in \cite{MR3981139}, Bardakov, Singh, and Singh initiated the systematic study of residual finiteness in the context of quandles, showing that free quandles and knot quandles possess this property. 
Building upon these results, the residual finiteness of link quandles was further established in  \cite{MR4075375}. As important consequences, it follows that link quandles are Hopfian and possess a solvable word problem.
\para

In an effort to incorporate ring-theoretic techniques into the study of quandles, Bardakov, Passi, and Singh in \cite{MR3977818} proposed a theory of quandle rings and algebras,  in analogy with the classical theory of group rings and algebras. This framework investigates the intricate relationships between quandles and their corresponding quandle rings and algebras, offering a new algebraic perspective on quandle structures. A notable feature of this theory is that quandle rings associated with non-trivial quandles are inherently non-associative. Moreover, it was shown in \cite{MR3915329} that these rings fail to be power-associative, placing them at the extreme end of the associativity spectrum. Further structural limitations were established in \cite{MR4450681}, where it was proved that quandle rings over base rings of characteristic greater than three cannot be alternative or Jordan algebras. Elhamdadi, Makhlouf, Silvestrov, and Zappala initiated a homological study of quandle rings in \cite{MR4443143}, providing a complete characterization of derivations for quandle algebras of dihedral quandles over fields of characteristic zero. Zero-divisors in quandle rings were investigated in \cite{MR4450681} using the idea of orderability of quandles. It was proved that quandle rings of left or right orderable quandles which are semi-latin have no zero-divisors. As a consequence, the integral quandle rings of free quandles are also free of zero-divisors. The orderability of quandles, including the knot quandles, was further explored in \cite{MR4330281}.  Units in group rings play a fundamental role in the structure theory of group rings. In contrast, it turns out that idempotents are the most natural objects in quandle rings, since each quandle element is, by definition, an idempotent of the quandle ring. In  \cite{MR4565221}, a detailed investigation of idempotents in quandle rings was carried out, with a focus on their precise computations. The concept of quandle coverings, as developed by Eisermann in \cite{MR1954330, MR3205568}, was utilized to determine the complete set of idempotents in the quandle ring of a finite-type quandle, which is a non-trivial covering over a nice base quandle. Additionally, idempotents in quandle rings of free products were also investigated, revealing an infinite family of quandles whose integral quandle rings possess only trivial idempotents.
\para

Over the past three decades, there has been an extraordinary amount of research focused on set-theoretic solutions to the Yang--Baxter equation. 
This monograph aims to provide a concise introduction to the algebraic framework underpinning these solutions. Our focus is centred on exploring fundamental algebraic structures such as skew braces, Rota--Baxter groups, quandles, racks, and their interconnectedness and applications to knot theory, while also shedding light on our own research contributions. It is important to acknowledge that this monograph is not intended to be a comprehensive reference for the continually expanding literature on these topics. Indeed, we believe that each of these topics merits a dedicated monograph to do full justice to their depth and complexity. Nonetheless, we have endeavoured to direct the reader to appropriate sources wherever a detailed treatment lies beyond the scope of this monograph. We have ensured that the proofs are provided in their entirety or accompanied by comprehensive references. One distinctive feature of this monograph is its progressive approach, commencing with basics and examples, enabling the reader to swiftly delve into meaningful material without requiring an extensive grasp of preparatory results. We have additionally posed several questions and problems for interested readers to explore. The target audience for this monograph includes graduate students and researchers engaged in the study of the Yang--Baxter equation, low dimensional topology, and related topics. We assume that the reader is  familiar with standard notions in algebra and topology, encompassing group theory, ring theory, homological algebra, algebraic topology, and knot theory.
\para

The monograph is structured into three main parts,  each addressing a distinct facet of the subject.  Part I is devoted to the algebraic theory of general set-theoretic solutions to the Yang--Baxter equation, with particular emphasis on skew left braces and Rota--Baxter groups. Part II presents a detailed treatment of the algebraic theory of racks and quandles. Part III, the most advanced of the three, is concerned with the homology and cohomology theories associated with solutions to the Yang--Baxter equation. From the point of view of logical dependency, Parts I and II are largely self-contained and may be read independently, while Part III builds upon foundational concepts introduced in the earlier parts.
\para

\bigskip
\bigskip

{\bf Acknowledgements.} It is appropriate to express our gratitude and extend acknowledgements to various individuals and institutions who have contributed to the development of this project. The initial seeds were planted during MS's visit to the Sobolev Institute of Mathematics in June 2019, and further progress was made during VB's visit to IISER Mohali in March 2020. Subsequently, from March to May 2022 and May to July 2023, MS dedicated time at the University of South Florida, with the visit being made possible through the generous support of the Fulbright-Nehru Academic and Professional Excellence Fellowship. Additionally, ME received support from the Simons Foundation Collaboration Grant 712462, while MS was fortunate to be supported by the SwarnaJayanti Fellowship grant of DST throughout the course of this project. We extend our sincere appreciation for the invaluable assistance and resources provided by these funding agencies and organizations. The authors are grateful to Tatyana Kozlovskaya for designing the figure on the cover page of the monograph. The authors thank their collaborators, students and friends, P. Belwal, O. V. Bryukhanov, N. Dhanwani, M. Kumar, S. V. Matveev, 
T. Nasybullov, M. V. Neshchadim, N. Rathee, H. Raundal, D. Saraf, M. Singh, D. V. Talalaev, and the late I. B. S. Passi for numerous mathematical discussions.
\para

\chapter*{Part I}

\begin{center}
\huge{General solutions, skew braces and Rota--Baxter groups} 
\end{center}

\chapter{Set-theoretic solutions to the Yang--Baxter equation}\label{chapter solution of YBE}

\begin{quote} 
In this chapter, we develop a comprehensive theory of set-theoretic solutions to the Yang--Baxter equation, building upon foundational contributions in the field. Our study includes numerous examples of such solutions and a thorough examination of the structure group associated with a solution, highlighting its key properties. Additionally, we explore braiding operators on groups and their significant implications. The chapter further presents results on indecomposable solutions, and investigates both linear and affine solutions within the framework of abelian groups.
\end{quote}
\bigskip

\section{Yang--Baxter equation and its solutions}

The theoretical physics and statistical mechanics witnessed the initial emergence of the quantum Yang--Baxter equation through the influential works of Yang \cite{MR0261870} and Baxter \cite{MR0290733, MR0690578}. Since then, this equation has sparked numerous intriguing applications in quantum groups, Hopf algebras, knot theory, tensor categories, and integrable systems, to name a few. 
\para 

A \index{monoidal category}{monoidal category} is a category $\mathcal{M}$ equipped with a functor $\otimes: \mathcal{M} \times \mathcal{M} \to \mathcal{M}$, called the {\it tensor product}, from the product category $\mathcal{M} \times \mathcal{M}$ to $\mathcal{M}$, such that $\mathcal{M}$ admits a unit object and the tensor product is associative. Given a monoidal category $\mathcal{M}$ and an object $V \in \mathcal{M}$, a quantum Yang--Baxter operator is a morphism $R: V \otimes V \to V \otimes V$ which satisfies the following \index{quantum Yang--Baxter equation}{\it quantum Yang--Baxter equation}
\begin{equation} \label{YBE}
R_{12} \,R_{13} \,R_{23} = R_{23} \,R_{13} \,R_{12}
\end{equation}
in $\End(V \otimes V \otimes V)$. Here, $R_{ij}$ denotes the morphism $V\otimes V \otimes V \to V\otimes V \otimes V$ that acts as $R$ on the $(i,j)$-th tensor factors and as the identity on the remaining third factor. The morphism $R$ is referred to as a solution to the quantum Yang--Baxter equation. When $\mathcal{M}$ is the category of vector spaces, then a solution to the quantum Yang--Baxter equation is also referred to as an \index{$R$-matrix} {\it $R$-matrix}.
\para 

Let $V$ be a vector space over a field. The linear map $\mathtt{t}: V \otimes V  \to V \otimes V$ given by $\mathtt{t}(v \otimes w)= w \otimes v$ is clearly a solution to the quantum Yang--Baxter equation. A direct check shows that a linear map $R: V \otimes V  \to V \otimes V$ is a solution to the quantum Yang--Baxter equation if and only if the linear map $S = \mathtt{t}~ R$ satisfies the following equation
\begin{equation} \label{BE}
S_{12} \,S_{23} \,S_{12} = S_{23} \,S_{12} \,S_{23}.
\end{equation}
In the literature, equation \eqref{BE} is commonly referred to as the \index{Yang--Baxter equation}{\it Yang--Baxter equation} or the \index{braid equation}{\it braid equation}. Further, $S$ is said to be a solution to the Yang--Baxter equation or the braid equation.  In this monograph, we will adopt the term ``Yang--Baxter equation'' for this  equation and focus on exploring set-theoretic solutions to it. Let $X$ be a basis of the vector space $V$. Then a map $r : X \times X \to X \times X$ satisfying the equation
\begin{equation}\label{set-theoretic solution of quantum YBE}
r_{12} \,r_{23}  \,r_{12} = r_{23}  \,r_{12} \, r_{23},
\end{equation}
where $r_{ij}: X \times X \times X \to X \times X \times X$ acts as $r$ on the $(i,j)$-th factors and as the identity on the remaining third factor, induces a solution to the Yang--Baxter equation on the vector space $V$. The pair $(X,r)$ is called a \index{set-theoretic solution}{\it set-theoretic solution} to the Yang--Baxter equation. Interest in the study of set-theoretic solutions to the quantum Yang--Baxter equation was intrigued by the paper \cite{MR1183474}  of Drinfel'd, published in 1992. Note that this is equivalent to studying set-theoretic solutions to the Yang--Baxter equation. However, such solutions had already appeared before Drinfel'd formulated his question. Notably, the work of Joyce \cite{MR2628474, MR0638121} and Matveev~\cite{MR0672410} introduced quandles as invariants of knots and links, which provide non-degenerate, bijective, set-theoretic solutions to the Yang--Baxter equation.
\para

In the interest of brevity, if $(X,r)$ is a set-theoretic solution to the Yang--Baxter equation, we commonly refer to it simply as a {\it solution to the Yang--Baxter equation} or a \index{braided pair}{\it braided pair}. Furthermore, the map $r$ is sometimes referred to as a \index{braiding }{\it braiding}.
\para

As the most elementary solution, if $X$ is any set, then the flip map $\mathtt{t}: X \times X \to X \times X$ given by $$\mathtt{t}(x, y)=(y, x)$$ is 
a solution to the Yang--Baxter equation, referred to as the \index{trivial solution}{\it trivial solution}. Now, let us explore some more intriguing examples.

\begin{example}\label{group commuting element example}
{\rm Let $G$ be a group and $a, b \in G$ two commuting elements. Then the map $r: G\times G \to G \times G$ given by
$$
r(x, y) = (ya, xb^{-1})
$$
is a solution to the Yang--Baxter equation.}
\end{example}

\begin{example}\label{self-dist example}
{\rm 
Let $X$ be a set with a binary operation $*$ such that $(x*y)*z=(x*z)*(y*z)$ holds for all $x, y, z \in X$. Then the map $r: X\times X \to X \times X$ given by
$$
r(x, y) = (y, x*y)
$$
is a solution to the Yang--Baxter equation.}
\end{example}

\begin{example} \label{YBEexamples} 
{\rm 
Let $\mathbb{k}$ be a commutative ring with unity $1$. Let $a,b \in \mathbb{k}$ be units  such that $(1-a)(1-b)=0$. 
Then a straightforward computation shows that if $r= \left[ \begin{array}{cc} 1-a & a \\ b & 1-b \end{array} \right]$, that is,
$$ r\big( (x, y) \big) = \big((1-a)x + ay, ~bx + (1-b)y \big)$$
for $x, y \in \mathbb{k},$ then $(\mathbb{k}, r)$ is a solution to the Yang--Baxter equation. 
\para
We now generalise the preceding example to matrices.  Let $r: \mathbb{k}^n \times \mathbb{k}^n \rightarrow  \mathbb{k}^n \times \mathbb{k}^n$ be defined by a matrix $\displaystyle r= \left[ \begin{array}{cc}I - Y & Y \\ Z & I-Z \end{array} \right]$, where $I$ is the identity matrix and $Y, Z$ are invertible  $n \times n$ matrices over $\mathbb{k}$ such that $YZ=ZY$ and $(I-Y)(I-Z)=0$. Then $(\mathbb{k}^n, r)$ is a solution to the Yang--Baxter equation. See \cite{MR2402508, MR2128040} for further matrix solutions to the Yang--Baxter equation.
} 
\end{example}
\para 

Let $X$ be a set and $r: X \times X \to X \times X$ a map. Then we can write $r$ in the form
$$
r(x, y) = \big(\sigma_x(y), \,\tau_y(x)\big),
$$
where $\sigma_x, \tau_y : X \to X$ are maps given by
$$
y \mapsto \sigma_x (y) \quad \textrm{and} \quad x \mapsto \tau_y(x)
$$
for $x, y \in X$.

\begin{proposition} \label{component conditions for a solution}
Let $X$ be a non-empty set and $r : X \times X \to X \times X$ a map given by
$$
r(x, y) = \big(\sigma_x(y), \,\tau_y(x)\big).
$$
Then $(X, r)$ is a solution to the Yang--Baxter equation if and only if the conditions
\begin{eqnarray}
\label{condition for solution 1} \sigma_{\sigma_x(y)} \big(\sigma_{\tau_y(x)} (z) \big) &=& \sigma_x \big( \sigma_y(z) \big),\\
\label{condition for solution 2} \tau_{\sigma_{\tau_y(x)}(z)}  \big( \sigma_x(y) \big) &=&\sigma_{\tau_{\sigma_y(z)}(x)}  \big( \tau_z(y) \big),\\
\label{condition for solution 3} \tau_{\tau_z(y)} \big(\tau_{\sigma_y(z)} (x) \big) &=& \tau_z \big( \tau_y(x) \big),
\end{eqnarray}
 hold for all $x, y, z \in X$.
\end{proposition}

\begin{proof}
Computing $r_{12}r_{23} r_{12}$, we obtain
\begin{eqnarray*}
r_{12}r_{23} r_{12}(x, y, z) &=& r_{12}r_{23} \big(\sigma_x(y), ~\tau_y(x),~ z \big)\\
&=& r_{12} \big(\sigma_x(y), ~\sigma_{\tau_y(x)}(z),~ \tau_z(\tau_y(x)) \big)\\
&=& \big(\sigma_{\sigma_x(y)}(\sigma_{\tau_y(x)}(z)), ~\tau_{\sigma_{\tau_y(x)}(z)}(\sigma_x(y)),~ \tau_z(\tau_y(x)) \big).
\end{eqnarray*}
Similarly, computing $r_{23} r_{12} r_{23}$ gives
\begin{eqnarray*}
r_{23} r_{12} r_{23}(x, y, z) &=& r_{23} r_{12} \big(x,~ \sigma_y(z),~ \tau_z(y) \big)\\
&=& r_{23} \big(\sigma_x(\sigma_y(z)), ~\tau_{\sigma_y(z)}(x), ~\tau_z(y)\big)\\
&=& \big(\sigma_x(\sigma_y(z)), ~\sigma_{\tau_{\sigma_y(z)}(x)}(\tau_z(y)), ~\tau_{\tau_z(y)}(\tau_{\sigma_y(z)}(x)) \big),
\end{eqnarray*}
which proves the assertion. $\blacksquare$
\end{proof}
\para

\begin{corollary}
Let $X$ be a non-empty set and $r : X \times X \to X \times X$ given by $r(x, y) = \big(\sigma_x(y), ~\tau_y(x) \big)$. Then the following  assertions hold:
\begin{enumerate}
\item If $\sigma_x = \id_X$ for all $x \in X$, then the pair $(X, r)$ is a solution to the Yang--Baxter equation if and only if 
$$
\tau_{\tau_z(y)} \big(\tau_{z} (x) \big) = \tau_z \big( \tau_y(x) \big)
$$
for all $x, y, z \in X$.
\item If $\tau_x = \id_X$ for all $x \in X$, then the pair $(X, r)$ is a solution to the Yang--Baxter equation if and only if
$$
\sigma_{\sigma_x(y)} \big(\sigma_{x} (z) \big) = \sigma_x \big( \sigma_y(z) \big)
$$
for all $x, y, z \in X$.
\end{enumerate}
\end{corollary}
\para

The subsequent classes of solutions hold specific significance within the subject.

\begin{definition}
A solution $(X, r)$ to the Yang--Baxter equation is said to be: 
\begin{enumerate}
\item bijective if $r : X\times X \to X\times X$ is a bijection;
\item \index{left  non-degenerate solution}{left  non-degenerate} if each $\tau_x$ is a bijection;
\item \index{right non-degenerate solution}{right non-degenerate} if each $\sigma_x$ is a bijection; 
\item \index{non-degenerate solution}{non-degenerate} if it is both left and right non-degenerate; 
\item \index{square-free solution}{square-free} if $r(x, x) = (x, x)$ for all $x \in X$; 
\item \index{involutive solution}{involutive}  if $r^2 = \id_{X\times X}$.
\item \index{finite solution}{finite}  if the set $X$ is finite.
\end{enumerate}
\end{definition}

An involutive solution $(X, r)$ to the Yang--Baxter equation is also called a \index{symmetric pair}{\it symmetric pair}.

\begin{example}\label{firs-examples}
{\rm Let us consider some examples:
\begin{enumerate}
\item The solution $(X,r)$, where $r(x, y)=(x, y)$, is bijective but neither left nor right non-degenerate.
\item The solutions in Examples \ref{group commuting element example} and \ref{self-dist example} are bijective as well as non-degenerate.
\item The solution $(X,r)$ with $r(x, y)=\big(y, ~y^{-1}x y\big)$, is bijective, square-free and non-degenerate.
\item The solution $(\mathbb{Z},r)$ with $r(a, b)=\big((-1)^ab, ~(-1)^ba\big)$, is involutive and non-degenerate.
\item Let $G$ be a non-trivial group with identity $1$. Then the solution $(X, r)$, where $r(x, y) = (1, x y)$, is non-bijective,  left  non-degenerate and  right degenerate.
\end{enumerate}
}
\end{example}

It is known that for finite left non-degenerate solutions, bijectivity and non-degeneracy are equivalent \cite[Theorem B]{MR4466104}.

\begin{theorem}
If $(X, r)$ is a finite left non-degenerate solution, then it is bijective if and only if it is right non-degenerate.
\end{theorem}

\begin{proposition}\label{involutive sol condition}
Let $(X, r)$ be a non-degenerate solution to the Yang--Baxter equation. Then $(X, r)$ is involutive if and only if
$$
\sigma_x(y) = \tau^{-1}_{\tau_y(x)}(y) \quad \textrm{and} \quad \tau_y(x) = \sigma^{-1}_{\sigma_x(y)}(x)
$$
for all $x, y \in X$. 
\end{proposition}

\begin{proof}
The proof follows from the identity $r^2 = \id_{X \times X}$ and  the non-degeneracy  of $\sigma_x$ and $\tau_y$.  
$\blacksquare$
\end{proof}

\begin{definition}
A solution $(X, r)$ to the Yang--Baxter equation with $r(x, y) = \big(\sigma_x(y),\, \tau_y(x)\big)$ is said to be an \index{elementary solution of first type}{elementary solution of first type} if $\tau_y = \id_X$ for all $y \in X$; and said to be an \index{elementary solution of second type}{elementary solution of second type} if $\sigma_x = \id_X$ for all $x \in X$.
\end{definition}

For instance, Example \ref{group commuting element example} gives an elementary solution of first type if we take $b=1$. Similarly, Example \ref{group commuting element example} and Example \ref{self-dist example} for $a=1$ give elementary solutions of second type.
\para

We now define a morphism between two solutions to the Yang--Baxter equation. 

\begin{definition}
Let $(X, r)$ and $(X', r')$ be solutions to the Yang--Baxter equation. A map $\varphi : X \to X'$ is called a \index{morphism of solutions}{morphism of solutions} if the following diagram commutes.
$$ \begin{CD}
X \times X @> r>> X \times X \\
@V{\varphi \times \varphi}VV @VV{\varphi \times \varphi}V \\
X' \times X' @> r' >> X' \times X'.
\end{CD}
$$
If $\varphi$ is bijective, then it is called an isomorphism of solutions.
\end{definition}

It is clear that the class of all set-theoretic solutions to the Yang--Baxter equation forms a category with the above definition of a morphism. Note that, if $r(x, y) = \big(\sigma_{x}(y), \,\tau_{y}(x)\big)$ and $r'(u, v) = \big(\sigma'_{u}(v), \,\tau'_{v}(u)\big)$, then 
$\varphi$ is a morphism of solutions if and only if 
$$\varphi \big(\sigma_{x}(y) \big) = \sigma'_{\varphi(x)} \big(\varphi(y) \big)\quad \textrm{and} \quad 
\varphi \big(\tau_{y}(x) \big) = \tau'_{\varphi(y)} \big(\varphi(x) \big)$$ for all $x, y \in X$.

\begin{definition}\label{conjugate sol definition}
Two solutions $(X,r)$ and $(X,r')$ to the Yang--Baxter equation are said to be \index{conjugate solutions}{conjugate} if there exists a bijection 
$T: X \times X \to X \times X$ such that $r' = T r T^{-1}$.
\end{definition}

Let $\mathtt{t}:X \times X \to X \times X$ be the flip map given by $\mathtt{t}(x, y)=(y, x)$. Then $(X, r)$ is an elementary solution of second type if and only if $ \mathtt{t} r \mathtt{t}$ is an elementary  solution of first type. Hence, there exists a one-to-one correspondence between elementary solutions of the first and second types.
\para

A set endowed with a single binary operation is referred to as a \index{groupoid}{\it groupoid}, while one equipped with two binary operations is termed a \index{bi-groupoid}{\it bi-groupoid}. There is an alternate approach to describe solutions to the  Yang--Baxter equation using bi-groupoids. A solution $(X, r)$ to the Yang--Baxter equation  defines a bi-groupoid $(X, \cdot, *)$ with binary operations
$$
y \cdot x =\sigma_x (y) \quad \textrm{and} \quad x * y = \tau_y(x)
$$
for $x, y \in X$. Thus, Proposition \ref{component conditions for a solution} takes the following form.

\begin{proposition} \label{f1}
Let $(X, \cdot, *)$ be a bi-groupoid  and $r : X \times X \to X \times X$ given by
$$
r(x, y) = (y \cdot x, ~x * y)
$$
for $x, y \in X$. Then, $(X, r)$ is a solution to the  Yang--Baxter equation if and only if the conditions
\begin{eqnarray*}
\big(z \cdot (x * y) \big)  \cdot (y \cdot x)&=&(z \cdot y) \cdot x,\\
(y\cdot x) * \big(z \cdot (x*y) \big)&=&(y*z) \cdot \big(x*(z \cdot y) \big),\\
(x*y)*z&= &\big(x*(z \cdot y) \big) *(y*z),
\end{eqnarray*}
hold for all $x, y, z \in X$.
\end{proposition}

In the general case, it is challenging to make definitive statements about bi-groupoids that satisfy these axioms. However, when we focus on elementary solutions, it is sufficient to examine certain groupoids.

\begin{corollary}
Let  $(X, \cdot, *)$ be a bi-groupoid such that $y \cdot x = y$ for all $x, y \in X$. Then the groupoid  $(X, *)$ defines a solution $r(x, y) = (y, x * y)$  if and only if the condition
\begin{equation}\label{right distributive condition}
(x * y) * z = (x*z) * (y * z)
\end{equation}
hold for all $x, y, z \in X$.
\end{corollary}

\begin{remark}
{\rm 
A groupoid satisfying condition \eqref{right distributive condition} is called a \index{distributive groupoid}{\it (right)~distributive groupoid} or a \index{shelf}{\it shelf}. A distributive groupoid $(X, *)$ in which the map $x \mapsto x*y$ is a bijection for each $y \in X$ is called a \index{rack}{\it rack}. If, in addition, $x*x=x$ for all $x \in X$, then $(X, *)$ is called a \index{quandle}{\it quandle}. In Part-II of the monograph, we will delve deeper into the intricate details of racks and quandles.}
\end{remark}

A non-degenerate elementary solution $(X,r)$ of the second type  with $r(x,y)=\big(y,\tau_y(x) \big)$ is said to be of \index{rack type solution}{\it rack type} if the groupoid $(X, *)$ with $x*y = \tau_y(x)$ is a rack. If $(X,*)$ is a quandle, then we say that the solution is of \index{quandle type solution}{\it quandle type}. The following result is immediate from the definitions.

\begin{proposition}
Let $X$ be a set and $r: X\times X\to X\times X$ a map given by $r(x,y) = \big(y, \tau_y(x)\big)$. Then the following  assertions hold:
\begin{enumerate}
\item $(X, r)$ is a~non-degenerate solution to the Yang--Baxter equation if and only if $(X,*)$ with $x*y = \tau_y(x)$ is a rack.
\item $(X, r)$ is a~non-degenerate square-free solution to the Yang--Baxter equation if and only if $(X,*)$ with $x*y = \tau_y(x)$ is a  quandle.
\end{enumerate}
\end{proposition}

In \cite[Theorem 2.3(iv)]{MR1809284}, Soloviev proved that every non-degenerate solution is conjugate to an elementary solution. This implies that instead of studying bi-groupoids corresponding to solutions, one can focus on studying groupoids.

\begin{proposition}\label{SolTh}
Let $(X, r)$ with $r(x,y) = \big(\sigma_x(y), \tau_y(x)\big)$ be a right non-degenerate solution to the Yang--Baxter equation. Then it is conjugate to the solution $(X, r^\star)$ with
$$
r^\star(x,y) = \left(y,~ \sigma_y(\tau_{\sigma^{-1}_x(y)}(x))\right).
$$
Further, the solution $(X,r^\star)$ is non-degenerate if and only if, for all $a, b \in X$, there exists a unique $x \in X$ such that
\begin{equation}\label{cond}
\sigma_a \big(\tau_{\sigma^{-1}_x(a)}(x) \big) = b.
\end{equation}
\end{proposition}

\begin{proof}
Since $(X, r)$ is right non-degenerate, $\sigma_x$ is bijective for each $x \in X$. Let $T_2: X \times X \to X \times X$ be the map defined by  $T_2(x, y) = \big(x, ~\sigma_x(y)\big)$ for $x, y \in X$. Then $T_2$ is invertible with $T_2^{-1}(x, ~y) = \big(x, ~\sigma^{-1}_x(y)\big)$ and 
$$
T_2  \,r  \,T_2^{-1}(x,~ y)
 = T_2 \,r \big(x, ~\sigma^{-1}_x(y) \big)
 = T_2 \big( \sigma_x(\sigma^{-1}_x(y)),~\tau_{\sigma^{-1}_x(y)}(x)\big)
 = \big(y, ~\sigma_y(\tau_{\sigma^{-1}_x(y)}(x)) \big)
 = r^\star(x,~y).
$$
Let us check  that $(X,r^\star)$ is a solution, that is,
$$
r_{12}^\star  \,r_{23}^\star  \,r_{12}^\star = r_{23}^\star  \,r_{12}^\star  \,r_{23}^\star
$$
where
$$
r_{12}^\star(x, y, z) = \left(y, ~\sigma_y(\tau_{\sigma^{-1}_x(y)}(x)), z\right) \quad \textrm{and}\quad r_{23}^\star(x, y, z) = \left(x,~z, ~\sigma_z(\tau_{\sigma^{-1}_y(z)}(y)) \right).
$$
Let $T_3 : X^3 \to X^3$ be given by 
$$T_3(x, y, z) = \big(x, ~\sigma_x(y),~ \sigma_x(\sigma_y(z)\big)$$
for $x, y, z \in X$. Then, its inverse is given by
$$
T_3^{-1} (x, y, z) = \left(x, ~\sigma^{-1}_x(y), ~\sigma^{-1}_{\sigma^{-1}_x(y)} (\sigma^{-1}_x(z)) \right).
$$
We claim that
$$
T_3  \,r_{12} \, T_3^{-1} = r_{12}^\star \quad \textrm{and}\quad T_3 \, r_{23}  \,T_3^{-1} = r_{23}^\star. 
$$
Since $(X, r)$ is a solution, this implies that $(X, r^\star)$ is also a solution.  We have 
$$
T_3  \,r_{12} (x, y, z) = T_3 \big(\sigma_x(y), \tau_y(x), z \big) =  \big(\sigma_x(y), ~ \sigma_{\sigma_x(y)}(\tau_y(x)), ~ \sigma_{\sigma_x(y)}(\sigma_{\tau_y(x)}(z)) \big)
$$
and
$$
 r_{12}^\star   \,T_3 (x, y, z) =  r_{12}^\star \big(x, \sigma_x(y), \sigma_x(\sigma_y(z))\big) =  \big(\sigma_x(y), ~ \sigma_{\sigma_x(y)}(\tau_y(x)), ~ \sigma_x (\sigma_y(z)) \big).
$$
By \eqref{condition for solution 1}, we have
$$ 
\sigma_{\sigma_x(y)} \big(\sigma_{\tau_y(x)} (z) \big) =  \sigma_x \big( \sigma_y(z) \big)
$$
for all $x, y, z \in X$, and hence the equality $T_3 r_{12} T_3^{-1} = r_{12}^\star$ holds. The proof of the equality $T_3 r_{23} T_3^{-1} = r_{23}^\star$ is similar, which establishes our claim. Further, the solution $(X, r^\star)$ is non-degenerate   if and only if, for every $a, b \in X$, there exists a~unique $x \in X$ such that
$\sigma_a \big(\tau_{\sigma^{-1}_x(a)}(x)\big) = b$. This completes the proof of the proposition. $\blacksquare$
\end{proof}

Next, we consider some more examples of solutions. 

\begin{example}\cite[Section 9]{MR1183474} 
{\rm   Let $X$ be a non-empty set and $r: X \times X \to X \times X$ given by $r(x, y) = \big(f(y), ~g(x)\big)$ for some maps $f, g : X \to X.$ 
Then the following  assertions hold:
\begin{enumerate}
\item $(X, r)$  is a solution to the Yang--Baxter equation if and only if $fg = gf$.
\item $(X, r)$ is non-degenerate if and only if  $f, g$ are bijective.
\item $(X, r)$ is involutive if and only if $g = f^{-1}$.
\end{enumerate}
\par

A solution of the preceding type is called a \index{permutation solution}{\it permutation solution}. If $f, g$ are cyclic permutations, then we  say that $(X, r)$ is a \index{cyclic permutation solution}{\it cyclic  permutation solution}. It is clear that two permutation solutions are isomorphic as solutions if and only if the corresponding permutations are conjugate.}
\end{example}

\begin{example}   \label{DirProd}
{\rm
If $(X, r)$ and $(Y, s)$ are two solutions, then  $(X \times Y, r  \times s)$  is a solution, called the \index{product solution}{\it product solution}.}
\end{example}

\begin{example} \label{OnGr}{\rm 
If  $G$ is a group, then we can define some elementary solutions on $G$ as follows:
\begin{enumerate}
\item $(x, y) \mapsto \big(y,~ y^{-k} x y^k \big)$ for $x, y \in G$ and fixed integer $k$.
\item $(x, y) \mapsto \big(y, ~y x^{-1} y \big)$ for $x, y \in G$.
\item $(x, y) \mapsto \big(y, ~\varphi(x y^{-1}) y \big)$ for $x, y \in G$ and fixed automorphism $\varphi \in \Aut(G)$.
\end{enumerate}
\para
One can check that all the preceding solutions are non-degenerate, elementary and square-free. In fact, all these solutions are of quandle type.}
\end{example}

The next example gives a family of non-elementary and non-square-free solutions.

\begin{example} \label{NonEl}
{\rm 
Let $G$ be a group and $\varphi:G \to G$ a homomorphism. Then $(G, r)$ with
$$
r(x, y) = \big(x y \varphi(x)^{-1}, ~\varphi(x)\big)
$$
is a solution to the Yang--Baxter equation. Further, $(G, r)$ is non-degenerate if and only if $\varphi \in \Aut(G)$.}
\end{example}

The next example is analogous to the preceding example. 

\begin{example} \label{NonEl2}
{\rm 
Let $G$ be a group and $\varphi:G \to G$ a homomorphism. Then $(G, r)$ with
$$
r(x, y) = \big(\varphi(y),~\varphi(y)^{-1} \, x y \big)
$$
is a solution to the Yang--Baxter equation. If $\varphi$ is the identity automorphism, then $r(x, y) = (y,~ y^{-1} x y)$ (see Example \ref{firs-examples}(3)). If $\varphi$ is the trivial homomorphism, then $r(x, y) = (1,~ x y)$ (see Example \ref{firs-examples}(5)). }
\end{example}

\begin{example}\cite[Example 3, p. 203]{MR1722951}
{\rm 
 Let $A$ be a ring with unity $1$ (not necessarily commutative) and $c\in A$ such that $1 + cx$ is invertible for any $x \in  A$. Define a new binary operation on $A$ by 
$$
x \circ y = x+y+xcy. 
$$
It is easy to check that the operation is associative and the zero element 0 of the ring $A$ acts as the identity element with respect to $\circ$. In fact, $(A, \circ)$ is a group, since for any $x \in A$, there exists
$$
x^{-1} = -x(1 + cx)^{-1}
$$
such that
$x\circ x^{-1} = x^{-1} \circ x = 0$.  Taking
$$
r(x, y) = \big(y(1 + cx + cxcy)^{-1}, ~x(1 + cy) \big),
$$
we see that $(A, r)$ is  a non-degenerate involutive solution to the Yang--Baxter equation. Further, $r$ defines a non-degenerate involutive solution on any right ideal of $A$.}
\end{example}

As we bring this section to a close, we examine how a solution to the Yang--Baxter equation  gives rise to  representations of braid groups. 

\begin{definition}
 The \index{braid group}{\it braid group} $\B_n$ on $n \ge 2$ strands is defined as the group generated by $n-1$ elements $\sigma_1,\ldots,\sigma_{n-1}$ with defining relations
\begin{align}
 \label{eq1}\sigma_i\sigma_{i+1}\sigma_i& = \sigma_{i+1}\sigma_i\sigma_{i+1} \quad \textrm{for}\quad 1 \le i \le n-2,\\
\label{eq2}\sigma_i\sigma_j &= \sigma_j\sigma_i \quad \textrm{for}\quad |i-j|\geq 2. 
\end{align}
\end{definition}
One can see that the braid equation \eqref{BE} resembles relations in the braid group $\B_n$.  Let $(X, r)$ be a bijective solution to the Yang--Baxter equation. Then, for each $1 \le i \le n-1$, we can define permutations $r_i : X^n \to X^n$ by setting
$$
r_i = \id_{X^{i-1}} \times r \times \id_{X^{n-i-1}}.
$$

The proof of the following proposition is immediate. 

\begin{proposition}\label{action of bn and sn via solution}
Let $X$ be a non-empty set and $r: X \times X \to X \times X$ a bijection. Then the following assertions  hold:
\begin{enumerate}
\item The assignment $\sigma_i \mapsto r_i$, for $1 \le i \le n-1$, extends to a representation $\B_n \to \Sigma_{X^n}$ if and only if $(X, r)$ is a solution to the Yang--Baxter equation.
\item The assignment $\sigma_i \mapsto r_i$, for $1 \le i \le n-1$, extends to a representation $\Sigma_n \to \Sigma_{X^n}$ if and only if $(X, r)$ is an involutive solution to the Yang--Baxter equation.
\end{enumerate}
\end{proposition}

\begin{exercise}{\rm 
Determine the representation of $\B_n$ corresponding to the solution $(\mathbb{Z}, r)$ with $r(a, b) = (b+1, a-1)$.}
\end{exercise}
\para

If $(X, r)$ is a solution to the Yang--Baxter equation such that $r(x,y) = \big(\sigma_x(y), \tau_y(x) \big)$ is bijective,  then we have a representation $\B_n \to \Sigma_{X^n}$  given by $\sigma_i \mapsto r_i$. It follows from Proposition \ref {SolTh} that $(X, r^\star)$, where  $$r^\star(x,y) = \left(y,~ \sigma_y(\tau_{\sigma^{-1}_x(y)}(x))\right),$$ is also a solution to the Yang--Baxter equation. Thus, we have another representation  $\B_n \to  \Sigma_{X^n}$ given by $\sigma_i \mapsto r_i^\star$, where
$$
r_i^\star = \id_{X^{i-i}} \times r^\star \times \id_{X^{n-i-i}}
$$
for each $1 \le i \le n-1$. Using the idea as in the proof of Proposition \ref {SolTh}, we can prove that these two representations of the braid group are conjugate  \cite[Theorem 2.3(iv)]{MR1809284}.

\begin{proposition} \label{SolRack} 
For each integer $n \geq 2$, there exists a bijection $T_n \colon X^n \to X^n$ such that
$$
T_n  \,r_i  \,T_n^{-1} = r_i^\star
$$
for all $1 \le i \le n-1$.
\end{proposition}
\para

Next, we consider an example from \cite{MR1918864} (see also \cite[Chapter 8]{MR1988550}), for which the induced action of the braid group has a nice geometrical meaning.

\begin{example}{\rm  
Consider the non-degenerate solution $(\mathbb{R}_+^{2}, r)$ to the Yang--Baxter equation, where $r$ is given by
$$
r(a, b, c, d) = \left( a+ ab + bc,~ \frac{bcd}{bc+a(1+b)(1+d)}, ~\frac{acd}{a+c+ad}, ~\frac{bc+a(1+b)(1+d)}{c} \right).
$$
To be precise, we have
$$
\sigma_{(a, b)}(c, d) = \left( a+ ab + bc,~ \frac{bcd}{bc+a(1+b)(1+d)} \right) \quad \textrm{and}\quad 
\tau_{(c, d)}(a, b) = \left(  \frac{acd}{a+c+ad}, ~\frac{bc+a(1+b)(1+d)}{c} \right).
$$

Replacing the triplet of operations $(\cdot, /, +)$ by $(+, -, \max)$, we obtain the non-degenerate solution $(\mathbb{R}^{2}, s)$ to the Yang--Baxter equation, where $s$ is given by
\begin{eqnarray*}
s(a, b, c, d) &=& \Big( \max\{a, a + b, b + c \}, ~d - \max \big\{0, a - b - c + \max\{0, b\} + \max\{0, d\} \big\},\\ 
&&  a + c + d - \max\{a, c, a + d\},~ \max \big\{b, a - c + \max\{0, b\} + \max\{0, d\} \big\}\Big).
\end{eqnarray*}
The two solutions give actions $\rho^r$ and $\rho^s$ of the braid group $\B_n$ on $\mathbb{R}_+^{2n}$ and $\mathbb{R}^{2n}$, respectively. It turns out that these actions have nice geometric meanings. To understand that, we restrict the action $\rho^s$ of $\B_n$ on the integer lattice $\mathbb Z^{2n}$ in $\mathbb R^{2n}$.  For each $x \in \mathbb Z$, let $x^+ = \max \{ 0, x \}$ and $x^- = \min \{ x, 0 \} $. We define right-actions
$$
\sigma, \, \sigma^{-1} : \mathbb Z^4 \to \mathbb Z^4
$$
on $(a,b,c,d) \in \mathbb Z^{4}$ as follows:
\begin{equation}
(a, b, c, d) \cdot  \sigma  = \big(a + b^+ + (d^+ - e)^+,  \, d - e^+,  \, c + d^- + (b^- + e)^-,  \, b + e^+ \big) \label{action-by-sigma}
\end{equation}
and
\begin{equation}
(a,b,c,d) \cdot  \sigma^{-1} = \big(a - b^+ - (d^+ + f)^+,  \, d + f^-,  \, c - d^-  - (b^- - f)^-,  \, b - f^- \big),  \label{action-by-sigma-inverse}
\end{equation}
where
\begin{equation}
e = a - b^- - c + d^+ \quad \textrm{and} \quad  f = a + b^- - c - d^+ . \label{definition-e-f}
\end{equation}
Let $(a_1, b_1, \ldots, a_n, b_n) \in \mathbb Z^{2n}$ and $\sigma_i^{\varepsilon} \in \B_{n}$, where $1 \le i \le n-1$ and $\varepsilon \in \{1, -1 \}$. We define the right-action of $\sigma_i^{\varepsilon} \in \B_{n}$ by
\begin{equation}
(a_1, b_1, \ldots, a_n, b_n) \cdot \sigma_i^{\varepsilon} = (a_1', b_1', \ldots, a_n', b_n'),  \label{action-by-generator-sigma-1}
\end{equation}
where $a_k' = a_k$ and $b_k' = b_k$ for $k \neq i, i+1$, and
\begin{equation}
(a_i', b_i', a_{i+1}', b_{i+1}') = \begin{cases}
\begin{array}{ll}
(a_i, b_i, a_{i+1}, b_{i+1}) \cdot \sigma & \mbox{if } \, \varepsilon = 1, \cr
(a_i, b_i, a_{i+1}, b_{i+1}) \cdot \sigma^{-1} & \mbox{if } \, \varepsilon = -1. \cr
\end{array}
\end{cases} \label{action-by-generator-sigma-2}
\end{equation}
For a word $w \in \B_n$ in the alphabet $\{ \sigma_1^{\pm 1}, \ldots, \sigma_{n-1}^{\pm 1}  \}$, we define the action of $w$ by
\begin{equation}
(a_1, b_1, \ldots,  a_n, b_n) \cdot w = \begin{cases}
\begin{array}{ll}
(a_1, b_1, \ldots, a_n, b_n) & \mbox{if } \, w= 1, \cr
((a_1, b_1, \ldots,  a_n, b_n)  \cdot \sigma_i^{\varepsilon}) \cdot w' & \mbox{if } \, w = \sigma_i^{\varepsilon} w'. \cr
\end{array}
\end{cases}
\end{equation}
It can be checked that the above action by $\B_n$ on $\mathbb Z^{2n}$ is well-defined. By  the \index{Dynnikov coordinates}{Dynnikov coordinates} of a braid $w \in \B_n$, we mean the vector $(0,1, \ldots, 0, 1) \cdot w$. The utilization of Dynnikov coordinates comes into  play when discerning the \index{positive braid}{positivity} or \index{negative braid}{negativity} of a braid \cite[Proposition 2]{MR1918864}, which is intimately related to ordering of braids \cite{MR1214782}.}
\end{example}
\bigskip
\bigskip


\section{Structure group of solutions}
Following \cite{MR1722951}, we now associate a group to a solution to the Yang--Baxter equation. The investigation of solutions will heavily rely on this group's pivotal role.

\begin{definition}
Let $(X,r)$ be a solution to the Yang--Baxter equation with $r(x, y)=\big(\sigma_x(y), \tau_y(x)\big)$. The \index{structure group}{\it  structure group} $G(X, r)$ of the solution $(X, r)$ is the group with the following presentation
$$
G(X, r)= \big\langle X  \, \mid \,  xy= \sigma_x(y) \tau_y(x)~\textrm{for all}~x, y \in X \big\rangle.
$$
\end{definition}
\para

\begin{example}{\rm 
Let $(X,r)$ be a solution to the Yang--Baxter equation.
\begin{enumerate}
\item If $r(x, y) = (x, y)$ is the identity map, then $G(X, r) \cong F(X)$, the free group with basis $X$.
\item If $r(x, y) = \mathtt{t}(x, y)=(y, x)$, then $G(X, r) \cong \mathbb{Z}^{X}$, the free abelian group of rank $|X|$.
\item If $X$ is a monoid with the identity element $e$, then $r(x, y) = (e, x y)$ is a solution and  $G(X, r) \cong F(X)$, the free group with basis $X$.
\end{enumerate}}
\end{example}

The following result gives a characterisation of a solution in terms of two compatible actions of the structure group on the underlying set of the solution \cite[Proposition 2.1]{MR1722951}.

\begin{proposition} \label{p122}
Let $X$ be a non-empty set and $r: X \times X \to X \times X$ a map given by $r(x, y)=\big(\sigma_x(y), \,\tau_y(x)\big)$ such that $\sigma_z$ and $ \tau_z$ are bijective for all $z\in X$. Then $(X, r)$  is a solution to the Yang--Baxter equation if and only if the following conditions are simultaneously satisfied:
\begin{enumerate}
\item The assignment $x \mapsto \tau_x$ is a right-action of $G(X, r)$ on  $X$.
\item The assignment $x \mapsto \sigma_x$ is a left-action of $G(X, r)$ on  $X$.
\item The linking relation
$$
\tau_{\sigma_{\tau_y(x)}(z)}  \big(\sigma_x(y)\big) = \sigma_{\tau_{\sigma_y(z)}(x)} \big(\tau_z(y)\big)
$$
holds for all $x, y, z \in X$.
\end{enumerate}
\end{proposition}

\begin{proof}
The assignment $x \mapsto \sigma_x$ defines a left-action of the group $G(X, r)$ on  $X$ if and only if 
$$\sigma_{\sigma_x(y)} \left(\sigma_{\tau_y(x)} (z) \right) = \sigma_x \big( \sigma_y(z) \big).$$
Similarly, $x \mapsto \tau_x$ defines a right-action of $G(X, r)$ on  $X$ if and only if
$$\tau_{\tau_z(y)} \big(\tau_{\sigma_y(z)} (x) \big) = \tau_z \big( \tau_y(x) \big).$$
Thus, by Proposition \ref{component conditions for a solution},  $(X, r)$  is a solution to the Yang--Baxter equation if and only if (1)-(3) hold simultaneously. $\blacksquare$
\end{proof}

Next, we have the following result \cite[Proposition 2.2]{MR1722951}. 

\begin{proposition} \label{p123}
The following  assertions hold:
\begin{enumerate}[(a)]
\item Let $X$ be a non-empty set and $r: X \times X \to X \times X$ with $r(x, y)=\big(\sigma_x(y), \,\tau_y(x)\big)$ an involution such that each $\tau_x$ is bijective and satisfy condition (1) in Proposition \ref{p122}. If $T : X \to X$ is the map given by $T(y) = \tau^{-1}_y(y)$, then 
\begin{equation}\label{equality of tau T and T sigma}
\tau_x^{-1}  \,T = T  \,\sigma_x.
\end{equation}
\item[] Suppose that, in addition, each $\sigma_x$ is bijective. Then we have the following:
\item The map $T$ is bijective. Consequently, the left-actions of $G(X, r)$ on $X$ given by $x \mapsto \tau^{-1}_x$ and $x \mapsto \sigma_x$  are isomorphic to each other.
\item  Condition (1) in Proposition \ref{p122} implies conditions (2) and (3) in Proposition \ref{p122}. Consequently, $(X, r)$ is a non-degenerate involutive solution if and only if the assignment $x \mapsto \tau_x$ is a right-action of $G(X, r)$ on $X$.
\end{enumerate}
\end{proposition}

\begin{proof}
Since $r$ is an involution, we have 
\begin{equation} \label{eqinvol}
\tau_{\tau_y(x)} \sigma_x(y) = y.
\end{equation} 
This together with condition (1) of Proposition \ref{p122} gives
$$
\tau_x^{-1} T(y) = \tau_x^{-1} \tau_y^{-1} (y) = \tau^{-1}_{\sigma_x(y)} \tau^{-1}_{\tau_y(x)} (y) = \tau^{-1}_{\sigma_x(y)} \sigma_{x} (y) = T \sigma_x(y),
$$
which is assertion (a).
\para 
Suppose that, in addition, each $\sigma_x$ is bijective. Taking $S(x) = \sigma_x^{-1}(x)$, we see that
\begin{equation}\label{equality of sigma S and S tau}
\sigma_x^{-1} S = S \tau_x.
\end{equation}
Now, using \eqref{equality of tau T and T sigma} and \eqref{equality of sigma S and S tau}, we obtain
$$TS(x)=T \sigma_x^{-1}(x)= \tau_x T (x)= \tau_x \tau^{-1}_x(x)=x=\sigma_x \sigma_x^{-1}(x)= \sigma_x S(x)=S\tau^{-1}_x(x) =ST(x),$$
which establishes assertion (b).
\para

The implication  $(1) \Rightarrow (2)$ follows from (a) and (b). Let us prove the implication $(1) \Rightarrow (3)$. Since $r$ is an involution, \eqref{eqinvol} holds and we can rewrite the linking relation in the form
$$
\tau_{\tau_{\tau_z \tau_y(x)}^{-1}(z)} \tau^{-1}_{\tau_y(x)}(y) = \tau^{-1}_{\tau_{\tau_z(y)} \tau_{\sigma_y(z)}(x)} \tau_z(y).
$$
Condition (1) of Proposition \ref{p122} gives $\tau_z \tau_y = \tau_{\tau_z(y)} \tau_{\sigma_y(z)}$. Setting $u = \tau_z \tau_y(x) = \tau_{\tau_z(y)} \tau_{\sigma_y(z)}(x)$, the preceding modified linking relation  takes the form
$$
\tau_{\tau_{u}^{-1}(z)} \tau^{-1}_{\tau^{-1}_z(u)} = \tau^{-1}_{u} \tau_{z},
$$
which proves assertion (c).
$\blacksquare$ 
\end{proof}
\para

Let $(X, r)$ be a non-degenerate bijective solution to the Yang--Baxter equation. We define another solution $(\overline{X}, \bar{r})$ as follows. For $x, y \in  X$, we write $x \sim y$ if $\sigma_x = \sigma_y$ and $\tau_x = \tau_y$. It is clear that $\sim$ is an equivalence relation on $X$. Let $\overline{X} = X /_\sim$ be the set of equivalence classes, and $\bar{x} \in \overline{X}$ denote the image of $x \in X$. A direct check yields the following result  \cite[Lemma 8.4]{MR3974961}.

\begin{proposition}
Let $(X, r)$ be a non-degenerate bijective solution to the Yang--Baxter equation. Then $(\overline{X}, \bar{r})$ is a non-degenerate bijective solution to the Yang--Baxter equation, where 
$$\bar{r}(\bar{x}, \bar{y}) = \big(\,\overline{\sigma_x(y)}, \, \overline{\tau_y(x)} \,\big)$$ 
for $\bar{x}, \bar{y} \in \overline{X}$. Further, the map $(X, r) \to (\overline{X}, \bar{r})$ given by $x\mapsto \bar{x}$ is a morphism of solutions.    
\end{proposition}

The solution $(\overline{X}, \bar{r})$ is called the \index{retraction}{\it retraction} of $(X, r)$, and is denoted by $\Ret(X, r)$. Setting $\Ret^1(X, r)=\Ret(X, r)$, we can repeat the process and define $\Ret^{n+1}(X, r) = \Ret \big(\Ret^{n}(X, r) \big)$ for each $n \ge 1$.
\para

\begin{remark}
{\rm The retraction of a solution was first defined for involutive solutions in \cite[Section 3.2]{MR1722951}. Note that, if $(X, r)$ is involutive, then by \eqref{equality of tau T and T sigma}, we have  $\tau_x = \tau_y$ iff $\sigma_x = \sigma_y$.}
\end{remark}
\para

\begin{definition}
A non-degenerate bijective solution  $(X, r)$ to the Yang--Baxter equation is said to be  a \index{multipermutation solution}{multipermutation solution} of level $n$ if $n$ is the smallest non-negative integer such that the underlying set of $\Ret^n(X, r)$ is a singleton set. Further,  $(X, r)$ is said to be  \index{irretractable solution}{irretractable} if $\Ret(X, r) = (X, r)$.
\end{definition}
\para

\begin{example}
A multipermutation solution of level zero is the trivial solution over a singleton set, while a multipermutation solution of level one is a permutation solution over a non-empty set. 
\end{example}

Several attempts have been made to classify solutions based on their multipermutation levels. Multipermutation solutions of level two can be characterized in a straightforward manner, thanks to results by Gateva-Ivanova  \cite[Proposition 4.7]{MR3861714}.

\begin{theorem}
Let $(X, r)$ be a non-degenerate involutive solution to the Yang--Baxter equation with $|X| \ge 2$. Then  $(X, r)$ is a multipermutation solution of level two if and only if 
$$\sigma_{\sigma_y(x)}=\sigma_{\sigma_z(x)}$$
for all $x, y, z \in X$.
\end{theorem}

Recall that, a group $G$ is called \index{poly-$\mathbb{Z}$ group}{\it poly-$\mathbb{Z}$} if it admits a subnormal series
$$1 = G_0  \triangleleft G_1 \triangleleft \cdots \triangleleft  G_n = G$$
such that $G_i/G_{i-1} \cong \mathbb{Z}$. Further, a group $G$ is \index{left-orderable group}{\it left-orderable} if there is a total order $<$ on $G$ such that 
$x < y$ implies $zx < zy$ for all $x, y, z \in G$. In general, finite non-degenerate involutive multipermutation solutions can be characterized in terms of their structure groups \cite[Theorem 2.1]{MR3815290}. 

\begin{theorem}
Let $(X, r)$ be a finite non-degenerate involutive solution to the Yang--Baxter equation. Then the following statements are equivalent:
\begin{enumerate}
\item $(X, r)$ is a multipermutation solution.
\item $G(X, r)$ is left-orderable.
\item $G(X, r)$ is poly-$\mathbb{Z}$.
\end{enumerate}
\end{theorem}

We note that the implication $(1) \Rightarrow (3)$ was already proved in  \cite{MR2189580}.

\begin{remark}
{\rm 
There is a continuously evolving  body of work on multipermutation solutions of finite levels:
\begin{enumerate}
\item In \cite[Theorem 4.21]{MR3861714}, several characterizations of solutions of finite multipermutation level are given.
\item In \cite{MR4123749}, a  classification of involutive solutions of multipermutation level two is given, with a bound on the number of such solutions. In the follow-up work \cite{MR4568107}, a complete characterization, up to isomorphism, of all indecomposable involutive solutions of multipermutation level two is given. 
\end{enumerate}}
\end{remark}
\para

Let $\Sigma_X$ be the group of permutations of $X$ and $\mathbb{Z}^X$ the free abelian group generated by $X$. We denote the generator of $\mathbb{Z}^X$ corresponding to $x \in X$ by $t_x$. Let $$M(X, r) = \Sigma_X \ltimes \mathbb{Z}^X$$  be the semi-direct product associated to the action of $\Sigma_X$ on $\mathbb{Z}^X$. For $\alpha \in \Sigma_X$ and $t_x \in \mathbb{Z}^X$, the relation $$\alpha t_x = t_{\alpha(x)} \alpha$$ holds in $M(X, r)$. With this set-up, we have the following result \cite[Proposition 2.3]{MR1722951}. 

\begin{proposition} \label{p129}
If $(X, r)$ is a non-degenerate involutive solution to the Yang--Baxter equation,  then the assignment $x \mapsto \tau^{-1}_{x} t_{x}$ extends to a group homomorphism $$\varphi_{\tau}: G(X, r) \to M(X, r).$$
\end{proposition}

\begin{proof}
The assignment $x \mapsto \tau^{-1}_{x} t_{x}$ gives a group homomorphism from the free group on $X$ to the group $M(X, r)$. For brevity, let us write   $r(x, y) = (u, v)$. Then we need to show that
 \begin{equation} \label{eq126}
\tau^{-1}_{x} t_{x} \tau^{-1}_{y} t_{y} = \tau^{-1}_{u} t_{u} \tau^{-1}_{v} t_{v}.
 \end{equation}
Using the relations in $M(X, r)$, we obtain
\begin{equation}\label{eq127}
\tau^{-1}_{x} t_{x} \tau^{-1}_{y} t_{y} = \tau^{-1}_{x} \tau^{-1}_{y}  t_{\tau_y(x)} t_{y}
\end{equation}
and
\begin{equation}\label{eq128}
\tau^{-1}_{u} t_{u} \tau^{-1}_{v} t_{v} = \tau^{-1}_{u} \tau^{-1}_{v}  t_{\tau_v(u)} t_{v}.
\end{equation}
Since $u = \sigma_x(y)$ and $v = \tau_y(x)$, the involutivity of $r$ gives $\tau_v(u) = y$. By Proposition \ref{p123}(3), the map $x \mapsto \tau_x$ is a right-action of $G(X, r)$ on $X$, which is equivalent to $\tau_y \tau_x = \tau_v \tau_u$. Using all this information on the right hand sides of \eqref{eq127} and \eqref{eq128} establishes \eqref{eq126}, which is desired. $\blacksquare$ 
\end{proof}

\begin{remark}{\rm 
If $(X, r)$ is a non-degenerate involutive solution to the Yang--Baxter equation, then we analogously have a homomorphism $$\varphi_{\sigma} : G(X, r) \to M(X, r)$$ defined by $x \mapsto  t_x \sigma_x$. Let $T:X \to X$ be the permutation given by $T(x)=\tau_x^{-1}(x)$ and $\hat{T}$ the automorphism of $M(X, r)$ induced by $T$. By Proposition \ref{p123}, we have $\tau_x^{-1} T = T \sigma_x$, which implies that $\hat{T}$ conjugates $\varphi_{\tau}$ to $\varphi_{\sigma}$. Hence, it is enough to study the properties of $\varphi_{\tau}$, and we denote it simply by $\varphi$ \cite[Section 2.2]{MR1722951}.}
\end{remark}

Let $\pi : G(X,r) \to \mathbb{Z}^X$  be the map defined by $\pi(g) = \pi_2\varphi(g)=t$, where $\varphi(g) = \alpha\, t$ for $\alpha \in \Sigma_X$ and $t \in \mathbb{Z}^X$, and $\pi_2: M(X, r)\to \mathbb{Z}^X$ is the projection onto the second factor. The assignment $x \mapsto \tau_x^{-1}$ gives a left-action of $G(X,r)$ on $X$, which further extends to a left $G(X,r)$-module structure on $\mathbb{Z}^X$. With this background, we have the following result \cite[Proposition 2.5]{MR1722951}.

\begin{proposition}  \label{p130}
Let $(X, r)$ be a non-degenerate involutive solution to the Yang--Baxter equation. Then the following  assertions hold:
\begin{enumerate}
\item $\pi : G(X,r) \to \mathbb{Z}^X$  is a group 1-cocycle of $G(X,r)$  with coefficients in the $ G(X,r)$-module $\mathbb{Z}^X$, that is,
$$\pi(g h) = \pi(g) + g \cdot \pi(h)$$
for any $g, h \in G$.
\item $\pi$ is bijective.
\end{enumerate}
\end{proposition}

\begin{proof}
Let $g, h \in G(X, r)$. Suppose that $\varphi(g)=\alpha \, s$ and $\varphi(h)=\beta \, t$ for $\alpha, \beta \in \Sigma_X$ and $s, t \in \mathbb{Z}^X$. Using the fact that $\varphi$ is a group homomorphism and the definition of a semi-direct product, we have 
$$\pi(gh)= \pi_2\varphi(gh)=\pi_2 \big(\varphi(g) \varphi(h) \big)= \pi_2 \big((\alpha \, s)(\beta \, t) \big)= \pi_2 \big((\alpha\beta) (s+ \alpha(t)) \big)= s+ \alpha(t)=\pi(g)+ g \cdot \pi(h),$$
which proves assertion (1). The proof of assertion (2) involves constructing an explicit inverse to the map $\pi$, and we refer the reader to \cite[Section 2.3]{MR1722951} for details. $\blacksquare$ 
\end{proof}

As an immediate consequence of Proposition \ref{p130}, we have the following result.

\begin{corollary}
Let $(X, r)$ be a non-degenerate involutive solution to the Yang--Baxter equation. Then the homomorphism $\varphi: G(X, r) \to M(X, r)$ is injective.
\end{corollary}
\medskip

\begin{remark}\label{remark on bij 1 cocycle}
{\rm 
Groups with bijective 1-cocycles have a nice geometrical interpretation. Namely, a bijective 1-cocycle on a simply connected Lie group $G$ with coefficients in some real representation is equivalent to a left-invariant affine structure on $G$ as a manifold, which identifies $G$ with an affine space. To admit such a structure, the group must be contractible as a topological space. It is known that many solvable Lie groups admit left-invariant affine structures as above \cite{MR1411303}. In fact, Milnor had conjectured that any solvable simply connected Lie group admit a left-invariant affine structure, but it was disproved by Burde \cite{MR1411303} by finding a nilpotent group that admit no such structure. Further, bijective cocycles on finite groups have been used by Etingof and Gelaki \cite{MR1653340} to construct new examples of semisimple Hopf algebras.}
\end{remark}

\begin{remark}\label{the g0 group}
{\rm 
Let $(X, r)$ be a non-degenerate involutive solution to the Yang--Baxter equation. Consider the group $$\Gamma =  G(X,r) \cap \mathbb{Z}^X,$$ where we consider $ G(X,r)$ as a subgroup of $M(X, r)$ via $\varphi$. Since $\Gamma$ is the kernel of the homomorphism $\rho :  G(X,r) \to \Sigma_X$, defined by the action of $ G(X,r)$ on $X$, it follows that $\Gamma$ is a normal subgroup of $ G(X,r)$. We set 
$$ G^0(X,r) =  G(X,r) / \Gamma \quad \textrm{and} \quad A = \mathbb{Z}^X / \Gamma.$$
Then we see that  $G^0(X,r)$ is the image of $G(X,r)$ in $\Sigma_X$. In particular, if the set $X$ is finite, then $G^0(X,r)$ is finite and $A$ is a finite abelian group with $|A| = |G^0(X,r)|$.}   
\end{remark}
\bigskip

Next, we associate another group $A(X, r)$ to each solution $(X, r)$ to the Yang--Baxter equation.

\begin{definition}
Let $(X,r)$ be a solution to the Yang--Baxter equation with $r(x, y)=\big(\sigma_x(y), \tau_y(x)\big)$. The \index{derived structure group}{\it derived structure group} $A(X, r)$ of the solution $(X, r)$ is a group with the following presentation
$$
A(X, r) = \big\langle X \, \mid \, x \; \sigma_x(y)  = \sigma_x(y) \; \sigma_{\sigma_x(y)}\left( \tau_y(x) \right)~ \textrm{for all}~ x, y \in X \big\rangle.
$$
\end{definition}
\para

\begin{example}\label{examples of derived str group}
{\rm 
Let $(X,r)$ be a solution to the Yang--Baxter equation.
\begin{enumerate}
\item If $r(x, y) = (x, y)$ is the identity map, then $A(X, r) \cong G(X, S) \cong F(X)$, the free group with basis $X$.
\item If $r(x, y) = (y, x)$, then $A(X, r) \cong G(X, r)\cong \mathbb{Z}^{X}$, the free abelian group generated by $X$.
\item If $(X, r)$ is an involutive solution and $r$ is not the identity map, then 
$$
r^2(x, y) = \big( \sigma_{\sigma_x(y)}\left( \tau_y(x) \right),  ~ \tau_{\tau_y(x)}\left( \sigma_x(y) \right) \big) = (x, y)
$$
and 
$$
A(X, r)= \big\langle X  \, \mid \, x \, \sigma_x(y)  = \sigma_x(y) \, x ~ \textrm{for all}~ x, y \in X \big\rangle.
$$
If $(X, r)$ is right non-degenerate, then  $A(X, r)\cong \mathbb{Z}^{X}$, the free abelian group of rank $|X|$.
\item If $X$ is a quandle and $r(x, y) = (y, \,x*^{-1}y)$, then $A(X, r) = G(X, r)$.
\end{enumerate}}
\end{example}
\para

Given a solution $(X, r)$ to the Yang--Baxter equation, let 
$$i_X : X \to G(X, r) \quad \textrm{and} \quad j_X : X \to A(X, r)$$ denote the two natural maps that send each element of $X$ to the corresponding generator of the group. We note that these maps are not always injective. 

\begin{example}{\rm 
Let $X = \mathbb{Z}$ and $r: X \times X \to X \times X$ be given by $r(i, j) = (j+1, i)$. Then $(X, r)$ is a non-degenerate bijective but non-involutive solution to the Yang--Baxter equation. Note that $G(X, r)$ is generated by elements $\{x_i \,\mid \; i \in \mathbb{Z} \}$
and has defining relations
$$
x _i \, x_j = x_{j+1}\,  x_i
$$
for $i,j \in \mathbb{Z}$. These relations imply  that $x_i = x_j$ for all $i, j \in \mathbb{Z}$, and hence $G(X, r) = \langle x_0 \rangle$
is the infinite cyclic group. Clearly, the natural map $i_X : X \to G(X, r)$ is not injective.}
\end{example}

In general, we have the following example (see \cite[Section 9]{MR1183474} and \cite[p. 171]{MR1722951}).

\begin{example}{\rm 
Let $X$ be a non-empty set, and $f$ and $g$ two commuting permutations of $X$. Then the map $$r(x,y) = \big(g(y),f(x) \big)$$ is a non-degenerate solution to the Yang--Baxter equation. In this case, we have the relation $xy = fg(x) \, fg(y)$ in $G(X,r)$ for all $x, y \in X$. Thus, if $fg$ has at least one fixed point in $X$, then we have relations $fg(y) = y$ in $G(X,r)$ for all $y \in X$, although it may not be so in $X$.}
\end{example}
\medskip

Example \ref{examples of derived str group}(4) supply infinitely many cases where $i_X$ fails to be injective, as well as infinitely many cases where $i_X$ is injective. See Theorem \ref{injectivity-dehn-eta} for a complete characterisation of quandles for which $i_X$ is injective. However, if the solution is involutive, then  $i_X: X \to G(X, r)$ is always injective (see Remark \ref{injectivity of eta remark}).
\medskip

\begin{example}{\rm 
Let $X$ be a monoid with the identity element $e$  and $r: X \times X \to X \times X$ be given by $r(x, y) = (e, x y)$. Then $(X, r)$ is a non-bijective, left  non-degenerate and right degenerate  solution to the Yang--Baxter equation. Note that $A(X, r)$ is the trivial group, and hence the natural map $j_X : X \to A(X, r)$ is not injective.}
\end{example}

For a non-degenerate solution $(X, r)$, let us write 
\begin{equation}\label{left-actions of the structure group}
x \circ y = \sigma_x(y) \quad \textrm{and} \quad  y * x = \tau_y^{-1}(x)
\end{equation}
for $x, y \in X$. By Proposition \ref{p122}, these operations can be extended to left-actions of the structure group $G(X, r)$ on $X$. We denote the action of an element $g \in  G(X,r)$ on an element $x \in  X$ by $g \circ x$ and $g * x$, respectively. 
Define the maps
$$
\phi: X \times X \to X \quad \textrm{by} \quad \phi(x, y)=y \circ \big( (x^{-1} \circ y)^{-1} * x \big)
$$
and
$$
r^\star: X \times X \to X \times X \quad \textrm{by} \quad r^\star(x, y)=\big(y, \phi(x, y)\big)
$$
for $x, y \in X$. Here, for $z \in X$, the element $z^{-1}$ is the inverse of $z$ when viewed as a generator of $G(X, r)$. Then we have the following result \cite[Theorem 2.3]{MR1809284}. 
\para

\begin{theorem} \label{deribed sol derived group}
Let $(X, r)$ be a non-degenerate solution to the Yang--Baxter equation. Then the following  assertions hold:
\begin{enumerate}
\item $\phi$ is $G(X, r)$-equivariant with respect to the $\circ$-action. In other words,
$$
g \circ \phi(y, z) = \phi(g \circ y, g \circ z)
$$
holds for all $g \in G(X, r)$ and $x, y \in X$.
\item  $(X, r^\star)$ is a non-degenerate solution to the Yang--Baxter equation (called the \index{derived solution}{\it derived solution} of $(X, r)$).
\item The structure group of the derived solution is the derived structure group.
\end{enumerate}
\end{theorem}

\begin{proof}
For assertion (1), it is enough to prove that
$$
t \circ \phi(y, z)  = \phi(t \circ y, t\circ z)
$$
for all $t, y, z \in X$. By definition of $\phi$, it is equivalent to proving that
\begin{equation} \label{eq-2.3}
t \circ \big( z \circ ( (y^{-1} \circ z)^{-1} * y ) \big)= 
(t \circ z) \circ  \big( ((t \circ y)^{-1} \circ (t \circ z) )^{-1} * (t \circ y) \big)
\end{equation}
for all $t, y, z \in X$. Since $(X, r)$ is a solution to the Yang--Baxter equation, the linking relation of Proposition \ref{p122} holds, that is, 
\begin{equation} \label{eq-2.4}
\big((y^{-1} * x) \circ t \big)^{-1} * (x \circ y) =  \big(( y \circ t)^{-1} * x \big) \circ (t^{-1} * y)
\end{equation}
for all $x, y, t \in X$. Putting $t = y^{-1} \circ z$ on the left hand side of  \eqref{eq-2.4} gives
$$
 \big((y^{-1} * x) \circ t \big)^{-1} * (x \circ y) =  \big((y^{-1} * x) \circ (y^{-1}x^{-1} \circ (x \circ z)) \big)^{-1} * (x \circ y),
$$
where $y^{-1}x^{-1}$ is the product  in the group $G(X, r)$. Since $$y^{-1}x^{-1} = \tau_y(x)^{-1}\sigma_x(y)^{-1}=(y^{-1} *x)^{-1} (x \circ y)^{-1}$$ in  $G(X, r)$, we have
$$
 \big((y^{-1} * x) \circ t \big)^{-1} * (x \circ y) =  \big((y^{-1} * x) \circ ((y^{-1} *x)^{-1} (x \circ y)^{-1} \circ (x \circ z)) \big)^{-1} * (x \circ y)= \big((x \circ y)^{-1} \circ (x \circ z)\big)^{-1} * (x \circ y).
$$
Hence, relation \eqref{eq-2.4} is equivalent to the relation
\begin{equation} \label{eq-2.5}
\big((x \circ y)^{-1} \circ (x \circ z) \big)^{-1} * (x \circ y) =(z^{-1} * x) \circ \big((y^{-1} \circ z)^{-1} * y \big)
\end{equation}
Using \eqref{eq-2.5} in the right hand side term of \eqref{eq-2.3} gives 
$$(t \circ z) \circ  \big( ((t \circ y)^{-1} \circ (t \circ z) )^{-1} * (t \circ y) \big)=(t \circ z) \circ \big( (z^{-1} * t) \circ ((y^{-1} \circ z)^{-1} * y)\big)=\big((t \circ z)(z^{-1} * t)\big) \circ \big((y^{-1} \circ z)^{-1} * y\big).$$
Using the relation
$$
tz = \sigma_t(z)\tau_z(t)= (t \circ z)(z^{-1} * t)
$$
in group $G(X,  r)$,  the preceding equation can be written as 
$$(t \circ z) \circ  \big( ((t \circ y)^{-1} \circ (t \circ z) )^{-1} * (t \circ y) \big)=t \circ \big( z \circ ((y^{-1} \circ z)^{-1} * y) \big),$$
which is desired.
\para 

Since $\phi(x, y)=\sigma_y \big(\tau_{\sigma^{-1}_x(y)}(x) \big)$, assertion (2) follows from Theorem  \ref{SolTh}. Finally, assertion (3) follows from the definitions of the structure group, the derived structure group and the map $\phi$. $\blacksquare$
\end{proof}

Let $\Aut_X \big(A(X, r) \big)$ denote the group of automorphisms of $A(X, r)$ that map the generating set $X$ onto itself. The next result shows that the natural map $j_X:X \to A(X, r)$ is $G(X, r)$-equivariant with respect to a suitable $G(X, r)$-action on $A(X, r)$ \cite[Theorem 2.4]{MR1809284}.

\begin{theorem} \label{action of str group on der str group}
Let $(X, r)$ be a non-degenerate solution to the Yang--Baxter equation. Then the group homomorphism  $\circ : G(X, r) \to \Sigma_X$  can be uniquely lifted to the group homomorphism  $\hat{\circ}: G(X, r) \to  \Aut_X \big(A(X, r)\big)$ such that $j_X:X \to A(X, r)$ is $G(X, r)$-equivariant. More precisely,
$$j_X(g \circ x) = g \,\hat{\circ}\, j_X(x)$$
for all  $g \in G(X, r)$ and $x \in X$.
\end{theorem}

\begin{proof}
It is clear that the homomorphism $\circ: G(X, r) \to \Sigma_X$ can be uniquely extended to a homomorphism from $G(X,r)$ to the group of automorphisms of the free group $F(X)$ on $X$. By Theorem \ref{deribed sol derived group}, each such automorphism of $F(X)$ arising from the $G(X, r)$-action  preserves relations in $A(X, r)$, and hence gives a homomorphism $\hat{\circ}: G(X, r) \to  \Aut_X \big(A(X, r)\big)$. The last assertion is immediate from the construction.  $\blacksquare$
\end{proof}

We have the following results (see \cite[Theorem 2.15]{MR1722951} and \cite[Proposition 6 and Proposition 10]{MR1769723}).

\begin{theorem} \label{CentSub} 
Let $(X, r)$ be a finite solution to the Yang--Baxter equation. Then the following assertions hold:
\begin{enumerate}
\item If  $(X, r)$ is non-degenerate or bijective, then $G(X, r)$ has a finitely generated abelian normal subgroup of finite index.
\item If $(X, r)$ is non-degenerate involutive, then $G(X, r)$ is solvable.
\end{enumerate}
\end{theorem}

Proposition \ref{deribed sol derived group}(3) and Theorem \ref{CentSub}(1) give the following result.

\begin{corollary}  \label{CentSubDer}
Let $(X, r)$ be a finite non-degenerate solution to the Yang--Baxter equation. Then $A(X, r)$ has a finitely generated abelian normal
subgroup of finite index. Moreover, $A(X,r)$ has a central subgroup of finite index.
\end{corollary} 

The next result gives equivalent conditions for the group $A(X, r)$ to be abelian (see \cite[Proposition 3.15]{MR3835326} and \cite[Proposition 4]{MR1769723}).

\begin{proposition} \label{AbelGr}
Let $(X, r)$ be a non-degenerate solution to the Yang--Baxter equation. Then the following statements are equivalent:
\begin{enumerate}
\item $A(X, r)$ is abelian.
\item $A(X, r)$ is free abelian with basis $X$.
\item $(X,r)$ is involutive.
\end{enumerate}
\end{proposition}

\begin{proof}
The implication $(2) \Rightarrow  (1)$ is obvious. Recall the presentation
$$A(X, r) = \big\langle X ~\mid ~x \; \sigma_x(y)  = \sigma_x(y) \; \sigma_{\sigma_x(y)}\left( \tau_y(x) \right)~ \textrm{for all}~ x, y \in X \big\rangle.$$
Further, note that $(X,r)$ is involutive if and only if 
$$\left( \sigma_{\sigma_x(y)}\left( \tau_y(x) \right),  ~ \tau_{\tau_y(x)}\left( \sigma_x(y) \right) \right) = (x, y)$$
for all $x, y \in X$. Since $(X, r)$ is a non-degenerate solution, the implication $(1) \Leftrightarrow (3)$ is now evident. Let us prove the implication $(1)  \Rightarrow  (2)$. Consider the map
$$
X \to H:=\big\langle j_X(X)~|~j_X(x) \, j_X(y) = j_X(y) \, j_X(x) ~\textrm{for all}~x, y \in X \big\rangle  \cong \mathbb{Z}^{j_X(X)}
$$
given by $x \mapsto j_X(x)$. We claim that it induces a group homomorphism $A(X, r) \to  H$. Since $A(X, r)$ is abelian, $(X,r)$ is involutive, which implies that
$$
j_X \left( \sigma_{\sigma_x(y)} (\tau_y(x))\right) = j_X(x)
$$
for all $x, y \in X$. Thus, if we apply $j_X$ to the  relation
$$
\sigma_x(y)^{-1} \, x^{-1} \, \sigma_x(y) \, \left( \sigma_{\sigma_x(y)} (\tau_y(x))\right) 
$$
of $A(X, r)$, then we obtain
$$
j_X \left( \sigma_x(y) \right)^{-1} \, j_X(x)^{-1} \, j_X(\sigma_x(y)) \, j_X\left( \sigma_{\sigma_x(y)} (\tau_y(x))\right) =
j_X \left( \sigma_x(y) \right)^{-1} \, j_X(x)^{-1} \, j_X(\sigma_x(y)) \, j_X (x) = 1.
$$
The homomorphism is in fact an isomorphism with its  inverse induced by the inclusion map $j_X(X) \to  A(X, r)$. $\blacksquare$ 
\end{proof}

For a finite solution $(X, r)$,  it is possible to characterize precisely when the group $A(X, r)$ is torsion-free \cite[Proposition 3.16]{MR3835326}.

\begin{proposition}
Let $(X, r)$ be a finite non-degenerate solution to the Yang--Baxter equation. Then $A(X, r)$ is torsion free if and only if $A(X, r)$ is abelian.
\end{proposition}

\begin{proof}
By Proposition \ref{AbelGr}, it is enough to show that $A(X, r)$ is torsion-free if and only if it is free abelian. Since $(X, r)$ is finite, by Corollary \ref{CentSubDer}, $A(X, r)$ is a finitely generated group containing an abelian normal subgroup of finite index. By \index{Schreier's lemma}{Schreier's lemma}, all finite index subgroups of finitely generated groups are also finitely generated \cite[Proposition 4.2]{MR1812024}. If $A(X, r)$ is torsion-free, we conclude that it contains a central finitely generated free abelian group of finite index. Hence, by \cite[Proposition 7.13]{MR1812024}, $A(X, r)$ is itself free abelian. $\blacksquare$
\end{proof}
\para

A group $G$ is called a \index{Bieberbach group}{\it Bieberbach group} if it is finitely generated, torsion free and abelian-by-finite. In other words, $G$ is an extension of a finitely generated free abelian group by a finite group. We conclude this section with the following result (see \cite[Theorem 1.6]{MR1637256} and \cite[Corollary 8.2.7]{MR2301033}).

\begin{theorem}
If $(X, r)$ is an involutive solution, then $G(X, r)$ is a Bieberbach group and is of $I$-type. 
\end{theorem}
\bigskip
\bigskip


\section{Solutions and braiding operators on groups}

In this section, our focus lies on exploring solutions to the Yang--Baxter equation within the context of groups. We begin with the following result \cite[Theorem 1]{MR1769723}.
      
\begin{theorem} \label{CompAct}
Let $G$ be a group. Let $\sigma$ and $\tau$ be left- and right-actions of $G$ on itself, denoted by $(x, y) \mapsto \sigma_x(y)$ and $(x, y) \mapsto \tau_y(x)$, respectively. If the two actions satisfy the compatibility condition
\begin{equation} \label{compat}
x y = \sigma_x(y) \tau_y(x)
\end{equation}
for all $x, y \in G$ and $r: G \times G \to  G \times G $ is the map given by
$$
r(x, y)= \big(\sigma_x(y), \, \tau_y(x)\big)
$$
for all $x, y \in G$, then $(G, r)$ is a non-degenerate bijective solution to the Yang--Baxter equation.
\end{theorem}

\begin{proof}
Recall that, by Proposition \ref{component conditions for a solution}, $(G, r)$ is a solution to the Yang--Baxter equation if and only if
\begin{eqnarray}
\sigma_{\sigma_x(y)} \big(\sigma_{\tau_y(x)} (z) \big) &=& \sigma_x \big( \sigma_y(z) \big), \label{CompAct (1)}\\
 \tau_{\sigma_{\tau_y(x)}(z)}  \big( \sigma_x(y) \big) &=&\sigma_{\tau_{\sigma_y(z)}(x)}  \big( \tau_z(y) \big), \label{CompAct (2)} \\
\tau_{\tau_z(y)} \big(\tau_{\sigma_y(z)} (x) \big) &=& \tau_z \big( \tau_y(x) \big), \label{CompAct (3)}
\end{eqnarray}
 hold for all $x, y, z \in X$. Using the fact that $(x, y) \mapsto \sigma_x(y)$ is a left-action and the compatibility condition \eqref{compat}, we have
$$ \sigma_{\sigma_x(y)} \big(\sigma_{\tau_y(x)} (z)\big)= \sigma_{\sigma_x(y)\tau_y(x)} (z)= \sigma_{xy} (z)= \sigma_{x} \sigma_y (z),$$
which establishes \eqref{CompAct (1)}. Similarly, using the fact that $(x, y) \mapsto \tau_y(x)$ is a right-action, we can prove \eqref{CompAct (3)}. A repeated use of the compatibility condition \eqref{compat} shows that the product of elements on the right hand side of \eqref{CompAct (1)}-\eqref{CompAct (3)} equals the element $xyz$, which further equals the product of elements on the left hand side of \eqref{CompAct (1)}-\eqref{CompAct (3)}. This shows that \eqref{CompAct (2)} also hold, and hence $(G,r)$ is a solution to the Yang--Baxter equation. 
\para
To see that $r$ is bijective, let $r(x, y) = \big(\sigma_x(y), \tau_y(x) \big)=(u, v)$. Then, the compatibility condition becomes $xy=uv$ and we get
$$
\left(\sigma_v(y^{-1}) \right) x = \sigma_v(y^{-1})  \tau_{y^{-1}}(v) = vy^{-1} = u^{-1} x.
$$
This implies that $y^{-1} =  \sigma_{v^{-1}}(u^{-1})$. Similarly, we have $x^{-1} = \tau_{u^{-1}}(v^{-1})$, and hence $r(v^{-1}, u^{-1})=(y^{-1}, x^{-1})$. It means that  if we put $s(x, y)=(y^{-1}, x^{-1})$, then
$(s r)^2 = \id_{X \times X}$, which implies that $r$ is bijective.
$\blacksquare$    \end{proof}          

Example \ref{OnGr} gives some elementary solutions on groups. Let us see which of these solutions satisfy the conditions of Theorem \ref{CompAct}.

\begin{example}{\rm 
Let $G$ be a group.
\begin{enumerate}
\item Consider the braiding given by $(x, y) \mapsto (y, y^{-k} x y^k)$. Then $\sigma_x = \id_G$ and $\tau_y(x) = y^{-k} x y^k$. We see that each $\tau_y$ is an automorphism of $G$. But,
$$
\tau_{y_1 y_2}(x) = (y_1 y_2)^{-k} x (y_1 y_2)^{k} \quad \textrm{and}\quad \tau_{y_2} \left(\tau_{y_1}(x) \right) = y_2^{-k} y_1^{-k} x y_1^{k} y_2^{k}.
$$
Hence, $\tau$ defines a right-action on any group if and only if $k=0, 1$. Further, the actions are compatible if and only if $k=1$.

\item For the braiding given by $(x, y) \mapsto (y, y x^{-1} y)$, we have $\sigma_x = \id_G$ and $\tau_y(x) =y x^{-1} y$. We see that in general $\tau_y$ is not an automorphism of $G$.
\item For the braiding given by $(x, y) \mapsto \big(y, \, \varphi(x y^{-1}) y \big)$, where $\varphi \in \Aut(G)$, we have $\sigma_x = \id_G$ and $\tau_y(x) =\varphi(x y^{-1}) y$. In this case, $ \tau_y$ is an automorphism of $G$ if and only if $y=1$.
\end{enumerate}}
\end{example} 

We now give an alternate description of solutions constructed in Theorem \ref{CompAct}.

\begin{definition}\label{def braiding operator}
 Let $G$ be a group with multiplication $m$. A \index{braiding operator}{\it braiding operator} on $G$ is a bijective map $r : G \times  G \to G \times G$ satisfying the following conditions:
\begin{enumerate}
\item For any $u, v, w \in G$,
\begin{equation} \label{braiding4}
r(u v, w)=(\id_G \times m) \, r_{12} \, r_{23}  \,(u, v, w), 
\end{equation}
\begin{equation} \label{braiding5}
r(u, v w)=(m \times \id_G)  \,r_{23}  \,r_{12} \, (u, v, w). 
\end{equation}

\item For any $u \in G$,
\begin{equation} \label{braiding6}
r(1, u)=(u, 1) \quad \textrm{and} \quad r(u, 1)=(1, u).
\end{equation}
\item For any $u, v \in G$,
\begin{equation} \label{braiding7}
m  \,r(u, v) = uv. 
\end{equation}
\end{enumerate}
\end{definition}

\begin{example}{\rm 
If $G$ is a group, then $r(u, v)=(v,  v^{-1} u v)$ is a braiding operator, called the \index{conjugate braiding}{\it conjugate braiding}. This braiding  comes from the Artin representation of the braid group $\B_n$ by automorphisms of the free group $F_n$ of rank $n$.}
\end{example}

\begin{definition}
Let $G$ and $A$ be groups. Let $G$ acts on $A$ on the left by automorphisms, denoted by $(g, a) \mapsto g \cdot a$ for $g \in G$ and $a \in A$. A bijective $1$-cocycle of $G$ with coefficient in $A$ is a bijection $\pi : G \to A$ such that 
$$
\pi (g h) = \pi(g) \, \big(g \cdot \pi(h) \big).
$$
holds for any $g, h \in G$.
\end{definition}

The following result gives a connection between braiding operators and bijective 1-cocycles \cite[Theorem 2]{MR1769723}.

\begin{theorem}\label{conjugate braiding thm}
The following statements are equivalent for any group $G$:
\begin{enumerate}
\item There exists a pair  $(\sigma, \tau)$ of compatible  left- and right-actions of $G$ on itself by automorphisms. 
\item There exists a braiding operator $r : G\times G \to G\times G$.
\item There exists a group $A$ on which $G$ acts on the left by  automorphisms and a bijective 1-cocycle $\pi : G \to A$. 
\end{enumerate}
\end{theorem}

\begin{proof} 
Let $(\sigma, \tau)$ be a pair of compatible  left- and right-actions of $G$ on itself by automorphisms. Let  $r : G\times G \to G\times G$ be given by
$$
r(u, v) = \big(\sigma_u(v), ~\tau_v(u) \big)
$$
for $u, v \in G$. To prove that $r$ is a braiding operator, we write
$$
r(u v, w) = (u_1, v_1) \quad \textrm{and} \quad (\id \times m)  \,r_{12} \, r_{23} \, (u, v, w) = (u_2, v_2).
$$
Direct calculations give
$$
u_1 = \sigma_{uv}(w)\quad \textrm{and} \quad u_2 = \sigma_{u} \sigma_{v}(w),
$$
which are equal since $\sigma$ is a left-action. Repeated use of the compatibility condition \eqref{compat} gives $uvw=u_1 v_1$ and $uvw=u_2 v_2$, which in turn implies that $v_1=v_2$. Similarly, one can prove that $r(u, v w)=(m \times \id_G)  \,r_{23}  \,r_{12} \, (u, v, w).$ Identities $r(1, u)=(u, 1)$ and $r(u, 1)=(1, u)$ are obvious, and  $m r(u, v) = uv$ is a direct  consequence of \eqref{compat}. Thus, assertion (1) implies assertion (2). 
\para
Conversely, suppose that $r$ is a braiding operator on $G$. We can write $r(x, y)= \big(\sigma_x(y), \tau_y(x)\big)$ for some maps $\sigma_x, \tau_y:G \to G$. It follows from \eqref{braiding4} that $\sigma_{uv}(w) = \sigma_{u} \sigma_{v}(w)$ for all $u, v,w \in G$. Moreover, \eqref{braiding6} implies that $\sigma_{1}(u) =u$ for all $u \in G$. Hence, $\sigma$ defines a left-action of $G$ on itself by automorphisms. Similarly, $\tau$ defines a right-action of $G$ on itself by automorphisms. The compatibility condition follows from \eqref{braiding7}. This proves that assertions (1) and (2) are equivalent.
\para

Let us prove that assertion  (1) implies assertion (3). Let $(\sigma, \tau)$ be a pair of  compatible actions of $G$ on itself by automorphisms. We take $A = G$ with the product
\begin{equation} \label{formula9}
u\star v = u \,\sigma_{u^{-1}}(v)
\end{equation}
for $u, v \in A$. By replacing $v$ in \eqref{formula9} with $\sigma_u (v)$, we obtain
$$
uv = u\star \sigma_u (v).
$$
Taking $\pi = \id_G : G \to A$, the preceding equation reads as
$$\pi(uv) = \pi(u)\star \big(u \cdot \pi(v) \big),$$
and hence $\pi$ is a bijective 1-cocycle. 
\para
It remains to show that $(A, \star)$ is indeed a group, and that $\sigma$ gives an action of $G$ on $A$ by automorphisms. 
Clearly, the identity 1 of the group $G$ is also the identity element of the groupoid $(A,\star)$. Further, $\sigma_{u} (u^{-1})$ is a right inverse of $u$ in $(A, \star)$. By the compatibility condition \eqref{compat}, we can rewrite $\star$ as
\begin{equation} \label{f10}
u \star v = u \big(\sigma_{u^{-1}}(v) \big) = v \big( \tau_v (u^{-1}) \big)^{-1}. 
\end{equation}
Then it is easy to see that $\left( \tau_{u^{-1}} (u)\right)^{-1}$ is a left inverse of $u$ in $(A, \star)$. But, by the compatibility condition \eqref{compat}, we have $\sigma_{u} (u^{-1})=\left( \tau_{u^{-1}} (u)\right)^{-1}$, and hence each element has an inverse in $(A, \star)$. A repeated use of the compatibility condition \eqref{compat} gives
\begin{equation} \label{f11}
\sigma_u (v x)  \tau_{v x} (u) =  u v x = \sigma_u (v)  \big(\tau_{v } (u)\big) x =
\big(\sigma_u (v) \big)  \big(\sigma_{\tau_{v } (u)} (x) \big)  \big(\tau_{x}(\tau_v (u)) \big).
\end{equation}
Since $\tau$ is a right-action, the rightmost factors on both sides are the same, and hence we  cancel this factor.  Taking  $x = \sigma_{v^{-1}}(w)$, we have
\begin{equation} \label{f12}
\sigma_u (v \star w) = \sigma_u (v x) \overset{(\ref{f11})}{=} \sigma_u (v)  \sigma_{\tau_v(u)}(x)\overset{(a)}{=} 
\sigma_u (v)  \sigma_{\tau_v(u)v^{-1}}(w)\overset{(\ref{compat})}{=} \sigma_u (v)  \sigma_{(\sigma_u(v))^{-1}u}(w)\overset{(a)}{=}
\sigma_u (v) \star \sigma_u (w),
\end{equation}
where $(a)$ means that we use the fact that $\sigma$ is a left-action. Further, we have
$$
u \star (v \star w) = u \big(\sigma_{u^{-1}}(v \star w) \big)  \overset{(\ref{f12})}{=} u \big(\sigma_{u^{-1}}(v)\big) 
\big( \sigma_{(\sigma_{u^{-1}}(v))^{-1}u^{-1}}(w) \big) = (u \star v) \big(\sigma_{(u \star v)^{-1}}(w)\big) =( u \star v ) \star w.
$$
Hence, $(A, \star)$ is a group and $\sigma$ is a left-action of $G$ on $(A, \star)$ by automorphisms.
\para 
To prove that assertion  (3) implies assertion (1), we use $\pi$ to identify $A$ with $G$. We take $\sigma$ to be the pullback of the action of $G$ on $A$ via the bijection $\pi$, defining an action of $G$ on itself. Let $\star$ denote  the pullback of the product of $A$ to $G$. Then $\sigma$ defines an action of $G$ on $(G, \star)$ as well. The 1-cocycle condition implies that \eqref{formula9} holds in $(G, \star)$. Defining $\tau$ as $\tau_v(u)=(\sigma_u(v))^{-1} u v$ according to the compatibility condition, we see  that \eqref{f10} also holds. It remains to show that $\tau$ defines a right-action of $G$ on itself. Since $\sigma$ is a left-action and $\pi(1)=1$, we have 
$$
1 \star u = u=u \star 1
$$
for all $u \in G$, and hence $1$ is also the identity of $(G,\star)$. Taking $v=1$ in \eqref{f10}, we get $\tau_1(u) = u$ for all $u \in G$. Note that, by the compatibility condition for $\sigma$ and $\tau$, we have the equality 
\eqref{f11}. Moreover, since $\sigma$ acts on $(G, \star)$, by taking  $x = \sigma_{v^{-1}}(w)$, we have
\begin{eqnarray*}
\sigma_{u}(v x) &=& \sigma_{u}(v \star w)\\
&=&  \sigma_{u}(v)  \star \sigma_{u}(w) \\
&=&  \sigma_{u}(v) \sigma_{(\sigma_{u}(v))^{-1}} \big(\sigma_{u}(w) \big)\\
&=&  \sigma_{u}(v) \sigma_{(\sigma_{u}(v))^{-1} u}(w)\\
&=& \sigma_{u}(v) \sigma_{(\tau_{v}(u)) v^{-1}}(w)\\
&=& \sigma_{u}(v)  \sigma_{\tau_v(u)}(x).
\end{eqnarray*}
It now follows from \eqref{f11} that $\tau_{v x} (u) = \tau_x \big(\tau_v(u) \big)$, and hence $\tau$ is a right-action of $G$ on itself.
$\blacksquare$   
 \end{proof}          

\begin{corollary}
Any braiding operator $r$ on a group $G$ defines a  non-degenerate bijective solution $(G, r)$ to the Yang--Baxter equation.
\end{corollary}

\begin{proposition}\label{commutative star product}
Let  $r$ be a braiding operator on a group $G$ such that $r(x, y)=\big(\sigma_x(y), \tau_y(x)\big)$.  Then  $r^2 =\id_{G \times G}$ if and only if
the product $\star$ given by
$$
u\star v = u \big(\sigma_{u^{-1}}(v) \big)
$$
 is commutative.
\end{proposition}

\begin{proof}
For each $u, v \in G$, we have
$$
r^{-1}(u, v) = \left( (\tau_{u^{-1}}(v^{-1}))^{-1},~  (\sigma_{v^{-1}}(u^{-1}))^{-1}\right)
$$
from the proof of bijectivity of $r$ in Theorem  \ref{CompAct}. Hence, the condition $r^2 = \id_{G \times G}$ is equivalent to
$$
\sigma_u(v) \tau_{u^{-1}}(v^{-1}) = 1.
$$
We have $u\star v = u \big(\sigma_{u^{-1}}(v) \big)$ and the compatibility condition for the pair $(\sigma, \tau)$ implies that $$v\star u = u  \big(\tau_{u}(v^{-1})\big)^{-1}.$$
Hence, the product $\star$ is commutative if and only if
$$
\sigma_{u^{-1}}(v) \tau_u(v^{-1}) = 1,
$$
which proves the proposition.
$\blacksquare$   
 \end{proof}

 \begin{proposition}\label{braidings on inverse}
 Let $r$ be a braiding operator on a group $G$. 
 If $ r(u, v) = (y, x)$, then 
 $$
 r(x^{-1}, y^{-1}) = (v^{-1}, u^{-1}),\quad r(u^{-1}, y) = (v, x^{-1}) \quad \textrm{and}\quad r(x, v^{-1}) = (y^{-1}, u).
 $$
 \end{proposition}

\begin{proof}
The proof of the equality $ r(x^{-1}, y^{-1}) = (v^{-1}, u^{-1})$  can be found in the proof of the invertibility of $r$ in Theorem \ref{CompAct}. To prove $ r(u^{-1}, y) = (v, x^{-1}),$ we note that $\sigma_{u^{-1}}(y) = v$ implies that $r(u^{-1}, y) = (v, z)$ for some $z$. Then from the compatibility condition \eqref{braiding7}, we have $u^{-1} y = v z$. On the other hand, applying the compatibility condition to
$r(u, v) = (y, x)$ gives $u v = y x$. Thus, we conclude that $z = x^{-1}$. The proof of the last equality is similar. $\blacksquare$   
 \end{proof}

Given a braiding operator $r$ on a group $G$, we constructed a new group $(G, \star)$ in the proof of Theorem \ref{conjugate braiding thm}. Define $r^\star: G \times G \to G \times G$  by setting
$$
r^{\star}(u, v) = (v, ~v^{\star(-1)} \star u \star v)
$$
for $u,v \in G$, where
$$
v^{\star(-1)} = \sigma_v(v^{-1}) = \big(\tau_{v^{-1}}(v) \big)^{-1}
$$
is the inverse of $v$ under the operation $\star$. Since $(G, \star)$ is a group, $r^\star$ is clearly a braiding operator (of quandle type) on $(G, \star)$, and called the \index{induced braiding operator}{\it induced braiding operator}. We now have the following result \cite[Theorem 6]{MR1769723}.

\begin{theorem}
Let $r$ be a braiding operator on a group $G$ and $r^{\star}$ the induced braiding operator on $(G, \star)$, where $r(u, v)=\big(\sigma_u(v),\, \tau_v(u)\big)$ for $u, v \in G$. Then the map $T_n: G^n \to G^n$ given by
$$
T_n(u_1, u_2, \ldots, u_{n-1}, u_n) = \big(u_1, ~\sigma_{u_1}(u_2),~ \ldots,~ \sigma_{u_1 \ldots u_{n-2}}(u_{n-1}), ~ \sigma_{u_1 \ldots u_{n-1}}(u_{n})\big)
$$
is an equivalence between the two representations $\B_n \to \Sigma_{G^n}$ of the braid group induced by $r$ and $r^{\star}$.
\end{theorem}

\begin{proof}
Since $r$ is a braiding operator on $G$, it follows from Theorem \ref{conjugate braiding thm} that $uv=\sigma_u(v) \tau_v(u)$ for all $u, v \in G$. Using this fact, the proof is now similar to that of Proposition \eqref{SolTh} and Proposition \ref{SolRack}. $\blacksquare$   
\end{proof}

\begin{remark}{\rm 
Although the braid group actions corresponding to the solutions $(G, r)$ and $(G, r^\star)$ to the Yang--Baxter equation are equivalent, the solutions themselves are not necessarily equivalent.}
\end{remark}

We conclude with the universal property satisfied by the structure group of a non-degenerate bijective solution to the  Yang--Baxter equation. The result is due to Lu, Yan, and Zhu \cite[Theorem 9]{MR1769723}, who also gave an elegant graphical interpretation of the maps involved.

\begin{theorem}\label{univ property non-deg solution st group}
Let $(X, r)$ be a  non-degenerate bijective solution to the Yang--Baxter equation. Let $G(X, r)$ be the structure group of the solution $(X, r)$ and $i_X: X \to G(X, r)$ the natural map. Then the following  assertions hold:
\begin{enumerate}
\item There is a unique braiding operator $r_G$ on $G(X, r)$ such that $r_G \,(i_X \times i_X) = (i_X\times i_X) \,r$, that is, $i_X$ is a morphism of solutions.
\item The group $G(X,r)$ and the braiding operator $r_G$ satisfy the following universal property: If $s$ is a braiding operator on a group $H$ and $\varphi: X \to H$ is a braiding-preserving map, then there is a unique braiding-preserving group homomorphism $\tilde{\varphi}:G(X, r) \to H$ such that $\varphi=\tilde{\varphi}\, i_X$.
\end{enumerate}
\end{theorem}

\begin{proof}
Let $(X, r)$ be a non-degenerate bijective solution to the Yang--Baxter equation. Denote
\begin{equation} \label{formula solution}
r(u, v) = (y, x)
\end{equation}
of $u, v, y, x \in X$. Since $r$ is  bijective, let $r^{-1} : X \times X \to X \times X$ be such that
$$
r^{-1} (y, x) = (u, v).
$$
It follows from the non-degeneracy of $(X,r)$ that for each $(u, y)$, there is a unique $(v, x)$ such that \eqref{formula solution} holds. Similarly, for each $(v, x)$, there is a unique $(u, y)$ such that \eqref{formula solution} holds. If we begin with a triple $(u, v, w)$, then
$$
r_{12}  \,r_{23}  \,r_{12}  \,(u, v, w) = (z, y, x)
=
r_{23}  \,r_{12} \, r_{23}  \,(u, v, w).
$$

We claim that  $r$ can be extended to a braiding operator $r_G:G(X, r)\times G(X, r) \to G(X, r) \times G(X, r)$.  Let $X'$ be another copy of $X$, with $x' \in X'$ denoting the element corresponding to the element $x \in X$. Let $\bar{X}$ be the disjoint union of $X$ and $X'$, and let
$$
U(X) = \bigsqcup_{i=0}^{\infty} \bar{X}^i = \{ e \} \bigsqcup \bar{X} \bigsqcup \bar{X} \times \bar{X} \bigsqcup \bar{X}  \times \bar{X}  \times  \bar{X} \bigsqcup \cdots
$$
be the free monoid generated by $\bar{X}$. Let $\sim $ be the equivalence relation on $U(X)$ generated by 
$$
g u u' h \sim  gh,~~ g u' u h \sim  gh,~~ g u v h \sim  g y x h~~\textrm{whenever}~~r(u, v) = (y, x),
$$
where $g, h \in  U(X)$ and $u, v, x, y \in X$. Then we have $G(X, r) = U(X) / \sim $.
\para

We begin the construction of $r_G$ by first extending $r$ to a solution $\bar{r} : \bar{X} \times \bar{X} \to \bar{X} \times \bar{X}$. Note that each $x' \in X'$ becomes $x^{-1}$ in $G(X, r)$.  Thus, in view of Proposition \ref{braidings on inverse}, the extension
$\bar{r}$ can be defined as follows. If $r(u, v) = (y, x)$ for $u, v, x, y \in X$, then we define
$$
\bar{r} (u, v) = (y, x),\quad \bar{r} (x', y') = (v', u'),\quad \bar{r} (u', y) = (v, x')\quad \textrm{and}\quad \bar{r} (x, v') = (y', u).
$$
It is not difficult to see that $(\bar{X},\bar{r})$ is a solution to the Yang--Baxter equation (see \cite[Proposition 8]{MR1769723} for a proof). 
\para 

Next, we extend $\bar{r}$ to a solution $r_U:U(X) \times U(X) \to U(X) \times U(X)$. The idea for the definition of $r_U$ comes from the formulas
$$
r(u v, w)=(\id_G \times m) \,r_{12} \,r_{23}  \,(u, v, w) 
\quad \textrm{and}\quad 
r(u, v w)=(m \times \id_G) \,r_{23} \,r_{12}  \,(u, v, w) 
$$
in Definition \ref{def braiding operator}. Note that the product of $u, v \in \bar{X}$ in $U(X)$ is the pair $(u, v) \in \bar{X}\times \bar{X}$. Hence, if $r_U$ satisfies the preceding two formulas, then we  have
$$
r_U \big((u, v), w \big) = \bar{r}_{12} \,\bar{r}_{23} \, (u, v, w) \quad \textrm{and}\quad r_U \big(u, (v, w)\big) = \bar{r}_{23} \,\bar{r}_{12}  \,(u, v, w),
$$
where the triples on the right sides are considered in $\bar{X} \times (\bar{X}\times \bar{X})$ and $(\bar{X}\times \bar{X}) \times \bar{X}$, respectively, since
$$
(\id_X \times m) (u, v, w) =  \big(u, (v, w)\big) \quad \textrm{and}\quad (m \times \id_X) (u, v, w) = \big((u, v), w\big) 
$$
in $U(X)$. This leads us to define $r_U \big((u_1, \ldots, u_m), (v_1, \ldots, v_n) \big)$. For instance,
$$
r_U \big((u_1, u_2, u_3), (v_1, v_2) \big) = \bar{r}_{23} \,\bar{r}_{12} \,\bar{r}_{34} \,\bar{r}_{23} \,\bar{r}_{45}\, \bar{r}_{34} \,\big((u_1, u_2, u_3), (v_1, v_2) \big).
$$
In general, we define $ r_U : U(X) \times U(X) \to U(X) \times U(X)$ 
by
$$r_U (e, u) = (u, e)\quad \textrm{and} \quad r_U (u, e) = (e, u),$$
and on $\bar{X}^m \times \bar{X}^n$ with $m, n \geq 1$ by
$$
r_U = (\bar{r}_{n,n+1} \ldots  \bar{r}_{23} \,\bar{r}_{12})  (\bar{r}_{n+1,n+2} \ldots  \bar{r}_{34}\, \bar{r}_{23}) \ldots  (\bar{r}_{n+m-1,n+m} \ldots  \bar{r}_{m+1,m+2} \,\bar{r}_{m,m+1}),
$$
where the image is considered to be in $\bar{X}^n \times \bar{X}^m$. It is easy to see  that the definition above is equivalent to the requirement that the first two conditions for $r_U$ to be a braiding operator on $U(X)$ (\eqref{braiding4}, \eqref{braiding5}, \eqref{braiding6}) are satisfied. Moreover, $\big(U(X),r_U\big)$ is a solution to the Yang--Baxter equation. 
\para 

It remains to prove that $r_U$ induces a braiding operator $r_G: G(X, r) \times G(X, r) \to G(X, r) \times G(X, r)$. We first consider the reduction with respect to the relations 
$$
g u u' h \sim  g h \quad \textrm{and} \quad g u' u h \sim  g h
$$
 in $U(X)$. Note that, if $r(u, v) = (y, x)$, then it is easy to see that 
$$
r_U \big((u, u'), y \big) = \big(y, (x, x') \big).
$$
Further,  $r_U(g u u' h, k)$  is obtained from $r_U(g h, k)$ by inserting some $x x'$ at appropriate places. Similarly, $r_U(g,  h u u' k)$, $r_U(g u' u h, k)$ and $r_U(g,  h u' u  k)$  are obtained from $r_U(g,  h k)$ and/or $r_U(g h, k)$ by inserting some $x x'$ and/or $x' x$ at appropriate places. This implies that $r_U$ is consistent with the relations
$$
g u u' h \sim  g h \quad \textrm{and} \quad g u' u h \sim  g h.
$$
\para 
Next, we consider the reduction with respect to the relation $g u v h \sim  g y x h$  whenever
$r(u, v) = (y, x)$. We first consider the simplest case. Suppose that $u, v, w \in X$, 
$ r(u, v) = (y, x)$,
$$
r_U \big((u, v), w \big) = (p_1, s_1, t_1) \quad \textrm{and} \quad r_U \big((y, x), w \big) = (p_2, s_2, t_2).
$$
Since $\big(U(X), r_U\big)$ is a solution to the Yang--Baxter equation, we conclude that $p_1 = p_2$ and $r(s_1, t_1) = (s_2, t_2)$.  Moreover, since $(\bar{X}, \bar{r})$ is a solution to the Yang--Baxter equation, the argument above also applies to the case when some of $u, v, w$ are in $X'$. We can show that $r_U(g y x h, k)$ is obtained from $r_U(g u v h, k)$ by applying $r$ at an appropriate adjacent
pair of coordinates. Similar statements can be proved regarding $r_U(g, h u v k)$ and $r_U(g, h y x k)$. This implies that $r_U$ is consistent with the relation  $g u v h \sim g y x h$. Hence, we have shown that $r_U$ can be reduced to an operator $r_G: G(X, r) \times G(X, r) \to G(X, r) \times G(X, r)$. Since $r_U$ satisfies the conditions \eqref{braiding4}, \eqref{braiding5} and \eqref{braiding6}, it follows that $r_G$ also satisfies the same conditions. 
\para 
To conclude that $r_G$ is a braiding operator on $G(X, r)$, it remains to show that $r_G$ satisfies the compatibility condition \eqref{braiding7}. By definition, the element $(g,h) \in  U(X) \times U(X)$ and $r_U (g, h)$ differ by successive applications of $r$ on some adjacent pairs of coordinates in the words (that is, ignoring the parentheses). Such applications of $r$ are considered as identity in $G(X, r)$. Thus, we conclude that $g h \in  G(X, r)$ is equal to the multiplication of the two coordinates in $r_G(g, h)$. This completes the construction of $r_G$. The uniqueness follows from the conditions  \eqref{braiding4}, \eqref{braiding5} and \eqref{braiding6}. Finally, the universal property is also obvious from the construction.
$\blacksquare$
\end{proof}

\begin{remark}{\rm 
It follows from Theorem \ref{univ property non-deg solution st group} that a solution $(X, r)$ to the Yang--Baxter equation is embedded in a group with a braiding operator if and only if the canonical map $i_X : X \to G(X,r)$ is an embedding. However, the canonical map is not always an embedding. In fact, we shall see in Part-II that solutions of quandle type do not always embed into their structure groups.}
\end{remark}

\begin{remark}\label{injectivity of eta remark}
{\rm 
In Theorem \ref{univ property non-deg solution st group}, if  we take $(X,r)$ to be a non-degenerate involutive solution, then it follows from the construction of $r_G$ that $\big(G(X, r), r_G\big)$ is also involutive. Further, by Proposition \ref{commutative star product}, the product $\star$ on $G(X,r)$ is commutative. Then it is not difficult to show that $\big(G(X,r),\star\big)$ is in fact a free abelian group with $X$ as a basis. This implies that the map $i_X$ is an embedding, recovering the classification theorem in \cite{MR1722951}.}
\end{remark}
\bigskip
\bigskip


\section{Indecomposable solutions}
To gain insight into the structure of solutions to the Yang--Baxter equation, we consider their invariant subsets.

\begin{definition}
Let $(X, r)$ be a solution to the Yang--Baxter equation. A subset $Y$ of $X$ is called \index{invariant subset}{invariant} if $r(Y\times Y) \subseteq Y \times Y$. Further, an invariant subset $Y$ of $X$ is said to be \index{non-degenerate invariant subset}{\it non-degenerate} if $(Y,~ r|_{Y \times Y})$ is a non-degenerate solution.
\end{definition}

The preceding definition leads to the following.

\begin{definition}
A solution $(X, r)$ to the Yang--Baxter equation is said to be \index{decomposable solution}{\it decomposable} if it is a union of two non-empty disjoint invariant subsets. Otherwise, $(X, r)$ is said to be \index{indecomposable solution}{\it indecomposable}.
\end{definition}

For example, a permutation solution to the Yang--Baxter equation is indecomposable if and only if it is cyclic.

\begin{remark}{\rm 
If $(X, r)$ is a finite non-degenerate involutive solution, then for an invariant subset $Y$ of $X$, the solution $(Y,~ r|_{Y \times Y})$ is also finite non-degenerate involutive. Indeed, since $r^2= \id_{X \times X}$, the map $r|_{Y\times Y}$ is bijective. Further, for each $y\in Y$, the maps $\sigma_y,\tau_y:Y \to Y$ are injective (by non-degeneracy of $X$), and hence bijective since $X$ is finite. Thus, $(Y,~ r|_{Y \times Y})$ is non-degenerate and involutive.}
\end{remark}

The next result gives a characterisation of indecomposability of a non-degenerate involutive solution \cite[Proposition 2.12]{MR1722951}.

\begin{proposition} \label{indecomposable vs transitivity}
A non-degenerate involutive solution $(X, r)$ is indecomposable if and only if $G(X, r)$ acts transitively on $X$.
\end{proposition} 

\begin{proof} 
Let $(X, r)$ be a non-degenerate involutive solution. We first claim that if $(X, r)$ is a union of two non-degenerate invariant subsets $X_1$ and $X_2$, then $r$ defines bijections $X_1 \times X_2 \to X_2 \times X_1$ and $X_2 \times X_1 \to  X_1 \times X_2$. Since $X_1$ and $X_2$ are invariant subsets, it follows that $r$ is a permutation of the set $X_1 \times X_2 \cup X_2 \times X_1$. Suppose that for $x_1 \in X_1$ and $x_2 \in X_2$, we have $r(x_1, x_2) = (y_1, y_2)$, where $y_1 \in X_1$ and $y_2 \in X_2$. Thus, $\sigma_{x_1}(x_2) = y_1$. By non-degeneracy of $X_1$, there exists $z_1 \in X_1$ such that $\sigma_{x_1}(z_1) = y_1$. But, this contradicts the non-degeneracy of $(X,r)$. Hence, our claim holds. 
\para 
It follows from the claim  that if $G(X, r)$ acts transitively on $X$, then $(X, r)$ is indecomposable. Conversely, if the action $x \mapsto \tau_x^{-1}$  is not transitive, then there exist two non-empty $G(X,r)$-invariant subsets $X_1$ and $X_2$ of $X$. Since the sets are invariant under each $\tau_x$, they are invariant under $T$, where $T(y) =  \tau_y^{-1}(y)$. Thus, by Proposition \ref{p123}(b), they are invariant under each $\sigma_x$. Hence, $X_1$ and $X_2$ are invariant subsets of $X$, which are clearly non-degenerate. This shows that  $(X, r)$ is decomposable.  $\blacksquare$ 
\end{proof}

Recall that, if a group $G$ acts on a set $X$ and $x\in X$, then the subgroup $\{g \in G \, \mid \, g\cdot x=x \}$ of $G$ is called the \index{stabilizer} \textit{stabilizer} of $x$ in $G$.  The next result  classifies finite indecomposable non-degenerate involutive solutions on a  set with $p$ elements, where $p$ is a prime \cite[Theorem 2.13]{MR1722951}.

\begin{theorem}
 Let $(X, r)$ be an indecomposable non-degenerate involutive solution  such that $|X| = p$, where $p$ is a prime. Then $(X, r)$ is isomorphic to the cyclic permutation solution $(\mathbb{Z}_p, r_0)$, where $r_0(x, y) = (y - 1,~ x + 1)$.
\end{theorem} 

\begin{proof}
Recall from Remark \ref{the g0 group} that $G^0(X, r)=G(X, r)/\Gamma$ and $A=\mathbb{Z}^X/\Gamma$, where $\Gamma$ is the kernel of the action of $G(X, r)$ on $X$. Then, by Proposition \ref{indecomposable vs transitivity},  the group $G^0(X, r)$ acts transitively on $X$ and its order is divisible by $p$. By Remark \ref{the g0 group}, $G^0(X, r)$ can be identified as a subgroup of $\Sigma_X$, and hence its order divides $p!$. This gives $|A|=|G^0(X,r)| = pn$, where $n$ is coprime to $p$, and  hence $A \cong \mathbb{Z}_ p  \oplus A_0$, where $|A_0|$ is coprime to $p$. Since $A_0$ is the group of all elements of $A$ whose order is not divisible by $p$, it follows that $A_0$ is invariant under the action of $G^0(X,r)$. The group 1-cocycle $\pi:G(X, r) \to \mathbb{Z}^X$ of Proposition \ref{p130} induces a group 1-cocycle $\bar{\pi} \colon G^0(X,r) \to A$. This further defines a group 1-cocycle $\pi' \colon G^0(X,r) \to  \mathbb{Z}_ p$. Set $H = (\pi')^{-1}(0) = \bar{\pi}^{-1}(A_0)$. It is easy to see that $H$ is a subgroup of $G^0(X, r)$ of order $n$. For an element $x \in  X$, let  $H_x$ be the stabilizer of $x$ in $G^0(X,r)$. Then $H_x$ is a subgroup of $G^0(X, r)$ of order $n$ and index $p$. 
\para 
We claim that  $H_x = H$ for all $x \in X$. Since $|H_x| = |H|$,  it is enough to show that $\pi'(H_x) = 0$. Let $\bar{t}_x \in  A / A_0$ be the image of the element $t_x \in \mathbb{Z}^X$. Then the set $\{ \bar{t}_x~\mid ~x \in X \}$ generates $A / A_0$ and each $\bar{t}_x$ is by definition fixed by $H_x$. Thus, $H_x$ acts trivially on $A/A_0=\mathbb{Z}_ p$. This implies that $\pi' |_{H_x}:H_x \to  \mathbb{Z}_ p$ is a group homomorphism. Since $|H_x| = n$ is coprime to $p$, this implies that $\pi'(H_x) =0$. It remains to prove that the image $\bar{t}_x \in  A / A_0$ of $t_x \in \mathbb{Z}^X$  is non-zero. Suppose  that $\bar{t}_x = 0$. Since $G(X,r)$ acts transitively on $X$, we get $\bar{t}_y = 0$ for all $y \in X$. This implies that the natural map $\mathbb{Z}^X \to  A / A_0$ is zero, which is a contradiction. Hence, $\bar{t}_x \neq 0$ for each $x \in X$ and the claim follows. It follows that $H$ acts trivially on $X$. This implies that $H=1$, and  hence $A_0 = 1$. Thus, $G^0(X,r) \cong A \cong  \mathbb{Z}_ p$, the action of $G^0(X,r)$ on $A$ is trivial and $\bar{\pi} = \id$.
\para
To conclude the proof, it is enough to observe that $\bar{t}_x = \bar{t}_y$ for any $x, y \in  X$. This follows from the fact that there exists $g \in  G^0(X,r)$ such that  $g \cdot x = y$, while the action of $G^0(X,r)$ on $A$ is trivial. Thus, $\bar{t}_x$ does not depend on $x$, and hence $\bar{\pi}^{-1} (\bar{t}_x) \in  G^0(X,r)$ does not depend on $X$. This proves that $(X, r)$ is a permutation solution.  $\blacksquare$
\end{proof}
\para

Let $(X, r_X)$ and $(Y, r_Y)$ be non-degenerate involutive solutions. Let $\Ext(X, Y)$ denote the set of all decomposable solutions $(Z, r)$ which are unions of $X$ and $Y$. We refer to elements of $\Ext(X, Y)$ as \index{extension of solutions}{\it extensions} of $X$ by $Y$. An element $(Z, r) \in \Ext(X, Y)$  is completely determined by the map
$$
r : X \times Y \to Y \times X.
$$
Let us write $r(x, y) = \big(\sigma_x(y), \tau_y(x)\big)$. Since $G(X, r_X)$ and $G(Y, r_Y)$ are subgroups of $G(Z, r)$, the following result is immediate  \cite[Proposition 2.17]{MR1722951}.
\para

\begin{proposition} 
If $(Z, r) \in \Ext(X, Y)$, then the assignments 
$$
x  \mapsto \sigma_x \quad \textrm{and} \quad  y  \mapsto \tau_y^{-1}
$$
are left-actions of  $(G, r_X)$ on $Y$ and of $G(Y, r_Y)$ on $X$, respectively.
\end{proposition} 
\para

\begin{definition}
An element $(Z, r) \in \Ext(X, Y)$ is called a \index{right-extension of solutions}{right-extension} of $X$ by $Y$ if $r(x, y) = \big(y, \tau_y(x)\big)$ for $x \in X$ and $y \in Y$. Similarly, $(Z, r)$ is called a \index{left-extension of solutions}{left-extension} of $X$ by $Y$ if $r(x, y) = \big(\sigma_x(y), x\big)$ for $x \in X$ and $y \in Y$.  The set of right- and left-extensions of $X$ by $Y$ will be denoted by $\Ext_+(X, Y)$ and $\Ext_-(X, Y)$, respectively.
\end{definition}
\para

It is clear that $\Ext_+(X, Y) =\Ext_-(Y, X)$. The following result gives a group-theoretic description of $\Ext_+(X, Y)$ \cite[Proposition 2.18]{MR1722951}.  
\para

\begin{proposition}
Let $Z= X \cup Y$ and $r: Z \times Z \to Z \times Z$ be such that $r(x, y) = \big(y, \tau_y(x)\big)$ for $x \in X$ and $ y \in Y$. Then, $(Z, r) \in \Ext_+(X, Y)$ if and only if  the assignment 
$$
y  \mapsto \tau_y^{-1}
$$
defines a left-action of  $G(Y, r_Y)$ on $X$.
\end{proposition}

\begin{proof} We have
$$
r_{12} r_{23} r_{12} (x_1, x_2, y) = \big(y, ~\tau_y \sigma_{x_1}(x_2),  ~ \tau_y  \tau_{x_2}(x_1) \big)
$$
and
$$
r_{23} r_{12} r_{23} (x_1, x_2, y) = \big(y, ~\sigma_{\tau_y(x_1)} \tau_{y}(x_2),  ~ \tau_{\tau_y(x_2)} \tau_y(x_1) \big)
$$
for $x_1, x_2 \in X$ and $y \in Y$. Equating the preceding equations shows that $\tau_y$ preserves $X$.
$\blacksquare$ 
\end{proof}
\para
\begin{corollary}
If $(Z,r) \in \Ext_+(X, Y)$,  then $$G(Z, r) \cong G(Y, r_Y) \ltimes G(X, r_X),$$ where $G(Y, r_Y)$-action on $G(X, r_X)$ is induced by the $G(Y, r_Y)$-action on $X$ via $y \mapsto \tau_y^{-1}$.
\end{corollary}
\para

We conclude this section with the following result \cite[Proposition 2.20]{MR1722951}.

\begin{proposition}  
Let $Z= X \cup Y$. Then the map $r(x, y) = \big(\sigma_x(y), \tau_y(x)\big)$ defines an element $(Z,r) \in \Ext(X, Y)$ if and only if the following conditions are simultaneously satisfied:
\begin{enumerate}
\item The assignments $x \mapsto\sigma_x$ and $y \mapsto \tau_y^{-1}$ are left-actions of $G(X, r_X)$ on $Y$ and of $G(Y, r_Y)$ on $X$.
\item The homomorphisms 
$$
\rho_X : G(X, r_X) \to \Sigma_Y \quad \textrm{and} \quad \rho_Y : G(Y, r_Y) \to \Sigma_X
$$
defined by assertion (1) are invariant under  $G(Y, r_Y)$ and $G(X, r_X)$, respectively. Here, the actions of  $G(X, r_X)$ on $\Sigma_X$ and  $G(Y, r_Y)$ on  $\Sigma_Y$ are trivial, and the actions of $G(X, r_X)$ and  $G(Y, r_Y)$ on each other are by $\rho_X$ and $\rho_Y$, respectively.
\end{enumerate}
\end{proposition} 
\bigskip 
\bigskip

\section{Linear  and affine  solutions on abelian groups}
In this section, we study non-degenerate solutions to the Yang--Baxter equation on abelian groups arising from their (affine) linear transformations
\para

We begin with the following definition.

\begin{definition}
Let $G$ be an abelian group. A solution $(G, r)$ to the Yang--Baxter equation is called a \index{linear solution}{\it linear solution} if $r$ has the form
\begin{equation} \label{linear solution expression}
r(x,y) = \big(ax + by,~ cx + dy \big)
\end{equation}
for some $a, b, c, d \in \End(G)$, where $\End(G)$ is the endomorphism ring of $G$.
\end{definition}

\begin{proposition}\label{nec suff condition for non-deg linear sol}
Let $G$ be an abelian group and $r : G \times G \to G \times G$ a map of the form \eqref{linear solution expression}. Then $(G, r)$ is a non-degenerate linear solution to the Yang--Baxter equation if and only if the following  conditions hold:
\begin{equation}\label{liner solution non-degen}
b, c \in \Aut(X),
\end{equation}
\begin{equation} \label{leq1}
a(1-a) = bac, \quad  d(1-d) = cdb,
\end{equation}
\begin{equation} \label{leq2}
ab = ba(1-d), \quad  ca = (1-d)ac, \quad  dc = cd(1-a),
\end{equation}
\begin{equation} \label{leq3}
bd = (1-a) db, \quad cb - bc = ada - dad.
\end{equation}
\end{proposition}

\begin{proof}
The proof follows immediately by evaluating $r_{12}r_{23} r_{12}(x, y, z)$ and $r_{23}r_{12} r_{23}(x, y, z)$, and using the condition of non-degeneracy. $\blacksquare$
\end{proof}

The following result determines non-degenerate linear solutions on abelian groups \cite[Lemma 8]{MR1809284}.

\begin{proposition} 
Let $G$ be an abelian group. Then non-degenerate linear solutions $(G, r)$ are in bijective correspondence with the quadruples $(a, b, d, s) \in \End(G)^4$ such that the following conditions hold:
\begin{enumerate}
\item $1-a$, $1-d$, $b$ and $1+s$ are invertible,
\item $s$ commutes with $a,b,d$ and $sa = sd = 0$,
\item $bdb^{-1} = (1-a)d$ and $b^{-1}ab = a(1-d)$. 
\end{enumerate}
In fact,  the bijective correspondence is given by $(a, b,c, d) \mapsto (a, b, d, s)$, where
\begin{equation}\label{linear and affine relation eq}
s=bc -(1-d+ad)(1-a).
\end{equation}
\end{proposition}

\begin{proof}
Let $(G, r)$ be a non-degenerate linear solution to the Yang--Baxter equation with
$$
r(x,y) = \big(ax + by,~ cx + dy \big).
$$
By Proposition \ref{nec suff condition for non-deg linear sol}, the quadruple $(a, b, c, d) \in \End(G)^4$ satisfies \eqref{liner solution non-degen}-\eqref{leq3}. The first equation of \eqref{leq3} implies that
$$
b d b^{-1}= (1-a) d,
$$
which further gives 
\begin{equation} \label{leq6a}
b (1 - d) b^{-1}= 1- d + a d.
\end{equation}
Multiplying the first equation of  \eqref{leq2} by $b^{-1}$ on the right gives $a = b a b^{-1} b (1 - d) b^{-1}$. Using \eqref{leq6a}, this further gives
\begin{equation} \label{leq6}
a=b a b^{-1} (1-d + ad).
\end{equation}
Similarly, using the second and the third equation of \eqref{leq2}, we obtain
\begin{equation} \label{leq7}
c a c^{-1}  = (1-d)a \quad \textrm{and} \quad c d c^{-1}  (1 - a + da) = d.
\end{equation}

Let $s$ be the endomorphism of $G$ defined by \eqref{linear and affine relation eq}. We claim that $s a = s d = 0$. Indeed, multiplying \eqref{linear and affine relation eq} by $a$ on the right and using the first relation of \eqref{leq1} and the second relation of \eqref{leq2}, we obtain
\begin{equation} \label{leq7a}
s a = b c a - (1 - d + ad) (1 - a) a = b (1 - d) a c - (1 - d + ad) b a c. 
\end{equation}
By \eqref{leq6a}, we have $b (1 - d) = (1 - d + a d)b$,  and using it in \eqref{leq7a} gives $sa=0$.
Similarly, we can show that $as = 0$.  The last equation of \eqref{leq3} implies that
\begin{equation} \label{leq8}
c b  =   (1 - a + da) (1 -  d) + s.
\end{equation}
Multiplying the preceding equation  by $d$  gives $sd = ds = 0$. We claim that the endomorphsims $1 - a$, $1 + s$ and $1 - d$ are invertible. Indeed, $b$ and $c$ are invertible by \eqref{liner solution non-degen}. Since
$$
b c = (1 -d + s + a d) (1 - a) = (1 - a) (1 -d + s +d a),
$$
$$
b c = (1 -d + a d) (1 - a + s) = (1 - a + s) (1 -d +  a d),
$$
$$
c b = (1 - a + s +  d a) (1 - d) = (1 - d) (1 - a + s + a d),
$$
we deduce that $1 - a$, $1 - a + s$ and $1 - d$  invertible. Considering the relation $1 - a + s = (1 - a) (1 + s)$,  it follows that $1 + s$ is also invertible. Our final claim is that $s$ commutes with $b$. We conjugate the equation \eqref{linear and affine relation eq} by $c$ and use equations \eqref{leq7} to conclude that
\begin{equation} \label{leq9}
c b  =   (1 - a + da) (1 -  d) + c s c^{-1}.
\end{equation}
Comparing \eqref{leq8} with \eqref{leq9} gives  $sc = c s$. Now, the equation \eqref{linear and affine relation eq}) implies that $s$ commutes with $b$ as well. Hence, starting with the quadruple $(a, b, c, d) \in \End(G)^4$ satisfying \eqref{liner solution non-degen}-\eqref{leq3}, we have constructed the quadruple $(a, b, d, s) \in \End(G)^4$ satisfying conditions (1)-(3) of the proposition.
\para 
Conversely, suppose that $(a, b, d, s) \in \End(G)^4$ satisfies conditions (1)-(3). Define $c$ by the formula 
$$c=b^{-1} \big(s+ (1-d+ad)(1-a) \big).$$
It is not difficult to check that the quadruple $(a, b, c, d) \in \End(G)^4$ satisfies the relations \eqref{liner solution non-degen}-\eqref{leq3}, which is desired. $\blacksquare$
\end{proof}

We can generalize linear solutions to affine solutions.

\begin{definition}
Let $G$ be an abelian group. A solution $(G, r)$ to the Yang--Baxter equation is called an \index{affine solution}{\it affine solution} if $r$ has the form
\begin{equation} \label{aff}
r(x,y) = \big(ax + by+z,~ cx + dy+t \big)
\end{equation}
for some $a, b, c, d \in \End(G)$ and $t, z \in G$.
\end{definition}

As in the linear case, the following result holds for affine solutions.

\begin{proposition}\label{prop affine solution condition}
Let $G$ be an abelian group and $r : G \times G \to G \times G$ a map of the form \eqref{aff}. Then $(G, r)$ is an affine solution to the Yang--Baxter equation  if and only if the conditions \eqref{liner solution non-degen}-\eqref{leq3} and 
\begin{equation} \label{affine solution condition}
cdz+dt = 0,\quad az+bat=0 \quad \textrm{and}\quad (c+d-ad-1)z + (da+1-a-b)t=0
\end{equation}
hold.
\end{proposition}

It follows from Proposition \ref{prop affine solution condition} that an affine solution $(G, r)$ gives rise to a linear solution $(G, r^*)$ with $r^*(x, y) = (ax + by, cx + dy)$. We refer to  $(G, r^*)$ as the linear part of the affine solution $(G, r)$. Similarly, we refer to $(G, r)$ as the affine solution associated with the linear solution $(G, r^*)$. We conclude with the following result \cite[Theorem 3.3 (1)]{MR1809284}.

\begin{proposition} 
Let $G$ be an abelian group and $(G, r^*)$ a non-degenerate linear solution with $r^*(x, y) = (ax + by, cx + dy)$. Then, $(G, r)$ given by 
$$
r(x, y)=\big(ax + by+z, \, cx + dy+t \big),
$$
where $a, b, c, d \in \End(G)$ and $z, t \in G$, is an affine solution associated with $(G, r^*)$ if and only if 
$$
ak = dk = 0 \quad \textrm{and}\quad (b-1)k = sz,
$$
where $s=bc -(1-d+ad)(1-a)$ and $k=c(1-a)^{-1}z +t $.
\end{proposition}

\begin{proof}
It follows from the first equation of \eqref{leq1} and the third equation of \eqref{leq2} that
$$
c d z + d t = d \big(c(1-a)^{-1} z + t \big) \quad \textrm{and} \quad az + bat = a \big(z + (1-a) c^{-1} t \big).
$$
Hence, we can rewrite the first two  equations of \eqref{affine solution condition} as
$$
d \big(c (1 - a)^{-1} z + t \big) = 0 \quad \textrm{and} \quad a \big(z +  (1 - a) c^{-1} t \big) = 0.
$$
If we define
\begin{equation} \label{leg13}
k=c (1 - a)^{-1} z +t,
\end{equation}
then it follows from \eqref{leq1} that $a k = d k = 0$. Now, \eqref{leq7} implies that
$$
c (1 - a)^{-1} c^{-1} = (1 - a + d a)^{-1},
$$
and we can modify \eqref{leg13} into
$$
(1 - a + d a) t = - c z + (1 - a + d a) k = - c z + k.
$$
We rewrite the last equation of  \eqref{affine solution condition} as 
\begin{equation} \label{leg14}
 (d - a d - 1) z - b t  + k = 0.
\end{equation}
If we substitute $t$ from \eqref{leg13} into \eqref{leg14} and use \eqref{linear and affine relation eq}, we obtain
$$
(b - 1) k = s z,
$$
which is the last desired equation. $\blacksquare$
\end{proof}
\bigskip
\bigskip

\chapter{Braces and skew braces}\label{chap braces and skew braces}

\begin{quote}
In this chapter, we develop the general algebraic theory of braces and skew braces. Braces were introduced by Rump to  provide non-degenerate involutive solutions to the Yang--Baxter equation. In an attempt to consider more general  non-degenerate bijective solutions, these structures were generalised to skew braces by Guarnieri and Vendramin. Beginning with basic definitions and plenty of examples, we discuss connections of braces with radical rings and pre-Lie algebras. We also give 
several constructions and properties of skew braces, and present the explicit form of the solution to the Yang–Baxter equation arising from skew braces.
\end{quote}
\bigskip

\section{Skew braces and radical rings}
In 2007, Rump \cite{MR2278047} introduced the notion of braces in the context of non-degenerate involutive solutions to the Yang--Baxter equation. Since then, this concept has garnered significant interest and has become a focal point of study. It is worth noting that braces had already been explored by Kurosh during his lectures \cite[Section 10]{MR0392756}. Following Rump's work, braces emerged as objects of fundamental significance. The notion was generalised by Guarnieri and Vendramin \cite{MR3647970} to skew braces, which give non-degenerate solutions to the Yang--Baxter equation that are not necessarily involutive.

\begin{definition}
A \index{skew left brace}{skew left brace} is a triple $(G, \cdot, \circ)$ such that $(G, \cdot)$ and $(G, \circ)$ are  groups and the compatibility condition
 \begin{equation}\label{slb}
a \circ (b \cdot c) =  (a \circ b) \cdot a^{-1} \cdot  (a \circ c)
 \end{equation}
holds for all $a, b, c \in G$, where $a^{-1}$ denotes the  inverse of $a$ in $(G, \cdot)$. We call $(G, \cdot)$ the \index{additive group}{additive group} and $(G, \circ)$ the \index{multiplicative group}{multiplicative  group} of the skew left brace $(G, \cdot, \circ)$. A skew left brace $(G, \cdot, \circ)$ is said to be a \index{left brace}{left brace} if $(G, \cdot)$ is an abelian group.  
\end{definition}
 
Similarly, a \index{skew right brace}{\it skew right brace} is a triple $(G, \cdot, \circ)$ such that $(G, \cdot)$ and $(G, \circ)$ are  groups and the compatibility condition
 \begin{equation}
(b \cdot c) \circ a=  (b \circ a) \cdot a^{-1} \cdot  (c \circ a)
 \end{equation}
holds for all $a, b, c \in G$.  If we change the multiplication in a skew left brace to the opposite multiplication, then we obtain a skew right brace, and vice-versa.  A skew left brace which is also a skew right brace is called a \index{skew two-sided  brace}{\it skew two-sided  brace}.
\para 

\begin{proposition}
If $(G, \cdot, \circ)$ is a skew left brace, then the identity element of both the groups $(G, \cdot)$ and $(G, \circ)$ is the same. 
\end{proposition}
\begin{proof}
Let $1$ be the identity of $(G, \cdot)$. We see that 
$$
a \circ 1 = a \circ (1 \cdot 1) = (a \circ 1) \cdot a^{-1} \cdot (a \circ 1)
$$
for all $a \in G$. Cancelling the element $(a \circ 1)$ with respect to the $\cdot$ operation gives $a = a \circ 1$ for all $a \in G$. Similarly, we can show that
$a = 1 \circ a$ for all $a \in G$. $\blacksquare$ 
 \end{proof}

\begin{remark}
{\rm The following conventions will be used throughout:
\begin{enumerate}
\item We denote the common identity element of the additive and the multiplicative group of a skew left brace by 1, and call it the identity of the skew left brace. 
\item When the additive group of a skew left brace is non-abelian, we sometimes use the notation $ab$ for $a \cdot b$. 
\item For a skew left brace $(G, \cdot, \circ)$, the inverse of an element $a \in G$ with respect to $\circ$ is denoted by $\bar a$ or $a^{\circ(-1)}$.
\item Sometimes we use the notation $G^{(\cdot)}$ to denote $(G, \cdot)$ and similarly  $G^{(\circ)}$ to denote $(G, \circ)$.
\item For a left brace, we denote the additive group operation by $+$. 
\end{enumerate}}
\end{remark} 
\para

\begin{example}\label{examples of solutions}
{\rm 
Let us consider some examples.
\begin{enumerate}
\item Any group $(G, \cdot)$ can be turned into a skew two-sided brace by putting $a \circ b= a \cdot b$ for all $a, b \in G$. Such a skew brace is called a \index{trivial brace}{\it trivial brace}. 
\item Any group $(G, \cdot)$ can be turned into a skew two-sided brace by putting $a \circ b = b \cdot a$ for all $a, b \in G$. If the group $(G, \cdot)$ is abelian, then we retrieve the skew brace of the preceding example. 
\item Let $G=\mathbb{Z}$ be the additive group of integers. Define $$n \circ m = n + {(-1)}^n m$$ for all $n, m \in \mathbb{Z}$. Then $(G, +, \circ)$ is a skew two-sided brace.
\item This example is due to Bachiller \cite[Example 2.2]{MR3835326}. Consider the group 
$$
G = \big\langle a, b \, \mid \, a^7 = b^3 = 1,~~b a b^{-1} = a^2 \big\rangle.
$$
It is easy to see that $G$ is a non-abelian group of order 21 and any element of $G$ can be presented in the form $a^i b^j$, where $0 \leq i \leq 2$ and $0 \leq j \leq 6$. Define the product
$$
(a^k b^l) \circ (a^i b^j) = a^{k + 4^l i} b^{l+j}.
$$
Then $(G, \cdot, \circ)$ is a skew left brace. Since
$$ (b \cdot b) \circ a = a^2 b^2 \quad \textrm{and} \quad (b \circ a) \cdot a^{-1} \cdot (b \circ a) = a^3 b^2,$$
this skew left brace is not two-sided.
\item Let $G$ and $K$ be groups and $\alpha : G \to \Aut(K)$  a group homomorphism. Then $G \times K$ is a skew left brace with
\begin{eqnarray*}
(a_1, b_1) \cdot (a_2, b_2) &=& (a_1 a_2, ~b_1 b_2),\\
(a_1, b_1)\circ (a_2, b_2) &=& \big(a_1 a_2, ~b_1 \alpha_{a_1}(b_2) \big),
\end{eqnarray*}
where $a_1, a_2 \in G$ and $b_1, b_2 \in K$.
\item If $(G_i, \cdot_i, \circ_i)$ is a family of skew left braces, then their direct sum $G = \oplus_i G_i$ is again a skew left brace in the usual way. More precisely, if $a = (a_i)$ and $b = (b_i) \in G$, then $a \cdot b = (a_i \cdot b_i)$ and $a \circ b = (a_i \circ b_i)$.
\end{enumerate}}
\end{example}

The following result gives some basic identities that hold in any skew left brace \cite[Lemma 1.7]{MR3647970}.

\begin{proposition}\label{idbr} 
The following  assertions hold for any skew left brace $(G, \cdot, \circ)$:
\begin{enumerate}
\item $a \circ (b^{-1} \cdot c) =  a \cdot (a \circ b)^{-1} \cdot  (a \circ c)$ for all $a, b, c \in G$.
\item $a \circ (b \cdot c^{-1}) =  (a \circ b) \cdot  (a \circ c)^{-1} \cdot a$ for all $a, b, c \in G$.
\end{enumerate}
\end{proposition}

\begin{proof}
Taking $c = b^{-1} \cdot d$, the skew left brace identity \eqref{slb} becomes
$$
a \circ d = (a \circ b) \cdot a^{-1}  \cdot \big(a \circ (b^{-1} \cdot d) \big) \Longleftrightarrow a \circ (b^{-1}  \cdot d) =  a \cdot (a \circ b)^{-1} \cdot  (a \circ d),
$$
which is assertion (1). Similarly, taking $b = d \cdot c^{-1}$ in \eqref{slb}, we get
$$
a \circ d = \big(a \circ (d \cdot c^{-1}) \big) \cdot a^{-1} \cdot (a \circ c)  \Longleftrightarrow a \circ (d  \cdot c^{-1} ) =  (a \circ d) \cdot  (a \circ c)^{-1} \cdot a,
$$
which is assertion (2). $\blacksquare$ 
\end{proof}

It is natural to expect that the underlying groups of a skew left brace are subject to certain constraints, as demonstrated by the following result. See \cite[Theorem 2.1]{MR2584610} and \cite[Theorem 2.15]{MR1722951}, as well as \cite[Theorem 5.2]{MR3824447} for a direct proof.

\begin{theorem}\label{multiplicative group slb solvable}
If $(G, +, \circ)$ is a finite left brace, then the multiplicative group $(G, \circ)$ is solvable.
\end{theorem}
\para

\begin{example}
{\rm  We give an example showing that the preceding theorem does not hold in general without the finiteness assumption \cite[Example 3.2]{MR4003478}.  For $n \ge 2$, let $U_n$ denote the set of strictly upper triangular $n \times n$ matrices over $\mathbb{Z}$. For $A, B \in U_n$, let $A + B$ denote the usual addition of matrices and let $A \circ B$ be defined by
$$A \circ B= (I_n + A)(I_n + B) - I_n,$$
where $I_n$ denote the identity matrix. Clearly, $(U_n, +)$ is isomorphic to the group  $\mathbb{Z}^{\frac{n(n-1}{2}}$ and $(U_n, \circ)$ is isomorphic to the group $UT(n, \mathbb{Z})$ of upper unitriangular matrices over $\mathbb{Z}$. For $A, B, C \in U_n$, we have
\begin{eqnarray*}
A \circ (B + C) &=& (I_n + A)(I_n + B + C) - I_n\\
&=& (I_n + A)(I_n + B) - I_n + (I_n + A)C\\
&=& (I_n + A)(I_n + B) - I_n + (I_n + A)(I_n + C) - (I_n + A)\\
&=& \big((I_n + A)(I_n + B) - I_n \big)- A + \big((I_n + A)(I_n + C) - I_n \big)\\
&=& (A \circ B)- A + (A \circ  C),
\end{eqnarray*}
and hence $A_n = (U_n, +, \circ)$ is a skew brace. Let $A=\oplus_{n \ge 2} A_n$ be the direct sum of $A_n$'s, equipped with the coordinate-wise operation. Then the additive group $(A, +)$ is isomorphic to the direct sum of infinitely many copies of $\mathbb{Z}$, and hence is abelian. The multiplicative group $(A, \circ)$ is isomorphic to the direct sum of groups $UT(n, \mathbb{Z})$ for $n \ge 2$. Since $UT(n, \mathbb{Z})$ is solvable of degree $\lceil \log_2(n) \rceil$, the group $(A, \circ)$ is not solvable. }
\end{example}

For two-sided skew braces, the following result is known \cite[Theorem 4.6]{MR4003478}.  

\begin{theorem}
Let $(G, \cdot, \circ)$ be a two-sided skew brace. If $(G,\circ)$ is nilpotent of nilpotency class $k$, then $(G, \cdot)$ is solvable of class at most $2k$.
\end{theorem}
\para

In \cite{MR3957824}, a two-sided brace is constructed whose multiplicative group contains a non-abelian free group, and hence is not solvable. The next result gives another such example \cite[Proposition 6.1]{MR4346001}.

\begin{proposition} \label{p7.1}
There is a two-sided brace $(A, +, \circ)$ such that its multiplicative group $(A,  \circ)$ contains a non-abelian free subgroup.
\end{proposition}

\begin{proof}
Let $\mathbb{Z}\langle\langle X_1, \ldots, X_n \rangle\rangle$ be the ring of formal power series in non-commuting variables $X_1, \ldots, X_n$ over $\mathbb{Z}$, where $n \ge 2$. Let $B$ be the two-sided ideal of $\mathbb{Z}\langle\langle X_1, \ldots, X_n \rangle\rangle$  generated by $X_1, \ldots, X_n$. Clearly, $(B, +)$ is an abelian group. Define a  binary operation $\circ$ on $B$ by setting
$$
a \circ b = ab + a + b
$$
for $a, b \in B$. We claim that $(B, +, \circ)$  is a left brace. Note that the binary operation $\circ$ is associative since the ring $\mathbb{Z}\langle\langle X_1, \ldots, X_n \rangle\rangle$ is so. Since
$$
0 \circ a = a \circ 0 = a
$$
for all $a \in B$, it follows that $0$ is the identity element of the semi-group $(B, \circ)$. Let $a, b \in B$ such that $a \circ b = 0$. Then, we have
$$
b = - a (1 + a)^{-1}  = - a + a^2 - a^3 + \cdots.
$$
Note that the element on the right side of the preceding equation  is well-defined and lies in $B$. Hence, $b$ is the right inverse of $a$, and it is easy to see that it is also the left inverse to $a$. This shows that $(B, \circ)$ is a group. Further, a direct check implies that $(B, +, \circ)$ is a two-sided brace.
\par
Let $F_n=\langle x_1, \ldots, x_n \rangle$ be the free group of rank $n$. It is well-known that there is an embedding $$
\mu : F_n \longrightarrow \mathbb{Z}\langle\langle X_1, \ldots, X_n \rangle\rangle,
$$
called the \index{Magnus embedding}{\it Magnus embedding} of the free group \cite[Theorem 5.6]{MR0207802}, which is defined on the generators and their inverses by 
$$
\mu(x_i) = 1 + X_i \quad \textrm{and}\quad \mu(x_i^{-1}) = 1 - X_i + X_i^2 - \cdots.
$$
It follows that $F_n$ is isomorphic to a subgroup of $(B, \circ)$, and hence $(B, \circ)$ contains a non-abelian free subgroup. $\blacksquare$ 
\end{proof}

\begin{remark}
{\rm
The examples of left braces constructed in  \cite{MR3957824} as well as Proposition \ref{p7.1} have infinitely generated additive groups.}
\end{remark}
\para

Let $(G, \cdot, \circ)$ be a skew left brace. Then, for each $a \in G$, there is a map $\lambda_a : G \to G$ defined by 
$$\lambda_a (b) = a^{-1 } \cdot (a \circ b)$$
for all $b \in G$. 

\begin{definition}\label{def skew brace homomorphism kernal ideal}
Let $(G, \cdot, \circ)$, $(G_1, \cdot_1, \circ_1)$ and $(G_2, \cdot_2, \circ_2)$ be skew left braces. 
\begin{enumerate}
\item A \index{sub skew left brace}{sub skew left brace} of $(G, \cdot, \circ )$ is a subset $H$ of $G$ such that  $(H, \cdot)$ is a subgroup of $(G,\cdot)$ and $(H, \circ)$ is a subgroup of $(G, \circ)$.
\item A homomorphism $f:(G_1, \cdot_1, \circ_1)\to (G_2,\cdot_2, \circ_2)$ of skew left braces is a map $f : G_1 \to G_2$ such that 
$$
f(a \cdot_1 b)  = f(a) \cdot_2  f(b) \quad \textrm{and}\quad f(a\circ_1 b)  = f(a) \circ_2 f(b)
$$
for all $a, b \in G$. 
\item The kernel of a homomorphism  $f:(G_1, \cdot_1, \circ_1)\to (G_2,\cdot_2, \circ_2)$ of skew left braces is defined as 
$$
\ker(f) = \big\{ a \in G_1\, \mid \, f(a) = 1 \big\},
$$
where $1$ denotes the identity of the skew left brace $(G_2,\cdot_2, \circ_2)$.
\item An ideal of $(G, \cdot, \circ)$ is a subset $I$ of $G$ such that $(I, \cdot)$ is a normal subgroup of $(G,\cdot)$, $(I, \circ)$ is a normal subgroup of $(G, \circ)$ and $\lambda_a(I) \subseteq I$ for any $a \in G$.  For example, the kernel of any skew left brace homomorphism is an ideal.
\end{enumerate}
\end{definition}

\begin{remark}
{\rm
A direct check shows that the set of all skew left braces with the above definition of a morphism forms a category, which we denote by $\mathcal{SLB}$. It is not difficult to see that the categories of skew left braces and skew right braces are equivalent, and hence we will consider only skew left braces throughout this monograph.}
\end{remark}

\begin{remark}
{\rm Describing the isomorphism types of skew  left braces of a given order is significantly more complex than describing all groups of the same order. For instance, while there exist 51 groups of order 32, there are a staggering 1223061 skew left braces and 25281 left braces of the same order  (see \cite{MR4113853}). Table \ref{number of skew left braces} borrowed from \cite{MR3647970} lists the number $c(n)$ of non-isomorphic skew left braces of size $n \le 30$. The task of comprehensively characterizing all skew  left braces, even with a known additive group, proves to be exceedingly challenging. As a result, there is a strong need for the construction of specific types of skew  left braces that can be more readily understood and analyzed.}
\end{remark}

\begin{table}[H]
\begin{center}
\begin{tabular}{|c|c|c|c|c|c|c|c|c|c|c|c|c|c|c|c|}
    \hline
$n$ & 1 & 2 & 3& 4 & 5 & 6 & 7 & 8 & 9 & 10 & 11 &  12 & 13 & 14 & 15 \\
\hline
$c(n)$ & 1 & 1 & 1 & 4 & 1 & 6 & 1 & 47 &4 & 6 & 1 & 38 &1 & 6 & 1 \\
\hline
\hline
$n$ & 16 & 17 & 18 & 19 & 20 & 21 & 22 & 23 & 24 & 25 & 26 & 27 & 28 & 29 & 30 \\
\hline
$c(n)$ & 1605 & 1 & 49 & 1 & 43 & 8 & 6 & 1 & 855 & 4 & 6 & 101 &29 & 1 & 36 \\
  \hline
\end{tabular}
\end{center}
\caption{The number of non-isomorphic skew left braces of order upto 30.} \label{number of skew left braces}
\end{table}

We conclude this section by briefly explaining the idea of Kurosh \cite[Section 10]{MR0080089}, which led to the concept of a brace. Let $(\mathbb{k}, +, \cdot)$ be an associative ring not necessarily with unity. Define an adjoint multiplication in $(\mathbb{k}, +, \cdot)$ by
$$
a \circ b = a + b - a \cdot b.
$$
Let $0$ be the zero element of the ring $(\mathbb{k}, +, \cdot)$. Since $ a \circ 0 = a= 0 \circ a,$ the element $0$ is the identity element with respect to the adjoint multiplication. Further, we see that
\begin{eqnarray*}
(a \circ b) \circ c 
&=& (a+b-a \cdot b) + c - (a + b - a \cdot b) \cdot c\\
&=& a + b + c -a \cdot b -a \cdot c-b \cdot c+(a \cdot b) \cdot c \\
&=& a+(b+c-b \cdot c)-a \cdot (b+c-b \cdot c)\\
&=& a \circ (b \circ c),
 \end{eqnarray*}
and hence the adjoint multiplication is associative. Thus, $(\mathbb{k}, \circ)$ is a semigroup with the identity element,  called the \index{adjoint semigroup of a ring}{\it adjoint semigroup} of the ring $(\mathbb{k}, +, \cdot)$.  The group of all invertible elements of $(\mathbb{k}, \circ)$ is called the \index{adjoint group of a ring}{\it adjoint group} of $(\mathbb{k}, +, \cdot)$. An associative ring $(\mathbb{k}, +, \cdot)$ is called a \index{radical ring}{\it radical ring} if $(\mathbb{k}, \circ)$ is a group. A ring is called \index{nilpotent ring}{\it nilpotent} if the product of arbitrary $n$ elements is zero for some positive integer $n$. 

\begin{proposition}
Every nilpotent ring is a radical ring.
\end{proposition}
\begin{proof}
Let $\mathbb{k}$ be a nilpotent ring and $a \in \mathbb{k}$ an arbitrary element. Consider the element
$$
x  = -a - a^2 - a^3 - \cdots,
$$
where $a^k= a \cdot a \cdot \cdots \cdot a$ ($k$-times). Since $\mathbb{k}$ is nilpotent, the element $x$ has only finitely many terms and 
$$
a \circ x =  a + x - a \cdot x = a+ (-a - a^2 - a^3 -\cdots )- (-a^2 - a^3 - a^4 -\cdots)=0.
$$
This shows that $x$ is the inverse of $a$ with respect to $\circ$, and hence  $\mathbb{k}$ is a radical ring. $\blacksquare$
\end{proof}

\begin{proposition}
If $(\mathbb{k}, +, \cdot)$ is a radical ring, then  $(\mathbb{k}, +, \circ)$ is a two-sided brace. Conversely, if $(\mathbb{k}, +, \circ)$ is a two-sided brace, then $(\mathbb{k}, +, \cdot)$ is a radical ring, where $a \cdot b= a +b - a \circ b$.
\end{proposition} 

\begin{proof}
If $(\mathbb{k}, +, \cdot)$ is a radical ring, then both $(\mathbb{k}, +)$ and $(\mathbb{k},\circ)$ are groups. Further, we have
$$ (a + b) \circ c = a \circ c + b \circ c - c \quad \textrm{and} \quad  c \circ (a + b) = c \circ a + c \circ b - c$$
for all $a, b, c \in \mathbb{k}$, and hence $(\mathbb{k}, +, \circ)$ is a two-sided brace.  
\para

Conversely, if $(\mathbb{k}, +, \circ)$ is a two-sided brace, then taking $a \cdot b= a +b - a \circ b$, a direct check shows that
$$(a+b) \cdot c= a \cdot c + b \cdot c, \quad  c \cdot (a + b)=  c \cdot a +  c \cdot b \quad \textrm{and} \quad  (a \cdot b) \cdot c= a \cdot (b \cdot c)$$
for all $a, b, c \in \mathbb{k}$. Hence,  $(\mathbb{k}, +, \cdot)$ is a radical ring. $\blacksquare$
\end{proof}
\bigskip
\bigskip

\section{Braces and pre-Lie algebras}
In the preceding section, we observed a relationship between braces and radical rings. In \cite{MR3291816}, Rump demonstrated a correspondence between left nilpotent right $\mathbb{R}$-braces and pre-Lie algebras. This correspondence was established through a geometric approach involving flat affine manifolds and affine torsors, and it operates at a local level. Smoktunowicz \cite{MR4391819} elucidated Rump's correspondence using purely algebraic formulas. In this section, we present some of these ideas.

\begin{definition}
A \index{pre-Lie algebra}{pre-Lie algebra} $(A,+, \cdot)$ is a vector space over a field $\mathbb{F}$ together with a $\mathbb{F}$-bilinear map $\cdot :A\times A\to A$ written as 
$(x, y) \mapsto  x \cdot y$ satisfying 
$$
(x \cdot y) \cdot z - x \cdot(y \cdot z) = (y\cdot x) \cdot z - y \cdot (x\cdot z)
$$
for all $x, y, z \in A$.
\end{definition}

Pre-Lie algebras are also called \index{left symmetric algebras}{\it left symmetric algebras} in the literature. Every associative algebra is clearly a pre-Lie algebra. Further, we can derive a  Lie algebra $L(A)$ from a pre-Lie algebra $(A,+, \cdot)$ by taking the Lie bracket as $$[a, b] = a \cdot b - b \cdot a.$$ It is known that pre-Lie algebras are in correspondence with the \'{e}tale affine representations of \index{nilpotent Lie algebra}{nilpotent Lie algebras} \cite{MR1486134} and Lie algebras with 1-cocycles \cite{MR4484613, MR2233854, MR3291816}.
\para 

\begin{definition}
A  pre-Lie algebra $(A,+, \cdot)$  is said to be  \index{nilpotent pre-Lie algebra}{\it nilpotent} if there exists some positive integer $n$ such that all products of $n$ elements in $A$ are zero. The smallest such integer is called the  \index{nilpotency index}{nilpotency index} of $(A,+, \cdot)$.
\end{definition}

Before proceeding further, we define a left ideal of a skew left brace.

\begin{definition}
A \index{left ideal}{\it left ideal} of a skew left brace $(G, \cdot, \circ)$  is a subgroup $(I, \cdot)$ of $(G, \cdot)$ such that $\lambda_a(I) \subseteq I$ for each $a \in G$.
\end{definition}

For a  left brace $(A, +, \circ)$ and $a, b \in A$, define $$a * b = a \circ b - a - b.$$
A direct check shows that $$a*(b+c)= a*b+a*c$$ for all $a, b, c \in A$. In \cite{MR2278047}, Rump introduced  the  left series 
$$A:=A^1  \supseteq A^2 \supseteq \cdots \supseteq A^n \supseteq \cdots $$
and the right series
$$A:=A^{(1)}  \supseteq A^{(2)} \supseteq \cdots \supseteq A^{(n)} \supseteq \cdots $$
of subsets of $(A, +, \circ)$ by defining
\begin{equation}\label{left and right series rump}
A^{n+1} = A * A^n \quad \textrm{and} \quad A^{(n+1)} = A^{(n)} * A
\end{equation}
for each $n \ge 1$. It follows from \cite{MR2278047} that the terms in the left series are left ideals of $(A, +, \circ)$, while those in the right series are ideals of $(A, +, \circ)$. 

\begin{definition}
A left brace $(A, +, \circ)$ is called \index{left nilpotent brace}{\it left nilpotent} if there exists some positive integer $n$ such that  $A^n = 0$. Similarly, a left brace $A$ is called \index{right nilpotent brace}{\it right nilpotent} if there exists some positive integer $n$  such that $A^{(n)} = 0$.
\end{definition}
\para

In \cite{MR3814340}, Smoktunowicz introduced the following chain of ideals for any left brace $(A, +, \circ)$. Set  $A^{[1]}= A$ and define
$$A^{[n+1]} = \sum_{i=1}^n A^{[i]} * A^{[n+1-i]}$$
for each $n \ge 1$. It is clear that 
$$ \cdots \subseteq A^{[n]} \subseteq  \cdots \subseteq A^{[2]} \subseteq A^{[1]}$$
 and each $A^{[n]}$ is an ideal of $(A, +, \circ)$.

\begin{definition}
A left brace $(A, +, \circ)$ is called \index{strongly nilpotent brace}{\it strongly nilpotent} if there exists a positive integer $n$ such that $A^{[n]} = 0$.
\end{definition}

The following result relates the two ideas of nilpotency for left braces \cite[Theorems 1.3. and 3.1]{MR3814340}.

\begin{theorem}\label{left and right nilpotent is strongly nilpotent}
A left brace is strongly nilpotent if and only if it is both left and right nilpotent.
\end{theorem}

\begin{proof}
 Let $(A, +, \circ)$ be a left brace such that $A^n = A^{(m)} = 0$ for some positive integers $n$ and $m$. We show that there is a positive integer $s$ such that $A^{[s]} = 0$. We proceed by induction on $m$. For $m = 2$, we have $0 = A^2 = A * A = A^{(2)} = A^{[2]}$ and the result holds. Suppose that there is a positive integer $s_{n,m}$ such that any left brace $(A, +, \circ)$ with $A^n =A^{(m)} = 0$ satisfies $A^{[s_{n,m}]} =  0$.  Further, suppose that  $A^{n} =  0$ and $A^{(m+1)} =  0$.  Let $p > n \, s_{n,m}$ and $a \in A^{[p]}$. Then we can write $a =\sum_i (a_i  * b_i)$ for some $a_i \in  A^{[p-q_i]}$ and $b_i \in  A^{[q_i]}$ for some $1 \le q_i \le p-1$. We see that if $q_i > s_{n,m}$, then $a_i \in  A^{(m)}$ (by the induction hypothesis applied to the left brace $A/A^{(m)}$). In this case, we get $a_i * b_i \in A^{(m)} * A = A^{(m+1)} = 0$. Hence, we can assume that each $q_i \leq s_{n,m}$. For each $i$, we can write $b_i = \sum_j (a_{i,j} * b_{i,j})$. By the same argument as before, we see that each $a_{i,j} \in  A^{[r_i]}$ for some $r_i \leq s_{n,m}$. Since $A$ is a left brace, we have
$$
a=\sum_i (a_i * b_i) = \sum_i  \left( a_i * \left(\sum_j a_{i,j} * b_{i,j} \right)\right)  = \sum_{i,j} a_i * (a_{i,j} * b_{i,j}).
$$
Continuing this way, we see that 
$$
a \in \sum_{c_1, \ldots, c_n \in A}  \big(c_n * (c_{n-1} * \cdots (c_2 * (c_1 * A)) \cdots )\big).
$$
Since $A^n =  0$, we get $a = 0$, and hence $A^{[p]} =  0$. The converse assertion is clear, and the proof is complete. $\blacksquare$
\end{proof}

The preceding theorem leads to the following result \cite[Theorem 3.2]{MR3814340}.

\begin{theorem}
Let $(A, +, \circ)$ be a left brace which is both left and right nilpotent. Then the multiplicative group $(A, \circ)$ is nilpotent.
\end{theorem}

\begin{proof}
We  construct a finite lower central series of $(A, \circ)$. For $a, b \in A$, we set $[a, b]= a \circ b \circ \bar{a} \circ \bar{b}$, where $\bar{a}, \bar{b}$ are the inverses of $a, b$ respectively in $(A, \circ)$. By Theorem \ref{left and right nilpotent is strongly nilpotent}, there is a positive integer $s$ such that $A^{[s]} = 0$. We proceed by induction on $s$. If $A^{[2]} = 0$, then $(A, \circ)$ is abelian, and the result holds. Suppose that the result holds for all positive integers smaller than $s$. By the induction hypothesis applied to $A/A^{[s-1]}$, we get 
$$\underbrace{[[[[A, A]A] \cdots]A]}_{m~{\rm brackets}} \in A^{[s-1]}$$ for some positive integer $m$.  Since $A^{[s-1]}$ is in the center of $A$, we get that 
$$\underbrace{[[[[A, A]A] \cdots]A]}_{(m + 1)~{\rm brackets}}= 0.$$ Hence, $(A, \circ)$ has a finite lower central series, and therefore is nilpotent. $\blacksquare$
\end{proof}

Let $(A, +)$ be a vector space over a field $\mathbb{F}$ such that $(A, + , \circ)$ is a left brace. Then $(A, + , \circ)$ is called a \index{left $\mathbb{F}$-brace}{\it left $\mathbb{F}$-brace} if $$a * (\alpha b) = \alpha (a * b)$$ for $a, b \in A$ and $\alpha \in \mathbb{F}$. Note that  $a * (\alpha b)= (\alpha a) * b$ for all $a, b \in A$ and $\alpha \in \mathbb{F}$.
\para 

We now associate a left $\mathbb{F}$-brace to a nilpotent pre-Lie algebra over $\mathbb{F}$ \cite{MR0579930, MR2839054, MR3291816}. 

\begin{theorem}\label{from pre lie algebra to brace}
Let $(A,+, \cdot)$ be a nilpotent pre-Lie algebra over a field $\mathbb{F}$ of characteristic zero (or of characteristic larger than the nilpotency index of $A$).
Then there is a left $\mathbb{F}$-brace $(A, +, \circ)$  with the same addition $+$ as that of  $(A,+, \cdot)$.
\end{theorem}

\begin{proof}
We define the multiplication $\circ$ as follows:
\begin{enumerate}
\item For each $a \in A$, let $L_a : A \to A$  denote the left multiplication $L_a(b) = a \cdot b$. Define 
$$
L_c \cdot L_b(a) = L_c \big(L_b(a)\big) = c \cdot (b \cdot a)
$$
and
$$
e^{L_a} (b) = b + a \cdot b + \frac{1}{2!} a \cdot (a \cdot b)  + \frac{1}{3!} a \cdot \big(a \cdot (a \cdot b)\big) + \cdots
$$
for $a, b, c \in A$.
\item We can formally consider element 1 such that $1\cdot  a = a \cdot 1 = a$ in the pre-Lie algebra and define
$$
W (a) = e^{L_a} (1) - 1 = a + \frac{1}{2!} a \cdot a   + \frac{1}{3!} a \cdot (a \cdot a) +  \cdots
$$
for $a \in A$. Since $A$ is a nilpotent pre-Lie algebra, the map $W : A \to A$ is bijective.

\item Let $\Omega : A \to A$  be the inverse of the map $W$, which is given by
$$
\Omega(a) = a -  \frac{1}{2} a \cdot a  + \frac{1}{4} (a \cdot a) \cdot a  + \frac{1}{12} a \cdot  (a \cdot a)  + \ldots
$$
for $a \in A$.
\item We define
$$
a \circ b = a + e^{L_{\Omega(a)}}(b)
$$
for $a, b \in A$. It follows from \cite{MR0579930} that $(A, \circ)$ is a group, called the \index{formal group of flows}{\it formal group of flows} of the pre-Lie algebra $(A, +, \cdot)$. Further, $(A, +, \circ)$ is a left $\mathbb{F}$-brace since
$$
a \circ (b + c) + a = a + e^{L_{\Omega(a)}}(b+c) + a = \big(a + e^{L_{\Omega(a)}}(b) \big) + \big(a + e^{L_{\Omega(a)}}(c)\big) = a \circ b + a \circ c
$$
and
$$a * (\alpha b) = \alpha (a * b)$$
for $a, b, c \in A$ and $\alpha \in \mathbb{F}$. Since $A$ is a nilpotent pre-Lie algebra, we have $A^n = 0$ for some positive integer $n$, and hence the above correspondence works globally. $\blacksquare$ 
\end{enumerate}
\end{proof}

The above formula can also be written using the \index{Baker--Campbell--Hausdorff formula}{Baker--Campbell--Hausdorff formula} and \index{Lazard's correspondence}{Lazard's correspondence} \cite{MR0579930, MR2839054}. In particular, in the left $\mathbb{F}$-brace $(A, +, \circ)$, we have
$$
W(a) \circ W(b) = W \big(C(a, b)\big),
$$
for all $a, b \in A$. Here, using the \index{Baker--Campbel--Hausdorf series}{Baker--Campbel--Hausdorf series} in the Lie algebra $L(A)$, we have
$$
C(a, b) = a + b +  \frac{1}{2} [a, b]  + \frac{1}{12} \big([a, [a, b]]  + [b, [b, a]] \big)  + \cdots
$$
for all $a, b \in A$.

\begin{example}{\rm 
Let $(A, +, \cdot)$  be a pre-Lie algebra over $\mathbb{F}$ such that $A^{[4]} = 0$. We calculate the formula for the multiplication $\circ$ in the corresponding left $\mathbb{F}$-brace $(A, +, \circ)$. Firstly, we have 
$$
\Omega(a) = a - \frac{1}{2} a \cdot a + c
$$
for some $c \in A^{[3]}$ and
$$
 e^{L_{\Omega(a)}}(b) = b + \Omega(a) \cdot b +  \frac{1}{2} \Omega(a) \big(\Omega(a) \cdot b\big).
$$
Thus, we have
$$
a \circ b = a + b + a \cdot b - \frac{1}{2}(a \cdot a) \cdot b + \frac{1}{2} a \cdot (a \cdot b),
$$
and hence
$$
a * b = a \cdot b - \frac{1}{2}(a \cdot a) \cdot b + \frac{1}{2} a \cdot (a \cdot b).
$$}
\end{example}

The preceding constructions lead to the following result \cite[Theorem 3.1]{MR4391819}.

\begin{theorem} \label{smokt1}
Let $(A, +, \cdot)$  be a nilpotent  pre-Lie algebra over a field $\mathbb{F}$ of characteristic zero and $(A, +, \circ)$ the corresponding left $\mathbb{F}$-brace. Then
$$
a * b = a \cdot b + \sum_{x \in B} \alpha_x x,
$$
where $\alpha_x \in \mathbb{F}$ and $B$  is the set of all products of elements $a$ and $b$ in $(A, \cdot)$ with $b$ appearing only at the end, and $a$ appearing at least two times in each product.
\end{theorem}

\begin{proof}
The proof follows from the construction of $\Omega(a)$, which is a sum of $a$ with a linear combination of all possible products of more than one element $a$ with any
distribution of brackets, which can be proved by induction.  $\blacksquare$    
\end{proof}

Next, we obtain a pre-Lie algebra from a left $\mathbb{F}$-brace \cite[Theorem 5.1 and Theorem 6.2]{MR4391819}.

\begin{theorem} \label{smokt2}
Let $(A, +, \circ)$ be a strongly nilpotent left $\mathbb{Q}$-brace. For $a, b \in A$, define 
$$
a \cdot b = \lim_{n \to \infty} 2^n \left( \frac{1}{2^n} a \right) * b.
$$
Then the following assertions hold:
\begin{enumerate}
\item $(A, +, \cdot)$ is a pre-Lie algebra over $\mathbb{Q}$.
\item  $(A,  \circ)$ is the formal group of flows of the pre-Lie algebra $(A, +, \cdot)$.
\item  $(A, +, \circ)$ can be obtained from the pre-Lie algebra $(A, +, \cdot)$ as in Theorem \ref{from pre lie algebra to brace}.
\end{enumerate}
\end{theorem}

We also have the following result \cite[Proposition 6.1]{MR4391819}.

\begin{proposition}
Let $(A, +, \cdot)$ be a nilpotent pre-Lie algebra over a field $\mathbb{F}$  of characteristic zero and  $(A, +, \circ)$ the corresponding left $\mathbb{F}$-brace. Suppose that $\mathbb{F} = \mathbb{R}$ or $\mathbb{Q}$. Then for every $a, b \in A$, the limit 
$$
\lim_{n \to \infty} 2^n \left( \frac{1}{2^n} a \right) * b
$$
exists and
$$
a \cdot b = \lim_{n \to \infty} 2^n \left( \frac{1}{2^n} a \right) * b.
$$
\end{proposition}

\begin{proof}
The result follows immediately from the fact that the multiplication in a pre-Lie algebra is bilinear, and from the fact that $a * b$ can be expressed as in Theorem \ref{smokt1}. $\blacksquare$ 
\end{proof}

The preceding results yield the following correspondence.

\begin{corollary}
 There is a one-to-one correspondence between the set of strongly nilpotent left $\mathbb{Q}$-braces and the set of nilpotent pre-Lie algebras over $\mathbb{Q}$.
\end{corollary}
\bigskip
\bigskip


\section{Properties and constructions of skew braces}
This section delves into the characteristics of skew  left braces, some of their constructions and maps linked to them.

\subsection{Maps associated to skew braces}

Let $(G, \cdot, \circ)$ be a skew left brace. For each $a \in G$, define $\lambda_a : G \to G$ by 
$$\lambda_a (b) = a^{-1 } \cdot (a \circ b)$$
for all $b \in G$. A direct check shows that each $\lambda_a$ is a bijection, with its inverse given by 
\begin{equation}\label{def lambda map}
\lambda_a^{-1} (b) = \bar{a} \circ(a \cdot b)
\end{equation}
for all $b \in G$.  This leads to the map $\lambda: G \to \Sigma_G$ given by $\lambda(a)= \lambda_a$. Equation \ref{def lambda map} allows us to write the two group operations in $(G, \cdot, \circ)$ in terms of each other, that is,
\begin{equation}\label{circ and dot by lambda}
a \circ b = a \cdot \lambda_a (b) \quad \textrm{and}\quad a \cdot b = a \circ \lambda_a^{-1} (b)    
\end{equation}
for $a, b \in G$.
\para

The following result characterises a skew left brace in terms of the map $\lambda$ \cite[Proposition 1.9]{MR3647970}.

\begin{proposition} \label{lambda-map}
Let $G$ be a set with two binary operations $\cdot$ and $\circ$ such that both $(G, \cdot)$ and $(G, \circ)$ are groups. Let $\lambda : G \to \Sigma_G$ be the map given by $\lambda(a) =  \lambda_a$, where $\lambda_a(b) =  a^{-1} \cdot (a \circ b)$ for $a, b \in G$. Then the following assertions are equivalent:
\begin{enumerate}
\item $(G, \cdot, \circ)$ is a skew left brace.
\item $\lambda_{a\circ b}(c) =\lambda_a \lambda_b(c)$ for all $a, b, c \in G$.
\item $\lambda_{a} (b  \cdot c) =\lambda_a (b) \cdot \lambda_a (c)$ for all $a, b, c \in G$.
\end{enumerate}
\end{proposition}

\begin{proof}
Let $(G, \cdot, \circ)$ be a skew left brace and $a, b, c \in G$. By Lemma \ref{idbr}(2), we have
$$
a \circ b^{-1} = a \cdot(a \circ b)^{-1}  \cdot a.
$$
This gives 
\begin{eqnarray*}
\lambda_a \lambda_b(c) 
&=& a^{-1}\cdot \big(a \circ \lambda_b(c) \big)\\
&=& a^{-1}  \cdot \big(a \circ (b^{-1} \cdot (b \circ c)) \big) \\
&=& a^{-1}\cdot  (a \circ b^{-1})\cdot  a^{-1} \cdot (a \circ b \circ c)\\
&=& (a \circ b)^{-1} \cdot (a \circ b \circ c)\\
&=& \lambda_{a\circ b}(c),
\end{eqnarray*}
which establishes the implication (1) $\Rightarrow$ (2).
\para 

Suppose that assertion (2) holds. Using the identity $ a \cdot b = a \circ \lambda_a^{-1} (b)$, we get
\begin{eqnarray*}
\lambda_{a} (b\cdot c) &=&  \lambda_{a} \big(b \circ \lambda_{b}^{-1} (c)\big) \\
&=&  a^{-1} \cdot \big(a \circ b  \circ \lambda_b^{-1}(c)\big) \\
&=&  a^{-1} \cdot (a \circ b)  \cdot (a \circ b)^{-1} \cdot \big(a \circ b \circ \lambda_b^{-1}(c)\big) \\
&=&  \lambda_a(b) \cdot \lambda_{a \circ b} \lambda_b^{-1}(c) \\
&=&  \lambda_a(b) \cdot \lambda_a \lambda_b \lambda_b^{-1}(c)\\
&=&  \lambda_a (b) \cdot \lambda_a (c),
\end{eqnarray*}
which proves the implication (2) $\Rightarrow$ (3).
\para 

Suppose that assertion (3) holds. Then, for $a, b, c \in G$, we have
$$
a^{-1} \cdot \big(a \circ (b \cdot c)\big) =  \lambda_a(b \cdot c)  =  \lambda_a(b) \cdot \lambda_a(c) = a^{-1} \cdot (a \circ b) \cdot a^{-1} \cdot (a \circ c),
$$
which proves the implication (3) $\Rightarrow$ (1). $\blacksquare$ 
\end{proof}

As a consequence, we have the following result.

\begin{corollary}
Let $(G, \cdot, \circ)$ be a skew left brace. Then  the map
$$
\lambda : (G, \circ) \to \Aut(G, \cdot)
$$
given by $a \mapsto \lambda_a$ is a group homomorphism.
\end{corollary}

Let $(G, \cdot, \circ)$ be a skew left brace. Then,  for each $b \in G$, define a map $\gamma_b : G \to G$  by 
$$\gamma_b(a) = \lambda^{-1}_{\lambda_a(b)} \left( (a \circ b)^{-1} \cdot a \cdot (a \circ b)\right)$$
for $a \in G$. These maps satisfy some interesting properties \cite[Lemma 2.4]{MR3835326}.

\begin{proposition}\label{gamma-map}
Let $(G, \cdot, \circ)$ be a skew left brace. Then, for any $a, b \in G$, the following  assertions hold:
\begin{enumerate}
\item $\gamma_b(a) = \big(\bar{b} \circ (\bar{a} \cdot b) \big)^{\circ(-1)} = \big((\bar{b} \circ \bar{a} ) \cdot (\bar{b})^{-1} \big)^{\circ(-1)}.
$
\item $\gamma_b$ is a bijection with $\gamma_{\bar{b}}$ as its inverse.
\item $\gamma_a \gamma_b = \gamma_{b \circ a}$.
\end{enumerate}
\end{proposition}

\begin{proof}
For $a, b \in G$, we have
$$
\lambda^{-1}_{\lambda_a(b)} \left( \lambda_a(b)^{-1} \cdot a \cdot \lambda_a(b) \right) = \lambda^{-1}_{\lambda_a(b)} \left( (a \circ b)^{-1} \cdot a \cdot a \cdot a^{-1} \cdot (a \circ b) \right) =\lambda^{-1}_{\lambda_a(b)} \left( (a \circ b)^{-1} \cdot a \cdot (a \circ b)\right) = \gamma_b(a).
$$
Writing $\gamma_b$ in terms of the operations of $G$, we see that
\begin{eqnarray*}
\gamma_b(a) &=& \lambda^{-1}_{\lambda_a(b)} \big( \lambda_a(b)^{-1} \cdot a \cdot \lambda_a(b) \big)\\
&=& \overline{\lambda_a(b)} \circ \big(\lambda_a(b) \cdot \lambda_a(b)^{-1} \cdot a \cdot \lambda_a(b) \big), \quad \textrm{using} ~\eqref{circ and dot by lambda}\\
&=& \overline{\lambda_a(b)} \circ \big(a \cdot \lambda_a(b)\big) \\
&=& \overline{\lambda_a(b)} \circ a \circ b \\
&=& \big(\bar{b} \circ \bar{a} \circ \lambda_a(b) \big)^{\circ(-1)}\\
&=& \big(\bar{b} \circ (\bar{a} \cdot b) \big)^{\circ(-1)},\quad \textrm{using}~ \eqref{circ and dot by lambda}\\
&=& \big((\bar{b} \circ \bar{a}) \cdot (\bar{b})^{-1} \big)^{\circ(-1)},\quad \textrm{using the defining condition of a skew left brace,}
\end{eqnarray*}
which proves assertion (1).
\para 

A direct check shows that $\gamma_b$ is invertible, with inverse given by
$$
a \mapsto \big(b \circ (\bar{a} \cdot \bar{b}) \big)^{\circ(-1)}, 
$$
which establishes assertion (2).
\para 

Now, we check that $\gamma: (G, \circ) \to \Sigma_G$  is an anti-homomorphism. For $a, b, c \in G$, we see that 
\begin{eqnarray*}
\gamma_{b \circ a}(c) &=& \big((\bar{a} \circ \bar{b} \circ \bar{c}) \cdot (\bar{a} \circ\bar{b})^{-1}\big)^{\circ(-1)}\\
&=& \big((\bar{a} \circ \bar{b} \circ \bar{c}) \cdot (\bar{a} \cdot \lambda_{\bar{a}}(\bar{b}))^{-1}\big)^{\circ(-1)}\\
&=&  \big((\bar{a} \circ \bar{b} \circ \bar{c}) \cdot \lambda_{\bar{a}}(\bar{b})^{-1} \cdot (\bar{a})^{-1}\big)^{\circ(-1)}\\
&=&  \big((\bar{a} \circ \bar{b} \circ \bar{c}) \cdot  \lambda_{\bar{a}}((\bar{b})^{-1}) \cdot (\bar{a})^{-1} \big)^{\circ(-1)} \\
&=& \big((\bar{a} \circ \bar{b} \circ \bar{c}) \cdot  (\bar{a})^{-1} \cdot  (\bar{a} \circ (\bar{b})^{-1}) \cdot (\bar{a})^{-1}\big)^{\circ(-1)}.
\end{eqnarray*}
On the other hand, we have
$$
\gamma_a \big(\gamma_b(c)\big) = \gamma_a\big(\big(\bar{b} \circ \bar{c}) \cdot (\bar{b})^{-1}\big)^{\circ(-1)}\big) = 
\big((\bar{a} \circ ((\bar{b} \circ \bar{c}) \cdot (\bar{b})^{-1})) \cdot (\bar{a})^{-1}\big)^{\circ(-1)} = 
\big((\bar{a} \circ \bar{b} \circ \bar{c})  \cdot (\bar{a})^{-1} \cdot (\bar{a} \circ (\bar{b})^{-1}) \cdot (\bar{a})^{-1}\big)^{\circ(-1)},
$$
which proves assertion (3). $\blacksquare$ 
\end{proof}

Given a skew left brace $(G, \cdot, \circ)$, we can define the semi-direct product $(G, \cdot) \rtimes (G, \circ)$ by the action $\lambda : (G, \circ) \to \Aut(G, \cdot)$. We have the following observation regarding this group \cite[Lemma 2.5]{MR3835326}.

\begin{proposition} \label{theta-map}
Let $(G, \cdot, \circ)$ be a skew left brace. Consider the map
$$
\theta : (G, \cdot) \rtimes (G, \circ) \to \Aut(G, \cdot)$$
given by $(a, b) \mapsto \theta_{(a, b)}$, where $\theta_{(a, b)}(c) = a \cdot \lambda_b(c) \cdot a^{-1}$ for $a, b, c \in G$. Then $\theta$ is a group homomorphism. 
\end{proposition}

\begin{proof}
Since $\theta_{(a, b)}$ is just the composition of $\lambda_b$ with the inner automorphism induced by $a$, it follows that $\theta_{(a, b)} \in \Aut(G, \cdot)$ for each $a, b \in G$. For $a_1, b_1, a_2, b_2, c \in G$, we have
\begin{eqnarray*}
\theta_{(a_1, b_1) (a_2, b_2)}(c) &=&  \theta_{(a_1 \cdot \lambda_{b_1}(a_2), \,b_1 \circ b_2)}(c)\\
&=& \big(a_1 \cdot \lambda_{b_1}(a_2) \big) \cdot \lambda_{b_1 \circ b_2}(c) \cdot  \big(a_1 \cdot \lambda_{b_1}(a_2)\big)^{-1}\\
&=& a_1 \cdot \lambda_{b_1}(a_2) \cdot \lambda_{b_1 \circ b_2}(c) \cdot \lambda_{b_1}(a_2)^{-1} \cdot a_1^{-1}\\
&=&  a_1 \cdot \lambda_{b_1} \big(a_2 \cdot \lambda_{b_2}(c) \cdot a_2^{-1} \big) \cdot a_1^{-1}\\
&=&  \theta_{(a_1, b_1)} \theta_{(a_2, b_2)}(c),
\end{eqnarray*}
and hence $\theta$ is a homomorphism. $\blacksquare$ 
\end{proof}

Recall from Definition \ref{def skew brace homomorphism kernal ideal} that, an ideal of a skew left brace $(G, \cdot, \circ)$ is a subset $I$ such that $(I, \cdot)$ is a~normal subgroup of  $(G, \cdot)$, $(I, \circ)$ is a normal subgroup of $(G, \circ)$ and $\lambda_a(b) \in I$ for any $a \in G$ and $b \in I$. Ideals of  skew left braces enjoy the following properties \cite[Lemma 2.3]{MR3647970}.

\begin{proposition} 
The following  assertions hold for an ideal~$I$ of a skew left brace $(G, \cdot, \circ)$:
\begin{enumerate}
\item  $a \circ I = a \cdot I$ for all $a \in G$.
\item  $I$ and $G/I$ are skew left braces.
\end{enumerate}
\end{proposition}

\begin{proof}
Assertion (1) is immediate from \eqref{circ and dot by lambda}. Clearly, $I$ is a skew left brace. By assertion (1),  the set of left cosets of $(I, \cdot)$ in $(G, \cdot)$ is the same as the set of left cosets of $(I, \circ)$ in $(G, \circ)$. Hence,  $G/I$ has the natural induced skew left brace structure. $\blacksquare$ 
\end{proof}
\para

\begin{definition}
The \index{socle}{socle} of a skew left brace $(G, \cdot, \circ)$ is defined as
$$
\Soc(G) = \{ a \in G \, \mid \,  \lambda_a = \id ~~ \textrm{and}~~ a \cdot b = b \cdot a~\mbox{for all}~ b \in G \} = \ker(\lambda) \cap \Z(G, \cdot).
$$
\end{definition}
\para

The following result gives some interesting properties of the socle of a skew left brace \cite[Proposition 2.8]{MR3835326}.

\begin{proposition} 
Let $(G, \cdot, \circ)$ be a skew left brace. Then the following  assertions hold:
\begin{enumerate}
\item $\Soc(G)$ is an ideal of $G$ contained in the center of $(G, \cdot)$. If $a \in G$ and $y \in \Soc(G)$, then $$\lambda_a(y) = a \circ y \circ \bar{a}.$$  Moreover, $\Soc(G)$ is a trivial brace, and consequently an abelian group.
\item $\Soc(G)$ coincides with the kernel of the homomorphism
$$
\varphi : (G, \circ) \to \Aut(G, \cdot) \times \Sigma_G
$$
given by
$$\varphi(a) = (\lambda_a, \gamma_a^{-1}).$$
\item $\Soc(G)$ coincides with the kernel of the homomorphism
$$
\Phi : (G, \circ) \to \Aut(G, \cdot) \times \Aut(G, \cdot)
$$
given by
$$\Phi(a) = (\lambda_a, h_a),$$
where $h_a(b) = a \cdot \lambda_a(b) \cdot a^{-1}$.
\end{enumerate}
\end{proposition}

\begin{proof}
For $x, y \in \Soc(G)$, we have 
$$
x \circ y = x \cdot \lambda_x(y) = x \cdot \id(y) = x \cdot y.
$$
Hence, the two group operations coincide on $\Soc(G)$. Let $x, y \in \Soc(G)$ and $a, b \in G$. We have $\lambda_{x \circ y} = \lambda_x \, \lambda_{y} = \id$ and
$$
a \cdot (x \circ y) \cdot a^{-1} = a \cdot x \cdot y \cdot a^{-1} = x \cdot y = x \circ y.
$$
This implies that $x \circ y \in \Soc(G)$. We also have $\lambda_{\bar{x}} = \lambda_{x}^{-1} = \id$ and
$$
a \cdot \bar{x} \cdot a^{-1} = a \cdot \lambda_{\bar{x}}(x)^{-1} \cdot a^{-1} = a \cdot x^{-1} \cdot a^{-1} = x^{-1} = \bar{x}.
$$
This shows that $\bar{x} \in \Soc(G)$, and hence $\Soc(G)$  is a subgroup of $(G, \circ)$. Moreover, we have
$$
\lambda_{b \circ x \circ \bar{b}} = \lambda_b \lambda_x \lambda_{\bar{b}} = \lambda_b \,  \id \, \lambda_{\bar{b}} = \id.
$$
Further, since
\begin{eqnarray*}
\bar{b} \circ \big(a \cdot  (b \circ x \circ \bar{b}) \cdot a^{-1} \big) &=& (\bar{b} \circ a) \cdot (\bar{b})^{-1} \cdot  (x \circ \bar{b})  \cdot (\bar{b})^{-1} \cdot (\bar{b} \circ a^{-1}) \\
&=& (\bar{b} \circ a) \cdot (\bar{b})^{-1} \cdot x \cdot (\bar{b} \cdot a^{-1}) \\
&=& x \cdot (\bar{b} \circ a) \cdot (\bar{b})^{-1} \cdot  (\bar{b} \circ a^{-1})\\
& =& x \cdot \bar{b} \circ (a \cdot a^{-1}) \\
&= & x \cdot \bar{b} \\
&=& x \circ \bar{b},
\end{eqnarray*}
we obtain
$$
a \cdot (b \circ x \circ \bar{b}) \cdot a^{-1} = b \circ x \circ \bar{b}.
$$
Hence, $\Soc(G)$ is a normal subgroup  of $(G, \circ)$. Further,
$$
\lambda_b(x) = b^{-1} \cdot  (b \circ x \circ \bar{b}) \cdot b =  b \circ x \circ \bar{b},
$$
and hence $\Soc(G)$  is invariant under $\lambda_b$ for each $b \in G$. Being closed under $\circ$ implies that $\Soc(G)$ is closed under $\cdot$, and hence it is a normal subgroup of $(G, \cdot)$. Finally, since the operations $\circ$ and $\cdot$  coincide in $\Soc(G)$, it  is a trivial brace, which proves assertion (1).
\para 

By Propositions \ref{lambda-map} and \ref{gamma-map},  the map $\varphi$ is well-defined and is a group homomorphism. Recall that
$\gamma_a(b) = \big((\bar{a} \circ \bar{b} \big) \cdot (\bar{a})^{-1} )^{\circ(-1)}$. Hence, if $\lambda_a = \id$, then 
$$
\gamma_a(b) = \big((\bar{a} \circ \bar{b} \big) \cdot (\bar{a})^{-1} \big)^{\circ(-1)} = 
\big((\bar{a} \cdot \bar{b} \big) \cdot (\bar{a})^{-1} \big)^{\circ(-1)} = \big(\bar{a} \cdot \bar{b}  \cdot (\bar{a})^{-1} \big)^{\circ(-1)}. 
$$
Thus, $\gamma_a(b) = b$ if and only $\bar{a} \cdot \bar{b} = \bar{b} \cdot \bar{a}$. Hence, we have
$$
\ker(\varphi) = \{ a \in G~|~\lambda_a = \gamma_a = \id \} = 
\{ a \in G~|~\lambda_a = \id~\textrm{and}~b \cdot a = a \cdot b~\mbox{for any}~ b \in G \} = \Soc(G),$$
proving assertion (2).
\para 

The proof of assertion (3) is analogous to (2) and follows from Proposition \ref{theta-map}.
$\blacksquare$ 
\end{proof}

\bigskip 
\bigskip


\subsection{Skew braces via exact factorisations of groups} \label{ExactFact}

We begin with the following definition.

\begin{definition}
A \index{factorization of a group}{factorization} of a group $G$ is the decomposition $$G = H L =\{ h l \, \mid \, h \in H ~~ \textrm{and}~ l \in L \}$$ of $G$ by two subgroups $H$ and $L$.  The factorization  is said to be \index{proper factorization}{proper} if $H$ and $L$ are proper subgroups of $G$, and is called \index{exact factorisation}{exact} if $H \cap L = 1$.
\end{definition}

For an exact factorization $G = H L$, it is not necessary that any of  the subgroups $H$ or $L$ is normal in $G$.  For example, the alternating group $A_5$ does not have proper normal subgroups, but admits a proper exact factorisation $A_5 = A_4 \big\langle (12345) \big\rangle$ by $A_4$ and a cyclic subgroup of order 5. 
\para

We have the following result \cite[Theorem 2.3]{MR3763907}.

\begin{theorem}  \label{skew left brace exact factor}
Let $(G, \cdot)$ be a group with an exact factorisation $G = H L$. Then $(G, \cdot, \circ)$ is a skew left brace, where
\begin{equation} \label{efeqn1}
a \circ b= h \cdot b \cdot l 
\end{equation}
for $a = h \cdot l, b\in HL$. Further, the group $(G, \circ)$ is isomorphic to the direct product  $H \times L$.
\end{theorem}

\begin{proof}
Let $\phi: H \times L \to G$ be the bijection given by $\phi(h, l)=h \cdot l^{-1}$. For elements $a=h \cdot l$ and $b=h_1 \cdot l_1$ of $G$, we see that
\begin{eqnarray*}
a \circ b &=& h \cdot b \cdot l\\
&=& h \cdot h_1\cdot l_1\cdot l\\
&=& \phi(h \cdot h_1, l^{-1} \cdot l_1^{-1})\\
&=& \phi \big((h, l^{-1})( h_1, l_1^{-1}) \big)\\
&=& \phi \big(\phi^{-1}(h\cdot l) \phi^{-1}(h_1\cdot l_1) \big)\\
&=& \phi \big(\phi^{-1}(a) \phi^{-1}(b) \big).
\end{eqnarray*}
Since $\phi$ is a bijection, it follows that $(G, \circ)$ is a group isomorphic to the direct product $ H \times L$. Let $a = h\cdot l \in HL$ and $b, c \in G$. Then we have
\begin{eqnarray*}
a \circ (b \cdot c) &=& h \cdot (b\cdot  c) \cdot l\\
&=& h \cdot b \cdot l \cdot (l^{-1} \cdot h^{-1}) \cdot h \cdot c \cdot l\\
&=& (h \cdot b \cdot l)\cdot  a^{-1} \cdot (h \cdot c \cdot l)\\
&=& (a \circ b) \cdot a^{-1} \cdot (a \circ c),
\end{eqnarray*}
and hence $(G, \cdot, \circ)$ is a skew left brace.  $\blacksquare$ 
\end{proof}

The following example gives an illustration of the preceding result \cite[Example 2.4]{MR3763907}.

\begin{example}
The general linear group $\GL(n, \mathbb{C})$ admits an exact factorization through subgroups $\Uu(n, \mathbb{C})$ and $\Tt(n, \mathbb{C})$, where 
$\Uu(n, \mathbb{C})$ is the unitary group and $\Tt(n, \mathbb{C})$ is the group of upper triangular matrices with positive diagonal entries. Thus, there exists a skew left brace $(G, \cdot, \circ)$ with
$(G, \cdot) \cong  \GL(n, \mathbb{C})$ and $(G, \circ) \cong \Uu(n, \mathbb{C}) \times \Tt(n, \mathbb{C})$.
\end{example}

\begin{proposition}\label{free product exact factorization}
A free product of non-trivial groups admits an exact factorization. 
\end{proposition}

\begin{proof}
Let $C$ and  $B$ be non-trivial groups and $G = C \ast B$  their free product. Let $D$ be the kernel of the natural homomorphism  $C \ast B \to C \times B$. Taking  $A=\langle  C, D  \rangle$, a direct check gives
$$
G=AB \quad \textrm{and}\quad A\cap B=1,
$$
and hence  $G$ admits an exact factorization. $\blacksquare$ 
\end{proof}
\para 

We have the following result  \cite[Proposition 5.11]{MR4346001}.

\begin{proposition}\label{exact factorisation of free group} 
There exists a skew  left brace $(G, \cdot, \circ)$ such that $(G, \cdot)$ is finitely generated, but $(G, \circ)$ is not.
\end{proposition}

\begin{proof}
Let $F_n=\langle x_1,\ldots,x_n \rangle$ be the free group of rank $n \ge 2$. Take  $C = F_{n-1} = \langle x_1, \ldots, x_{n-1} \rangle$ and $B = \langle  x_n \rangle$. Then $F_n = C \ast B$ and $A = \big\langle C, [x_n, F_n] \big\rangle$. By Proposition \ref{free product exact factorization}, we have
$$
F_n=AB \quad \textrm{and}\quad A \cap B=1,
$$
and hence $F_n$ admits an exact factorization. It follows from Theorem \ref{skew left brace exact factor} that $(F_n, \cdot, \circ)$ is a skew left brace with  $\circ$ defined by \eqref{efeqn1} and $(F_n, \circ) \cong A \times B$. Since the operations $\circ$ and $\cdot$ coincide on $A$ as well as on  $B$, it follows that $A\cong F_{\infty}$, the free group of countably infinite rank.  $\blacksquare$ 
\end{proof}

\begin{definition}
A \index{matched pair}{\it matched pair} of groups is a quadruple $(H, L, \rightharpoonup, \leftharpoonup)$, where $H, L$ are groups and $\rightharpoonup, \leftharpoonup$ are two actions
$$
L \overset{\leftharpoonup}{\longleftarrow} L \times H  \quad \textrm{and} \quad  L \times H  \overset{\rightharpoonup}{\longrightarrow} H
$$
satisfying the conditions
$$
l \rightharpoonup (h h') = (l \rightharpoonup  h) \big( (l \leftharpoonup h) \rightharpoonup  h' \big) \quad \textrm{and} \quad
(l l') \leftharpoonup h = \big(l \leftharpoonup (l' \rightharpoonup h) \big)(l' \leftharpoonup h)
$$
for all $h, h' \in H$ and $l, l' \in L$. Here, the notation $l \rightharpoonup h$ means that $l$ acts on $h$ from the left, and similarly, $l \leftharpoonup h$ means that $h$ acts on $l$ from the right.
\end{definition}

If $(H, L, \rightharpoonup, \leftharpoonup)$ is a matched pair of groups, then $H \times L$ is a group with multiplication
$$
(h, l) (h', l') = \big(h(l \rightharpoonup h'), (l \leftharpoonup h')l'\big)
$$
for $h, h' \in H$ and $l,  l' \in L$. The inverse of $(h, l)$ in this group is
$$
(h, l)^{-1} = \big(l^{-1} \rightharpoonup h^{-1}, (l \leftharpoonup (l^{-1} \rightharpoonup h^{-1}))^{-1}\big).
$$
This group is denoted by $H \bowtie L$, and known as the \index{biproduct}{\it biproduct} of $H$ and $L$. Note that the biproduct $H \bowtie L$  admits an exact factorization through its subgroups $H \bowtie 1 \cong H$ and $1 \bowtie L \cong L$. As a consequence of Theorem \ref{skew left brace exact factor},  we obtain the following result \cite[Corollary 2.7]{MR3763907}.

\begin{corollary}\label{skew left brace match product}
Let $(H, L, \rightharpoonup, \leftharpoonup)$ be a matched pair of groups. Then the biproduct $H \bowtie L$ is a skew left brace with
$$
(h, l) (h', l') = \big(h (l \rightharpoonup h'), (l \leftharpoonup h') l' \big) \quad \textrm{and} \quad  (h, l) \circ (h', l') = (h h', l l')
$$
for $h, h' \in H$ and $l,  l' \in L$.
\end{corollary}

We use exact factorization to prove that the knot group of a knot admits a non-trivial skew left brace structure. 

\begin{proposition}
 Let $K$ be a  knot  in $\mathbb{S}^3$ and $G(K) = \pi_1(\mathbb{S}^3 \setminus K)$ its knot group.  Then there is a non-trivial skew left brace structure on $G(K)$.
\end{proposition}

\begin{proof}
We set $G = G(K)$. Since $G / [G, G] \cong \mathbb{Z}$,  there exists  a factorization $G = H L$ with $H=[G, G]$, $L \cong \mathbb{Z}$ and $H \cap L = 1$. Hence, by Theorem \ref{skew left brace exact factor},  we can define a binary operation $\circ$ on $G$ such that $(G, \cdot, \circ)$ is a non-trivial skew left brace.
$\blacksquare$ 
\end{proof}

Another application of exact factorization  yields a  skew left brace whose additive group is solvable but not nilpotent, while its multiplicative group is abelian \cite[Proposition 5.13]{MR4346001}.

\begin{proposition}
There is a skew left brace $(G, \cdot, \circ)$ such that $(G, \cdot)$ is a metabelian group which is not nilpotent, while $(G,  \circ)$ is a free abelian group of countably infinite rank.
\end{proposition}

\begin{proof}
Let $G$ be the group with the following presentation
$$
G=\big\langle\,  x, y_i, \, i\in \mathbb{Z}   \, \mid \,    xy_i x^{-1}=y_{i+1} ~~ \textrm{and} ~~  [y_i,y_j]=1  ~~ \textrm{for all} ~~ i,j\in \mathbb{Z}\,\big\rangle.
$$
It follows that $G$ is the \index{wreath product}{wreath product}  $\mathbb{Z} \wr \mathbb{Z}$ of two infinite cyclic groups, and hence is  metabelian. Consider the subgroups
$$
H = \big\langle\,  y_i, \,i\in \mathbb{Z}   \, \mid \,     [y_i,y_j]=1 ~~ \textrm{for all} ~~ i,j\in \mathbb{Z}\,\big\rangle \cong \mathbb{Z}^{\infty}   ~~ \textrm{and} ~~ L = \langle x \rangle \cong \mathbb{Z}
$$
of $G$. Then $G = H \rtimes L$ is an exact factorization of $G$. Hence, Theorem \ref{skew left brace exact factor} gives a binary operation $\circ$ on $G$ such that $(G, \cdot, \circ)$ is a  skew left brace and
$$
(G, \circ) \cong H \times L \cong \mathbb{Z}^{\infty}  \times \mathbb{Z} \cong \mathbb{Z}^{\infty}.
$$
Further, $\gamma_2(G, \cdot)=\gamma_3(G, \cdot)$ is non-trivial abelian, and hence $(G, \cdot)$ is metabelian but not nilpotent. $\blacksquare$
\end{proof}
\bigskip 
\bigskip


\subsection{Skew braces, bijective 1-cocycles and regular subgroups}
Let $G$ and $\Gamma$ be groups such that the map $\Gamma \times G \to G$ given by $(\gamma, a) \mapsto \gamma \sbullet a$,
is a left-action of $\Gamma$ on $G$ by automorphisms. A bijective \index{1-cocycle}{1-cocycle} or \index{crossed homomorphism}{a crossed homomorphism} is a bijective map $\pi : \Gamma \to G$ such that
$$\pi(x y) = \pi(x) \, \big(x \sbullet \pi(y) \big)$$
for all $x,y\in \Gamma$. 
\para

The following result gives an intimate connection between bijective 1-cocycles and skew braces on a given group \cite[Proposition 1.11]{MR3647970}.

\begin{proposition}\label{bijective 1-cocycle and skew brace}
Let $(G, \cdot)$ be a group. Then there is a one-to-one correspondence between bijective 1-cocycles with values in $(G, \cdot)$ and skew left braces $(G, \cdot, \circ)$.
\end{proposition}

\begin{proof}
Let $(G, \cdot)$ and $\Gamma$ be groups such that $\Gamma \times G \to G$ given by $(\gamma, a) \mapsto \gamma \sbullet a$,
is a left-action of $\Gamma$ on $G$. Let $\pi : \Gamma \to G$ be a bijective 1-cocycle with respect to this action. Then the binary operation $\circ$ by
\begin{equation} \label{1-cocycle condtion}
a \circ b = \pi\big(\pi^{-1}(a) \,\pi^{-1}(b) \big)
\end{equation}
for $a,b \in G$, gives rise to a group $(G, \circ)$. Further, for $a, b, c \in G$, we have
\begin{eqnarray*}
a \circ (b \cdot c) &=& \pi \big(\pi^{-1}(a)\, \pi^{-1}(b \cdot c)\big)\\
&=& a \cdot \big(\pi^{-1}(a) \sbullet (b \cdot c) \big)\\
&=& a \cdot \big((\pi^{-1}(a) \sbullet b) \cdot (\pi^{-1}(a)\sbullet c)\big)\\
&=& a \cdot (\pi^{-1}(a) \sbullet b) \cdot a^{-1} \cdot a \cdot (\pi^{-1}(a)\sbullet c)\\
&=& \pi \big(\pi^{-1}(a) \pi^{-1}(b) \big) \cdot a^{-1} \cdot \pi \big(\pi^{-1}(a) \pi^{-1}(c) \big)\\
&=& (a \circ b) \cdot a^{-1} \cdot (a \circ c),
\end{eqnarray*}
and hence $(G, \cdot, \circ)$ is a~skew left brace.
\para
Conversely, let $(G, \cdot, \circ)$  be a skew left brace. Taking $\Gamma = (G, \circ)$, the map $\lambda:\Gamma \to \Aut(G, \cdot)$ is an action of $\Gamma$ on $(G, \cdot)$ by automorphisms, which we write as $a \sbullet b= \lambda_a(b)$ for $a \in \Gamma$ and $b \in G$. Taking $\pi:\Gamma \to G$ to be the identity map, we see that $\pi (a \circ b)= a \circ b= a \cdot \lambda_a(b) =a \cdot (a \sbullet b )= \pi(a) \cdot \big(a \sbullet \pi(b) \big)$ for all $a, b \in \Gamma$. Hence, $\pi : \Gamma \to G$ is a bijective 1-cocycle. $\blacksquare$ 
\end{proof}
\bigskip

The \index{holomorph}{\it holomorph} of a group $G$ is defined to be the group $\Hol(G) = \Aut(G) \ltimes G$, in which the product is given by
$$
(f,a)(g,b)=\big(fg,~af(b) \big)
$$
for all $a,b \in G$ and $f,g \in \Aut(G)$. Each subgroup $H$ of $\Hol(G)$ admits a left-action on $G$ given by
$$
(f,a) \sbullet b = af(b)
$$
for $a,b \in G$ and $f \in \Aut(G)$. A subgroup $H$ of $\Hol(G)$ is said to be \index{regular subgroup}{\it regular} if the action of $H$ on $G$ is free and transitive. Regularity is equivalent to the fact that for each $a \in G$ there exists a unique $(f,x) \in H$ such that $xf(a)=1$.  Let $\pi_2 : \Hol(G) \to G$ denote the projection onto the second factor, that is, $\pi_2((f, a) )=a$ for all $(f, a) \in \Hol(G)$. Then we have the following result.

\begin{proposition}\label{reg sub to group bijection}
Let $G$ be a group and $H$ a regular subgroup of $\Hol(G)$. Then the restriction $\pi_2|_H : H \to G$ is bijective.
\end{proposition}

\begin{proof}
Let $(f, a),(g, b) \in H$ be such that $\pi_2 \big((f, a) \big) = \pi_2 \big((g, b) \big)$. This gives $a = b$. Since $H$ is a subgroup, both the elements $(f, a)^{-1}= \big(f^{-1}, f^{-1}(a^{-1}) \big)$ and $(g, b)^{-1}= \big(g^{-1}, g^{-1}(b^{-1}) \big)$ lie in $H$. Since $f^{-1}(a^{-1})f^{-1}(a)=g^{-1}(b^{-1})g^{-1}(b)=1$ and $H$ is a regular subgroup, it follows that $f=g$, and hence $\pi_2$ is injective. For surjectivity, let $a \in G$. The regularity of $H$ implies
the existence of an automorphism $f \in \Aut(G)$ such that $\big(f,f(a^{-1}) \big) \in H$. Then $(f^{-1}, a) \in H$ and the proposition  follows. $\blacksquare$ 
\end{proof}

The following theorem  provides a connection between skew left braces and regular subgroups of the holomorph \cite[Theorem 4.2]{MR3647970}. The result goes back to Bachiller for left braces \cite{MR3465351} and also generalises a result of Catino and Rizzo \cite{MR2486886}.

\begin{theorem}\label{gv2017}
Let $(G, \cdot)$ be a group. Then the following assertions  hold:
\begin{enumerate}
\item  If $(G, \cdot, \circ)$ is a skew left brace, then $ \{(\lambda_a, a) \, \mid \, a \in G  \}$
is a regular subgroup of $\Hol(G)$, where $\lambda_a(b) = a^{-1} \cdot (a \circ b)$ for all $a, b \in G$.
\item If $H$ is a regular subgroup of $\Hol (G, \cdot)$, then $(G, \cdot, \circ)$ is a skew left brace with $(G, \circ) \cong H$, where $a \circ b=a f(b)$ for some $f \in \Aut(G, \cdot)$ such that  $(\pi_2|_H)^{-1}(a)=(f, a)$.
\end{enumerate}
\end{theorem}

\begin{proof}
Let $(G, \cdot, \circ)$ be a skew left brace. Since $\lambda: (G, \circ) \to \Aut(G, \cdot)$ is a group homomorphism and $a \cdot \lambda_a(b)=a\circ b$ for all $a,b\in G$, it follows that $\{(\lambda_a,a)  \, \mid \,  a\in G\}$ is a subgroup of $\Hol(G,\cdot)$. Further, since $(G,\circ)$ is a group, the regularlity of $H$ also follows.
\para

Now, let $H$ be a regular subgroup of $\Hol(G, \cdot)$. By Proposition \ref{reg sub to group bijection}, $\pi_2|_H:H \to G$ is a bijection. Hence, we can transport the group structure on $H$ onto $G$  by defining 
$$ a\circ b=\pi_2|_H\big((\pi_2|_H)^{-1}(a)\, (\pi_2|_H)^{-1}(b)\big)=a \cdot f(b),$$
where $a,b\in G$ and $(\pi_2|_H)^{-1}(a)=(f,a)\in H$. Since
$$ a\circ (b \cdot c)=a \cdot f(b \cdot c)=a \cdot f(b) \cdot f(c)=a\cdot f(b)\cdot a^{-1}\cdot a\cdot f(c)=(a\circ b)\cdot a^{-1}\cdot (a\circ c)$$
for all $a,b,c\in G$, it follows that $(G, \cdot, \circ)$ is a skew left brace. $\blacksquare$ 
    \end{proof}

\begin{proposition}
Let $(G, \cdot)$ be a group.  Then there is a bijective correspondence between skew left brace structures on $(G, \cdot)$ and regular subgroups of
	$\Hol(G, \cdot)$. Moreover, isomorphic skew left braces structures on $(G, \cdot)$ correspond to conjugate subgroups of $\Hol(G, \cdot)$ by elements of $\Aut(G, \cdot)$.
\end{proposition}

\begin{proof}
Suppose that the group $(G, \cdot)$ has two skew left brace structures, for which the multiplicative groups are given by $(a,b)\mapsto a\circ b$ and $(a,b)\mapsto a\times b$. Let $\phi\in\Aut(G, \cdot)$ such that $\phi(a\circ b)=\phi(a)\times\phi(b)$ for	 all $a,b\in G$. Let $\lambda_a(b)=a^{-1}\cdot (a\circ b)$ and $\mu_a(b)=a^{-1}\cdot (a\times b)$ for all $a,b\in G$. We claim that the regular subgroups $\{(\lambda_a,a) \, \mid \, a\in G$\} and 	$\{(\mu_a,a) \, \mid\, a\in G\}$ are conjugate by the element $\phi$. Since
$$    \phi\lambda_a\phi^{-1}(b)=\phi \big(a^{-1}\cdot (a\circ\phi^{-1}(b))\big)=\phi(a)^{-1}\cdot \big(\phi(a)\times b \big)=\mu_{\phi(a)}(b), $$
we see that $(\phi, 1)(\lambda_a,a)(\phi, 1)^{-1}=\big(\mu_{\phi(a)},\phi(a)\big)$, and hence the claim follows.
\para

Conversely, let $H$ and $K$ be regular subgroups of $\Hol(G, \cdot)$ such that there exists $\phi\in\Aut(G, \cdot)$ with $(\phi, 1)H(\phi, 1)^{-1}=K$. For $a \in G$, let $(f,a)=(\pi_2|_H)^{-1}(a)\in H$ and $(g,a)=(\pi_2|_K)^{-1}(a)\in K$. Write  $a\circ b=a \cdot f(b)$ and 
$a\times b=a \cdot g(b)$ for $a, b \in G$. It follows from Theorem \ref{gv2017} that both $(G, \cdot, \circ)$ and $(G, \cdot, \times)$ are skew left braces.
Since  $(\phi, 1)(f,a)(\phi, 1)^{-1}= \big(\phi f\phi^{-1},\phi(a) \big)\in K$, it follows that  $ (\pi_2|_K)^{-1}(\phi(a))= \big(\phi f\phi^{-1}, \phi(a)\big)$. This gives
$$ \phi(a)\times\phi(b)=\phi(a) \cdot  (\phi f\phi^{-1})\big(\phi(b)\big)=\phi(a) \cdot \phi \big(f(b)\big)=\phi \big(a \cdot  f(b)\big)=\phi(a\circ b),$$
and hence the skew left braces corresponding to $H$ and $K$ are isomorphic. $\blacksquare$ 
\end{proof}
\bigskip
\bigskip


\section{Solutions arising from skew  braces}

In this section, we study solutions to the Yang--Baxter equation that arise from skew left braces. Recall from Definition \ref{def braiding operator} the notion of a braiding operator on a group. According to Theorem \ref{conjugate braiding thm}, every braiding operator yields a non-degenerate bijective solution to the Yang--Baxter equation. The following theorem is due to Guarnieri and Vendramin \cite[Theorem 3.1]{MR3647970}. 

\begin{theorem} \label{solution from a skew brace}
Let $(G,\cdot, \circ)$ be a~skew left brace and $r_G:G \times G \to G \times G$ given by
\begin{equation} \label{BraceToYBE}
r_G(a, b) = \left(\lambda_a(b), ~\lambda^{-1}_{\lambda_a(b)}\big((a \circ b)^{-1} \cdot a\cdot (a \circ b)\big) \right)= \big(\lambda_a(b), ~ \overline{\lambda_a(b)} \circ (a \circ b)\big)
\end{equation}
for $a, b \in G$. Then $(G, r_G)$ is a non-degenerate bijective solution to the Yang--Baxter equation.
\para 
If $(G,\cdot, \circ)$ is a left brace (that is, $(G, \cdot)$ is abelian), then
\begin{equation}
r_G(a, b) = \left(\lambda_a(b), ~\lambda^{-1}_{\lambda_a(b)}(a) \right).
\end{equation}
\para 
Moreover, $r_G$ is involutive if and only if $(G, \cdot)$ is abelian.
\end{theorem}

\begin{proof} 
 For convenience of notation, we denote $r_G$ by $r$.  By Theorem \ref{conjugate braiding thm}, every braiding operator on a group is a non-degenerate bijective solution to the Yang--Baxter equation. Thus, it is enough to prove that $r$ is a braiding operator on $(G, \circ)$. Since $\lambda_a^{-1}(b) = \bar{a} \circ (a \cdot b)$ for all $a, b \in G$, we have 
$$
\lambda_{\lambda_a(b)}^{-1}\big((a \circ b)^{-1} \cdot a \cdot (a \circ b)\big) = \overline{\lambda_a(b)} \circ \big(\lambda_a(b) \cdot (a \circ b)^{-1} \cdot a\cdot  (a \circ b)\big) = \overline{\lambda_a(b)} \circ (a \circ b)
$$
 for all $a, b \in G$. Thus, if $m$ is the multiplication in $(G, \circ)$, then $m \, r(a, b) = a \circ b$  for all $a, b \in G$. Clearly, $r(a, 1) = (1, a)$ and $r(1, a) = (a, 1)$ for all $a \in G$. Since $\lambda : (G, \circ) \to \Aut(G, \cdot)$ is a group homomorphism, we have
\begin{eqnarray*}
(\id \times m)\, r_{12} \,r_{23} \,(a, b, c) &=& (\id \times m) \, r_{12}\, \big(a, ~\lambda_b(c),  ~\overline{\lambda_b(c)}  \circ b \circ c \big)\\
&=& (\id \times m) \, \big(\lambda_a\lambda_b(c), ~\overline{\lambda_a\lambda_b(c)}  \circ a \circ \lambda_b(c),~  \overline{\lambda_b(c)}  \circ b \circ c \big)\\ 
&=& \big(\lambda_a\lambda_b(c),~ \overline{\lambda_a\lambda_b(c)}  \circ a \circ b \circ c \big)\\
&=& r(a \circ b, ~c).
\end{eqnarray*}
For $a, b, c \in G$, since $b \circ c=b \cdot \lambda_b(c)$, we have
$$
\lambda_a (b \circ c) = \lambda_a (b) \cdot \lambda_{a \circ b} (c).
$$
This gives
$$\lambda_a (b) \circ \lambda_{\overline{\lambda_a(b)} \circ a \circ b} (c) = \lambda_a (b) \circ \lambda^{-1}_{\lambda_{a}(b)} \lambda_{a}\lambda_{b}(c)= \lambda_a (b) \cdot \lambda_{a \circ b} (c)=\lambda_{a}(b  \circ c)$$
for $a, b, c \in G$. Thus, 
\begin{eqnarray*}
(m \times \id) \,r_{23}\, r_{12} \,(a, b, c) 
&=&  (m \times \id)\, r_{23}\, \big(\lambda_a (b),~ \overline{\lambda_a (b)} \circ a \circ b, ~c \big) \\
&=&   (m \times \id) \, \big(\lambda_a (b),  ~\lambda_{\overline{\lambda_a (b)} \circ a \circ b} (c), ~\overline{\lambda_{\overline{\lambda_a (b)} \circ a \circ b}(c)} \circ \overline{\lambda_a (b)} \circ a \circ b \circ c \big) \\
&=&  \big(\lambda_a (b) \circ  \lambda_{\overline{\lambda_a (b)} \circ a \circ b} (c), ~\overline{\lambda_{\overline{\lambda_a (b)} \circ a \circ b}(c)} \circ
\overline{\lambda_a (b)} \circ a \circ b \circ c \big)\\
&=&  \big(\lambda_a (b\circ c), ~\overline{\lambda_a (b\circ c)} \circ a \circ b \circ c \big) \\
&=&  r(a, ~b \circ c)
\end{eqnarray*}
for all $a, b, c \in G$. Finally, we prove that $r$ is bijective. Suppose that  $r(a, b)=r(c, d)$. This gives $\lambda_a(b)=  \lambda_c(d)$ and $\overline{\lambda_a(b)} \circ (a \circ b)= \overline{\lambda_c(d)} \circ (c \circ d)$, from which we deduce that $a \circ b= c \circ d$. Now, $a^{-1} \cdot (a \circ b)= \lambda_a(b)=  \lambda_c(d)= c^{-1} \cdot (c \circ d)$ implies that $a=c$. Plugging this back into  $\lambda_a(b)=  \lambda_c(d)$ implies that $b=d$, and hence $r$ is injective. For surjectivity of $r$, if $(x, y) \in X \times X$, then a direct check  shows that 
$$r\big((x \circ y) \cdot x^{-1},~ (\,\overline{(x \circ y) \cdot x^{-1})} \circ x \circ y \big)= (x, y),$$
and hence $r$  is a braiding operator on $(G, \circ)$. The remaining two assertions are immediate, which completes the proof. $\blacksquare$ 
\end{proof} 

\begin{remark}
{\rm 
Bachiller, Ced\'{o}, and Jespers showed in \cite{MR3527540} that non-degenerate involutive solutions can be explicitly generated using left braces. Subsequently, Bachiller further established in \cite{MR3835326} that non-degenerate bijective solutions can be derived from skew left braces. This reduces the problem of classifying all non-degenerate bijective solutions to that of classifying skew left braces.
}
\end{remark}

\begin{example}
{\rm 
Let $(G, \cdot, \circ)$ be a skew left brace such that the map $\lambda$ is trivial, that is, $\lambda_a = \id$ for all $a \in G$. The identity $\lambda_a(b) = a^{-1} \cdot (a \circ b) = b$ implies that $a \circ b = a \cdot b$, and hence $(G, \cdot, \circ)$ is a trivial skew left brace. Then the formula (\ref{BraceToYBE}) has the form
\begin{equation} \label{trivBr}
s_G(a, b) = (b, ~b^{-1} \cdot a \cdot  b)
\end{equation}
for $a, b \in G$. Note that this is a solution of quandle type on the group $(G, \cdot)$.}
\end{example}

\begin{example}{\rm 
Let $G$ and $K$ be groups and $\alpha : G \to \Aut(K)$  a group homomorphism.  Consider the skew left brace structure on $G \times K$ as in Example \eqref{examples of solutions}(5). Then the solution to the Yang--Baxter equation corresponding to the skew left brace $(G \times K, \cdot, \circ)$ has braiding
$$
r : (G \times K) \times  (G \times K)  \to (G \times K) \times  (G \times K)
$$
is given by
$$
r\big((a_1, b_1), ~(a_2, b_2)\big) = \left( \big(a_2, \alpha_{a_1}(b_2) \big), ~\big(a_2^{-1} a_1 a_2, b_2^{-1} \alpha_{a_2^{-1}}(b_1 \alpha_{a_1}(b_2) ) \big) \right)
$$
for $a_1, a_2 \in G$ and $b_1, b_2 \in K$.}
\end{example}
\para 

Next, we show that the solution to the Yang--Baxter equation arising from a skew left brace as in \eqref{BraceToYBE}  is conjugate to the solution of the form \eqref{trivBr} \cite[Proposition 3.7]{MR3647970}.

\begin{proposition} \label{GV}
Let $(G,\cdot,\circ)$ be a skew left brace. Then the solution to the Yang--Baxter equation given by \eqref{BraceToYBE} is conjugate to the solution of quandle type on the group $(G, \cdot)$ given by \eqref{trivBr}.
\end{proposition}

\begin{proof}
Let $T : G\times G  \to G\times G $ be defined by
$$
T(a, b) = \big(a, ~\lambda_a(b) \big)
$$
for $a, b \in G$. Then $T$ is invertible with inverse
$$
T^{-1}(a, b) = \big(a, ~\lambda^{-1}_a(b) \big)
$$
for $a, b \in G$. Let $r_G$ be the solution given by ~\eqref{BraceToYBE} and $s_G$ the solution of quandle type given by \eqref{trivBr}.
A direct calculation gives
\begin{eqnarray*}
T \, r_G \, (a, b) &=& T \,\big(\lambda_a(b), ~\lambda^{-1}_{\lambda_a(b)}((a \circ b)^{-1} \cdot a\cdot (a \circ b))\big)\\
&=& \big(\lambda_a(b),~ (a \circ b)^{-1} \cdot a\cdot (a \circ b)\big)\\
&=& \big(\lambda_a(b), ~\lambda_a(b)^{-1} \cdot a \cdot \lambda_a(b) \big)\\ 
&=& s_G \, \big(a, ~\lambda_a(b)\big)\\ 
 &=& s_G \, T \,(a, b)
\end{eqnarray*}
for all $a, b \in G$, which proves the proposition. $\blacksquare$  
\end{proof}          

\begin{remark}
{\rm The preceding proposition also follows from Proposition \ref{SolTh} by taking $\sigma_a(b)=\lambda_a(b)$ and $\tau_b(a)=\overline{\lambda_a(b)} \circ (a \circ b)$ for $a,  b \in G$.
}    
\end{remark}
\para

Recall that, to each solution $(X,r)$ to the Yang--Baxter equation, we can associate its structure group $G(X, r)$ and its derived structure group $A(X, r)$. Also, we have natural maps $i_X: X \to G(X, r)$ and $j_X: X \to A(X, r)$. By Theorem \ref{action of str group on der str group}, if $(X, r)$ is a non-degenerate solution, then the left-action $G(X, r) \to \Sigma_X$, given by $x \mapsto \sigma_x$, lifts uniquely to a left-action $$\lambda:G(X, r) \to  \Aut\big(A(X, r)\big).$$ In other words, 
$$\lambda_{x}(y)=\sigma_x(y)$$
for all $x, y \in X$. The next result gives the structure  of a skew left brace on the structure group $G(X, r)$ of a non-degenerate solution $(X, r)$ (see \cite[Theorem 3.7]{MR3835326} and \cite[Theorem 9]{MR1769723}).

\begin{theorem} \label{skew brace on structure group}
Let $(X, r)$ be a  non-degenerate bijective solution to the Yang--Baxter equation and $G(X, r)$ its structure group with group operation $\circ$. Then there is another group operation $\cdot$ on the set $G(X, r)$ and a group isomorphism $\varphi : \big(G(X, r), \cdot \big) \to A(X, r)$ with $\varphi(x) = x$ for all $x \in X$, such that $\big(G(X, r), \cdot, \circ \big)$ is a skew left brace, where 
\begin{equation}\label{bijective 1-cocycle inverse}
a \cdot b = a \circ \varphi^{-1} \big(\lambda_a^{-1}(\varphi(b)) \big)    
\end{equation}
for $a, b\in G(X, r)$.
\end{theorem}

Note that \eqref{bijective 1-cocycle inverse} is equivalent to the map $$\varphi : \big(G(X, r), \circ \big) \to A(X, r)$$ being a bijective 1-cocycle. By Theorem \ref{univ property non-deg solution st group}, the structure group of a non-degenerate bijective solution to the Yang--Baxter equation satisfies a universal property, which can be rephrased in terms of skew left braces as follows  \cite[Theorem 3.9]{MR3647970}.

\begin{theorem}
Let $(X,r)$ be a non-degenerate bijective solution to the Yang--Baxter equation and $i_X  : X \to G(X, r)$ the natural map. Then the following  assertions hold:
\begin{enumerate}
\item There exists a unique skew left brace structure on the group $G(X, r)$ such that $i_X$ becomes a morphism of solutions. In other words, if $r_G$ is the solution corresponding to the skew left brace on the group $G(X, r)$, then
$$
(i_X \times i_X) \, r= r_G \,(i_X \times i_X).
$$
\item If $(A, \cdot, \circ)$ is any skew left brace with corresponding solution $r_A$ and $\varphi : X \to A$ is a morphism of solutions, that is,
$$
(\varphi \times \varphi)\, r = r_A \, (\varphi \times \varphi),
$$
then there exists a unique skew brace homomorphism $\tilde{\varphi} : G(X, r) \to A$ such that $\varphi = \tilde{\varphi} \,i_X$ and
$$
(\tilde{\varphi} \times \tilde{\varphi}) \,r_G = r_A \,(\tilde{\varphi} \times \tilde{\varphi}).
$$
\end{enumerate}
\end{theorem} 

\begin{remark}
{\rm 
It follows from Proposition \ref{involutive sol condition} that if  $(X, r)$ is an involutive solution, then
$$
\sigma_{\sigma_x(y)} \big(\tau_y(x) \big) = x \quad \textrm{and}\quad \tau_{\tau_y(x)} \big(\sigma_x(y) \big) = y
$$ 
for all $x, y \in X$. In addition, if $(X, r)$ is non-degenerate, then  $\tau_y(x) = \sigma^{-1}_{\sigma_x(y)}(x)$ and we can write
$$
r(x, y) = \big(\sigma_x(y),~ \sigma^{-1}_{\sigma_x(y)}(x)\big).
$$
Under these conditions, the derived structure group has the presentation 
$$
A(X, r) = \big\langle X  \, \mid \, x  \sigma_x(y) = \sigma_x(y) \sigma_{\sigma_x(y)}(\tau_y(x)) \big\rangle = \big\langle X  \, \mid \, x \sigma_x(y) = \sigma_x(y) x \big\rangle =
 \big\langle X  \, \mid \, x  y= y  x  \big\rangle, 
$$
and hence is a free abelian group generated by $X$. This implies that the natural map $j_X : X \to A(X, r)$ is injective. Theorem \ref{skew brace on structure group} gives an isomophism $\varphi : \big(G(X, r), \cdot \big) \to A(X, r)$ with $\varphi(x) = x$ for all $x \in X$. This implies  that for non-degenerate involutive solutions $(X, r)$, the map $i_X:X \to G(X, r)$ is injective. In other words, no two elements of $X$ coincide as generators of $G(X, r)$.
}    
\end{remark}
\bigskip
\bigskip


\chapter{Skew braces and associated actions} \label{chapter skew braces and solution of YBE}

\begin{quote} 
The structure of a skew left brace $(G, \cdot, \circ)$ is governed by its associated left-action $\lambda: (G, \circ) \to \Aut(G, \cdot)$. In this chapter, we focus on skew left braces for which the map $\lambda: (G, \cdot) \to \Aut(G, \cdot)$ is also a left-action and whose image is cyclic. These additional constraints have significant implications on the structure of the skew left brace.  We present several constructions of such skew left braces on free abelian groups and investigate the associated solutions to the Yang--Baxter equation that emerge from these structures.
\end{quote}
\bigskip

\section{Lambda-homomorphic skew braces} \label{lambda}

We begin with the following definition.

\begin{definition}
A skew left brace $(G, \cdot, \circ)$ is said to be \index{$\lambda$-homomorphic skew left brace}{$\lambda$-homomorphic} if the map $\lambda : (G, \cdot) \to \Aut(G, \cdot)$, defined  by  $\lambda_a(b) = a^{-1} \cdot  (a \circ b)$ for $a, b \in G$, is a group homomorphism.
\end{definition}

Given a group  $(G, \cdot)$, we desire to define a homomorphism $\lambda : (G, \cdot) \to \Aut(G, \cdot)$ such that $(G, \cdot, \circ)$ is a skew left brace, where  $a \circ b = a \cdot \lambda_a (b)$ for all $a, b \in G$. In this regard, we prove the following result \cite[Theorem 2.3]{MR4346001}.

\begin{theorem}  \label{t1}
Let $(G, \cdot)$ be a group and $\lambda : (G, \cdot) \to \Aut(G, \cdot)$ a group homomorphism. Then $ H_{\lambda} = \{ (\lambda_a, a) \,\mid \, a \in G \}$ is a subgroup of $\Hol(G)$ if and only if
\begin{equation} \label{inc}
[G, \lambda(G)] = \big\{ b^{-1} \cdot \lambda_a (b) \, \mid \, a, b \in G \big\} \subseteq \ker(\lambda).
\end{equation}
Moreover, if  $H_{\lambda}$ is a subgroup of $\Hol(G)$, then it is regular.
\end{theorem}

\begin{proof}
Suppose that  $H_{\lambda}$ is a subgroup of  $\Hol(G)$. Then, for $a, b \in G $, we have
$$ (\lambda_a, a)(\lambda_b, b)=\big(\lambda_a \lambda_b, ~a \cdot \lambda_{a} (b)\big) =\big(\lambda_{a \cdot b}, ~a \cdot \lambda_a (b)\big) \in H_{\lambda}.$$
This implies that
\begin{equation}
a \cdot \lambda_a(b)=a\cdot b,
\end{equation}
and hence
\begin{equation} \label{2}
\lambda \big(a \cdot \lambda_a(b) \big)=\lambda(a\cdot b).
\end{equation}
Since $\lambda$ is a group homomorphism, we get
\begin{equation} \label{3}
\lambda \big(\lambda_a(b)\big)=\lambda(b),
\end{equation}
which is equivalent to
\begin{equation} \label{4}
\lambda \big(b^{-1} \cdot \lambda_a(b)\big) = \id.
\end{equation}
This proves that $[G, \lambda(G)]  \subseteq \ker(\lambda)$.
\para

Conversely, suppose that $[G, \lambda(G)] \subseteq \ker(\lambda)$. Then  (\ref{4}), (\ref{3}) and (\ref{2}) hold, which implies that
$H_{\lambda}$ is closed under the product in $\Hol(G)$. It remains  to check that  $H_{\lambda}$ is closed under inversion in $\Hol(G)$. Note that, for   $(f, x)\in \Hol(G)$, we have  $(f, x)^{-1}= \big(f^{-1}, ~f^{-1}(x^{-1})\big)$, and hence
$$
(\lambda_a, a)^{-1}=\left( \lambda_a^{-1},~ \lambda_a^{-1}(a^{-1})    \right).
$$
By \eqref{3}, we have   $\lambda\left( \lambda_{a^{-1}}(a^{-1})  \right)=\lambda_{a^{-1}}$.
Since $\lambda$ is a homomorphism, we get
$$
(\lambda_a, a)^{-1}=\left( \lambda_{a^{-1}}, ~\lambda_{a^{-1}}(a^{-1})  \right)=
\left( \lambda ( \lambda_{a^{-1}}(a^{-1}) ), ~\lambda_{a^{-1}}(a^{-1})  \right) \in H_{\lambda},
$$
which proves that $H_{\lambda}$ is a subgroup of $\Hol(G)$.
\para 

Now, suppose that  $H_{\lambda}$ is a subgroup of $\Hol(G)$. By definition, the projection $\pi_2 : H_{\lambda} \to G$
given by $\pi_2 \big((\lambda_a, a)\big) = a$ is a bijection. Given $a \in G$, we see that $(\lambda_a, a)$ is the unique element such that  $(\lambda_a, a)\sbullet 1= a \cdot \lambda_a(1)= a$. Hence, the action of $H_{\lambda}$ on  $G$ is free and  transitive, and it is a regular subgroup of $\Hol(G)$. $\blacksquare$ 
\end{proof}

\begin{corollary} \label{c1}
Let $(G, \cdot)$ be a group and $\lambda : (G, \cdot) \to \Aut(G, \cdot)$ a group homomorphism. If  $G$ is non-trivial and $\lambda$ is injective, then $H_{\lambda}$ is not a group.
\end{corollary}

\begin{proof}
If $H_{\lambda}$ is a group, then injectivity of $\lambda$ together with \eqref{inc} gives  $\lambda_a = \id$ for all $a \in G$. But, this implies that   $\ker(\lambda)=G$,  a contradiction.  $\blacksquare$ 
\end{proof}
\para

 The following example shows that $H_{\lambda}$ can be a subgroup of $\Hol(G)$ for an injective anti-homomorphism $\lambda : G \to \Aut(G)$.

\begin{example}{\rm 
For a group $G$ with trivial center, consider the subset $H_\lambda= \{(\lambda_{a}, a) \, \mid \,  a \in   G \}$  of $\Hol(G)$, where $\lambda_{a}(b)=a^{-1}ba$ for all $b \in G$. Note that $\lambda : G \to \Aut(G)$ given by $a \mapsto \lambda_a$, is an injective anti-homomorphism. Since
$$
(\lambda_{a},a)(\lambda_{b},b)=\big(\lambda_{a}\lambda_{b},a \, \lambda_{a}(b)\big)=
(\lambda_{ba},ba),
$$
it follows that  $H$ is a subgroup of $\Hol(G)$. Further, the group operation  $\circ$ on the set  $G$ defined by
$$
a\circ b=a  \, \lambda_{a}(b)=ba,
$$
gives a  skew  left brace $(G, \cdot, \circ)$.}
\end{example}
\para

We introduce the following definition to describe the structure of a $\lambda$-homomorphic skew left brace.

\begin{definition}
A skew left brace $(G, \cdot, \circ )$ is said to be  \index{meta-trivial skew left brace}{\it meta-trivial}  if the following conditions hold:
\begin{enumerate}
\item There exists a trivial sub skew left brace $(H, \cdot, \circ )$ of $(G, \cdot, \circ )$.
\item $(H, \cdot)$ is normal in $(G,\cdot)$ and $(H, \circ)$ is normal in $(G, \circ)$.
\item $H \cdot a = H \circ a$ for all $a \in G$.
\item The quotient $G/H$  is a trivial skew left brace.
\end{enumerate}
\end{definition}

We have the following result \cite[Theorem 2.12]{MR4346001}.

\begin{theorem} 
Any $\lambda$-homomorphic skew left brace is meta-trivial.
\end{theorem}

\begin{proof}
Let $(G, \cdot, \circ )$ be a $\lambda$-homomorphic skew left brace, that is, $\lambda : (G, \cdot) \to \Aut(G, \cdot)$ is a group homomorphism. Recall that  $\lambda : (G, \circ) \to \Aut(G, \cdot)$ is always a group homomorphism. Note that
$$\ker (\lambda) = \{ a \in G \, \mid \, a \circ b = a \cdot b \, \mbox{  for all } b \in G \}.$$
Setting $H = \ker (\lambda)$, we see that $(H, \cdot)=(H, \circ)$ is a normal subgroup of both $(G, \cdot)$ as well as $(G, \circ)$. Clearly, $(H, \cdot, \circ)$  is a trivial sub skew left brace of $(G, \cdot, \circ )$. Set $(\overline{G}, \cdot) = (G, \cdot)/(H, \cdot)$, the set of right cosets of $(H, \cdot)$ in $(G, \cdot)$. Then
we see that 
$$H\cdot b = \{a \cdot b \, \mid \, a \in H\} = \{a \circ b \, \mid \, a \in H\}  = H \circ b$$
for all $b \in G$. Hence, we also have the group $(\overline{G}, \circ) = (G, \circ)/(H, \circ)$. Clearly, we have
$$ (H \cdot a) \cdot (H \cdot b)= H \cdot(a \cdot b )$$
and 
$$ (H \circ a) \circ (H \circ b) =  H \circ (a \circ b)$$
for all $a, b \in G$. By Theorem \ref{t1}, we have
$$ b^{-1} \cdot \lambda_a(b) \in \ker(\lambda)=H$$
for all $a, b \in G$. Since $a \circ b= \big(a \cdot  \lambda_a(b) \cdot b^{-1} \cdot a^{-1} \big) \cdot a \cdot b$ and $(H, \cdot)$ is normal in $(G, \cdot)$, it follows that 
$$ (H \circ a) \circ (H \circ b) =  H \circ (a \circ b)=H \cdot (a \circ b) = H \cdot \big(a \cdot  \lambda_a(b) \cdot b^{-1} \cdot a^{-1}\big) \cdot a \cdot b= H \cdot  (a \cdot b)= (H \cdot  a) \cdot (H \cdot b).$$
 This proves that $(\overline{G}, \cdot, \circ)$ is a trivial skew brace, which completes the proof.
$\blacksquare$ 
\end{proof}

\begin{remark}
{\rm 
It is natural to ask whether every meta-trivial skew left brace is $\lambda$-homomorphic. It turns out to be false in general. For example, consider the symmetric group $\Sigma_3$ on three letters with group operations $a \cdot b=ab$ and $a\circ b=ba$. Then $(\Sigma_3, \cdot, \circ)$ is a meta-trivial skew left brace that is not $\lambda$-homomorphic. 
}
\end{remark}

The following result gives a reduction argument, which allows us to verify  condition of Theorem \ref{t1} on generators of a given group and on the images of the generators under the map $\lambda$ \cite[Proposition 2.14]{MR4346001}. 

\begin{proposition}\label{redprop}
Let $G$ be a group generated by  $X=\left\{ x_i \,\mid\, i\in I\right\}$ and $\lambda : G \to \Aut (G)$, given by $a \mapsto \lambda_a$,  a group homomorphism.
If $\lambda_{x_i} = \varphi_i$ and $x_i^{-1} \varphi_j(x_i) \in \ker(\lambda)$ for $i, j \in I$, then
\begin{equation*} \label{7}
b^{-1} \lambda_a (b) \in \ker(\lambda)
\end{equation*}
for all $a, b \in G$.
\end{proposition}

\begin{proof}
Suppose that $x_i^{-1} \varphi_j(x_i) \in \ker(\lambda)$ for all $i, j \in I$. Let $\alpha \in \Aut(G)$.
Since $x\alpha(x^{-1}) = x \big(\alpha(x) \big)^{-1}=x \big(x^{-1}\alpha(x) \big)^{-1} x^{-1}$ for all $x \in G$, it follows that
\begin{equation}\label{eqn1sec2}
x^{-1}\alpha(x) \in \ker (\lambda) \Longleftrightarrow x\alpha(x^{-1})  \in \ker(\lambda).
\end{equation}
 Next, since
$x^{-1}\alpha^{-1}(x)= \left( \alpha^{-1}(x^{-1}) \alpha (\alpha^{-1}(x)) \right)^{-1}$ for all $x \in G$, we have
\begin{equation}\label{eqn2sec2}
z^{-1}\alpha(z)  \in \ker(\lambda) \Longleftrightarrow x^{-1}\alpha^{-1}(x)  \in \ker(\lambda),
\end{equation}
where $z = \alpha^{-1}(x)$.
\para
We claim that $x_i^{-1}\varphi_j^{-1}(x_i) \in \ker(\lambda)$ for all $i, j \in I$. Let us fix a pair $i, j \in I$. In view of \eqref{eqn2sec2}, we only need to show that $z_1^{-1}\varphi_j(z_1) \in \ker(\lambda)$, where $z_1 = \varphi_j^{-1}(x_i)$. Let $z_1 = x_{i_{1}}^{e_1} \cdots  x_{i_{t}}^{e_t}$ for some positive integer $t$, where $i_1,\ldots,i_t \in I$ and $e_1, \ldots, e_t  \in \{1, -1 \}$. If $t = 1$, then we are done by the given hypothesis along with \eqref{eqn1sec2}. We now apply induction on $t$. Note that 
$$(xy)^{-1}\alpha(xy) = \big(y^{-1}(x^{-1}\alpha(x))y \big)  \big(y^{-1}\alpha(y) \big)$$
 for all $x, y \in G$ and $\alpha \in \Aut(G)$. As a consequence,  if both $x^{-1}\alpha(x)$ and $y^{-1}\alpha(y)$ lie in  $\ker(\lambda)$ for  given $x, y \in G$, then $(xy)^{-1}\alpha(xy)  \in \ker(\lambda)$. Using this fact for $\alpha = \varphi_j$, $x = x_{i_{1}}^{e_1}$ and $y = x_{i_{2}}^{e_2} \cdots  x_{i_{t}}^{e_t}$, our claim follows by induction.
\para
Let $a=x_{i_{1}}^{\epsilon_1} \cdots  x_{i_{k}}^{\epsilon_k}$ and $b=x_{j_{1}}^{\varepsilon_1} \cdots  x_{j_{m}}^{\varepsilon_m}$ be non-trivial elements of  $G$, where  $i_1,\ldots,i_k,j_1,\ldots,j_m \in I$ and $\epsilon_1,\ldots,\epsilon_k,\varepsilon_1,\ldots,\varepsilon_m \in \{1, -1 \}$. Then, we have
$$
b^{-1}\lambda_a(b)=x_{j_{m}}^{-\varepsilon_m} \cdots  x_{j_{1}}^{-\varepsilon_1} \,
\varphi_{i_{1}}^{\epsilon_1} \cdots  \varphi_{i_{k}}^{\epsilon_k} \, 
x_{j_{1}}^{\varepsilon_1}\cdots  x_{j_{m}}^{\varepsilon_m}.
$$
\para
If $k=m=1$, then the assertion holds by the given hypothesis, \eqref{eqn1sec2} and the preceding claim. Assume that at least one of $k$ or $m$ is larger than $1$. As explained above,  if both $x^{-1}\alpha(x)$ and $y^{-1}\alpha(y)$ lie in  $\ker(\lambda)$ for  given $x, y \in G$, then $(xy)^{-1}\alpha(xy)  \in \ker(\lambda)$. Consequently, we can assume, without loss of generality, that $m=1$. Let $\beta \in \Aut(G)$ be another automorphism. Note that 
$$x^{-1} (\alpha\beta)(x)= \big(x^{-1}\beta(x) \big) \big(\beta(x)^{-1} (\alpha \beta)(x) \big)$$
for all $x \in G.$  Thus, if both $x^{-1}\beta(x)$ and $\beta(x)^{-1}\alpha(\beta(x))$  lie in $\ker(\lambda)$  for a given $x \in G$, then $x^{-1} (\alpha \beta)(x)  \in \ker(\lambda)$. Take $x = b = x_{j_{1}}^{\varepsilon_1}$, $\alpha = \varphi_{i_{1}}^{\epsilon_1}$ and $\beta = \varphi_{i_{2}}^{\epsilon_2} \cdots  \varphi_{i_{k}}^{\epsilon_k}$. Using the arguments mentioned above, it is easy to conclude that  $\beta(x)^{-1}\alpha(\beta(x)) \in \ker(\lambda)$.  It now remains to prove that $x^{-1}\beta(x) \in \ker(\lambda)$. But this follows by the  iterative process replacing, successively, $\alpha$ by $\varphi_{i_{r}}^{\epsilon_r}$ and $\beta$ by  $\varphi_{i_{r+1}}^{\epsilon_{r+1}} \cdots  \varphi_{i_{k}}^{\epsilon_k}$ for $2 \le r \le k-1$, and using the given hypothesis, \eqref{eqn1sec2}  and the above claim.  This completes the proof.
$\blacksquare$ 
\end{proof}

We construct skew  left braces on  groups admitting subgroups of index two. 

\begin{proposition} \label{abgr}
Let $(G, \cdot)$ be a group,  $B$  an index two subgroup of $G$ and $\phi \in \Aut(G, \cdot)$ an involution such that $\phi(B) = B$.  Then $(G, \cdot , \circ)$ is a skew  left brace, where 
$$
a \circ b =
\left\{
\begin{array}{ll}
a \cdot b \quad \quad  \textrm{if}~ a \in B~\textrm{and}~ b \in G,\\
a \cdot \phi(b) \quad   \textrm{if}~a \not\in B~\textrm{and}~ b \in G. &
\end{array}
\right.
$$
\end{proposition}

\begin{proof}
Consider the map $\lambda : (G, \cdot) \to \Aut(G, \cdot)$ given by
$$
\lambda_a =
\left\{
\begin{array}{ll}
\id \quad \textrm{if}~ a \in B, \\
\phi \quad \textrm{if}~ a \not\in B. &
\end{array}
\right.
$$
Note that $\lambda$ is a homomorphism. A straightforward calculation shows that $b^{-1} \cdot \lambda_a(b) \in \ker (\lambda)$ for all $a, b \in G$. Hence, by Theorem \ref{t1}, it follows that $H_{\lambda}$ is a regular subgroup of $\Hol (G, \cdot)$. Consequently, by Theorem \ref{gv2017}, $(G, \cdot, \circ)$ is a skew left brace with $\circ$ as given in the statement.
$\blacksquare$ 
\end{proof}
\para

Next, we construct some non-trivial $\lambda$-homomorphic left braces on free abelian groups. 

\begin{example}
{\rm The infinite cyclic group $(\mathbb{Z}, +)$ admits a non-trivial left brace structure with its multiplicative group structure given by
$$
m \circ n = m + (-1)^m n
$$
for $m, n \in \mathbb{Z}$  \cite{MR3957824, MR4003478}. Note that $\Aut(\mathbb{Z}, +) \cong \{ 1, -1\}$, where $1$ denotes the identity automorphism and $-1$ denotes multiplication by $-1$. Further,
$$
\lambda : (\mathbb{Z}, +) \to \Aut(\mathbb{Z}, +)
$$
is given by $\lambda_m = (-1)^m$, and  is a group homomorphism. It follows from Theorem \ref{t1} that $(\mathbb{Z}, +, \circ)$ is a $\lambda$-homomorphic left brace.}
\end{example}
\para

\begin{example}\label{fag}
{\rm We now construct  $\lambda$-homomorphic left braces on higher rank free abelian groups. Let $A  \cong \mathbb{Z}^n$ be the free abelian group of rank $n$ with free basis  $\{ x_1, \ldots, x_n \}$. For each $1 \le i \le n-1$, define $\varphi_i \in \Aut(A)$ by
$$
\varphi_i :
\left\{
\begin{array}{ll}
x_i \to x_i + x_n, &  \\
x_j \to x_j & ~\mbox{for} ~j \not= i,
\end{array}
\right.
$$
and let
$$\varphi_n = \id.$$
Note that, these automorphisms  are pairwise commuting, and hence the map $\lambda  : (A, +) \to \Aut(A, +)$ given by $\lambda_{x_i}=\varphi_i$ for each $i$, is a group homomorphism. Further, we have
$$
\lambda \big(\varphi_i(x_j)-x_j \big)=\id
$$
for all $1 \le i, j \le n$. Hence, it follows from Proposition \ref{redprop} and  Theorem \ref{t1} that
$$
H_{\lambda }= \big\{  (\lambda_a, a) \, \mid \, a \in A   \big\}=
 \left\{ \, \left(\prod_{i=1}^{n-1} \varphi_i^{\alpha_i}, ~\sum_{i=1}^n \alpha_i x_i \right)  \, \mid \,  \alpha_i \in \mathbb{Z}  \, \right\}
$$
is a regular subgroup of  $\Hol(A, +)$, where $a = \sum_{i=1}^n \alpha_i x_i$. Thus, by Theorem \ref{gv2017}, we obtain a $\lambda$-homomorphic left brace $(A, + , \circ )$, where $a \circ b = a+ \lambda_a(b)$ for $a, b \in A$.
\para 
Next, we investigate  the structure of the multiplicative group  $(A, \circ )$. Note that the two operations  $+$ and  $\circ$ agree on $ \ker(\lambda) =\left\langle x_n\right\rangle \cong \mathbb{Z}$. Further, for elements $a =  \sum_{i=1}^n \alpha_i x_i$ and $b = \sum_{i=1}^n \beta_i x_i$ of $(A, +)$, we have
\begin{eqnarray*}
a \circ b &=& a + \lambda_a(b) \\
&=& (\alpha_1 + \beta_1) x_1 + (\alpha_2 + \beta_2) x_2 + \cdots +(\alpha_{n-1} + \beta_{n-1}) x_{n-1} +\\
& & (\alpha_n + \beta_n + \alpha_1 \beta_1 + \alpha_2 \beta_2 + \cdots + \alpha_{n-1} \beta_{n-1}) x_n\\
& =&  b \circ a,
\end{eqnarray*}
which implies that $(A, \circ)$ is an abelian group. Moreover, the powers of the generators $x_i$ with respect to $\circ$ are given by
$$
x_i^{\circ k} = k x_i + \frac{k(k-1)}{2}x_n
$$
for each $1 \le i \le n-1$ and
$$
x_n^{\circ k} = k x_n
$$
for each $k \in \mathbb{Z}$. Thus, any element  $b = \sum_{i=1}^n \beta_i x_i$ of $(A,+)$ has the unique form
$$
b = x_1^{\circ \alpha_1} \circ \cdots \circ x_n^{\circ \alpha_n},
$$
where $\alpha_i = \beta_i$ for $ 1 \le  i \le n-1$ and
$$
\alpha_n = \beta_n - \sum_{j=1}^{n-1} \frac{\beta_j (\beta_j - 1 )}{2}.
$$
Hence,  we have $(A, \circ) \cong (A, +)$, and  there is a short exact sequence
$$
1 \to \ker(\lambda) \to (A, \circ) \to \mathbb{Z}^{n-1} \to 1.
$$}
\end{example}
\para

We now construct some $\lambda$-homomorphic skew left braces on non-abelian free groups.  

\begin{example}\label{homomorphic skew brace on free group}
{\rm 
Let $F_n = \langle x_1,\ldots,x_n \rangle$  be the free group of rank $n$. Let  $\varphi_1, \ldots, \varphi_n \in \Aut(F_n)$ be a set of pairwise commuting automorphisms such that
$$
x_i^{-1}\varphi_j(x_i)\in F_n'
$$
for all $1 \le i,j \le n.$ Let  $\lambda : F_n \to \Aut(F_n)$ be the homomorphism defined on the generators by
$$
\lambda_{x_i} = \varphi_i
$$
for  $1 \le i  \le n$. Clearly, $\ker(\lambda)$ contains the commutator subgroup  $[F_n, F_n]$ of $F_n$, and hence
$$
x_i^{-1}\varphi_j(x_i)\in \ker(\lambda)
$$
for all $1 \le i,j \le n$. It follows from Proposition \ref{redprop} that $b^{-1} \lambda_a (b)\in \ker(\lambda)$ for all  $a, b \in F_n$. Hence, by Theorem \ref{t1}, $H_{\lambda}= \{ (\lambda_a, a) \, \mid \, a \in F_n \}$ is a regular subgroup of  $\Hol(F_n)$. Thus, $( F_n, \cdot, \circ  )$ is a $\lambda$-homomorphic skew left brace with $a \circ b=a \cdot \lambda_a(b)$ for $a, b\in F_n$, where $\cdot$ is the product in $F_n$.}
\end{example}
\para

An automorphism $\alpha$ of a group $G$ is called an \index{IA-automorphism}{\it IA-automorphism} if $a^{-1}\alpha(a) \in [G, G]$ for all $a \in G$. The set $\IA(G)$ of all IA-automorphisms forms a subgroup of $\Aut(G)$.  Thus, the preceding construction  is valid for  any choice of pairwise commuting automorphisms $\varphi_1, \ldots, \varphi_n\in \IA(F_n)$. We now carry-out  this construction concretely for $n=4$.

\para

\begin{example} \label{fIA}
{\rm 
Let $F_4 = \langle x_1, x_2, x_3, x_4 \rangle$ be the free group of rank 4. For $u, v \in [F_4, F_4]$, let $\varphi_1, \varphi_2, \varphi_3, \varphi_4$ be automorphisms of  $F_4$ defined on free generators as 
$$
\varphi_1=\varphi_2 : \left\{
\begin{array}{ll}
x_1 \to x_1 u, &  \\
x_i \to x_i & \mbox{for} ~i \not= 1,
\end{array}
\right.
\quad \textrm{and}\quad
\varphi_3=\varphi_4 : \left\{
\begin{array}{ll}
x_i \to x_i & \mbox{for} ~i \not= 4.\\
x_4 \to x_4 v. &  \\
\end{array}
\right.
$$
It is not difficult to check that $\varphi_1 \varphi_3 = \varphi_3 \varphi_1$. Hence, by Example \ref{homomorphic skew brace on free group}, we get a skew left brace $(F_4, \cdot \; , \circ)$, for which the operation $\circ$ on free generators is given by
\begin{eqnarray*}
& x_1 \circ x_1=x_1^2u, \quad x_1 \circ x_2=x_1x_2,\quad x_1 \circ x_3=x_1x_3,\quad x_1 \circ x_4=x_1x_4,\\
& x_2 \circ x_1=x_2x_1u,\quad x_2 \circ x_2=x_2^2,\quad x_2 \circ x_3=x_2x_3,\quad x_2 \circ x_4=x_2x_4,\\
&x_3 \circ x_1=x_3x_1,\quad x_3 \circ x_2=x_3x_2,\quad x_3 \circ x_3=x_3^2,\quad x_3 \circ x_4=x_3x_4v,\\
&x_4 \circ x_1=x_4x_1,\quad x_4 \circ x_2=x_4x_2,\quad x_4 \circ x_3=x_4x_3,\quad x_4 \circ x_4=x_4^2v.
\end{eqnarray*}
The operations  $\cdot$ and  $\circ$ coincide  on  $[F_4, F_4]$. Note that the inverse $\bar x$ of an element $x$ with respect to $\circ$ is given by $\bar x =\lambda_x^{-1}(x^{-1})$. Hence, we have
$$
\bar x_1 = ux_1^{-1},\quad
\bar x_2 = x_2^{-1},\quad
\bar x_3 = x_3^{-1},\quad
\bar x_4 = v x_4^{-1}.
$$
We can write
$$
F_4= \sqcup~ x_1^{a_1}x_2^{a_2}x_3^{a_3}x_4^{a_4} \, [F_4, F_4],
$$
where  $a_1, a_2, a_3, a_4 \in \mathbb{Z}$ and
$$
x_1^{a_1}\circ x_2^{a_2}\circ x_3^{a_3}\circ x_4^{a_4}\equiv
x_1^{a_1}x_2^{a_2}x_3^{a_3}x_4^{a_4}\, \mod \, [F_4, F_4].
$$
Hence, it follows that
$$
(F_4, \circ)= \sqcup~ x_1^{a_1}\circ x_2^{a_2}\circ x_3^{a_3}\circ x_4^{a_4} \,[(F_4, \circ), (F_4, \circ)]
$$
is a decomposition of  $(F_4, \circ)$ into disjoint union of left cosets of $[(F_4, \circ), (F_4, \circ)]$, which gives
$$
(F_4, \circ)/[(F_4, \circ), (F_4, \circ)] \cong \mathbb{Z}^4.
$$}
\end{example}
\bigskip 
\bigskip


\section{Lambda-cyclic skew braces} \label{lambda-cyc}

We begin this section with the following definition.

\begin{definition}
A $\lambda$-homomorphic skew  left brace $(G, \cdot, \circ)$ is said to be \index{$\lambda$-cyclic skew brace}{$\lambda$-cyclic} if the image of $\lambda$ is a cyclic subgroup of $\Aut (G, \cdot)$.
\end{definition}

Let $G$ be a group with a presentation $G = \langle S\mid R \rangle,$ where $S= \{ g_1, g_2,g_3, \ldots \}$ is the set of generators and $R = \{ r_1, r_2, r_3, \ldots \}$ is the set of defining relations. Then any element  $w\in G$ can be written as
\begin{equation} \label{w}
w = g_{i_1}^{\alpha_1} g_{i_2}^{\alpha_2} \ldots g_{i_s}^{\alpha_s}
\end{equation}
for some $g_{i_j} \in S$ and $\alpha_j \in \mathbb{Z}$. We define the \index{logarithm of a word}{\it logarithm} $\mathfrak{L}(w)$ of $w$ as 
$$
\mathfrak{L}(w) = \sum_{j=1}^s \alpha_j.
$$

A presentation $G = \langle S\mid R \rangle$ is called a \index{homogeneous presentation}{\it homogeneous presentation} if  $\mathfrak{L}(r_i) = 0$ for all $r_i \in R$. A group is said to be a \index{homogeneous group}{\it homogeneous group} if it admits a homogeneous presentation. It is not difficult to see that for a homogeneous group $G$, the logarithm $$\mathfrak{L} : G \to \mathbb{Z}$$ is well-defined. In fact, the map $\mathfrak{L} : G \to \mathbb{Z}$ is well-defined for each group $G$ which is free in some variety. Note that a homogeneous group can have torsion. For example, the relation $(g_i g_j^{-1})^2 = 1$ has zero logarithm,  but the element $g_i g_j^{-1}$ has finite order.

\begin{example}{\rm 
Consider the following examples:
\begin{enumerate}
\item Any free group, free solvable group and free nilpotent group is homogeneous.
\item If $L$ is a link in  $\mathbb{S}^3$, then the Wirtinger presentation of the link group $G(L) = \pi_1(\mathbb{S}^3 \setminus L)$ is homogeneous.
\end{enumerate}}
\end{example}

\begin{lemma} \label{l1}
Let $G$ be a homogeneous group admitting a homogeneous presentation $G=\langle S\mid R \rangle$. Let $\varphi \in \Aut (G)$ be such that
$\mathfrak{L}(\varphi (g_i)) = 1$ for all $g_i \in S.$ Then the following  assertions hold for all $u, v, w \in G$:
\begin{enumerate}
\item  $\mathfrak{L}(u v) = \mathfrak{L}(u) + \mathfrak{L}(v)$.
\item  $\mathfrak{L}(u^{-1}) = -\mathfrak{L}(u)$.
\item  $\mathfrak{L} (\varphi^k (w)) = \mathfrak{L}(w)$ for any integer $k$.
\end{enumerate}
\end{lemma}

\begin{proof}
The first and second assertions follow directly from the definition of the logarithm. Clearly, it is sufficient to prove the  third assertion for  $k=1$. Suppose that $w = g_{i_1}^{\alpha_1} g_{i_2}^{\alpha_2} \ldots g_{i_s}^{\alpha_s}$ for some $g_{i_j} \in S$ and $\alpha_j \in \mathbb{Z}$. Since 
$$\varphi(w) = \varphi(g_{i_1})^{\alpha_1} \varphi(g_{i_2})^{\alpha_2} \cdots \varphi(g_{i_s})^{\alpha_s},$$
using the first two assertions and the property of  $\varphi$, we get
\begin{eqnarray*}
\mathfrak{L} \big(\varphi(w) \big) &=& \alpha_1 \mathfrak{L} \big(\varphi(g_{i_1})\big) + \alpha_2 \mathfrak{L} \big(\varphi(g_{i_2})\big) + \cdots + \alpha_s \mathfrak{L} \big(\varphi(g_{i_s})\big) \\
&=& \alpha_1  + \alpha_2  + \cdots + \alpha_s\\
& =& \mathfrak{L}(w),
\end{eqnarray*}
which is desired. $\blacksquare$ 
\end{proof}

\begin{proposition} \label{lcb}
Let $G$ be a homogeneous group admitting a homogeneous presentation $G= \langle S\mid R\rangle$. Let $\varphi \in \Aut (G)$ be such that $\mathfrak{L}(\varphi (g_i)) = 1$ for all $g_i \in S.$ Then the set
$
H_{\lambda} = \big\{ (\varphi^{\mathfrak{L}(w)}, w)  \, \mid \, w \in G \big\}
$
is a regular subgroup of $\Hol (G)$.
\end{proposition}

\begin{proof}
By Lemma \ref{l1}(1), the map
$$
\lambda : G \to \Aut (G),
$$
given  by  $\lambda(w) = \varphi^{\mathfrak{L}(w)}$, is a group homomorphism. By Lemma \ref{l1}(3), we  get   $\lambda \big(w^{-1}\cdot  \lambda_u(w)\big) = \id$ for all $u, w \in G$.  It now follows from  Theorem \ref{t1} that   $H_{\lambda} = \big\{(\varphi^{\mathfrak{L}(w)}, w) \, \mid \, w \in G \big\}$ is a regular subgroup of
$\Hol (G)$. $\blacksquare$ 
\end{proof}

Let us see an example satisfying the hypothesis of Proposition \ref{lcb}.

\begin{example}{\rm 
The free group $F_n=\langle x_1, \ldots, x_n \rangle$ of rank $n$ is homogeneous. Let $\theta$ be an $\IA$-automorphism of  $F_n$. It follows that $\mathfrak{L}(\theta(x_i)) = 1$ for all $i$. Let $\lambda : F_n \to \Aut(F_n)$ be the homomorphism defined on the free generators as
 $$
\lambda(x_i)=\theta
$$
for all $1 \le i  \le n$. Hence, by Proposition \ref{lcb}, the set $H_{\lambda}$ is a regular subgroup of  $\Hol(F_n)$.}
\end{example}

The following example shows that Proposition \ref{lcb} does not work for all automorphisms of $F_n$.

\begin{example}{\rm 
Let   $\theta$ be the automorphism of $F_n=\langle x_1, \ldots, x_n \rangle$ defined on its free generators as
$$
\theta(x_i)= 
\left\{
\begin{array}{ll}
x_1x_2 & ~\mbox{if} ~i=1, \\
x_i & ~\mbox{if} ~i \neq 1.
\end{array}
\right.
$$
Then, we have
$$
\lambda \big(x_1^{-1}\lambda_{x_2}(x_1)\big) = \lambda \big(x_1^{-1} \theta(x_1)\big) = \lambda(x_2) = \theta \ne \id.
$$
Hence, by Theorem \ref{t1}, $H_{\lambda}$ is not a subgroup of $\Hol( F_n)$. Note that, in this case,  $\mathfrak{L}(\theta(x_1)) = 2$.}
\end{example}

The subsequent result is straightforward.

\begin{proposition}\label{impprop}
Let   $\theta$ be the automorphism of $F_n=\langle x_1, \ldots, x_n \rangle$ defined on its free generators as
\begin{eqnarray*}
\theta(x_i) &=&  x_{\sigma(i)} u_i\,  \mbox{for some } u_i \in F_n'\; \mbox{and some permutation } \sigma \in \Sigma_n\\
\textrm{or} &  &\\
\theta(x_i) &=&  x_{\sigma(i)}^{-1} u_i\,  \mbox{for some } u_i \in F_n'\; \mbox{and some permutation } \sigma \in \Sigma_n \; \mbox{and } \theta^2 =1.
\end{eqnarray*}
Let $\lambda : F_n\to \Aut(F_n)$ be the homomorphism given by $\lambda(x_i)=\theta$ for all $i$. Then $H_\lambda$ is a regular subgroup of $\Hol(F_n)$.
\end{proposition}

\begin{proof}
Let $\theta(x_i) =  x_{\sigma(i)} u_i$ for some $u_i \in F_n'$ and some permutation $\sigma \in \Sigma_n$. Let $\lambda : F_n \to \Aut (F_n)$
be the homomorphism given by $\lambda(w) = \theta^{\mathfrak{L}(w)}$. If $\lambda_a(b) = \theta^{\mathfrak{L}(a)}(b) = c$, then we see that $\mathfrak{L}(c) = \mathfrak{L}(b)$. This gives
$$
\lambda \big(b^{-1}\lambda_a(b) \big)= \lambda (b^{-1} c) = \theta^{\mathfrak{L}(b^{-1} c)}  =\theta^{-\mathfrak{L}(b)} \, \theta^{\mathfrak{L}(c)} = \theta^{-\mathfrak{L}(b)}\,  \theta^{\mathfrak{L}(b)}
= 1.
$$
Hence, by Theorem \ref{t1}, $H_{\lambda}$ is a regular subgroup of $\Hol( F_n)$. The second case is similar.  $\blacksquare$ 
\end{proof}
\para

The rest of the section is focussed on constructing Lambda-cyclic skew left braces on the free group of rank two with respect to specific automorphisms. We begin with the following result \cite[Proposition 3.10]{MR4346001}.

\begin{proposition} \label{fp}
Let $F_2=\langle x, y \rangle$ be the free group of rank two and $\theta\in \Aut(F_2)$  such that  $\theta(x)=y$ and $\theta(y)=x$. Let $H_0$ be the kernel of the homomorphism  $F_2 \to  \Aut(F_2)$ given by $x,y \mapsto \theta$. If $\cdot$ denote the product in $F_2$ and
$$
a\circ b=
\left\{
\begin{array}{ll}
  ab & \textrm{if}~ a\in H_0~\textrm{and}~b \in F_2,\\
  a\theta(b) & \textrm{if}~ a\not\in H_0~\textrm{and}~b \in F_2, \\
\end{array}
\right.
$$
then $(F_2, \cdot, \circ )$ is a skew left brace with $(F_2, \circ) \cong (F_2, \cdot)$.
\end{proposition}

\begin{proof}
Let $\lambda : (F_2, \cdot)\to \Aut(F_2, \cdot)$ be the homomorphism given by $\lambda(x)=\lambda(y)=\theta$. Then $\theta$ satisfies the hypothesis of Proposition \ref{impprop}, and therefore $H_{\lambda}$ is a regular subgroup of  $\Hol(F_2, \cdot)$.  By Reidemeister-Schreier method \cite[Theorem 2.6]{MR0207802}, it follows that $H_0 = \langle xy, yx, x^2, y^2 \rangle$. Note that  $(H_0, \cdot)$ is a free group with free generators $\{xy, x^2,xy^{-1} \}$, and has index two in $(F_2, \cdot)$. Since $(F_2, \cdot)$ is a normal subgroup of $\Hol(F_2, \cdot)$, it follows that $H_{\lambda}$ and $F_2$ are not conjugate in $\Hol(F_2, \cdot)$. Further, since $\lambda_a = \id$ for $a \in H_0$, we see that $( F_2, \cdot, \circ )$ is a skew  left brace with
$$
a\circ b=
\left\{
\begin{array}{ll}
  ab & \textrm{if}~ a\in H_0 ~\textrm{and}~b \in F_2,\\
  a\theta(b) & \textrm{if}~ a\not\in H_0~\textrm{and}~b \in F_2. \\
\end{array}
\right.
$$
Note that the binary operations $\circ$ and $\cdot$ coincide on $H_0$, and hence $(H_0, \circ)$ is a free group of rank three. We claim that $(F_2, \circ)$ is torsion-free. Note that  $H_0$  is a normal and free subgroup of  $(F_2,  \circ)$ of index two. Thus, if $(F_2,  \circ)$ has an element of finite order, then it must be an element of order two. Let $1 \ne a \in F_2$ such that $a\circ a= 1$.
Then $a=xa_0$ for some $a_0 \in H_0$  and $xa_0y \theta(a_0)= a \circ a = 1$. Since $\mathfrak{L} \big(xa_0y \theta(a_0)\big)=\mathfrak{L}(1)$, we have
$$
2+ \mathfrak{L}(a_0) + \mathfrak{L} \big(\theta(a_0)\big)=0.
$$
Noting that $\mathfrak{L}(a_0)+ \mathfrak{L} \big(\theta(a_0)\big)\equiv 0 \mod 4$, we get $2\equiv 0 \mod 4$, which is absurd. Hence, $(F_2,  \circ)$ does not admit any torsion element. It is a well-known result that a finitely generated torsion free group admitting a free subgroup of finite index is itself free (see \cite{MR0228573} and \cite[Theorem B]{MR0240177}). Hence, $(F_2, \circ)$ is a free group. Using the facts that $\mathrm{rank}(H_0, \circ)=3$, $|(F_2,  \circ): (H_0, \circ)|=2$ and the well-known formula
$$
\mathrm{rank}(H_0, \circ)=1+\big(\mathrm{rank}(F_2,  \circ)-1\big)~ |(F_2,  \circ): (H_0, \circ)|,
$$
we see that  $\mathrm{rank}(F_2,  \circ)=2$. This proves that $(F_2, \circ) \cong (F_2, \cdot)$.
$\blacksquare$ 
\end{proof}

\begin{proposition} \label{fi}
Let  $F_2 = \langle x, y \rangle$ be the free group of rank two and  $\theta \in \Aut(F_2)$  such that $\theta(x)= x^{-1}$ and $\theta(y)= y^{-1}$. Let
$H_0$ be  the kernel of the homomorphism  $F_2 \to \Aut(F_2)$ given by $x, y \mapsto \theta$. If $\cdot$ denote the product in $F_2$ and
$$
a\circ b=
\left\{
\begin{array}{ll}
  ab &    \textrm{if}~ a\in H_0~\textrm{and}~b \in F_2,\\
  a\theta(b) & \textrm{if}~ a\not\in H_0~\textrm{and}~b \in F_2, \\
\end{array}
\right.
$$
then $(F_2, \cdot, \circ )$ is a skew left brace. Moreover, 
$$
(F_2, \circ) \cong F_3 \rtimes \mathbb{Z}_2 \cong \left\langle\, p,q,r,s  \, \mid \,  s^2=1,\,\, sps=p^{-1},\,\, sqs=q^{-1},\,\, srs=pr^{-1}p^{-1}\, \right\rangle.
$$
\end{proposition}

\begin{proof}
Let $\lambda :  (F_2, \cdot) \to \Aut(F_2, \cdot)$ be the homomorphism given by $\lambda(x) = \lambda(y)=\theta$. Then, by Proposition \ref{impprop}, $H_{\lambda}$ is a regular subgroup of  $\Hol(F_2, \cdot)$. Note that the subgroup   $(H_0, \cdot)$ consists of words of even length. Since $(F_2, \cdot)$ is a normal subgroup of $\Hol(F_2, \cdot)$, it follows that  $H_{\lambda}$ and $F_2$ are not conjugate in $\Hol(F_2, \cdot)$. Since $\lambda_a = \id$ for all $a \in H_0$, it follows that $( F_2, \cdot, \circ )$ is a skew left brace with
$$
a\circ b=
\left\{
\begin{array}{ll}
  ab &  \textrm{if}~ a\in H_0~\textrm{and}~b \in F_2,\\
  a\theta(b) &   \textrm{if} ~ a\not\in H_0~\textrm{and}~b \in F_2. \\
\end{array}
\right.
$$
We now determine  the structure of  $(F_2, \circ)$. Note  that $x \circ x = 1= y \circ y$ and $(F_2, \circ)$ is generated by $\{H_0, x\}$. Since $x \circ x=1$ and $H_0$ is normal in  $(F_2, \circ)$, it follows that
$$
(F_2, \circ) \cong H_0 \rtimes \mathbb{Z}_2,
$$
where $\mathbb{Z}_2 \cong \left\langle s \right\rangle$ for $s \not\in H_0$. Note that  $H_0$ is freely generated by $\{p, q, r \}$, where
$$
p= xy,\quad q=x^2\quad \textrm{and}\quad r=y^2.
$$
To avoid confusion, let us denote $x$ in $(F_2, \circ)$ by $s$. Then, we have the conjugation relations
$$
sps=p^{-1},\quad sqs=q^{-1}\quad \textrm{and}\quad  srs=pr^{-1}p^{-1}.
$$
Hence,  $(F_2, \circ)$ has the  presentation
$$
(F_2, \circ)  \cong \left\langle\, p,q,r,s  \, \mid \,  s^2=1,\,\, sps=p^{-1},\,\, sqs=q^{-1},\,\, srs=pr^{-1}p^{-1}\, \right\rangle,
$$
and the proof is complete. $\blacksquare$ 
\end{proof}

\begin{remark}
{\rm 
Propositions \ref{fp} and \ref{fi} provide examples satisfying the hypothesis of Proposition \ref{abgr}. Further, the $\lambda$-cyclic skew  left braces constructed in  Propositions \ref{fp} and \ref{fi} are not isomorphic. In fact, in the first case $(F_2, \circ) \cong F_2$, while in the second case $(F_2, \circ) \cong F_3 \rtimes \mathbb{Z}_2$. But the socles and the quotients by the socles for these skew  left braces are isomorphic.
}
\end{remark}

For abelian groups, we have the following result.

\begin{proposition} \label{abgr1}
Let $(A, + )$ be an abelian group admitting a subgroup $B$ of index two such that $A=B \cup (a_0 + B)$
for some  $a_0 \in A$.  Then $(A, + , \circ)$ is a left brace, where
$$
a \circ b =
\left\{
\begin{array}{ll}
a + b & \textrm{if}~ a \in B~\textrm{and}~b \in A,\\
a - b &   \textrm{if}~ a \not\in B~\textrm{and}~b \in A. 
\end{array}
\right.
$$
Moreover,  $(A, \circ) \cong (B, +) \rtimes \mathbb{Z}_2$.
\end{proposition}

\begin{proof}
Let $( A, +)$ be an abelian group and  $\theta$ the automorphism given by  $\theta (a)=-a$ for all $a \in A$. Since $\theta^2 = \id$  and $\theta(B) = B$, it follows from Proposition \ref{abgr} that $(A, + , \circ)$ is a left brace, where $\circ$ is as given in the statement of the proposition. Note that $a_0 \circ a_0 = 0$, where $0$ is the identity element of $(A, +)$. Thus, we have $\bar a_0 = a_0$, and therefore $\bar a_0 \circ b \circ a_0 = -b$ for all $b \in B$. Since the operations $\circ$ and $+$ coincide on $B$, it follows that $(A, \circ) \cong (B, +) \rtimes \mathbb{Z}_2$.
$\blacksquare$ 
\end{proof}

Let $A = \mathbb{Z}_2 \times \mathbb{Z}_4$, where $\mathbb{Z}_2= \langle a \rangle$ and  $\mathbb{Z}_4=\langle b \rangle$. Let $B_1 = \langle b \rangle = \mathbb{Z}_4$ and $B_2 = \langle a \rangle \times \langle b^2 \rangle = \mathbb{Z}_2 \times \mathbb{Z}_2$. Obviously, $B_1$ and $B_2$ are not isomorphic, and hence  the left braces constructed from $B_1$ and $B_2$ using Proposition \ref{abgr1} are not isomorphic. We conclude this section with the following questions.

\begin{question}
Under what conditions are the left braces, constructed from two subgroups of index two in an abelian group, isomorphic?
\end{question}

\begin{question} 
Is it possible to construct a skew left brace merely with the information that its additive group admits  a subgroup of index two?
\end{question}
\bigskip
\bigskip

\section{Lambda-cyclic braces on free abelian groups} \label{cyc-ab}

Consider the free abelian group
$$
\mathbb{Z}^n=\big\{ \, m_1x_1+\cdots + m_nx_n  \, \mid \,  m_1, \ldots , m_n \in   \mathbb{Z}    \big\}
$$
 of rank  $n$ with free basis  $\{x_1,\ldots, x_n\}$. For $a = m_1x_1+ \cdots + m_n x_n \in \mathbb{Z}^n$, we define 
 $$\mathfrak{L}(a) = \sum_{i=1}^{n}m_i.$$ 
Let $\varphi$ be an automorphism of $\mathbb{Z}^n$ and $\lambda: \mathbb{Z}^n \to \Aut(\mathbb{Z}^n)$ the map given by $a \mapsto \lambda_a= \varphi^{\mathfrak{L}(a)}$. Then $\lambda_{x_i} = \varphi$ for each $1 \le i \le n$. Let $H_{\lambda} = \{\,(\lambda_a, a) \, \mid \, a \in \mathbb{Z}^n \,\}$ and 
$$\mathbb{Z}_0^n = \big\{ a \in \mathbb{Z}^n \, \mid \, \mathfrak{L}(a) = 0  \big\}.$$
Note that $\mathbb{Z}_0^n$ is a subgroup of $\mathbb{Z}_0$.
\para

\begin{proposition} \label{cycprop1}
Let $\mathbb{Z}^n$ be the free abelian group of rank $n$ with free basis  $\{x_1,\ldots, x_n\}$ and $\varphi \in \Aut(\mathbb{Z}^n)$. If $\varphi(x_i) \equiv x_i \mod \mathbb{Z}^n_0$ for each $1 \le i \le n$, then  $H_{\lambda}$ is a regular subgroup of  $\Hol(\mathbb{Z}^n)$. Conversely, if $\varphi$ is of infinite order such that   $H_{\lambda}$ is a regular subgroup of  $\Hol(\mathbb{Z}^n)$, then $\varphi(x_i) \equiv x_i \mod \mathbb{Z}^n_0$ for each $1 \le i \le n$.
\end{proposition}

\begin{proof}
If $\varphi(x_i) \equiv x_i \mod \mathbb{Z}^n_0$ for each $1 \le i \le n$, then it follows from Proposition \ref{lcb} that $H_{\lambda}$ is a regular subgroup of  $\Hol(\mathbb{Z}^n)$.
\para
Suppose that $H_{\lambda}$ is a subgroup of $\Hol(\mathbb{Z}^n)$. Then
$$
(\varphi, x_i) (\varphi, x_j)=\big(\varphi^2,~ x_i+\varphi(x_j) \big)\in H_{\lambda}
$$
for each $1 \le i, j \le n$, which implies that $\mathfrak{L}(x_i+\varphi(x_i)) = 2$ for each $1 \le i \le n$. Hence, if
\begin{equation}\label{varphix expression}
\varphi(x_i)=m_{i1}x_1 + \cdots + m_{in}x_n   
\end{equation}
for some $m_{i1},\ldots , m_{in}\in \mathbb{Z}$, then  $m_{i1} + \cdots + m_{in}=1$ for each $1 \le i \le n$.  Subtracting $x_i$ on both the sides of \eqref{varphix expression} gives
$$
\varphi(x_i)\equiv x_i \mod \mathbb{Z}^n_0
$$
for each $1 \le i \le n$, which proves the converse. $\blacksquare$ 
\end{proof}

\begin{example}{\rm 
The converse statement in Proposition \ref{cycprop1} does not hold if $\varphi$ is of finite order. Let $\mathbb{Z}^2$ be the free abelian group with free basis $\{x_1, x_2 \}$. Consider the automorphism $\psi$ of $\mathbb{Z}^2$ given by $\psi(x_1) = x_1$ and $\psi(x_2) = -x_2.$ Then $H_{\lambda}$ is a regular subgroup of  $\Hol(\mathbb{Z}^2)$, but $\psi(x_2) \not\equiv x_2 \mod \mathbb{Z}^2_0$.}
\end{example}
\para

Our next result is the following theorem \cite[Theorem 4.10]{MR4346001}.

\begin{theorem} 
Let $A \cong (\mathbb{Z}^n, +)$ be the free abelian group of rank $n\geq 2$ with free basis $\{x_1,  \ldots,x_n\}$ and  $\psi \in \Aut(A)$ the  automorphism defined by
$$
\psi :
\left\{
\begin{array}{ll}
x_i \to x_{i+1} &    \textrm{if}~ i \ne n,\\
x_n \to x_1. &
\end{array}
\right.
$$
Let $\lambda : A \to \Aut(A)$ be the group homomorphism given by  $a \mapsto \lambda_a= \psi^{\mathfrak{L}(a)}$ for all $a \in A$.  Then $(A, +, \circ)$ is a $\lambda$-cyclic left brace with
$$
a\circ b=a+\psi^{\mathfrak{L}(a)} (b)
$$
for $a , b \in A.$ Moreover,
$$
(A, \circ) = \big\langle \, x,\,z  \, \mid \,  [x^{\circ n},z] = 1,\,(z \circ x)^{\circ n}=x^{\circ n}, \,
[z,x^{\circ k}, z]=1 ~~ \textrm{for all}~~  1 \le k \le n-2 \big\rangle,
$$
where  $x=x_1$ and $z=x_2-x_1$.
\end{theorem}

\begin{proof}
Note that $A$ is a homogeneous group such that $\mathfrak{L} \big(\psi (x_i) \big) = 1$ for each $1 \le i \le n$. Then, by Proposition \ref{lcb}, $H_{\lambda}$ is a regular subgroup of $\Hol(A)$.  Thus, we obtain a $\lambda$-cyclic left brace $(A, +, \circ)$, where 
$$
a \circ b = a+\psi^{\mathfrak{L}(a)}(b)
$$
for all $a , b  \in A$.
\para 
We now determine the structure of the multiplicative group $(A, \circ)$. Note that
$$
x_1 \circ x_2=x_1+x_3 \;\;\mbox{ and } \;\; x_2 \circ x_1=2x_2,
$$
and hence $(A, \circ)$ is a non-abelian group. Consider the subgroup
$$
A_0 =
\left\{ \, \sum\limits_{i=1}^{n} \alpha_i x_i \in A  \, \mid \, 
\sum\limits_{i=1}^{n} \alpha_i \equiv 0 \mod n  \, \right\}
$$
of $A$. It is not difficult to see that $A_0 = \ker( \lambda)$, and hence the two operations $+$ and $\circ$ coincide on $A_0$. We have the decomposition 
$$
A =A_0 \cup (x_1 + A_0) \cup (2x_1 + A_0) \cup \cdots \cup \big((n-1) x_1 + A_0 \big)
$$
of  $A$ as a union of cosets of $A_0$. Since $\psi(A_0)=A_0$, it follows that
$$
(kx_1)\circ A_0 =(kx_1) + A_0
$$
for each $k\in \mathbb{Z}$. Hence, we obtain the decomposition 
$$
(A, \circ)= A_0 \cup (x_1 \circ A_0) \cup (2x_1 \circ A_0) \cup \cdots \cup \big((n-1)x_1 \circ A_0 \big)
$$
of   $(A, \circ)$  as a union of cosets  of $A_0$. Observe that
$$
kx_1\equiv x_1^{\circ k} \mod A_0
$$
for each $1 \le k \le n$. Hence, we can write
$$
(A, \circ) = A_0 \cup (x_1\circ A_0) \cup (x_1^{\circ 2} \circ A_0) \cup \cdots \cup (x_1^{\circ (n-1)} \circ A_0),
$$
and  therefore $(A, \circ)$ is generated by the element  $x_1$ and the subgroup  $(A_0, \circ)$. We take
$$
\big\{z_1:=x_1+\cdots+x_n, ~~ z_2:=x_2-x_1, \ldots,  z_n:=x_n-x_{n-1} \big\}
$$
as a set of free generators of  $(A_0, \circ)$. Let $[a, b]$ denote the element $\bar a \circ \bar b \circ a \circ b$ for all $a , b \in (A, \circ)$. A direct check shows that we have the relations
\begin{eqnarray*}
x_1^{\circ n}  & =& z_1,\\
\left[ z_i, z_j \right] & =& 1 \quad \textrm{for}~~ 1 \le i, j \le n,\\
x_1\circ z_k & =& x_1\circ (x_k-x_{k-1})=x_1+ x_{k+1}-x_k =z_{k+1}+x_1=z_{k+1}\circ x_1 \quad \textrm{for}~~2 \le k \le n-1,\\
x_1\circ z_n & =& x_1\circ (x_n-x_{n-1})=x_1+x_1-x_n=(x_1-x_n) \circ x_1= (\overline{z_2+z_3+\cdots+z_n}) \circ x_1= \bar z_2\circ \bar z_3\circ \cdots \circ \bar z_n \circ x_1.
\end{eqnarray*}
This shows that  $(A, \circ)$ has a presentation
\begin{eqnarray*}
(A, \circ) &= & \Big\langle x_1,z_1,\ldots,z_n  \, \mid \,  z_1=x_1^{\circ n},~~x_1 \circ z_n  \circ \bar x_1 = \bar z_2 \circ \bar z_3 \circ \cdots \circ \bar z_n,~~ [z_i,z_j]=1 ~~\textrm{for}~~1 \le i,j \le n,\\
& &  x_1 \circ z_k \circ \bar x_1=z_{k+1}
~~\textrm{for}~~2 \le k \le n-1 \Big\rangle.
\end{eqnarray*}
To simplify the presentation of $(A, \circ)$, we note that
\begin{eqnarray*}
(A, \circ) &=& \Big\langle \, x_1,z_2,\ldots,z_n  \, \mid \,  [x_1^{\circ n},z_2]= \cdots  = [x_1^{\circ n},z_n]=1,\,\, [z_i,z_j]=1 ~~\textrm{for}~~2 \le i,j \le n, \\
& & z_3=x_1 \circ z_2 \circ \bar x_1, \ldots, z_n=x_1 \circ z_{n-1} \circ \bar x_1,\,\,
x_1 \circ z_n \circ \bar x_1 = \bar z_2 \circ \bar z_3 \circ \cdots \bar \circ z_n \Big\rangle\\
&=& \left\langle \, x_1,\,\,z_2 ~ \, \mid \, ~
[x_1^{\circ n},z_2]=[x_1^{\circ n},x_1\circ z_2 \circ \bar x_1]= \cdots = [x_1^{\circ n},x_1^{\circ(n-2)} \circ z_2 \circ \bar x_1^{\circ (n-2)}]=1, \right.\\
& &
[z_2,x_1 \circ z_2 \circ \bar x_1]= \cdots = [z_2,x_1^{\circ (n-2)} \circ z_2 \circ \bar x_1^{\circ (n-2)}]=1, \;\ldots,\\
& &[x_1^{\circ(n-3)} \circ z_2 \circ \bar x_1^{\circ(n-3)},x_1^{\circ (n-2)}\circ z_2\circ \bar x_1^{\circ(n-2)}]=1,\\
& & \left. x_1^{\circ (n-1)}\circ z_n \circ \bar x_1^{\circ (n-1)}=\overline{z_2 \circ x_1\circ z_2 \circ \bar x_1 \circ \cdots \circ x_1^{\circ (n-2)}\circ z_2 \circ \bar x_1^{\circ (n-2)}}
\,\right\rangle.
\end{eqnarray*}
Setting  $x = x_1$ and $z = z_2$, we get
\begin{eqnarray*}
(A, \circ) &=& \left\langle \, x,\, z  \, \mid \,  [x^{\circ n},z]=1,\, [z, x\circ z \circ \bar x]= \cdots =[z, x^{\circ (n-2)}\circ z \circ \bar x^{\circ (n-2)}]=1, \right.\\
& & \left. \bar x \circ z \circ  x= \overline{(z \circ x)^{\circ (n-1)} \circ \bar x^{\circ (n-1)}}
\,\right\rangle\\
&=& \Big\langle \, x,\, z  \, \mid \,  [x^{\circ n},z]=1,\, [z, x^{\circ k}\circ z \circ \bar x^{\circ k}]=1 ~~\textrm{for}~~1 \le k \le n-2, \\
& &\bar x \circ z \circ x = x^{\circ (n-1)} \circ (\overline{z \circ x})^{\circ (n-1)} \Big\rangle\\
& = &
\big\langle \, x,\,z  \, \mid \,  [x^{\circ n}, z]=1,\, [z, x^{\circ k}\circ z \circ \bar x^{\circ k}]=1 ~~\textrm{for}~~1 \le k \le n-2,~~ (\overline{z \circ x})^{\circ n}=\overline{x}^{\circ n} \big\rangle.
\end{eqnarray*}
Observing that
\begin{eqnarray*}
[z,x^{\circ k}\circ z\circ \bar x^{\circ k}] = x^{\circ k} \circ [\bar x^{\circ k}\circ z\circ x^{\circ k}, z]\circ \bar x^{\circ k} = x^{\circ k} \circ [z, [z, x^{\circ k}], z] \circ \bar x^{\circ k}
=x^{\circ k}\circ [z,x^{\circ k}, z] \circ \bar x^{\circ k},
\end{eqnarray*}
 we obtain 
$$
(A, \circ) = \big\langle \, x,\,z  \, \mid \,  [x^{\circ n},z] = 1,\,(z \circ x)^{\circ n}=x^{\circ n},\,
[z,x^{\circ k}, z]=1 ~~\textrm{for all}~~1 \le k \le n-2 \big\rangle.
$$
This completes the proof of the theorem. $\blacksquare$ 
\end{proof}
\para

We conclude this section by  constructing some non-isomorphic $\lambda$-cyclic left braces on $\mathbb{Z}^2$ \cite[Theorem 4.9]{MR4346001}.

\begin{theorem}\label{cyclic left braces on Z2}
Let $(\mathbb{Z}^2, +)$ be the free abelian group with free basis  $\{x_1,x_2\}$ and $\varphi \in \Aut (\mathbb{Z}^2, +)$ such that
\begin{equation}\label{varphi x cong x condition}
\varphi(x_1) \equiv x_1 \mod \mathbb{Z}^2_0 \quad \textrm{and} \quad \varphi(x_2) \equiv x_2 \mod \mathbb{Z}^2_0.
\end{equation}
Let  $\lambda :\mathbb{Z}^2 \to \Aut (\mathbb{Z}^2)$ be the group homomorphism defined by $a \mapsto \lambda_a = \varphi^{\mathfrak{L}(a)}$ for all $a \in \mathbb{Z}^2$. Then $(\mathbb{Z}^2, +, \circ)$ is a $\lambda$-cyclic left brace such that either $(\mathbb{Z}^2, \circ) \cong (\mathbb{Z}^2 , +)$ or $(\mathbb{Z}^2, \circ)$ is isomorphic to the fundamental group of the Klein bottle. 
\par
Moreover, for an automorphism of $(\mathbb{Z}^2, +)$ satisfying \eqref{varphi x cong x condition}, these are the only possibilities for the multiplicative group of a $\lambda$-cyclic left brace $(\mathbb{Z}^2, +, \circ)$.
\end{theorem}

\begin{proof}
Observe that  $\mathbb{Z}^2_0$ is generated by $x_1-x_2$. Since $\varphi$ satisfies
$$
\varphi(x_1) \equiv x_1 \mod \, \mathbb{Z}^2_0 \quad   \textrm{and} \quad \varphi(x_2) \equiv x_2 \mod \, \mathbb{Z}^2_0,
$$
we can write
$$
\varphi(x_1)=x_1+p(x_1-x_2) \quad   \textrm{and} \quad \varphi(x_2)=x_2+q(x_1-x_2)
$$
for some integers  $p$ and $q$. Since $\varphi$ is an automorphism, we have
$$
\det \big([\varphi]\big)=p-q+1=\pm 1,
$$
where $[\varphi]$ is the matrix of $\varphi$ with respect to  the basis $\{x_1, x_2\}$. Hence, we have either $q=p$ or $q=p+2$, and therefore
$$
[\varphi]=\left(%
\begin{array}{cc}
  1+p & -p \\
  p & 1-p
\end{array}%
\right)
\quad
\mbox{or}
\quad
[\varphi]=\left(%
\begin{array}{cc}
  1+p & -p \\
  p+2 & -1-p
\end{array}%
\right).
$$
\para

In the first case,  $\varphi$ has infinite order and takes the form
\begin{equation}\label{sec4eqn1}
\varphi : \left\{
\begin{array}{l}
x_1 \mapsto (1+p) x_1 - p x_2,\\
x_2 \mapsto p x_1 + (1-p) x_2.\\
\end{array}
\right.
\end{equation}
Using induction on $k$, we can prove that
$$
\varphi^k : \left\{
\begin{array}{l}
x_1 \mapsto (1 + k p) x_1 - k p x_2,\\
x_2 \mapsto k p x_1 + (1 - k p) x_2,\\
\end{array}
\right.
$$
and hence
$$
[\varphi]^k = \left(%
\begin{array}{cc}
  1 + k p & -k p \\
  k p & 1 - k p \\
\end{array}%
\right).
$$
\para

Let $a=\alpha_1x_1+\alpha_2x_2$ and  $b=\beta_1x_1+\beta_2x_2$ be elements in $\mathbb{Z}^2$. Then, we have
$$
a\circ b=a+\varphi^{\alpha_1+\alpha_2}(b) = a + b + p (\alpha_1 + \alpha_2)(\beta_1 + \beta_2)  (x_1 - x_2).
$$
We see that if $\alpha_1+\alpha_2=0$, then $a\circ b=a+b$. Note that  ${\bf 0} = 0x_1 + 0x_2$ is the identity element of $(\mathbb{Z}^2, \circ)$. Suppose that $a \circ b = {\bf 0}$. Then we have the equations
$$ \alpha_1 + \beta_1 +  p (\alpha_1+\alpha_2) (\beta_1 + \beta_2) = 0 \quad \textrm{and} \quad  \alpha_2 + \beta_2 -  p (\alpha_1+\alpha_2) (\beta_1 + \beta_2) = 0.$$
This implies that the right inverse $\bar a$ of $a$ has the form
$$
\bar a = -a + p (\alpha_1+\alpha_2)^2 (x_1 - x_2).
$$
Since $\bar a$ is also the left inverse of  $a$, it is the inverse of $a$. 
\para

We now determine the structure of $(\mathbb{Z}^2, \circ)$. Note that
$$
(\alpha_1 x_1 + \alpha_2 x_2) \circ (m x_1) = \big(\alpha_1 + m + m p(\alpha_1 + \alpha_2) \big) x_1 + \big(\alpha_2 - m p (\alpha_1 + \alpha_2) \big) x_2.
$$
Thus, if $\alpha_1 + \alpha_2 = 0$, then
$$
(\alpha_1 x_1 + \alpha_2 x_2) \circ (m x_1) = (\alpha_1 + m) x_1 + \alpha_2 x_2.
$$
This gives a decomposition
$$
\mathbb{Z}^2=\sqcup_{m\in \mathbb{Z}} ~ \mathbb{Z}^2_0 \circ (m x_1)
$$
of $\mathbb{Z}^2$ as a union of cosets of $\mathbb{Z}^2_0$. Note that the $m$-th power of  $x_1$ in $(\mathbb{Z}^2, \circ)$ has a form
$$
x_1^{\circ m}=x_1\circ \cdots \circ x_1 = \big( \id + \varphi + \cdots + \varphi^{m-1} \big)(x_1).
$$
Since
$$
I + [\varphi]+ \cdots + [\varphi]^{m-1}=
\left(%
\begin{array}{cc}
  m+\frac{m(m-1)}{2} p & -\frac{m(m-1)}{2} p \\
  & \\
  \frac{m(m-1)}{2} p & m-\frac{m(m-1)}{2}p \\
\end{array}%
\right),
$$
we obtain
$$
x_1^{\circ m}= \left(m + \frac{m(m-1)}{2} p \right) x_1  - \frac{m(m-1)}{2} p x_2,
$$
which further gives
$$
x_1^{\circ m}-mx_1 = \frac{m(m-1)}{2} p  x_1  - \frac{m(m-1)}{2} p x_2 = \frac{m(m-1)}{2} p (x_1 - x_1).
$$
We also have
\begin{equation}\label{x1 and x1-x2 relation}
 x_1 \circ (x_1 - x_2) = 2 x_1 - x_2 =  (x_1 - x_2) \circ  x_1.
 \end{equation}
This gives a decomposition
$$
(\mathbb{Z}^2, \circ) = \sqcup_{m\in \mathbb{Z}}  ~\mathbb{Z}^2_0 \circ (m x_1)
$$
 of $(\mathbb{Z}^2, \circ)$ as a union of cosets of $\mathbb{Z}^2_0$. In particular,
$(\mathbb{Z}^2, \circ)$ is generated by $\{x_1, x_1 - x_2 \}$. In view of \eqref{x1 and x1-x2 relation}, it follows that $(\mathbb{Z}^2, \circ)$ is the free abelian group of rank two.
\para

In the second  case, $\varphi$ is of order two and takes the form
\begin{equation}\label{sec4eqn2}
\varphi : \left\{
\begin{array}{l}
x_1 \mapsto (1+p) x_1 -p x_2,\\
x_2 \mapsto (2+p) x_1 - (1+p) x_2.\\
\end{array}
\right.
\end{equation}

Let  $a=\alpha_1x_1+\alpha_2x_2$ and $b=\beta_1x_1+\beta_2x_2$ be elements of $\mathbb{Z}^2$. Then, we have
$$
a\circ b = a + \varphi^{\alpha_1+\alpha_2}(b)  =
\left\{\begin{array}{ll}
  a+b & \textrm{if}~ \alpha_1+\alpha_2 \equiv 0 \mod 2,  \\
  a + b + \big(2 \beta_2 + p (\beta_1 + \beta_2)\big) (x_1 - x_2) & \textrm{if}~ \alpha_1+\alpha_2 \equiv 1 \mod 2. \\
\end{array}
\right.
$$
A direct computation along the lines of the first case shows that the inverse of $a=\alpha_1x_1+\alpha_2x_2$ in $(\mathbb{Z}^2, \circ)$ has the form
$$
\bar a =
\left\{\begin{array}{ll}
-a & \textrm{if} ~  \alpha_1+\alpha_2 \equiv 0 \mod 2 ,  \\
  -a  +  \big(2 \alpha_2 + p (\alpha_1 + \alpha_2) \big) (x_2 - x_1) & \textrm{if} ~ \alpha_1+\alpha_2\equiv 1 \mod 2. \\
\end{array}
\right.
$$
\para

Given $a=\alpha_1x_1+\alpha_2x_2 \in \mathbb{Z}^2_0$, we can write it as
$$
a= \frac{\alpha_1+\alpha_2}{2}(x_1+x_2)+\frac{\alpha_1-\alpha_2}{2}(x_1-x_2).
$$
This shows that $\mathbb{Z}^2_0$ is generated by the elements $z_1 = x_1 + x_2$ and $z_2 = x_1 - x_2$. The group $(\mathbb{Z}^2, +)$ has the coset decomposition
$$
(\mathbb{Z}^2, +)=\mathbb{Z}^2_{0} ~\sqcup~(x_1+ \mathbb{Z}^2_{0}).
$$
Since $\varphi(\mathbb{Z}^2_{0})=\mathbb{Z}^2_{0}$, we get
$$
x_1+ \mathbb{Z}^2_{0}=x_1+\varphi(\mathbb{Z}^2_{0})=x_1\circ \mathbb{Z}^2_{0}.
$$
Hence, we have the coset decomposition
$$
(\mathbb{Z}^2, \circ) =\mathbb{Z}^2_{0} ~\sqcup~(x_1\circ \mathbb{Z}^2_{0})
$$
of $(\mathbb{Z}^2, \circ )$.
\para

To determine the structure of $(\mathbb{Z}^2, \circ )$, observe that it is generated by $\{x_1, z_1, z_2\}$ and has relations
$$
z_1 \circ z_2 = z_2 \circ z_1 \quad \textrm{and}\quad x_1 \circ x_1 = z_1 + (p+1) z_2.
$$
Since $\mathbb{Z}^2_{0}$ is normal in $(\mathbb{Z}^2, \circ )$, the conjugation relations
$$
\bar x_1 \circ z_1 \circ x_1 =  z_1 + 2(1 + p) z_2~~\textrm{and}~~~\bar x_1 \circ z_2 \circ x_1 = - z_2
$$
hold. Hence, the group $(\mathbb{Z}^2, \circ)$ has the presentation
\begin{eqnarray*}
(\mathbb{Z}^2, \circ) &= &\Big\langle x_1, z_1, z_2  \, \mid \,  z_1 \circ z_2 = z_2 \circ z_1, \quad x_1 \circ x_1 = z_1 \circ z_2^{\circ (p+1)},\\
& &  \;\; \bar x_1 \circ z_1 \circ x_1 =  z_1 + 2(1 + p) z_2, \quad \bar x_1 \circ z_2 \circ x_1 = \bar z_2 \Big\rangle.
\end{eqnarray*}
Since  $x_1^{\circ 2} \circ {\bar z_2}^{\circ (p+1)} = z_1$, we obtain
$$
(\mathbb{Z}^2, \circ) = \big\langle x_1, z_2  \, \mid \,  \bar x_1 \circ z_2 \circ x_1 = \bar z_2  \big\rangle,
$$
which  is the fundamental group of the Klein bottle. $\blacksquare$
\end{proof}

\begin{remark}
{\rm 
It is conceivable that two distinct tuples of integers could result in isomorphic left braces. It would be intriguing to ascertain the conditions under which this occurs.}    
\end{remark}
\bigskip
\bigskip


\section{Solutions arising from Lambda-homomorphic skew braces} 
In this section, we discuss solutions to the Yang--Baxter equation arising from Lambda-homomorphic skew braces constructed in the preceding sections.
Our first result is the following theorem  \cite[Theorem 7.4]{MR4346001}.

\begin{theorem}\label{lambda homomorphic solution 1}
Let $A \cong (\mathbb{Z}^n, +)$ be the free abelian group with free basis $\{x_1, \ldots, x_n \}$. For elements $a = \sum_{i=1}^n \alpha_i x_i$ and $b = \sum_{j=1}^n \beta_j x_j$ of $A$, let  $\langle a, b \rangle = \sum_{i=1}^{n-1} \alpha_i \beta_i$ denote their restricted scalar product and $r: A \times A \to A \times A$ the map given by
$$
r(a, b) = \big(b~ + \langle a, b \rangle x_n, ~a~ -  \langle a, b \rangle x_n \big).
$$
Then $(A, r)$ is a non-degenerate involutive solution to the Yang--Baxter equation.
\end{theorem}

\begin{proof}
Let $\varphi_1,  \ldots, \varphi_n \in \Aut( A)$ be pairwise commuting automorphisms of $A$ given by
 $$
\varphi_i : \left\{
  \begin{array}{ll}
x_i \to x_i + x_n, & \\
x_j \to x_j & \textrm{if} \quad j \not= i
  \end{array}
\right.
$$
for $1 \le i \le n-1$ and $\varphi_n = \id$. Then the map  $\lambda : A \to \Aut( A)$ defined on the generators by
$$
\lambda_{x_i}  = \varphi_i
$$
for $1 \le i \le n$, gives a group homomorphism. Let $a = \sum_{i=1}^n \alpha_i x_i$ and $b = \sum_{j=1}^n \beta_j x_j$ be elements of $A$. Then we have 
$$
\lambda_a = \varphi_1^{\alpha_1} \varphi_2^{\alpha_2} \ldots \varphi_{n-1}^{\alpha_{n-1}} :
\left\{
  \begin{array}{ll}
x_i \to x_i + \alpha_i x_n & \textrm{if} \quad i \ne n,\\
x_n \to x_n. &
  \end{array}
\right.
$$
Thus, we can write $\lambda_a(b) = b + \langle a, b \rangle x_n$. It follows from Proposition \ref{redprop} and  Theorem \ref{t1} that
$$
a \circ b = a + \lambda_a(b) =  a + b~ + \langle a, b \rangle x_n
$$
defines a $\lambda$-homomorphic left brace $(A, +, \circ)$. By Theorem \ref{solution from a skew brace}, the left brace gives a non-degenerate involutive solution to the Yang--Baxter equation with braiding
$$
r(a, b) = \big( \lambda_a(b), ~\lambda^{-1}_{\lambda_a(b)}(a)\big).
$$
It is not difficult to see that $\lambda^{-1}_{c} = \lambda_{-c}$ for all $c \in A$. Using this, we obtain 
$$
r(a, b) = \big(b~ +  \langle a, b \rangle x_n, ~a~ -  \langle a, b \rangle x_n \big),
$$
which is desired. $\blacksquare$    
\end{proof}          
\para

Let $A \cong (\mathbb{Z}^n, +)$ be the free abelian group with free basis $\{x_1, \ldots, x_n \}$.  Consider the subgroup
$$
B = \Big\{ \sum_{i=1}^n \alpha_i x_i \in A  \, \mid \, \sum_{i=1}^n \alpha_i  = 0 \Big\}
$$
of $A$ and the automorphism $\varphi \in \Aut(A)$ such that
$$
\varphi (x_i) \equiv x_i \mod B
$$
for each $1 \le i \le n$. Then the binary operation
$$
a \circ b = a + \varphi^{\mathfrak{L}(a)}(b),
$$
for $a, b \in A$, gives a left brace $(A, +, \circ)$. By Theorem \ref{solution from a skew brace}, the braiding given by
$$
r(a, b) = \big(\varphi^{\mathfrak{L}(a)}(b),~ \varphi^{-\mathfrak{L}(b)}(a)\big),
$$
gives a non-degenerate involutive solution $(A, r)$ to the Yang--Baxter equation.  The next result extends the preceding observation \cite[Proposition 7.8]{MR4346001}. 

\begin{proposition}
Let $G$ be an abelian group admitting a homogeneous presentation $G= \langle S\mid R\rangle$. Let $\varphi \in \Aut (G)$ be such that $\mathfrak{L}(\varphi (g)) = 1$ for all $g \in S$. Then the map
$r : G \times G \to G \times G$ given by
$$
r(a, b) = \big(\varphi^{\mathfrak{L}(a)} (b), ~\varphi^{\mathfrak{L}(\varphi^{\mathfrak{L}(a)}(b^{-1}))} (a)\big)
$$
gives a non-degenerate involutive solution $(G, r)$ to the Yang--Baxter equation.
\para 
In addition, if $\varphi$ satisfies $\mathfrak{L} \big(\varphi(b)\big) = \mathfrak{L}(b)$ for each $b \in G$, then
$$
r(a, b) = \big(\varphi^{\mathfrak{L}(a)} (b), ~\varphi^{-\mathfrak{L}(b)} (a)\big).
$$
\end{proposition}

\begin{proof}
Since $G$ is homogeneous and $\mathfrak{L}(\varphi (g_i)) = 1$ for all $g_i \in S$, it follows from Proposition \ref{lcb} that $(G, \cdot, \circ)$ 
is a $\lambda$-cyclic skew left brace. Here, the map $\lambda : G \to \Aut(G^{(\cdot)})$ is given by $\lambda_a (b) = \varphi^{\mathfrak{L}(a)} (b)$. Hence, the first component of the induced braiding has the desired form. For the second component, note that
$$
\lambda_a^{-1} (b) = \varphi^{-\mathfrak{L}(a)} (b) = \varphi^{\mathfrak{L}(a^{-1})} (b) = \lambda_{a^{-1}} (b).
$$
Thus, we have
$$
\lambda_{\lambda_a (b)}^{-1} (a) = \lambda_{\varphi^{\mathfrak{L}(a)} (b)}^{-1} (a) =
\lambda_{\varphi^{\mathfrak{L}(a)} (b^{-1})} (a) =
\varphi^{\mathfrak{L}(\varphi^{\mathfrak{L}(a)}(b^{-1}))} (a),
$$
which is desired. The last assertion is immediate. $\blacksquare$    
\end{proof}          
\para

Next, we construct some specific solutions on $\mathbb{Z}^2$. 

\begin{proposition}
Let $A \cong (\mathbb{Z}^2, +)$, $p$ an integer and $r:A \times A \to A \times A$ the map given by
$$r(a, b) = \big(b + p \mathfrak{L}(a) \mathfrak{L}(b) (x_1 - x_2),~ a - p \mathfrak{L}(a) \mathfrak{L}(b) (x_1 - x_2)\big).$$
Then  $(A, r)$ is a non-degenerate involutive solution  to the Yang--Baxter equation.
\end{proposition}

\begin{proof}
Consider the automorphism $\varphi$ of $A$ given by 
$$
\varphi : \left\{
  \begin{array}{l}
x_1 \to (1+p) x_1 - p x_2, \\
x_2 \to p x_1 + (1-p) x_2.
  \end{array}
\right.
$$
By Theorem \ref{cyclic left braces on Z2}, $(A, +, \circ)$ is a $\lambda$-cyclic left brace, where  $\lambda :A \to \Aut (A)$ is the group homomorphism given by $a \mapsto \lambda_a = \varphi^{\mathfrak{L}(a)}$ for all $a \in A$. By Theorem \ref{solution from a skew brace}, the left brace gives a non-degenerate involutive solution to the Yang--Baxter equation with braiding
$$
r(a, b) = \big( \lambda_a(b), ~\lambda^{-1}_{\lambda_a(b)}(a)\big).
$$
For elements $a = \alpha_1 x_1 + \alpha_2 x_2$ and $b = \beta_1 x_1 + \beta_2 x_2$ of $A$, we have
$$ \lambda_a(b) = \varphi^{\alpha_1 + \alpha_2}(b) = b + p (\alpha_1  + \alpha_2) (\beta_1 + \beta_2)(x_1 - x_2)=b + p \mathfrak{L}(a) \mathfrak{L}(b) (x_1 - x_2) $$
and 
$$ \lambda_{\lambda_a (b)}^{-1} (a) =\varphi^{\mathfrak{L}(\varphi^{\mathfrak{L}(a)}(b^{-1}))} (a)=a - p (\alpha_1  + \alpha_2) (\beta_1 + \beta_2) (x_1 - x_2) =a - p \mathfrak{L}(a) \mathfrak{L}(b) (x_1 - x_2).$$
This completes the proof.  $\blacksquare$
\end{proof}

Considering the automorphism $\varphi$ of $A \cong (\mathbb{Z}^2, +)$ given by 
$$
\varphi : \left\{
  \begin{array}{l}
x_1 \to (1+p) x_1 - p x_2, \\
x_2 \to (2 + p) x_1 - (1+p) x_2,
  \end{array}
\right.
$$
gives the following result.

\begin{proposition}
Let $A \cong (\mathbb{Z}^2, +)$ and $p$ an integer. For elements $a = \alpha_1 x_1 + \alpha_2 x_2$ and $b = \beta_1 x_1 + \beta_2 x_2$ of $A$, let $r:A \times A \to A \times A$ be the map given by
$$
r(a, b) =  \left\{
  \begin{array}{ll}
(b, a) & \textrm{if}~ \mathfrak{L}(a) \equiv \mathfrak{L}(b) \equiv 0 \mod~2, \\
\big(b + (2 \beta_2 + p \mathfrak{L}(b)) (x_1 - x_2),~ a \big)  &\textrm{if} ~ \mathfrak{L}(a) \equiv 1 \mod~2~\textrm{and}~\mathfrak{L}(b) \equiv 0 \mod~2,\\
\big(b, ~a + (2 \alpha_2 - p \mathfrak{L}(b)) (x_1 - x_2) \big)  &\textrm{if} ~ \mathfrak{L}(a) \equiv 0 \mod~2~\textrm{and}~\mathfrak{L}(b) \equiv 1 \mod~2,\\
\big(b + (2 \beta_2 + p \mathfrak{L}(b)) (x_1 - x_2), ~a + (2 \alpha_2 - p \mathfrak{L}(b)) (x_1 - x_2) \big)  &\textrm{if} ~ \mathfrak{L}(a) \equiv \mathfrak{L}(b) \equiv 1 \mod~2.\\
  \end{array}
\right.
$$
Then  $(A, r)$ is a non-degenerate involutive solution  to the Yang--Baxter equation.
\end{proposition}
\para

Similarly, considering the automorphism $\varphi$ of $A \cong (\mathbb{Z}^2, +)$ given by 
$$
\psi : \left\{
  \begin{array}{l}
x_1 \to  x_{2},  \\
x_2 \to  x_1,
  \end{array}
\right.
$$
gives the following result.

\begin{proposition}
Let $A \cong (\mathbb{Z}^2, +)$. For elements $a = \alpha_1 x_1 + \alpha_2 x_2$ and $b = \beta_1 x_1 + \beta_2 x_2$ of $A$, let $r:A \times A \to A \times A$ be the map given by
$$
r(a, b) =  \left\{
  \begin{array}{ll}
(b, a) & \textrm{if}~ \mathfrak{L}(a) \equiv \mathfrak{L}(b) \equiv 0 \mod~2, \\
\big(b + (\beta_1 - \beta_2) (x_1 - x_2), ~a\big)  &\textrm{if}~ \mathfrak{L}(a) \equiv 1 \mod~2~\textrm{and}~\mathfrak{L}(b) \equiv 0 \mod~2,\\
\big(b, ~a + (\alpha_1 + \alpha_2) (x_1 - x_2)\big)   &\textrm{if}~ \mathfrak{L}(a) \equiv 0 \mod~2~\textrm{and}~\mathfrak{L}(b) \equiv 1 \mod~2,\\
\big(b + (\beta_1 - \beta_2) (x_1 - x_2), ~a + (\alpha_1 + \alpha_2) (x_1 - x_2)\big)  &\textrm{if}~ \mathfrak{L}(a) \equiv \mathfrak{L}(b) \equiv 1 \mod~2.\\
  \end{array}
\right.
$$
Then  $(A, r)$ is a non-degenerate involutive solution  to the Yang--Baxter equation.
\end{proposition}
\para

As a final example,  we consider the solution arising from a left brace with non-solvable multiplicative group \cite[Proposition 7.16]{MR4346001}.

\begin{proposition}
Let $\mathbb{Z}\langle\langle X_1, \ldots, X_n \rangle\rangle$ be the ring of formal power series in non-commuting variables $X_1, \ldots, X_n$ over $\mathbb{Z}$, where $n \ge 2$. Let $B$ be the two-sided ideal of $\mathbb{Z}\langle\langle X_1, \ldots, X_n \rangle\rangle$  generated by $X_1, \ldots, X_n$. Let $r:B \times B \to B \times B$ be the map given by
$$
r(a, b) = \big(b + a b,~ a + (\overline{b+ba}) a \big)
$$
for $a, b \in B$, where $\overline{c}= -c + c^2 - c^3 + \cdots$. Then $(B, r)$ is a non-degenerate involutive solution to the Yang--Baxter equation.
\end{proposition}

\begin{proof}
It follows from the proof of Proposition \ref{p7.1} that  $(B, +, \circ)$ is a two-sided brace with
$$
a \circ b = ab + a + b
$$
for $a, b \in B$. We know that the left brace gives a non-degenerate involutive solution with braiding
$$
r(a, b) = \big(\lambda_a (b), ~\lambda^{-1}_{\lambda_a (b)} (a)\big)
$$
for all $a, b \in B$. In this case, we have
$$
\lambda_a (b) = -a + a \circ b = b + ab.
$$
Further, $\lambda_a^{-1} = \lambda_{\overline{a}}$, where
$$
\overline{a} = -a + a^2 - a^3 + \cdots.
$$
 Hence, we have
$$
r(a, b) = \big(\lambda_a (b), ~\lambda^{-1}_{\lambda_a (b)} (a)\big)  = \big(b + a b,~ a + (\overline{b+ba}) a \big)
$$
for all $a, b \in B$, which is desired. $\blacksquare$
\end{proof}          
\bigskip


\chapter{Rota--Baxter operators on groups }\label{chap RBO on groups}

\begin{quote}
This chapter develops the foundational theory of Rota--Baxter operators on groups. We begin by introducing the definition, key examples, and basic properties of Rota--Baxter groups. Various constructions of these structures are then explored in detail.  Furthermore, we examine the problem of extending set-theoretic maps to Rota--Baxter operators on groups. We also investigate extensions of Rota--Baxter groups and examine the relationship between Rota--Baxter operators on groups and their corresponding Lie rings.
\end{quote}
\bigskip

\section{Properties of Rota--Baxter groups} 
The first appearance of Rota--Baxter operators on commutative algebras can be traced back to Baxter's seminal paper \cite{MR0119224} in 1960. Since then, the theory of Rota--Baxter operators has undergone extensive development by various authors, with comprehensive details available in the monograph \cite{MR3025028}. These operators play a pivotal role due to their connections with several mathematical concepts, including the Yang--Baxter equation \cite{MR0674005,MR0725413}, Loday algebras \cite{MR3021790}, double Poisson algebras \cite{MR2568415}, and others. Rota--Baxter operators on Lie algebras were first studied by Belavin, Drinfel'd, and Semenov-Tian-Shansky as operator forms of the classical Yang--Baxter equation.  In their work \cite{MR4271483}, Guo, Lang and Sheng  introduced the notion of a Rota--Baxter operator on a group.  Integrating the Rota--Baxter operators on Lie algebras, they introduced the notion of Rota--Baxter operators on Lie groups and more generally on (abstract) groups.
\para

Let $A$ be an algebra over a field~$\Bbbk$. A linear operator $R:A \to A$ is called
a Rota--Baxter operator of weight~$\lambda\in \Bbbk$ if
\begin{equation}\label{RBAlgebra}
R(x)R(y) = R\big( R(x)y + xR(y) + \lambda xy \big)
\end{equation}
for all $x,y\in A$. An algebra endowed with a Rota--Baxter operator is called a \index{Rota--Baxter algebra}{\it Rota--Baxter algebra}. Following \cite{MR4271483}, we consider an analogue of a Rota--Baxter operator of weight~$\pm1$ on a group.

\begin{definition}
Let $G$ be a group.
\begin{enumerate}
\item A map $B: G\to G$ is called a \index{Rota--Baxter operator}{Rota--Baxter operator} of weight~1 if
\begin{equation}\label{RB}
B(g)B(h) = B\big(g B(g) h B(g)^{-1} \big)
\end{equation}
for all $g,h\in G$.
\item A map $C: G\to G$ is called a {\it Rota--Baxter operator} of weight~$-1$ if
$$
C(g) C(h) = C\big( C(g) h C(g)^{-1} g \big)
$$
for all $g,h\in G$.
\end{enumerate}
\end{definition}

A group endowed with a Rota--Baxter operator is called a \index{Rota--Baxter group}{\it Rota--Baxter group}. 

\begin{example}{\rm 
Let $G$ be a group.
\begin{enumerate}
\item The map given by $B_0(g) = 1$ is a Rota--Baxter operator on $G$.
\item The map given by $B_{-1}(g) = g^{-1}$ is a  Rota--Baxter operator on $G$.
\end{enumerate}}
\end{example}

The Rota--Baxter operators $B_0$ and $B_{-1}$ are referred to as \index{elementary Rota--Baxter operator}{\it elementary Rota--Baxter operators}. A~group $G$ is called \index{Rota--Baxter elementary group}{\it Rota--Baxter elementary} if any Rota--Baxter operator on $G$ is elementary.
\para

The following result provides a method for constructing new Rota--Baxter operators from existing ones \cite[Proposition 2.4]{MR4271483}.

\begin{proposition}\label{RBO from old ones}
Let $G$ be a group, $B: G \to G$ a Rota--Baxter operator of weight~1 and $C: G \to G$ a Rota--Baxter operator of weight $-1$. Then the following  assertions hold:
\begin{enumerate}
\item The map $\widetilde{C}:G \to G$ given by $\widetilde{C}(g) = g C(g^{-1})$ is a Rota--Baxter operator of weight $-1$.
\item The map $\widetilde{B}:G \to G$ given by $\widetilde{B}(g) = g^{-1}B(g^{-1})$
is a Rota--Baxter operator of weight $1$.
\item The map $\bar{B}:G \to G$  given by $\bar{B}(g) = B(g^{-1})$ is a Rota--Baxter operator of weight $-1$.
\item The map $B_{+}:G \to G$  given by $B_{+}(g) = g B(g)$ is a Rota--Baxter operator of weight $-1$.
\end{enumerate}
\end{proposition}

\begin{proof}
All the assertions follow by direct calculations. For instance, we have
$$
B_{+}(g)B_{+}(h)= g B(g) h B(h)
$$
and
$$B_{+}\big( B_{+}(g) h B_{+}(g)^{-1} g \big)=B_{+}\big( gB(g) h B(g)^{-1} \big)= \big( gB(g) h B(g)^{-1} \big)B\big(gB(g) h B(g)^{-1} \big)$$
for all $g, h \in G$. Since $B$ is a Rota--Baxter operator of weight 1, we have $B_{+}(g)B_{+}(h)= B_{+}\big( B_{+}(g) h B_{+}(g)^{-1} g \big)$,
which is assertion (4).
$\blacksquare$ 
\end{proof}

Proposition \ref{RBO from old ones} shows that there is a bijection between Rota--Baxter operators of weight 1 and $-1$ on a given group. Thus, we consider only Rota--Baxter operators of weight 1, and refer them simply as Rota--Baxter operators.
\para

Given a group $G$ and subgroups $H, L$ of $G$,  let $[H,L]$ be the subgroup of $G$ generated by the set
$$\big\{[h,l] = h^{-1}l^{-1}hl  \, \mid \,  h\in H~\textrm{and}~l\in L \big\}.$$
Next, we observe some elementary facts about Rota--Baxter operators on groups.

\begin{lemma}\label{lem:elementary}
Let $B$ be a Rota--Baxter operator on a group $G$. Then the following  assertions hold:
\begin{enumerate}
\item $B(1) = 1$.
\item $B(g)B(g^{-1}) = B\big([g^{-1},B(g)^{-1}] \big)$ for all $g \in G$.
\item $B(g)B(B(g)) =  B\big(g B(g)\big)$  for all $g \in G$.
\item  $B(g)^{-1} = B \big(B(g)^{-1}g^{-1}B(g)\big)$  for all $g \in G$. 
\end{enumerate}
\end{lemma}

\begin{proof}
Assertion (1) follows from~\eqref{RB} by taking $g = h = 1$. Assertion (2) follows from~\eqref{RB} by taking $h = g^{-1}$, and assertion (3) follows by taking $h = B(g)$. Finally, since
$$B(g)B \big(B(g)^{-1}g^{-1}B(g)\big)= B\big(g B(g)B(g)^{-1}g^{-1}B(g)B(g)^{-1}\big) = B(1) = 1,$$
we obtain assertion (4).
$\blacksquare$ 
\end{proof}

The following result is a simple yet important observation \cite[Lemma 3.2]{MR4271483}.

\begin{lemma}\label{KerImSub}
Let $B$ be a Rota--Baxter operator on a group~$G$. Then the following  assertions hold:
\begin{enumerate}
\item $\ker(B)$ is a subgroup of $G$.
\item $\im (B)$ is a subgroup of $G$.
\end{enumerate}
\end{lemma}

\begin{proof}
Recall from Lemma~\ref{lem:elementary}(1) that $B(1) = 1$. For $g,h\in \ker(B)$, by~\eqref{RB}, we have
$$
B(g^{-1}) = B(g)B(g^{-1}) = B \big(gB(g)g^{-1}B(g)^{-1}\big) = B(g g^{-1})=B(1) = 1$$
and
$$B(gh)= B \big(gB(g)hB(g)^{-1}\big)= B(g)B(h) =  1,
$$
which proves assertion (1).
\para 
It is clear from \eqref{RB} that if $g,h\in \im(B)$, then $gh\in \im (B)$. Further, if $g\in \im (B)$, then by Lemma \ref{lem:elementary}(4), we have $g^{-1}\in \im (B)$. This proves assertion (2).
$\blacksquare$ 
\end{proof}

A direct check yields the following observation.

\begin{lemma} \label{kernelCosets}
Let $B$ be a Rota--Baxter operator on a group~$G$. If $B(g) =  1$ for some $g\in G$, then $B(h) = B(gh)$ for any $h\in G$. In particular, if $G = \sqcup_{i \in I}~ \ker(B)g_i$ is the decomposition of $G$ as right cosets of $\ker(B)$, then $B(x) = B(y)$ if $x$ and $y$ lie in the same coset.
\end{lemma}

The next result is crucial from the point of  view of classification of Rota--Baxter operators on a given group \cite[Lemma 9 and Lemma 10]{MR4556953}.

\begin{lemma}\label{lem:Aut}
Let $B$ be a Rota--Baxter operator on a group~$G$ and $\varphi \in \Aut(G)$ an automorphism. Then the following  assertions hold:
\begin{enumerate}
\item $B^{(\varphi)} = \varphi^{-1}B\varphi$
is a~Rota--Baxter operator on $G$.
\item $(\widetilde{B})^{(\varphi)} = \widetilde{B^{(\varphi)}}$.
\end{enumerate}

\end{lemma}

\begin{proof}
Direct computations give
\begin{eqnarray*}
\varphi \big( B^{(\varphi)}(g) B^{(\varphi)}(h) \big)
 &=& \varphi \big( \varphi^{-1}(B(\varphi(g))) \,\varphi^{-1}(B(\varphi(h))) \big)\\
 &=&  B \big(\varphi(g)\big) \, B \big(\varphi(h)\big) \\
 &=&  B\big(\varphi(g)B(\varphi(g))\varphi(h)B(\varphi(g))^{-1} \big)\\
 &=&  \varphi\varphi^{-1}B\varphi \big( g\varphi^{-1}B(\varphi(g))h\varphi^{-1}(B(\varphi(g))^{-1}) \big) \\
 &=&  \varphi\big( B^{(\varphi)} \big(gB^{(\varphi)}(g)h(B^{(\varphi)}(g))^{-1} \big) \big)
\end{eqnarray*}
for all $g, h \in G$, which is assertion (1).
\para 
Further, we see that
$$
(\widetilde{B})^{(\varphi)} (g)
 = \varphi^{-1} \big( \widetilde{B} (\varphi(g)) \big)
 = \varphi^{-1} \big( \varphi(g)^{-1}B( \varphi(g)^{-1} ) \big)
 = g^{-1}B^{(\varphi)}(g^{-1})
 = \widetilde{B^{(\varphi)}}(g), 
$$
which is assertion (2). $\blacksquare$ 
\end{proof}

The following result shows how a Rota--Baxter operator gives rise to an alternative group structure  \cite[Section 3]{MR4271483}.

\begin{proposition}\label{prop:Derived}
Let $B$ be a Rota--Baxter operator on a group~$(G, \cdot)$. Then the following  assertions hold:
\begin{enumerate}
\item The set $G$ with the binary operation
\begin{equation}\label{R-product}
g\circ h = gB(g)hB(g)^{-1}
\end{equation}
for $g,h\in G$, is a group.
\item The map $B: (G,\circ) \to (G,\cdot)$ is a homomorphism of groups.
\item The map $B$ is a Rota--Baxter operator on the group $(G,\circ)$.
\end{enumerate}
\end{proposition}

\begin{proof}
Let $g, h, k \in G$. Using \eqref{RB}, we see that
$$(g \circ h) \circ k=  \big(gB(g)hB(g)^{-1} \big) \circ k= g B(g) h B(h) k B(h)^{-1} B(g)^{-1} =g \circ \big(h B(h) k B(h)^{-1} \big)= g \circ (h \circ k).$$
Further, if 1 is the identity element of $(G, \cdot)$, then we see that $g \circ 1=g= 1 \circ g$ and $g^{\circ (-1)}= B(g)^{-1} g^{-1} B(g)$ for all $g \in G$. Thus, $(G, \circ)$ is a group, which proves assertion (1).
\para
If $g, h \in G$, then using \eqref{RB}, we have
$$B(g \circ h)=B \big(gB(g)hB(g)^{-1}\big)= B(g) B(h).$$
Thus, $B: (G, \circ) \to (G, \cdot)$ is a group homomorphism, which proves assertion (2). 
\para 
Using assertion (2), we see that $$B(g) \circ B(h)= B(g) B \big(B(g)\big) B(h) B \big(B(g)\big)^{-1}= B \big(g \circ B(g) \circ h \circ B(g)^{\circ (-1)}\big),$$
which is assertion (3). $\blacksquare$ 
\end{proof}

We shall denote the group $(G,\circ)$ by $G_B$. Given a Rota--Baxter group $(G,B)$,  consider the map $B_+:G \to G$ defined as $B_+(g) = gB(g)$ for $g \in G$. Then $B_+$ is a group homomorphism from $G_B$ to $G$. It follows from Lemma~\ref{lem:elementary}(3)  that $B_+B = BB_+$.

\begin{exercise}{\rm 
Let $B$ be a Rota--Baxter operator on a group $(G, \cdot)$. Then the following  assertions hold:
\begin{enumerate}
\item  If $B = B_0$, then
$$
g\circ h = gB_0(g)hB_0(g)^{-1} = g h.
$$
Hence, the group operations $\circ$ and $\cdot$ coincide.
\item  If $B = B_{-1}$, then
$$
g\circ h = gB_{-1}(g) h B_{-1}(g)^{-1} = h g.
$$
Hence, $(G, \circ)$ is the opposite group of $(G, \cdot)$.
\item  If $G$ is abelian, then $g\circ h = g\cdot h$ for all $g,h\in G$.    
\end{enumerate}}
\end{exercise}

The following result of Guo, Lang, and Sheng \cite[Proposition 3.3]{MR4271483} is analogous
to the one established by Semenov-Tyan-Shanskii for Lie algebras \cite{MR0725413}.

\begin{proposition} Let $(G, B)$ be a Rota--Baxter group. Then the following  assertions hold: 
\begin{enumerate}
\item $\ker(B)$ and $\ker(B_+)$ are normal subgroups of $G_B$.
\item $\ker(B)$ is a normal subgroup of $\im(B_+)$.
\item $\ker(B_+)$ is a normal subgroup of $\im(B)$.
\item  There is an isomorphism of groups
$$ \im(B_+)/\ker(B)\cong \im(B)/\ker(B_+). $$
\item There is a factorization of $G$ as
\begin{equation}\label{ImageFactorization}
G = \im(B_+)\im(B).
\end{equation}
\end{enumerate}
\end{proposition}
\para

Let $(G, B)$ be a Rota--Baxter group and $a \in G$. Given an integer~$k$, let $a^{\circ(k)}$ denote the $k$-th power of $a$ with respect to $\circ$. The following result gives a precise formula for a word written in terms of the group operation in $G_B$ \cite[Proposition 14]{MR4556953}. 

\begin{proposition} \label{for} Let $(G, B)$ be a Rota--Baxter group and $A$ a subset of $G$. Let
$$
w = a_{i_1}^{\circ(k_1)} \circ a_{i_2}^{\circ(k_2)} \circ \cdots
    \circ a_{i_s}^{\circ(k_s)},
$$
where $a_{i_j} \in A$ and $k_j \in \mathbb{Z}$, be an element presented with respect to the group operation $\circ$. Then
$$
w = \big(a_{i_1}B(a_{i_1})\big)^{k_1}
    \big(a_{i_2}B(a_{i_2})\big)^{k_2} \ldots
    \big(a_{i_s}B(a_{i_s})\big)^{k_s}
    B(a_{i_s})^{-k_s}
    B(a_{i_{s-1}})^{-k_{s-1}} \ldots
    B(a_{i_1})^{-k_1}.
$$
In particular, $a^{\circ(-1)} = B(a)^{-1} a^{-1} B(a)$ and $a^{\circ(k)} = B_+(a)^k B(a)^{-k}$ for each $k \in \mathbb{Z}$.
\end{proposition}

\begin{proof}
We first show using induction on $n$ that $a^{\circ(n)} = B_+(a)^n B(a)^{-n}$. Indeed, for $n = 0$, the assertion is vacuous. The induction step follows from the equalities
$$
a^{\circ(n)}
 = a\circ (a^{\circ(n-1)})
 = aB(a)a^{\circ(n-1)}B(a)^{-1}
 = \big(aB(a)\big)^n B(a)^{-n}=B_+(a)^n B(a)^{-n}.
$$
We can check directly that $B(a)^{-1}a^{-1}B(a)$
is the inverse element to~$a$ in the group $G_B$.
Again, using induction on $n$ and
Lemma~\ref{lem:elementary}(4), we obtain
\begin{eqnarray*}
a^{\circ(-n)}
 &=& a^{\circ(-1)}\circ(a^{\circ(-(n-1))}) \\
&=& B(a)^{-1}a^{-1} B(a) \, B\big( B(a)^{-1}a^{-1}B(a) \big) \big(aB(a)\big)^{-(n-1)} B(a)^{n-1} \,\big(B( B(a)^{-1}a^{-1}B(a) ) \big)^{-1} \\
&=& \big(aB(a)\big)^{-1} \big(aB(a)\big)^{-(n-1)} B(a)^{n-1}B(a)\\
 &=& \big(aB(a)\big)^{-n} B(a)^{n}\\
&=& B_+(a)^{-n} B(a)^{n}.
\end{eqnarray*}
We write $w' = a_{i_2}^{\circ(k_2)} \circ \cdots \circ a_{i_s}^{\circ(k_s)}$.
Using the fact that $B: G_B \to G$ is a group homomorphism, we get
$$
w = a_{i_1}^{\circ(k_1)}\circ w'
  = \big(B_+(a_{i_1}) \big)^{k_1} \big(B(a_{i_1})\big)^{-k_1}
  B(a_{i_1}^{\circ(k_1)})w'B(a_{i_1}^{\circ(k_1)})^{-1}
  = \big(B_+(a_{i_1})\big)^{k_1}w' \big(B(a_{i_1})\big)^{-k_1}.
$$
Thus, using induction on~$s$, we derive the desired expression.
$\blacksquare$ 
\end{proof}
\bigskip
\bigskip

\section{Constructions of Rota--Baxter operators on groups}

In this section, we present some constructions of Rota--Baxter operators on groups. 

\subsection{Constructions via factorizations}
The following example allows us to construct non-elementary Rota--Baxter operators.

\begin{example}\label{exm:split} {\rm 
Let $G$ be a group with an exact factorization $G = HL$. Then the map $B: G\to G$ defined as
$$
B(hl) = l^{-1}
$$
is a Rota--Baxter operator on~$G$.}
\end{example}

We refer to such Rota--Baxter operators as ~\index{splitting Rota--Baxter operator}{\it splitting Rota--Baxter operators}. In this case, it follows from the definition of the group operation in $G_B$ that  $G_B\cong H\times L^{\rm op}$, where $L^{\rm op}$ is the opposite group of $L$. We now have the following result \cite[Proposition 16]{MR4556953}.

\begin{proposition}\label{SplittingCond}
Let $B$ be a Rota--Baxter operator on a group $G$. Then $B$ is a splitting Rota--Baxter operator on~$G$
if and only if $B \big(gB(g) \big) = 1$ for all $g\in G$.
\end{proposition}

\begin{proof}
Suppose that $B$ is a splitting Rota--Baxter operator on~$G$. Then there exist subgroups $H, L$ of $G$ giving an exact factorisation $G=HL$.
Let $g = hl \in G$, where $h\in H$ and $l\in L$. Then we have
$$
B \big(gB(g) \big) = B \big(hlB(hl) \big) = B(hll^{-1}) = B(h) = 1.
$$

Conversely, suppose that $B \big(gB(g)\big) = 1$ for all $g\in G$. We claim that $G = \ker(B)\im(B)$ is an exact factorization. Firstly, let $a\in \ker(B)\cap\im(B)$. Then $B(a) = 1$ and there exists $g\in G$ such that $a = B(g)$.
By the hypothesis, we have $1 = B \big(gB(g)\big) = B(ga)$. Since $a, ga\in \ker(B)$, by Lemma~\ref{KerImSub}(1) we conclude that $g\in \ker(B)$. Thus, we have $a = 1$. Secondly, let $g\in G$. Then the equality $B \big(gB(g)\big) = 1$ implies that $gB(g) \in\ker(B)$.
Thus, by Lemma~\ref{KerImSub}(2), we can write $g = ab$ for some $a \in \ker(B)$ and $b = B(g)^{-1}\in \im(B)$. $\blacksquare$ 
\end{proof}

\begin{corollary}
Let $G = H\times H\times L$ be a direct product of groups, where $H$ is a non-trivial group. Then the map $B: G\to G$ defined by
$B \big((h_1,h_2,l)\big) = (1,h_1,1)$ is a non-splitting Rota--Baxter operator on $G$.
\end{corollary}

\begin{proof}
A direct check shows that $B$ is a Rota--Baxter operator on $G$. Let $g = (h_1,h_2,l) \in G$. Since $B \big(gB(g)\big) = (1,h_1,1) \ne (1, 1, 1)$ for $h_1\neq 1$, it follows from Proposition~\ref{SplittingCond} that  $B$ is non-splitting. $\blacksquare$ 
\end{proof}

\begin{remark}{\rm 
By Lemma~\ref{lem:elementary}, the condition $B\big(gB(g)\big) = 1$ is equivalent to the relation $B(g) B^2 (g) = 1$ or $B^2 (g) = B(g)^{-1}$. It means that if $B$ is a splitting Rota--Baxter operator, then it inverts all elements of $\im(B)$.}
\end{remark}

The following result extends Example~\ref{exm:split} to  triple factorizations \cite[Proposition 18]{MR4556953}.

\begin{proposition}\label{triangular}
Let $G$ be a group such that $G = HLM$, where $H$, $L$ and $M$ are subgroups of $G$ with pairwise trivial intersections.
Let $C$ be a Rota--Baxter operator on~$L$ and suppose that $[H,L] = [C(L),M] = 1$.
Then the map $B: G\to G$ defined by
$$
B(hlm) = C(l)m^{-1}
$$
is a Rota--Baxter operator on~$G$.
\end{proposition}

\begin{proof}
Let $h,h'\in H$, $l,l'\in L$ and $m,m'\in M$. Then we have
\begin{eqnarray*}
B(hlm)B(h'l'm') &=& C(l)m^{-1}C(l')m'^{-1}\\
&=& C(l)C(l')m^{-1}m'^{-1}\\
&=& C \big(lC(l)l'C(l)^{-1} \big)(m'm)^{-1}, \quad \textrm{since $C$ is a Rota--Baxter operator}\\
&=& B \big(hh' l C(l)l'C(l)^{-1}m'm \big), \quad \textrm{by definition of $B$}\\
 &=& B \big(hlm C(l) m^{-1} h'l'm' m C(l)^{-1} \big), \quad \textrm{since $[H,L] = [C(L),M] = 1$}\\
 &=& B \big(hlm B(hlm) (h'l'm')B(hlm)^{-1} \big),
\end{eqnarray*}
and hence $B$ is a Rota--Baxter operator on $G$.
$\blacksquare$ 
\end{proof}

Note that, in the case of Proposition \ref{triangular}, we have 
$$G_B\cong H\times L_C\times M^{\rm op},$$ where $M^{\rm op}$ is the opposite group of $M$.

\begin{corollary}
Let $G = \langle a,b \rangle$ be a 2-generated 2-step nilpotent group and let $c = [b,a]$. Then the following assertions hold:
\begin{enumerate}
\item  $G = \langle a \rangle \langle c \rangle \langle b \rangle$.
\item For any integer~$k$, the map $B:G \to G$ given by
$$
B(a^{\alpha} c^{\beta} b^{\gamma})
 = c^{k\beta} b^{-\gamma}
$$
for $\alpha,\beta,\gamma \in \mathbb{Z}$, is a Rota--Baxter operator on~$G$.
\end{enumerate}
\end{corollary}

Given a~splitting Rota--Baxter operator $B$ on a~group~$G$, we have an exact factorization
\begin{equation}\label{KerImFactor}
G = \ker(B)\im(B).
\end{equation}
It is easy to see that if $G$~is abelian, then
there exists a non-splitting Rota--Baxter operator~$B$ on~$G$ such that~\eqref{KerImFactor} holds. Indeed, if $G = H \times L$, where $H, K$ are abelian groups, then any map $hl\mapsto \varphi(l)$, where $\varphi$~is an endomorphism of~$L$, defines a Rota--Baxter operator on~$G$ satisfying~\eqref{KerImFactor}. The following result generalizes this observation \cite[Proposition 20]{MR4556953}.

\begin{proposition}\label{semi-direct}
Let $G = H\rtimes L$ be a semi-direct product and let $C$ be a~Rota--Baxter operator on $L$. Then the map $B: G\to G$ defined by  $B(hl) = C(l)$, where $h\in H$ and $l\in L$, is a~Rota--Baxter operator on $G$.
\end{proposition}

\begin{proof}
Let $h,h'\in H$ and $l,l'\in L$. Then we have
\begin{eqnarray*}
B(hl)B(h'l') &=& C(l)C(l')\\
&=& C \big(lC(l)l'C(l)^{-1} \big)\\
&=& B \big(hh'^{C(l)^{-1}l^{-1}}l C(l)l'C(l)^{-1} \big) \\
&=& B \big(hl C(l)h'l'C(l)^{-1} \big)\\
&=& B \big(hlB(hl)h'l'B(hl)^{-1} \big),
\end{eqnarray*}
which is desired. $\blacksquare$ 
\end{proof}

Note that, in the case of Proposition \ref{semi-direct}, we have $G_B\cong H\rtimes L_C$.
\bigskip
\bigskip

\subsection{Constructions via homomorphisms and exact formulas}

Note that, if $G$ is an abelian group, then each endomorphism of $G$ is a Rota--Baxter operator on~$G$.

\begin{proposition}\label{prop:RB-Hom}
The following assertions hold:
\begin{enumerate}
 \item If $G$ is a group and $H$ is an abelian subgroup of $G$, then any homomorphism (or anti-homomorphism) $B: G \to H$ is a Rota--Baxter operator on $G$.
 \item If $B$ is a Rota--Baxter operator on $G$ which is also an automorphism of $G$, then  $G$ is abelian.
\end{enumerate}
\end{proposition}

\begin{proof}
If $H$ is abelian and $B: G \to H$ is a homomorphism, then
$$
B \big(gB(g)hB(g)^{-1} \big)= B(g)B^2(g) B(h) \big(B^2(g)\big)^{-1}= B(g)B(h).  
$$
The case of anti-homomorphism is similar, which proves assertion (1).
\para

If $B$ is a Rota--Baxter operator on $G$ which is also an endomorphism, then
$$
B(g)B(h)
 = B \big(gB(g)hB(g)^{-1} \big)
 =B(g)B^2(g) B(h) \big(B^2(g)\big)^{-1}
$$
for all $g,h \in G$. It is equivalent to the relation
\begin{equation}\label{RB-Hom}
[B^2(g), B(h)] = 1.
\end{equation}
for all $g,h \in G$. Thus, if $B$ is an automorphism, then $G$ is abelian, which proves assertion (2).  $\blacksquare$ 
\end{proof}

\begin{corollary}
Let $G$ be a group and $a,b\in G$. Then the map $B: G\to G$ defined as 
$B(g) = agb$, is a Rota--Baxter operator on $G$
if and only if $G$ is abelian and $b = a^{-1}$.
\end{corollary}

\begin{proof}
If $G$ is abelian, then every endomorphism of $G$ is a Rota--Baxter operator on $G$, including the identity map. Conversely, suppose that $B:G \to G$ defined as
$B(g) = agb$ is a~Rota--Baxter operator on $G$.
Since $B(1) = 1$, we obtain $b = a^{-1}$. Thus, $B$ is an (inner) automorphism of~$G$. It follows from Proposition~\ref{prop:RB-Hom}(2) that $G$ is abelian, and hence $b = a^{-1}$.
$\blacksquare$ 
\end{proof}

Given an integer~$k$, a group $G$ is called \index{$k$-abelian group}{$k$-abelian} if $(gh)^k = g^k h^k$ for all $g,h\in G$.

\begin{proposition}\label{prop:n-abelian}
Let $G$ be a group and $n \ge 1$ an integer. Then the map $B: G\to G$ defined by
$B(g) = g^n$ is a Rota--Baxter operator on $G$ if and only if $G$ is $(n+1)$-abelian.
\end{proposition}

\begin{proof}
Given $g,h\in G$, the Rota--Baxter defining identity~\eqref{RB} for $B$ has the form
$$
g^n h^n
 = (g^{n+1}hg^{-n})^n
 = g^{n+1}hg^{-n} \, g^{n+1}hg^{-n} \cdots g^{n+1}hg^{-n}
 = g^n(gh)^n g^{-n}.
$$
This gives
$(gh)^n = h^n g^n$, which is further equivalent to
$$
(hg)^{n+1} = h(gh)^n g = h^{n+1}g^{n+1}.
$$
Thus, $B$ is a Rota--Baxter operator on $G$ if and only if $G$ is $(n+1)$-abelian.
$\blacksquare$ 
\end{proof}

\begin{remark}{\rm 
Since every $k$-abelian group is also $(1-k)$-abelian, the group operation in $G_B$ is given by
$$
g\circ h
 = gg^nhg^{-n}
 = g^{n+1}h^{n+1}h^{-n}g^{-n}
 = (gh)^{n+1}(hg)^{-n}
$$
for all $g, h \in G$.}
\end{remark}

\begin{proposition}\cite[Proposition 24]{MR4556953}\label{prop:center-conj}
Let $G$ be a group and $g \in G$. Then the map $B_g(x) = g^{-1} x^{-1} g$ is a Rota--Baxter operator on $G$ if and only if $[g,G]\subseteq \Z(G)$.
\end{proposition}

\begin{proof}
Suppose that $B_g$ is a Rota--Baxter operator on $G$. This gives
$$
g^{-1} x^{-1}y^{-1} g
 = g^{-1}( x g^{-1}x^{-1}g y g^{-1} x g )^{-1}g
 = g^{-1}g^{-1}x^{-1}gy^{-1}g^{-1}xg x^{-1}g
$$
for all $x,y\in G$. The preceding identity is equivalent to
$$
x^{-1}y^{-1}
 = [g,x]x^{-1}y^{-1}[g,x^{-1}].
$$
Writing $s = x^{-1}y^{-1}$, we have
$[x,g]^s = [g,x^{-1}]$ for all $x,s\in G$.
Fixing $x\in G$, we conclude that $[g,x]\in \Z(G)$ for each $x\in G$.
\para 
Conversely, suppose that $g\in G$ such that $[g,x]\in \Z(G)$ for each $x\in G$. By the preceding argument, it is enough to prove that $[x,g] = [g,x^{-1}]$ for all $x\in G$.
Since
$$
[g,x^{-1}]
 = g^{-1}xgx^{-1}
 = x^{-1}g^{-1}xg[g^{-1}xg,x^{-1}]
 = [x,g]\, [x^g,x^{-1}],
$$
we need to prove that $[x^g,x^{-1}] = 1$. Writing $c = [g,x^{-1}]$, we see that
$$
[x^g,x^{-1}]
 = [[g,x^{-1}]x,x^{-1}]
 = [cx,x^{-1}]
 = x^{-1}c^{-1}xcxx^{-1}
 = 1, 
$$
which is desired. $\blacksquare$ 
\end{proof}

\begin{corollary}\label{coro:2step}
Let $G$ be a 2-step nilpotent group and $g \in G$. Then the map $B_g(x) = g^{-1} x^{-1} g$ is a Rota--Baxter operator.
\end{corollary}

\begin{remark}
{\rm 
Let $G$ be a group and $g \in G$ such that $[g,G]\subset \Z(G)$. If $B_g$ is the Rota--Baxter operator as defined in Proposition \ref{prop:center-conj}, then it is evident that $B_g$ is a bijection. Further, computing the group operation in $G_B$, we see that
$$
x\circ y
 = xB(x)yB(x)^{-1}
 = xg^{-1} x^{-1} gyg^{-1} xg
 = yxg^{-1} x^{-1} gg^{-1} xg
 = yx.
$$
Thus, the group structure induced by $B_g$ coincides with the group structure induced by the Rota--Baxter operator $B_{-1}$.}
\end{remark}

Let $F_n = \langle x_1, \ldots, x_n  \rangle$ be the free group of rank $n$ and $\pi_1 : F_n \to \langle x_1 \rangle$ be the projection given by
$$
\pi_1(x_1) = x_1 \quad \textrm{and} \quad \pi_1(x_i) = 1\quad \textrm{for} \quad i\neq1.
$$
Then we can write $F_n = \ker(\pi_1) \langle x_1 \rangle$. Further, for any integer $k$, the map $B: F_n \to F_n$ given by $B(g) = \pi_1(g)^k$ is a Rota--Baxter operator on $F_n$. The following result generalizes this observation.

\begin{proposition}\label{fr}
Let $G = K * L$ be a free product of two groups, where $L$ is abelian, and let $\pi_1: G \to L$ be the projection onto $L$. Then for any endomorphism $\chi$ of $L$, the map $B:G \to G$
 given by $B(g) = \chi(\pi_1(g))$ is a~Rota--Baxter operator on $G$.
\end{proposition}

We conclude this subsection with the following question.

\begin{question}
Let $G$ be a group and $H$ its proper subgroup. Under what conditions does there exist a Rota--Baxter operator $B$ on $G$ such that $B(G) = H$? In particular, if $F$ is a non-abelian free group, does there exist a Rota--Baxter operator on $F$ whose image is the commutator subgroup of $F$?
\end{question}
\bigskip
\bigskip


\subsection{Constructions via direct products}

Rota--Baxter operators on direct products of groups can be constructed from Rota--Baxter operators on the component groups \cite[Proposition 30]{MR4556953}.

\begin{proposition} \label{direct}
Let $G = H\times L$ be a direct product of groups $H$ and $L$. Then the following assertions hold: 
\begin{enumerate}
\item Let $B_H$ be a Rota--Baxter operator on $H$ and $B_L$ a Rota--Baxter operator on $L$.
Then the map $B: G\to G$ defined by
$B(hl) = B_H(h)B_L(l)$, where $h\in H$ and $l\in L$, is a Rota--Baxter operator.
\item If $L$ is abelian and $|L|>2$, then
there exists a non-splitting Rota--Baxter operator on $G$.

\end{enumerate}
\end{proposition}

\begin{proof}
The proof of assertion (1) is straightforward. For assertion (2), suppose that every Rota--Baxter operator on $G$ is splitting. Let $\psi$ be an automorphism of $L$. Define $B_\psi : G \to G$ by $B_\psi \big((h,l)\big) = \psi(l)$. Then, by Proposition~\ref{prop:RB-Hom}(1), $B_\psi$ is a~Rota--Baxter operator on~$G$. Further, by Proposition~\ref{SplittingCond}, we should have
$B_\psi \big(lB_\psi(l) \big) = \psi(l)\psi^2(l) = 1$ for all $l\in L$. Taking $\psi = \id$, we have $l^2 = 1$ for all $l\in L$. Since $L$ is abelian, we see that $L$ is a~direct sum of some copies of $\mathbb{Z}_2$. The condition $|L|>2$ implies that we can find two distinct non-trivial elements $l_1,l_2\in L$. Consider the automorphism~$\psi$ of~$L$ which interchanges $l_1$ and $l_2$ and fixes all other elements. Then, it is easy to see that $B_\psi$ is non-splitting, which is a~contradiction.
$\blacksquare$ 
\end{proof}

\begin{example} \label{exm:direct}
{\rm
Let $G$ be a group and $G^n = G\times \cdots \times G$. Then the maps 
$$ B \big((g_1,\ldots,g_n) \big) = \big(1,g_1,g_2g_1,g_3g_2g_1,\ldots,g_{n-1}g_{n-2}\ldots g_1 \big)$$
and
$$\widetilde{B} \big((g_1,\ldots,g_n) \big) = \big(g_1^{-1},g_2^{-1}g_1^{-1},g_3^{-1}g_2^{-1}g_1^{-1},\ldots, g_n^{-1}g_{n-1}^{-1}g_{n-2}^{-1}\ldots g_1^{-1} \big)$$
are Rota--Baxter operators on $G^n$.}
\end{example}

Next, consider the following example from 
 \cite{MR2450698,MR0448155}. Let $\Bbbk$ be a field and $\Bbbk^n = \Bbbk e_1 \oplus \cdots \oplus \Bbbk e_n$, where $e_ie_j = \delta_{ij}e_i$, viewed as an associative $\Bbbk$-algebra. A~linear operator $R:\Bbbk \to \Bbbk$ given by 
$R(e_i) = \sum\limits_{k=1}^n r_{ik}e_k$, where 
$r_{ik}\in \Bbbk$, is a Rota--Baxter operator of weight~1 on~$\Bbbk^n$
if and only if the following conditions are satisfied:
\begin{enumerate}
\item $r_{ii} = 0$ and $r_{ik}\in\{0,1\}$
or $r_{ii} = -1$ and $r_{ik}\in\{0,-1\}$ for all $k\neq i$.
\item If $r_{ik} = r_{ki} = 0$ for $i\neq k$,
then $r_{il}r_{kl} = 0$ for all $l\not\in\{i,k\}$.
\item If $r_{ik}\neq0$ for $i\neq k$,
then $r_{ki} = 0$ and
$r_{kl} = 0$ or $r_{il} = r_{ik}$ for all $l\not\in\{i,k\}$.
\end{enumerate}
We can assume that the matrix of $R$ is upper-triangular~\cite{MR0448155,MR4120095}. We use the preceding example to construct a Rota--Baxter operator on $G^n$ for any group $G$ \cite[Theorem 32]{MR4556953}.

\begin{theorem} \label{theo:directProduct}
Let $G$ be a group, $\Bbbk$ a field and $n \ge 1$ an integer. Let $R$ be a~Rota--Baxter operator of weight~1 on $\Bbbk^n$, given by $R(e_i) = \sum\limits_{i=1}^n r_{ik}e_k$, such that the matrix of $R$ is upper-triangular. Define a~map $B: G^n\to G^n$ by
\begin{equation} \label{RBOnDirectProduct}
B \big((g_1,\ldots,g_n) \big) = (t_1,\ldots,t_n),~~ \textrm{where}~~
t_i = g_i^{r_{ii}}g_{i-1}^{r_{i-1\,i}}\ldots g_1^{r_{1i}}.
\end{equation}
Then $B$ is a Rota--Baxter operator on $G^n$.
\end{theorem}

\begin{proof}
We write  the left and the right hand sides of the Rota--Baxter identity \eqref{RB} for~$B$. Let 
$g = (g_1,\ldots,g_n)$ and $h = (h_1,\ldots,h_n)$ be elements of $G^n$. Then we have
\begin{small}
\begin{eqnarray}
B(g)B(h)
& = &\big(g_1^{r_{11}}h_1^{r_{11}},g_2^{r_{22}}g_1^{r_{12}}h_2^{r_{22}}h_1^{r_{12}},\ldots,
 g_n^{r_{nn}}g_{i-1}^{r_{n-1\,n}}\ldots g_1^{r_{1n}}
 h_n^{r_{nn}}h_{i-1}^{r_{n-1\,n}}\ldots h_1^{r_{1n}}\big), \label{LeftRBSumFields} \\
B\big(gB(g)hB(g)^{-1}\big)
 &= & \big(t_1^{r_{11}},t_2^{r_{22}}t_1^{r_{12}},\ldots,
 t_n^{r_{nn}}t_{n-1}^{r_{n-1\,n}}\ldots t_1^{r_{1n}}\big), \label{RightRBSumFields}
\end{eqnarray}
\end{small}
where
\begin{eqnarray*}
t_1 &=& g_1g_1^{r_{11}}h_1g_1^{-r_{11}},\\
t_2 &=& g_2g_2^{r_{22}}g_1^{r_{12}}h_2g_1^{-r_{12}}g_2^{r_{22}},\\
&& \vdots \\
t_n  &=& g_n g_n^{r_{nn}}g_{n-1}^{r_{n-1\,n}}\ldots g_1^{r_{1n}}h_n g_1^{-r_{1n}}\ldots
 g_{n-1}^{-r_{n-1\,n}}g_n^{-r_{nn}}.
\end{eqnarray*}

We prove by induction on $i$ that the $i$-th coordinates of~\eqref{LeftRBSumFields} and~\eqref{RightRBSumFields} are equal.
For $i = 1$, we desire that
$$
g_1^{r_{11}}h_1^{r_{11}}
 = \big(g_1g_1^{r_{11}}h_1g_1^{-r_{11}}\big)^{r_{11}}.
$$
When $r_{11} = 0$, the equality becomes $1=1$,
and when $r_{11} = -1$, then we get $g_1^{-1}h_1^{-1} = (h_1g_1)^{-1}$. This proves the base step of the induction. Suppose that the $i$-coordinates
of~\eqref{LeftRBSumFields} and~\eqref{RightRBSumFields} are equal
for all $i<k$, where $2\leq k\leq n$. We prove the equality
$$
g_i^{r_{ii}}g_{i-1}^{r_{i-1\,i}}\ldots g_1^{r_{1i}}
h_i^{r_{ii}}h_{i-1}^{r_{i-1\,i}}\ldots h_1^{r_{1i}}
 = t_i^{r_{ii}}t_{i-1}^{r_{i-1\,i}}\ldots t_1^{r_{1i}}
$$
for $i = k$. If $r_{ii} = 0$, then we apply the induction hypothesis for $i-1$ and
for the coefficients $r'_{s\,i-1} = r_{si}$.
If $r_{ii} = -1$, then we apply the induction hypothesis for $i-1$ and
for the coefficients $r'_{s\,i-1} = r_{si}$ to derive
\begin{eqnarray*}
t_i^{r_{ii}}t_{i-1}^{r_{i-1\,i}}\ldots t_1^{r_{1i}}
 &=& (\underline{g_i g_i^{-1}}g_{i-1}^{r_{i-1\,i}}\ldots g_1^{r_{1i}}h_i g_1^{-r_{1i}}
 \ldots g_{i-1}^{-r_{i-1\,i}}g_i)^{-1}
 t_{i-1}^{r_{i-1\,i}}\ldots t_1^{r_{1i}} \\
 &=& g_i^{-1}g_{i-1}^{r_{i-1\,i}}\ldots g_1^{r_{1i}}h_i^{-1}\underline{g_1^{-r_{1i}}
 \ldots g_{i-1}^{-r_{i-1\,i}}
 g_{i-1}^{r_{i-1\,i}}\ldots g_1^{r_{1i}}}h_{i-1}^{r_{i-1\,i}}\ldots h_1^{r_{1i}} \\
 &=& g_i^{r_{ii}}g_{i-1}^{r_{i-1\,i}}\ldots g_1^{r_{1i}}
h_i^{r_{ii}}h_{i-1}^{r_{i-1\,i}}\ldots h_1^{r_{1i}}. 
\end{eqnarray*}
This proves the theorem.
$\blacksquare$ 
\end{proof}

It is not difficult to see that $G^n_B\cong G^n$ in Theorem \ref{theo:directProduct}.

\begin{corollary}
Let $G$ be a group and $\psi_2,\ldots,\psi_n\in\Aut(G)$. Then the map $P: G^n\to G^n$
given by 
$$ P \big((g_1,\ldots,g_n) \big) = (t_1,\ldots,t_n),$$
where 
$$t_1 = g_1^{r_{11}}\quad \textrm{and} \quad t_i = g_i^{r_{ii}}\psi_i\big(g_{i-1}^{r_{i-1\,i}}\psi_{i-1}\big(g_{i-2}^{r_{i-2\,i}}
 \ldots \psi_2\big(g_1^{r_{1i}}\big)\big)\big)$$
 for $i\geq2$, is a Rota--Baxter operator on $G^n$.
\end{corollary}

\begin{proof}
Define an automorphism $\varphi$ of $G^n$ by
$$
\varphi \big((g_1,\ldots,g_n)\big)
 = \big(g_1,\psi_2^{-1}(g_2),\psi_2^{-1}\psi_3^{-1}(g_3),\ldots,\psi_2^{-1}\ldots \psi_n^{-1}(g_n)\big).
$$
It follows from Lemma~\ref{lem:Aut}(1) that $B^{(\varphi)}$ is a Rota--Baxter operator on $G^n$, where $B$ is defined by \eqref{RBOnDirectProduct}. Since $P = B^{(\varphi)}$, the assertion follows.
$\blacksquare$ 
\end{proof}
\bigskip
\bigskip


\section{Extending maps to Rota--Baxter operators on groups}

Let $(G, B)$ be a Rota--Baxter group.  A natural problem is to identify a minimal subset $A$ of $G$ such that the values of $B$ on $A$ uniquely determine the values of $B$ on the entirety of $G$. Given a subset $A$ of $G$, let $\langle A \rangle_{B}$ denote the subgroup of $G_B$ generated by $A$. The following result addresses the preceding problem \cite[Theorem 36]{MR4556953}.

\begin{theorem}  \label{t}
Let $G$ be a group, $A$ a subset of $G$ and $B$ a Rota--Baxter operator on $G$. Then the following assertions hold:
\begin{enumerate}
\item The values of $B$ on elements of $A$ determine the values of $B$ on elements of $\langle A \rangle_{B}$.  In particular,  $B(\langle A \rangle_{B})$ is generated by $B(A)$.
\item  If $A$ is a generating subset of $G_B$, then the values of $B$ on elements of $A$ uniquely determine the values of~$B$ on the entirety of $G$.
\end{enumerate}
\end{theorem}

\begin{proof}
Let $A = \{a_i \,\mid \, i \in I\}$ and  $B(a_i) = u_i$ for each $i\in I$. Suppose that $g \in \langle A \rangle_{B}$, which is written as
$$
g = a_{i_1}^{\circ(k_1)} \circ a_{i_2}^{\circ(k_2)} \circ \cdots \circ a_{i_s}^{\circ(k_s)},
$$
where $a_{i_j}\in A$ and $k_{i_j} \in \mathbb{Z}$. Since $B$ is a group homomorphism from $G_B$ to $G$, we have
$$
B(g) = B\big(a_{i_1}^{\circ(k_1)}\big) B\big(a_{i_2}^{\circ(k_2)}\big)
 \cdots B\big(a_{i_s}^{\circ(k_s)}\big) = u_{i_1}^{k_1} u_{i_2}^{k_2} \ldots u_{i_s}^{k_s}.
$$
Hence, we can determine the values of $B$ on the subgroup $\langle A\rangle_B $, which proves assertion (1). In particular, if $\langle A \rangle_{B}$ equals $G_B$, then we can also determine values of $B$ on elements of $G$, which establishes assertion (2).
$\blacksquare$ 
\end{proof}

Next, we consider the general situation. Let $G$ be a~group generated by a~set $A = \{ a_i \mid i \in I \}$. Let $U = \{u_i \, \mid \, i \in I \}$ be a subset of~$G$ with the same cardinality as that of $A$ and $\beta : A \to U$ the map defined by $\beta(a_i) = u_i$ for each $i \in I$. Consider a set $\bar{A} = \{\bar{a_i} \,\mid \, i \in I \}$ which is disjoint from $A$ and has the same cardinality as that of $A$.
Let $F(\bar{A})$ be the free group with basis $\bar{A}$, and denote the group operation in $F(\bar{A})$ by~$\circ$. An element $w$ of $F(\bar{A})$ can be uniquely presented by a~reduced word as
$$
w = \bar{a}_{i_1}^{\circ(k_1)} \circ \bar{a}_{i_2}^{\circ(k_2)}
 \circ \cdots \circ \bar{a}_{i_s}^{\circ(k_s)},
$$
where $\bar{a}_{i_j} \in \bar{A}$ and $k_{i_j} \in \mathbb{Z} \setminus \{0\}$.
Define a~map $\bar{\beta} : F(\bar{A}) \to G$ by
$$
\bar{\beta}(w) = \beta(a_{i_1})^{k_1} \beta(a_{i_2})^{k_2} \cdots \beta(a_{i_s})^{k_s}
 = u_{i_1}^{k_1}  u_{i_2}^{k_2}  \cdots u_{i_s}^{k_s}.
$$
Then $\bar{\beta}$ is a group homomorphism with image $\langle U \rangle$. Define a map $\pi : F(\bar{A}) \to G$ by 
$$
\pi(w) = \big(a_{i_1}\beta(a_{i_1})\big)^{k_1}
    \big(a_{i_2}\beta(a_{i_2}) \big)^{k_2} \cdots
    \big(a_{i_s}\beta(a_{i_s}) \big)^{k_s}
    \beta(a_{i_s})^{-k_s}
    \beta(a_{i_{s-1}})^{-k_{s-1}} \cdots
    \beta(a_{i_1})^{-k_1}.
    $$
The following result clarifies the connection between the maps $\pi$ and $\bar{\beta}$.

\begin{lemma} \label{l7.3}
Let $w, w' \in F(\bar{A})$. Then the following  assertions hold:
\begin{enumerate}
\item $\pi(w^{\circ({-1})}) = \bar{\beta}(w)^{-1} \pi(w)^{-1} \bar{\beta}(w)$.
\item $\pi(w \circ w' ) =  \pi(w) \bar{\beta}(w) \pi(w')   \bar{\beta}(w)^{-1}$.
\item $\pi \big((w')^{\circ(-1)} \circ w \circ w' \big) = \bar{\beta}(w')^{-1} \pi(w')^{-1} \pi(w)
 \bar{\beta}(w) \pi(w') \bar{\beta}(w)^{-1} \bar{\beta}(w')$.
\end{enumerate}
\end{lemma}

\begin{proof}
If $w = \bar{a}_{i_1}^{\circ(k_1)} \circ \bar{a}_{i_2}^{\circ(k_2)}
 \circ \cdots \circ \bar{a}_{i_s}^{\circ(k_s)},
$ then its inverse is
$$
w^{\circ({-1})} = \bar{a}_{i_s}^{\circ(-k_s)}
 \circ \bar{a}_{i_{s-1}}^{\circ(-{k_{s-1}})}
 \circ \cdots \circ \bar{a}_{i_1}^{\circ(-k_1)}
$$
and
$$
\pi(w^{\circ({-1})})
= \big(a_{i_s}\beta(a_{i_s})\big)^{-k_s}
    \big(a_{i_{s-1}}\beta(a_{i_{s-1}})\big)^{-k_{s-1}} \cdots
    \big(a_{i_1}\beta(a_{i_1})\big)^{-k_1}
    \beta(a_{i_1})^{k_1}
    \beta(a_{i_{2}})^{k_{2}} \cdots
    \beta(a_{i_s})^{k_s}.
$$
This gives
$$
\pi(w^{\circ({-1})})
 = \beta(a_{i_s})^{-k_s} \cdots \beta(a_{i_1})^{-k_1} \pi(w)^{-1}
 \beta(a_{i_1})^{k_1} \cdots \beta(a_{i_s})^{k_s}=\bar{\beta}(w)^{-1} \pi(w)^{-1} \bar{\beta}(w),
$$
which establishes assertion (1).
\para 
Given $w' = \bar{a}_{j_1}^{\circ(l_1)} \circ \bar{a}_{j_2}^{\circ(l_2)}
   \circ \cdots \circ \bar{a}_{j_t}^{\circ(l_t)},
$
we have
$$
w \circ w'
 = \bar{a}_{i_1}^{\circ(k_1)}  \circ \cdots \circ \bar{a}_{i_s}^{\circ(k_s)}
  \circ \bar{a}_{j_1}^{\circ(l_1)}  \circ \cdots \circ \bar{a}_{j_t}^{\circ(l_t)}
$$
and
\begin{eqnarray*}
\pi(w \circ w')
 &=& \big(a_{i_1}\beta(a_{i_1})\big)^{k_1} \cdots \big(a_{i_s}\beta(a_{i_s})\big)^{k_s} \,
  \big(a_{j_1}\beta(a_{j_1}) \big)^{l_1} \cdots  \big(a_{j_t}\beta(a_{j_t})\big)^{l_t} \\
 &&\beta(a_{j_t})^{-l_t} \cdots \beta(a_{j_1})^{-l_1}\, \beta(a_{i_s})^{-k_s} \cdots \beta(a_{i_1})^{-k_1}.
\end{eqnarray*}
Hence, we have
$$
\pi(w \circ w' )
 = \pi(w) \beta(a_{i_1})^{k_1} \cdots \beta(a_{i_s})^{k_s}
   \pi(w') \, \beta(a_{i_s})^{-k_s} \cdots \beta(a_{i_1})^{-k_1}
 = \pi(w) \bar{\beta}(w) \pi(w') \bar{\beta}(w)^{-1},
$$
which proves assertion (2).
\para 

Considering the product $(w')^{\circ(-1)} \circ w \circ w'$, we see that
\begin{eqnarray*}
&& \pi \big((w')^{\circ(-1)} \circ w \circ w' \big)\\
 &=& \beta(a_{j_l})^{-l_t} \cdots \beta(a_{j_1})^{-l_1}
 \big( \pi(w')^{-1} \pi(w) \beta(a_{i_1})^{k_1} \cdots \beta(a_{i_s})^{k_s} \pi(w')
  \beta(a_{i_s})^{-k_s} \cdots \beta(a_{i_1})^{-k_1} \big)
 \beta(a_{j_1})^{l_1} \cdots \beta(a_{j_t})^{l_t} \\
 &=& \bar{\beta}(w')^{-1} \pi(w')^{-1} \pi(w) \bar{\beta}(w)
 \pi(w') \bar{\beta}(w)^{-1} \bar{\beta}(w'),
\end{eqnarray*}
which is assertion (3).
$\blacksquare$ 
\end{proof}

Lemma \ref{l7.3} immediately yields the following result.

\begin{lemma}
Let $R = \{ w \in F(\bar{A}) \, \mid \, \pi(w) = 1 \}$.  
\begin{enumerate}
\item $R$ is a~subgroup of $F(\bar{A})$. 
\item Suppose that $\pi$ and $\bar{\beta}$ satisfy the following implication:
\begin{equation} \label{condition ab}
\pi(w) = \pi(w') \Rightarrow \bar{\beta}(w) = \bar{\beta}(w').
\end{equation}
Then $R$~is a normal subgroup of $F(\bar{A})$.
\end{enumerate}
\end{lemma}

\begin{lemma}\label{lem:RB-extCond}
Let $w, w' \in F(\bar{A})$ be such that $\pi(w) = \pi(w')$ and $\bar{\beta}(w) \neq \bar{\beta}(w')$. Then the map $\beta$~has no extension to a Rota--Baxter operator on $G$.
\end{lemma}

\begin{proof}
Suppose that $B$~is an extension of~$\beta$ to a Rota--Baxter operator~on~$G$. Let
$$
w = \bar{a}_{i_1}^{\circ(k_1)} \circ \bar{a}_{i_2}^{\circ(k_2)}
 \circ \cdots \circ \bar{a}_{i_s}^{\circ(k_s)}\quad \textrm{and}\quad
w' = \bar{a}_{j_1}^{\circ(l_1)} \circ \bar{a}_{j_2}^{\circ(l_2)}
 \circ \cdots \circ \bar{a}_{j_t}^{\circ(l_t)}.
$$
Consider the elements
$$
u = a_{i_1}^{\circ(k_1)} \circ a_{i_2}^{\circ(k_2)}
 \circ \cdots \circ a_{i_s}^{\circ(k_s)} \quad \textrm{and}\quad
u' = a_{j_1}^{\circ(l_1)} \circ a_{j_2}^{\circ(l_2)}
 \circ \cdots \circ a_{j_t}^{\circ(l_t)}
$$
of $G_B$. Then we have
$$
B(u) = \beta(a_{i_1})^{k_1}\cdots \beta(a_{i_s})^{k_s}
 = \bar{\beta}(w)
 \neq \bar{\beta}(w')
 = \beta(a_{j_1})^{l_1}\cdots \beta(a_{j_t})^{l_t}
 = B(u').
$$
On the other hand, in view of Proposition \ref{for}, we have
$$
B(u)
 = B \big(\pi(w)\big)
 = B \big(\pi(w')\big)
 = B(u'),
$$
which is a contradiction.
$\blacksquare$ 
\end{proof}

If $R$ is a~normal subgroup of $F(\bar{A})$,
then we can define the group $$\bar{G}_{\beta} = F(\bar{A}) / R.$$
Let $[w]$ denote the element of $\bar{G}_{\beta}$ represented by an element $w \in F(\bar{A})$. Then we write the product in $\bar{G}_{\beta}$ as $[w] \circ [w']  = [w \circ w']$. Note that $\bar{G}_{\beta}$ is generated by $\big\{ [\bar{a_i}] \, \mid \,  i\in I \big\}$ and has defining relations $[w] = [w']$ for all words
$$
w = \bar{a}_{i_1}^{\circ(k_1)} \circ \bar{a}_{i_2}^{\circ(k_2)}
 \circ \cdots \circ \bar{a}_{i_s}^{\circ(k_s)} \quad \textrm{and}\quad
w' = \bar{a}_{j_1}^{\circ(l_1)} \circ \bar{a}_{j_2}^{\circ(l_2)}
 \circ \cdots \circ \bar{a}_{j_t}^{\circ(l_t)}
$$
such that $\pi(w) = \pi(w')$. In other words, $[w] = [w']$ if and only if the relation
\begin{eqnarray*}
&& \big(a_{i_1}\beta(a_{i_1})\big)^{k_1}
    \big(a_{i_2}\beta(a_{i_2})\big)^{k_2} \ldots
    \big(a_{i_s}\beta(a_{i_s})\big)^{k_s}
    \beta(a_{i_s})^{-k_s}
    \beta(a_{i_{s-1}})^{-k_{s-1}} \ldots
    \beta(a_{i_1})^{-k_1} \\
& =&  \big(a_{j_1}\beta(a_{j_1})\big)^{l_1}
    \big(a_{j_2}\beta(a_{j_2})\big)^{l_2} \ldots
    \big(a_{j_t}\beta(a_{j_t})\big)^{l_t}
    \beta(a_{j_l})^{-l_t}
    \beta(a_{j_{t-1}})^{-l_{t-1}} \ldots
    \beta(a_{j_1})^{-l_1}
\end{eqnarray*}
holds in  $G$.
\para
If condition~\eqref{condition ab} is satisfied, then the homomorphism
$\bar{\beta} : F(\bar{A}) \to \langle U \rangle \leq G$
induces a homomorphism
$\bar{G}_{\beta} \to \langle U \rangle$, which we again denote by $\bar{\beta}$.
Also, the map $\pi : F(\bar{A}) \to G$ induces a map
$\bar{\pi}:\bar{G}_{\beta} \to G$. It is easy to see that $\im(\pi) = \im(\bar{\pi})$. Let us consider an illustration of Lemma~\ref{lem:RB-extCond}.

\begin{example}{\rm 
Consider the symmetric group
$$
\Sigma_3 = \big\langle s_1, s_2  \, \mid \,  s_1^2 = s_2^2 = 1~~\textrm{and} ~~s_1 s_2 s_1 = s_2 s_1 s_2 \big\rangle
$$
on three symbols. 
\begin{enumerate}
\item If we take the map $\beta : \{ s_1, s_2 \} \to \{ s_1, s_2 \}$
such that $\beta(s_1)= s_1$ and $\beta(s_2)= s_2$, then it has an extension to the Rota--Baxter operator
$B_{-1}(g) = g^{-1}$ on $\Sigma_3$. Let us determine the group $\bar{G}_{\beta}$.
Let $F(\bar{A})$ be the free group with basis $\bar{A} = \{ t_1, t_2 \}$. We have the map $\pi : F(\bar{A}) \to G$ defined by 
$$
\pi \big(t_1^{\circ(k_1)} \circ t_2^{\circ(l_1)} \circ
 \cdots \circ t_1^{\circ(k_s)} \circ t_2^{\circ(l_s)}\big)
 = s_2^{-l_s} s_1^{-k_s} \cdots s_2^{-l_1} s_1^{-k_1}.
$$
In this case, $R = \{ w \in F(\bar{A}) \, \mid \, \pi (w) = 1 \}$ is normal in $F(\bar{A})$ and the product in $\bar{G}_{\beta}$ is opposite to the product in~$G$.
\item  Note that $\Sigma_3$ is also generated by the elements $\tau_1 = s_1$,
$\tau_2 = s_2$ and $\tau_3 = s_1 s_2 s_1$. Suppose that
$$
\beta(\tau_1) = \tau_1,\quad
\beta(\tau_2) = \tau_2 \quad \textrm{and} \quad
\beta(\tau_3) = \tau_2 \tau_1.
$$
Let  $F(\bar{A})$ be the free group with basis $\bar{A} =\{t_1,t_2,t_3\}$.
Then we have
\begin{gather*}
\pi(t_1 \circ t_1) = \pi(t_2 \circ t_2) = \pi(t_3 \circ t_1) = 1, \quad
\pi(t_1 \circ t_2 \circ t_1) = \pi(t_2 \circ t_1 \circ t_2) = \tau_1 \tau_2 \tau_1, \\
\pi(t_1 \circ t_2) = \pi(t_2 \circ t_3) = \pi(t_3 \circ t_3) =  \tau_2 \tau_1 \quad \textrm{and} \quad
\pi(t_2 \circ t_1) = \pi(t_3 \circ t_2) = \pi(t_1 \circ t_3) = \tau_1 \tau_2.
\end{gather*}
Hence, we have $\pi(F(\bar{A}))  = \Sigma_3$. On the other hand, we have
$$
1 = \bar{\beta}(t_1 \circ t_1)  \not= \bar{\beta}(t_3 \circ t_1) = t_2.
$$
Thus,  by Lemma~\ref{lem:RB-extCond}, the map $\beta$ has no extension to a~Rota--Baxter operator on $\Sigma_3$.
\end{enumerate}}
\end{example}

We now formulate the main result of this subsection \cite[Theorem 41]{MR4556953}.

\begin{theorem} \label{thm:RBextension}
Let $G$ be a~group generated by a~set $A = \{ a_i \, \mid\, i \in I \}$. Let $U = \{u_i \,\mid \, i \in I\}$ be a subset of~$G$ with the same cardinality as that of $A$ and $\beta : A \to U$ the map given by $\beta(a_i) = u_i$ for each $i \in I$. If condition~\eqref{condition ab} is satisfied and the map $\bar{\pi} : \bar{G}_{\beta} \to G$ is bijective, then $\beta$~can be extended to a Rota--Baxter operator on~$G$ with $G_B\cong \bar{G}_{\beta}$.
\end{theorem}

\begin{proof}
Let $g \in G$ be an element. Since $\bar{\pi}$ is invertible, there exists a~unique element
$\bar{g} \in \bar{G}_{\beta}$ such that $\bar{\pi}(\bar{g}) = g$. Suppose that $\bar{g} = [w]$, where
$$
w = \bar{a}_{i_1}^{\circ(k_1)} \circ \bar{a}_{i_2}^{\circ(k_2)}
 \circ \cdots \circ \bar{a}_{i_s}^{\circ(k_s)}
$$
is a~word in the free group $F(\bar{A})$. We define $B:G \to G$ by
$$
B(g)
 = B \big(\bar{\pi}(\bar{g})\big)
 = B \big(\pi(w)\big)
 = \bar{\beta}(w)
 =  u_{i_1}^{k_1} u_{i_2}^{k_2} \cdots  u_{i_s}^{k_s}.
$$
Condition~\eqref{condition ab} implies that $B$ does not depend on the choice of~$w$.
\para 
Next we prove that~$B$ is a~Rota--Baxter operator on~$G$. For $g' \in G$, there exists
$$
w' = \bar{a}_{j_1}^{\circ(l_1)} \circ \bar{a}_{j_2}^{\circ(l_2)}
 \circ \cdots \circ \bar{a}_{j_t}^{\circ(l_t)}
$$
such that $\pi(w') = \bar{\pi}([w']) = g'$. Thus, we have $B(g') = u_{j_1}^{l_1} u_{j_2}^{l_2} \ldots  u_{j_t}^{l_t}$. On the other hand, by Lemma \ref{l7.3}(2), we get
$$
\pi(w \circ w') = \pi(w) \bar{\beta}(w) \pi(w') \bar{\beta}(w)^{-1},
$$
and consequently
\begin{eqnarray*}
B\big( gB(g)g'B(g)^{-1} \big)
 &=& B \big(\pi(w \circ w') \big)\\
 &=& \bar{\beta}(w \circ w')\\
 &=& \bar{\beta}(w) \bar{\beta}(w') \\
 &=& u_{i_1}^{k_1} u_{i_2}^{k_2} \cdots  u_{i_s}^{k_s}\cdot u_{j_1}^{l_1} u_{j_2}^{l_2} \ldots  u_{j_t}^{l_t}\\
 &=& B(g)B(g').
\end{eqnarray*}
Hence, $B$ is indeed a Rota--Baxter operator on $G$, which extends $\beta$.
$\blacksquare$ 
\end{proof}

The next result shows that condition~\eqref{condition ab} is not sufficient for extending a given map to a  Rota--Baxter operator.

\begin{proposition}
Let $\Sigma_3$ be the symmetric group on 3 symbols. Then the following  assertions hold:
\begin{enumerate}
\item The map $\beta$ given by $\beta(s_1)= s_1$ and  $\beta(s_2)= 1$ defines the group
$\bar{G}_{\beta}$, which is isomorphic to $\mathbb{Z}_2\times \mathbb{Z}_2$.
\item There is no~Rota--Baxter operator on $\Sigma_3$ extending the map $\beta$.
\end{enumerate}
\end{proposition}

\begin{proof}
Consider the free group $\langle t_1, t_2 \rangle$ of rank two. Then, we have
$$
\pi(t_1 \circ t_1) = \pi(t_2\circ t_2) = 1 \quad \textrm{and}\quad \pi(t_1 \circ t_2) = \pi(t_2\circ t_1) = s_2 s_1.
$$
It follows from the equalities
$$
\bar{\beta}(t_1 \circ t_1) = \bar{\beta}(t_2\circ t_2) = 1 \quad \textrm{and}\quad
\bar{\beta}(t_1 \circ t_2) = \bar{\beta}(t_2\circ t_1) =  s_1
$$
that \eqref{condition ab} holds, and 
$$
\bar{G}_{\beta}  = \big\langle [t_1], [t_2]  \, \mid \,  [t_1] \circ [t_1] = [t_2] \circ [t_2] = 1~\textrm{and}~[t_1] \circ [t_2] = [t_2]\circ [t_1] \big\rangle \cong \mathbb{Z}_2 \times \mathbb{Z}_2.
$$
\para 
Suppose that there is a Rota--Baxter operator $B$ extending the map $\beta$. Since $B(s_2) = 1$, by Lemma \ref{kernelCosets}, we have
$$
B(s_1 s_2) = B(s_2s_1s_2) = B(s_1 s_2 s_1).
$$
Noting that
$$
s_1 B(s_1s_2)
 = B(s_1)B(s_1s_2)
 = B(s_1 s_1 s_1s_2 s_1)
 = B(s_1 s_2 s_1)
 = B(s_1s_2), 
$$
we arrive at a contradiction. $\blacksquare$ 
\end{proof}

The following result gives some elementary Rota--Baxter operators on free groups \cite[Proposition 43]{MR4556953}. 

\begin{proposition} 
Let $F(A)$ be a free group with a basis $A = \{a_i \, \mid \, i \in I\}$.  Then the following  assertions hold:
\begin{enumerate}
\item The map $\beta$, given by $\beta(a_i) = 1$ for all $i \in I$, has the unique extension to the trivial Rota--Baxter operator $B_0$ on $F(A)$. 
\item The map $\beta$, given by $\beta(a_i)=  a_i^{-1}$ for all $i \in I$, has the unique extension to the Rota--Baxter operator $B_{-1}$ on $F(A)$. 
\end{enumerate}
\end{proposition}

\begin{proof}
Let $B$ be an extension of $\beta$ to a Rota--Baxter operator on $F(A)$, where $\beta(a_i)=1$ for all $i \in I$. We use $\circ$ to denote the product  given by \eqref{R-product}. In this case, $a_i^{\circ (k)} = a_i^k$, and hence $B(a_i^{\circ (k)}) = 1$
for all $i \in I$ and $k \in \mathbb{Z}$. If
$$
w = a_{i_1}^{\circ(k_1)} \circ a_{i_2}^{\circ(k_2)} \circ \cdots
 \circ a_{i_s}^{\circ(k_s)},
$$
where $a_{i_j} \in A$ and $k_{i_j} \in \mathbb{Z}$, is an arbitrary element, then we have
$$
B(w)
 = B\big(a_{i_1}^{\circ(k_1)}\big)
   B\big(a_{i_2}^{\circ(k_2)}\big)
   \cdots B\big(a_{i_s}^{\circ(k_s)}\big)
 = 1,
$$
which proves assertion (1). The proof of assertion (2) is similar. $\blacksquare$ 
\end{proof}

The following problem seems natural at this point.

\begin{problem}
Determine all Rota--Baxter operators on non-abelian free groups.
\end{problem}
\para

The following example shows that the values of~$\beta$ on the generators of $F_n$ do not determine the values of~$\beta$ on all elements of $F_n$.

\begin{example} \label{e7.8}
{\rm 
Let $F_2 = \langle a, b \rangle$ be the free group of rank two. We take $\beta(a) = a$ and $\beta(b) = 1$.
By Proposition \ref{fr}, there exists a Rota--Baxter operator on $F_2$ which is extension of the map~$\beta$. In fact, the Rota--Baxter operator $B$ is the endomorphism of $F_2$ given by $B(a) = \beta(a) = a$ and $B(b) = \beta(b) = 1$. As usual, we write the product in $F_2$ as the juxtaposition and the product in $(F_2)_B$ as $\circ$. We attempt do define the binary operation~$\circ$ on~$F_2$ by applying only the map~$\beta$.
We have
\begin{eqnarray*}
a^{\circ(k)} &=& a^k,\\
b^{\circ(l)} &=& b^l,\\
\beta(a^{\circ(k)}) &=& \beta(a)^k = a^k,\\
\beta(b^{\circ(l)}) &=&  \beta(b)^l = 1
\end{eqnarray*}
for all $k, l \in \mathbb{Z}$. For words of syllable length two, we have
$$ a^{\circ(k)} \circ b^{\circ(l)} = a^{2k} b^l a^{-k}\quad \textrm{and} \quad b^{\circ(l)} \circ a^{\circ(k)} =  b^l a^{k}$$
for all $k,l \in \mathbb{Z}$. By Proposition~\ref{for}, for a word of arbitrary length, we have
$$
w = a^{\circ(k_1)} \circ b^{\circ(l_1)} \circ a^{\circ(k_2)} \circ b^{\circ(l_2)}
  \circ \cdots \circ a^{\circ(k_s)} \circ b^{\circ(l_s)}
  = a^{2k_1} b^{l_1} a^{2k_2} b^{l_2}\cdots a^{2k_s} b^{l_s} a^{-\sum k_i}.
$$
Further, we see that
$$
\beta(w) = a^{\sum k_i}.
$$
Hence, we can construct an extension of~$\beta$ to the subgroup~$S = \langle a,b  \rangle_B$ of $(F_2)_B$ generated by $a$ and~$b$. Thus, $\ker(\beta) = \{ w  \, \mid \,  \sum k_i = 0 \}$
is the normal closure of~$b$ in $S$. Hence, $S \neq F_2$ and $ab \not\in S$. Let us add the element $ab$ to the generating set $\{a, b\}$ and put
$$
\beta_1(a) = \beta_1(ab) = a \quad \textrm{and} \quad \beta_1(b) =  1.
$$
Then, by Proposition~\ref{for}, for each $m \in \mathbb{Z}$, we have
$$
(ab)^{\circ (m)} = (aba)^m a^{-m} 
$$
in the group $T = \langle a, b, ab \rangle_B$. In particular, $ab, a^{-1} b^{-1}\in T$, and $T$ also consists of elements
\begin{eqnarray*}
 && a^{\circ(k_1)} \circ b^{\circ(l_1)} \circ (ab)^{\circ(m_1)}
 \circ a^{\circ(k_2)} \circ b^{\circ(l_2)} \circ  (ab)^{\circ(m_2)} \circ \cdots
  \circ a^{\circ(k_s)} \circ b^{\circ(l_s)} \circ (ab)^{\circ(m_s)} \\
  &=& a^{2k_1}  b^{l_1} (aba)^{m_1} a^{2k_2} b^{l_2} (aba)^{m_2}
  \cdots  a^{2k_s}  b^{l_s} (aba)^{m_s} a^{-\sum (k_i+m_1)}
\end{eqnarray*}
for integers $k_i, l_i, m_i$. It is not difficult to see that $T\neq F_2$.}
\end{example}

\bigskip
\bigskip


\section{Extensions of Rota--Baxter groups}
The set of all Rota--Baxter groups can be turned into a category $\mathcal{RBG}$ with the following definition of a morphism in this category.

\begin{definition}\label{definition of hom of RBG}
Let $(G_1, B_1)$ and $(G_2, B_2)$ be Rota--Baxter groups. A \index{morphism of Rota--Baxter groups}{morphism of Rota--Baxter groups} 
 $\psi:(G_1, B_1) \to (G_2, B_2)$ is a group homomorphism $\psi:G_1 \to G_2$ such that $B_2\, \psi= \psi \,B_1$.
\para 
A morphism between Rota--Baxter groups is said to be an isomorphism or an automorphism of Rota--Baxter groups if it is an isomorphism or an automorphism of the respective groups.
\end{definition}

\begin{definition}
Let $(G, B)$ be a Rota--Baxter group. A \index{Rota--Baxter subgroup}{Rota--Baxter subgroup} of $(G, B)$ is a subgroup $H$ of $G$  such that $B(H) \subseteq H$. Moreover, if $H$ is normal in $G$, then we refer to $H$ as a Rota--Baxter normal subgroup of $(G, B)$.
\end{definition}
\para 

It follows that if $H$ is a Rota--Baxter subgroup of $(G, B)$, then $(H, B|_H)$  is itself a Rota--Baxter group, and the inclusion map $H \hookrightarrow G$ is a morphism of Rota--Baxter groups. Further, if $H$ is a Rota--Baxter normal subgroup of $(G, B)$, then the map 
$\overline{R}:G/H \to G/H$, defined by  $\overline{R}(\overline{g})= \overline{R(g)}$, is a Rota--Baxter operator on $G/H$, where 
$\overline{g} \in G/H$ denote the image of  $g \in G$. It is not difficult to see that the kernel and the image of a morphism between Rota--Baxter groups is a Rota--Baxter normal subgroup and a Rota--Baxter subgroup of their respective groups. This allows us to define an extension  of Rota--Baxter groups.

\begin{definition}
An \index{extension of Rota--Baxter groups}{extension of Rota--Baxter groups} is a sequence 
\begin{equation} \label{ses}
1 \to (H, B_H) \stackrel{i}{\to}(G, B_G) \stackrel{\pi}{\to} (L, B_L) \to 1
\end{equation}
of Rota--Baxter groups and morphisms such that $i$ is injective, $\pi$ is surjective and 
$\im(i)= \ker(\pi)$.
\end{definition}

Analogous to the groups, the following problem can be formulated.

\begin{problem}\label{Rota--Baxter extension problem}
Let $(H, B_H)$ and $(L, H_L)$ be Rota--Baxter groups. Determine all Rota--Baxter groups $(G, B_G)$  which fits into the exact sequence~\eqref{ses}.
\end{problem}

\begin{remark}
{\rm
Proposition~\ref{direct} shows that Rota--Baxter groups as desired in Problem \ref{Rota--Baxter extension problem} do exist. In fact, \cite[Theorem 4.4]{MR4644858} provides an answer for abelian extensions of Rota--Baxter groups in terms of some appropriately defined second cohomology group.}
\end{remark}
\para 

Next, we recall the construction of wreath product of groups. Let $H$ and $L$ be groups. Let $\Map(L, H)$ and $\map(L, H)$ be the Cartesian product and the direct sums of copies $H$ indexed by elements of $L$, respectively. In other words, $\Map(L, H)$ is the group of all maps $L \to H$
and $\map(L, H)$ is its subgroup consisting of maps with a finite support. For $f \in \Map(L, H)$ and $l \in L$, consider the map $f^l$ defined by 
$f^l(x) = f(l x)$ for all $x \in L$. Then the map
$$
\hat{l} : \Map(L, H) \to \Map(L, H),
$$
given by $f \mapsto f^l$, is an automorphism of $\Map(L, H)$, which sends $\map(L, H)$ to itself. Further, the maps
$$
L \to \Aut \big( \Map(L, H) \big)~~\textrm{and}~~ L \to \Aut \big( \map(L, H) \big)
$$
given by $l \mapsto \hat{l}$ are embeddings of groups. The Cartesian wreath product $H \bar{\wr} L$ is defined as the semi-direct product $L  \ltimes \Map(L, H)$, where the group  operation is given by
$$
(l f) \, (l' f') = (l \, l') (f^{l'} f').
$$
Similarly, the direct wreath product is defined as the group $H \wr L = L \ltimes \map(L, H)$.
\para

The following result gives some Rota--Baxter operators on Cartesian and  direct wreath products of groups.

\begin{proposition}
Let $G$ be a Cartesian or a direct wreath product of groups $H$ and $L$. Then the following  assertions hold:
\begin{enumerate}
\item The map $B(lf) = f^{-1}$ is a Rota--Baxter operator on $G$.
\item If $L$ is abelian and $\varphi:L \to L$ is a group homomorphism,
then the map $B(l f) = \varphi(l)$ is a Rota--Baxter operator on $G$.
\item If the action of $L$ on $\Map(L, H)$
is trivial, and $B_L$ and $B_H$ are Rota--Baxter operators on $L$ and on $\Map(L, H)$, then the map $B(l f) = B_L(l) B_H(f)$
is a Rota--Baxter operator on~$G$.
\item  If the action of $L$ on $\map(L, H)$ 
is trivial, and $B_L$ and $B_H$ are Rota--Baxter operators on $L$ and on $\map(L, H)$, then the map $B(l f) = B_L(l) B_H(f)$ is a Rota--Baxter operator on~$G$.
\end{enumerate}
\end{proposition}

\begin{proof}
Since $G$ has an exact factorization $G = L \ltimes \Map(L, H)$ or $G = L \ltimes \map(L, H)$,  assertion (1) follows from Example~\ref{exm:split}. Assertion (2) follows from Proposition~\ref{prop:RB-Hom}(1). Finally, assertions (3) and (4) follow from Proposition~\ref{direct}.
$\blacksquare$ 
\end{proof}

We conclude with the following questions concerning extensions of Rota--Baxter groups.

\begin{question}
Let  $1 \to H \to G \to L \to 1$ be a short exact sequence of groups.
\begin{enumerate}
\item Let $(H, B_H)$ be a Rota--Baxter group. Under what conditions does there exists a Rota--Baxter group $(G, B_G)$
such that $(H, B_H)$ is its Rota--Baxter-subgroup?
\item Let $(L, B_L)$ be a Rota--Baxter group. Under what conditions does there exist a Rota--Baxter group $(G, B_G)$ such that $B_L$ is induced by $B_G$?
\item It is well-known that any extension of a group~$H$ by a group~$L$ can be embedded into the Cartesian wreath product $H \bar{\wr} L$,
known as the \index{Frobenius embedding}{\it Frobenius embedding}.  Does every Rota--Baxter-extension of $(H, B_H)$ by $(L, B_L)$  embed into $(H \bar{\wr} L, B)$ for some Rota--Baxter operator $B$ on $H \bar{\wr} L$?
\end{enumerate}
\end{question}

\begin{remark}
We expect that the obstruction to the extension and the lifting problem lie in some appropriately defined cohomology of Rota--Baxter groups.  The reader may refer to \cite{MR4644858} for some related work. The answer is well-known for the case of group automorphisms \cite{MR3887655, MR0272898}.
\end{remark}
\bigskip
\bigskip


\section{Rota--Baxter operators on groups and associated Lie rings}

In \cite{MR4271483}, motivated by Rota--Baxter operators on Lie algebras arising from Lie groups,  Guo et al. introduced Rota--Baxter operators on (abstract) groups. We see that this connection can be extended to the \index{Lie ring}{Lie ring} of any group.

\begin{definition}
Let $G$ be a group and $\{\gamma_n(G) \}_{n \ge 1}$ be its lower central series, where $\gamma_1(G) = G$ and $\gamma_{n+1}(G) = [G,\gamma_n(G)]$ for each $n \ge 1$. The Lie ring $L(G)$ of $G$
is the associated graded abelian group
$$L(G) = \bigoplus_{n\geq1} \gamma_n(G)/ \gamma_{n+1}(G)$$
equipped with the Lie bracket
$[x \gamma_{i+1}(G), \, y \gamma_{j+1}(G)] = [x,y] \gamma_{i+j+1}(G)$, where $[x,y]$ is the usual commutator in $G$. 
\end{definition}

We have the following result  \cite[Theorem 46]{MR4556953}. 

\begin{proposition} \label{prop:LieRing}
Let $B$ be a Rota--Baxter operator on a group $G$ such that $B\big(\gamma_n(G)\big)\subseteq \gamma_n(G)$ for all $n\geq1$.
Then the map $R:L(G) \to L(G)$ defined by $R \big(x \gamma_{i+1}(G)\big) = B(x) \gamma_{i+1}(G)$ is a Rota--Baxter operator
of weight~1 on the Lie ring $L(G)$.
\end{proposition}

\begin{proof}
Given $h\in \gamma_i(G)$ and $g\in \gamma_{i+1}(G)$, it follows from Lemma~\ref{lem:elementary}(4) that
$B(h)^{-1}B(hg)\in \gamma_{i+1}(G)$, and hence $R$ is well-defined. For $x\in \gamma_i(G)$ and $y\in \gamma_j(G)$, we need to prove that
$$
\big[R(x \gamma_{i+1}(G)),R(y \gamma_{j+1}(G)) \big]
 = R\big( [R(x \gamma_{i+1}(G)),y \gamma_{j+1}(G)] + [x \gamma_{i+1}(G), R(y \gamma_{j+1}(G))] + [x \gamma_{i+1}(G),y \gamma_{j+1}(G)] \big),
$$
which in terms of $B$ is equivalent to 
\begin{equation}\label{LieRing}
B(x)^{-1}B(y)^{-1}B(x)B(y) \gamma_{i+j+1}(G)
 = B\big( [B(x),y] + [x,B(y)] + [x,y] \big)\gamma_{i+j+1}(G).
\end{equation}
Applying the commutator properties and the fact that each $\gamma_{i+j}(G)/\gamma_{i+j+1}(G)$ is abelian,
we rewrite the expression inside the brackets on the right-hand side of~\eqref{LieRing}  as
\begin{eqnarray*}
&& [B(x),y] + [x,B(y)] + [x,y] \\
 &=&  [B(x),y]^{B(y)}  + [x,B(y)]^{B(x)} + [x,y]^{B(y)B(x)} 
 +  [B(x),B(y)] + [B(y),B(x)] \\
 &=& [x,yB(y)]^{B(x)}  +  [B(x),yB(y)] + [B(y),B(x)]\\
 &=& [xB(x),yB(y)] + [B(y),B(x)] \quad \mod \gamma_{i+j+1}(G).
\end{eqnarray*}
Finally, equation \eqref{LieRing} is a consequence of the Rota--Baxter identities satisfied in $G$, that is,
\begin{eqnarray*}
&& B(y)B(x) ~ B\big([xB(x),yB(y)] \, [B(y),B(x)]\big) \\
 &=& B \big(yB(y)xB(y)^{-1}\big) ~ B \big([xB(x),yB(y)] \, [B(y),B(x)] \big) \\
 &=& B\Big( yB(y)x B(y)^{-1}~B(y) B(x)~(xB(x))^{-1}(yB(y))^{-1}xB(x)y B(y)~B(y)^{-1} B(x)^{-1}B(y)B(x)~B(x)^{-1}B(y)^{-1} \Big) \\
 &= & B \big( yB(y)xB(x)B(x)^{-1}x^{-1} B(y)^{-1}y^{-1}xB(x)yB(x)^{-1}B(y)B(y)^{-1} \big) \\
 &=& B \big(yB(y)B(y)^{-1}y^{-1}xB(x)yB(x)^{-1} \big)\\
 &=& B(x)B(y),
\end{eqnarray*}
and the proof is complete.
$\blacksquare$ 
\end{proof}

\begin{example}{\rm 
Let $G$ be a 2-step nilpotent group and $B$ a Rota--Baxter operator on~$G$ defined as in Proposition~\ref{prop:center-conj}.
Then by Proposition~\ref{prop:LieRing}, we get a Rota--Baxter operator
$R: L(G)\to L(G)$, where $L(G) = G/[G,G]\oplus [G,G]$. A direct check shows that $R = -\id$.}
\end{example}
\para

An automorphism of a group  is called \index{fixed-point free automorphism}{\it fixed-point free} if it does not fix any non-trivial element. We now state a group theoretic analogue of \cite[Proposition 2.21]{MR3977738}, which was proven therein for Lie algebras (see also~\cite{MR0674005}).

\begin{theorem} \cite[Theorem 50]{MR4556953} \label{thm:invertibleRB}
Let $G$ be a finite non-abelian simple group and $B$ a Rota--Baxter operator on~$G$ such that $\ker(B) = 1$.
Then $B(g) = g^{-1}$ for all $g \in G$.
\end{theorem}

\begin{proof}
Since $B$ is a group homomorphism from $G_B$ to $G$ with trivial kernel, we conclude that $G_B\cong G$, and hence $G_B$ is also simple.
Suppose that $\ker(B_+)\neq G_B$, where $B_+(g) = gB(g)$ is a~homomorphism from $G_B$ to $G$. It follows from the simplicity of $G_B$ that $\ker(B_+) =1$. Thus, we have an automorphism $\varphi = B^{-1}B_+$ of the simple group $G_B$. Note that $\varphi$ is fixed-point free. Indeed, $g = \varphi(g) = B^{-1} \big(gB(g)\big)$ would imply that $B(g) = gB(g)$, that is, $g = 1$ By \cite{MR1334233}, a finite group admitting a fixed-point free automorphism is solvable, which is a contradiction. Hence, $\ker(B_+)= G_B$, that is, $B(g) = g^{-1}$ for all $g \in G$. $\blacksquare$ 
\end{proof}

\begin{corollary}
Let $G$ be a finite group and $B$ a Rota--Baxter operator on~$G$. If  $G_B$ is a non-abelian simple group, then $B$ is elementary and $G\cong G_B$.
\end{corollary}

\begin{proof}
Since $\ker(B)$ is a normal subgroup of $G_B$ and $G_B$ is simple, we have either $\ker(B) = G$ or $\ker(B)=1$. If $\ker(B) = G$, then $B$ is elementary and $G\cong G_B$. If $\ker(B)=1$, then $G\cong G_B$ by Theorem~\ref{thm:invertibleRB}.
$\blacksquare$ 
\end{proof}

\begin{corollary}\label{coro:invertibleRB}
Let $G$ be a finite non-abelian simple group. If $G$ is not factorizable, then $G$ is Rota--Baxter elementary.
\end{corollary}

\begin{proof}
Let $B$ be a non-elementary Rota--Baxter operator on $G$. By Theorem~\ref{thm:invertibleRB}, we have a proper factorization~\eqref{ImageFactorization},
which is a contradiction.
$\blacksquare$ 
\end{proof}


\chapter{Rota--Baxter groups and skew braces }\label{chap RBG and YBE}

\begin{quote}
This chapter investigates the relationship between Rota--Baxter groups and skew left braces. We show that every Rota--Baxter group naturally gives rise to a skew left brace, and conversely, any skew left brace with a complete additive group can be constructed in this way. The chapter also examines set-theoretic solutions to the Yang--Baxter equation that emerge from Rota--Baxter groups, along with the structure and characteristics of symmetric skew left braces. Additionally, we show that the obstruction to induction of a skew left brace from a Rota--Baxter operator lies in certain second group cohomology.
\end{quote}
\bigskip

\section{Relations between Rota--Baxter groups and skew braces}

We begin with the following result that gives a connection between Rota--Baxter groups and skew left braces \cite[Proposition 3.1]{MR4370524}.

\begin{proposition} \label{RBToSLB}
Let  $B$ be a Rota--Baxter operator on a group $(G,\cdot)$. If we set $x \circ_B y = xB(x)yB(x)^{-1}$ for all $x, y \in G$, then $(G, \cdot, \circ_B)$ is a skew left brace.
\end{proposition}

\begin{proof}
By~Proposition~\ref{prop:Derived}(1), $(G,\circ_B)$ is a group. Writing the product in $(G, \cdot)$ as juxtaposition, we see that
$$
(x\circ_B y) x^{-1}(x\circ_B z)
 = xB(x) y \, B(x)^{-1} x^{-1}xB(x) z B(x)^{-1}
 = xB(x)y z B(x)^{-1}
 = x\circ_B(yz)
$$
for all $x, y, z \in G$, which is desired. $\blacksquare$ 
\end{proof}

For brevity,  the skew left brace $(G, \cdot, \circ_B)$ induced by the Rota--Baxter group $(G, B)$ will be denoted by $G(B)$.

\begin{corollary}
Let  $B$ be a Rota--Baxter operator on a group $(G,\cdot)$.  If $\im (B)\subseteq \Z(G)$, then the skew left brace $(G, \cdot, \circ_B)$ is trivial.
\end{corollary}

\begin{example}{\rm 
Let $G$ be a group and $B_0(g)=1$ and $B_{-1}(g)=g^{-1}$ for all $g \in G$. Then $G(B_0)$ is a trivial skew left brace. Similarly, $G(B_{-1})$ is a trivial skew left brace, since $x \circ y = y \cdot x$ for all $x, y \in G$.}
\end{example}

For Rota--Baxter operators of weight $-1$, we have the following analogue of Proposition~\ref{RBToSLB}.

\begin{proposition}
Let $C$ be a~Rota--Baxter operator of weight $-1$ on a group $(G,\cdot)$. If we set $x \circ_C y = C(x)yC(x)^{-1}x$ for all $x,y \in G$, then $(G, \cdot, \circ_C)$ is a~skew left brace. If $(G, \cdot)$ is an abelian group, then
$(G, \cdot, \circ_C)$ is a trivial skew left brace.
\end{proposition}

\begin{proposition}
Let  $B$ be a Rota--Baxter operator on a group $(G,\cdot)$. Then the group $(G, \circ_B)$ is abelian if and only if 
\begin{equation} \label{RBlsbAbelianCirc}
[x^{-1}, B(y)] \,[B(x), y^{-1}] = [x^{-1}, y^{-1}]
\end{equation}
 for all $x, y \in G$.
\end{proposition}

\begin{proof}
By the definition of $\circ_B$, the identity
$x\circ_B y = y\circ_B x$ is equivalent to the identity
$$
x B(x) y B(x)^{-1} = y B(y) x B(y)^{-1}.
$$
If $[a, b]$ denotes $aba^{-1}b^{-1}$, then we rewrite the preceding identity as
$$
x y [y^{-1}, B(x)] = y x [x^{-1},B(y)]
$$
Since, $yx = xy[y^{-1}, x^{-1}]$,  we get
$$
[y^{-1}, B(x)] = [y^{-1}, x^{-1}] \,[x^{-1}, B(y)],
$$
which is equivalent to~\eqref{RBlsbAbelianCirc}. $\blacksquare$ 
\end{proof}

\begin{proposition}
Let $B$ be a Rota--Baxter operator on a group $(G,\cdot)$. Then $(G,\circ_B)\cong (G,\circ_{\widetilde{B}})\cong (G,\circ_{B^{(\varphi)}})$ for each $\varphi\in\Aut(G, \cdot)$.
\end{proposition}

\begin{proof} 
For $g, h \in G$, we have
$$
g\circ_{\widetilde{B}}h
 = g \big(g^{-1}B(g^{-1}) \big) h \big(g^{-1}B(g^{-1})\big)^{-1}
 = B(g^{-1})h B(g^{-1})^{-1}g
 = (g^{-1}\circ_B h^{-1})^{-1},
$$
and hence $x\to x^{-1}$ is an isomorphism of groups $(G,\circ_B)$ and $(G,\circ_{\widetilde{B}})$.
Similarly,
\begin{eqnarray*}
g\circ_{B^{(\varphi)}}h &=& g\varphi^{-1} \big(B(\varphi(g))\big) h \big(\varphi^{-1} \big(B(\varphi(g))) \big)^{-1}\\
 &=& \varphi^{-1} \big( \varphi(g)B(\varphi(g))\varphi(h)B(\varphi(g))^{-1} \big) \\
 &=& \varphi^{-1} \big(\varphi(g)\circ_B\varphi(h)\big),
\end{eqnarray*}
and hence $\varphi(g\circ_{B^{(\varphi)}}h) = \varphi(g)\circ_B \varphi(h)$,
and the corresponding groups are isomorphic. In fact, $\varphi$ is an isomorphism of skew left braces
$(G,\cdot,\circ_B)$ and $(G,\cdot,\circ_{B^{(\varphi)}})$. $\blacksquare$ 
\end{proof}

The next result relates any skew left brace to a Rota--Baxter group \cite[Theorem 3.6]{MR4370524}.

\begin{theorem}  \label{Embedding}
Every skew left brace can be embedded into a Rota--Baxter group.
\end{theorem}

\begin{proof}
Let $(G,\cdot,\circ)$ be a skew left brace.
Consider the semi-direct product $\widetilde{G} = G\ltimes G$ with the group operation
\begin{equation}\label{CircProduct}
(x,y)*(z,t) = \big(x\circ z, ~y\lambda_x(t)\big).
\end{equation}
Consider the subgroups
$H = \{(g,g) \, \mid \, g\in G\}$ and $L = \{(g,1) \,\mid \, g\in G\}$ of $\widetilde{G}$.
Since $(x,y) = (y,y)*(y^{\circ(-1)}\circ x,1)$, we can write $\widetilde{G} = H*L$.  We define a splitting Rota--Baxter operator $B: \widetilde{G} \to \widetilde{G}$ by $B(h*l) = l^{*(-1)}$ for $h \in H$ and $l \in L$. Thus, $B \big((x,y)\big) = (x^{\circ(-1)}\circ y,1)$, and  $\widetilde{G}$ is a Rota--Baxter group.
\para

Define $\psi: G \to \widetilde{G}$ by $\psi(g)= (1,g)$. We claim that $\psi$ is an isomorphism of skew left braces $(G, \cdot, \circ)$ and $\im(\psi)$, where $\im(\psi)$ is considered
as a~subbrace of
$\widetilde{G}(B) = (\widetilde{G},*,\circ_B)$. Indeed, we have
$$
\psi(g)*\psi(h)
 = (1,g)*(1,h)
 = (1,g\cdot h)
 = \psi(g\cdot h)
$$
and
\begin{eqnarray*}
\psi(g) \circ_B \psi(h)
 &=& (1,g)*B \big((1,g)\big)*(1,h) * \big(B((1,g)) \big)^{-1} \\
 &=& (1,g)*(g,1)*(1,h)*(g^{\circ(-1)},1)\\
 &=& (g,g)*(g^{\circ(-1)},h)\\
 &=& \big(1,g\lambda_g(h)\big)\\
 &=& (1,g\circ h)\\
 &=& \psi(g\circ h)
\end{eqnarray*}
for all $g, h \in G$, which is desired. $\blacksquare$ 
\end{proof}

The following example illustrates Theorem \ref{Embedding}.

\begin{example}{\rm 
Let $G = (\mathbb{Z}, +, \circ)$ be the left brace, where $+$ is the usual addition of integers and
$$a \circ b = a + (-1)^a b$$ for  $a, b \in \mathbb{Z}$. Note that the inverse of an element $a$ under $\circ$ is given by $a^{\circ(-1)} = (-1)^{a+1} a$.
We now construct a Rota--Baxter group containing the left brace $G$. Define a group operation on the set
$\widetilde{G} = \mathbb{Z} \times \mathbb{Z}$ by setting
$$
(a, b) * (c, d) = \big(a+(-1)^a c,~ b +(-1)^a d \big)
$$
for $(a, b),(c, d) \in \widetilde{G}$. For each $a \in \mathbb{Z}$, the inverse of $(a, 0)$ in $\widetilde{G}$ is
$$
(a, 0)^{*(-1)} = \big((-1)^{a+1} a, 0 \big) = (a^{\circ(-1)}, 0).
$$
Consider the subgroups $H = \{(g,g) \,\mid \, g\in \mathbb{Z}\}$ and $L = \{(g,0) \,\mid \, g\in \mathbb{Z}\}$ of $\widetilde{G}$.
Since
$$
(x,y)
 = (y,y)*(y^{\circ(-1)}\circ x,0)
 = (y, y)* \big((-1)^{y+1}(y-x), 0\big),
$$
we can write $\widetilde{G} = H*L$. Define a Rota--Baxter operator $B: \widetilde{G} \to \widetilde{G}$ by
$$
B \big((g, h)\big)
 = \big((-1)^{h+1}(h-g),0\big)^{*(-1)}
 = \big((-1)^{g+1}(h-g),0\big)
$$
for $(g,h)\in G$. Hence, $G$ embeds into the Rota--Baxter group $(\widetilde{G}, *,B)$ via the map $g\mapsto(0,g)$.}
\end{example}

\begin{remark} \label{tilde{G}Nilpotent}
{\rm
The construction of the group $\widetilde{G}$ appeared first in~\cite[Section 4]{MR3957824} in the context of nilpotent skew left braces. More precisely, we can define a new binary operation~$*$ on a~skew left brace~$(G, \cdot, \circ)$ by
$$g* h = g^{-1} \cdot (g\circ h) \cdot h^{-1}.$$
Note that, in terms of the group~$\widetilde{G}$, we have
$$
[(1,h),(g,1)]
 = (1,g* h).
$$
This observation helped Ced\'{o}, Smoktunowicz, and Vendramin \cite{MR3957824} to clarify the construction of left and right series of~$G$, and hence led to the definition of left (right) nilpotency of~$G$. Note that, if the skew left brace arise from a Rota--Baxter group $(G,B)$, then $g* h = [B(g)^{-1},h^{-1}]$.}
\end{remark}

In \cite[Problem 19.90(d)]{Kourovka}, the question was posed whether a finite skew left brace can exist with a non-solvable additive group and a nilpotent multiplicative group. This question was negatively answered by Tsang and Qin in \cite{MR4077413}. An alternative proof of this result can be obtained using Theorem \ref{Embedding}.

\begin{corollary}\label{coro:Kegel}
Let $(G,\cdot,\circ)$ be a finite skew left brace. If $(G, \circ)$ is nilpotent, then $(G, \cdot)$ is solvable.
\end{corollary}

\begin{proof}
Given a finite skew left brace~$(G,\cdot,\circ)$, we construct $(\widetilde{G}, *)$ as in the proof of Theorem~\ref{Embedding}.
We also have the decomposition $\widetilde{G} = H*L$, where
$H = \{(g,g) \,\mid\,  g\in G\}$ and $L = \{(g,1)\,\mid\, g\in G\}$.
Note that $H\cong L\cong (G,\circ)$. Since $(G, \circ)$ is nilpotent, it follows from Kegel's theorem \cite{MR0133365} that $(\widetilde{G}, *)$ is solvable. Thus, its subgroup $\{(1,g) \,\mid \, g\in G\}\cong (G, \cdot)$ is also solvable.
$\blacksquare$ 
\end{proof}

Another interesting application of Theorem~\ref{Embedding} is the following result \cite{MR4003478,MR4077413}.

\begin{corollary}
Let $(G,\cdot,\circ)$ be a skew left brace. If $(G, \circ)$ is abelian, then $(G, \cdot)$ is metabelian.
\end{corollary}

\begin{proof}
We have the factorization $\widetilde{G} = H*L$. By Ito's theorem~\cite{MR0071426}, $(\widetilde{G}, *)$ is metabelian.
Thus, its subgroup $\{(1,g)\, \mid \, g\in G\}\cong (G, \cdot)$ is also metabelian.
$\blacksquare$ 
\end{proof}

A group $G$ is called \index{complete group}{\it complete} if both its center $\Z(G)$ and its outer automorphism group $\Out(G)$ are trivial.

\begin{proposition}\label{complete SLB is RBO}
Let $(G,\cdot,\circ)$ be a skew left brace such that $(G, \cdot)$ is complete. Then there exists a Rota--Baxter operator~$B$ on $(G, \cdot)$ inducing $(G,\cdot,\circ)$.
\end{proposition}

\begin{proof}
Since $\Out(G,\cdot)$ is trivial, for each $g \in G$, there exists $\tilde{g}\in G$ such that
$\lambda_g(x) = \tilde{g} x {\tilde{g}}^{-1}$ for all $x \in G$. Since $\Z(G)$ is trivial, such an element $\tilde{g}$ is unique. Define $B:G \to G$ by $$B(g)=\tilde{g}$$ for all $g \in G$. We claim that $B$ is the desired Rota--Baxter operator on $(G, \cdot)$. For $g, h\in G$, we see that 
$$\widetilde{(g \circ h)} \, x \,{\widetilde{(g \circ h)}}^{-1}= \lambda _{g \circ h} (x)=\lambda _g \lambda_h (x)= \lambda_g(\tilde{h} x \tilde{h}^{-1})=\tilde{g} (\tilde{h} x \tilde{h}^{-1} )\tilde{g}^{-1}= (\tilde{g} \tilde{h}) x (\tilde{g}\tilde{h})^{-1}$$
for all $x \in G$. Since $\Z(G)$ is trivial, it follows that $\widetilde{g \circ h} = \tilde{g} \tilde{h}$. This means that  $B(g \circ h)= B(g) B(h)$ for all $g, h \in G$. Further, $g \circ h= g  \lambda_g(h)=g (\tilde{g} h\tilde{g} ^{-1})= g B(g) h B(g)^{-1}= g \circ_B h$ for all $g, h \in G$. This completes the proof. $\blacksquare$ 
\end{proof}

\begin{remark}\label{remark lambda RBO}
{\rm 
If $B$ is a Rota--Baxter operator on a group $(G,\cdot)$. Then the map $\lambda: (G, \circ_B) \to \Aut(G, \cdot)$ is given by $\lambda_g(h) = B(g)hB(g)^{-1}$ for all $g, h \in G$. Thus, for a skew left brace arising from a Rota--Baxter group, the map $\lambda$ takes values in the  inner automorphism group.}
\end{remark}

In~\cite[Problem~12]{MR3974481}, Vendramin  posed the problem of describing a free skew left brace. The following result provides some insight into the structure of a skew left brace arising from a free Rota--Baxter group \cite[Proposition 3.14]{MR4370524}.

\begin{proposition} 
Let $(G,B)$ be a free Rota--Baxter group generated by a~set~$X$. Let $S$ denote the skew left subbrace of $G(B)$ generated by~$X$.
Then $S$~is a~free skew left brace generated by~$X$.
\end{proposition}

\begin{proof}
Let $(A,\cdot,\circ)$ be a~skew left brace generated by~$X$.
By Theorem~\ref{Embedding}, $A$~embeds into the Rota--Baxter group~$\tilde{A}$,
which itself is generated by the set~$X$ as a Rota--Baxter group.
Thus, there exists a surjective homomorphism $\psi: G \to \tilde{A}$ of Rota--Baxter groups.
Hence, $\psi$~is a~surjective homomorphism from $S$ onto $A$
considered as a~subbrace of~$\tilde{A}$.
Since $A$ is an arbitrary skew left brace generated by~$X$, $S$ is also generated by~$X$ and $A$ is a homomorphic image of~$S$,
we conclude that $S$~is the free skew left brace generated by~$X$.
$\blacksquare$ 
\end{proof}

In the preceding chapter, we gave some constructions of Rota--Baxter groups. We now analyse skew left braces induced by these Rota--Baxter groups. 

\begin{example}
{\rm 
Let $G = HL$ be an exact factorization. Let $B$ be the splitting Rota--Baxter operator on $G$, that is, $B(hl) = l^{-1}$ for all $h \in H$ and $l \in L$. Then the multiplicative group structure of the skew left brace $G(B)$ is given by
$$
x \circ_B y
 = (hl)\circ_B y
 = hlB(hl)yB(hl)^{-1}
 = hyl.
$$ 
Such skew left braces appeared in several works, for instance, \cite[Example~1.6]{MR3647970},
~\cite[Theorem~2.3]{MR3763907} and~\cite[p. 19]{MR4346001}.}
\end{example}

\begin{example}
{\rm 
Recall the Rota--Baxter operator given by Proposition \ref{triangular}. Let $G$ be a group such that $G = HLM$, where $H$, $L$ and $M$ are subgroups of $G$ with pairwise trivial intersections. Let $C$ be a Rota--Baxter operator on~$L$ and  $[H,L] = [C(L),M] = 1$. Then the map $B: G\to G$ defined by 
$$
B(hlm) = C(l)m^{-1}
$$
is a Rota--Baxter operator on~$G$. For $x=hlm$ and $y=h'l'm'$, we have
\begin{eqnarray*}
x\circ_B y  &=& (hlm)\circ_B (h'l'm')\\
 &=& (hlm)B(hlm)h'l'm'B(hlm)^{-1}\\
 &=& hlmC(l)m^{-1}h'l'm'mC(l)^{-1} \\
 &=& hl mm^{-1}C(l)h'l'm'mC(l)^{-1}\\
 &=& hh'lC(l)l'm'mC(l)^{-1}\\
 &=& hh'(l\circ_C l') m'm.
\end{eqnarray*}
Hence, we have $(G,\circ_B)\cong H\times (L, \circ_C)\times M^{\rm op}$, where $M^{\rm op}$ is the opposite group of $M$.}
\end{example}

\begin{example}
{\rm 
Consider the Rota--Baxter operator given by Proposition \ref{semi-direct}. Let $G = H\rtimes L$ be a semi-direct product and $C$ a~Rota--Baxter operator on $L$. Then the map $B: G\to G$ defined by  $B(hl) = C(l)$, where $h\in H$ and $l\in L$, is a~Rota--Baxter operator on $G$.
For $x=hl$ and $y=h'l'$, we have
\begin{eqnarray*}
x\circ_B y &=& (hl)\circ_B (h'l')\\
 &=& hlB(hl)h'l'B(hl)^{-1}\\
 &=& hlC(l)h'l'C(l)^{-1}\\
 &=& hh'^{C(l)^{-1}l^{-1}}lC(l)l'C(l)^{-1} \\
 &=& hh'^{C(l)^{-1}l^{-1}}(l\circ_C l'),
\end{eqnarray*}
and hence $(G, \circ_B)\cong H\rtimes (L, \circ_C)$.}
\end{example}

\begin{example}\label{iso slb induced by differnt rbo}
{\rm 
Let $\Sigma_3$ be the symmetric group generated by two transpositions $s_1$ and $s_2$.
\begin{enumerate}
\item Consider the  factorisation $\Sigma_3 = \langle s_2 \rangle A_3$ and the splitting Rota--Baxter operator $B_1 : \Sigma_3 \to \Sigma_3$ given by $B_1(c a) = a^{-1}$ for $c\in \langle s_2 \rangle$ and $a \in A_3$. More precisely, we have
$$
B_1(s_1) = s_1 s_2,\quad
B_1(s_2) = 1,\quad
B_1(s_1 s_2) = s_2 s_1, \quad
B_1(s_2 s_1) = s_1s_2 \quad \textrm{and} \quad
B_1(s_1 s_2 s_1) = s_2 s_1.
$$
If $(\Sigma_3, \cdot, \circ)$ is the skew left brace  induced by $B_1$, then $(\Sigma_3, \circ)$ is the cyclic group of order six generated by $s_1$ such that
$$
s_1^{\circ 2} = s_1 \circ s_1 = s_1 s_2,\quad
s_1^{\circ 3} = s_2, \quad
s_1^{\circ 4} = s_2 s_1, \quad
s_1^{\circ 5} = s_1 s_2 s_1  \quad \textrm{and} \quad
s_1^{\circ 6} = 1.
$$
\item Consider the Rota--Baxter operator $B_2 : \Sigma_3 \to \Sigma_3$ given by
$$
B_2(s_1) = s_1, \quad
B_2(s_2) = s_1, \quad
B_2(s_1 s_2 s_1) = s_1, \quad
B_2(s_2 s_1) = 1 \quad \textrm{and} \quad
B_2(s_1 s_2) = 1.
$$
If $(\Sigma_3, \cdot, \circ)$ is the skew left brace  induced by $B_2$, then we see that the group $(\Sigma_3, \circ)$ is the cyclic group of order six generated by $s_2$ such that
$$
s_2^{\circ 2} = s_1 s_2, \quad
s_3^{\circ 3} = s_1, \quad
s_2^{\circ 4} = s_2 s_1, \quad
s_5^{\circ 5} = s_1 s_2 s_1 \quad \textrm{and} \quad
s_2^{\circ 6} = 1.
$$
\item[]  It is not difficult to see that the Rota--Baxter groups considered above are non-isomorphic, yet they induce isomorphic skew left braces.
\end{enumerate}}
\end{example}

Recall the Definition \ref{definition of hom of RBG} of a homomorphism of Rota--Baxter groups.

\begin{proposition}
Let $\varphi: (G_1, B_1) \to (G_2, B_2)$ be a homomorphism of Rota--Baxter groups with underlying groups $(G_1, \cdot_1)$ and $(G_2, \cdot_2)$. Then $\varphi: (G_1, \cdot_1, \circ_{B_1}) \to (G_2, \cdot_2, \circ_{B_2})$ is a homomorphism of induced skew left braces. In particular, an isomorphism of Rota--Baxter groups gives rise to an isomorphism of induced skew left braces.
\end{proposition}

\begin{proof}
It only needs to be shown that $\varphi$ preserves the group structure induced by  the Rota--Baxter operators. Indeed, for $g, h \in G_1$, we have
\begin{eqnarray*}
\varphi (g \circ_{B_1} h) &=& \varphi \big(g \cdot_1{B_1}(g) \cdot_1h \cdot_1 {B_1}(g)^{-1} \big)\\
&=& \varphi (g) \cdot_2 \varphi \big(B_1(g)\big)\cdot_2 \varphi(h)\cdot_2 \varphi \big(B_1(g)^{-1}\big)\\
&=& \varphi (g) \cdot_2 {B_2} \big(\varphi(g)\big) \cdot_2 \varphi(h) \cdot_2 {B_2} \big(\varphi(g)\big)^{-1}\\
&=& \varphi (g) \circ_{B_2} \varphi (h),
\end{eqnarray*}
which is desired. The second assertion is immediate.
$\blacksquare$ 
\end{proof}

\begin{remark}
{\rm 
Mapping a Rota--Baxter group to its induced skew left brace gives a functor from the category $\mathcal{RBG}$ of Rota--Baxter groups to the category $\mathcal{SLB}$ of skew left braces.}
\end{remark}

We now examine the solution to the Yang--Baxter equation derived from a Rota--Baxter group \cite[Theorem 6.1]{MR4370524}.

\begin{theorem} 
Let $(G, B)$ be a Rota--Baxter group and $\lambda_a(b) = B(a) b B(a)^{-1}$ for $a, b \in G$. Then $r : G \times G \to G\times G$ given by
\begin{equation} \label{RBE}
r(a, b) = \big(\lambda_a(b), ~a^{\lambda_a(b)B(\lambda_a(b))}\big),
\end{equation}
gives a~non-degenerate bijective solution $(G, r)$ to the Yang--Baxter equation. Moreover, $S$ is involutive if and only if $G$ is abelian.
\end{theorem}

\begin{proof}
By Proposition \ref{RBToSLB}, taking $a \circ_B b = a B(a) b B(a)^{-1}$, we obtain the skew left brace $(G, \cdot, \circ_B)$, where $\cdot$ is the group operation in $G$. In this case, the map $\lambda: (G, \circ_B) \to \Aut(G, \cdot)$ is given by $\lambda_a(b) = B(a) b B(a)^{-1}$ for $a, b \in G$. We see that
$$
(a \circ_B b)^{-1} a (a \circ_B b)
 = B(a) b^{-1} B(a)^{-1} a B(a) b B(a)^{-1}
 = a^{B(a) b B(a)^{-1}}
$$
and
$$ \lambda^{-1}_{\lambda_a(b)} \big((a \circ_B b)^{-1} a (a \circ_B b) \big) = B \big(\lambda_a(b)\big)^{-1}a^{B(a) b B(a)^{-1}} B \big(\lambda_a(b)\big)= a^{B(a) b B(a)^{-1} B(\lambda_a(b))}= a^{\lambda_a(b) B(\lambda_a(b))}.$$
By Theorem \ref{solution from a skew brace}, for a skew left brace $(G,\cdot, \circ_B)$, the map
\begin{equation}
r(a, b) = \big(\lambda_a(b), ~\lambda^{-1}_{\lambda_a(b)}((a \circ_B b)^{-1}  a  (a \circ_B b))\big)
\end{equation}
gives a~non-degenerate bijective solution $(G, r)$ to the Yang--Baxter equation. Moreover, $r$~is involutory if and only if $a b = b  a$ for all $a, b \in G$. This completes the proof. $\blacksquare$   
\end{proof}

We conclude with the following examples.

\begin{example}
{\rm 
If a group $G$ has an exact factorization $G = HL$, then the binary operation
$$(a_1 b_1) \circ (a_2 b_2) = a_1 a_2 b_2 b_1 $$
for $a_i \in H$ and $b_i \in L$, gives a skew left brace $(G, \cdot, \circ)$. Further, the map 
$B : G \to G$ given by $B(ab)= b^{-1}$ is a Rota--Baxter operator on $G$. Hence, we have a non-degenerate bijective solution to the Yang--Baxter equation with braiding
$$
r(a_1b_1,~a_2b_2)=
\big(b_1^{-1}a_2b_2b_1,~ b_2(a_2b_2b_1)^{-1}(b_1a_1)(a_2b_2b_1)b_2^{-1}\big) = \big( (a_2 b_2)^{b_1},~ (b_1 a_1)^{a_2 b_2 b_1 b_2^{-1}} \big).
$$}
\end{example}

\begin{example}
{\rm 
Recall from Proposition \ref{exact factorisation of free group} that the free group $F_n = \langle x_1, \ldots, x_n \rangle$ of rank $n \ge 2$ has an exact factorization $F_n = AB$, where
$$
A =  \big\langle x_1, x_2, \ldots, x_{n-1}, [x_n, F_n] \rangle \quad \textrm{and}\quad B =  \langle x_n \big\rangle.
$$
For a reduced word $w$ representing an element of $F_n$ and $1 \le i \le n$, let $\mathfrak{L}_i(w)$ be the sum of the powers  of $x_i$ in $w$. Each element $u \in F_n$ can be written in the form $u=(u x_n^{-\mathfrak{L}_n(u)}) x_n^{\mathfrak{L}_n(u)}$, where $u x_n^{-\mathfrak{L}_n(u)} \in A$ and $x_n^{\mathfrak{L}_n(u)} \in B$. Let $r: F_n \times F_n \to F_n \times F_n$ be the map given by
$$
r(u, v) = \big(v^{\pi(u)}, ~u^{\pi(u)^{-1} v \pi(u) \pi(v)^{-1}} \big),
$$
for $u, v \in F_n$, where $\pi(u) = x_n^{\mathfrak{L}_n(u)}$ and $\pi(v) = x_n^{\mathfrak{L}_n(v)}$. Then $(F_n, r)$ is a non-degenerate bijective solution to the Yang--Baxter equation.}
\end{example}
\bigskip
\bigskip

\section{Structure of skew braces  arising from Rota--Baxter groups}\label{TermsViaRB}

Recall from Proposition \ref{bijective 1-cocycle and skew brace} that, there is a one-to-one correspondence between  skew  left braces and bijective  1-cocycles.
\para

 If $(G, \cdot)$ is a~group and $\pi : \Gamma \to G$ is a bijective group 1-cocycle, then the binary operation given by
\begin{equation} \label{1-cocycle}
a \circ b = \pi \big(\pi^{-1}(a) \,\pi^{-1}(b) \big)
\end{equation}
for $a,b \in G$, gives a skew left brace $(G, \cdot, \circ)$. Conversely, suppose that $(G, \cdot, \circ)$ is a~skew left brace. Taking $\Gamma = G$ with group operation $(a,b)\mapsto a\circ b$ and the action $a \mapsto \lambda_a$ of $\Gamma$ on $G$,  we see that the identity map $\id:\Gamma \to G$ is a bijective 1-cocycle. In particular, when $G$ acts on itself by conjugation, the inverse of a bijective 1-cocycle given by \eqref{1-cocycle} is simply a Rota--Baxter operator on~$G$~\cite{MR4271483}.

\begin{proposition} \cite[Proposition 5.1]{MR4370524}
Let $B$ be a Rota--Baxter operator on a group $(G,\cdot)$ and  $(G,\cdot, \circ_B)$ the induced skew left brace. Then $(G,\cdot,\circ_B)$ is a two-sided skew  left brace
if and only if the map $\psi_g: G^{\rm op} \to G^{\rm op}$
given by $\psi_g(x)= [B(x)^{-1},g]$ is a~1-cocycle for each $g\in G$.
\end{proposition}

\begin{proof}
By~Proposition~\ref{RBToSLB}, $G$ is a skew left brace. Rewriting the axiom 
$$
(b \cdot c) \circ a  =  (b \circ a) \cdot a^{-1} \cdot (c\circ a)
$$
of a  skew  right  brace in terms of $B$, we get
$$
abB(ab)cB(ab)^{-1}
 = aB(a)cB(a)^{-1}c^{-1}bB(b)cB(b)^{-1},
$$
which is further equivalent to
$$
b[B(ab)^{-1},c^{-1}]
 = [B(a)^{-1},c^{-1}]b[B(b)^{-1},c^{-1}].
$$
Thus, $\psi_{c^{-1}}(ab) = \psi_{c^{-1}}(a)^b \psi_{c^{-1}}(b)$, and hence $\psi_g(b*a) = \psi_g(b)* \big(b\cdot \psi_g(a)\big)$ for all $a,b,g\in G^{\rm op}$.
Here, $b*a = ab$ is the opposite group operation in~$G$. $\blacksquare$ 
\end{proof}

The following result gives a criteria for a Rota--Baxter group induced skew left brace to be $\lambda$-homomorphic \cite[Proposition 5.14]{MR4370524}.

\begin{proposition}
Let $B$ be a Rota--Baxter operator on a group $(G,\cdot)$ and  $(G,\cdot, \circ_B)$ the induced skew left brace. Then $(G,\cdot,\circ_B)$ is a $\lambda$-homomorphic skew left brace if and only if $B(a \cdot b)^{-1} \cdot B(a) \cdot B(b)\in \Z(G,\cdot)$ for all $a,b\in G$.
\end{proposition}

\begin{proof}
Note that $\lambda_{a \cdot c}(b) = B(a \cdot c) \cdot b \cdot B(a \cdot c)^{-1}$ for all $a, b, c \in G$. On the other hand, we have
$$
\lambda_{a} \lambda_{c}(b) = \lambda_{a} \big( B(c) \cdot b \cdot B(c)^{-1} \big) = B(a) \cdot B(c) \cdot b \cdot B(c)^{-1} \cdot B(a)^{-1}
$$
for all $a, b, c \in G$. Hence, $\lambda: (G, \cdot) \to \Aut (G, \cdot)$ is a  homomorphism if and only if $b\cdot B(a\cdot c)^{-1}\cdot B(a)\cdot B(c) = B(a \cdot c)^{-1}\cdot B(a)\cdot B(c) \cdot b$ for all $a,b,c \in G$. $\blacksquare$ 
\end{proof}

Recall that, an ideal of a skew left brace $(G, \cdot, \circ)$ is a subset $I$ of $G$ such that $(I, \cdot)$ is a normal subgroup of $(G,\cdot)$, $(I, \circ)$ is a normal subgroup of $(G, \circ)$ and $\lambda_a(I) \subseteq I$ for any $a \in G$.  The condition $\lambda_a(I) \subseteq I$ is equivalent to the equality $a\circ I = a \cdot I$ of cosets for all $a\in G$. Thus, the notion of an ideal allows us to consider the quotient skew left brace $G/I$. 

\begin{definition}
A \index{left ideal}{\it left ideal} of a skew left brace $(G, \cdot, \circ)$  is a subgroup $(L, \cdot)$ of $(G, \cdot)$ such that $\lambda_a(L) \subseteq L$ for each $a \in G$.   A \index{strong left ideal}{\it strong left ideal} of $(G, \cdot, \circ)$ is a normal subgroup  $(L, \cdot)$ of $(G, \cdot)$ such that $\lambda_a(L) \subseteq L$ for each $a \in G$. 
\end{definition}

Clearly, every ideal of a skew left brace is a left as well as a strong left ideal. If $L$ is a left ideal of a skew left brace $(G, \cdot, \circ)$, then we have
$$
a \circ b = a \cdot \lambda_a(b) \in L \quad \textrm{and} \quad \bar{a} = \big(\lambda_{\bar{a}}(a) \big)^{-1} \in L
$$
for all $a, b \in L$. It follows that a left ideal of a skew left brace is always a subbrace. Strong left ideals provide decompositions of the corresponding solutions to the Yang--Baxter equation~\cite{MR4023387}, and are also connected with intermediate fields of Galois extensions~\cite{MR4033084}.
\para 

Inspired by Theorem~\ref{Embedding}, Remark~\ref{tilde{G}Nilpotent}, and the analogy with post-algebras, we apply the map $\psi:G \to \widetilde{G}= G \ltimes G$, given by  $\psi(g) = (1, g)$, to define some classes of skew left braces \cite{MR3957824}.

\begin{definition}\label{DefViaDelta}
Let~$(G, \cdot, \circ)$ be a~skew left brace and $\mathcal{X}$ a class of groups.
\begin{enumerate}
\item  A~non-empty subset $I$ of $G$ is said to be of type~$\mathcal{X}$
if $I$ is of type $\mathcal{X}$ in the group~$(G,\cdot)$.
\item  A~non-empty subset $I$ of $G$ is said to be of strong type $\mathcal{X}$
if $\psi(I)$ is of type $\mathcal{X}$ in the group~$\widetilde{G}$.
\end{enumerate}
\end{definition}

In \cite{MR4427114}, many notions of skew left braces $(G, \cdot, \circ)$ have been interpreted in terms of the group~$\widetilde{G}$. 

\begin{proposition} \label{IdealCriterion}
Let $(G, \cdot, \circ)$ be a skew left brace. A~non-empty subset~$I$ is a~strong left ideal of~ $(G, \cdot, \circ)$ if and only if
$\psi(I)$ is a normal subgroup of~$\widetilde{G}$.
\end{proposition}

\begin{proof}
Recall that $\psi:G \to \widetilde{G}$ is given by $\psi(g)=(1, g)$. Let $a\in G$ and $i\in I$. The condition that $\psi(I)$ is a normal subgroup of~$\widetilde{G}$ is equivalent to the inclusions
$$
(1,a)^{*(-1)}*(1,i)*(1,a)\in \{1\} \times I\quad \textrm{and}\quad (a,1)^{*(-1)}*(1,i)*(a,1)\in \{1\} \times I.
$$
By~\eqref{CircProduct}, we obtain  $a^{-1}ia\in I$ and $\lambda_{a^{\circ(-1)}}(i)\in I$. Since $a$ and $i$ are arbitrary, this implies that $\psi(I)$ is a normal subgroup of~$\widetilde{G}$
if and only if~$I$ is a~normal subgroup of $(G,\cdot)$ and $\lambda_x(I)\subseteq I$ for all $x\in G$. $\blacksquare$ 
\end{proof}

\begin{proposition}
Let $B$ be a Rota--Baxter operator on a group $(G,\cdot)$ and $I\subseteq G$. Then the following  assertions hold:
\begin{enumerate}
\item $I$~is an ideal of~$(G, \cdot, \circ_B)$ if and only if $I$~is a~normal subgroup of both
$(G, \circ_B)$ and $(G, \cdot)$.
\item $I$~is a left ideal of $(G, \cdot, \circ_B)$ if and only if $I$ is a subgroup of $(G, \cdot)$ normalized by $\im(B)$.
\end{enumerate}
\end{proposition}

\begin{proof}
If $I$ is an ideal of~$(G, \cdot, \circ_B)$, then $I$~is a~normal subgroup of both $(G, \cdot)$ and $(G, \circ_B)$ and $\lambda_a(I) \subseteq I$ for all $a \in G$. The last condition means that $B(a)iB(a)^{-1}\in I$ for all $a\in G$ and $i \in I$, which follows from the normality of~$I$ in $(G, \cdot)$, thereby establishing assertion (1). The proof of  assertion (2) is similar. $\blacksquare$ 
\end{proof}

Next, we introduce the notion of the left center of a~skew left brace.

\begin{definition}
Let $(G, \cdot, \circ)$ be a~skew left brace. The \index{left center of skew brace}{\it left center} of~$G$ is defined as 
$$
 \Z_l(G)
 = \{c\in G \, \mid \, c\in \Z(G,\cdot) ~\textrm{and}~
 g \cdot c = g\circ c\mbox{ for all }g\in G\}.
$$
\end{definition}

\begin{proposition}\label{strong let center prop}
The following assertions hold:
\begin{enumerate}
\item  If $(G, \cdot, \circ)$ is a skew left brace, then $\psi \big( \Z_l(G)\big) = \Z(\widetilde{G})\cap \psi(G)$.
\item If $G(B) = (G,\cdot,\circ)$ is a~skew left brace induced by a Rota--Baxter operator~$B$ on $(G, \cdot)$, then
$ \Z_l \big(G(B)\big) = \Z(G,\cdot)$.
\end{enumerate}
\end{proposition}

\begin{proof}
An element~$(1,s)$ lies in $\Z(\widetilde{G})$ if and only if $(a,b)*(1,s) = (1,s)*(a,b)$ for all $a,b\in G$. This means that $b \cdot \lambda_a(s) = s\cdot \lambda_1(b) = s\cdot b$ or equivalently $\lambda_a(s) = b^{-1}\cdot s\cdot b$. Taking~$a = 1$, we get $s\in \Z(G,\cdot)$. Thus, $s = \lambda_a(s) = a^{-1}\cdot (a\circ s)$ for all $a\in G$. This gives $a\cdot s = a\circ s$, which proves assertion (1).
\para 
If $c\in \Z_l \big(G(B)\big)$, then $c\in \Z(G, \cdot)$. Conversely, if $c\in \Z(G, \cdot)$, then $g\circ c = g \cdot B(g)\cdot c \cdot B(g)^{-1} = g \cdot c$ for all $g\in G$, which proves assertion (2). $\blacksquare$ 
\end{proof}

It follows from  Proposition~\ref{IdealCriterion} and Proposition \ref{strong let center prop} that $ \Z_l(G)$ is a strong left ideal of $(G, \cdot, \circ)$. In general, the left center is not an ideal of a skew left brace.

\begin{example}{\rm 
Let $A = \mathbb{Z}_3 = \langle a\rangle$ and 
$B = \mathbb{Z}_2 = \langle b\rangle$ be trivial braces.  Consider the homomorphism $\beta : B \to \Aut(A)$ defined such that $\beta(b):A \to A$ acts as $x\to x^2$. Then we have the semi-direct product $G = A\rtimes B$ with $(G,\cdot)\cong \mathbb{Z}_3\times\mathbb{Z}_2$ and $(G,\circ)\cong \Sigma_3$.  It is not difficult to see that $ \Z_l(G) = \{1\}\rtimes B$, which is not a normal subgroup of $(G,\circ)$.}
\end{example}

Recall that the socle of a skew left brace $(G, \cdot, \circ)$ is defined as
$$
\Soc(G)
 = \big\{a\in G \, \mid \, a\in \Z(G,\cdot)~\textrm{and}~ a\circ b = a \cdot b ~\mbox{ for all }b\in G\big\}.
$$

In analogy with the concept of the left center, we refer to $\Soc(G)$ as the right center of $(G, \cdot, \circ)$. The following result provides a description of the socles of skew left braces induced by Rota--Baxter operators.

\begin{proposition}
Let $B$ be a Rota--Baxter operator on a group $(G,\cdot)$ and $G(B)$ its induced skew left brace. Then
$$
\Soc(G(B)) = \big\{ a \in G \, \mid \, a,B(a) \in \Z(G, \cdot)\big\}
 = \Z(G, \cdot)\cap B^{-1}\big(\Z(G, \cdot)\big).
$$
\end{proposition}

\begin{proof}
The equality $a\circ b = a \cdot b$ is equivalent to $B(a)\cdot b\cdot B(a)^{-1} = b$, that is, $B(a)\in \Z(G, \cdot)$. $\blacksquare$ 
\end{proof}

In~\cite{MR3917122}, the \index{annihilator}{\it annihilator} of a~skew left brace~$(G, \cdot, \circ)$ is defined as $$\Ann(G) = \Z(G,\circ)\cap\Soc(G).$$ 
Thus, we can redefine the annihilator as $\Ann(G) = \Z_l(G)\cap\Soc(G)$. In~\cite{Bonatto}, $\Ann(G)$ is referred to as the {\it center} and has been applied to develop central nilpotency of skew left braces.
\para 

We define the \index{upper central series}{\it upper central series} of a skew left brace $(G, \cdot, \circ)$ by setting $\zeta_1(G) = \Z_l(G)$ and $\zeta_{n+1}(G)$ to be the strong left ideal of~$G$ such that
$$
\zeta_{n+1}(G)/\zeta_n(G) = \psi(G)/\zeta_n(G)\cap \Z \big(\widetilde{G}/\zeta_n(G)\big).
$$
Here, we identify $\zeta_i(G)$ with its image in $\widetilde{G}$ under the  map $\psi:G \to \widetilde{G}$.

\begin{definition}
A~skew left brace  $(G, \cdot, \circ)$ is called {\it strong left-nilpotent} if there exists an integer $m \ge 1$ such that $\zeta_m(G) = G$.
\end{definition}

Extending the idea of nilpotency of left braces by Rump \cite{MR2278047},  left- and right-nilpotency of skew left braces was defined in \cite{MR3957824}. The left series of a skew left brace $(G, \cdot, \circ)$ is the sequence of subsets $$G^1 \supseteq G^2 \supseteq G^3 \supseteq  \cdots \supseteq G^n \supseteq \cdots,$$ 
where  $G^1=G$, $G^{n+1} = G * G^n$ for each $n \ge 1$ and $$g* h = g^{-1} \cdot (g\circ h) \cdot h^{-1}.$$
A skew left brace $(G, \cdot, \circ)$ is said to be {\it left-nilpotent} if there exists $n \ge 1$ such that $G^n$ is trivial. 

\begin{proposition}
Let $(G, \cdot, \circ)$ be a strong left-nilpotent skew left brace. Then $(G, \cdot, \circ)$ is left-nilpotent and $(G,\cdot)$ is nilpotent.
\end{proposition}

\begin{proof}
Since $ \Z_l \big(\widetilde{G}/\zeta_n(G) \big) \subseteq \Z \big(\psi(G)/\zeta_n(G) \big)$,  we see that $\zeta_1(G),\zeta_2(G),\ldots$ is a central series of the group~$(G,\cdot)$, and hence  $(G,\cdot)$ is nilpotent. Suppose that $\zeta_m(G) = G$ for some $m \ge 1$. Then by definition of $\zeta_i(G)$, we have $G* \zeta_i(G)\subseteq \zeta_{i+1}(G)$. Hence, $G^{m+1}$ is trivial and $(G, \cdot, \circ)$ is left-nilpotent. $\blacksquare$ 
\end{proof}

We conclude the section with the following result \cite[Theorem 5.11]{MR4370524}.

\begin{theorem}
Let $(G, \cdot, \circ)$ be a strong left-nilpotent skew left brace and $I$ a non-trivial strong left ideal of $(G, \cdot, \circ)$. Then $ \Z_l(G)\cap I$ is non-trivial.
\end{theorem}

\begin{proof}
Suppose that $ \Z_l(G)\cap I$ is trivial. Let  $k\geq1$ be such that $\zeta_k(G)\cap I$ is trivial, but $\zeta_{k+1}(G)\cap I$ is non-trivial. Consider a non-trivial element $ i\in \zeta_{k+1}(G)\cap I$. By  definition of $\zeta_{k+1}(G)$, we have $i\in \Z \big(\widetilde{G}/\zeta_k(G)\big)$. Thus, $[i,g]\in \zeta_k(G)$ for all $g\in \widetilde{G}$. Since $I$ is a normal subgroup of~$\widetilde{G}$, we also have $[i,g]\in I$. By our assumption,  $\zeta_k(G)\cap I$ is trivial, and hence $[i,g] = 1$ for all $g\in\widetilde{G}$. Finally, we conclude that $i\in \Z(\widetilde{G})$, and hence $i\in \zeta_1(G)\cap I$, which is a~contradiction. $\blacksquare$ 
\end{proof}
\bigskip
\bigskip


\section{Symmetric skew braces} \label{sym}

In this section, we discuss symmetric skew left braces.

\begin{definition}
A skew  left brace $(G, \cdot, \circ)$ is said to be \index{symmetric skew  left brace}{symmetric} if  $(G, \circ, \cdot)$ is also a skew  left brace.
\end{definition}

These braces were introduced by Childs \cite{MR3982254} under the name bi-skew braces and studied subsequently by Caranti \cite{MR4130907}. The next result gives a characterisation of a symmetric skew  left brace \cite[Proposition 5.2]{MR4346001}.

\begin{proposition}  \label{symprop1}
A skew  left brace  $(G, \cdot, \circ)$ is symmetric if and only if
\begin{eqnarray} \label{inc1}
\overline{b}\circ (a\cdot b)\circ \overline{a} \in \ker(\lambda)
\end{eqnarray}
for all $a,b \in G$, where $\bar a$ denotes the inverse of $a$ in $(G,\circ)$.
\end{proposition}

\begin{proof}
Let $(G, \cdot, \circ)$ be a skew  left brace. Recall that $\lambda : (G, \circ) \to \Aut (G, \cdot)$ is a group homomorphism, where
\begin{equation}\label{eqn1sec5}
\lambda_a(b)=a^{-1} \cdot (a\circ b)
\end{equation}
 for all $a, b \in G$. Now, $( G, \circ, \cdot)$ is a symmetric skew  left brace if and only if
\begin{equation}\label{eqn2sec5}
a\cdot(b\circ c)=(a\cdot b)\circ \overline{a}\circ (a\cdot c)
\end{equation}
for all $a, b, c \in G$. Using \eqref{eqn1sec5}, we see that \eqref{eqn2sec5} is equivalent to
$$  
a \cdot b \cdot \lambda_b(c)=a \cdot b \cdot  \lambda_{a \cdot b} \big(\overline{a} \cdot \lambda_{\overline{a}} (a \cdot c) \big),$$
which is further equivalent to
$$
\lambda_b(c)=\lambda_{a \cdot b} \big(\overline{a} \cdot \lambda_{\overline{a}} (a) \cdot \lambda_{\overline{a}} (c)\big).
$$
Since $\lambda_{\overline{a}} (a)={\overline{a}}^{-1}\cdot (\overline{a}  \circ a)={\overline{a}}^{-1}$, we obtain
$$
\lambda_b(c)=\lambda_{a \cdot b} \big(\lambda_{\overline{a}} (c) \big) \Longleftrightarrow
\lambda_b(c)=\lambda_{(a \cdot b)\circ \overline{a}} (c) \Longleftrightarrow
\lambda_{\overline{b} \circ (a \cdot b)\circ \overline{a}} (c)=1.
$$
Since $a, b, c \in G$ are arbitrary elements, the assertion follows.
$\blacksquare$ 
\end{proof}
\para

\begin{remark}\label{symcor1}
{\rm 
The expression \eqref{inc1} is equivalent to $\lambda_{\overline{b} \circ (a \cdot b)\circ \overline{a}} =\id$, which can also be rewritten in the forms
$$
\lambda_{\overline{a}\circ\overline{b} \circ (a \cdot b) }=\id,\quad
\lambda_{b\circ a  }=\lambda_{a \cdot b } \quad \textrm{or} \quad
\lambda_{ (\overline{a \cdot b}) \circ (b \circ a)}=\id.
$$}
\end{remark}
\para

Let us consider some examples.

\begin{example}{\rm 
Consider the  left  brace $( \mathbb{Z}, +, \circ)$, where  $+$  is the usual addition and
$$
m \circ n=m+(-1)^{m}n
$$
for $m,n \in \mathbb{Z}$. In this case, $\lambda: (\mathbb{Z}, \circ) \to \Aut(\mathbb{Z}, +)$ is the homomorphism such that $\lambda(1) = -\id$ and
$\ker(\lambda)=2\mathbb{Z}$. Further, for $m,n \in \mathbb{Z}$, we have
$$
m\circ n\equiv m + n \mod\, 2 \quad \textrm{and}\quad  \overline{m}= (-1)^{m+1}m \equiv m \mod\, 2.
$$
Thus, we have
$$
\overline{m} \circ (m+n)\circ \overline{n} \equiv 2m+2n \equiv 0\mod\, 2.
$$
Hence, by Proposition \ref{symprop1}, $(\mathbb{Z}, \circ, + )$ is a  symmetric left  brace.}
\end{example}
\para

\begin{example}{\rm 
Let $\mathbb{Z}^n=\left\langle x_1,\ldots,x_n\right\rangle$ be the  free abelian group  of rank $n$  and 
$$
\mathbb{Z}^n_0=
\left\{\, \sum\limits_{i=1}^{n} a_ix_i \in  \mathbb{Z}^n \,  \mid \,  \sum a_i=0 \,\right\}.
$$
Let $\varphi\in \Aut (\mathbb{Z}^n, +)$ be an automorphism such that $\varphi(x_i) \equiv x_i\mod \mathbb{Z}^n_0$ for each $1 \le i  \le n$. Consider the left brace $(\mathbb{Z}^n, +, \circ)$, where $a \circ b=a+\varphi^{\sum a_i}(b)$. Note that 
$\mathbb{Z}^n_0=\ker(\lambda)$. Since
$$
a \circ b \equiv a+b \mod \mathbb{Z}^n_0 \quad \textrm{and} \quad \overline{a}=-\varphi^{-\sum a_i}(a)\equiv -a \mod \mathbb{Z}^n_0,
$$
we obtain
$$
\overline{b} \circ (a+b)\circ \overline{a} \equiv -b+a+b-a \equiv 0 \mod \mathbb{Z}^n_0.
$$
Hence, by Proposition \ref{symprop1}, $(\mathbb{Z}^n, \circ, +)$ is a symmetric skew  left brace.}
\end{example}

The next result gives more examples of symmetric skew left braces \cite[Theorem 5.8]{MR4346001}.

\begin{theorem} \label{t6.6}
Every $\lambda$-cyclic skew  left brace is symmetric.
\end{theorem}

\begin{proof}
Let   $(G, \cdot, \circ)$ be a $\lambda$-cyclic skew  left brace. Then  $\lambda : (G, \cdot) \to \Aut (G, \cdot)$ is a group homomorphism and $H_\lambda$ is a regular subgroup of $\Hol (G, \cdot)$. Thus, by Theorem \ref{t1},  we have
$$
\lambda \big(b^{-1} \cdot \lambda_a(b) \big)=\id,
$$
which is equivalent to
\begin{equation}\label{cyclic skew is symmetric}
\lambda \big((a \cdot b)^{-1} \cdot (a\circ b) \big)=\id.
\end{equation}
Since the image of $\lambda$ is a cyclic subgroup of  $\Aut (G, \cdot)$, it follows that  $\ker(\lambda)$ contains the commutator subgroup  of $(G, \cdot)$. Thus, \eqref{cyclic skew is symmetric} gives
$$
 \lambda \big((b \cdot a)^{-1} \cdot (a\circ b) \big)=\id,$$
 which further yields $\lambda(b \cdot a) = \lambda(a\circ b)$. Thus, by Proposition \ref{symprop1} and Remark \ref{symcor1}, it follows that  $(G, \cdot, \circ)$ is a symmetric skew  left brace. This completes the proof.
$\blacksquare$ 
\end{proof}

The next result shows that skew  left braces constructed from some exact factorizations are symmetric  \cite[Theorem 5.14]{MR4346001}.

\begin{theorem}
Let $G = A B$ be an exact factorization of a group $G$ such that $A$ is a normal subgroup of $G$. Then the skew  left brace obtained from this exact factorization is symmetric.
\end{theorem}

\begin{proof}
In view of Theorem \ref{skew left brace exact factor}, given the exact factorisation $G = AB$, we have the skew left brace $(G, \cdot, \circ)$, where
\begin{equation} \label{efeqn1}
x \circ y= a  y  b 
\end{equation}
for $x, y\in G$, where $x = a  b$. Further, for $x=ab \in G$, we have $\bar x = a^{-1}b^{-1}$. Hence, for $x_1=a_1b_1$ and $x_2=a_2b_2 \in G$, we have
$$
\bar x_1 \circ \bar x_2 \circ (x_1x_2)=a_1^{-1}a_2^{-1}a_1b_1a_2b_1^{-1}.
$$
Note that  $\lambda_x$ is the inner automorphism of $G$ induced by  $b$, that is,
$$
\lambda_x (y)=b^{-1}yb
$$
for all $y \in G$. Thus, if $x \in A$, then $\lambda_x = \id$.
Since  $A$ is a normal subgroup of  $G$, we get
$$
a_1^{-1}a_2^{-1}a_1b_1a_2b_1^{-1} \in A
$$
for all $a_1, a_2 \in A$ and $ b_1 \in B $. Hence, we have $\bar x_1 \circ \bar x_2 \circ (x_1x_2) \in \ker(\lambda)$. Thus, by  Proposition \ref{symprop1} and Remark \ref{symcor1},   $(G, \cdot, \circ)$ is a symmetric skew left brace, which completes the proof.
$\blacksquare$ 
\end{proof}

We conclude this section by introducing a generalization of a skew left brace, which might be of independent interest.

\begin{definition}
Let $k$ be a positive integer. A $(k+1)$-groupoid $(G, \circ_0, \circ_1, \circ_2, \ldots,\circ_k)$ is called a \index{skew left $k$-brace}{skew left $k$-brace} if the following conditions hold:
\begin{enumerate}
\item $(G, \circ_i)$ is a~group for each $0 \le i \le k$.
\item  For each $0 < i < k$ and $a, b, c \in G$, we have
$$
a \circ_{i} (b \circ_{i-1} c)
 = (a \circ_{i} b) \circ_{i-1} a^{\circ_{i-1}(-1)} \circ_{i-1} (a \circ_{i} c),
$$
where $a^{\circ_{i-1}(-1)}$ is the inverse to $a$ in the group $(G, \circ_{i-1})$. Note that a skew left 1-brace is simply a skew left brace.
\end{enumerate}
\end{definition}

\begin{proposition}
Let $B$ be a  Rota--Baxter operator on a group $(G, \cdot)$ and $k$ a positive integer. Define binary operations $\circ_0, \circ_1,\circ_2,\ldots,\circ_k$ on $G$ by
$$
x \circ_{i+1} y = x \circ_i B(x) \circ_i y \circ_i \big(B(x)\big)^{\circ_i(-1)}
$$
where $\circ_0 = \cdot$. Then
$(G,\circ_0, \circ_1, \circ_2,\ldots,\circ_k)$ is a~skew left $k$-brace.
\end{proposition}

\begin{proof}
The assertion follows directly from Proposition~\ref{prop:Derived}. $\blacksquare$ 
\end{proof}
\bigskip
\bigskip


\section{Obstruction to induction of skew braces from Rota--Baxter groups} \label{obstruction to SLB arising from RBG}
In view of Remark \ref{remark lambda RBO}, if $B$ is a Rota--Baxter operator on a group $(G,\cdot)$ and $(G, \cdot, \circ_B)$ is the induced skew left brace, then $$\lambda_g(h) = B(g)hB(g)^{-1}$$ for all $g, h \in G$. Thus, for a skew left brace induced by a Rota--Baxter group, the map $\lambda$ takes values in the inner automorphism group of $(G, \cdot)$.
\para

Let $(G, \cdot, \circ)$ be a skew left brace such that  $\lambda$ takes values in the inner automorphism group of $(G, \cdot)$.  It is natural to ask whether the skew left brace $(G, \cdot, \circ)$ (equivalently, the map $\lambda$)  is induced by a Rota--Baxter operator on $(G, \cdot)$. In view of Proposition \ref{complete SLB is RBO}, the answer is affirmative when the group $(G, \cdot)$ is complete. 
In general, it was shown in \cite{MR4531710} that the obstruction to the existence of a Rota--Baxter operator is captured by a certain second group cohomology. In this section, we explore this idea in detail.
\para

We begin by recalling the connection  between group extensions and the second group cohomology \cite{MR0672956}.  Let $G$ be a group and $A$ an abelian group admitting a right action of $G$ by automorphisms. In this case, we say that $A$ is a right $G$-module. The \index{group cohomology}{group cohomology} $\Ho^\ast_{\rm Grp}(G;A)$ of  $G$ with coefficients in $A$ is the cohomology of the cochain complex $C^*(G;A)=\{C^n(G;A),\partial^n\}$ with cochain groups $C^n(G;A) = \Map(G^n,A)$ and coboundary maps $\partial^n: C^n(G;A) \to C^{n+1}(G;A)$ given by
\begin{eqnarray}\label{group coboundary 1}
\partial^n f (g_1,\ldots,g_{n+1}) &=& f(g_2,\ldots,g_{n+1}) + (-1)^{n+1} f(g_1,\ldots,g_{n}) \cdot g_{n+1}\\
\nonumber  &&+ \sum_{i=1}^n (-1)^{i} f(g_1,\ldots, g_{i-1}, g_i g_{i+1}, g_{i+2}, \ldots,  g_{n+1})
\end{eqnarray}
for all $f \in C^n(G;A)$ and $(g_1,\ldots,g_{n+1}) \in G^{n+1}$. If the coefficients are in a left $G$-module~$A$, then the coboundary maps $\partial^n: C^n(G;A) \to C^{n+1}(G;A)$ are given by
\begin{eqnarray}\label{group coboundary 2}
\partial^n f (g_1,\ldots,g_{n+1}) &=& g_1 \cdot f(g_2,\ldots,g_{n+1}) + (-1)^{n+1} f(g_1,\ldots,g_{n})\\
 \nonumber &&+ \sum_{i=1}^n (-1)^{i} f(g_1,\ldots, g_{i-1}, g_i g_{i+1}, g_{i+2}, \ldots,  g_{n+1}). 
\end{eqnarray}
\para

For each $n \ge 0$, we write $\Ho^n_{\rm Grp}(G;A)= Z^n(G;A)/ B^n(G;A)$, where $Z^n(G;A)=\ker(\partial^n)$ is the group of $n$-cocycles and  $B^n(G;A)=\im(\partial^{n-1})$ is the group of $n$-coboundaries. In particular, for $n=2$ we have $\Ho^2_{\rm Grp}(G; A)=Z^2(G; A)/B^2(G; A)$, where
$$Z^2(G; A) =\big\{\theta:G \times G \to A  \, \mid \,   \big(x \cdot \theta(y,z) \big)\theta(x, yz)= \theta(xy, z)\theta(x, y)~ \textrm{for all}\ x, y, z\in G\big\}$$
and
$$B^2(G;A) = \big\{\theta:G \times G \to A  \, \mid \,  \textrm{there exists}~\lambda:G \to A~\textrm{with}~ \theta(x, y)=~ \big(x \cdot \lambda(y) \big)\lambda(xy)^{-1} \lambda(x)~\textrm{for all}~x, y\in G \big\}.$$
\para

Let  $G$ be  a group and $A$ a (left) $G$-module. We denote the action of $G$ on $A$ by $(g, a) \mapsto g \cdot a$ for $g \in G$ and $a \in A$. Let $E$ be an extension  of $G$  by $A$,  that is,  there is a short exact sequence of groups
$$ 1 \to A \to E \xrightarrow{\pi} G\to 1$$
such that the conjugation action of $G$  on $A$ coincides with the given action, where we identify $A$ with  its image in $E$. Let  $s : G  \to  E$  be  a normalised  \index{set-theoretic section}{set-theoretic section},  that  is,  $s(1)=1$ and $\pi \big(s(x)\big) = x$ for all $x  \in G$.  Then  there exists  a  map $\theta : G \times  G \to  A$ such  that 
\begin{equation}
  \label{equation theta}
  s(x) s(y)  s(x y)^{-1}= \theta(x, y)
\end{equation}
for  all $x,  y \in G$. Associativity in $E$ implies that $\theta$ satisfy
\begin{equation*}
\big( x \cdot {\theta(y,z)} \big) \, \theta(x,yz)=\theta(xy,z) \, \theta(x,y)
\end{equation*}
for all $x,y,z\in G$,  and hence $\theta$ is  a group $2$-cocycle. Conversely, if $\theta: G \times
G \to A$ is a group $2$-cocycle, then the set $E = A \times G$ can be endowed with the group operation
\begin{equation}
  \label{equation E operation}
  (a_1, x_1) (a_2, x_2)  =  \big(a_1(x_1 \cdot a_2) \theta(x_1, x_2), \, x_1 x_2 \big)
\end{equation}
for all $a_1, a_2 \in A$ and $x_1, x_2 \in G$. The preceding observations lead to the following well-known result \cite[Chapter IV,  Section  3]{MR0672956}. 

\begin{proposition} \label{prop:rob}
Let $G$ be a group and $A$  a $G$-module. Then the following assertions hold:
  \begin{enumerate}
  \item  Given an extension $E$ of $G$  by $A$ and a set-theoretic section $s$, the map $\theta$ as defined in \eqref{equation theta}  is a  group $2$-cocycle whose cohomology class  in $\Ho^2_{\rm Grp}(G; A)$ does not  depend on the choice of the set-theoretic section.
  \item   If  $\theta: G \times  G \to  A$ is  a group $2$-cocycle, then there  exists an extension $E$ of $G$  by $A$ and a set-theoretic section $s$ such
    that \eqref{equation theta} holds. 
  \item    The extension $E$ of $G$ by $A$ splits if and only if the cohomology class of $\theta$ in  $\Ho^2_{\rm Grp}(G; A)$ is trivial.
  \end{enumerate}
\end{proposition}

\begin{remark}
{\rm  If $A$ is  a trivial $G$-module, then the corresponding  group extension is central, that is, $A \le \Z(E)$. Note that, if the central extension $E$ of $G$ by $A$ splits,  then $E\cong G\times A$.}
\end{remark}

 Let $U$ and $V$  be groups. Let $A$ be an  abelian normal subgroup of $V$ and $\psi: U \to V/A$  a group homomorphism. Note that $A$ is a $V/A$-module  under  conjugation, and  hence  a  $U$-module via the homomorphism $\psi$. Let $C: U \to V$ be a lifting of $\psi$, that is, $C$ is a set-theoretic map such that $\psi(u) = C(u) A$ for all $u \in
U$. Since $\psi$ is a group homomorphism, there exists a map  $\kappa : U \times U \to A$ such that
\begin{equation}
  \label{eq:enter-coc}
  C(u_{1}) C(u_{2}) = \kappa(u_{1}, u_{2}) C(u_{1} u_{2})
\end{equation}
for all $u_{1}, u_{2} \in U$. Expanding the expression $C((u_{1} u_{2}) u_{3})  = C(u_{1} (u_{2} u_{3}))$ for $u_1, u_2, u_3 \in U$,  it follows  that   $\kappa$    is   a group $2$-cocycle. If $\sigma: U \to A$ is a map and $C'(u) = \sigma(u) C(u)$, then the group 2-cocycle $\kappa'$ associated to $C'$ is given by
\begin{equation*}
  \kappa'(u_{1}, u_{2})
  =
  \kappa(u_{1}, u_{2}) \sigma(u_{1})
  \big(u_1 \cdot \sigma(u_{2})  \big) \sigma(u_{1} u_{2})^{-1}
\end{equation*}
for all $u_{1}, u_{2} \in U$. Thus,  $\kappa$ and $\kappa'$ differ by a group $2$-coboundary.  Proposition~\ref{prop:rob} together with \cite[Chapter  IV, Section   3,  Exercise 1(a)]{MR0672956} yield the following result.

\begin{proposition}  \label{prop:MOF}
  Let $U$  and $V$ be groups, $A$ an abelian normal subgroup  of $V$ and   $\psi : U \to V/A$  a group homomorphism. Let $C: U \to V$ be a
  lifting of $\psi$   and $\kappa:U \times U \to A$ the map defined  by~\eqref{eq:enter-coc}. Then the following assertions hold:
  \begin{enumerate}
  \item    $\kappa$ is a group $2$-cocycle whose cohomology class in $\Ho^2_{\rm Grp}(U; A)$ does not depend on the choice of the lifting.
  \item     There exists a group homomorphism $\phi : U \to V$ which is a lifting of $\psi$ if and only if the cohomology class of $\kappa$ in $\Ho^2_{\rm Grp}(U; A)$ is trivial.
  \item    Two group homomorphisms $\phi_{1}, \phi_{2}$ are liftings of $\psi$ if  and only if there exists a group $1$-cocycle $\zeta : U \to A$ such that
$ \phi_{2}(u) =  \zeta(u) \phi_{1}(u)$  for all $u \in U$.
  \end{enumerate}
\end{proposition}

Specialising to the set-up of a skew left brace, Proposition \ref{prop:MOF} leads to the following result \cite[Theorem 3.2]{MR4531710}.

\begin{theorem}  \label{thm:thm}
Let $(G, \cdot, \circ)$ be a skew left brace such that the homomorphism $\lambda: (G, \circ) \to \Aut(G, \cdot)$ takes values in $\Inn(G, \cdot)$.  Let $C  :  G \to  G$  be a  map such  that  $\lambda(g) = \iota_{C(g)}$, the inner automorphism induced by $C(g)$, for all $g  \in G$. Let $\kappa :  (G, \circ)  \times (G, \circ)  \to \Z(G, \cdot)$  be  the map defined by
  \begin{equation}
    \label{equation coc}
    C(g)\cdot C(h)
    =
    \kappa(g, h) \cdot  C(g \circ h)
  \end{equation}
  for all $g, h \in G$.
Then the  following assertions hold:
  \begin{enumerate}
  \item    $\kappa$ is a group $2$-cocycle whose cohomology class in $\Ho^2_{\rm Grp} \big((G, \circ); \Z(G, \cdot) \big)$   does not  depend on the choice  of $C$.
  \item      The map $\lambda$ is induced by a Rota--Baxter operator on $(G, \cdot)$ if and only if  the cohomology class of $\kappa$ in $\Ho^2_{\rm Grp} \big((G, \circ); \Z(G, \cdot) \big)$ is trivial. 
  \item    Two Rota--Baxter operators $B_1, B_2$ on $(G, \cdot)$ induce the same map $\lambda$  if and only if there exists a group homomorphism $\zeta : (G, \circ)
    \to \Z(G, \cdot)$ such that $  B_{2}(g) = \zeta(g) \cdot B_{1}(g) $
    for all $g \in G$.
  \end{enumerate}
\end{theorem}

\begin{proof}
Let $(G, \cdot, \circ)$ be a skew left brace such that the homomorphism $\lambda: (G, \circ) \to \Aut(G, \cdot)$ takes values in $\Inn(G, \cdot)$. Then there exists a map $C  :  G \to  G$  such  that  $\lambda(g) = \iota_{C(g)}$ for all $g  \in G$. Composing  $\lambda$  with   the  natural
isomorphism  $\Inn(G, \cdot)  \to  (G,  \cdot)/\Z(G, \cdot)$, we  obtain a group homomorphism $\psi\colon (G, \circ) \to (G, \cdot)/\Z(G, \cdot)$. Note that a group homomorphism $B : (G, \circ) \to (G, \cdot)$ is a lifting  of $\psi$ if and only if $B$ is a Rota--Baxter operator on $(G, \cdot)$ which induces $\lambda$.  Note that $\Z(G, \cdot)$ is  a trivial $(G,\circ)$-module via $\psi$. Thus, a group 1-cocycle $(G, \circ)    \to \Z(G, \cdot)$  is precisely a group homomorphism. The theorem now follows directly from Proposition~\ref{prop:MOF}.  $\blacksquare$ 
\end{proof}

\begin{remark}
{\rm 
We now see how we can explicitly determine the Rota--Baxter operator associated with the given skew left brace $(G, \cdot, \circ)$. Let $C  :  G \to  G$  be a  map such  that $\lambda(g)=\iota_{C(g)}$ for each $g \in G$. Let $\kappa$ be the group 2-cocycle defined as in~\eqref{equation coc} and suppose that the cohomology class of $\kappa$ in $\Ho^2_{\rm Grp} \big((G, \circ); \Z(G, \cdot) \big)$ is trivial. Thus, there exists a map $\sigma:  (G, \circ) \to   \Z(G, \cdot)$  such  that  
$$  \kappa(g,h)=\sigma(g)^{-1}\cdot \sigma(h)^{-1}\cdot\sigma(g\circ h)$$
for all $g, h \in G$. If  we  define $B(g)= \sigma(g) \cdot C(g)  $  for   all   $g\in  G$,   then a direct check shows that $B$ is  the  Rota--Baxter operator on $(G, \cdot)$ which induces the map $\lambda$.}
\end{remark}

\begin{example}
{\rm 
We give an example of  a skew left  brace $(G,  \cdot, \circ)$ for which the map $\lambda$ takes values in the inner automorphism group of $(G, \cdot)$, but it is not induced by a Rota--Baxter operator. In fact, the example is a simplified version  of~\cite[Example 5.4]{MR4473766}. Let $p$ be an odd prime and $H$ the Heisenberg group of order
$p^{3}$ given by the presentation
$$
  H  = \big\langle u, v, k \mid u^{p}, v^{p}, k^{p},    [u, v] = k, [u, k], [v, k] \big\rangle.
$$
Note that every element of $H$ can be written uniquely as $u^i v^j k^{q}$ for some $0\le i,j,q<p$. Let $S = \langle x,  y \rangle$ be the elementary  abelian $p$-group  of  order  $p^2$ and take  $$(G,\cdot)  = S \times  H.$$  We write $\Z(G)=\Z(G,\cdot)$ and   $K =  \langle k \rangle=\Z(H) \le
\Z(G)$. Since $H/K = \langle u K,  v K \rangle$ is the elementary abelian $p$-group of order
$p^{2}$, the assignment
$$
 \psi (x K) = u K,  \quad   \psi(y K) = v K  \quad \textrm{and} \quad  \psi(h K) = K
$$
for all $h\in H$, uniquely defines an  endomorphism $\psi:G/K \to G/K$,
whose kernel and  image are $H/K$.
\para

Let $C : G \to G$ be a map such that 
\begin{equation}
  \label{equation of lift}
  C(g) K = \psi(g K)
\end{equation}
for all $g  \in  G$. Note that $C(G) \subseteq H$ and $C(H)\subseteq K$. It is not difficult to see that  the map $\lambda : G \to \Inn(G)$  defined by $\lambda_g = \iota_{C(g)}$ depends only on $\psi$ and is independent of the choice of $C$. Further, $\lambda$  satisfies 
$$\lambda \big(g\cdot \lambda_g( h) \big) = \lambda (g) \lambda (h)$$
for all $g, h \in G$. This defines a  skew left brace  $(G, \cdot, \circ)$,  where the group operation $\circ$ is given by
$$
  g \circ h  =  g\cdot C(g)\cdot h\cdot C(g)^{-1}
$$
for all $g,h\in G$. Note that $[a, b] = 1$  for all $a \in S$  and $b \in  C(G) \subseteq H$. Thus,  the group operations $\cdot$  and $\circ$ coincide on $S$, and hence $S$ is a subgroup  of $(G, \circ)$.
\para

We now compute the associated group $2$-cocycle $\kappa$. We can  choose any map $C$ satisfying~\eqref{equation of lift}. In particular, we take the map $C : G \to G$ defined by
\begin{equation*}
  C(x^{i} y^{j} c)
  =
  u^{i} v^{j}
  \in H
\end{equation*}
for all $0 \le i, j < p$ and $c \in H$. For all  $0 \le i, j, m, n < p$ and $c, d \in H$, we have
$$ C \big((x^{i} y^{j} c) \cdot  C(x^{i} y^{j} c) \cdot (x^{m} y^{n} d)  \cdot C(x^{i} y^{j} c)^{-1} \big)=	 C(x^{i+m} y^{j+n} e)= u^{i+m} v^{j+n},$$
for some $e \in H$.
We also have
$$  C(x^{i} y^{j} c) \cdot C(x^{m} y^{n} d)=  u^{i} v^{j} u^{m} v^{n}=  u^{i+m} v^{j+n} [v^{-j}, u^{-m}]= u^{i+m} v^{j+n} k^{-j m}.$$
Thus, the group $2$-cocycle $\kappa$ is given by
\begin{equation}
  \label{equation of kappa}
  \kappa( x^{i} y^{j} c, x^{m} y^{n} d)  =  k^{-j m}
\end{equation}
and has its image in $K = \Z(H) \le \Z(G) = S \times K$.
\para

We claim that the cohomology class of $\kappa$ in $\Ho^2_{\rm Grp} \big( (G, \circ); \Z(G, \cdot) \big)$ is non-trivial.   Let $E = (S \times K) \times (G, \circ)$. Consider the central
  extension  
  \begin{equation}
    \label{eq:alt}
    1 \to S \times K \to E \to (G, \circ) \to 1
  \end{equation}
  defined by $\kappa$  and the set-theoretic section given by $s(g) =  (1, g)$. The group operation in $E$  is given as in ~\eqref{equation E operation}. It  follows from the expression \eqref{equation of kappa} for   $\kappa$ that
  \begin{eqnarray*}
    \big[(1, x), (1, y)     \big]    &= &    (1, x) (1, y) \big((1, y) (1, x) \big)^{-1}\\
&=&    (1, x y) \big(\kappa(x, y), 1\big) \big( (\kappa(y, x), 1) (1, y x) \big)^{-1}\\
&=&     (1, x y) (1, y x)^{-1} \big(\kappa(y, x), 1\big)^{-1}\\
&=&    \big(\kappa(y, x)^{-1}, 1\big)\\
&=&     (k, 1),
  \end{eqnarray*}
  so that $(1 \times K, 1)$ is contained in the  commutator subgroup    $[E, E]$ of $E$.  Suppose  that the exact sequence~\eqref{eq:alt}
  splits, and let $T$ be a complement  to $S \times K$ in $E$. Then $M  = (S  \times 1, 1)  T$ is a maximal  subgroup of
  the  finite  $p$-group $E$,  which  does  not contain  $(1  \times K, 1)$, a non-trivial subgroup of $[E, E]$.  This contradiction shows that the sequence \eqref{eq:alt} does
  not split. Hence, by Proposition~\ref{prop:rob}, the cohomology class of $\kappa$ in $\Ho^2_{\rm Grp} \big((G, \circ); \Z(G) \big)$ is  non-trivial. }
\end{example}


\chapter*{Part II}
\begin{center}
{\huge{Racks and quandles}}
\end{center}

\chapter{Racks and quandles}\label{chapter preliminaries on racks and quandles}

\begin{quote}
Racks and quandles offer non-degenerate bijective solutions to the Yang--Baxter equation, and share a close connection with low-dimensional topology. This chapter develops the fundamental algebraic theory of racks and quandles, building upon the pioneering works of Joyce and Matveev. We showcase a range of examples of quandles, which extend across multiple mathematical disciplines, and investigate key constructions and properties that are intrinsic to racks and quandles.
\end{quote}
\bigskip

\section{Properties and examples}
A quandle is a non-empty set with a binary operation that satisfies three axioms which are algebraic formulations of the three Reidemeister moves of planar diagrams of knots in the 3-space. Ignoring the first Reidemeister move gives rise to a weaker structure called a rack. While racks and quandles are abundant in the vast realm of mathematics, their study gained momentum only after fundamental works of Joyce \cite{MR2628474, MR0638121} and Matveev \cite{MR0672410}, who independently associated a quandle to each link in the 3-space and proved it to be an almost complete invariant of knots. Although knot quandles are strong invariants of knots, it is usually difficult to distinguish two quandles from their presentations. The algebraic theory of racks and quandles has seen numerous captivating developments as a result of this.

\para

We begin with the definition of a quandle.

\begin{definition}
A \index{quandle}{\it quandle} is a non-empty set $X$ with a binary operation $(x,y) \mapsto x *y$ satisfying the following axioms:
\begin{enumerate}
\item[(Q1)] $x*x=x$ for all $x \in X$.
\item[(Q2)] For each $x, y \in X$, there is a unique $z \in X$ such that $x=z*y$.
\item[(Q3)] $(x*y)*z=(x*z) * (y*z)$ for all $x, y, z \in X$.
\end{enumerate}
An algebraic structure satisfying (Q2) and (Q3) is called a \index{rack}{\it rack}, whereas the one satisfying (Q3) is called a \index{shelf} {\it shelf} or a \index{distributive groupoid}{\it distributive groupoid}.
\end{definition}

\begin{remark}{\rm 
It is worth noting that Burstin and Mayer \cite{MR1581179} initially investigated distinct classes of finite distributive groupoids. Further, given a shelf $(X, *)$, the map $r: X \times X \to X\times X$ given by $r(x, y)=(y, x*y)$ gives a solution to the Yang--Baxter equation. Additionally, if $(X, *)$ is a rack, then the solution $(X, r)$ is bijective and non-degenerate. It was already observed by Brieskorn \cite{MR0975077} that racks provide solutions to the Yang--Baxter equation.}
\end{remark}

If $(X, *)$ is a rack, then it follows from the second quandle axiom that there exists another binary operation $*^{-1} : X \times X \to X$ such that
$$
z = x *^{-1} y  \quad \textrm{if and only if} \quad x = z * y  
$$
for all $x, y, z \in X$. This equivalent to
$$
(x * y) *^{-1} y = x = (x *^{-1} y) * y
$$
for all $x, y \in X$. The preceding identity together with the third quandle axiom implies that
\begin{eqnarray*}
(x * y) *^{-1} z &=&  (x *^{-1} z) * (y *^{-1} z),\\
(x *^{-1} y) * z &=& (x * z) *^{-1} (y * z) \quad \textrm{and}\\
(x *^{-1} y) *^{-1} z &=& (x *^{-1} z) *^{-1} (y *^{-1} z)
\end{eqnarray*}
for all $x, y, z \in X$.
\para 

If $(X,*)$ and $(Y, \circ)$ are two quandles, then a map $f:X \to Y$ is called a \index{quandle homomorphism}{\it quandle homomorphism} if
$$f(x*y)=f(x) \circ f(y)$$
for all $x, y \in X$. The set of all quandles forms a category $\mathcal{Q}$ with objects as quandles and quandle homomorphisms as morphisms in this category. We write $\Hom_{\mathcal{Q}} (X, Y)$ to denote the set of all quandle homomorphisms from $X$ to $Y$. As usual, a quandle isomorphism is a simply a bijective homomorphism of quandles. A non-empty subset $Y$ of a quandle $(X, *)$ is called a \index{subquandle}{\it subquandle} of $X$ if $Y$ forms a quandle under the operation $*$. The essential axiom to be checked in $Y$ is the second quandle axiom. Morphisms and \index{subrack} subracks of racks are defined analogously. Throughout the text, by \index{order}{\it order} of a rack or a quandle $X$, we mean the cardinality of the set $X$.
\para

If $X$ is a rack, then for each $x \in X$, we can define a map $S_x : X \longrightarrow X$ by $$S_x(y)= y*x$$ for each $y \in X$.
\para

\begin{proposition}
If $X$ is a rack, then $S_x$ is an automorphism of $X$ for each $x \in X$.
\end{proposition}

\begin{proof}
Note that the third quandle axiom can be rewritten as $S_z(x*y)=S_z(x)*S_z(y)$ for all $x, y, z \in X$. hence, each $S_x$ is a homomorphism of $X$. Further, the second quandle axiom is equivalent to each $S_x$ being a bijection. $\blacksquare$
\end{proof}

The map $S_x$ is referred to as the \index{inner automorphism}{\it inner automorphism} of $X$ induced by $x$, or  the {\it symmetry} at $x$. Let $\Aut(X)$ denote the group of all automorphisms of $X$. Then the subgroup $\Inn(X)$ of $\Aut(X)$ generated by the set $\{S_x \mid x \in X \}$ is called the group of \index{inner automorphism group}{inner automorphisms} of $X$. It is evident that $\Inn(X)$ is a normal subgroup of $\Aut(X)$.
\para 

\begin{proposition}\label{S quandle homomorphism}
If $X$ is a quandle, then the map $S: X \to \Conj(\Inn(X))$ given by $S(x)=S_x$ is a quandle homomorphism.
\end{proposition}

\begin{proof}
For $x, y, z \in X$, we have
\begin{eqnarray*}
S_{x*y}(z) &=& z * (x*y)\\
&=& \big((z*^{-1} y)*y \big)* (x*y)\\
&=& \big(((z*^{-1} y)*x)* y \big)\\
&=&  S_yS_x S_y^{-1} (z)\\
&=&  (S_x *S_y) (z),
\end{eqnarray*}
and hence $S$ is a quandle homomorphism. $\blacksquare$
\end{proof}
\para

A quandle $X$ is called \index{trivial quandle}{\it trivial} if $x*y=x$ for all $x,y \in X$. Trivial quandles are the simplest example of quandles. Note that a trivial quandle can contain arbitrary number of elements, and that any quandle on at most two elements is trivial. We denote the trivial quandle with $n$ elements by $\T_n$ and the trivial quandle with countably infinite number of elements by $\T_{\infty}$.
\medskip

It is interesting to understand similarities between quandles and groups.

\begin{proposition}
The following statements hold:
\begin{enumerate}
\item A quandle is associative if and only if it is trivial.
\item A quandle contains an identity element if and only if it is the trivial quandle with one element.
\end{enumerate}
\end{proposition}

\begin{proof}
Suppose that
$$
(x * y)  * z = x * (y * z).
$$
for all $x, y, z \in X$. Taking $z = y$, we get $(x * y) * y = x * y$. It follows from the second quandle axiom that $x*y=x$, and hence $X$ is a trivial quandle. The converse is obvious, and hence assertion (1) holds.
\para
If $X$ contains an identity element, say $e$, then $e * x = x=x*x$ for all $x \in X$. It now follows from the second quandle axiom that $x = e$. The converse is obvious, and hence we have assertion (2). $\blacksquare$
\end{proof}

Let us now see some non-trivial examples of quandles arising in different contexts. Groups turn out to be a large source of quandles. We shall often denote a finite quandle simply by its multiplication table. For this, we first order the elements of $X$ as $\{x_1, \ldots, x_n \}$. Then the $(i,j)$-th entry in the table represents the element $x_i * x_j$ of the quandle.

\begin{example}{\rm 
If $G$ is a group and $n$ an integer, then the set $G$ equipped with the binary operation
$$
a*b= b^{n} a b^{-n}
$$ for $a, b \in G$, becomes a quandle called the {\it $n$-conjugation quandle} of $G$ and is denoted by $\Conj_n(G)$. If $n=1$, then we  refer this quandle simply as the \index{conjugation quandle}{\it conjugation quandle} and denote it by $\Conj(G)$. Obviously, if $G$ is an abelian group, then $\Conj(G)$ is a trivial quandle.
\para
 For instance, if
$$
\Sigma_3 = \big\{ e, (12), (23), (13), (123), (132) \big\}
$$
is the symmetric group on three symbols, then the  conjugation quandle $\Conj(\Sigma_3)$ has the following multiplication table.
$$
\begin{tabular}{|c||c|c|c|c|c|c|}
    \hline
$*$ & $e$ & (12) & (23) & (13) & (123) & (132) \\
  \hline \hline
$e$ & $e$ & $e$ & $e$ & $e$ & $e$ & $e$ \\
\hline
(12) & (12) & (12) & (13) & (23) & (23) & (13) \\
\hline
(23) & (23) & (13) & (23) & (12) & (13) & (12) \\
\hline
(13) & (13) & (23) & (12) & (13) & (12) & (23) \\
\hline
(123) & (123) & (132) & (132) & (132) & (123) & (123) \\
\hline
(132) & (132) & (123) & (123) & (123) & (132) & (132) \\
  \hline
\end{tabular}
$$

The table shows that $\Conj(\Sigma_3)$ contains the subquandle  $\{ e, (12), (23), (13) \}$
with 4 elements. Hence, we do not have an analog of the Lagrange's theorem for finite quandles.}
\end{example}

\begin{example}{\rm 
If $G$ is a group, then the set $G$ equipped with the binary operation
$$
a*b= b a^{-1} b
$$ forms a quandle called the \index{core quandle}{\it core quandle} of $G$, which we denote by $\Core(G)$. For example, $\Core(\Sigma_3)$ has the following multiplication table.
$$
\begin{tabular}{|c||c|c|c|c|c|c|}
    \hline
$*$ & $e$ & (12) & (23) & (13) & (123) & (132) \\
  \hline \hline
$e$ & $e$ & $e$ & $e$ & $e$ & (132) & (123) \\
\hline
(12) & (12) & (12) & (13) & (23) & (12) & (12) \\
\hline
(23) & (23) & (13) & (23) & (12) & (23) & (23) \\
\hline
(13) & (13) & (23) & (12) & (13) & (13) & (13) \\
\hline
(123) & (132) & (123) & (123) & (123) & (123) & e \\
\hline
(132) & (123) & (132) & (132) & (132) & e & (132) \\
  \hline
\end{tabular}
$$}
\end{example}

\begin{example}{\rm 
If $G$ is an  abelian group, then the core quandle of $G$ is referred in the literature as the \index{Takasaki quandle}{\it Takasaki quandle} of $G$, and is denoted by $\T(G)$. Note that, for an additively written abelian group, the binary operation in $\T(G)$ is given as $$a*b= 2b-a.$$ Such quandles (earlier called \index{Kei}{\it Kei}) appeared first in the work of Takasaki \cite{MR0021002} in the context of finite geometry. 
\para

The Takasaki quandle of the cyclic group $\mathbb{Z}_n$ of integers modulo $n$ is usually referred to as the \index{dihedral quandle}{\it dihedral quandle} and is denoted by $\R_n$. We shall write $\R_n= \{0, 1, 2, \ldots, n-1 \}$, where $$i*j= 2j-i \mod n.$$ Interestingly, the dihedral quandle $\R_n$ has a nice geometrical interpretation (see Figure \ref{regular pentagon dihedral} for $\R_5$ and $\R_6$). Consider a regular $n$-gon and  denote its vertices as $a_0, a_1, \ldots, a_{n-1}$, along the counterclockwise direction.  The binary operation in $\R_n$ can be viewed as the reflection of this  $n$-gon. More precisely, $a_{i * j}$ is the vertex of the $n$-gon which one gets after reflecting $a_i$ across the line which goes through $a_j$ and is an axis of symmetry of the $n$-gon.

\begin{figure}[hbtp]
\centering
\includegraphics[height=6cm]{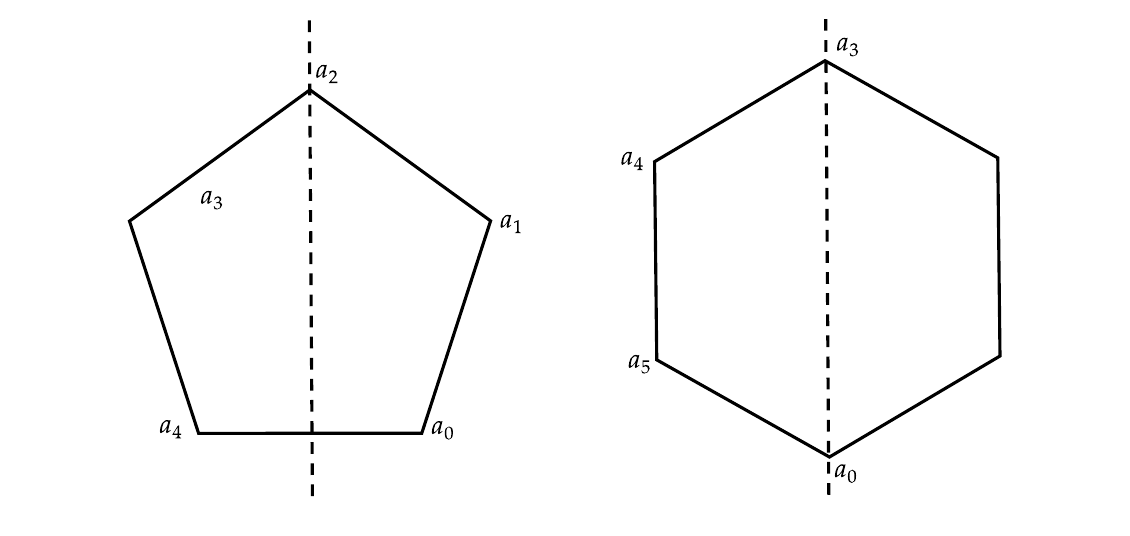}
\caption{Geometric interpretation of binary operation in dihedral quandles.}
\label{regular pentagon dihedral}
\end{figure}}
\end{example}

\begin{example}\label{example generalised Alexander quandle}
{\rm 
If $G$ is a group and $\varphi \in \Aut(G)$ an automorphism of $G$, then the binary operation $$a*b=\varphi(ab^{-1})b$$ gives a quandle structure on $G$, which we denote by $\Alex(G, \varphi)$. This quandle is referred to as the \index{generalised Alexander quandle}{\it generalised Alexander quandle}.
\para
If $G$ is an additive abelian group and  $\varphi=-\id$, the inversion in $G$, then $a*b=2b-a$. Thus,  $\Alex(G, \varphi)=\T(G)$, the Takasaki quandle of $G$. On the other hand, for arbitrary $\varphi \in \Aut(G)$, we have $$a*b=\varphi(a) + (\id-\varphi)(b),$$ and  $\Alex(G, \varphi)$ is the \index{Alexander quandle}{\it Alexander quandle} of $G$.
\para

If $M$ is a module over the ring $\mathbb {Z} \left[t,\,t^{-1}\right]$ of Laurent polynomials over integers, then the right multiplication by $t$ is an automorphism of the additive abelian group of $M$. This gives the Alexander quandle structure on $M$ with
    $$
    a* b=ta+(1-t)b.
    $$}
\end{example}

\begin{example}\label{example spherical quandle}
{\rm 
Let $\mathbb{k}$ be a  field of characteristic other than 2 and $\langle ~,~ \rangle : \mathbb{k}^{n+1} \times \mathbb{k}^{n+1} \to \mathbb{k}$ the standard symmetric bilinear form. Let
$$
S(\mathbb{k}^{n+1}) = \big\{ x \in \mathbb{k}^{n+1} \, \mid \,  \langle x, x \rangle = 1 \big\}
$$
be the unit sphere in $\mathbb{k}^{n+1}$ with respect to the form $\langle ~,~ \rangle$. Then the binary operation 
$$x * y = 2 \langle x, y \rangle y -x$$
gives a quandle structure on $S(\mathbb{k}^{n+1})$, and is referred to as the \index{spherical quandle}{\it spherical quandle}.}
\end{example}

\begin{example}\label{coset quandle with auto}
{\rm 
Let $G$ be a group and $\phi$ an automorphism of $G$. Let $H$ be a subgroup of $G$ such that $H \le \Fix(\phi)= \{x \in G \mid \phi(x)=x \}$. Then the set $G/H$ of right cosets of $H$ in $G$ has a quandle structure given by 
$$Hx*Hy= H\phi(xy^{-1})y$$
for $x, y \in G$. We denote this quandle by $\Alex(G, H, \varphi)$. Taking $H$ to be the trivial subgroup recovers Example \ref{example generalised Alexander quandle}. Further, by taking $G=\Oo(n + 1, \mathbb{F})$, $H = \Oo(n,\mathbb{F})$ viewed as a subgroup of $G$ and $\phi$ to be the inner automorphism induced by the diagonal matrix $\diag(1, \ldots, 1, -1)$, we recover Example \ref{example spherical quandle}.}
\end{example} 

\begin{example}{\rm 
Let $\mathbb{k}$ be any commutative ring and $X$ a free $\mathbb{k}$-module equipped with an antisymmetric bilinear form $\langle-, - \rangle : X \times X \to \mathbb{k}$. If the additive group of the ring $\mathbb{k}$ has any two-torsion, then we require alternating rather than just antisymmetric. Then $X$ can be turned into a quandle when equipped with the binary operation
$$x * y = x + \langle x, y \rangle y
$$
for all $x, y \in X$. This quandle first appeared in \cite[Example 4]{MR1985908} where it is referred to as an \index{alternating quandle}{\it alternating quandle}. Later, from the work \cite{MR2493966} onwards, such quandles are referred to as \index{symplectic quandle}{\it symplectic quandles}.
\para

Let us consider special symplectic quandles arising from orientable surfaces. Let $\mathcal{S}_g$ be a closed orientable surface of genus $g \ge 1$. For isotopy classes $a$ and $b$ of transverse, oriented, simple closed curves in $\mathcal{S}_g$, the {\it algebraic intersection number} $\hat{i}(a,b)$ is defined as the sum of the indices of the intersection points of $a$ and $b$, where an intersection point is of index +1 when the orientation of the intersection agrees with the orientation of the surface, and is -1 otherwise. It turns out that $\hat{i}(a,b)$ depends only on the homology classes (and hence isotopy classes) of $a$ and $b$. Further, it is known \cite[Section 6.1.2]{MR2850125} that $\hat{i}(-, -)$ gives a skew-symmetric (in fact, symplectic) bilinear form on the $\mathbb{Z}$-module $\mathrm{H}_1^{\rm Cell}(\mathcal{S}_g; \mathbb{Z}) \cong \mathbb{Z}^{2g}$, where $\mathrm{H}_1^{\rm Cell}(\mathcal{S}_g; \mathbb{Z})$ is the cellular homology of $\mathcal{S}_g$. Then the binary operation given by
\begin{equation}
x*y=x + \hat{i}(x,y) y
\end{equation}
for $x, y \in \mathrm{H}_1^{\rm Cell}(\mathcal{S}_g; \mathbb{Z})$, gives a symplectic quandle structure on $\mathrm{H}_1^{\rm Cell}(\mathcal{S}_g; \mathbb{Z})$.}
\end{example}

The next example shows that quandles appear naturally in algebraic geometry as well. We refer the reader to \cite{MR0463157} for the necessary background on algebraic geometry.

\begin{example}{\rm 
Let $\Bbbk$ be an algebraically closed field of an arbitrary characteristic. A \index{quandle variety}{\it quandle variety} over $\Bbbk$ is a reduced and irreducible algebraic variety $X$ over $\Bbbk$ equipped with a quandle operation $*$ such that the map $X \times X \to X \times X$ given by 
$$(x, y) \mapsto (x, x *y)$$ is an automorphism of varieties. 
\para
Let $X$ be an algebraic variety, $A$ a connected commutative \index{algebraic group}{\it algebraic group} and $f: X \times X \to A$ a regular map with  $f(x, x) = 0$ for all $x$.  Then the set $X\times A$ has a quandle variety structure with quandle operation given by
$$(x, a) * (y, b) = \big(x, \, a + f(x, y)\big).$$ 
\para

If $G$ is a connected algebraic group, then $\Conj(G)$ is a quandle variety. Similarly, if $G$ is a connected algebraic group and $\varphi$ is an algebraic automorphism, then $\Alex(G, \varphi)$  is a quandle variety. Quandle varieties were introduced and investigated by Takahashi in \cite{MR3492048}.}
\end{example}

\begin{example}{\rm 
Example \ref{coset quandle with auto} arise naturally in the setting of Lie groups. A \index{differentiable reflection space}{\it differentiable reflection space} is a differentiable manifold with a differentiable binary operation satisfying the quandle axioms \cite{MR0239005}. For example, if $G$ is a connected Lie group with an involutive automorphism $\varphi$ and $H$ a subgroup of the fixed-point subgroup containing the identity component, then the homogeneous space $G/H$ acquires the structure of differentiable reflection space with the binary operation $$Hx *Hy=H\varphi(x y^{-1})x$$ for all $x, y \in G$. Conversely, any differentiable reflexion space can be obtained from a symmetric space \cite{MR0239005}. }
\end{example}

Examples of the preceding kind have led to transfer of ideas (for example, flatness) from the theory of Riemannian symmetric spaces to the theory of quandles \cite{MR3544543, MR3576755}. Next, we consider some quandles arising from surfaces.

\begin{example}{\rm 
Let $\mathcal{S}_g$ be a closed orientable surface of genus $g \ge 1$, $\mathcal{M}_g$ the mapping class group of $\mathcal{S}_g$ and $\mathcal{D}_g$ the set of isotopy classes of simple closed curves in $\mathcal{S}_g$. It is well-known that $\mathcal{M}_g$ is generated by Dehn twists along finitely many simple closed curves from $\mathcal{D}_g$ \cite[Theorem 4.1]{MR2850125}. The binary operation
$$\alpha * \beta= T_\beta(\alpha)$$ for $\alpha, \beta \in \mathcal{D}_g$, where $T_\beta$ is the Dehn twist along $\beta$ (See Figure \ref{dehn quandle operation}), turns $\mathcal{D}_g$ into a quandle called the \index{Dehn quandle of surface}{\it Dehn quandle} of the surface $\mathcal{S}_g$.

\begin{figure}[hbtp]
\centering
\includegraphics[height=3.8cm]{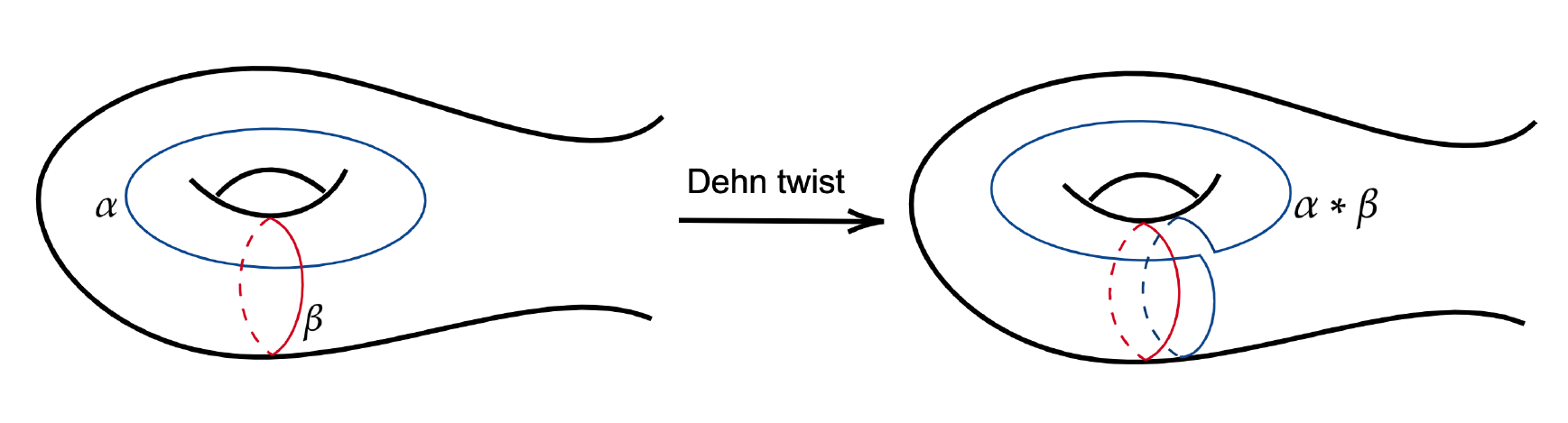}
\caption{Dehn twist.}
\label{dehn quandle operation}
\end{figure}

By identifying the isotopy class of a simple closed curve with the corresponding Dehn twist, it turns out that the quandle $\mathcal{D}_g$ is a subquandle of the conjugation quandle $\Conj(\mathcal{M}_g)$ of the mapping class group. These quandles originally appeared in the work of Zablow \cite{MR2699808, MR1967241}. Further, \cite{MR1865704,Yetter2002, MR1985908} considered a quandle structure on the set of isotopy classes of simple closed arcs in orientable surfaces with at least two punctures, and called it quandle of cords. A presentation for the quandle of cords of the disk and the 2-sphere has been given in \cite{MR1865704}.}
\end{example}

More generally, we have the following example.

\begin{example}\label{example Dehn quandle}
{\rm 
Let $G$ be a group, $A$ a non-empty subset of $G$ and $A^G$ the set of all conjugates of elements of $A$ in $G$. The \index{Dehn quandle}{\it Dehn quandle} $\mathcal{D}(A^G)$ of $G$ with respect to $A$ is defined as the set $A^G$ equipped with the binary operation of conjugation, that is, $$x*y=yxy^{-1}$$ for all $x, y \in A^G$.
\para

Clearly, $\mathcal{D}(A^G)$ is a subquandle of $\Conj(G)$ for each subset $A$ of $G$. We shall prove later in Theorem \ref{injectivity-dehn-eta} that every subquandle of $\Conj(G)$ is a Dehn quandle. Dehn quandles of groups generalise many well-known constructions of quandles from groups. 
\begin{enumerate}
\item If $A=G$, then $\mathcal{D}(G^G)$ is the conjugation quandle $\Conj(G)$. 
\item If $F(S)$ is the free group on the set $S$, then $\mathcal{D}(S^{F(S)})$ is the \index{free quandle}free quandle on $S$ \cite{Kamada2012, MR3588325, MR3729413}. 
\item If $\mathcal{W}$ is a \index{Coxeter group}Coxeter group with Coxeter generating set $S$, then $\mathcal{D}(S^{\mathcal{W}})$ is the so called \index{Coxeter quandle}Coxeter quandle \cite{MR4175808, MR3821082}. Similarly, if $\mathcal{A}$ is an \index{Artin group}Artin group with Artin generating set $S$, then $\mathcal{D}(S^{\mathcal{A}})$ is the \index{Artin quandle}Artin quandle.  
\item Let $S_{g}$ be a closed orientable surface of genus $g$ and $\mathcal{M}_{g}$ its mapping class group. If $S$ is the set of Dehn twists about essential simple closed curves, then $\mathcal{D}(S^{\mathcal{M}_{g}})$ is the Dehn quandle of the surface $\mathcal{S}_g$. 
\end{enumerate}}
\end{example}

Dehn quandles have been investigated in detail in \cite{MR4669143}. Further, in \cite{MR4642200}, two approaches to write efficient presentations for such quandles have been given. Next, we give examples of quandles arising from knots and links in the 3-sphere.  Numerous captivating texts on knot theory are readily accessible, such as \cite{MR0808776, MR0146828, MR0907872, MR1417494, MR1472978, MR1253070, MR1391727, MR0515288}.

\begin{example}\label{construction link quandle}
{\rm 
Let $L$ be a link in the 3-sphere. Joyce \cite{MR2628474, MR0638121} and Matveev \cite{MR0672410} associated a quandle $Q(L)$ to $L$ called the \index{link quandle}{\it link quandle} of the link $L$. We fix a diagram $D(L)$ of $L$ and label its arcs. Then the link quandle $Q(L)$ is the quandle generated by labelings of arcs of $D(L)$ with a defining relation at each crossing in $D(L)$ given as shown in Figure \ref{Link quandle relation}.

\begin{figure}[hbtp]
\centering
\includegraphics[height=3cm]{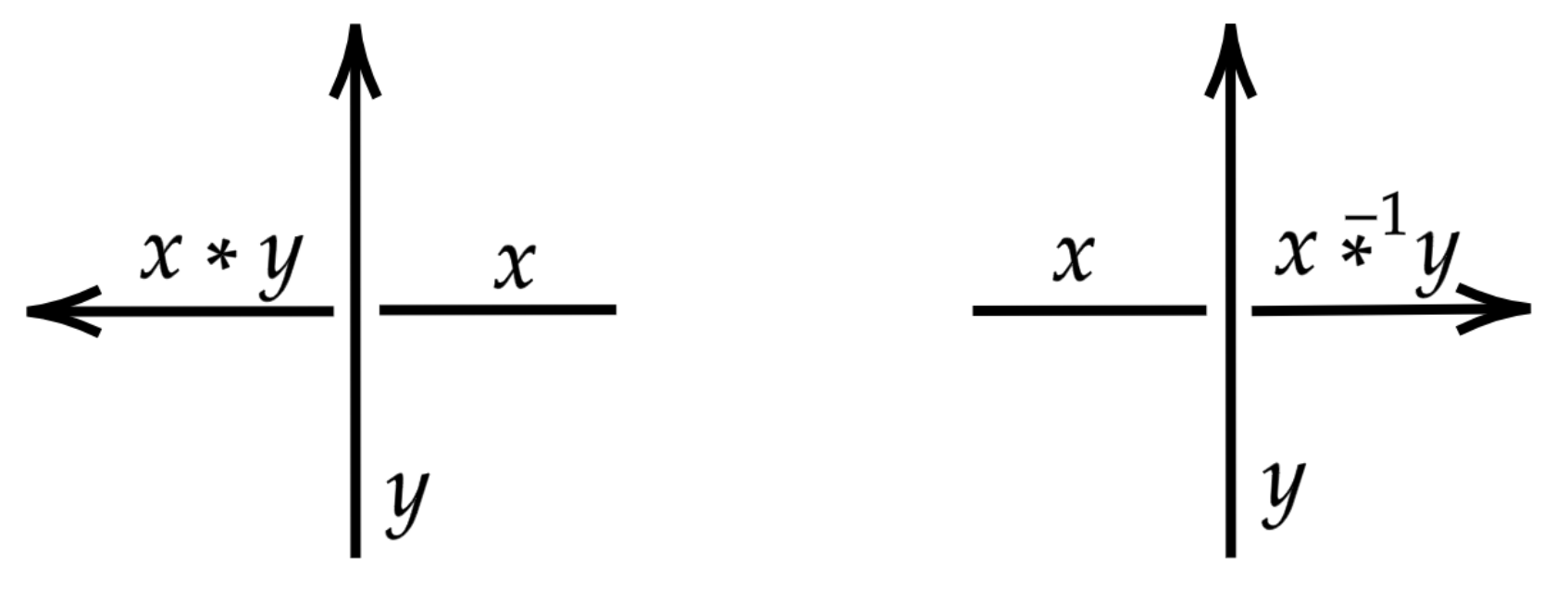}
\caption{Link quandle relations.}
\label{Link quandle relation}
\end{figure}

One can check that the three quandle axioms are equivalent to the three Reidemeister moves of planar diagrams of links as shown in Figure \ref{reidemeister moves and quandle axioms}.

\begin{figure}[hbtp]
\centering
\includegraphics[height=4.4cm]{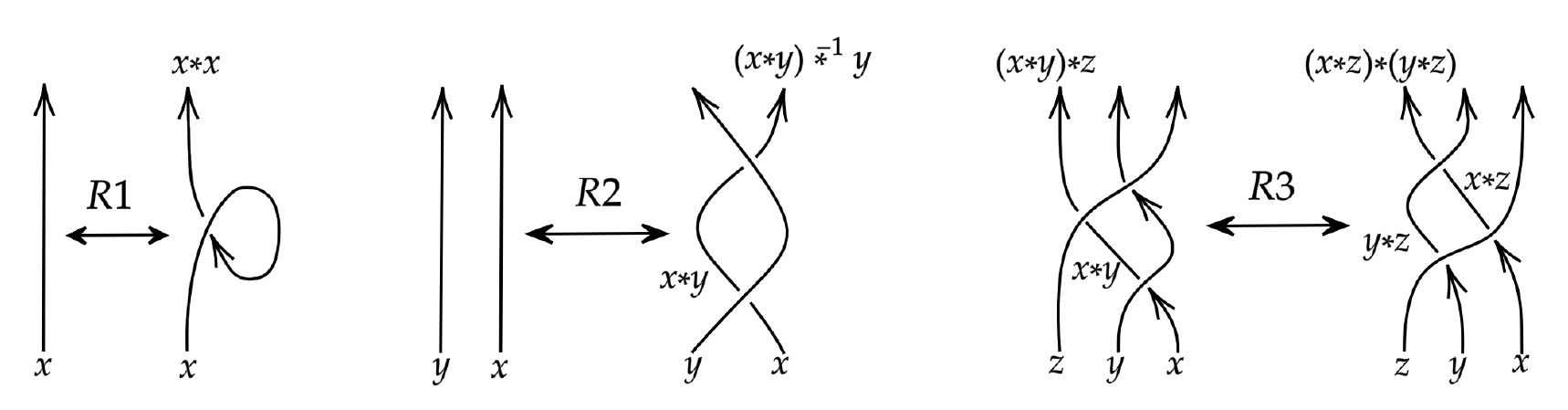}
\caption{Reidemeister moves and quandle axioms.}
\label{reidemeister moves and quandle axioms}
\end{figure}

It turns out that $Q(L)$ is independent of the diagram of $L$, that is, the quandles obtained from any two diagrams of $L$ are isomorphic. We refer the reader to the original sources 
\cite{MR2628474, MR0638121, MR0672410} for details. In a rather curious connection of knot quandle of knots with Dehn quandles of surfaces, Niebrzydowski and Przytycki \cite{MR2583322} proved that the Dehn quandle of the torus is isomorphic to the knot quandle of the trefoil knot.}
\end{example}

Using homomorphisms of quandles we can prove that the trefoil knot is non-trivial.

\begin{example}\label{trefoil-coloring-r3}
{\rm
Let $K$ be the trefoil knot (see Figure \ref{Trefoil knot}).  Then its knot quandle has a presentation
$$
Q(K) = \big\langle x, y, z   \, \mid \,  x * y = z, ~y * z = x,~ z * x = y \big\rangle.
$$

\begin{figure}[hbtp]
\centering
\includegraphics[height=4cm]{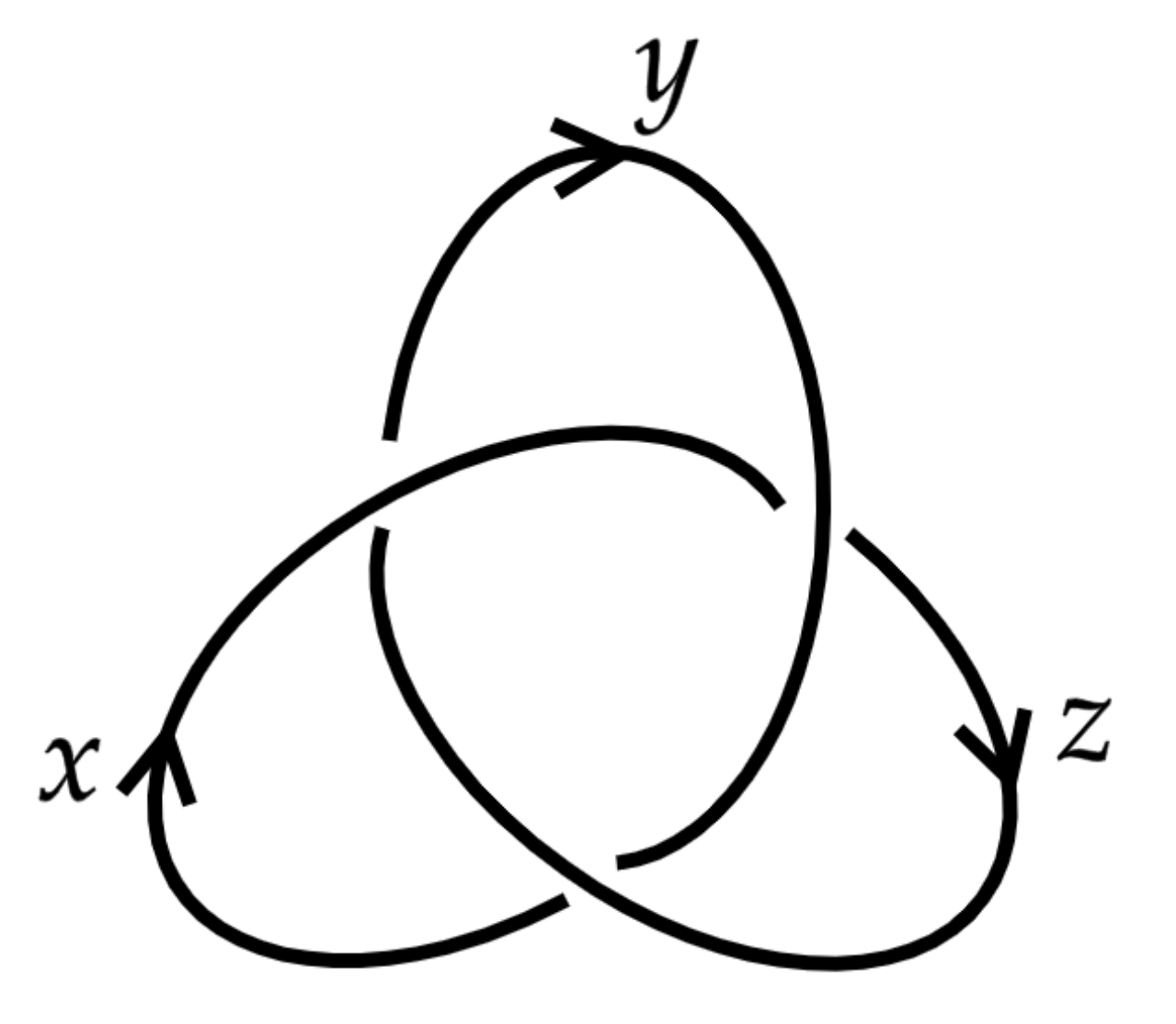}
\caption{Trefoil knot.}
\label{Trefoil knot}
\end{figure}

We see that the map $x \mapsto 0,~~y \mapsto 1,~~z \mapsto 2$ defines a quandle epimorphism from the knot quandle $Q(K)$ to the dihedral quandle $\R_3 = \{ 0, 1, 2 \}$. Since the fundamental knot of the trivial knot is the 1-element trivial quandle, it follows that $K$ is a non-trivial knot.}
\end{example}

We introduce a class of quandles in analogy to Fibonacci groups.

\begin{definition}
Let $n \geq 2$ be an integer. A \index{Fibonacci  quandle} {\it Fibonacci  quandle} $\Fib(n)$ is a quandle with a presentation
$$
\Fib(n)=\big\langle x_0, x_1, \ldots, x_{n-1} \, \mid \, x_0 * x_1 = x_2, \ldots, x_{n-2} * x_{n-1} = x_0, ~~x_{n-1} * x_0 = x_1 \big\rangle.
$$
\end{definition}

It follows from Example \ref{trefoil-coloring-r3} that the knot quandle of the trefoil knot is isomorphic to $\Fib(3)$. Thus, the following question seems natural.

\begin{question}
Which Fibonacci quandles can be realised as link quandles? More generally, which Fibonacci quandles admit non-trivial colorings by link quandles?
\end{question}

\begin{example}{\rm 
Let $\mathbb{k}$ be any ring and $r \in \mathbb{k}$ a non-zero element. Then the binary operation given by $$x * y = x+r$$ for $x, y \in \mathbb{k}$, gives a rack structure on $\mathbb{k}$ which is not a quandle.}
\end{example}

\begin{example}{\rm 
Let $X$ be any set and $f \in \Sigma_X$ a permutation of $X$. Then the binary operation given by  $x*y=f(x)$ gives a rack structure on $X$,  called a \index{permutation rack }{\it permutation rack}. A permutation rack is a quandle if and only if $f=\id_X$.}
\end{example}

\begin{exercise}{\rm 
Check that the quandle axioms are satisfied in each of the preceding examples.}
\end{exercise}

\begin{exercise}\label{all three element quandles}
{\rm 
Prove that the 3-element set $\{ 1, 2, 3 \}$ has precisely three non-isomorphic quandle structures with multiplication tables  as follows:
$$
\T_3=\begin{tabular}{|c||c|c|c|}
    \hline
$*$ & 1 & 2 & 3  \\
  \hline \hline
1 & 1 & 1 & 1  \\
\hline
2 & 2 & 2 & 2  \\
\hline
3 & 3 & 3 & 3  \\
  \hline
\end{tabular}, 
\quad
\R_3=\begin{tabular}{|c||c|c|c|}
    \hline
$*$ & 1 & 2 & 3  \\
  \hline \hline
1 & 1 & 3 & 2  \\
\hline
2 & 3 & 2 & 1  \\
\hline
3 & 2 & 1 & 3  \\
  \hline
\end{tabular},
\quad
\J_3=\begin{tabular}{|c||c|c|c|}
    \hline
$*$ & 1 & 2 & 3  \\
  \hline \hline
1 & 1 & 1 & 1  \\
\hline
2 & 3 & 2 & 2  \\
\hline
3 & 2 & 3 & 3  \\
  \hline
\end{tabular},
$$
where $\T_3$ is the trivial quandle, $\R_3$ is the dihedral quandle and $\J_3$ is a union of two trivial quandles. See Example \ref{union quandle example} for the general definition of union of two quandles.}
\end{exercise}

\begin{remark}{\rm 
Classification of finite quandles is a major theme in the subject. It is known that there are 7 non-isomorphic 4-element quandles, 22 non-isomorphic 5-element quandles, 73 non-isomorphic 6-element quandles and 298 non-isomorphic 7-element quandles. Table \ref{number of racks and quandles}
borrowed from \cite{MR3957904} lists the number $r(n)$ of non-isomorphic racks and the number $q(n)$ of non-isomorphic quandles of size $n \le 13$.  We refer the reader to \cite{MR3665565, MR2414453, MR2900878, MR3543136, MR3685034, MR2926571, MR3957904} and references therein for more related works.}
\end{remark}

\begin{table}[H]
\begin{center}
\begin{tabular}{|c|c|c|c|c|c|c|c|c|c|c|c|c|c|c|c|}
    \hline
$n$ & 1 & 2 & 3& 4 & 5 & 6 & 7 & 8 & 9 & 10 & 11 &  12 & 13 \\
\hline
$q(n)$ & 1 & 1 & 3 & 7 & 22 & 73 & 298 & 1581 & 11079 & 102771 & 1275419 & 21101335 & 469250886 \\
\hline
\hline
$n$ & 1 & 2 & 3 & 4 & 5 & 6 & 7 & 8 & 9 & 10 & 11 & 12 & 13\\
\hline
$r(n)$ & 1 & 2 & 6 & 19 & 74 & 353 & 2080 & 16023 & 159526 & 2093244 & 36265070 & 836395102 & 25794670618 \\
  \hline
\end{tabular}
\end{center}
\caption{The number of non-isomorphic racks and quandles.} \label{number of racks and quandles}
\end{table}

Though we do not have associativity in (non-trivial) quandles, we can use the defining axioms to write a typical element in a standard form.

\begin{proposition}\label{left associated form in quandles}
Any element in a quandle $X$ can be written in the form
\begin{equation*}
\big(\left(\cdots\left(\left(x_0*^{\epsilon_1} x_1\right)*^{\epsilon_2} x_2\right)*^{\epsilon_3}\cdots\right)*^{\epsilon_{n-1}} x_{n-1}\big)*^{\epsilon_n} x_n,
\end{equation*}
where $x_i \in X$ and $\epsilon_i \in \{ 1, -1 \}$ are such that $x_0\neq x_1$, and if $x_i=x_{i+1}$ for any $1 \le i \le n-1$, then $\epsilon_i= \epsilon_{i+1}$.
\end{proposition}

\begin{proof}
Using the third quandle axiom in $X$, we see that
\begin{equation}\label{left association identity}
x*^\epsilon\left(y*^\delta z\right)=\big((x*^{-\delta} z)*^\epsilon y\big)*^\delta z
\end{equation}
for all $x,y,z\in X$ and $\epsilon,\delta\in\{1,-1\}$.  By a repeated use of \eqref{left association identity}, any element of $X$ can be written in the form
\begin{equation}\label{left association identity 2}
\big(\left(\cdots\left(\left(x_0*^{\epsilon_1} x_1\right)*^{\epsilon_2} x_2\right)*^{\epsilon_3}\cdots\right)*^{\epsilon_{n-1}} x_{n-1}\big)*^{\epsilon_n} x_n,
\end{equation}
where $x_i \in X$ and $\epsilon_i \in \{ 1, -1 \}$. If $x_0= x_1$, then the expression \eqref{left association identity 2} reduces to
$$ \big(\cdots((x_0 *^{\epsilon_2} x_2 )*^{\epsilon_3}\cdots )*^{\epsilon_{n-1}} x_{n-1} \big)*^{\epsilon_n} x_n.$$
And, if $x_i=x_{i+1}$ for any $1 \le i \le n-1$ and $\epsilon_i= -\epsilon_{i+1}$, then the expression reduces to
$$
\big((\cdots (( x_0*^{\epsilon_1} x_1 ) *^{\epsilon_2}\cdots *^{\epsilon_{i-1}} x_{i-1})*^{\epsilon_{i+2}} x_{i+2}  \cdots )*^{\epsilon_{n-1}} x_{n-1} \big)*^{\epsilon_n} x_n.
$$

After all the cancellations have been incorporated in \eqref{left association identity 2}, we can write each element of $X$ in the desired form. $\blacksquare$
\end{proof}

For simplicity, we write the expression \eqref{left association identity 2} as
\begin{equation*}
x_0*^{\epsilon_1} x_1*^{\epsilon_2}\cdots*^{\epsilon_n} x_n.
\end{equation*}
and refer it as the \index{left-associated product}{\it left-associated product}. We will use this convention throughout the monograph.
\para

Proposition \ref{left associated form in quandles} yields the following result \cite[Lemma 4.4.8]{MR2634013}.

\begin{proposition}\label{Lem:the canonical left associated form}
Let $X$ be a quandle. Then the product 
\begin{equation*}
(x_0*^{\epsilon_1}x_1*^{\epsilon_2}\cdots*^{\epsilon_m}x_m)*^{\mu_0}\left(y_0*^{\mu_1}y_1*^{\mu_2}\cdots*^{\mu_n}y_n\right)
\end{equation*}
of two left-associated expressions in $X$ is the left-associated expression
\begin{equation*}
x_0*^{\epsilon_1}x_1*^{\epsilon_2}\cdots*^{\epsilon_m}x_m*^{-\mu_n}y_n*^{-\mu_{n-1}}y_{n-1}*^{-\mu_{n-2}}\cdots*^{-\mu_1}y_1*^{\mu_0}y_0*^{\mu_1}y_1*^{\mu_2}\cdots*^{\mu_n}y_n.
\end{equation*}
\end{proposition}
\bigskip
\bigskip

\section{Special classes  of quandles}\label{special classes of quandles}

In this section, we introduce some special classes of quandles. Following \cite[Definition 1.11]{MR1994219}, we have the following definition.

\begin{definition}
A quandle is said to be \index{faithful quandle} {\it faithful} if the natural map $S: X \to \Inn(X)$ given by $S(x)=S_x$ is injective.    
\end{definition}

For instance, dihedral quandles of odd order are faithful, while those of even order are not. Similarly, trivial quandles of order more than one are not faithful.

\begin{definition}
A quandle $X$ is said to be \index{medial quandle}{\it medial} if
\begin{equation}\label{medial identity}
  (w * x) * (y * z) = (w * y) * (x * z)
\end{equation}
for all $x, y, z, w \in X$.
\end{definition}

We note that medial quandles are referred to as abelian quandles in \cite[p. 38]{MR0638121}, and the term `medial' for quandles satisying \eqref{medial identity} was first used in \cite{MR3400403}.  Note that a medial quandle is both right-distributive and left-distributive. For example, if $G$ is an abelian group and $\phi$ an automorphism of $G$, then the generalised Alexander quandle structure given by
$$
x *y=  \phi (x) +(\id-\phi)(y)
$$
is medial.
\para

The following are two related definitions.

\begin{definition}
A quandle $X$ is called \index{abelian quandle}{ \it abelian} if 
$$(x * y) * z = (x * z) * y$$
for all $x, y, z \in X$. 
\end{definition}

\begin{proposition}
Every abelian quandle is medial.
\end{proposition}

\begin{proof}
Let $X$ be an abelian quandle and $x, y, z, w \in X$. Note that $X$ being abelian implies that $(x*y)*^{-1}z=(x*^{-1}z)*y$. We compute
\begin{eqnarray*}
(w * x) * (y * z) &=& \big(((w * x) *^{-1} z )* y \big) * z\\
&=& \big(((w * x) *y ) *^{-1} z \big) * z\\
&=& \big((w * x) *y \big)\\
&=& \big((w * y) * x \big)\\
&=& \big(((w * y) * x) *^{-1} z \big) * z\\ 
&=& \big(((w * y) *^{-1} z) * x \big)* z\\ 
&=& (w * y) * (x * z), 
\end{eqnarray*}
and hence $X$ is medial. $\blacksquare$
\end{proof}
 
 Lebed and Mortier considered abelian quandles in \cite{MR4116819}, wherein they showed that quandles with abelian structure groups are necessarily abelian.

\begin{definition}
A quandle $X$ is said to be \index{commutative quandle}{ \it commutative} if
$$x *y = y *x$$
for all $x, y \in X$. 
\end{definition}

Note that, unlike in group theory, being abelian and being commutative do not mean the same for quandles. For example, the quandle with the following multiplication table is abelian (and hence medial) but not commutative.
$$
\begin{tabular}{|c||c|c|c|c|}
    \hline
$*$ & 1 & 2 & 3 &4 \\
  \hline \hline
1 & 1 & 3 & 4 & 2  \\
\hline
2 & 4 & 2 & 1 & 3  \\
\hline
3 & 2 & 4 & 3 & 1 \\
\hline
4 & 3 & 1 & 2 & 4 \\
  \hline
\end{tabular}
$$
In fact, any trivial quandle with more than one element is abelian but not commutative. On the other hand, the dihedral quandle $\R_3$ on three elements is both abelian as well as commutative.
\para

\begin{proposition}
Let $G$ be a group and $\varphi \in \Aut(G)$. If $\Alex(G, \varphi)$ is commutative, then $\varphi(a^2)=a$ for all $a \in G$. Further, the converse holds if and only if $G$ is abelian.
\end{proposition}

\begin{proof}
If $\Alex(G, \varphi)$ is commutative, then for $a, b \in G$, we have $\varphi(a)\varphi(b)^{-1}b=\varphi(b)\varphi(a)^{-1}a$. Taking $b=1$ yields $\varphi(a^2)=a$. Conversely, if $\varphi(a^2)=a$ for all $a \in G$, then $a*b=b*a$ if and only if $\varphi(ab)= \varphi(ba)$. This holds if and only if $G$ is abelian.  $\blacksquare$
\end{proof}

If $G$ is an additive abelian group, then $\Alex(G, \varphi)$ is the usual Alexander quandle. As a consequence, we recover the following result of Bae and Choi \cite[Theorem 2.8]{MR3265406}.

\begin{corollary}
Let $G$ be an additive abelian group and $\varphi \in \Aut(G)$. Then the  Alexander quandle $\Alex(G, \varphi)$ is commutative if and only if $2\varphi = \id_G$.
\end{corollary}
\para

\begin{definition}
A quandle $X$ is called simple if every homomorphism from $X$ to any quandle $Y$ is either constant or injective.
\end{definition}

It turns out that simple quandles embed in conjugation quandles of groups.

\begin{proposition} 
Let $X$ be a simple quandle of order more than 2. Then the following  assertions hold:
\begin{enumerate}
\item The map $S : X \to \Inn(X)$ is injective.
\item $X$ is connected.
\item The image $S(X)$ of $X$ is a generating conjugacy class in $\Inn(X)$.
\item  The center $\Z(\Inn(X))$ of $\Inn(X)$ is trivial.
\end{enumerate}
\end{proposition}

\begin{proof} 
If $S$ is not injective, then there exist $x \neq y \in X$ such that $x * y =x$ and $y*x=y$. But, the only simple quandle satisfying $x * y =x$ has no more than two elements. This proves assertion (1).
\para
For assertion (2), it is enough to consider the quotient of $X$ whose elements are its connected components. Since $X$ is simple, it follows that it must be connected.
\para

By definition of $\Inn(X)$, the image $S(X)$ of $X$ generates $\Inn(X)$. Since $X$ is connected, it follows that $S(X)$ is a conjugacy class in $\Inn(X)$, which proves assertion (3). 
\para

Let $f$ be an element of the center of $\Inn(X)$. Then, for each $x \in X$, we have
$$
S_x = f S_x f^{-1} = S_{f(x)}.
$$
It follows from the injectivity of $S$ that $x = f(x)$, and hence $f$ is the identity automorphism on $X$. $\blacksquare$
\end{proof}

\begin{definition}
Let $n \ge 1$ be an integer. A quandle $X$ is called an \index{$n$-quandle}{\it $n$-quandle} if 
$$x*^n y:=x* \underbrace{y *y*\cdots *y}_{n ~\mathrm{times}}=x$$
for all $x, y \in X$.
\end{definition}

Equivalently, a quandle is an $n$-quandle if and only if each inner automorphism $S_x$ has order dividing $n$. For example, if $G$ is a finite group of order $n$ or of exponent $n$, then $\Conj(G)$ is an $n$-quandle. Note that a $1$-quandle is simply a trivial quandle. In the literature, 2-quandles are referred to as \index{involutory quandle}{\it involutory} quandles. For example, $\Core(G)$ is an involutory quandle for any group $G$. If $X$ is an involutory quandle, then $$(x*y)*y=x$$ for all $x, y \in X$. Equivalently,
$$
x * y = x *^{-1} y
$$
for all $x, y \in X$.

\begin{exercise}{\rm 
Determine which quandles from the preceding section are commutative, abelian or involutory.}
\end{exercise}
\medskip

Motivated by \index{Burnside group}{Burnside groups} and the preceding definition, we now introduce the following.

\begin{definition}
Let  $m, n \geq 1$ be integers. The \index{Burnside quandle} {\it Burnside quandle} $\BQ(m, n)$ is the $n$-quandle generated by $m$ elements, that is,
$$
\BQ(m, n) = \big\langle x_1, \ldots, x_m   \, \mid \,  x *^n y = x ~\mbox{for all} ~x, y \big\rangle.
$$
\end{definition}

We formulate the following quandle analogue of the well-known Burnside problem for groups.

\begin{problem} (Burnside problem for quandles)
Determine pairs $(m, n)$ for which the  Burnside quandle $\BQ(m, n)$ is finite.
\end{problem}

\begin{remark}{\rm 
One can see that for each $m, n \ge 1$, $\BQ(m, 1)$ is simply the trivial quandle with $m$ elements and $\BQ(1, n)$ is the trivial quandle with one element. Further, if $m\ge 2$, then $\BQ(m, 2)$ is the free involutory quandle generated by $m$ elements, which is left-orderable by Theorem \ref{free invol left orderable}. But, any non-trivial left-orderable quandle is infinite  by Proposition \ref{orderable quandle infinite}.}
\end{remark}
\medskip

Next, we consider quandles determined by words in free groups. Note that the conjugation quandle and the core quandle of a group are defined by words $w_1(a, b) = b a b^{-1}$  and $w_2(a, b) = b a^{-1} b$, respectively. In general, suppose that $w = w(x, y)$ is a reduced word in the free group $F_2$ on two generators $x$ and $y$. For each group $G$, we define the binary operation $g *_w h = w(g, h)$. If the algebraic structure $(G,*_w)$ is a rack, then we call it a \index{verbal rack}{\it verbal rack}  defined by the word $w$. Similarly, if   $(G,*_w)$ is a quandle, then we call it a \index{verbal quandle}{\it verbal quandle}. The following result classifies all words in $F_2$ which gives a rack or  a quandle structure on groups \cite[Proposition 3.1]{MR4406425}.

\begin{proposition}\label{verbal rack quandle prop}
Let $w=w(x,y)\in F(x,y)$ be a reduced word such that $(G, *_{w})$ is a rack for every group $G$. Then $w(x,y) = x^{\pm 1}$, $w(x,y) = y x^{-1} y$ or $w(x,y) = y^n x y^{-n}$ for some $n \in \mathbb{Z}$. Further, in the latter two cases, $(G, *_{w})$ is a quandle.
\end{proposition}

\begin{proof}Let $w = x^{\alpha_1} y^{\beta_1} \cdots x^{\alpha_k} y^{\beta_k}$ for $\alpha_i, \beta_i \in \mathbb{Z}$ be a reduced word, where all $\alpha_i, \beta_i$ are non-zero with possible exceptions for $\alpha_1$ and $\beta_k$. Since $(G,*_w)$ is a rack, it follows from the second quandle axiom that for every $a, b \in G$ there exists an element $c \in G$ such that 
$$c^{\alpha_1} a^{\beta_1} \cdots c^{\alpha_k} a^{\beta_k}=c *_w a = b.$$ 
But this is possible if and only if $w =  y^{\alpha}  x^{\varepsilon} y^{\beta}$ for $\alpha, \beta \in \mathbb{Z}$ and $\varepsilon \in \{\pm 1 \}$. It follows from the third quandle axiom that 
$(x *_w y) *_w z = (x *_w z) *_w(y *_w z)$  for all $x,y,z\in G$. This can be rewritten as
\begin{equation}\label{a3}
z^{\alpha} (y^{\alpha} x^{\varepsilon} y^{\beta})^{\varepsilon} z^{\beta} = (z^{\alpha} y^{\varepsilon} z^{\beta})^{\alpha} (z^{\alpha} x^{\varepsilon} z^{\beta})^{\varepsilon} (z^{\alpha} y^{\varepsilon} z^{\beta})^{\beta}.
\end{equation}
Putting $x=y=1$ in \eqref{a3}, we obtain
$$
z^{(\alpha + \beta)} = z^{(\alpha + \beta) (\alpha + \beta + \varepsilon)}.
$$
Thus,  $\alpha + \beta = 0$ and $\varepsilon = \pm 1$ or $\alpha + \beta \not = 0$ and $1 = \alpha + \beta + \varepsilon$. Thus, in this case, we have $\varepsilon = -1$. We now proceed case by case as follows.
\para
\begin{itemize}
\item[Case 1:] $\alpha + \beta = 0$. In this case \eqref{a3} can be rewritten in the form
$$
z^{\alpha} y^{\alpha} x y^{-\alpha} z^{-\alpha} = z^{\alpha} y^{\alpha \varepsilon} x y^{-\alpha \varepsilon} z^{-\alpha},
$$
which is further equivalent to 
$$y^{\alpha} x y^{-\alpha} =  y^{\alpha \varepsilon} x y^{-\alpha \varepsilon} .$$
Since this equality holds, in particular, in the free group $F(x, y)$ on two generators, we have $\alpha(1-\varepsilon)=0$. Thus, either $\alpha=0$ and $\varepsilon = \pm 1$ or $\alpha\neq 0$ and $\varepsilon = 1$. Hence, in this case, we have either 
$w=(x, y)=x^\varepsilon$ or  $w(x, y) = y^{\alpha} x y^{-\alpha}$.

\item[Case 2:]  $\alpha + \beta = 2$. In this case $\varepsilon = -1$, $\beta = 2 - \alpha$, and equality \eqref{a3} can be rewritten in the form
\begin{equation}\label{case2verbal}
z^{\alpha} y^{\alpha-2} x y^{-\alpha} z^{2-\alpha} = (z^{\alpha} y^{-1} z^{2-\alpha})^{\alpha}
z^{\alpha-2} x z^{-\alpha} (z^{\alpha} y^{-1} z^{2-\alpha})^{2-\alpha}.
\end{equation}
Since this equality must hold in every group for all $x,y,z$, in particular, it holds in the free group $F(x, y, z)$ on three generators. Depending on $\alpha$, we have the following subcases.
\para
\item[Case 2.1:]   $\alpha<0$. In this case \eqref{case2verbal} implies that
$$
z^{\alpha} y^{\alpha-2} x y^{-\alpha} z^{2-\alpha} = \underbrace{(z^{\alpha-2} y z^{-\alpha}) \ldots (z^{\alpha-2} y z^{-\alpha})}_{-\alpha~ \mbox{\tiny times}} z^{\alpha-2}  x z^{-\alpha}
\underbrace{(z^{\alpha} y^{-1} z^{2-\alpha}) \ldots (z^{\alpha} y^{-1} z^{2-\alpha})}_{2-\alpha~ \mbox{\tiny times}}.  
$$
The left side of this equality is reduced, while the word on the right side of this equality after cancellations starts with $z^{\alpha-2} y z^{-2}$. This equality cannot hold in $F(x,y,z)$, and hence this case is not possible.

\item[Case 2.2:]   $\alpha = 0$. In this case equality \eqref{case2verbal} implies that $y^{-2} x = z^{-2} x y^{-1} z^{2} y^{-1}$. But it is not true in $F(x,y,z)$, and hence this case is not possible.

\item[Case 2.3:] $\alpha = 1$. In this case $w = y x^{-1} y$ and equality \eqref{case2verbal} obviously holds.

\item[Case 2.4:] $\alpha = 2$. In this case \eqref{case2verbal} implies that $x y^{-2} = y^{-1} z^{2} y^{-1} x z^{-2}$. But it is not true in $F(x,y,z)$, and hence this case is not possible.

\item[Case 2.5:] $\alpha >2$. In this case $2-\alpha < 0$, and equality \eqref{case2verbal} implies that
$$
z^{\alpha} y^{\alpha-2} x y^{-\alpha} z^{2-\alpha} = \underbrace{(z^{\alpha} y^{-1} z^{2-\alpha}) \ldots (z^{\alpha} y^{-1} z^{2-\alpha})}_{\alpha~ \mbox{\tiny times}}
z^{\alpha-2}  x z^{-\alpha}
\underbrace{(z^{\alpha-2} y z^{-\alpha}) \ldots (z^{\alpha-2} y z^{-\alpha})}_{\alpha-2~ \mbox{\tiny times}}.
$$
The left side of this equality is reduced, while the word on the right side of this equality after cancellations starts by $z^{\alpha} y^{-1} z^{2}$. This equality cannot hold in $F(x,y,z)$, and hence this case is also not possible.
\end{itemize}

Thus, we have proven that a word $w$ defines a rack $(G, *_w)$ on arbitrary group $G$ if and only if $w(x,y)=x^{\pm 1}$, $w(x,y) = y^n x y^{-n}$ for some $n \in \mathbb{Z}$ or $w(x,y) = y x^{-1} y$. In the latter two cases the rack $(G, *_w)$ is, in fact,  a quandle. $\blacksquare$
\end{proof}

\begin{remark}{\rm 
A similar question for groups was formulated by Cooper \cite[Problem 6.47]{Kourovka}. Let $G$ be a group and $w$ a word in two variables such that the operation $x *y = w(x, y)$ defines a new group structure $G_w = (G, *)$ on the set $G$. Does $G_w$ always lie in the variety generated by $G$? The problem seems to be open in general. }
\end{remark}

We have encountered examples of quandles defined by words over certain extended alphabets. This gives rise to the following natural problem.

\begin{problem}
Let $F_n$ be the free group of rank $n$ with a free basis $\{x_1, \ldots, x_n\}$. Determine all words in the extended alphabets $\{ x_1^{\pm 1}, \ldots, x_n^{\pm 1}, \varphi_1^{\pm 1}, \ldots, \varphi_m^{\pm 1} \}$, where $\varphi_i \in \Aut(F_n)$, which define quandle structures on $F_n$.
\end{problem}

Concerning the preceding problem, the paper \cite{Markhinina-Nasybullov-2022} determines all words $w(x,y,z)$ in the free group $F(x,y,z)$, such that for every group $G$ and an element $c \in G$, the set $G$ with the binary operation  given by $a*b=w(a,b,c)$ for $a,b\in G$, is a quandle.
\bigskip
\bigskip


\section{Free quandles and free racks}\label{section free quandles}
Free quandles are free objects in the category of quandles.  The following construction of a model of a free rack is due to Fenn and Rourke \cite[p. 351]{MR1194995}. Let $S$ be a set and $F(S)$ the free group on $S$. Define
 $$FR(S):=S \times F(S) = \big\{ a^{w}:= (a,w)  \, \mid \,  a \in S, ~w \in F(S) \big\} $$ with the operation defined as $$a^{w} \ast b^{u}:= a^{u bu^{-1}w}.$$
It can be seen that $FR(S)$ is a \index{free rack}{\it free rack} on the set $S$. 

A model of the free quandle on the set $S$ is due to Kamada \cite{Kamada2012, MR3588325}, who defined the \index{free quandle}{\it free quandle} $FQ(S)$ on $S$ as a quotient of $FR(S)$ modulo the equivalence relation generated by $$a^{w}= a^{wa}$$ for $a \in S$ and $w \in F(S)$. It is not difficult to check that $FQ(S)$ is quandle satisfying the universal property that for any quandle $X$ and a set-theoretic map $\psi: S \to X$, there exists a unique quandle homomorphism $\widetilde{\psi}: FQ(S) \to X$ extending $\phi$. For convenience, we use the same notation to denote elements of both $FQ(S)$ and $FR(S)$. If $S$ is a finite set with $n$ elements, then we denote the free quandle $FQ(S)$ by $FQ_n$.

\para

There is another model of free quandle on a set $S$ \cite[Example 2.16]{MR3729413}, which is defined as the subquandle of $\Conj\big(F(S) \big)$ consisting of all conjugates of elements of $S$ in the free group $F(S)$, that is,
$$FQ(S)= \bigcup_{x\in S} \big\{w x w^{-1}  \, \mid \,  w \in F(S) \big\}.$$

 Below is an explicit isomorphism between the two models of a free quandle.

\begin{proposition}\label{equivalence-two-models}
The map $\Phi: FQ(S) \rightarrow \Conj \big(F(S)\big)$ given by $\Phi (a^{w})=w a w^{-1} $ is an embedding of quandles.
 \end{proposition}

\begin{proof}
Let $a_{1}^{w_{1}},a_{2}^{w_{2}} \in FQ(S)$. Then $\Phi(a_{1}^{w_{1}})=w_{1} a_{1}w_{1}^{-1}$, $\Phi(a_{2}^{w_{2}})=w_{2} a_{2}w_{2}^{-1}$ and $a_{1}^{w_{1}} \ast a_{2}^{w_{2}}= a_{1}^{w_{2}a_{2}w_{2}^{-1}w_{1}}$. Further,
\begin{align*}
\Phi(a_{1}^{w_{1}} \ast a_{2}^{w_{2}})&=\Phi(a_{1}^{w_{2}a_{2}w_{2}^{-1}w_{1}})\\
&=(w_{2}a_{2}w_{2}^{-1}w_{1}) a_1 (w_{1}^{-1}w_{2}a_{2}^{-1}w_{2}^{-1})\\
&=(w_{2}a_{2}w_{2}^{-1})(w_{1} a_1 w_{1}^{-1})(w_{2}a_{2}w_{2}^{-1})^{-1}\\
&= \Phi(a_2^{w_2}) \Phi(a_1^{w_1}) \Phi(a_2^{w_2})^{-1}\\
&=\Phi(a_1^{w_1}) \ast \Phi(a_2^{w_2}),
\end{align*}
and hence $\Phi$ is a quandle homomorphism.  Let $a_1^{w_1},  a_2^{w_2} \in FQ(S)$ such that $a_1^{w_1}$  $ \neq $  $a_2^{w_2}$.
\para
\begin{enumerate}
\item[Case 1:] Suppose $ a_1 $ $ \neq $ $ a_2 $. If $\Phi( a_1^{w_1} )= \Phi( a_2^{w_2} )$, then $ w_1 a_1 w_1^{-1} = w_2 a_2 w_2^{-1}$, which contradicts the fact that $F(S)$ is a free group. Hence, we have  $\Phi( a_1^{w_1} ) \neq \Phi( a_2^{w_2} )$.
\para
\item[Case 2:] Suppose $a_1 = a_2=a$. If $\Phi( a^{w_1} )= \Phi( a^{w_2} )$, then $ w_1 a w_1^{-1} = w_2 a w_2^{-1}$, which further implies that $ w_1^{-1} w_2$ commutes with $a$ in $F(S)$. Since $F(S)$ is a free group, only powers of $a$ can commute with $a$, and hence $w_1^{-1} w_2$ = $a^{i}$ for some integer $i$. Thus, we have $w_2$= $w_1a^{i}$, which implies that $ a^{w_2}$ = $ a^{w_1a^{i}} = a^{w_1}$ in $ FQ(S)$, a contradiction. Hence, $\Phi( a_1^{w_1} ) \neq \Phi( a_2^{w_2} )$ and $\Phi$ is an embedding of quandles.
\end{enumerate}
This completes the proof. $\blacksquare$
\end{proof}
    
The classical Nielsen--Schreier theorem states that every subgroup of a free group is free \cite[Chapter 3, Theorem 3.3]{MR1812024}. We now present an analogue of this theorem for quandles due to Ivanov, Kadantsev, and Kuznetsov \cite{IvanovKadantsevKuznetsov}.

\begin{lemma}\label{lemma_if_S_ind}
Let $Y$ be a subset of the free 	quandle $FQ(S)$ such that $Y$ is a basis of a free subgroup of the free group $F(S)$. Then the subquandle of $FQ(S)$ generated by  $Y$ is free.
\end{lemma}

\begin{proof}
Consider the group homomorphism $F(Y)\to F(S)$ induced by the embedding $Y\hookrightarrow F(S)$. Since $Y$ is a basis of a free subgroup of $F(S)$, the homomorphism is injective. Then the restriction to the free quandles $FQ(Y)\to FQ(S)$ is also injective. The subquandle generated by $Y$ in $FQ(S)$ is equal to the image of the injective quandle morphism $FQ(Y)\to FQ(S),$ and hence it  is free.  $\blacksquare$
\end{proof}

\begin{theorem} \cite{IvanovKadantsevKuznetsov}
\label{IvanovKadantsevKuznetsov theorem}
 A subquandle of a free quandle is free.
\end{theorem}

\begin{proof}
Let $Q$ be a subquandle of the free quandle $FQ(S)$ on the set $S$. We are going to construct a subset $T(Q)$ of $Q$ and prove that it is a free basis of $Q$. We view an element $w=w_1w_2\dots w_n$ of the free group $F(S)$ as a reduced word, where $w_i\in S\cup S^{-1}$ and $w_i \neq w_{i-1}^{-1}$ for all $i$.
The length of the word $w$ is denoted by $|w|:=n$. First, for an element $x\in S$ we consider the subset
$$ P_x:= \big\{ w\in F(S)  \, \mid \,  x^w\in Q, ~\textrm{where}~ w_n\ne x^{\pm 1} \big\},$$
consisting of all reduced words $w$ such that $x^w\in Q$ and whose last letter, if it exists, differs from $x$ and $x^{-1}.$  Then we define the set
$$R_x := \big\{ w \in P_x  \, \mid \,   |q^\varepsilon w| > |w| ~\textrm{for all}~ q\in Q~\textrm{and all}~ \varepsilon\in \{1,-1\}  \big\}.$$
In some sense, $R_x$ consists of ``non-shrinkable'' elements of $P_x.$  Finally, we consider the set 
$$T :=  \bigcup_{x \in S} \big\{x^w  \, \mid \,  w \in R_x \big\},$$
which is obviously a subset of  $Q$.
\para
We first claim that the set $T$ generates the quandle $Q$. Let $\langle T \rangle$ denote the subquandle of $Q$ generated by $T$. It suffices to prove that $Q \subseteq \langle T\rangle$. Note that for any element $q\in Q$, there is an element $x\in S$ and $w\in P_x$ such that $q=x^w.$ We proceed by induction on the length $n$ of such words $w=w_1w_2\dots w_n$.  If $n=0$, then $w=1$, and hence $1\in P_x$ and $1\in R_x$. Thus, $x \in T \subseteq \langle T \rangle$, which proves the base step.  
\para 
Now, we prove the induction step. Assume that $|w|=n$, where $x\in S$ and $x^w\in Q.$ If $w\in R_x,$ then the statement is obvious. So, we can assume that $w\notin R_x.$ In this case there exists $q\in Q$ such that $|q^\varepsilon w|\leq |w|$ for some $\varepsilon\in \{1,-1 \}.$ Note that  elements from $FQ(S)$ have odd lengths. Hence,  $|q^\varepsilon w|$ and $|w|$ have different parity. Thus, they cannot be equal, and consequently $|q^\varepsilon w|<n$. We set $w':=q^\varepsilon w.$ Then $x^{w'}= (q^\varepsilon w) x (q^\varepsilon w)^{-1}= x^w *^{\varepsilon} q$, and hence $x^{w'}\in Q.$ By induction hypothesis, we have $x^{w'}\in \langle T \rangle.$ Since $q\in Q,$ there exists $y\in S$ and $u\in P_y$ such that $q=y^u.$ Note that the last letter of $u$ is not equal to $y,$ and hence there are no cancelations in the product $u y^\varepsilon u^{-1}.$  Since the length of $w'=u y^\varepsilon u^{-1} w $ is less than the length of $w,$ we obtain that $ u^{-1}$ completely cancels in the product $u^{-1} w$.  Thus, we have $|u|<|w|=n.$ Again, by the induction hypothesis, we obtain $q=y^u\in \langle T\rangle.$ Combining this with the fact that $x^{w'}\in \langle T \rangle$  and the equation $x^w=x^{w'} *^{-\varepsilon} q,$ we obtain $x^w\in \langle T \rangle.$
\para

Next, we claim that the set $T$ is a basis of a free subgroup of the free group $F(S)$. Following Hall \cite[\S \, 7.2]{MR0103215} we say that a subset $Y$ of the free group $F(S)$ such that $Y\cap Y^{-1}=\emptyset$ {\it posses significant factors} if there is a collection of indexes $\{i(w)\}_{w\in Y\cup Y^{-1}}$ such that the following conditions are satisfied:
\begin{enumerate}[(i)]
\item $0\leq i(w)\leq |w|$,
\item   $i(w^{-1})=|w|+1-i(w)$,
\item for each $w,v\in Y\cup Y^{-1}$ with $w v\ne 1$,  the cancelation in the product $wv$ does not reach the factors $w_{i(w)}$ and $v_{i(v)}$.  
\end{enumerate}
It is known that if a subset possess significant factors, then it is a basis of a free subgroup  \cite[Theorem 7.2.2]{MR0103215}. Thus, it suffices to show that $T$ possess significant factors. For each $x^w\in T$ we choose the central factor $x$ as a significant factor and set $i(x^w):=|w|+1.$ Then we only need to prove that the central factors $x^\varepsilon,y^\delta$  of the words $x^{\varepsilon w} $ and $y^{\delta v}$ do not cancel in the product $x^{\varepsilon w}y^{ \delta v}$ for any $x,y\in S$, $w\in R_x$, $ v\in R_y$ and $\varepsilon,\delta\in \{1,-1\}$ such that $x^{\varepsilon w}y^{\delta v} \neq 1$.
\para
Suppose on the contrary that one of the factors $x^\varepsilon$ or $y^\delta $ cancels in the product $x^{\varepsilon w}y^{ \delta v}=w x^\varepsilon w^{-1} v y^\delta v^{-1}$. Then one of the following  assertions holds: 
\begin{enumerate}[(i)]
\item $v$ can be presented as a product without cancellations  $v=w x^{-\varepsilon}u$ for some $u$.
\item  $w^{-1}$ can be presented as a product without cancellations $w^{-1}=uy^{-\delta}v^{-1}$ for some $u.$
\end{enumerate} 
In the first case we have $x^{\varepsilon w} v=wu,$ and hence $|x^{\varepsilon w} v|<|v|,$ which contradicts to the fact that $v\in R_y.$ In the second case we have $y^{-\delta v} w=v u^{-1},$ and hence $|y^{-\delta v} w|<|w|,$ which contradicts the fact that $ w\in R_x.$
\para

We have proved that $Q$ is generated by the set $T$, which is a basis of a free subgroup of $F(S)$. The result now follows from Lemma \ref{lemma_if_S_ind}.  $\blacksquare$
\end{proof}

As expected, free quandles are universal in the category of quandles, and every quandle is quotient of some free quandle. Thus, every quandle $X$ has a presentation $X= \langle S \mid R \rangle$, where $S$ is a set of generators and $R$ is a set of defining relations. We will not delve into presentations of quandles; instead, we refer the reader to 
 \cite[Section 4.2]{MR2634013} for details. 
\para 

We can define free products of quandles analogous to other algebraic structures. Let
$$
X = \langle S_1 \mid R_1 \rangle ~\textrm{and}~Y = \langle S_2\mid R_2 \rangle
$$
be two quandles with non-intersecting sets of generators. Then their \index{free product of quandles}{\it free product} $X \star Y$ is a quandle that is defined by the presentation
$$
X \star Y = \langle S_1 \sqcup S_2 \mid R_1 \sqcup R_2\rangle.
$$
The definition can be extended to arbitrary free product of quandles. For example, if  $FQ_n$ is the free quandle of rank $n$, then one can see that
$$
FQ_n = \underbrace{T_1 \star T_1 \star \cdots \star T_1}_{n~\textrm{copies}},
$$
the free product of $n$ copies of trivial one element quandles. Free product of racks can be defined analogously.  We refer the reader to \cite[Section 7]{MR4129183} for more on free products of racks and quandles. 
\bigskip
\bigskip


\section{Connected and flat quandles}\label{section connected and flat quandles}

The natural action of the inner automorphism group $\Inn(X)$ on a quandle $X$ via evaluation yields orbits that are central to the structural analysis of 
$X$.

\begin{definition}
A quandle $X$ is said to be \index{connected quandle}{\it connected} if the inner automorphism  group $\Inn(X)$ acts transitively on $X$. 
\end{definition}

In other words, a quandle $X$ is connected if and only if, for each pair  $x, y$ of elements in $X$, there exist elements $z_1,\ldots, z_n \in X$ and $\epsilon_1, \ldots, \epsilon_n \in \{1,-1\}$ such that
$$
x *^{\epsilon_1} z_1 *^{\epsilon_2}  \ldots *^{\epsilon_n} z_n = y.
$$
\para

By an \index{orbit}{\it orbit} of a quandle, we mean an orbit under the action of its inner automorphism group. An orbit is sometimes referred to as a \index{connected component}{\it connected component}. It is clear that every trivial quandle with more than one element is disconnected; in fact, each orbit in this case is a singleton.
\para

\begin{exercise}{\rm 
Show that a dihedral quandle of odd order is connected, while a dihedral quandle of even order is disconnected, consisting of exactly two orbits.}
\end{exercise}

\begin{proposition}
Homomorphic image of a connected quandle is connected.
\end{proposition}

\begin{proof}
The proof follows immediately from the connectedness of the quandle and the surjectivity of the homomorphism. $\blacksquare$
\end{proof}

\begin{definition}
A connected $3$-manifold $M$ is said to be \index{irreducible 3-manifold}{\it irreducible} if every embedded $2$-sphere in $M$ bounds a $3$-ball in $M$. A 3-manifold that is not irreducible is called \index{reducible}{\it reducible}. 
\end{definition}

We have the following definition concerning links in the 3-sphere.

\begin{definition}
A link in the 3-sphere $\mathbb{S}^3$ is said to be \index{split link}{\it split} if the complement of the link in $\mathbb{S}^3$ is reducible. If the complement is irreducible, we say that the link is \index{non-split}{\it non-split}.
\end{definition}

The following result is known due to Joyce \cite[Corollary 15.3]{MR0638121} and Matveev \cite{MR0672410}.

\begin{theorem}
The link quandle of a link is connected if and only if the link is non-split.
\end{theorem}

\begin{definition}
A quandle in which left multiplication by each element is a bijection is called a \index{latin quandle}{\it latin} quandle. In particular, the multiplication table of a finite latin quandle is a latin square. 
\end{definition}

It is evident that every latin quandle is connected, but the converse is not true in general. For instance, if $X = \{ (12),~ (13),~(14),~~ (23),~  (24),~ (34) \}$ is the subquandle of $\Conj(\Sigma_4)$ consisting of all transpositions, then $X$ is connected but not latin.
\para

For a finite connected quandle, left multiplications by any two of its elements have the same list (possibly with repeats) of lengths of their disjoint cycles.  A criterion for a finite connected quandle to be latin has been given by Lages and Lopes \cite[Theorem 2.2]{LagesLopes2021}. 

\begin{theorem}
If  $X$ is a finite connected quandle for which the list of lengths of disjoint cycles of its left multiplications by elements has no repeats, then $X$ is latin. 
\end{theorem}

For instance, cyclic quandles satisfy the hypothesis of the preceding theorem, where a finite quandle $X$ is called \index{cyclic quandle} \textit{cyclic} if each $S_x$ is a cycle of length $|X|-1$.
\bigskip

Given a quandle $X$, we consider two subgroups of $\Inn(X)$. First being the \index{transvection group}{\it transvection group} $$\Trans(X):=\langle S_xS_y^{-1}  \, \mid \,  x,y\in X \rangle$$ considered first by Joyce \cite{MR0638121}, and the second being the \index{displacement group}{\it displacement group} $$\Dis(X):=\langle S_xS_y  \, \mid \,  x,y\in X\rangle,$$ which appeared in the work \cite[Section 4]{MR3544543}. The latter definition is motivated by Riemannian symmetric spaces, where it is known that a connected Riemannian symmetric space is flat (that is, the Riemannian curvature vanishes identically) if and only if its group of displacements is abelian. 
\para

\begin{definition}
A quandle is called \index{flat quandle}{\it flat} if its group of displacements is abelian.     
\end{definition}

Clearly, trivial quandles and Takasaki quandles of abelian groups are flat. Note that, for any $x, y \in X$, we have 
$$S_x S_y^{-1}=(S_x S_x) (S_y S_x)^{-1} \in \Dis(X),$$ and hence $\Trans(X) \le \Dis(X)$. It turns out that just like $\Inn(Q)$, both $\Trans(X)$ and $\Dis(X)$ are also normal subgroups of $\Aut(X)$. Hence, there is a normal series
$$[\Inn(X), \Inn(X)] \trianglelefteq \Trans(X) \trianglelefteq \Dis(X) \trianglelefteq \Inn(X) \trianglelefteq \Aut(X).$$

\begin{proposition}
The following statements hold for any quandle $X$:
\begin{enumerate}
\item  $[\Inn(X), \Inn(X)] \le \Trans(X)$. Further, if $X$ is connected, then  $[\Inn(X),\Inn(X)]  = \Trans(X)$.
\item  $\Inn(X) / \Trans(X)$ is a cyclic group.  
\end{enumerate}
\end{proposition}

\begin{proof}
For any $x, y \in X$, since $S_xS_y^{-1} \in \Trans(X)$, we have the equality $S_x\Trans(X)= S_y\Trans(X)$ of cosets. This implies that $\Inn(X)/\Trans(X)$ is a cyclic group, and hence $[\Inn(X),\Inn(X)] \le \Trans(X)$. Suppose that $X$ is connected. If $x, y \in X$, then there exists an element $f \in \Inn(X)$ such that $f(y)=x$. Since $S_{f(y)}= f S_y f^{-1} $, we have  $S_xS_y^{-1}=S_{f(y)}S_y^{-1}= f S_y f^{-1} S_y^{-1} \in [\Inn(X),\Inn(X)]$, and hence the equality holds. This proves assertion (1). Assertion (2) is immediate from the definition of $\Trans(X)$. $\blacksquare$
\end{proof}
\para

It turns out that connectedness of a quandle $X$ can be characterized in terms of the action of $\Dis(X)$ \cite[Lemma 4.3]{MR3544543}.

\begin{proposition}\label{Isihara Tamaru prop}
A quandle $X$ is connected if and only if $\Dis(X)$ acts transitively on $X$.
\end{proposition}

\begin{proof}
Observe that, for any $x, y \in X$, we have
$$S_x^{-1}  S_y= (S_y  S_x)^{-1}  (S_y  S_y) \in \Dis(X)$$
and
$$S_x  S_y^{-1}= (S_x  S_x)  (S_y  S_x)^{-1} \in \Dis(X).$$
Let $X$ be connected and  $x,y \in X$. Then there exists $f \in \Inn(X)$ such that $f(x) = y $. We write $f= S_{x_k}^{\pm1}  \cdots  S_{x_1}^{\pm1}$ for some $x_i \in X$. If $k$ is even, then we can write
$$f= (S_{x_k}^{\pm1}  S_{x_{k-1}}^{\pm1})  \cdots  (S_{x_2}^{\pm1}  S_{x_1}^{\pm1}),$$
and hence $f \in \Dis(X)$ by the preceding observations. And, if $k$ is odd, then by the same argument $f  S_x \in \Dis(X)$ and $f  S_x (x) = f(x) = y$. Hence, $\Dis(X)$ acts transitively on $X$. $\blacksquare$
\end{proof}
\para

Next, we provide a complete characterisation of flat connected quandles. Before that, we recall Example \ref{coset quandle with auto}.  Let $G$ be a group, $\phi \in \Aut(G)$, and $H$ a subgroup of $G$ such that $H \le \Fix(\phi)= \{x \in G \mid \phi(x)=x \}$. Then the set of right cosets of $H$ in $G$ forms a quandle $\Alex(G, H, \phi)$, with its binary operation given by 
$$Hx*Hy= H\phi(xy^{-1})y$$
for all $x, y \in G$. For brevity, we denote the coset of an element $x$ by $\overline{x}$. Then we have the following result \cite[Proposition 3.7]{MR3544543}.

\begin{proposition}
Let $G_1, \ldots, G_n$ be groups, where $n \ge 2$. For each $i$, let $\phi_i \in \Aut(G_i)$ and $H_i$ a subgroup of $G_i$ such that $H_i \le \Fix(\phi_i)= \{x \in G_i \mid \phi_i(x)=x \}$. Then 
$$\Alex\Big(\prod_{i=1}^n G_i, \prod_{i=1}^n H_i, \prod_{i=1}^n \varphi_i \Big) \cong \prod_{i=1}^n \Alex(G_i, H_i, \varphi_i).$$
\end{proposition}

\begin{proof}
It suffices to consider the case $n=2$. Consider the map $\Phi: (G_1 \times G_2)/ (H_1 \times H_2) \to (G_1/ H_1) \times (G_2/ H_2)$ given by $$\Phi\big(\,\overline{(x_1, x_2)} \,\big)= \big(\overline{x_1}, \overline{x_2}\big)$$ for $(x_1, x_2) \in G_1 \times G_2$. Clearly, $\Phi$ is well-defined and bijective. It only remains to show that $\Phi$ is a quandle homomorphism. For $(x_1, x_2), (y_1, y_2) \in G_1 \times G_2$, we have
\begin{eqnarray*}
\Phi\big(\, \overline{(x_1, x_2)}*\overline{(y_1, y_2)} \, \big) &=& \Phi\big(\, \overline{(\varphi_1 \times \varphi_2) (x_1y_1^{-1}, x_2y_2^{-1})(y_1, y_2)} \, \big)\\
&=& \Phi\big(\, \overline{(\varphi_1(x_1y_1^{-1})y_1, \varphi_2(x_2y_2^{-1})y_2)} \,\big)\\
&=& \big(\, \overline{\varphi_1(x_1y_1^{-1})y_1}, \overline{\varphi_2(x_2y_2^{-1})y_2}\, \big)\\
&=& \big(\overline{x_1}*\overline{y_1},  \overline{x_2}*\overline{y_2} \big)\\
&=& \big(\overline{x_1},  \overline{x_2} \big)*\big(\overline{y_1},  \overline{y_2} \big)\\
&=& \Phi\big(\,\overline{(x_1, x_2)} \,\big) * \Phi\big(\,\overline{(y_1, y_2)} \,\big),
\end{eqnarray*} 
which is desired. $\blacksquare$
\end{proof}

As a consequence, we have the following result for Takasaki quandles.

\begin{corollary}\label{product of takasaki quandles}
If $G_1, \ldots, G_n$ are abelian groups, then 
$$\T(G_1 \times \cdots \times G_n) \cong \T(G_1)  \times \cdots \times \T(G_n)$$
\end{corollary}
\para

The next result shows that flat connected quandles are particular types of Alexander quandles \cite[Proposition 5.2]{MR3544543}. 

\begin{theorem}\label{Ishihara and tamaru theorem}
Let $X$ be a flat connected quandle. Then $X \cong  \Alex(A,\varphi)$, where $A$ is an abelian group and $\varphi \in \Aut(A)$ is an involution.
\end{theorem}

\begin{proof}
Let $x \in X$ be a fixed element and $\Dis(X)_x$ the stabilizer at $x$ under the natural action of $\Dis(X)_x$ on $X$. We claim that  $\Dis(X)_x$ is trivial. Let $\alpha \in \Dis(X)_x$ and $y \in X$ be arbitrary elements. Since $X$ is connected, by Proposition \ref{Isihara Tamaru prop}, there exists $\beta \in \Dis(X)$ such that $\beta (x)=y$. Since $\Dis(X)$ is abelian, we obtain 
$$\alpha(y)= \alpha \beta (x)= \beta \alpha(x)=\beta(x)=y,$$
and hence $\Dis(X)_x$ is trivial. Thus, the orbit map $\Psi: \Dis(X) \to X$ given by $\alpha \mapsto \alpha(x)$ is a bijection. In fact, it turns out that $\Phi: \Alex \big(\Dis(X), \iota_{S_x} \big) \cong X$ is an isomorphism, where $\iota_{S_x}$ is the inner automorphism of $\Dis(X)$ induced by $S_x$. Furthermore, since $S_x^2 \in \Dis(X)$, we have $S_x^2 \, \alpha= \alpha \, S_x^2$ for all $\alpha \in \Dis(X)$. Hence, $\iota_{S_x}$ is an order two automorphism of $\Dis(X)$. $\blacksquare$
\end{proof}
\para

Note that if $G$ is an abelian group and $\varphi \in \Aut(G)$, then the map $\tilde{\varphi}:G \to G$ defined by $$\tilde{\varphi}(a)=a^{-1}\varphi(a)$$ is a group homomorphism. 

\begin{proposition}\label{alexander connnected-abelian}
Let $G$ be an abelian group and $\varphi \in \Aut(G)$ an involution. Then $\Alex(G, \varphi)$ is connected if and only if $\tilde{\varphi}$ is surjective.
\end{proposition}

\begin{proof}
For $g, a \in G$, we have
$$S_g  S_g (a)=S_g \big (\varphi(a) \varphi(g^{-1}) g\big)=\varphi^2(a) \varphi^2(g^{-1}) \varphi(g) \varphi(g^{-1}) g=a.$$
Thus, $\Alex(G, \varphi)$ is involutory, and equivalently $S_g=S_g^{-1}$ for all $g \in G$. Let $\Alex(G, \varphi)$ be connected and  $b \in G$. Then there exists a $f \in \Inn(G)$ such that $b^{-1}= f(1)$, where $f= S_{g_k}  \cdots S_{g_1}$ for some $g_i \in G$. Then we have
\begin{eqnarray*}
b^{-1} & = & S_{g_k}  \cdots  S_{g_1}(1)\\
&=& \big(\cdots ((1 * g_1) * g_2) *\cdots \big) *g_k\\
&=& \big(\cdots((\varphi(g_1^{-1}) g_1)* g_2 \big) *\cdots ) *g_k\\
&=& \varphi^k(g_1^{-1}) \varphi^{k-1}(g_1) \varphi^{k-1}(g_2^{-1}) \varphi^{k-2}(g_2) \varphi^{k-2}(g_3^{-1}) \cdots \varphi^2(g_{k-1}^{-1}) \varphi(g_{k-1}) \varphi(g_k^{-1}) g_k\\
&=& \varphi^k(g_1^{-1}) \varphi^{k-1}(g_1) \varphi^{k-1}(g_2^{-1}) \varphi^{k-2}(g_2) \varphi^{k-2}(g_3^{-1}) \cdots \varphi^2(g_{k-1}^{-1})   \varphi(g_k^{-1}) g_k \varphi(g_{k-1}),\\
&=& \varphi^k(g_1^{-1}) \varphi^{k-1}(g_1) \varphi^{k-1}(g_2^{-1}) \varphi^{k-2}(g_2) \varphi^{k-2}(g_3^{-1}) \cdots \varphi \big(\varphi(g_{k-1}^{-1})  g_k^{-1} \big) \big(\varphi(g_{k-1}^{-1}) g_k^{-1} \big)^{-1}.
\end{eqnarray*}
We see that
$$b^{-1} =\varphi(a)a^{-1},~ \textrm{where}~ a=\varphi^{k-1}(g_1^{-1})\varphi^{k-2}(g_2^{-1}) \cdots \varphi(g_{k-1}^{-1})g_k^{-1},$$
which gives $b= a \varphi(a^{-1})=\tilde{\varphi}(a^{-1}) $ and hence $\tilde{\varphi}$ is surjective.
 \para
Conversely, suppose that $\tilde{\varphi}$ is surjective.  If $a, b \in G$,  then there exists $c \in G$ such that $\tilde{\varphi}(c)= \varphi(a)b^{-1}$. This implies that $\varphi(c)c^{-1} = \varphi(a)b^{-1}$. Equivalently $b=a*c=S_c(a)$, and $\Alex(G, \varphi)$ is connected. $\blacksquare$
\end{proof}

The following result gives a complete characterisation of flat connected quandles \cite[Theorem 4.6]{MR3576755}.

\begin{theorem}\label{flat-connected}
A quandle  $X$ is flat connected if and only if $X \cong  \T(A)$ for a 2-divisible group $A$.
\end{theorem}

\begin{proof}
If $X$ is flat and connected, then Theorem \ref{Ishihara and tamaru theorem} gives $X \cong  \Alex(A,\varphi)$, where $A$ is an (additive) abelian group and $\varphi \in \Aut(A)$ is an involution. Further, since $X$ is connected, by Proposition \ref{alexander connnected-abelian}, the map $\tilde{\varphi}$ is surjective. Thus, given $x \in A$ there exists $y \in A$ such that $x = \tilde{\varphi}(y)= \varphi(y)-y$. This implies that $\varphi(x) = \varphi^2(y)- \varphi(y) = y-\varphi(y) = -x$, and hence $\varphi$ is simply the inversion. Thus, $\Alex(A,\varphi)= \T(A)$ and $A = 2A$, that is, $A$ is a 2-divisible group.
\para
For the converse, let $A$ be an additive 2-divisible group. If $x, y, x', y', g \in A$, then we have
$$S_x  S_y (g) = S_x \big(2y-g)= 2x-2y+g= S_{x-y}(-g).$$
Further, we see that
$$ \big((S_x  S_y )  (S_{x'}  S_{y'}) \big) (g)  =  (S_x  S_y) \big(S_{x'-y'}(-g)\big)= S_x \big( S_{y-x'+y'}(g)\big)= S_{x+x'-y-y'}(-g)$$
and
$$ \big((S_{x'}  S_{y'})  (S_x  S_y) \big) (g) = (S_{x'}  S_{y'}) \big(S_{x-y}(-g)\big)= S_{x'} \big( S_{y'-x+y}(g)\big)= S_{x+y'-x-y'}(-g).$$
This shows that $\T(A)$ is flat. Further, since $A$ is 2-divisible, $\T(A)$ is connected by Proposition \ref{alexander connnected-abelian}. $\blacksquare$
\end{proof}

As a consequence, we recover the main result of \cite[Theorem 1.1]{MR3544543} stated in another form.

\begin{corollary}
A quandle $X$ is finite, flat and connected if and only if $X$ is the Takasaki quandle of a finite abelian group of odd order.
\end{corollary}

\begin{remark}
{\rm 
In view of the structure theorem for finite abelian groups, $A \cong \mathbb{Z}_{n_1} \times \cdots \times \mathbb{Z}_{n_k}$, where each $n_i$ is an odd prime power.  Hence, by Corollary \ref{product of takasaki quandles}, a quandle is finite, flat and connected if and only if  it is isomorphic to a direct product of dihedral quandles of odd prime power orders.
\par
It is not difficult to construct examples of quandles which are infinite, flat and connected. For example, if $A$ is a divisible abelian group, then $\T(A)$ is flat and connected. In fact, divisible groups are not even finitely generated. One can construct many more examples by taking direct sums and quotients of divisible abelian groups.}
\end{remark}

Since all flat connected quandles arise from abelian groups, it is worth investigating flatness of other quandles arising from groups. We present the following result for conjugation quandles of groups.

\begin{proposition}
Let $G$ be a group. Then $\Conj(G)$ is flat if and only if $G$ is nilpotent of class 2.
\end{proposition}

\begin{proof}
If $a, b, a', b' \in G$, then $S_a  S_b (g) =  S_a \big(b g b^{-1} \big)=S_{ab}(g)$ and $S_{a'}  S_{b'} = S_{a'b'}$ for all $g \in G$.  Now, $\Conj(G)$ is flat if and only if
$$S_{(ab)(a'b')}=S_{ab} \, S_{a'b'} =S_{a'b'} \, S_{ab}=S_{(a'b')(ab)}.$$
This is further equivalent to $[ab,a'b']  \in \Z(G)$. Taking $b=b'=1$, we get $[a, a' ] \in \Z(G)$, and hence $[G,G] \leq \Z(G)$. Thus, $\Conj(G)$ is flat if and  only if $G$ is nilpotent of class 2.
$\blacksquare$
\end{proof}
\para

Next, we consider the connectivity of generalized Alexander quandles.  Before doing so, we introduce the following definition.

\begin{definition}
An automorphism $\varphi$ of a group $G$ is called \index{central automorphism}{\it central} if $a^{-1}\varphi(a) \in \Z(G)$ for all $a \in G$.    
\end{definition}

Note that these are precisely the automorphisms of $G$ which induce the identity map on the central quotient $G/\Z(G)$. The group of all central automorphisms is a subgroup of $\Aut(G)$, and denoted by $\Autcent(G)$. Note that $\Autcent(G)=\Aut(G)$ for an abelian group $G$. However, there exist non-abelian groups $G$ for which $\Autcent(G)=\Aut(G)$. See Curran and McCaughan \cite{MR1837963} for detailed references. 

\begin{lemma}\label{central-auto}
Let $G$ be a group and $\varphi \in \Autcent(G)$. Then the following  assertions hold:
\begin{enumerate}
\item The map $\tilde{\varphi}:G \to \Z(G)$ defined by $\tilde{\varphi}(a)=a^{-1}\varphi(a)$ is a group homomorphism. 
\item The map $\varphi \mapsto \tilde{\varphi}$ gives an injection $\Autcent(G) \to \Hom_{\mathcal{G}}\big(G, \Z(G)\big)$.
\item If $\varphi$ is fixed-point free, then $G$ is abelian.
\end{enumerate}
\end{lemma}

\begin{proof}
For $a, b \in G$, we have
\begin{eqnarray*}
\tilde{\varphi}(ab) & = & (ab)^{-1}\varphi(ab)\\
&=& b^{-1}a^{-1}\varphi(a) \varphi(b)\\
&=& a^{-1}\varphi(a) b^{-1} \varphi(b),~\textrm{since}~a^{-1}\varphi(a) \in \Z(G)\\
&=& \tilde{\varphi}(a)\tilde{\varphi}(b),
\end{eqnarray*}
and hence assertion (1) holds. Clearly, the map $\varphi \mapsto \tilde{\varphi}$ is injective, which is assertion (2). Further, note that $\varphi$ is fixed-point free if and only if  $\tilde{\varphi}$ is injective. Thus, $\tilde{\varphi}:G \to \Z(G)$ is an embedding and $G$ is abelian, which proves asertion (3). $\blacksquare$
 \end{proof}      

We have the following result concerning the connectedness of generalized Alexander quandles.

\begin{theorem}\label{connnected-abelian}
Let $G$ be a group and $\varphi \in \Autcent(G)$ an involution. If $\Alex(G, \varphi)$ is connected, then $G$ is abelian.
\end{theorem}

\begin{proof}
Since $\varphi^2= \id_G$, for all $g, a \in G$, we have
$$S_g  S_g (a)=S_g \big (\varphi(a) \varphi(g^{-1}) g\big)=\varphi^2(a) \varphi^2(g^{-1}) \varphi(g) \varphi(g^{-1}) g=a.$$
Hence, $S_g=S_g^{-1}$ for all $g \in G$. Let $\Alex(G, \varphi)$ be connected and  $b \in G$. Then there exists a $f \in \Inn(G)$ such that $b^{-1}= f(1)$, where $f= S_{g_k}  \cdots  S_{g_1}$ for some $g_i \in G$. We compute
\begin{eqnarray*}
b^{-1} & = & S_{g_k}  \cdots   S_{g_1}(1)\\
&=& \big(\cdots ((1 * g_1) * g_2) *\cdots \big) *g_k\\
&=& \big(\cdots((\varphi(g_1^{-1}) g_1)* g_2) *\cdots \big) *g_k\\
&=& \varphi^k(g_1^{-1}) \varphi^{k-1}(g_1) \varphi^{k-1}(g_2^{-1}) \varphi^{k-2}(g_2) \varphi^{k-2}(g_3^{-1}) \cdots \varphi^2(g_{k-1}^{-1}) \varphi(g_{k-1}) \varphi(g_k^{-1}) g_k\\
&=& \varphi^k(g_1^{-1}) \varphi^{k-1}(g_1) \varphi^{k-1}(g_2^{-1}) \varphi^{k-2}(g_2) \varphi^{k-2}(g_3^{-1}) \cdots \varphi^2(g_{k-1}^{-1})   \varphi(g_k^{-1}) g_k \varphi(g_{k-1}),\\
&& \textrm{since}~a^{-1}\varphi(a) \in \Z(G)~\textrm{for all}~ a \in G\\
&=& \varphi^k(g_1^{-1}) \varphi^{k-1}(g_1) \varphi^{k-1}(g_2^{-1}) \varphi^{k-2}(g_2) \varphi^{k-2}(g_3^{-1}) \cdots \varphi \big(\varphi(g_{k-1}^{-1})  g_k^{-1} \big) \big(\varphi(g_{k-1}^{-1}) g_k^{-1} \big)^{-1}.
\end{eqnarray*}
By repeatedly using the preceding property of $\varphi$, we get 
$$b^{-1} =\varphi(a)a^{-1},~ \textrm{where}~ a=\varphi^{k-1}(g_1^{-1})\varphi^{k-2}(g_2^{-1}) \cdots \varphi(g_{k-1}^{-1})g_k^{-1}.$$
 This gives $b= a \varphi(a^{-1})=\tilde{\varphi}(a^{-1}) \in \Z(G)$, and hence  $G=\Z(G)$.
$\blacksquare$  
\end{proof}          
\para

\begin{proposition}
Let $G$ be a group and $\varphi \in \Aut(G)$ a fixed-point free automorphism. The natural map $S:\Alex(G, \varphi) \hookrightarrow \Conj \big(\Inn(\Alex(G, \varphi))\big)$ given by $S(x)=S_x$ is an embedding of quandles.
\end{proposition}

\begin{proof}
By Proposition \ref{S quandle homomorphism}, $S$ is a quandle homomorphism. Further, $S(x)=S(y)$ together with the fixed-point freeness of $\varphi$ implies that $x=y$. Hence, $S$ is an embedding of quandles.  $\blacksquare$
\end{proof}
\para

The following result gives a characterisation of finite connected Alexander quandles (\cite[Theorem 2.5]{MR3265406} or \cite[Corollary 7.2]{MR3399387}).

\begin{theorem}\label{Bae-Choe-theorem}
Let $G$ be a finite abelian group and $\varphi \in \Aut(G)$. Then the following assertions are equivalent:
\begin{enumerate}
\item $\tilde{\varphi} \in \Aut(G)$.
\item $\varphi$ is fixed-point free.
\item  $\Alex(G, \varphi)$ is connected.
\end{enumerate}
\end{theorem}

\begin{proof}
Since $G$ is abelian, it follows that $\tilde{\varphi}$ is a group homomorphism. Further, since $G$ is finite, $\tilde{\varphi}$ is an automorphism if and only if $\ker(\tilde{\varphi})$ is trivial. But, $\ker(\tilde{\varphi})$ is trivial if and only if $\varphi$ is fixed-point free. This proves the equivalence of (1) and (2).
\para
Let $x, y \in G$. Assuming that $\tilde{\varphi} \in \Aut(G)$, there exists $z \in G$ such that $\tilde{\varphi}(z)=y^{-1}\varphi(x)$. Rewriting the expression gives  $y=\varphi(x)\varphi(z)^{-1} z= x*z$. Thus, assertion (1) implies assertion (3). Now, suppose that $\Alex(G, \varphi)$ is connected. Since $G$ is abelian, $\tilde{\varphi}$ is already a group homomorphism. Let $y \in G$ be an element. Since the quandle is connected, there exist $x_1, \ldots, x_k \in G$ and $\epsilon_1, \ldots, \epsilon_k \in \{1, -1 \}$ such that $y=S_{x_k}^{\epsilon_k} \cdots S_{x_1}^{\epsilon_1}(e)$, where $e$ is the identity element of the group. Rewriting the equation gives $y=e*^{\epsilon_1} x_1 *^{\epsilon_2}  \cdots *^{\epsilon_k} x_k$, which is equivalent to
\begin{equation}\label{phi tilde and phi}
y=\varphi^{\epsilon_2 + \cdots +\epsilon_k} \big(\varphi^{\epsilon_1}(x_1^{-1}) x_1 \big)~ \varphi^{\epsilon_3 + \cdots +\epsilon_k} \big(\varphi^{\epsilon_2}(x_2^{-1}) x_2 \big)  \cdots \varphi^{\epsilon_k}(x_k^{-1}) x_k.
\end{equation}
Note that, if $\epsilon_i=1$, then $\varphi^{\epsilon_i}(x_i^{-1}) x_i= \tilde{\varphi}(x_i^{-1})$. Similarly, if  $\epsilon_i=-1$, then $\varphi^{\epsilon_i}(x_i^{-1}) x_i= \varphi^{-1}\big(\tilde{\varphi}(x_i)\big)$. Since, $\tilde{\varphi} \varphi= \varphi \tilde{\varphi}$, it follows from \eqref{phi tilde and phi} that $y=\tilde{\varphi}(z)$ for some $z \in G$. This proves the equivalence of (1) and (3).
\end{proof}

An immediate consequence of  Theorem  \ref{Bae-Choe-theorem} is the following result.

\begin{corollary}
The dihedral quandle $\R_n$ is connected if and only if $n$ is odd.
\end{corollary}

\begin{remark}\label{purely-non-ab-def}
{\rm 
A group is said to be \index{purely non-abelian group} {\it purely non-abelian} if it does not have any non-trivial abelian direct factor. Adney and Yen \cite[Theorem 1]{MR0171845} proved that if $G$ is a finite purely non-abelian group, then the map $\varphi \mapsto \tilde{\varphi}$ is a bijection. Thus, for a finite purely non-abelian group $G$, in view of Theorem \ref{connnected-abelian}, the set $\Hom_{\mathcal{G}}\big(G, \Z(G)\big)=\Hom_{\mathcal{G}}\big(G/[G,G], \Z(G)\big)$ gives central automorphisms for which the generalised Alexander quandle structures on $G$ are not connected. For example, if $Q_8=\langle i, j, k  \, \mid \,  ijk=i^2=j^2=k^2=-1 \rangle $ is the quaternion group of order 8,  then $|\Autcent(Q_8)|=|\Hom_{\mathcal{G}}\big(Q_8, \mathbb{Z}_2\big)|=|\Hom_{\mathcal{G}}\big(\mathbb{Z}_2 \oplus \mathbb{Z}_2, \mathbb{Z}_2\big)|=4$.}
\end{remark}
\bigskip
\bigskip

\section{Construction of quandles via automorphisms}\label{sec9}
In this section, we describe some general approaches to constructing quandles using their automorphisms. For any group $G$ and an element $x \in G$, let $\iota_x$ denotes the inner automorphism of $G$ given by $\iota_x(y)=xy x^{-1}$. If $G$ is a group, then a map $\phi:G \to \Aut(G)$ is said to be \textit{compatible} if for every element $x \in G$ the following diagram commutes
$$ \begin{CD}
G @>\phi>> \Aut(G) \\
@V{\phi(x)}VV @VV{\iota_{\phi(x)}}V \\
G @>>\phi> \Aut(G).
\end{CD}
$$
For example, the maps $x\mapsto \id_G$ and $x\mapsto \iota_x$ are both compatible maps from $G$ to $\Aut(G)$.

\begin{proposition}\label{new1}
Let $G$ be a group and $\phi: G \to \Aut(G)$ be a compatible map such that $\phi(x)(x)=x$ for all $x \in G$. Then the set $G$ with the operation  $x*y=\phi(y)(x)$ is a quandle. Moreover, if $\phi$ is an injective map which satisfies $\phi(xy)=\phi(y)\phi(x)$ for all $x,y\in G$, then  $x*y=y x  y^{-1}$.
\end{proposition}

\begin{proof}
The equality $x*x=\phi(x)(x)=x$ for all $x \in G$ is given as the condition of the proposition. Thus, the first quandle axiom holds. For an element $x \in G$, the map $S_x:y\mapsto y*x=\phi(x)(y)$ is an automorphism of $G$, and hence is a bijection, which is the second quandle axiom.
Since $\phi$ is a compatible map, for $y,z\in G$, we have $\phi(z)\phi(y)=\phi\big(\phi(z)(y)\big)\phi(z)$. This yields
\begin{small}
$$(x*y)*z=\phi(z)(x*y)=\phi(z)\phi(y)(x)=\phi \big(\phi(z)(y) \big)\phi(z)(x)=\phi(z)(x)*\phi(z)(y)=(x*z)*(y*z),$$
\end{small}
and the set $G$ with the binary operation $x*y=\varphi(y)(x)$ is a quandle. 
\para 
If $\phi$ is a compatible map which satisfies the equality $\phi(xy)=\phi(y)\phi(x)$ for all $x,y\in G$, then  $\phi(x)\phi(y)=\phi\big(\phi(x)(y)\big)\phi(x)$. Hence, $\phi \big(\phi(x)(y)\big)=\phi(xyx^{-1})$, and if $\phi$ is injective, then $y*x={\phi(x)(y)}=xyx^{-1}$. $\blacksquare$
\end{proof}

For the compatible map $G\to \Aut(G)$ given by $x\mapsto \id_G$, the quandle constructed in Proposition \ref{new1} is trivial. Further, for the compatible map of the form $x \mapsto \iota_x$, the quandle constructed in Proposition \ref{new1} is the conjugation quandle $\Conj(G)$.

\begin{definition}
Let $X$ be a quandle and $S$ a set. Then a \index{quandle action}{\it quandle action} of $X$ on $S$ is a quandle homomorphism
$$\phi: X \to \Conj(\Sigma_S),$$
where $\Sigma_S$ is the symmetric group on the set $S$. 
\end{definition}

We can also define an action of a quandle $X$ on another quandle $Y$ by replacing $\Sigma_Y$ by $\Aut(Y)$, the group of quandle automorphisms of $Y$. For $x \in X$ and $z \in S$, we write $$z \cdot x := \phi(x) (z),$$ and we say that $X$ acts on $S$ from the right. Then $\phi$ being a quandle homomorphism is equivalent to the equation
$$(z \cdot x_2) \cdot (x_1 * x_2)= (z \cdot x_1 ) \cdot x_2$$
for $x_1, x_2 \in X$ and $z \in S$. In particular, if $S=X$ and $\phi(x)=S_x$ for each $x \in X$, then we recover the usual action of a quandle on itself by inner automorphisms. 
\para

As another example, if $G$ is a group acting on a set $S$, then the induced map $\Conj(G) \to \Conj(\Sigma_S)$ is a quandle action of  $\Conj(G)$ on the set $S$.
\para

The following result gives another interesting construction of quandles \cite[Proposition 11]{MR3948284}.

\begin{proposition}\label{new2}
Let $(X, *)$ and $(Y,  \circ)$ be quandles. Let $\sigma: X \to  \Conj \big(\Aut(Y) \big)$ and $\tau: Y \to  \Conj \big(\Aut(X) \big)$ be quandle homomorphisms. Then the set $Z=X \sqcup Y$ with the binary operation
$$
x\star y=\begin{cases}
x*y& ~\textit{if}~ x, y \in X, \\
x\circ y  & ~\textit{if}~ x, y \in Y, \\
{\tau(y)}(x)  & ~\textit{if}~x \in X, ~y \in Y, \\
{\sigma(y)}(x) & ~\textit{if}~x \in Y, ~y \in X,
\end{cases} 
$$
is a quandle if and only if the following conditions hold:
\begin{enumerate}
\item $\tau(z)(x)* y=\tau\big(\sigma(y)(z)\big)(x* y)$ for $x, y \in X$ and $z \in Y$,
\item $\sigma(z)(x)\circ y=\sigma\big(\tau(y)(z)\big)(x\circ y)$ for $x, y \in Y$ and $z \in X$.
\end{enumerate}
\end{proposition}

\begin{proof}
For  each $x\in Z$, we have $x\star x=x*x=x$ or $x\star x=x\circ x=x$. Thus, the first quandle axiom holds. For each $x\in X$, the map $S_x:y\mapsto y\star x$ acts as the map $y\mapsto y*x$ on $X$ and acts as the map $y\mapsto \sigma(x)(y)$ on $Y$. Thus, the map $S_x$ is a bijection on $Z$. Similarly, for each $x\in Y$, the map $S_x$ is a bijection on $Z$, and the second quandle axiom holds. Let $x,y,z\in Z$. If $x,y,z\in X$ or $x,y,z\in Y$, then the third quandle axiom holds obviously. Suppose that $x, y \in X$ and $z \in Y$. In this case, using (1), we have the equalities
\begin{align}
\notag(x\star y)\star z&=\tau(z)(x* y)=\tau(z)(x)* \tau(z)(y)=(x\star z)\star (y\star z),\\
\notag(z\star x)\star y &=  \sigma(x)(z)\star y=  \sigma(y) \sigma(x)(z)= \sigma(y) \sigma(x)\sigma(y)^{-1} \big(\sigma(y)(z) \big)\\
\notag&=  \sigma(x*y) \big(\sigma(y)(z)\big)= \sigma(y)(z)\star(x\star y)= (z\star y)\star(x\star y),\\
\notag(x\star z)\star y&=\tau(z)(x)* y=\tau\big(\sigma(y)(z)\big)(x* y)=(x* y)\star \sigma(y)(z)=(x\star y)\star(z\star y),
\end{align}
and hence the third quandle axiom holds. Similarly, using (2), we see that the third quandle axiom holds for $x, y \in Y$ and $z \in X$. $\blacksquare$
\end{proof}

\begin{example}\label{union quandle example}
{\rm 
Let us see some special cases of the preceding construction.
\begin{enumerate}
\item If $\sigma:X\to\{\id_{Y}\}$ and $\tau:Y\to\{\id_{X}\}$ are the trivial maps, then conditions (1) and (2) of Proposition \ref{new2} are satisfied. This gives a quandle structure on $X\sqcup Y$, which we refer to as the union of quandles. 

\item Let $X$ be an \textit{involutory quandle}, that is, $$(x*y)*y=x$$ for all $x,y\in X$. Let $Q_1$ and $Q_2$ be quandles such that there are isomorphisms $\varphi:X\to Q_1$ and $\psi:X\to Q_2$. Direct computations show that the maps $\sigma:\varphi(x)\mapsto S_{\psi(x)}$ and $\tau:\psi(x)\mapsto S_{\varphi(x)}$ satisfy conditions (1) and (2) of Proposition \ref{new2}. Thus, we can define a quandle structure on $Q_1\sqcup Q_2$ for each involutory quandle $X$.

\item Consider the quandle $\J_3=\{1,2,3 \}$ of Exercise \ref{all three element quandles} with multiplication table as follows.
$$
\begin{tabular}{|c||c|c|c|}
    \hline
$*$ & 1 & 2 & 3  \\
  \hline \hline
1 & 1 & 1 & 1  \\
\hline
2 & 3 & 2 & 2  \\
\hline
3 & 2 & 3 & 3  \\
  \hline
\end{tabular}
$$
Then, $\J_3$ can be thought of as obtained from two trivial quandles $X=\{2,3\}$ and $Y=\{1\}$ using the procedure of Proposition \ref{new2}, where $\sigma(2)=\sigma(3)=\id_{Y}$ and $\tau(1)$ is the automorphism of $X$ which permutes $2$ and $3$.
\end{enumerate}}
\end{example}


\chapter{Adjoint groups of quandles}\label{chapter-adjoint-groups-quandles}

\begin{quote}
The adjoint group of a rack or a quandle is simply the structure group of the solution to the Yang--Baxter equation associated to it. In this chapter, we discuss characteristics of adjoint groups within these frameworks. We explore augmented quandles and those that embed into their adjoint groups, thereby becoming unions of conjugacy classes. Additionally, we provide a classification of finite quandles whose adjoint groups are abelian.
\end{quote}
\bigskip

\section{Properties of adjoint groups}
In this section, we associate a group to a quandle that satisfies some universal property, and plays a crucial role in the structure theory of quandles.
\para

\begin{definition}
The  \index{adjoint group} {\it adjoint group} of a quandle $X$ is the group with the following presentation
$$\Adj(X)=\big\langle \textswab{a}_x, ~x \in X  \, \mid \,  \textswab{a}_{x*y}=\textswab{a}_y \textswab{a}_x \textswab{a}_y^{-1} ~ \textrm{for all} ~ x, y \in X \big\rangle.$$     
\end{definition}

\begin{remark}
{\rm 
Recall that a quandle $(X, *)$ gives rise to a solution $(X, r)$ to the Yang--Baxter equation via the braiding $r(x, y)=(y, x*y)$. It follows from the definitions that both the structure group and the derived structure group of the solution $(X, r)$ are simply the group $\Adj(X)$.  Following \cite[Section 6]{MR0638121},  we refer the  group as the adjoint group, but it is also known in the literature as  the {\it associated group} \cite[Section 2]{MR1194995}, the {\it enveloping group} \cite[Definition 1.5]{MR1994219} and the {\it structure group} \cite[Section 2]{MR1722951}.
}    
\end{remark}

\begin{example}
{\rm 
If $\T_n$ is the trivial quandle on $n$ elements, then $$\Adj(\T_n)=  \big\langle \textswab{a}_x, ~x \in \T_n  \, \mid \,  \textswab{a}_y \textswab{a}_x= \textswab{a}_x \textswab{a}_y ~ \textrm{for all} ~ x, y \in \T_n \big\rangle,$$ and hence $\Adj(\T_n)$ is the free abelian group of rank $n$.
}    
\end{example}

Let $\mathcal{Q}$ and $\mathcal{G}$ denote the categories of quandles and groups, respectively.  Taking the adjoint group of a quandle is a functor $$\Adj: \mathcal{Q} \to \mathcal{G}.$$ Further, $\Adj$ is left adjoint of the functor $$\Conj: \mathcal{G} \to \mathcal{Q}.$$ More precisely, there is a bijection
$$\Hom_{\mathcal{Q}} \big(X, \Conj(G) \big) \longleftrightarrow \Hom_{\mathcal{G}} \big(\Adj(X), G \big)$$
for each $X \in \mathcal{Q}$ and $G \in \mathcal{G}$. By functoriality, any quandle homomorphism $f: X \rightarrow Y$ induces a group homomorphism $f_\#: \Adj(X) \rightarrow \Adj(Y)$ defined as $$f_\#(\textswab{a}_x)=\textswab{a}_{f(x)}$$ for any $x \in X$. 
\para

There is a natural map $i_X: X \rightarrow \Adj(X)$ defined as $$i_X(x)=\textswab{a}_x$$ for all $x \in X$. The following result of Winker \cite[Theorem 5.1.7]{MR2634013} gives a smaller presentation of the adjoint group $\Adj(X)$ using a given presentation of the quandle $X$.

\begin{theorem}\label{presentation-of-adjoint-group}
Let $X$ be a quandle with a presentation $\langle S \mid R \rangle$. Then the adjoint group $\Adj(X)$ has a presentation $\langle \overline{S} \mid \overline{R} \rangle$, where 
$ \overline{S} = \{\textswab{a}_x  \, \mid \,  x \in S \}$ and $\overline{R}$ consists of relations in $R$ with each expression $x*y$ replaced by $\textswab{a}_y \textswab{a}_x \textswab{a}_y^{-1}$ and the expression $x *^{-1} y$ replaced by $\textswab{a}_y^{-1} \textswab{a}_x \textswab{a}_y$.
\end{theorem}

\begin{proof}
It follows from Proposition \ref{left associated form in quandles} that any relation in $R$ is of the form 
$$x_0 *^{\epsilon_1} x_1  *^{\epsilon_2} \cdots  *^{\epsilon_n} x_n=  y_0 *^{\varepsilon_1} y_1  *^{\varepsilon_2} \cdots  *^{\varepsilon_m} y_m$$
for some $x_i, y_j \in X$ and $\epsilon_i, \varepsilon_j \in \{1, -1\}$. Setting $G=\langle \overline{S}  \, \mid \,  \overline{R} \rangle$, it suffices to show that $G \cong \Adj(X)$. The map $\phi: \overline{S} \to \Adj(X)$ given  by $\phi(\textswab{a}_x)=\textswab{a}_x$ for $x \in S$ defines a group homomorphism from the free group $F(\overline{S})$ to $\Adj(X)$, which we denote by $\phi$ itself. If $x \in X$, then by Proposition \ref{left associated form in quandles}, we can write
$x=x_0 *^{\epsilon_1} x_1  *^{\epsilon_2} \cdots  *^{\epsilon_n} x_n$ for some $x_i\in S$ and $\epsilon_i \in \{1, -1\}$. Thus,
\begin{eqnarray*}
    \textswab{a}_x &=& \textswab{a}_{x_0 *^{\epsilon_1} x_1  *^{\epsilon_2} \cdots  *^{\epsilon_n} x_n}\\
    &=& \textswab{a}_{x_n}^{\epsilon_n} \textswab{a}_{x_{n-1}}^{\epsilon_{n-1}} \cdots \textswab{a}_{x_1}^{\epsilon_1} \textswab{a}_{x_0}\textswab{a}_{x_1}^{-\epsilon_1} \cdots  \textswab{a}_{x_{n-1}}^{-\epsilon_{n-1}} \textswab{a}_{x_n}^{-\epsilon_n}\\
    &=& \phi\big( \textswab{a}_{x_n}^{\epsilon_n} \textswab{a}_{x_{n-1}}^{\epsilon_{n-1}} \cdots \textswab{a}_{x_1}^{\epsilon_1} \textswab{a}_{x_0}\textswab{a}_{x_1}^{-\epsilon_1} \cdots  \textswab{a}_{x_{n-1}}^{-\epsilon_{n-1}} \textswab{a}_{x_n}^{-\epsilon_n}\big),
\end{eqnarray*}
and hence $\phi$ is surjective. By definition of $\Adj(X)$, if $\overline{r}=\overline{s}$ is any relation in $\overline{R}$, then $\phi(\overline{r})=\overline{r}=\overline{s}=\phi(\overline{s})$ in $\Adj(X)$ as well, and hence $\phi$ induces a surjective group homomorphism $\overline{\phi}: G \to \Adj   (X)$. Finally, injectivity of $\overline{\phi}$ follows from the fact that relations in the definition of $\Adj(X)$ can be derived from relations in $\overline{R}$. $\blacksquare$
\end{proof}
\para

Let  $X \star Y$ denote the free product of quandles $X$ and $Y$. For the sake of clarity, we denote by $G \hexstar H$ the free product of groups $G$ and $H$.

\begin{proposition}\label{presentation-of-adjoint-group-of-free-product-of-quandles}
If $X$ and $Y$ are quandles, then $\Adj(X \star Y) \cong  \Adj(X) \hexstar \Adj(Y).$
\end{proposition}

\begin{proof}
If $X$ and $Y$ have presentations $X = \langle S_1 \,\mid \, R_1 \rangle$ and $Y = \langle S_2  \, \mid \, R_2 \rangle$, then $X \star Y = \langle S_1 \sqcup S_2 \, \mid \, R_1 \sqcup R_2 \rangle$. Now, by Theorem \ref{presentation-of-adjoint-group}, we have
\begin{eqnarray*}
\Adj(X \star Y) & \cong & \langle \textswab{a}_x, ~x \in  S_1 \sqcup S_2 \, \mid \, \overline{R}_1 \sqcup \overline{R}_2 \rangle \\
& \cong &  \langle \textswab{a}_x, ~x \in S_1  \, \mid \, \overline{R}_1 \rangle \ast  \langle \textswab{a}_x,~x \in S_2  \, \mid \, \overline{R}_2 \rangle\\
& \cong &  \Adj(X) \hexstar \Adj(Y),
\end{eqnarray*}
which completes the proof. $\blacksquare$   
 \end{proof}

\begin{proposition}\label{rank-of-free-quandle}
Let $FQ(S)$ and $ FQ(T)$ be free quandles on sets $S$ and $T$, respectively. If $FQ(S)\cong FQ(T)$, then $|S|=|T|$.
\end{proposition}

\begin{proof}
It follows from Theorem \ref{presentation-of-adjoint-group} that for free quandles $FQ(S)$ and $FQ(T)$, we have $\Adj \big(FQ(S) \big)\cong F(S)$ and $\Adj\big(FQ(T) \big) \cong F(T)$, the free groups on the sets $S$ and $T$, respectively. Since $FQ(S)\cong FQ(T)$, we must have $\Adj \big(FQ(S) \big) \cong \Adj \big(FQ(T)\big)$, and hence $|S|=|T|$.
$\blacksquare$    \end{proof}

In view of Proposition \ref{rank-of-free-quandle}, we can define the \index{rank of free quandle}{\it rank of a free quandle} as the cardinality of its any free generating set.
\para

The adjoint group of a quandle satisfies the following properties \cite[Section 6]{MR0638121}.

\begin{proposition}\label{Prop: natural map from quandle to its adjoint group}
The following statements hold for any quandle $X$.
\begin{enumerate}
\item The natural map $i_X: X \rightarrow \Conj \big(\Adj(X) \big)$ is a quandle homomorphism.
\item For any group $G$ and a quandle homomorphism $f: X \rightarrow \Conj(G)$, there exists a unique group homomorphism $\bar{f}: \Adj(X) \rightarrow G$ such that the diagram 
	
		$$\begin{tikzcd}
		X \arrow[d, "f"'] \arrow[r, "i_X"] & \Conj \big(\Adj(X) \big) \arrow[d, "\bar{f}"] \\
		\Conj(G) \arrow[r, "\id_G"']           & G                          
		\end{tikzcd}$$
	\end{enumerate}
	commutes. Here, $\bar{f}$ is also viewed as a map $\Conj \big(\Adj(X) \big) \rightarrow G$.
\end{proposition}

\begin{proof}
If $x, y \in X$, then $i_X(x*y)=\textswab{a}_{x*y}= \textswab{a}_y \textswab{a}_x \textswab{a}_y^{-1}= \textswab{a}_x*\textswab{a}_y= i_X(x)* i_X(y)$. Thus, $i_X$ is a quandle homomorphism, which is assertion (1).
\para

Given a quandle homomorphism $f: X \to \Conj(G)$, we define $\bar{f}: \Adj(X) \rightarrow G$ by $\bar{f}(\textswab{a}_x)= f(x)$ for all $x \in X$. Then, for $x, y \in X$, we have
$$\bar{f} (\textswab{a}_{x*y})= f(x*y)=f(x)*f(y)=f(y) f(x) f(y)^{-1}= \bar{f}(\textswab{a}_y) \bar{f}(\textswab{a}_x) \bar{f}(\textswab{a}_y)^{-1}= \bar{f}(\textswab{a}_y \textswab{a}_x \textswab{a}_y^{-1}),$$
and hence $\bar{f}$ is a group homomorphism. Further, $\bar{f} i_X (x)= \bar{f}(\textswab{a}_x)= f(x)$ for all $x \in X$, which proves assertion (2). $\blacksquare$
\end{proof}

\begin{remark} \label{image-under-eta}
{\rm
If $X$ is a quandle and $f : X \rightarrow \{a\}$ the trivial quandle homomorphism, then the  induced group homomorphism $f_{\#}: \Adj(X) \rightarrow \Adj(\{a\}) \cong \mathbb{Z}$ is given by $f_{\#}(\textswab{a}_x) = 1$ for all $x \in X$. Thus, $i_X: X \rightarrow \Adj(X)$ does not map any element of $X$ onto the identity element of $\Adj(X)$.}
\end{remark}

\begin{remark} \label{eta not injective}
{\rm
The natural map $i_X: X \rightarrow \Adj(X)$ is not injective in general. For example, consider the quandle
$$
\J_3=\begin{tabular}{|c||c|c|c|}
    \hline
$*$ & 1 & 2 & 3  \\
  \hline \hline
1 & 1 & 1 & 1  \\
\hline
2 & 3 & 2 & 2  \\
\hline
3 & 2 & 3 & 3  \\
  \hline
\end{tabular}
$$
of Example \ref{all three element quandles}. Since $1=1*2$, we have $\textswab{a}_{1}=\textswab{a}_{1*2}=\textswab{a}_2 \textswab{a}_1 \textswab{a}_2^{-1}$, and hence $\textswab{a}_2\textswab{a}_1=\textswab{a}_1 \textswab{a}_2$. Now, $3=2*1$ implies that $$i_X(3)=\textswab{a}_3=\textswab{a}_{2*1}= \textswab{a}_1 \textswab{a}_2 \textswab{a}_1^{-1}=\textswab{a}_2=i_X(2),$$ and hence $i_X$ is not injective.}
\end{remark}

Given a quandle $X$, there is a group homomorphism $$\zeta: \Adj(X) \rightarrow \Inn(X)$$ defined on generators by $\zeta(\textswab{a}_x)=S_x$ for each $x \in X$. 

\begin{proposition}\label{adj inn central extension}
If $X$ is a quandle, then
$$
1 \to \ker(\zeta) \to \Adj(X) \to \Inn(X) \to 1
$$
is a central extension of groups.
\end{proposition}

\begin{proof}
The map $\zeta$ is clearly a surjective group homomorphism.  Further, the homomorphism $\zeta$ induces a left-action of $\Adj(X)$ on $X$ by $\textswab{a}_y \cdot x=S_y(x)=x*y$ for $x, y \in X$. If $x \in X$ and $g=\textswab{a}_{x_1}^{\epsilon_1}\textswab{a}_{x_2}^{\epsilon_2} \cdots \textswab{a}_{x_k}^{\epsilon_k} \in \Adj(X)$ for some $x_i \in X$ and $\epsilon_i \in \{1, -1\}$, then the definition of left group action implies that
\begin{eqnarray*}
\textswab{a}_{g \cdot x} &=& \textswab{a}_{(\textswab{a}_{x_1}^{\epsilon_1}\textswab{a}_{x_2}^{\epsilon_2} \cdots \textswab{a}_{x_k}^{\epsilon_k}) \cdot x}\\
&=& \textswab{a}_{\textswab{a}_{x_1}^{\epsilon_1} \cdot ( \cdots (\textswab{a}_{x_{k-1}}^{\epsilon_{k-1}} \cdot (\textswab{a}_{x_k}^{\epsilon_k} \cdot x)))}\\
&=& \textswab{a}_{(((x *^{\epsilon_k} x_k) *^{\epsilon_{k-1}} x_{k-1}) \cdots ) *^{\epsilon_1} x_1}\\
&=& \textswab{a}_{x_1}^{\epsilon_1} \textswab{a}_{x_2}^{\epsilon_2}  \cdots \textswab{a}_{x_k}^{\epsilon_k} \textswab{a}_x  \textswab{a}_{x_k}^{-\epsilon_k} \cdots \textswab{a}_{x_2}^{-\epsilon_2} \textswab{a}_{x_1}^{-\epsilon_1}\\
&=& g \textswab{a}_x g^{-1}.
\end{eqnarray*} 
Now, if  $g \in \ker (\zeta)$, then $g \cdot x=x$ for all $x \in X$. Thus, the preceding identity gives $\textswab{a}_x= g \textswab{a}_x g^{-1}$ for all $x \in X$, and hence $\ker(\zeta)$ lies in the center of $\Adj(X)$. $\blacksquare$
\end{proof}

One of the most important results concerning link quandles was independently established by Joyce \cite{MR0638121} and Matveev \cite{MR0672410}.

\begin{theorem}\label{Thm:MJ-2}
If $L$ is a link in $\mathbb{S}^3$, then $\Adj \big(Q(L) \big)\cong G(L)=\pi_1(\mathbb{S}^3 \setminus L)$.
\end{theorem}

\begin{example}\label{trefoil adjoint group b3}
{\rm 
Let $K$ be the trefoil knot. As in Example \ref{trefoil-coloring-r3}, the knot quandle of $K$ has a presentation
$$
Q(K) = \langle x, y, z   \, \mid \,  x * y = z,~ y * z = x,~ z * x = y \rangle.
$$
By Theorem \ref{presentation-of-adjoint-group}, a presentation of the adjoint group is 
$$\Adj \big(Q(K)\big)= \big\langle \textswab{a}_x, \textswab{a}_y, \textswab{a}_z  \, \mid \,  \textswab{a}_y \textswab{a}_x \textswab{a}_y^{-1} = \textswab{a}_z,~ \textswab{a}_z \textswab{a}_y \textswab{a}_z^{-1} = \textswab{a}_x, ~  \textswab{a}_x \textswab{a}_z \textswab{a}_x^{-1} = \textswab{a}_y \big\rangle.$$
We can remove the generator $\textswab{a}_z$ and the presentation takes the form
$$\Adj \big(Q(K)\big)= \big\langle \textswab{a}_x, \textswab{a}_y \, \mid \, \textswab{a}_y \textswab{a}_x \textswab{a}_y = \textswab{a}_x \textswab{a}_y \textswab{a}_x \big\rangle,$$ 
which is isomorphic to the braid group $\B_3$ on three strands.}
\end{example}

\begin{example}\label{dihh}
{\rm 
The multiplication table for the dihedral quandle $\R_3$ is as follows.
\begin{center}
\begin{tabular}{|c||c|c|c|}
    \hline
$*$ & 0 & 1 & 2  \\
  \hline \hline
0 & 0 & 2 & 1  \\
\hline
1 & 2 & 1 & 0  \\
\hline
2 & 1 & 0 & 2  \\
  \hline
\end{tabular}
\end{center}
The adjoint group $\Adj(\R_3)$ has three generators $\textswab{a}_0, \textswab{a}_1, \textswab{a}_2$ and following defining relations
$$\textswab{a}_0 \textswab{a}_1 \textswab{a}_0^{-1}=\textswab{a}_2=\textswab{a}_1 \textswab{a}_0 \textswab{a}_1^{-1},$$
$$\textswab{a}_1 \textswab{a}_2 \textswab{a}_1^{-1} =\textswab{a}_0=\textswab{a}_2 \textswab{a}_1 \textswab{a}_2^{-1},$$
$$\textswab{a}_0 \textswab{a}_2 \textswab{a}_0^{-1} =\textswab{a}_1=\textswab{a}_2 \textswab{a}_0 \textswab{a}_2^{-1}.$$
We remove the generator $\textswab{a}_2$ using the first relation and obtain the following relations
$$\textswab{a}_0 \textswab{a}_1 \textswab{a}_0^{-1}= \textswab{a}_1 \textswab{a}_0 \textswab{a}_1^{-1}, \quad \textswab{a}_0 \textswab{a}_1 \textswab{a}_0^{-1}= \textswab{a}_1^{-1} \textswab{a}_0 \textswab{a}_1, \quad \textswab{a}_1 \textswab{a}_0 \textswab{a}_1^{-1} = \textswab{a}_0^{-1} \textswab{a}_1 \textswab{a}_0.$$
The first and the second relation  imply that $\textswab{a}_1^2 \textswab{a}_0 \textswab{a}_1^{-2}= \textswab{a}_0$, that is, $\textswab{a}_1^2$ belongs to the center of $\Adj(\R_3)$. Similarly, the first and the third relation imply that $\textswab{a}_0^2$ belongs to the center of $\Adj(\R_3)$. The square of the first relation gives $\textswab{a}_0^2=\textswab{a}_1^2$, and we obtain $\textswab{a}_0 \textswab{a}_1\textswab{a}_0=\textswab{a}_1\textswab{a}_0\textswab{a}_1$. Hence, we have 
$$\Adj(\R_3)=\big\langle \textswab{a}_0,\, \textswab{a}_1  \, \mid \,  \textswab{a}_1 \textswab{a}_0 \textswab{a}_1 = \textswab{a}_0 \textswab{a}_1 \textswab{a}_0 ~\textrm{and}~ \textswab{a}_0^2 = \textswab{a}_1^2 \big\rangle.$$ In particular, $\Adj(\R_3)$ is a homomorphic image of the braid group $\B_3$  on three strands.
}
\end{example}
\para

For a quandle $X$, let $$\delta:\Adj(X) \to \mathbb{Z}$$ be the unique surjective homomorphism of groups given by $\delta(\textswab{a}_x)=1$ for all $x\in X$.  The following result gives a decomposition of adjoint groups of quandles \cite[Lemma 2.3]{MR3671570}.

\begin{proposition}\label{adjoint group as kernel}
The following statements hold for any quandle $X$.
\begin{enumerate}
\item $\Adj(X) \cong \ker (\delta) \rtimes \mathbb{Z}$.
\item If $X$ is connected, then $\ker (\delta) = [\Adj(X), \Adj(X)]$.
\end{enumerate}
\end{proposition}

\begin{proof}
One could think of $\delta$ as the group homomorphism of adjoint groups induced by the quandle homomorphism of $X$ onto the trivial one element quandle. Let $x \in X$ an element. Since $\ker (\delta)$ is a normal subgroup of $\Adj(X)$, it follows that the product $ \ker (\delta) \langle \textswab{a}_x\rangle$ is a subgroup of $\Adj(X)$. If $y \in X$, then $\textswab{a}_y= (\textswab{a}_y\textswab{a}_x^{-1}) \textswab{a}_x \in  \ker (\delta) \langle \textswab{a}_x\rangle$, and hence  $\Adj(X)=\ker (\delta) \langle \textswab{a}_x\rangle$. Since $\ker (\delta) \cap \langle \textswab{a}_x\rangle=1$ by definition of $\delta$, it follows that
 $\Adj(X)=\ker (\delta) \rtimes \langle \textswab{a}_x  \rangle \cong \ker (\delta) \rtimes \mathbb{Z}$. This proves assertion (1).
\para
It is clear that $[\Adj(X), \Adj(X)] \subseteq\ker (\delta)$.  Consider the length function $\ell:\Adj(X)\to \mathbb{Z}$ defined as $\ell(g)=n$ if $g=\textswab{a}_{x_1}^{\epsilon_1}\dots \textswab{a}_{x_n}^{\epsilon_n}$ is a  reduced expression of $g$ in terms of generators of $\Adj(X)$, where $x_i \in X$ and $\epsilon_i\in\{\pm1\}$. We show that $\ker (\delta) \subseteq [\Adj(X), \Adj(X)]$ by induction on $\ell(g)$ for $g\in\ker (\delta)$.  If $\ell(g)=2$, then $g=\textswab{a}_{x_1}^{\pm1} \textswab{a}_{x_2}^{\mp1}$ for some $x_1, x_2 \in X$. We can assume that $g=\textswab{a}_{x_1} \textswab{a}_{x_2}^{-1}$.  Since $X$ is connected, there exists $h \in \Adj(X)$ such
that $h\cdot \textswab{a}_{x_2}=\textswab{a}_{x_1}$. Hence, $g=h \textswab{a}_{x_2} h^{-1} \textswab{a}_{x_2}^{-1} = [h, \textswab{a}_{x_2}] \in [\Adj(X), \Adj(X)]$, and we are done. Now, suppose that $\ell(g)>2$. Then there is a reduced expression of $g$ (or $g^{-1}$) such that $g=g_1 \textswab{a}_{x_1} \textswab{a}_{x_2}^{-1}g_2$ for some $x_1, x_2 \in X$ and $g_1,g_2\in \Adj(X)$.  Since $g \in \ker(\delta)$, we have 
$$0=\delta (g)=\delta(g_1)+1-1+\delta(g_2)=\delta(g_1g_2),$$
and hence $g_1g_2 \in \ker (\delta)$. As $\ell(g_1g_2)<\ell(g)$, by induction hypothesis we have $g_1g_2\in [\Adj(X), \Adj(X)]$. Writing $g=(g_1  \textswab{a}_{x_1} \textswab{a}_{x_2}^{-1}g_1^{-1})(g_1g_2)$, we see that $g\in [\Adj(X), \Adj(X)]$, which proves assertion (2).  $\blacksquare$
\end{proof}

\begin{proposition}\label{fingen} 
Let $X$ be a quandle. Then the following assertions hold:
\begin{enumerate}
\item If $X$ has finitely many orbits, then the abelianization of the adjoint group $\Adj(X)$ is isomorphic to the free abelian group $\mathbb{Z}^{\left| \mathcal{O}(X)\right|}$, where $\mathcal{O}(X)$ is the set of all orbits of $X$.
\item If $X$ is finite, then the commutator subgroup of the adjoint group $\Adj(X)$  is finitely generated.
\end{enumerate}
\end{proposition}
\begin{proof}
Recall that $\Adj(X)$ is generated by elements $\{ \textswab{a}_x \, \mid \, x \in X \}$ subject to the defining relations $\textswab{a}_{x*y}= \textswab{a}_y \textswab{a}_x \textswab{a}_y^{-1}$. Thus, modulo the commutator subgroup $[\Adj(X), \Adj(X)]$ of $\Adj(X)$, the relation $\textswab{a}_{x*y}= \textswab{a}_y \textswab{a}_x \textswab{a}_y^{-1}$ has the form $\textswab{a}_{x*y}=\textswab{a}_x$ for all $x, y \in X$. Hence, if $X$ is finite, then we have an isomorphism between the abelianization of $\Adj(X)$ and $\mathbb{Z}^{|\mathcal{O}(X)|}$, which is assertion (1).
\para

It follows from the defining relations of $\Adj(X)$ that $\textswab{a}_{x *^{-1} y}= \textswab{a}_y^{-1} \textswab{a}_x \textswab{a}_y$. Further, we have
$$
\textswab{a}_{z} (\textswab{a}_x \textswab{a}_y \textswab{a}_x^{-1} \textswab{a}_y^{-1} ) \textswab{a}_z^{-1} = \textswab{a}_{z} [\textswab{a}_x,  \textswab{a}_y]
 \textswab{a}_z^{-1} = [\textswab{a}_{x*z},  \textswab{a}_{y*z}] \quad \textrm{and} \quad 
 \textswab{a}_{z}^{-1}  [\textswab{a}_x,  \textswab{a}_y]  \textswab{a}_z = [\textswab{a}_{x *^{-1} z},  \textswab{a}_{y *^{-1}  z}].
$$
Using these identities and the  commutator identities
$$
[a b, c] = a [b, c] a^{-1} [a, c] \quad \textrm{and} \quad [c, a b] = [c, a] a [c, b] a^{-1},
$$
we can write any commutator from  $[\Adj(X), \Adj(X)]$ as a product of commutators $ [\textswab{a}_x,  \textswab{a}_y]$ for $x, y \in X$. Since the number of such commutators is finite, assertion (2) follows. $\blacksquare$
\end{proof}

Observe that the dihedral quandle $\R_3$ is connected. Thus, by Proposition~\ref{fingen}(1), the abelianization of the adjoint group $\Adj(\R_3)$ is isomorphic to $\mathbb{Z}$. The same conclusion also follows from Example \ref{dihh}, where we proved that $\Adj(\R_3)=\langle \textswab{a}_0, \textswab{a}_1  \, \mid \, \textswab{a}_1 \textswab{a}_0 \textswab{a}_1 = \textswab{a}_0 \textswab{a}_1 \textswab{a}_0~\textrm{and}~\textswab{a}_0^2 = \textswab{a}_1^2\rangle.$
\para

Recall from Remark \ref{eta not injective} that the natural map  $i_X: X\to \Adj(X)$ is not injective in general. The following result shows a relation between the automorphism group of $X$ and the automorphism group of $\Adj(X)$ when $i_X$ is injective \cite[Proposition 4]{MR3948284}.

\begin{proposition}\label{inclus} 
Let $X$ be a quandle such that the natural map $i_X: X\to \Adj(X)$ is injective. If $\Aut_X \big(\Adj(X)\big)$ denotes the subgroup of $\Aut \big(\Adj(X) \big)$ consisting of automorphisms which keep $i_X(X)$ invariant, then $\Aut_X \big(\Adj(X) \big)\cong  \Aut(X)$.
\end{proposition}

\begin{proof}
Since $i_X$ is injective, for each $f\in \Aut_X \big(\Adj(X) \big)$, the map $f': X \to X$ given by $f'(x)= i_X^{-1}(f(\textswab{a}_x))$ is a well-defined bijective map. Moreover, if $x, y \in X$, then 
\begin{eqnarray*}
f'(x*y) &=& i_X^{-1} \big(f(\textswab{a}_{x*y})\big)\\
&=& i_X^{-1} \big(f(\textswab{a}_y) f(\textswab{a}_x) f(\textswab{a}_y)^{-1} \big)\\
&=& i_X^{-1} \big(f(\textswab{a}_x) * f(\textswab{a}_y) \big)\\
&=& i_X^{-1} \big(f(\textswab{a}_x) \big) * i_X^{-1} \big(f(\textswab{a}_y) \big)\\
&=& f'(x)*f'(y),
\end{eqnarray*}
 and hence $f' \in \Aut(X)$.  This defines a map $\phi:\Aut_X \big(\Adj(X\big)) \to \Aut(X)$ given by $\phi(f)=f'$, which is clearly a group homomorphism. We claim that $\phi$ is a bijection.
\para
If $\phi(f)=\phi(g)$ for $f,g\in \Aut \big(\Adj(X) \big)$, then $f'(x)=g'(x)$ for all $x\in X$. This implies that $f=g$ since $\Adj(X)$ is generated by $i_X(X)=\{\textswab{a}_x \,\mid \, x \in X \}$. For surjectivity of $\phi$, let $\alpha \in \Aut(X)$. Then the injective composition $i_X \alpha \, i_X^{-1}:i_X(X) \to \Adj(X)$ can be extended to a group homomorphism $f:F(X) \to \Adj(X)$, where $F(X)$ is the free group on $i_X(X)$. Further, for elements $x, y \in X$, we have
\begin{eqnarray*}
f(\textswab{a}_{x*y}) &=& i_X \alpha \, i_X^{-1}(\textswab{a}_{x*y})\\
&=& i_X \alpha(x*y)\\
&=& i_X \alpha(x)*i_X\alpha(y)\\
&=& i_X\alpha(y) ~i_X \alpha(x) ~i_X\alpha(y)^{-1}\\
&=& f(\textswab{a}_y) ~ f(\textswab{a}_x) ~f(\textswab{a}_y)^{-1}\\
&=& f(\textswab{a}_y \textswab{a}_x \textswab{a}_y^{-1}),
\end{eqnarray*}
and hence $f$ can be viewed as a group homomorphism $f: \Adj(X) \to \Adj(X)$. It is easy to see that $f$ is an automorphism of $\Adj(X)$ which keeps $X$ invariant such that $\phi(f)= \alpha$.  This completes the proof of the proposition. $\blacksquare$
\end{proof}

\begin{remark}{\rm 
If $i_X:X \to \Adj(X)$ is injective,  then by Proposition \ref{inclus}, the group $\Aut(X)$ is a subgroup of $\Aut \big(\Adj(X) \big)$, which is not necessarily normal. Indeed, if $\T_n=\{x_1,\dots,x_n\}$ is the trivial quandle, then the natural map $i_X$ is injective and $\Adj(\T_n)$ is the free abelian group of rank $n$. We have $\Aut(\T_n) \cong \Sigma_{n}$, which consists of all permutations of elements of $\T_n$. On the other hand, $\Aut \big(\Adj(\T_n) \big) \cong \GL(n, \mathbb{Z})$, which acts naturally on $\Adj(\T_n)\cong \mathbb{Z}^{n}$.  If $n>1$, then consider the automorphisms
$$
\varphi :\begin{cases}x_1\mapsto x_2,&\\
x_2\mapsto x_1,&\\
x_i\mapsto x_i& ~\textrm{for}~~i>2,
\end{cases}
\quad \textrm{and } \quad 
\psi: \begin{cases}
\textswab{a}_{x_1}\mapsto \textswab{a}_{x_1}+ \textswab{a}_{x_2},&\\
\textswab{a}_{x_i}\mapsto \textswab{a}_{x_i}&~\textrm{for}~~ i>1,
\end{cases}
$$
in the groups $\Aut(\T_n)$ and $\Aut \big(\Adj(\T_n) \big)$, respectively. Viewing $\Aut(\T_n)$ as a subgroup of $\Aut \big(\Adj(\T_n)\big)$, we see that the automorphism $\psi^{-1}\varphi\psi$ of $\Adj(\T_n)$ has the form
$$
\psi^{-1}\varphi\psi : \begin{cases}
\textswab{a}_{x_1} \mapsto \textswab{a}_{x_1},& \\
\textswab{a}_{x_2} \mapsto \textswab{a}_{x_1} - \textswab{a}_{x_2},& \\
\textswab{a}_{x_i} \mapsto \textswab{a}_{x_i}& ~\textrm{for}~~i > 2,
\end{cases}
$$
and hence does not lie in $\Aut(\T_n)$. Thus, $\Aut(\T_n)$ is not a normal subgroup of $\Aut \big(\Adj(\T_n) \big)$.}
\end{remark}

Next, we show that every quandle admits a homogeneous representation in terms of its adjoint group. For this, we first recall the quandle of Example \ref{coset quandle with auto}. The construction presented below generalizes this example.

\begin{proposition}
Let $G$ be a group, $\{z_i \mid i \in I\}$ a set of elements of $G$ and $H_i \le \C_G(z_i)$ a family of subgroups of $G$. Then the disjoint union $\bigsqcup_{i \in I} G/H_i$ of right cosets forms a quandle with respect to the binary operation
$$
H_i x * H_j y = H_i z_i^{-1} x y^{-1} z_j y
$$
for $x, y \in G$. We denote this quandle by $\bigsqcup_{i \in I}(G, H_i, z_i)$.
\end{proposition}

\begin{proof}
The proof is a direct check, and hence omitted. $\blacksquare$
\end{proof}

Observe that a quandle $X$ is connected if and only if $\Adj(X)$ acts transitively on $X$. The next result shows that every quandle admits a homogeneous representation \cite{MR0638121, MR0672410}.

\begin{theorem}\label{thm connected quandle is coset quandle}
The following statements hold:
\begin{enumerate}
\item  Let $X$ be a connected quandle,  $x_0 \in X$ a fixed element and $z = \textswab{a}_{x_0} \in\Adj(X)$. Let $H$ be the stabilizer of $x_0$ under the action of $\Adj(X)$.
Then the orbit map $\xi : \Adj(X) \to X$ given by $\xi(w)= x_0 \cdot w$ induces a quandle isomorphism $$ \big(\Adj(X), H,  z \big) \cong X.$$
\item Every quandle $X$ is representable as $\bigsqcup_{i \in I}(G, H_i, z_i)$ for some group $G$, a set $\{z_i \mid i \in I\}$ of elements of $G$ and a family of subgroups $H_i \le \C_G(z_i)$ of $G$.
\end{enumerate}
\end{theorem}

\begin{proof} For assertion (1), since $X$ is connected, the map $\xi$ descends to a bijection $\overline{\xi}: \Adj(X)/H\to X$. Hence, it is enough to show that $\xi$ is a quandle homomorphism. Indeed, we see that
\begin{eqnarray*}
\overline{\xi}(Hw_1 * Hw_2) &=& \overline{\xi} (H z^{-1} w_1 w_2^{-1} z w_2)\\
&=& x_0 \cdot (z^{-1} w_1 w_2^{-1} z w_2)\\
&=& \big(((x_0 \cdot z^{-1}) \cdot w_1w_2^{-1}) \cdot z \big) \cdot w_2\\
&=& \big((x_0 \cdot w_1w_2^{-1})* x_0 \big) \cdot w_2\\
&=& \big((x_0 \cdot w_1w_2^{-1})\cdot w_2 \big) * (x_0 \cdot w_2)\\
&=& (x_0 \cdot w_1) * (x_0 \cdot w_2)\\
&=& \overline{\xi}(Hw_1) * \overline{\xi}(Hw_2)
\end{eqnarray*}
for any $w_1, w_2 \in  \Adj(X)$.
\para

Take $G = \Adj(X)$ and write $X = \bigsqcup_{i \in I} X_i$ as a disjoint union of orbits, where $X_i$ is the orbit of the element $x_i \in X$. For each $i$, let  $z_i=\textswab{a}_{x_i} \in G$ and $H_i$ the stabilizer of $x_i$, giving rise to the quandle $(G, H_i, z_i)$. For each $i$, define a map $\xi_i : G \to X_i$ by $\xi(w) = x_i \cdot w$, which induce a bijection $\overline{\xi}_i:G/H_i \to X_i$. As in  assertion (1), the map $\bigsqcup_{i \in I} \overline{\xi}_i: \bigsqcup_{i \in I}(G, H_i, z_i) \to X$ is a quandle isomorphism, which proves assertion (2). $\blacksquare$
\end{proof}

\begin{remark}{\rm  
The homogeneous representation of a connected quandle $X$ as in Theorem \ref{thm connected quandle is coset quandle} is not unique. For example, we can replace $\Adj(X)$ by groups $\Inn(X)$ or $\Aut(X)$.}
\end{remark}

\bigskip
\bigskip


\section{Adjoint functors between categories of quandles and groups}\label{adjoint-functor-section}
Let $\mathcal{Q}$ and $\mathcal{G}$ denote the categories of quandles and groups, respectively.  Let $w = w(x,y)$ be a word in the free group $F(x,y)$ on generators $x$ and $y$. For any group $G$, the word $w$ defines a map $G\times G\to G$ given by $(g, h) \mapsto w(g, h)$  for $g, h \in G$. Setting $$g * h = w(g, h)$$ gives an algebraic structure $(Q_w(G), *)$. If $w(x,y) = y^{n} x y^{-n}$ for some integer $n$, then $Q_w(G)$ is the $n$-th conjugation quandle $\Conj_n(G)$, and if $w(x,y) = y x^{-1} y$, then $Q_w(G)$ is the core quandle $\Core(G)$. In fact, Proposition \ref{verbal rack quandle prop} shows that these are all the possibilities for $(Q_w(G), *)$ to be a quandle.
\para
For a word $w = w(x,y)\in F(x,y)$ and a quandle $X \in \mathcal{Q}$, we define the group
$$\Adj_w(X)= \big\langle \textswab{a}_x, ~x \in X \, \mid \, \textswab{a}_{x*y}=w(\textswab{a}_x, \textswab{a}_y),~~x, y \in X \big\rangle.$$
We note that if $w(x,y)=yxy^{-1}$, then $\Adj_w(X)$ is simply the adjoint group $\Adj(X)$ of $X$.

\begin{proposition}\label{verbal-adjoint}
Let $w(x,y) = y x^{-1} y$ or $y^{n} x y^{-n}$ for some $n \in \mathbb{Z}$. Then $\Adj_w: \mathcal{Q} \to \mathcal{G}$ is a functor that is left adjoint to the functor $Q_w:\mathcal{G} \to \mathcal{Q}$.
\end{proposition}

\begin{proof}
Let $X, Y \in \mathcal{Q}$ and $\phi \in \Hom_{\mathcal{Q}}(X, Y)$ a quandle homomorphism. Then we obtain a group homomorphism
$\Adj_w(\phi): F(X) \to \Adj_w(Y)$ given by $\Adj_w(\phi)(\textswab{a}_x)= \textswab{a}_{\phi(x)}$, where $F(X)$ is the free group on the set $\{\textswab{a}_x \, \mid \, x \in X\}$. Further, for $x, y \in X$, we have
\begin{eqnarray*}
\Adj_w(\phi)(\textswab{a}_{x*y})& =& \textswab{a}_{\phi(x)*\phi(y)}\\
& =& w(\textswab{a}_{\phi(x)}, \textswab{a}_{\phi(y)})\\
& =& \Adj_w(\phi) \big(w(\textswab{a}_x, \textswab{a}_y)\big),
\end{eqnarray*}
and hence we obtain a group homomorphism $\Adj_w(\phi): \Adj_w(X) \to \Adj_w(Y).$
\para

Conversely, let $G, H \in \mathcal{G}$ and $f \in \Hom_{\mathcal{G}} (G, H)$. Defining $Q_w(f): Q_w(G) \to Q_w(H)$ by $Q_w(f)(a)=f(a)$, we see that
\begin{eqnarray*}
Q_w(f)(a*b)& =& f(a*b)\\
& =& f\big(w(a, b)\big)\\
& =& w\big(f(a), f(b)\big)\\
& =& f(a)*f(b)\\
& =& Q_w(f)(a)*Q_w(f)(b)
\end{eqnarray*}
for all $a, b \in G$. Thus, $Q_w(f) \in \Hom_{\mathcal{Q}} \big(Q_w(G), Q_w(H)\big)$. A direct check shows that both $\Adj_w$ and $Q_w$ are functors.
\para

Let $X \in \mathcal{Q}$, $G \in \mathcal{G}$ and $\phi \in \Hom_{\mathcal{Q}} \big(X, Q_w(G)\big)$ a quandle homomorphism. Define $\widetilde{\phi}: F(X) \to G$ by setting $\widetilde{\phi}(\textswab{a}_x)= \phi(x)$. Then, for $x, y \in X$, we have
\begin{eqnarray*}
\widetilde{\phi}(\textswab{a}_{x*y})& =& \phi(x*y)\\
& =& \phi(x)*\phi(y)\\
& =& w\big(\phi(x), \phi(y)\big)\\
& =& \widetilde{\phi}\big(w(\textswab{a}_{x}, \textswab{a}_{y})\big),
\end{eqnarray*}
and hence $\widetilde{\phi} \in \Hom_{\mathcal{G}} \big(\Adj_w(X), G\big)$. Similarly, let $f \in \Hom_{\mathcal{G}} \big(\Adj_w(X), G\big)$ be a group homomorphism. Define $\widehat{f}: X \to Q_w(G)$ by setting $\widehat{f}(x)= f(\textswab{a}_x)$. Then, for $x, y \in X$, we have
\begin{eqnarray*}
\widehat{f}(x*y)& =& f(\textswab{a}_{x*y})\\
& =& f\big(w(\textswab{a}_x, \textswab{a}_y)\big)\\
& =& w\big(f(\textswab{a}_x), f(\textswab{a}_y)\big)\\
& =& w\big(\widehat{f}(x), \widehat{f}(y)\big)\\
& =& \widehat{f}(x) * \widehat{f}(y)
\end{eqnarray*}
and hence $\widehat{f} \in \Hom_{\mathcal{Q}} \big(X, Q_w(G)\big)$. It follows easily that $\widehat{\widetilde{\phi}}=\phi$ and $\widetilde{\widehat{f}}=f$, which estabishes that $\Adj_w$ is left adjoint to $Q_w$.  $\blacksquare$
\end{proof}

\begin{corollary}
For any quandle $X$ and any group $G$, there is a natural bijection
$$
\Hom_{\mathcal{Q}} \big(X, \Conj(G) \big) \cong \Hom_{\mathcal{G}} \big(\Adj(X), G \big).  
$$
\end{corollary}
\para

Let $\mathcal{Q}'$ denote the category of pairs $(X, \phi)$, where $X \in \mathcal{Q}$ and $\phi \in \Aut(X)$. If $(X_1, \phi_1), (X_2, \phi_2) \in \mathcal{Q}'$, then a morphism from $(X_1, \phi_1)$ to $(X_2, \phi_2)$ is a quandle homomorphism $\psi:X_1 \to X_2$ such that $$\phi_2 \, \psi= \psi \, \phi_1.$$ In case of groups, the category $\mathcal{G}'$ is defined analogously. Each object $(G, f) \in \mathcal{G}'$ gives an object $\big(\Alex_f(G),  \Alex(f) \big)$, where $\Alex(f)$ acts as $f$. On the other hand, for each $(X, \phi) \in \mathcal{Q}'$, we define the group
$$\Adj_\phi(X)= \big\langle \textswab{a}_x, ~x \in X \, \mid \, \textswab{a}_{x*y}=\textswab{a}_{\phi(x)}\textswab{a}_{\phi(y)}^{-1}\textswab{a}_{y}\quad \textrm{for all} \quad x, y \in X \big\rangle.$$
Further, the map $\phi$ extends to an automorphism $\Adj(\phi)$ of the free group $F(X)$. Moreover, for $x, y \in X$, we have 
$$ \Adj(\phi)(\textswab{a}_{x*y}) = \textswab{a}_{\phi(x)*\phi(y)}= \textswab{a}_{\phi^2(x)}\textswab{a}_{\phi^2(y)}^{-1}\textswab{a}_{\phi(y)}=
 \Adj(\phi) \big(\textswab{a}_{\phi(x)}\textswab{a}_{\phi(y)}^{-1}\textswab{a}_{y} \big),
$$
and hence $\big(\Adj_\phi(X), \Adj(\phi)\big) \in \mathcal{G}'$. It follows easily that both $\Alex$ and $\Adj$ are functors. Imitating the proof of Proposition \ref{verbal-adjoint}, we obtain the following result.

\begin{proposition}\label{alexander-adjoint}
The functor $\Adj: \mathcal{Q}' \to \mathcal{G}'$ is left adjoint to the functor $\Alex:\mathcal{G}' \to \mathcal{Q}'$.
\end{proposition}
\bigskip
\bigskip


\section{Quandles that embed into their adjoint groups}
This section is motivated by the question of which quandles can be embedded into the conjugation quandles of groups. We begin with a generalisation of the construction of free quandles and free racks. 
\para

Let $G$ be a group and $A$ a subset of $G$. On the set $A\times G$, we define the binary operation
$$
(a, u) * (b, v) = (a, v b v^{-1} u)
$$
for $a, b \in A$ and $u, v \in G$. It is not difficult to check that the algebraic structure $(A\times G, *) $ forms a rack, which we denote by $R(G,A)$ and refer as a {\it $(G,A)$-rack}. Let $Q(G,A)$ be the quotient of $R(G,A)$ by the equivalence relation generated by $(a,uw)\sim(a,u)$ for all $a \in A$, $u \in G$ and $w\in \C_G(a)=\{x\in G \, \mid \, xa=ax\}$. Then the binary operation
$$
[(a, u)] * [(b, v)] = [(a, v b v^{-1} u)],
$$
where $[(a,u)]$ denotes the equivalence class of $(a,u)$, turns $Q(G,A)$ into a quandle referred to as a {\it $(G,A)$-quandle}. For the sake of simplicity, we write $(a,u)$ instead of $[(a,u)]$ to denote an element of  $Q(G,A)$. Note that, if $F(X)$ is the free group with free basis $X$, then $Q \big(F(X),X  \big)=FQ(X)$, the free quandle on $X$. The following observation follows from the definition.
\para

A direct computation yields the following result.

\begin{lemma}\label{simple2} 
If $x,y\in Q(G,A)$ are such that $x*y=x$, then $y*x=y$.
\end{lemma}

It seems that, similar to free quandles, the quandle $Q(G,A)$ is isomorphic to the subquandle of $\Conj(G)$ consisting of union of conjugacy classes of elements from $A$. However, it is not true, in general. For example, if $G=F_2=\langle x,y\rangle$ is a free group and $A=\{x,y,y^{-1}xy\}$, then the quandle $Q(G,A)$ has three orbits under the action of $\Inn  \big(Q(G,A) \big)$, but the subquandle of $\Conj(G)$ consisting of union of conjugacy classes of $x$ and $y$ has only two orbits under the action of its group of inner automorphisms. However, if the elements of $A$ are pairwise non-conjugate, then we have the following result \cite[Proposition 4.3]{MR4129183}.

\begin{proposition} \label{ga subquandle of conj}
Let $G$ be a group and $A$ a subset of $G$ such that the elements of $A$ are pairwise non-conjugate in $G$. Then $Q(G,A)$ is isomorphic to the subquandle of $\Conj(G)$ consisting of complete conjugacy classes of elements from $A$.
\end{proposition}

\begin{proof} 
Let $X$ denote the subquandle of $\Conj(G)$ consisting of union of conjugacy classes of elements from $A$ and consider the map $\varphi:Q(G,A)\to X$ given by $\varphi\big((a,u)\big)=u au^{-1}$. If $w\in X$, then there exist elements $b\in A$ and $x\in G$ such that $w=x bx^{-1}$. Thus, $\varphi\big((b,x)\big)=w$, and $\varphi$ is surjective. Suppose that $u au^{-1}=\varphi \big((a,u)\big)=\varphi\big((b,v)\big)=v bv^{-1}$ for $a,b\in A$ and $u,v\in G$. Since elements from $A$ are pairwise non-conjugate, it follows that $a=b$ and $vu^{-1}\in \C_G(a)$. Thus, $(a,u)=(b,v)$ in $Q(G,A)$ and $\varphi$ is injective. Further, direct computations give
\begin{eqnarray*}
\varphi \big((a,u)*(b,v)\big)&=& \varphi(a,v b v^{-1}u)\\
&=&(v b v^{-1}u) a (v b v^{-1}u)^{-1}\\
&=& (v b v^{-1}) (u a u^{-1}) (v b v^{-1})^{-1}\\
&=& \big(\varphi(b,v)\big) \big(\varphi(a,u)\big) \big(\varphi(b,v)\big)^{-1}\\
&=& \varphi \big((a,u)\big)* \varphi\big((b,v)\big),
\end{eqnarray*}
and hence $\varphi$ is an isomorphism. $\blacksquare$
\end{proof}

\begin{corollary} 
Let $G$ be a group and $A$ a set of representatives of conjugacy classes of $G$. Then $Q(G,A) \cong \Conj(G)$.
\end{corollary}

The following result demonstrates that a wide variety of quandles can be realized as $(G,A)$-quandles \cite[Theorem 4.5]{MR4129183}.

\begin{theorem}\label{theorem quandle is ga quandle}
Let $X$ be a quandle such that the natural map $i_X:X\to \Adj(X)$ is injective, and $A$ be a complete set of representatives of orbits of $X$. Then $X \cong Q  \big(\Adj(X),i_X(A)  \big)$.
\end{theorem}

\begin{proof} 
If $a\in X$, then there exist an element $x\in A$ and an inner automorphism $\varphi \in \Inn(X)$ such that $a=\varphi(x)$. Note that $\varphi$ can be written in the form $\varphi=S_{x_n}^{\varepsilon_n}S_{x_{n-1}}^{\varepsilon_{n-1}}\cdots S_{x_1}^{\varepsilon_1}$ for some $x_i\in X$ and $\varepsilon_i \in\{\pm1\}$. Thus, we can write
\begin{align}\label{canon}
a=x*^{\varepsilon_1}x_1*^{\varepsilon_2}\dots *^{\varepsilon_n}x_n.
\end{align}
Define $\psi:X\to Q  \big(\Adj(X),i_X(A) \big)$ by
$$\psi(a)= (\textswab{a}_x, \,\textswab{a}_{x_n}^{\varepsilon_n} \cdots \textswab{a}_{x_1}^{\varepsilon_1}).$$
We claim that $\psi$ is an isomorphism of quandles. First, we prove that $\psi$ is well-defined. Let $a \in X$ be such that
\begin{equation}\label{korre} 
x*^{\varepsilon_1}x_1*^{\varepsilon_2}\dots *^{\varepsilon_n}x_n=a=y*^{\epsilon_1}y_1*^{\epsilon_2}\dots *^{\epsilon_k}y_k
\end{equation}
for some $x,y\in A$, $x_i,y_j\in X$ and $\varepsilon_i,\epsilon_j \in\{\pm 1\}$. Since all the elements from $A$ belong to different orbits, we must have $x=y$. The equality \eqref{korre} implies that
$$
x*^{\varepsilon_1}x_1*^{\varepsilon_2}\dots *^{\varepsilon_n}x_n*^{-\epsilon_k}y_k*^{-\epsilon_{k-1}}\dots*^{-\epsilon_1}y_1=x.
$$
Applying the map $i_X$, we see that the element $w=\textswab{a}_{y_1}^{-\epsilon_1}\textswab{a}_{y_2}^{-\epsilon_2} \cdots \textswab{a}_{y_{k-1}}^{-\epsilon_{k-1}}\textswab{a}_{y_k}^{\epsilon_k}\textswab{a}_{x_n}^{\varepsilon_n}\textswab{a}_{x_{n-1}}^{\varepsilon_{n-1}} \cdots \textswab{a}_{x_2}^{\varepsilon_2}\textswab{a}_{x_1}^{\varepsilon_1}$ belongs to the centraliser $\C_{\Adj(X)}(\textswab{a}_x)$. Thus, we have
$$\psi(a)=(\textswab{a}_y,\textswab{a}_{y_k}^{\epsilon_k}\cdots \textswab{a}_{y_1}^{\epsilon_1})=(\textswab{a}_x,\textswab{a}_{y_k}^{\epsilon_k}\cdots \textswab{a}_{y_1}^{\epsilon_1})=(\textswab{a}_x, \textswab{a}_{y_k}^{\epsilon_k}\cdots \textswab{a}_{y_1}^{\epsilon_1}w)=(\textswab{a}_x, \textswab{a}_{x_n}^{\varepsilon_n} \cdots \textswab{a}_{x_1}^{\varepsilon_1})=\psi(a),$$
and the map $\psi$ is well-defined.
\para
Next, we prove that $\psi$ is bijective. Let $(\textswab{a}_x,w) \in Q  \big(\Adj(X),i_X(A) \big)$, where $\textswab{a}_x\in i_X(A)$ and $w\in \Adj(X)$. Since $i_X(X)=\{ \textswab{a}_x \, \mid \, x \in X\}$ generates $\Adj(X)$, we have $w=\textswab{a}_{x_n}^{\varepsilon_n}\cdots \textswab{a}_{x_1}^{\varepsilon_1}$ for some $x_i\in X$ and $\varepsilon_i \in\{\pm1\}$. Setting $a=x*^{\varepsilon_1}x_1*^{\varepsilon_2}\dots *^{\varepsilon_n}x_n$, we see that $\psi(a)=(\textswab{a}_x, \textswab{a}_{x_n}^{\varepsilon_n} \cdots \textswab{a}_{x_1}^{\varepsilon_1})$, and $\psi$ is surjective. For injectivity, consider two elements
$$ a=x*^{\varepsilon_1}x_1*^{\varepsilon_2}\dots *^{\varepsilon_n}x_n \quad \textrm{and} \quad  b=y*^{\epsilon_1}y_1*^{\epsilon_2}\dots *^{\epsilon_k}y_k,$$
where $x,y\in A$, $x_i, y_j\in X$ and $\varepsilon_i,\epsilon_i\in\{\pm 1\}$, such that
\begin{equation}\label{1inject}
(\textswab{a}_x, \textswab{a}_{x_n}^{\varepsilon_n} \cdots \textswab{a}_{x_1}^{\varepsilon_1})=\psi(a)=\psi(b)=(\textswab{a}_y, \textswab{a}_{y_k}^{\epsilon_k} \cdots \textswab{a}_{y_1}^{\epsilon_1}).
\end{equation}
It follows from the construction of $Q(G,A)$ that $\textswab{a}_x=\textswab{a}_y$. Since $i_X$ is injective, we get $x=y$. Equation \eqref{1inject} implies that the element $w=\textswab{a}_{y_1}^{-\epsilon_1} \textswab{a}_{y_2}^{-\epsilon_2} \cdots \textswab{a}_{y_{k-1}}^{-\epsilon_{k-1}} \textswab{a}_{y_k}^{\epsilon_k} \textswab{a}_{x_n}^{\varepsilon_n} \textswab{a}_{x_{n-1}}^{\varepsilon_{n-1}} \cdots \textswab{a}_{x_2}^{\varepsilon_2} \textswab{a}_{x_1}^{\varepsilon_1}$ belongs to $\C_{\Adj(X)}(\textswab{a}_x)$. Consider the element
$$c=x*^{\varepsilon_1}x_1*^{\varepsilon_2}\dots *^{\varepsilon_n}x_n*^{-\epsilon_k}y_k*^{-\epsilon_{k-1}}\dots *^{-\epsilon_1}y_1.$$
Applying $i_X$ on both the sides give $i_X(c)=w i_X(x) w^{-1}=i_X(x)$. Since $i_X$ is injective, we get $c=x=y$. This gives
$$b=y*^{\epsilon_1}y_1*^{\epsilon_2}\dots *^{\epsilon_k}y_k=c*^{\epsilon_1}y_1*^{\epsilon_2}\dots *^{\epsilon_k}y_k=x*^{\varepsilon_1}x_1*^{\varepsilon_2}\dots *^{\varepsilon_n}x_n=a$$
and $\psi$ is injective.
\para
Finally, we prove that $\psi$ is a homomorphism. Let $a=x*^{\varepsilon_1}x_1*^{\varepsilon_2}\dots *^{\varepsilon_n}x_n$ and $b=y*^{\epsilon_1}y_1*^{\epsilon_2}\dots *^{\epsilon_k}y_k$ be elements of $X$. Then $\psi(a)=(\textswab{a}_x, \textswab{a}_{x_n}^{\varepsilon_n} \cdots \textswab{a}_{x_1}^{\varepsilon_1})$, $\psi(b)=(\textswab{a}_y, \textswab{a}_{y_k}^{\epsilon_k} \cdots \textswab{a}_{y_1}^{\epsilon_1})$ and 
\begin{equation}\label{oneside}
\psi(a)*\psi(b)=(\textswab{a}_x, \textswab{a}_{y_k}^{\epsilon_k} \cdots \textswab{a}_{y_1}^{\epsilon_1} \textswab{a}_y \textswab{a}_{y_1}^{-\epsilon_1} \cdots \textswab{a}_{y_k}^{-\epsilon_k} \textswab{a}_{x_n}^{\varepsilon_n} \cdots \textswab{a}_{x_1}^{\varepsilon_1}).
\end{equation}

Note that the right-distributivity axiom can be rewritten as
\begin{equation}\label{hope}
x*^{-1}z*z*y*z=(x*z)*(y*z)
\end{equation}
for all $x, y, z \in X$. Applying \eqref{hope} repeatedly, we obtain
\begin{eqnarray*}
a*b&= & a*(y*^{\epsilon_1}y_1*^{\epsilon_2}\dots *^{\epsilon_k}y_k)\\
 &= & \big((a*^{-\epsilon_k}y_k)*(y*^{\epsilon_1}y_1*^{\epsilon_2}\dots *^{\epsilon_{k-1}}y_{k-1})\big)*^{\epsilon_k} y_k\\
 &\vdots &\\
&=&a*^{-\epsilon_k}y_k*^{-\epsilon_{k-1}}\dots *^{-\epsilon_1}y_1*y*^{\epsilon_1}y_1*^{\epsilon_2}\dots*^{\epsilon_k}y_k.
\end{eqnarray*}
Applying $i_X$ to the preceding equality and comparing with \eqref{oneside} proves that $\psi$ is a homomorphism.  $\blacksquare$
\end{proof}

\begin{remark}{\rm 
Note that Example \ref{example Dehn quandle} gives the following similar but different construction of quandles from groups. Let $G$ be a group and $A$ a non-empty subset of $G$. Let $A^G$ be the set of all conjugates of elements of $A$ in $G$. Then the {\it Dehn quandle} $\mathcal{D}(A^G)$ of $G$ with respect to $A$ is defined as the set $A^G$ equipped with the binary operation  $x*y=yxy^{-1}$ of conjugation.
\para
By Proposition \ref{ga subquandle of conj}, the quandles $Q(G, A)$ and $\mathcal{D}(A^G)$ are the same if and only if elements of $A$ are pairwise non-conjugate in $G$. Thus, for a fairly large class of groups including knot groups, mapping class groups, Artin groups and Coxeter groups, the construction in Example \ref{example Dehn quandle} differs from the one presented in this section.}
\end{remark}
\para

The next result gives a complete answer to the question posed in the beginning of this section \cite[Proposition 3.9]{MR4669143}.

\begin{theorem}\label{injectivity-dehn-eta}
The following statements are equivalent for any quandle $X$:
\begin{enumerate}
\item $X$ embeds in $\Conj(H)$ for some group $H$.
\item The natural map $i_X: X \to \Adj(X)$ is injective.
\item $X \cong \mathcal{D}(S^G)$ for some group $G$ and a generating set $S$ of $G$.
\end{enumerate}
\end{theorem}

\begin{proof}
For (1) $\Rightarrow$ (2), suppose that $f: X \hookrightarrow \Conj(H)$ is an embedding for some group $H$. Then, by Proposition \ref{Prop: natural map from quandle to its adjoint group}(2), there is a unique group homomorphism $\bar{f}: \Adj(X) \to H$ such that $f=\bar{f} ~i_X$. Thus, $i_X$ is injective. The implication (2) $\Rightarrow$ (1) is obvious.
\para
For (2) $\Rightarrow$ (3), suppose that $X$ is a quandle for which the map $i_X$ is injective. Let $A$ be the set of representatives of orbits (connected components) of $X$. We claim that elements of $i_X(A)$ are pairwise non-conjugate in $\Adj(X)$. Let $x$ and $y$ be distinct elements of  $A$ such that $\textswab{a}_x=g \textswab{a}_y g^{-1}$, where $g= \textswab{a}_{x_1}^{\epsilon_1} \textswab{a}_{x_2}^{\epsilon_2} \cdots \textswab{a}_{x_k}^{\epsilon_k}$ for $x_i \in X$ and $\epsilon_i \in \{ 1, -1\}$. Then we have
$$\textswab{a}_x=\textswab{a}_{x_1}^{\epsilon_1} \textswab{a}_{x_2}^{\epsilon_2} \cdots \textswab{a}_{x_k}^{\epsilon_k} \textswab{a}_y \textswab{a}_{x_k}^{-\epsilon_k} \cdots \textswab{a}_{x_2}^{-\epsilon_2} \textswab{a}_{x_1}^{-\epsilon_1} = \textswab{a}_{y *^{\epsilon_k} x_k *^{\epsilon_{k-1}} x_{k-1} *^{\epsilon_{k-2}} \cdots *^{\epsilon_1} x_1}.$$
Since $i_X$ is injective, it follows that $x=y *^{\epsilon_k} x_k *^{\epsilon_{k-1}} x_{k-1}*^{\epsilon_{k-2}} \cdots *^{\epsilon_1} x_1$, that is, $x$ and $y$ are in the same orbit, a contradiction. It now follows from Proposition \ref{ga subquandle of conj} and Theorem \ref{theorem quandle is ga quandle} that $X\cong \mathcal{D} \big(i_X(A)^{\Adj(X)}\big)$. If $S$ is a generating set for $X$, then $i_X(S)$ is a generating set for $\Adj(X)$. Further, $S$ intersects every orbit of $X$, and therefore we can choose a representative of an orbit from $S$. In other words, the set $A$ can be chosen to be a subset of $S$. For this choice of $A$, elements of $i_X(S)$ are conjugates of elements of $i_X(A)$ in $\Adj(X)$. Thus, $i_X(S)^{\Adj(X)}=i_X(A)^{\Adj(X)}$ as sets, and hence $X\cong \mathcal{D}\big(i_X(S)^{\Adj(X)}\big)$.
\para
For (3) $\Rightarrow$ (2), suppose that $X \cong \mathcal{D}(S^G)$ for some group $G$ and a generating set $S$ of $G$. Let $\tau: \mathcal{D}(S^G) \to \Conj(G)$ be the natural embedding and  $i_X: \mathcal{D}(S^G) \to \Conj \big(\Adj(\mathcal{D}(S^G))\big)$ the natural quandle homomorphism. Define $\phi: \Adj \big(\mathcal{D}(S^G)\big) \to G$ by $\phi(\textswab{a}_x)= x$ for $x \in \mathcal{D}(S^G)$. Since $$\phi(\textswab{a}_{x*y})=\phi(\textswab{a}_{yxy^{-1}})=yxy^{-1}=\phi(\textswab{a}_y\textswab{a}_x\textswab{a}_y^{-1})$$ for all $x, y \in \mathcal{D}(S^G)$, it follows that $\phi$ is a group homomorphism, and hence induces homomorphism of conjugation quandles. Since $\phi ~i_X=\tau$, it follows that $i_X$ is injective.  $\blacksquare$
\end{proof}

\begin{remark}{\rm 
Given a quandle $X$, there are two natural maps $i_X: X \to \Adj(X)$ and $S: X \to \Inn(X)$. The universal property of $\Adj(X)$ implies that if $S$ is injective, then $i_X$ is also injective. However, the converse is not true. For instance, if $X$ is a trivial quandle with more than one element, then $i_X$ is injective, but $S$ is not.}
\end{remark}

A knot is said to be \index{prime knot}{\it prime} if it cannot be written as a connected sum of two non-trivial knots. The following result of Ryder provides an elegant characterisation of prime knots \cite[Corollary 3.6]{MR1388194}.
 
\begin{theorem}\label{prime-quandles-injective-in-adjoint-group}
A knot $K$ in $\mathbb{S}^3$ is prime if and only if its knot quandle $Q(K)$ embeds into its adjoint group $\Adj \big(Q(K)\big)=G(K)$.
\end{theorem}
\para

Constructing families of quandles that embed into their adjoint groups is of particular interest. In this context,  Bardakov and Nasybullov proved the following result \cite[Lemma 7.1]{MR4129183}.

\begin{proposition}\label{embedding-of-quandles-into-adjoint-group}
Let $X$ and  $Y$  be quandles. If the natural maps $X \rightarrow \Adj(X)$ and $Y \rightarrow \Adj(Y)$ are injective, then the natural map $X \star Y \rightarrow \Adj(X) \star \Adj(Y)$ is injective.
\end{proposition}
 
\begin{proof}
By Proposition \ref{presentation-of-adjoint-group-of-free-product-of-quandles}, we have $\Adj(X \star Y) \cong  \Adj(X) \hexstar \Adj(Y)$. Let $i_X: X \to \Adj(X)$, $i_Y: Y \to \Adj(Y)$ and $i_{X\star Y}: X \star Y \to \Adj(X) \hexstar \Adj(Y)$  be natural maps.  Let $a, b\in X\star Y$ be such that $i_{X\star Y}(a)=i_{X\star Y}(b)$. The elements $a,b$ can be written in the forms $a=a_1*^{\varepsilon_2}\dots *^{\varepsilon_n} a_n$ and $b=b_1*^{\mu_2}\dots *^{\mu_k} b_k$, where each $a_i,b_j$ belongs to either $X$ or $Y$ and $\varepsilon_2,\dots,\varepsilon_n,\mu_2,\dots,\mu_k\in \{1, -1\}$. Since there is no non-trivial element in $\Adj(X)$ that is conjugate to an element in $\Adj(Y)$, it follows that $a_1,b_1$ both belong to either $X$ or $Y$. Without loss of generality, we can assume that $a_1,b_1\in X$. Setting
$$c:=a_1*^{\varepsilon_2}\dots *^{\varepsilon_n} a_n*^{-\mu_k}b_k*^{-\mu_{k-1}}\dots *^{-\mu_2} b_2,$$
we see that $i_{X\star Y}(c)=i_{X\star Y}(b_1)$.  By definition of $i_{X\star Y}$, we have
$$(\textswab{a}_{b_2}^{-\mu_2} \cdots \textswab{a}_{b_k}^{-\mu_k} \textswab{a}_{a_n}^{\varepsilon_n}\cdots \textswab{a}_{a_3}^{\varepsilon_3} \textswab{a}_{a_2}^{\varepsilon_2}) \textswab{a}_{a_1} (\textswab{a}_{b_2}^{-\mu_2} \cdots \textswab{a}_{b_k}^{-\mu_k} \textswab{a}_{a_n}^{\varepsilon_n}\cdots \textswab{a}_{a_3}^{\varepsilon_3} \textswab{a}_{a_2}^{\varepsilon_2})^{-1}=i_{X\star Y}(c)=i_{X\star Y}(b_1)=\textswab{a}_{b_1}.$$
It follows that the element $\textswab{a}_{b_2}^{-\mu_2} \cdots \textswab{a}_{b_k}^{-\mu_k} \textswab{a}_{a_n}^{\varepsilon_n}\cdots \textswab{a}_{a_3}^{\varepsilon_3} \textswab{a}_{a_2}^{\varepsilon_2}$ belongs to $\Adj(X)$, since two elements from $\Adj(X)$ are conjugated in $\Adj(X)\ast \Adj(Y)$ if and only if they are conjugated in $\Adj(X)$. Thus, the element can be written as
$$\textswab{a}_{b_2}^{-\mu_2} \cdots \textswab{a}_{b_k}^{-\mu_k} \textswab{a}_{a_n}^{\varepsilon_n}\cdots \textswab{a}_{a_3}^{\varepsilon_3} \textswab{a}_{a_2}^{\varepsilon_2}=\textswab{a}_{x_1}^{\xi_1}\dots \textswab{a}_{x_t}^{\xi_t}$$
for some $x_1,\dots, x_t\in X$ and $\xi_1,\dots,\xi_t\in\{1, -1\}$. Hence, we have the equality
$$c=a_1*^{\varepsilon_2}\dots *^{\varepsilon_n} a_n*^{-\mu_k}b_k*^{-\mu_{k-1}}\dots *^{-\mu_2} b_2=a_1*^{\xi_t}x_t*^{\xi_{t-1}}\dots*^{\xi_1}x_1,$$
which implies that $c\in X$. Since $i_X(c)=i_X(b_1)$ and $i_X$ is injective, we conclude that $c=b_1$, and hence $a=b$. $\blacksquare$
\end{proof}
\bigskip
\bigskip


\section{Quandles with abelian adjoint groups}

The purpose of this section is to classify finite quandles with abelian adjoint groups. We begin with the following result of Lebed and Vendramin \cite[Theorem 7.15]{MR3974961}. 

\begin{theorem}\label{thm Lebed Vendramin}
The following statements are equivalent for a finite quandle $X$:
\begin{enumerate}
\item $\Adj(X)$ has no torsion.
\item $\Adj(X)$ is free abelian.
\item $\Adj(X)$ is abelian.
\end{enumerate}

\end{theorem}
\begin{proof}
Let $(X, *)$ be the given quandle. The implications $(2) \Rightarrow (1)$  and $(2) \Rightarrow (3)$ are immediate. Assume that $\Adj(X)$ has no torsion. By Proposition \ref{adj inn central extension}, we have the central extension of groups
$$
1 \to \ker(\zeta) \to \Adj(X) \to \Inn(X) \to 1.
$$
Since $X$ is finite, it follows that the center of $\Adj(X)$ is of finite index. Then, by Schur's theorem \cite[Theorem 5.32]{MR1307623}, the commutator subgroup of $\Adj(X)$ is finite. Since $\Adj(X)$ has no torsion, this subgroup is trivial, and hence $\Adj(X)$ is free abelian. This proves  $(1) \Rightarrow (2)$.
\para
Now, assume that $\Adj(X)$ is abelian. Define an equivalence relation $\sim$ on $X$ by setting $x \sim y$ if and only if $\textswab{a}_x=\textswab{a}_y$ in $\Adj(X)$. If $X/\!\sim$ is the set of equivalence classes, then the binary operation $[x]*' [y]:=[x*y]$ turns $X/\!\sim$ into a quandle. Further, by the construction of  $X/\!\sim$, it follows that $\Adj(X/\!\sim) \cong \Adj(X)$. For any $x,y\in X$, we have $\textswab{a}_x \textswab{a}_y=\textswab{a}_y \textswab{a}_x=\textswab{a}_{x * y} \textswab{a}_y$ in $\Adj(X)$, and hence $x \sim x * y$. Thus, the induced  quandle $(X/\!\sim,*')$ is trivial. Since  $\Adj(X/\!\sim) \cong \Adj(X)$, the triviality of $X/\!\sim$ implies that $\Adj(X/\!\sim)$, and hence $\Adj(X)$ is free abelian. This proves $(3) \Rightarrow (2)$. $\blacksquare$
\end{proof}

Theorem \ref{thm Lebed Vendramin} shows that adjoint groups of finite quandles exhibit the following dichotomy:
\begin{enumerate}
\item either they are free abelian of rank $r = |\mathcal{O}(X)|$, where $\mathcal{O}(X)$ is the set of orbits of $X$,
\item or they are non-abelian and have torsion.
\end{enumerate}

It is natural to ask which quandles lie in the first class. The condition $\Adj(X)\cong \mathbb{Z}$ is easily seen to be equivalent to $X$ being a one element quandle. Quandles with $\Adj(X) \cong \mathbb{Z}^2$ were first characterised in \cite{MR4129183}. This work was extended in \cite{MR4116819} to characterise all finite quandles $X$ with $\Adj(X) \cong \mathbb{Z}^r$ for arbitrary $r$, which we present in this section. We first show that quandles with abelian adjoint groups are necessarily abelian, that is, they satisfy the condition
\begin{equation}\label{equation abelian quandle}
(a * b) * c = (a * c) * b.
\end{equation}
\para

We fix a positive integer $r \geq 2$, and consider a collection 
\begin{equation}\label{parameter matrix}
M=(m_{i,j})_{1 \leq j \leq i <r}
\end{equation}
of $\frac{r(r-1)}{2}$ integers, where $1 \leq m_{i,i}$ and  $0 \leq m_{j,i} < m_{i,i}$ for $i < j$. We think of $M$ as a lower triangular matrix of size $r-1$, and associate an abelian group $G(M)$ by defining
$$G(M) = \big\langle x_1, \ldots, x_{r-1} ~|~ x_ix_j=x_jx_i, \quad x_1^{m_{i,1}}x_2^{m_{i,2}}\cdots x_i^{m_{i,i}} =1 \quad \textrm{for} \quad 1 \le i, j \le r-1 \big\rangle.$$
Using the notations $x_0:=1$ and $m_i:=m_{i,i}$, we see that the group $G(M)$ is finite abelian of order $m_{1}m_{2}\cdots m_{r-1}$.
\para

For instance, $r=2$ gives the cyclic group of order $m_{1}$. For $r=4$ and 
$M=\begin{psmallmatrix}
m_1& 0& 0\\ 
m_{2,1}&m_{2}&0&\\
m_{3,1}&m_{3,2}&m_{3}
\end{psmallmatrix}$, 
we get commuting generators $x_1, x_2, x_3$ subject to relations
$$ x_1^{m_{1}}=1, \quad x_1^{m_{2,1}}x_2^{m_{2}}=1, \quad x_1^{m_{3,1}}x_2^{m_{3,2}}x_3^{m_{3}}=1.
$$

We take $r$ collections $M^{(1)}$, \ldots, $M^{(r)}$ as above, and consider the disjoint union
$$Q(M^{(1)}, \ldots, M^{(r)}):=G(M^{(1)})\sqcup \cdots \sqcup G(M^{(r)}).$$
We denote the generator $x_i$ of $G(M^{(j)})$ by $x_i^{(j)}$, and equip $Q(M^{(1)}, \ldots, M^{(r)})$ with a binary operation $*$ as follows. For any $a^{(i)} \in G(M^{(i)})$ and $b^{(i+k)}\in G(M^{(i+k)})$, where $0 \leq k < r$ and the sum $i+k$ is considered modulo $r$, we put
$$a^{(i)} * b^{(i+k)} = a^{(i)}x_k^{(i)} \in G(M^{(i)}).$$
Since $x_0^{(i)}=1$, we have $a^{(i)} * b^{(i)} = a^{(i)}$ for each $i$. It is not difficult to prove the following result  \cite[Proposition 2.1]{MR4116819}.

\begin{proposition} \label{abelian quandle from parameters}
The quandle $\big(Q(M^{(1)}, \ldots, M^{(r)}), * \big)$ is abelian with $r$ components $G(M^{(1)}), \ldots, G(M^{(r)})$ as its orbits. 
\end{proposition}

\begin{proof}
The first two quandle axioms, and the assertion about the orbits, are clear from the construction. Moreover, the groups  $G(M^{(i)})$ are abelian, and hence all right $*$-actions commute. Consequently, we obtain the condition  \eqref{equation abelian quandle}. By construction, all elements from the same orbit $G(M^{(i)})$ of $Q(M^{(1)}, \ldots, M^{(r)})$ have the same right $*$-actions. Hence, we obtain
$$(a * b) * c = (a * c) * b = (a * c) * (b * c)$$
 for all $a,b,c \in Q(M^{(1)}, \ldots, M^{(r)})$, which is the third quandle axiom.
  $\blacksquare$
\end{proof}

\begin{definition}
A quandle $X$ such that $X \cong Q(M^{(1)}, \ldots, M^{(r)})$ is called a \index{filtered-permutation quandle} {\it filtered-permutation quandle}, and $M^{(1)}, \ldots, M^{(r)}$ are called the {\it parameters} of $X$.
\end{definition}

\begin{remark}
{\rm For $r=2$,  by assuming $m^{(1)}_1\leq m^{(2)}_1$, we recover the quandle $U_{m_1^{(1)},m_1^{(2)}}$ considered in \cite{MR4129183}. These quandles are parametrised by coprime tuples $(m,n)$ with $m\leq n$, and are presented as 
\begin{equation}
U_{m,n} = \big\{x_0,x_1, \ldots, x_{m-1}, y_0,y_1,\ldots, y_{n-1} \big\}
\end{equation}
with
$$x_i * x_j = x_i, \quad y_k * y_l = y_k,\quad  x_i * y_k = x_{i+1}, \quad  y_k * x_i = y_{k+1},$$
where $x_{m}=x_0$, $y_{n}=y_0$, $0 \leq i,j \leq m-1$ and $0 \leq k,l \leq n-1$.}
\end{remark}
\medskip

Using structural results on medial quandles, that is, quandles satisfying the condition $$(a * b) * (c * d)=(a * c) * (b * d),$$
the following result is derived in \cite{MR3400403}. The proof presented below is from \cite[Theorem 2.3]{MR4116819}.

\begin{theorem}\label{abelian filtered per module}
The following statements are equivalent for a finite quandle $X$.
\begin{enumerate}
\item $X$ is abelian.
\item $X$ is isomorphic to a filtered-permutation quandle.
\end{enumerate}
Moreover, two filtered-permutation quandles with $r$ ordered orbits are isomorphic if and only if they have the same parameters $M^{(1)}, \ldots, M^{(r)}$.
\end{theorem}

\begin{proof} 
The implication (2) $\Rightarrow$ (1) follows from Proposition~\ref{abelian quandle from parameters}. For (1) $\Rightarrow$ (2), first note that if $X$ is an abelian quandle, then
$$(a * b) * c = (a * c) * b = (a * b) * (c * b)$$
for all $a,b,c \in X.$ Since the right translation $- * b$ is bijective, we have
$$a * c = a * (c * b)$$
for all $a,b,c \in X.$ Further, since $X$ is finite abelian, we have $a * a' = a * a = a$ for $a$ and $a'$ from the same orbit. Let $O_1,\ldots,O_r$ be the orbits of $X$ under the action of $\Inn(X)$. Further, we have permutations $f_{i,j} \in \Sigma_{O_i}$, $1 \leq i, j \leq r$ such that
\begin{equation}
a * b =f_{i,j}(a)
\end{equation}
for all $a \in O_i, b \in O_j$. These permutations satisfy the following conditions:
\begin{enumerate}[(i)]
\item Commutativity: $f_{i,j}f_{i,k}=f_{i,k}f_{i,j}$.
\item Transitivity: the $f_{i,j}$, $1 \leq j \leq r$, generate a transitive subgroup $G_i$ of $\Sigma_{O_i}$.
\item Freeness: $a \cdot g = a$ for some $a \in O_i$ and $g \in G_i$ implies $a' \cdot g = a'$ for all $a' \in O_i$.
\end{enumerate}
Indeed, (i) follows from abelianity condition, and (ii) from the definition of orbits and the finiteness of~$X$. For (iii), using transitivity, write $a'=a\cdot h$ for some $h \in G_i$ to get
$$a' \cdot g = (a\cdot h) \cdot g = a\cdot (hg)= a\cdot (gh) = (a\cdot g) \cdot h = a\cdot h = a'.$$

Now, fix an index $i$. All the indices below are considered modulo~$r$. By the freeness, the permutation $f_{i,i+1}$ consists of cycles of the same length; denote this length by $m^{(i)}_1$. Further, take an $a \in O_i$; the permutation $f_{i,i+2}$ will send $a$ to a possibly different $f_{i,i+1}$-cycle, but after $m^{(i)}_2$ iterations will bring it back to the original $f_{i,i+1}$-cycle for the first time. This yields a condition $f_{i,i+2}^{m^{(i)}_2}(a)=f_{i,i+1}^{-m^{(i)}_{2,1}}(a)$ for some $0 \leq m^{(i)}_{2,1} < m^{(i)}_1$. Once again, the freeness yields the relation $f_{i,i+1}^{m^{(i)}_{2,1}} f_{i,i+2}^{m^{(i)}_2}=1$ in $\Sigma_{O_i}$. Similarly, by looking when  $f_{i,i+3}$ brings $a$ back to its original $\left\langle f_{i,i+1},f_{i,i+2}\right\rangle$-orbit (where we are considering the subgroup of $\Sigma_{O_i}$ generated by $f_{i,i+1}$ and $f_{i,i+2}$), one finds a relation $f_{i,i+1}^{m^{(i)}_{3,1}} f_{i,i+2}^{m^{(i)}_{3,2}} f_{i,i+3}^{m^{(i)}_3}=1$ in $\Sigma_{O_i}$, with $0 \leq m^{(i)}_{3,1} < m^{(i)}_1$ and $0 \leq m^{(i)}_{3,2} < m^{(i)}_2$.
\para 
Iterating this argument, we obtain a parameter collection $M^{(i)}$ of the form~\eqref{parameter matrix}, and a transitive action of the group $G(M^{(i)})$ on $O_i$.	To be precise, the generator $x^{(i)}_k$ of $G(M^{(i)})$ act by $f_{i,i+k}$. If the action is not free, then we have a relation $f_{i,i+1}^{m'_{k,1}}f_{i,i+2}^{m'_{k,2}}\cdots  f_{i,i+k-1}^{m'_{k,k-1}}f_{i,i+k}^{m'_k}=1$ in $\Sigma_{O_i}$, with $1 \leq k \leq r-1$, and $0<m'_k<m^{(i)}_k$. But this contradicts the minimality in the choice of $m^{(i)}_k$. Hence, the action must be free. 
\para
Now, choosing an $a_i \in \Sigma_{O_i}$ for each $i$, we get the bijections
$$ G(M^{(i)}) \longleftrightarrow O_i$$
given by
$$(x_1^{(i)})^{n_1}(x_2^{(i)})^{n_2}\ldots (x^{(i)}_{r-1})^{n_{r-1}} \longleftrightarrow f_{i,i+1}^{n_1}f_{i,i+2}^{n_2}\ldots f_{i,i+r-1}^{n_{r-1}}(a_i).$$
Since the action of~$f_{i,i+k}$ on~$O_i$ corresponds to multiplying by $x^{(i)}_k$ in $G(M^{(i)})$, we obtain a quandle isomorphism $Q(M^{(1)}, \ldots, M^{(r)}) \cong X$.

Finally, by the freeness of the $G_i$-action on~$O_i$, the parameter collection $M^{(i)}$ is independent of the choice of the orbit representative~$a_i$, and hence is uniquely determined by the isomorphism class of $X$, where we require isomorphisms to preserve a chosen order of orbits. $\blacksquare$
\end{proof}

\begin{remark}
{\rm To avoid ambiguity, we consider finite quandles with ordered orbits. The parameters $M^{(1)}, \ldots, M^{(r)}$  uniquely determine finite abelian quandles up to a reordering of orbits, that is, up to the natural action of the symmetric group $\Sigma_r$. }
\end{remark}
\para

\begin{definition}
Let $X \cong Q(M^{(1)}, \ldots, M^{(r)})$ be a finite abelian quandle with parameters $M^{(1)}, \ldots, M^{(r)}$. Then the group
\begin{equation}
G'(X)= {\prod_{i=1}^r G(M^{(i)})}/ {\left\langle x^{(i)}_{j-i}x^{(j)}_{i-j} \,,\, 1\leq i < j \leq r \right\rangle}
\end{equation}
is called the parameter group of $X$.
\end{definition}

For example, for $r=2$, we have
\begin{eqnarray*}
G' (X) &=& {G(M^{(1)})\times G(M^{(2)})}/{\left\langle x^{(1)}_{1}x^{(2)}_{1} \right\rangle} \\
&\cong&  \left\langle x^{(1)}_{1},x^{(2)}_{1} \,\mid\, (x^{(1)}_{1})^{m^{(1)}_1}, \, (x^{(2)}_{1})^{m^{(2)}_1}, \, x^{(1)}_{1}x^{(2)}_{1} \right\rangle\\
&\cong&  \mathbb{Z}_{\gcd(m^{(1)}_1,m^{(2)}_1)}.
\end{eqnarray*}

Next, we have the following result \cite[Theorem 3.2]{MR4116819}.

\begin{theorem}  \label{ajoint central parameter}
Let $X$ be a finite abelian quandle with $r$ orbits. Then the following assertions hold:
\begin{enumerate}
\item $\Adj(X)$ is a central extension of $\mathbb{Z}^r$ by its parameter group $G'(X)$.
\item $G'(X)$ is a finite abelian group, and is isomorphic to the commutator subgroup of $\Adj(X)$.
\end{enumerate}
\end{theorem}

\begin{proof}
Recall that
$$\Adj(X) = \big\langle \textswab{a}_a, \, a \in X \,\mid\, \textswab{a}_b \textswab{a}_a \textswab{a}_b^{-1} =  \textswab{a}_{a* b} ~~ \textrm{for all}~~  a,b \in X \big\rangle.$$
By Theorem~\ref{abelian filtered per module}, we can assume that $X \cong Q(M^{(1)}, \ldots, M^{(r)})$. The defining relations of $\Adj(X)$ are
\begin{equation}\label{adjoint defining parameter}
\textswab{a}_b \textswab{a}_a = \textswab{a}_{ax_{j-i}} \textswab{a}_b 
\end{equation} 
for $a \in G(M^{(i)})$ and $b \in G(M^{(j)})$. As usual, the index $(j-i)$ is taken modulo $r$. We omit the indexing ${}^{(i)}$ when the context is clear.
\para

Let $G_i$ denote the subgroup of $\Adj(X)$ generated by $\{\textswab{a}_a \, \mid \,a \in G(M^{(i)}) \}$. Since $x^{(i)}_0=1$, it is abelian. Further, we can rewrite \eqref{adjoint defining parameter} as
\begin{equation}\label{adjoint defining parameter3}
\textswab{a}_b \textswab{a}_a \textswab{a}_b^{-1} \textswab{a}_a^{-1} = \textswab{a}_{ax_{j-i}} \textswab{a}_a^{-1} \in G_i.
\end{equation}
In particular, this expression is independent of~$b$. Interchanging the roles of $a$ and $b$ gives
$$\textswab{a}_a \textswab{a}_b \textswab{a}_a^{-1} \textswab{a}_b^{-1}= \textswab{a}_{bx_{i-j}} \textswab{a}_b^{-1}  \in G_j,$$
which is independent of~$a$, and is the inverse of the preceding expression. Denoting both sides of~\eqref{adjoint defining parameter3} by $g_{i,j}$, we  obtain elements $g_{i,j} \in G_i \cap G_j$, which allow us to break the relations \eqref{adjoint defining parameter3} into two parts, namely,
\begin{eqnarray}
\textswab{a}_b \textswab{a}_a  &=& g_{i,j} \textswab{a}_a \textswab{a}_b  \; \text{ for } a \in G(M^{(i)}) \text{ and } b\in G(M^{(j)}),\label{adjoint defining parameter2}\\
\textswab{a}_{ax_{j-i}} &=& g_{i,j} \textswab{a}_a \; \text{ for } a \in G(M^{(i)}).\label{adjoint defining parameter4}
\end{eqnarray} 
Moreover, $g_{i,j}$ satisfy the conditions
\begin{eqnarray}
g_{i,j}g_{j,i} &=& 1 \; \text{ for } 1 \leq i < j \leq r,\label{equation of gij}\\
g_{i,i}&=& 1 \; \text{ for } 1 \leq i \leq r.\label{equation of gij2}
\end{eqnarray}

We claim that
\begin{equation}\label{equation of gij3}
g_{i,j} \text{ lies in the center of } \Adj(X) ~\text{for all}~1 \leq i, j \leq r.
\end{equation}
Indeed, $g_{i,j} $ can be written as $\textswab{a}_{a'} \textswab{a}_a^{-1}$, with $a$ and $a'$ from the same orbit $G(M^{(i)})$. Taking $b \in G(M^{(j)})$, we see that 

\begin{equation}\label{E:Compute}
\textswab{a}_b \textswab{a}_{a'} {\textswab{a}_a}^{-1}=\textswab{a}_{a'} \textswab{a}_b g_{i,j} {\textswab{a}_a}^{-1}= \textswab{a}_{a'} \textswab{a}_b {\textswab{a}_a}^{-1}g_{i,j}  =\textswab{a}_{a'} {\textswab{a}_a}^{-1} \textswab{a}_b g_{j,i} g_{i,j}=\textswab{a}_{a'} {\textswab{a}_a}^{-1} \textswab{a}_b,
\end{equation}
since $g_{i,j}$ commutes with $G_i$. For each $1 \le i \le r$, set $h_i:=\textswab{a}_{1^{(i)}} \in G_i$. It follows from \eqref{adjoint defining parameter4} that the set 
$$ \big\{g_{i,j} \, \mid \,  1 \le i,j \le r \big\} \cup \big\{h_i \,  \mid  \,  1 \le i \le r \big\}$$
generates $\Adj(X)$. Indeed, we can put
$$g_{(x_1^{(i)})^{n_1}(x_2^{(i)})^{n_2}\ldots (x^{(i)}_{r-1})^{n_{r-1}}} = g_{i,i+1}^{n_1}g_{i,i+2}^{n_2}\ldots g_{i,i+r-1}^{n_{r-1}}h_i.$$
This is well-defined if and only if
\begin{equation} \label{equation of gij4}
 g_{i,i+1}^{m^{(i)}_{j,1}}g_{i,i+2}^{m^{(i)}_{j,2}}\ldots g_{i,i+j}^{m^{(i)}_{j,j}}=1
\end{equation}
for $1 \leq i \leq r$ and  $1 \leq  j < r$.  If we assume these conditions, then the relations~\eqref{adjoint defining parameter4} become redundant. Finally, since each $g_{i,j}$ is central, it is  sufficient to check relations~\eqref{adjoint defining parameter2} for the generators $h_i$ only, that is,
\begin{equation} \label{equation of gij5}
h_j h_i  = g_{i,j}  h_i h_j.  
\end{equation}
This gives a new presentation of the adjoint group as
\begin{equation}\label{new adjoint grp present}
\Adj(X) \cong \Big\langle g_{i,j}, \, 1 \leq i,j \leq r \, ; \, h_i, \, 1 \leq i \leq r \;\mid\; \eqref{equation of gij}, \eqref{equation of gij2}, \eqref{equation of gij3}, \eqref{equation of gij4}, \eqref{equation of gij5} \Big\rangle.
\end{equation}

It follows from the presentation~\eqref{new adjoint grp present} that the commutator subgroup $[\Adj(X), \Adj(X)]$ of ~$\Adj(X)$ is generated by $\{ g_{i,j} \mid 1 \le i, j \le r \}$. It is well-known that the abelianisation of~$\Adj(X)$ is $\mathbb{Z}^r$, and hence we obtain the short exact sequence
$$0 \to [\Adj(X), \Adj(X)] \to \Adj(X) \to \mathbb{Z}^r \to 0$$
of groups. Since each $g_{i,j}$ is central in~$\Adj(X)$, it follows that $\Adj(X)$ is a central extension of $\mathbb{Z}^r$ by $[\Adj(X), \Adj(X)]$. 
\para 
It remains to prove that $G'(X) \cong [\Adj(X), \Adj(X)]$. Relations \eqref{equation of gij}, \eqref{equation of gij2} and \eqref{equation of gij4} give a surjective  homomorphism  $$\psi : G'(X)\to [\Adj(X), \Adj(X)]$$ defined by $\psi(x^{(i)}_j) =  g_{i,i+j}$. For injectivity, we construct a set-theoretic map $\pi : \Adj(X) \to G'(X)$ whose restriction to $[\Adj(X), \Adj(X)]$ is the inverse of $\psi$. Take an element $g \in \Adj(X)$ written in terms of the generators $g_{i,j}$ and $h_i$. Move all the occurrences of $h_1^{\pm 1}$ to the left using the centrality of the $g_{i,j}$ and the twisted commutativity~\eqref{equation of gij5} of the $h_i$. Similarly, move all the occurrences of $h_2^{\pm 1}$ right after the $h_1^{\pm 1}$, and so on. Using the trivial relations $h_i h_i^{-1}= h_i^{-1} h_i = 1$, we obtain a word of the form $h_1^{k_1}\ldots h_r^{k_r}g''$, where $k_1,\ldots,k_r \in \mathbb{Z}$ and $g''$ is a product of the $g_{i,j}^{\pm 1}$. Next, in $g''$ replace each generator $g_{i,j}$ by $x^{(i)}_{j-i}$ and denote by $g'$ the word obtained. Considering it as an element of~$G'(X)$, we set $\pi(g)=g'$. The map is well-defined, since relations \eqref{equation of gij}, \eqref{equation of gij2}, \eqref{equation of gij4} and $g_{i,j}^{\pm 1}g_{i,j}^{\mp 1}=1$ have counterparts in~$G'(X)$. Further, the relation~\eqref{equation of gij3} does not change the result by construction. Similarly, the relations \eqref{equation of gij5} and $h_i^{\pm 1}h_i^{\mp 1}=1$ do not change the result as shown by a computation similar to~\eqref{E:Compute} together with~\eqref{equation of gij}. Consider the restriction $$\varphi: [\Adj(X), \Adj(X)] \to G'(X)$$
of $\pi$ to $[\Adj(X), \Adj(X)]$. It simply replaces each $g_{i,j}$ by $x^{(i)}_{j-i}$ in any representative of an element of~$[\Adj(X), \Adj(X)]$, and hence is the desired inverse of~$\psi$. Finally, since the groups $G(M^{(i)})$ are finite abelian groups, it follows that $G'(X)$ is a finite abelian group.  $\blacksquare$
\end{proof}

\begin{definition}
Let $X$ be a finite abelian quandle with parameters $M^{(1)}, \ldots, M^{(r)}$. The parameter matrix $\mathcal{M}(X)$ of $X$ is defined as follows: Its columns are indexed by tuples $(i,j)$ with $1 \leq i < j \leq r$. Its rows are indexed by tuples $(i,j)$ with $1 \leq i \leq r$ and $1 \leq j < r$. All tuples are ordered lexicographically here. The row $(i,j)$ corresponds to the $j$-th row of $M^{(i)}$; for all $1 \le k \le \min \{j, r-i\}$, it contains $m^{(i)}_{j,k}$ in the column $(i,i+k)$, and for all $r-i < k \leq j$, it contains $-m^{(i)}_{j,k}$ in the column $(i+k-r,i)$. 
\end{definition}

Let us see two special cases. For $r=2$, $M^{(1)} = \begin{psmallmatrix} m^{(1)}_1 \end{psmallmatrix}$ and $M^{(2)} = \begin{psmallmatrix} m^{(2)}_1 \end{psmallmatrix}$, we have
$$\mathcal{M}(X)=\begin{psmallmatrix}
m^{(1)}_1 \\ -m^{(2)}_1
\end{psmallmatrix}.$$

For $r=3$, there are three columns $(1,2)$, $(1,3)$, $(2,3)$, and 
$$\mathcal{M}(X)=\begin{psmallmatrix}
m^{(1)}_1 & 0 & 0\\
m^{(1)}_{2,1} & m^{(1)}_{2} & 0\\ 
0 & 0 & m^{(2)}_1\\
-m^{(2)}_{2} & 0 & m^{(2)}_{2,1}\\
0 & -m^{(3)}_1 & 0\\
0 & -m^{(3)}_{2,1} & -m^{(3)}_{2}\\
\end{psmallmatrix}.$$

The following result classifies finite quandles with abelian adjoint groups \cite[Theorem 4.2]{MR4116819}.

\begin{theorem}\label{classify quandle abelian adjoint}  
Let $X$ be a finite quandle. Then the following statements are equivalent:
\begin{enumerate}
\item $\Adj(X)$ is abelian.
\item $X$ is abelian and its parameter group $G'(X)$ is trivial.
\item $X$ is abelian and the maximal minors of its parameter matrix $\mathcal{M}(X)$ are globally coprime.
\end{enumerate}
\end{theorem}

\begin{proof}
There is a natural  left-action of $\Adj(X)$ on~$X$ given by $\textswab{a}_b \cdot a =a * b$ for $a, b \in X$. If $\Adj(X)$ is abelian, then
$$(a * b) * c = \textswab{a}_c \cdot (\textswab{a}_b \cdot a) = (\textswab{a}_c \textswab{a}_b) \cdot a = (\textswab{a}_b \textswab{a}_c) \cdot a=\textswab{a}_b \cdot  (\textswab{a}_c  \cdot a)= (a * c) * b.$$
Hence, if $\Adj(X)$ is abelian, then $X$ is abelian. Thus, we only need to understand which finite abelian quandles have abelian adjoint groups. By Theorem~\ref{ajoint central parameter}, this happens if and only the parameter group $G'(X)$ is trivial. Furthermore, if $G'(X)$ is trivial, then $\Adj(X) \cong \mathbb{Z}^r$. And, if $\Adj(X)$ is abelian, then by Theorem \ref{thm Lebed Vendramin}, it is free abelian. Hence, the only possibility for its finite subgroup~$G'(X)$ is that it must be trivial. Thus, we have proved the equivalence of (1) and (2).
\para
Assume that $X$ is abelian and has parameters $M^{(1)}, \ldots, M^{(r)}$. The group~$G'(X)$ admits as generators the elements $x^{(i)}_{j}$ for $i+j \leq r$, since for $i+j > r$ one has $x^{(i)}_{j}=(x^{(j+i-r)}_{r-j})^{-1}$. With these $n:=\frac{r(r-1)}{2}$ generators, $G'(X)$ is isomorphic to the quotient of $\mathbb{Z}^n$ by the row space of the matrix $\mathcal{M}(X)$. Indeed, the rows of~$\mathcal{M}$ encode the defining relations of the components $G(M^{(i)})$ of~$\Adj(X)$, taking into account the identification $x^{(i)}_{j}=(x^{(j+i-r)}_{r-j})^{-1}$. We claim that triviality of $G'(X)$ is equivalent to the maximal minors of $\mathcal{M}(X)$ being globally coprime. To see this, note that for a given finitely generated abelian group, both its isomorphism class and the greatest common divisor of the maximal minors of its presentation matrix as above are invariant under elementary row and column operations. Further, for a matrix in Smith normal form, the triviality of the group and the minors condition are both equivalent to the matrix being of maximal rank with all diagonal entries equal to $1$. This establishes the equivalence of (2) and (3). 
$\blacksquare$
\end{proof}

For $r=2$, the coprimality condition in Theorem \ref{classify quandle abelian adjoint} becomes $\gcd \big(m^{(1)}_1,~m^{(2)}_1 \big)=1$, and recovers the main result of Bardakov and Nasybullov \cite{MR4129183}. The next result provides a large supply of quandles satisfying the conditions of the preceding theorem  \cite[Proposition 4.3]{MR4116819}.

\begin{proposition} 
Let $X$ be a finite abelian quandle with parameters $M^{(1)}, \ldots, M^{(r)}$ such that $m^{(i)}_{j,k}$ vanish whenever $k<j$. Then the following conditions are equivalent:
\begin{enumerate}
\item  $\Adj(X) \cong \mathbb{Z}^r$.
\item $\gcd \big(m^{(i)}_{j},~m^{(j+i-r)}_{r-j} \big)=1$ whenever $i+j>r$.
\end{enumerate}
\end{proposition}

\begin{proof}
We check the triviality of the abelian group $G'(X)$ and use the equivalence (1) $\Longleftrightarrow$ (2) of Theorem~\ref{classify quandle abelian adjoint}. In our situation, $G'(X)$ has the presentation
$$G'(X) \cong \Big\langle \, x^{(i)}_j \,\mid\, \left( x^{(i)}_j\right)^{m^{(i)}_{j}} =1,\, \, x^{(l)}_{k-l}x^{(k)}_{l-k} =1 \, \Big\rangle,$$
where $1\leq i \leq r$, $1\leq j < r$ and $1\leq l < k \leq r$. The last condition means that the generators $x^{(l)}_{k-l}$ and $x^{(k)}_{l-k}$ are mutually inverses. The above presentation can be rewritten as
$$G'(X) \cong \Big\langle \, x^{(i)}_j, \, i+j>r \,\mid \, \left( x^{(i)}_j\right)^{\gcd (m^{(i)}_{j},m^{(j+i-r)}_{r-j})} =1~~\textrm{for}~ i+j>r \, \Big\rangle,$$
where $1\leq i \leq r$ and $1\leq j < r$. But this is the direct product of cyclic groups of orders $\gcd \big(m^{(i)}_{j},~m^{(j+i-r)}_{r-j} \big)$, where $1\leq i \leq r$, $1\leq j < r$ and $ i+j>r$, and the proof is complete. $\blacksquare$
\end{proof}

\begin{remark}
{\rm
The quandles from the preceding proposition have $r(r-1)$ non-zero parameters $m^{(i)}_{j}$, and condition (2) divides them into $\frac{r(r-1)}{2}$ coprime pairs. Thus, one obtains, for each~$r,$ an infinite family of quandles with adjoint group $\mathbb{Z}^r$.
}
\end{remark}
\bigskip
\bigskip


\section{Augmented  quandles}

In this section, we consider quandles that arise from an action of a group on a set. The construction is attributed to Joyce \cite[Section 9]{MR0638121}.

\begin{definition} 
 An  \index{augmented quandle}{\it augmented  quandle} is a pair  $(X, G)$, where $X$ is a set equipped with a right-action of a group $G$ such that there is a map $\varepsilon : X \to G$,  called the \index{augmentation  map}{\it augmentation  map},  which satisfies the following conditions:
\begin{enumerate}
\item $x \cdot \varepsilon(x) = x  ~\mbox{for all}~x \in X$.
\item $\varepsilon(x \cdot g) = g \varepsilon(x) g^{-1}$ for all $x \in X$ and $g \in G$.
\end{enumerate}
\end{definition}

Given an augmented quandle $(X, G)$, we can associate a quandle structure on the set $X$ by setting
$$
x * y = x \cdot \varepsilon(y)
$$
for $x, y \in X$. Then the action of $G$ on $X$ is by quandle automorphisms, and the augmentation map $\varepsilon  :  X \to G$ is a quandle homomorphism,  when $G$ is viewed as the conjugation quandle.
\para

If $X$ is also a group and $\varepsilon$ is a group homomorphism, the pair $(X, G)$ is called a \index{crossed module}{\it crossed module}. We can associate an augmented quandle to each quandle. More precisely, for a quandle $X$, we have the augmented quandle  $\big(X, \Adj(X)\big)$, where  the augmentation map is given by $x \mapsto \textswab{a}_x$. Similarly, $\big(X, \Aut(X)\big)$ and $\big(X, \Inn(X)\big)$ are also augmented quandles, where the augmentation maps are $x \mapsto S_x$ in each case.

\begin{definition} A  morphism of augmented quandles from $(X, G)$ to  $(Y, H)$ consists
of a group homomorphism  $g : G \to H$ and a map  $f :  X \to Y$ such that the diagram
$$
\begin{array}{ccccc}
  X \times G & \longrightarrow & X & \overset{\varepsilon_X}{\longrightarrow} & G \\
{\scriptstyle f\times g} \downarrow &  & {\scriptstyle f} \downarrow &  & {\scriptstyle g} \downarrow \\
  Y \times H & \longrightarrow & Y & \overset{\varepsilon_Y}{\longrightarrow}  & H
\end{array}
$$
commutes.
\end{definition}

 It follows that $f$ is a homomorphism of associated quandles.
 
\begin{example}{\rm   If $X$ is a quandle, then $\big(X, \Adj(X)\big)$, $\big(X, \Aut(X)\big)$ and $\big(X, \Inn(X)\big)$ are augmented quandles with associated quandle as $X$. The set of all augmentations of $X$ forms a  category in which $\big(X, \Aut(X)\big)$ is the final object.  In other words, for each augmentation  $(X, G)$, there is a unique group
homomorphism  $f : G \to \Aut(Q)$ such that the diagram
$$
\begin{array}{ccccc}
  X \times G & \longrightarrow & X & \overset{\varepsilon}{\longrightarrow} & G \\
{\scriptstyle \id\times f} \downarrow &  & {\scriptstyle \id} \downarrow &  & {\scriptstyle f} \downarrow \\
 X \times \Aut(X) & \longrightarrow & X & \overset{S}{\longrightarrow}  & \Aut(X)
\end{array}
$$
commutes. The map  $f$  is readily defined from the action $X \times G \to X$ as $f(g)(x)=x \cdot g$.}
\end{example}
\bigskip

Let $(X, G)$ be an augmented quandle with augmentation map $\varepsilon:X \to G$. Given a group homomorphism $f : G \to H$, we define another augmented quandle $(X \otimes_G H, H)$ as follows. Note that the cartesian product $X \times H$ can be turned into a right $H$-set with the action
$$
(x, h) \cdot k := (x, h k)
$$
for $x \in X$ and $h, k \in H$. Define a $H$-equivalence  on $X \times H$ by
$$
(x, h) \sim (y, k) \Longleftrightarrow h k^{-1} = f (g) ~\mbox{and}~ y = x \cdot g
$$
for some $g \in G$. Let $X \otimes_G H$ denote the set of congruence classes, and $x \otimes h$ denote the congruence class of the element $(x,h) \in X \times H$.
 Define a map
$$
\varepsilon' : X \otimes_G H \to  H 
 $$
 given by
$$\varepsilon'(x\otimes h) = h^{-1} \big(f (\varepsilon(x)) \big) h.$$
Then $\varepsilon'$ is well-defined and the pair $(X \otimes_G H, H)$ is an augmented quandle. Using this construction, we can define quotients of quandles as follows.

\begin{definition} 
Let $X$ be a quandle and $N$ a normal subgroup of $G = \Adj(X)$. Consider the augmented quandle $(X, G)$. The quotient of $X$ by $N$ is the set $X \otimes_G G / N$ equipped with the quandle operation
$$
[p] * [q] := \big[p \cdot \varepsilon'(q) \big]
$$
for $p, q \in X \otimes_G G / N$. We denote the quotient quandle by $X / N$.
\end{definition}

There is a natural quandle epimorphism $X \to  X / N$. Conversely,  given a quandle epimorphism, we can construct a corresponding quotient of the quandle.

\begin{proposition} 
Let $f : X \to Y$ be a quandle epimorphism, which induces a bijection between the orbit sets $\mathcal{O}(X)$ and $\mathcal{O}(Y)$. Then $Y$ is
isomorphic to the quandle from the augmented data $\big(X \otimes_{\Adj(X)} \Adj(Y), \Adj(Y) \big)$.
\end{proposition}

The preceding proposition is analogous to the first isomorphism theorem in group theory. In fact, the group $\Adj(Y)$ is isomorphic to $\Adj(X) / \ker(\Adj(f))$, where $\Adj(f):\Adj(X) \to \Adj(Y)$ is the induced homomorphism.
\para 

We know from Proposition \ref{equivalence-two-models} that the free quandle $FQ(S)$ on a set $S$ is isomorphic to the conjugacy classes of $S$ in the free group $F(S)$. By definition,  $FQ(S)$ has exactly $|S|$ many connected components and $\Adj \big(FQ(S)\big)=F(S)$. Furthermore, the natural inclusion $\varepsilon_S : FQ(S) \hookrightarrow F(S)$ gives an augmented quandle $\big(FQ(S), F(S)\big)$.  Note that, for any normal subgroup $N$ of $F(S)$, we can consider the quotient quandle $FQ(S) / N$. 
\bigskip

Next, we use the idea of augmented quandles to describe specific quotients  of quandles, in particular,  the largest medial  and the largest involutory quotients. Let $(X, G)$ be an augmented quandle with augmentation map $\varepsilon:X \to G$ and $N$  a normal subgroup of $G$. Elements of the quotient quandle $X/N$ correspond  to equivalence classes
$$
\overline{q}= \big\{ q\cdot n \in X   \, \mid \,   n \in N \big\}
$$
for $q\in X$. We denote elements of $G/N$ as $\overline{x} = x N$ for $x\in G$. Then the right-action
$$
X/N \times  G/N \to X/N
$$
is given by
$$
\overline{q} \cdot \overline{x} = \overline{q \cdot x},
$$
and the augmentation  $\varepsilon : X/N \to G/N$ is given by $$\varepsilon(\overline{q}) = \overline{\varepsilon(q)}.$$
\para

Let us first consider the largest medial quotient of a quandle. Let  $(X, G)$ be  an augmented quandle.  In order that $X$ be medial, we require that
$$
  (p * q) * (r * s) = (p * r) * (q * s),
$$
which is equivalent to
$$
\varepsilon(q) \varepsilon \big(r\varepsilon(s) \big)  = \varepsilon(r) \varepsilon \big(q\varepsilon(s) \big).
$$
In other words, every element of the form
\begin{equation} \label{m1}
\varepsilon(q) \varepsilon(s)^{-1} \varepsilon(r) \varepsilon(q)^{-1} \varepsilon(s) \varepsilon(r)^{-1}
\end{equation}
equals to 1. Let $N$ be the normal subgroup of $G$ generated by elements of the
form \eqref{m1}. Then the augmented quandle  $(X/N, G/N)$, which is the quotient of the augmented quandle $(X, G)$, is  medial. It is
evident that $(X/N, G/N)$  has the universal property that each morphism $(X, G) \to (Y, H)$ of augmented quandles
factors uniquely through $(X/N, G/N)$ whenever $Y$ is a medial quandle. Moreover,
if $\varepsilon(X)$ generates $G$, then $X/N$ is the largest medial quotient  of $X$.  The following result of Joyce is well-known \cite[Theorem 10.1]{MR0638121}.

\begin{theorem}  Let $A$ be a set and $G$ the group generated by $A$ modulo the
relations $a b^{-1} c = c b^{-1} a$ for conjugates $a, b, c$ of the generators of $G$. Then the free
medial quandle on $A$ appears as the conjugates of the generators of $G$ as a
subquandle of $\Conj (G)$.
\end{theorem}

This method works for any family of quandles determined by equations of the form
$$
p  *^{\varepsilon_1} \varphi_1 *^{\varepsilon_2} \ldots *^{\varepsilon_n} \varphi_n = p  *^{\mu_1} \psi_1 *^{\mu_2} \ldots *^{\mu_m} \psi_m,
$$
where $\varphi_i$ and $\psi_j$ are expressions not involving $p$. For example, the equation
$$
p *^n q = p
$$
 describes $n$-quandles. The following result is known for $n$-quandles \cite[Theorem 10.2]{MR0638121}.

\begin{theorem}
Let $(X, G)$ be an augmented quandle such that $\varepsilon(X)$ generates $G$, and $n$ a positive integer. Let $N$ be the normal subgroup of $G$ generated by elements
of the form $\varepsilon(x)^n$, where $x \in X$. Then $X/N$ is the largest quotient of $X$ which is an $n$-quandle.
\end{theorem}

\begin{corollary}  Let $X$ be a set and $G$ be the group with presentation
$$
G = \big\langle  X   \, \mid \,   x^n = 1 ~\textit{for all}~ x \in X \big\rangle.
$$
Then the free $n$-quandle on $X$ consists of conjugates of generators of $G$.
\end{corollary}

\begin{corollary}  
The free  involutory quandle  on two elements is isomorphic  to
$\Core(\mathbb{Z})$  with generators 0 and 1.
\end{corollary}

\begin{proof}  Let $G =  \langle a, b \, \mid \,  a^2 = b^2 =  1 \rangle \cong \mathbb{Z}_2 * \mathbb{Z}_2$
be the infinite dihedral group. The quandle of conjugates of  $a$  and
$b$  in $G$ is
$$
Q =  \big\{ a x^k    \, \mid \,   k \in \mathbb{ Z} \big\}.
$$
The equation $ a x^n * a x^m = a x^{2n-m}$ implies that the map
$$
f :  Q \to \Core(\mathbb{Z}), 
$$
given by $f(a x^n) = n$, is an isomorphism of quandles.  $\blacksquare$
\end{proof}

We conclude with the following result \cite[Theorem 10.5]{MR0638121}.

\begin{theorem} The free medial involutory quandle on $(n + 1)$  generators can be identified with the subquandle
$$
P = \big\{ (k_1, \ldots,  k_n) \in \mathbb{Z}^n  \, \mid \,  \mbox{at most one}~ k_i~ \mbox{is odd} \big\}
$$
of $\Core(\mathbb{Z}^n)$. Further, generators of $P$ are
$$
e_0 = (0, 0, \ldots, 0), \, e_1 = (1, 0, \ldots, 0), \, \ldots, \, e_n = (0, 0, \ldots, 0, 1).
$$
\end{theorem}

\begin{proof} Let $G$ be the group presented as
$$
G=\big\langle a_0, \ldots, a_n   \, \mid \,  a_i^2 =1~\textrm{and}~ a_i a_j a_k  = a_k a_j a_i \, \textrm{for all} \,   i,  j,  k \big\rangle,
$$
and $Q$ be the set of conjugates of generators of $G$. Then $Q$ is the free medial involutory quandle on  $S = \{ a_0, \ldots, a_n \}$.  Since  $P$  is a medial involutory quandle, there is a unique map  $h : Q \to P$  such that  $h(a_i) = e_i$ for each  $i$.  We  claim that $h$  is an isomorphism. Setting  $t_0 = a_0 a_0$ and $t_i=a_0 a_i$  for $1 \le i \le n$, we see that $t_j t_k  =  t_k t_j$.  The  conjugates  of $a_i$  are  of the form
$$
a_i * a_{j_1} * \cdots  * a_{j_r} = a_i * a_{j_1} a_{j_2} \cdots   a_{j_r}.
$$
Since $a_i * a_i = a_i$, we can take $r$ to be even. Further, since
$$
a_i * a_{j_1} * \cdots  * a_{j_r} = a_i * t_{j_1}^{-1} t_{j_2} \cdots t_{j_{r-1}}^{-1} t_{j_r},
$$
we see that $Q$ consists of the elements of the form $ a_i * t_{1}^{k_1} \cdots t_{n}^{k_n},$
where $k_i \in \mathbb{Z}$ and $1 \le i \le n$. We now have 
$$
h(a_i * t_{1}^{k_1} \cdots t_{n}^{k_n}) = e_i + (2 k_1, \ldots, 2 k_n),
$$
which proves the bijectivity of $h$. $\blacksquare$
\end{proof}

\begin{remark}
{\rm
Alternatively,  we can describe the free medial involutory  quandle on $(n + 1)$ generators  as the subquandle
$$
\big\{ (k_0, \ldots , k_n) \in \mathbb{Z}^{n+1}  \, \mid \,   \mbox{exactly  one}~ k_i~ \mbox{is odd} \big\}
$$
 of $\Core(\mathbb{Z}^{n+1})$.}
\end{remark}
\bigskip
\bigskip


\chapter{Automorphisms of quandles}\label{chapter-auto-quandles}

\begin{quote}
A comprehensive understanding of automorphisms of algebraic structures is crucial for their systematic classification. In this chapter, we explore automorphisms of Takasaki quandles and conjugation quandles of groups. Additionally, we study automorphisms of quandle extensions, as well as the subgroup of quasi-inner automorphisms of these structures. As we progress, we will also touch upon topics such as connectivity and commutativity in Alexander quandles, and the multi-transitive actions of these automorphisms groups.
\end{quote}
\bigskip

\section{Automorphisms of Takasaki quandles}\label{section4}

Automorphisms play a crucial role in the structure theory of quandles. We begin by making some basic observations about automorphisms of quandles. 
\para
It is readily apparent that the set of all inner automorphisms of a group is a group. But, this is not the case for quandles. The following example shows, in general, that the group $\Inn(X)$ is bigger than the set $\{ S_x \, \mid \, x \in X \}$ of all symmetries.

\begin{example}{\rm 
Consider the dihedral quandle $\R_3 = \{ 0, 1, 2 \}$, where $i * j = 2j-i \mod 3$. Then the set of symmetries  contains only three elements, but $\Inn(\R_3)=\Aut(\R_3) \cong \Sigma_3$, the permutation group on three symbols.}
\end{example}

The following example shows that, in general, $\Aut(X) \neq \Inn(X)$.

\begin{example}{\rm 
Let $G = \{ \pm 1, \pm i, \pm j, \pm k \}$ be the quaternion group of order 8 and $\Conj(G)$ its conjugation quandle.
Let $X$ be the subquandle of $\Conj(G)$ consisting of conjugacy classes of elements from the set $\{\pm i, \pm j \}$. Then we see that $\Inn(X) = \langle S_i, S_j \rangle \cong \mathbb{Z}_2 \times  \mathbb{Z}_2$ and $\Aut(X) = \langle S_i, S_j, \varphi \rangle$ has order 8, where $\varphi(\pm i) = \pm j$ and $\varphi(\pm j) = \pm i$.}
\end{example}

First, we present the following slightly more general result \cite[Proposition 4.1]{MR3718201}.

\begin{proposition}\label{joyce-embedding}
Let $G$ be a group and $\varphi \in \Aut(G)$. Then there is an embedding $\Z(G) \rtimes \C_{\Aut(G)}(\varphi)  \hookrightarrow \Aut\big(\Alex(G, \varphi)\big)$, where $\C_{\Aut(G)}(\varphi)$ is the centraliser of $\varphi$ in $\Aut(G)$ and $\Z(G)$ is the center of $G$.
\end{proposition}

\begin{proof}
Let $f \in \C_{\Aut(G)}(\varphi)$. Then for $a, b \in G$, we see that
\begin{eqnarray*}
f(a*b) &=& f\big( \varphi(a)\varphi(b)^{-1}b\big)\\
& = & f\big(\varphi(a)\big) f\big(\varphi(b)\big)^{-1}f(b)\\
& = & \varphi\big(f(a)\big) \varphi\big(f(b)\big)^{-1}f(b)\\
& = & f(a)*f(b).
\end{eqnarray*}
 Hence, we have $\C_{\Aut(G)}(\varphi) \leq  \Aut\big(\Alex(G, \varphi)\big)$.
\para 
Now, for $a \in \Z(G)$, let $t_a: G \to G$ be given by $t_a(b)=ba$ for all $b \in G$. Then $t_a \in \Aut\big(\Alex(G, \varphi)\big)$ and the map $a \mapsto t_a$ gives an embedding $\Z(G) \hookrightarrow \Aut\big(\Alex(G, \varphi)\big)$.

Finally, we show that the map $\Phi:\Z(G) \rtimes \C_{\Aut(G)}(\varphi) \to \Aut\big(\Alex(G, \varphi)\big)$ given by $\Phi(a,f) = t_a  f$ is an embedding. The map is clearly injective. Further, for $(a_1,f_1), (a_2, f_2) \in  \Z(G) \rtimes \C_{\Aut(G)}(\varphi)$, we have
\begin{eqnarray*}
\Phi\big((a_1,f_1)(a_2, f_2) \big) & = & \Phi\big((a_1 f_1(a_2), f_1  f_2) \big)\\
& = & t_{a_1 f_1(a_2)}  (f_1  f_2)\\
& = & (t_{a_1}  f_1)  (t_{a_2}  f_2)\\
& = & \Phi\big((a_1,f_1)\big)  \Phi\big((a_2, f_2) \big).
\end{eqnarray*}
This completes the proof of the proposition. $\blacksquare$
\end{proof}

We now take $G$ to be an additive abelian group and $\varphi=-\id_G$. In this case $\Alex(G, \varphi)=\T(G)$ and  $\C_{\Aut(G)}(\varphi)= \Aut(G)$. The following result gives a complete description of the  automorphism and the inner automorphism groups of some Takasaki quandles \cite[Theorem 4.2]{MR3718201}. 

\begin{theorem}\label{main-theorem}
Let $G$ be an additive abelian group without 2-torsion. Then the following  assertions hold:
\begin{enumerate}
\item $\Aut\big(\T(G)\big) \cong G \rtimes \Aut(G)$.
\item $\Inn\big(\T(G)\big) \cong 2G \rtimes \mathbb{Z}_2$.
\end{enumerate}
\end{theorem}

\begin{proof}
We already proved in Proposition \ref{joyce-embedding} that there is an embedding $\Phi:G \rtimes \Aut(G)  \hookrightarrow \Aut\big(\T(G)\big)$.  Let $f \in \Aut\big(\T(G)\big)$. Then $f(a*b)=f(a)* f(b)$ implies that $f(2b-a)=2f(b)-f(a)$ for all $a, b \in \T(G)$. Taking $a=2b$, we have  $$f(2b)=2f(b)-f(0)$$ for all $b \in G$, where $0$ is the identity element  of the group $G$.
\para
Define $h: G \to G$ by $h(a)=f(a)-f(0)$ for $a \in G$. We claim that $h \in \Aut(G)$. Since $f$ is a bijection, it follows that $h$ is also a bijection.  To show that $h$ is a group homomorphism, it suffices to show that $f(a+b)=f(a)+f(b)-f(0)$ for all $a, b \in G$. We can write
$$f(a+b)=f\big(2a-(a-b)\big)=2f(a)-f(a-b)$$
and
$$f(a+b)=f\big(2b-(b-a)\big)=2f(b)-f(b-a).$$
 Adding the two identities yields $2 f(a+b)= 2f(a)+ 2f(b)-f(a-b)-f(b-a)$. But, $f(a-b)=f\big(0-(b-a)\big)= 2f(0)-f(b-a)$. This gives $2 f(a+b)= 2f(a)+ 2f(b)- 2 f(0)$. Since $G$ has no 2-torsion, our claim follows. By definition of $h$, we can write $f=t_{f(0)}  h \in \Phi \big(G \rtimes \Aut(G) \big)$, and the proof of assertion (1) is complete.
\para
Let $r \in \Aut(G)$ be the automorphism given by $r(a) = -a$ for all $a \in G$. Then $\langle r \rangle \cong \mathbb{Z}_2$. Note that $\Phi(2a, r) = t_{2a}  r=S_a$, and hence $\Phi$ embeds $2G \rtimes \mathbb{Z}_2$ into $\Inn\big(\T(G)\big)$. Conversely, any $S_a \in \Inn\big(\T(G)\big)$ can be written as $S_a= t_{2a}  r$, and the proof of assertion (2) is complete.  $\blacksquare$
\end{proof}

As an immediate consequence, we obtain the following result.

\begin{corollary}
$\Aut\big( \T(\mathbb{Z}^n)\big) \cong \mathbb{Z}^n \rtimes \GL(n, \mathbb{Z})$ and  $\Inn\big( \T(\mathbb{Z}^n)\big) \cong (2\mathbb{Z})^n \rtimes \mathbb{Z}_2$ for all $n \geq 1$.
\end{corollary}

\begin{remark}{\rm 
In \cite[Proposition 2.1]{MR2831947}, Hou proved a general result  analogous to Theorem \ref{main-theorem} for Alexander quandle of a group $G$ under the condition that $G= \{\varphi(a)^{-1}a \, \mid \, a \in G \}$. Note that this set is precisely the \index{twisted conjugacy class} {\it $\varphi$-twisted conjugacy class} of the identity element of $G$ and is of independent interest. For instance, see \cite{MR3089322} for some related results. When $G$ is a finite group, this condition is equivalent to the automorphism $\varphi$ being fixed-point free. In particular, when $G$ is a  finite additive abelian group and $\varphi = -\id_G$, then this condition is equivalent to $G$ being without 2-torsion. However, when $G$ is infinite, then the two conditions are different. For example, $G=\mathbb{Q}/\mathbb{Z}$ satisfy $G=2G$, but has 2-torsion. On the other hand, $G=\mathbb{Z}$ is without 2-torsion, but $G \neq 2G$. Thus, Theorem \ref{main-theorem}, though restricted to Takasaki quandles of abelian groups, can be considered as complementary to that of Hou's result. Later, in Theorem \ref{fixed-point-free-theorem}, we prove a stronger version of Theorem \ref{main-theorem} for finite abelian groups.}
\end{remark}

The following result provides a complete characterization of the automorphism and the inner automorphism groups of dihedral quandles \cite[Theorems 2.1 and 2.2]{MR2900878}. Observe that when $n$ is odd, the result also follows from Theorem \ref{main-theorem}.

\begin{theorem}\label{Aut}
Let $\R_n$ be the dihedral quandle of order $n$. Then the following assertions hold:
\begin{enumerate}
\item $\Aut(\R_n) \cong \mathbb{Z}_n \rtimes \mathbb{Z}_n^\times$, where $\mathbb{Z}_n^\times$ is the group of units of $\mathbb{Z}_n$ when viewed as a ring.
\item   $\Inn(\R_n) \cong D_m$, the dihedral group of order $m$, where $m$ is the least common multiple of $n$ and $2$.
\end{enumerate}
\end{theorem}

\begin{proof}
Note that $\R_n= \T(\mathbb{Z}_n)$, where $\mathbb{Z}_n$ is the cyclic group of order $n$. By Proposition \ref{joyce-embedding}, there is an embedding $\mathbb{Z}_n \rtimes \Aut(\mathbb{Z}_n)\hookrightarrow \Aut(\R_n)$. Conversely, suppose that  $f \in \Aut(\R_n)$. Define $g:\R_n \rightarrow \R_n$ by $g(i)=f(i)-f(0)$. Since $f$ is a bijection, it follows that $g$ is also a bijection. Since $f(2j-i)=2f(j)-f(i)$ for all $i, j \in \R_n$, 
it follows that 
\begin{equation}\label{formula for g}
g(2j-i)= f(2j-i)-f(0)= 2f(j)-f(i)-f(0)=2g(j)-g(i).  
\end{equation}
Note that $g(0)=0$. Taking $j=0$ in \eqref{formula for g} gives $g(-i)=-g(i)$ for all $i$. Similarly, taking $i=0$ in \eqref{formula for g} gives $g(2j)=2g(j)$ for all $j$. Iterating the procedure gives $g(kj)= k g(j)$ for all even $k$. For odd multiples, we have  $g \big((2k+1)j \big)=g \big(2kj -(-j) \big)=2g(kj)-g(-j)=2kg(j)+g(j)=(2k+1)g(j)$. Hence, $g \in \Aut(\mathbb{Z}_n)$, and therefore $f=t_{f(0)}\, g$. Since $\Aut(\mathbb{Z}_n)\cong \mathbb{Z}_n^\times$, the proof of assertion (1) is complete.
\para

For $i \in \R_n$,  we have $S_i(j)= 2i-j \mod n$. Further, $S_i$ can be thought of as a reflection of a regular $n$-gon.  If $n$ is odd, the axis of symmetry of $S_i$ connects the vertex $i$ to the mid-point of the side opposite to $i$. If $n=2m$ is even, the axis of symmetry of $S_i$ passes through the opposite vertices $i$ and $(i+m) \mod 2m$. The result now follows from these observations. Observe that, if $n$ is odd, then Theorem \ref{main-theorem}(2) also shows that
$\Inn( \R_n) \cong \mathbb{Z}_n \rtimes \mathbb{Z}_2$, which establishes assertion (2).
$\blacksquare$
\end{proof}

We conclude the section with an embedding result for automorphism groups of core quandles.

\begin{proposition}\label{coreaut}
Let $G$ be a group with center $\Z(G)$. Then there is an embedding $\Z(G)\rtimes\Aut(G) \hookrightarrow \Aut(\Core(G))$.
\end{proposition}

\begin{proof} 
Recall that, in $\Core(G)$, we have $x*y=yx^{-1}y$ for all $x,y\in G$. If $\varphi \in \Aut(G)$, then by functoriality, $\varphi$ can be viewed as an automorphism of $\Core(G)$. Let $H_1=\Aut(G)$ viewed as automorphisms of $\Core(G)$. For an element $a\in \Z(G)$, let $t_a$ denote the bijection of $\Core(G)$ of the form $t_a(x)= ax$, which is clearly an automorphism of $\Core(G)$.  Let $H_2$ denote the subgroup of $\Aut(\Core(G))$ generated by $t_a$ for all $a\in \Z(G)$. Since  $\varphi(1)=1$ for each automorphism $\varphi\in H_1$, it follows that $H_1 \cap H_2= 1$. For $\varphi\in H_1$ and $t_a\in H_2$, we have $\varphi t_a\varphi^{-1}=t_{\varphi(a)}$, and hence $H_2$ is the normal subgroup of the group $H$ generated by $\{H_1, H_2\}$. Thus, $H$ is the semi-direct product $H_2\rtimes H_1$ of $H_2 \cong \Z(G)$ and $H_1 \cong \Aut(G)$.
$\blacksquare$    \end{proof}          
\para

Recall that an automorphism $\varphi$ of a group $G$ is called \index{quasi-inner automorphism} {\it quasi-inner} if for every $g \in G$, there exists $h \in G$ such that $\varphi(g) = h g h^{-1}$. In other words, $\varphi$ is  quasi-inner if for every $g \in G$, there exists an inner automorphism $\iota_h$ of $G$ such that $\varphi(g) = \iota_h(g)$. Every inner automorphism of $G$ is obviously quasi-inner.
\para

In 1913, Burnside posed the question of whether every quasi-inner automorphism of a group must necessarily be inner. He later provided a negative answer in \cite{MR1577234}, presenting an example of a finite group that admits a quasi-inner automorphism which is not inner. In the setting of quandles, this concept leads to two distinct definitions.

\begin{definition}
An automorphism $\varphi$ of a quandle $X$ is called a {\it quasi-inner automorphism in strong sense} if for every $x \in X$ there exists $y\in X$ such that $\varphi(x)=x*y$. Similarly, an automorphism $\varphi$ of a quandle $X$ is called a {\it quasi-inner automorphism in weak sense} if for every $x \in X$ there exists an inner automorphism $S \in \Inn(X)$ such that $\varphi(x)=S(x)$.
\end{definition}

Every quasi-inner automorphism in strong sense is obviously quasi-inner in weak sense. However, in general, the converse is not true since  the group of inner automorphisms $\Inn(X)=\langle S_x \,\mid \, x\in X \rangle$ does not have to coincide with the set $\{ S_x \,\mid \, x \in X \}$. If $X=\Conj(G)$ for some group $G$, then the two definitions of a quasi-inner automorphism are the same. Further, if $G$ is a $2$-step nilpotent group, then the two definitions of a quasi-inner automorphism are the same for the quandles $X=\Conj_n(G)$ for any $n$.
\para

We denote by $\QInn(X)$ the set of all automorphisms of the quandle $X$ that are quasi-inner in the weak sense. It is clear that $\QInn(X)$ forms a subgroup of $\Aut(X)$ and contains $\Inn(X)$ as a subgroup. 

\begin{lemma}\label{qaut}
If $X$ is a connected quandle, then $\Aut(X)=\QInn(X)$.
\end{lemma} 

\begin{proof}
Let $\phi$ be an automorphism of $X$. Since $X$ is connected, for each element $\phi(x) \in X$, there exists an inner automorphism $S$ such that $\phi(x)=S(x)$. Thus, we have $\Aut(X)=\QInn(X)$.
$\blacksquare$    \end{proof}

If $X$ is a trivial quandle, then  both $\Inn(X)$ and $\QInn(X)$ are trivial. An analogue of \index{Burnside's problem}{Burnside's problem} for quandles is the question of whether there exists a quandle $X$ such that $\Inn(X) \neq \QInn(X)$. Let $G$ be the group constructed by Burnside in \cite{Burnside1913} that admits a non-inner quasi-inner automorphism $\varphi$. Then the automorphism of $\Conj(G)$ induced by $\varphi$ is non-inner quasi-inner in a strong sense. This gives a negative answer to the preceding question. The following result provides an example of a quandle that is not a conjugation quandle and admits a non-inner quasi-inner automorphism.

\begin{proposition}\label{dihburn}
If $n\geq5$ is odd, then $\Inn(\R_n) \neq \QInn(\R_n)$.
\end{proposition}

\begin{proof}
By Theorem \ref{main-theorem}, if $n$ is  odd, then $\Aut( \R_n) \cong \mathbb{Z}_n\rtimes \mathbb{Z}_n^\times$, where $\mathbb{Z}_n^\times$ is the multiplicative group of the ring $\mathbb{Z}_n$. Again, by Theorem~\ref{main-theorem}, if $n$ is odd, then  $\Inn(\R_n) \cong  \mathbb{Z}_n \rtimes \mathbb{Z}_2$. Thus, if $n\ge 5$ is odd, then the groups $\Inn(\R_n)$ and $\Aut( \R_n)$ do not coincide. Since for odd $n$, the quandle $\R_n$ is connected, by Lemma \ref{qaut}, we have $\QInn(\R_n)=\Aut(\R_n)\neq \Inn(Q)$.
$\blacksquare$    \end{proof} 
\bigskip
\bigskip


\section{Automorphisms of conjugation quandles}\label{sec4}
Let $G$ be a group and $\varphi$ an automorphism of $G$. By functoriality, an automorphism of $G$ induces an automorphism of $\Conj_n(G)$ for each $n$, and hence $\Aut(G) \leq \Aut\big(\Conj_n (G)\big)$ for each $n$. The group $\Aut(G)$ is not necessarily normal in  $\Aut\left(\Conj_n (G)\right)$. For example, if $G$ is a cyclic group of prime order $p\geq5$, then the order of $\Aut(G)$ is equal to $(p-1)$. The quandle $\Conj(G)$ is trivial and $\Aut\big(\Conj(G)\big)=\Sigma_p$. The only proper normal subgroup of $\Sigma_p$ for $p\geq5$ is the alternating group $A_p$, which has order $\frac{p!}{2}>p-1$. Thus, $\Aut(G)$ cannot be normal in $\Aut \big(\Conj(G)\big)$.
\para

We shall show that $\Aut(G)=\Aut\big(\Conj(G) \big)$ if and only if $\Z(G)=1$. Throughout the section, $x^y$ means $yxy^{-1}$. First, note that 
$$\Inn \big(\Conj(G)\big) \cong \Inn(G) \cong G/\Z(G)$$ for any group $G$. Our first result is the following \cite[Proposition 4.7]{MR3718201}.

\begin{proposition}\label{aut-conj}
Let $G$ be a group with center $\Z(G)$. Then there is an embedding of groups $\Z(G) \rtimes \Aut(G) \hookrightarrow \Aut\big(\Conj(G)\big)$.
\end{proposition}

\begin{proof}
For each $a \in \Z(G)$, let $t_a: G \to G$ be given by $t_a(b)=ba$ for all $b \in G$. Then $t_a \in \Aut\big(\Conj(G))$ and the map $a \mapsto t_a$ gives an embedding $G \hookrightarrow \Aut\big(\Conj(G)\big)$.  Obviously, $\Aut(G)\subseteq  \Aut\big(\Conj(G)\big)$.
\para
It remains to show that the map $\Phi: \Z(G) \rtimes \Aut(G) \to \Aut\big(\Conj(G)\big)$ given by $\Phi(a,f) = t_a  f$ is an embedding. The map is clearly injective. Further, for $(a_1,f_1), (a_2, f_2) \in  \Z(G) \rtimes \Aut(A)$, we have
\begin{eqnarray*}
\Phi\big((a_1,f_1)(a_2, f_2) \big) & = & \Phi\big((a_1f_1(a_2), f_1  f_2) \big)\\
& = & t_{a_1f_1(a_2)}  (f_1  f_2)\\
& = & (t_{a_1}  f_1)  (t_{a_2}  f_2)\\
& = & \Phi\big((a_1,f_1)\big)  \Phi\big((a_2, f_2) \big).
\end{eqnarray*}
This completes the proof of the proposition.  $\blacksquare$
\end{proof}

It is known that $\Aut\big(\Conj(\Sigma_3)\big) \cong \Inn\big(\Conj(\Sigma_3)\big) \cong \Sigma_3$ \cite{MR2900878}. Also, it is well-known that if $n \ge 3$ and $n \neq 6$, then $\Aut\big(\Sigma_n \big) \cong \Inn\big(\Sigma_n  \big) \cong \Sigma_n$. In view of Proposition \ref{aut-conj} and the preceding discussion, it is a natural problem to classify groups $G$ for which $\Aut\big(\Conj(G)\big) \cong \Z(G) \rtimes \Aut(G)$. 
The main result of this section gives such a classification. First, we present the following result.

\begin{lemma}
Let $G$ be a group, $X=\Conj(G)$ and $\varphi \in \Aut(X)$. Then for $x,y\in X$, there exists an element $z\in \Z(G)$ such that $\varphi(xy)=\varphi(x)\varphi(y)z$.
\end{lemma}

\begin{proof} For $x,y,z\in X$, we have $$\varphi(x^{yz})=\varphi(x*{yz})=\varphi(x)*{\varphi(yz)}=\varphi(x)^{\varphi(yz)}.$$
On the other hand, we have 
\begin{multline*}\varphi(x^{yz})=\varphi \big((x^y)^z \big)=\varphi(x^y*z)=\varphi(x*y)*{\varphi(z)}=\big(\varphi(x)*\varphi(y) \big)*{\varphi(z)}=\varphi(x)^{\varphi(y)\varphi(z)}
\end{multline*}
and $\varphi(x)^{\varphi(yz)}=\varphi(x)^{\varphi(y)\varphi(z)}$. Since $x \in G$ is arbitrary, the element $\varphi(yz) \big(\varphi(y)\varphi(z)\big)^{-1}$ lie in $\Z(G)$.
$\blacksquare$    \end{proof}          

As an immediate consequence, we have the following result.

\begin{corollary}\label{aut conj trivial center}
If $G$ is a group with trivial center, then $\Aut(G)=\Aut \big(\Conj(G)\big)$.
\end{corollary}

Next, we prove the following \cite[Proposition 5]{MR3948284}.

\begin{proposition}\label{left}
If $G$ is a group with non-trivial center, then $\Aut(G)\neq \Aut \big(\Conj (G)\big)$.
\end{proposition}

\begin{proof}
Let $x\in \Z(G)$ be a non-trivial element. Let $\varphi$ denote the bijection of $X=\Conj(G)$ such that  $\varphi(1)=x$, $\varphi(x)=1$ and it 
fixes all other elements of $X$. We claim that $\varphi\in \Aut(X)$. Let $y,z \in X$. If $y\in\{1,x\}$, then $\varphi(y)\in\Z(G)$, and hence
$\varphi(y*z)=\varphi(y^z)=\varphi(y)=\varphi(y)*{\varphi(z)}$. If $y\notin\{1,x\}$, then $y^z\notin\{1,x\}$, and hence  $\varphi(y*z)=\varphi(y^z)=y^z=\varphi(y)^{\varphi(z)}=\varphi(y)*\varphi(z)$. This shows that $\varphi$ is an automorphism of $X$. Since $\varphi(1)\neq1$, the map $\varphi$ does not belong to $\Aut (G)$, and hence $\Aut(G)\neq \Aut \big(\Conj (G)\big)$.
$\blacksquare$    \end{proof}          

\begin{corollary}\label{equal}
Let $G$ be a group. Then $\Aut(G)=\Aut \big(\Conj (G)\big)$ if and only if $\Z(G)=1$.
\end{corollary}

By Proposition \ref{aut-conj},  for each group $G$, there is an embedding $\Z(G)\rtimes \Aut(G)\hookrightarrow \Aut\big(\Conj (G)\big)$. The following result shows that if $G$ is a non-abelian group, then $\Aut (G) \times  \Sigma_{\Z(G)}$ also embeds in $\Aut\big(\Conj (G)\big)$ \cite[Proposition 6]{MR3948284}.

\begin{proposition}\label{direct aut conj}
Let $G$ be a non-abelian group. Then the direct product $\Aut (G) \times  \Sigma_{\Z(G)}$ embeds in $\Aut \big(\Conj(G) \big)$.
\end{proposition}

\begin{proof} 
Consider the map $f:\Aut(G)\to\Aut \big(\Conj(G) \big)$ which sends an automorphism $\varphi \in \Aut(G)$ to the map $f(\varphi):\Conj (G) \to  \Conj(G)$, defined as
$$f(\varphi)(x)=\begin{cases}
\varphi(x)&~\textrm{if}~ x\notin \Z(G),\\
x& ~\textrm{if}~x\in \Z(G).
\end{cases}$$
Since $\varphi\in \Aut (G)$, it maps $\Z(G)$ to $\Z(G)$ and maps $G\setminus \Z(G)$ to $G\setminus \Z(G)$. Hence, $f(\varphi)$ is a bijection on $\Conj (G)$. A direct check shows that $f(\varphi)$ preserves the quandle operation in $\Conj (G)$, and hence  $f(\varphi) \in \Aut \big(\Conj (G)\big)$.
\para 
For automorphisms $\varphi,\psi\in \Aut(G)$ and an element $x\in \Z(G)$, we have $f(\varphi\psi)(x)=x=f(\varphi)f(\psi)(x)$. If $x\notin \Z(G)$, then $\psi(x)\notin \Z(G)$, and hence $f(\varphi\psi)(x)=\varphi \big(\psi(x)\big)=f (\varphi) \big(f(\psi)(x)\big)=f(\varphi)f(\psi)(x)$. Thus, $f$ is a group homomorphism. If $f(\varphi)=f(\psi)$, then for all $x\notin \Z(G)$, we have $\varphi(x)=f(\varphi)(x)=f(\psi)(x)=\psi(x)$. In particular, if $x \in G\setminus \Z(G)$ and $y\in \Z(G)$, then $\varphi(x)=\psi(x)$ and $\varphi(xy)=\psi(xy)$. This gives $\varphi(y)=\psi(y)$, and hence $\varphi=\psi$. Hence, $f: \Aut(G) \to \Aut \big(\Conj(G)\big)$ is injective.
\para
Consider the map $g:\Sigma_{\Z(G)}\to \Aut \big(\Conj(G) \big)$ which sends a permutation $\sigma$ of $\Z(G)$ to the map $g(\sigma): \Conj (G) \to  \Conj(G)$, defined as
$$g(\sigma)(x)=\begin{cases}
x&~\textrm{if}~   x\notin \Z(G),\\
\sigma(x)&~\textrm{if}~ x\in \Z(G).
\end{cases}$$
As in the case of $f$, it is straightforward to verify that $g:\Sigma_{\Z(G)} \to \Aut \big(\Conj (G) \big)$ is an injective group homomorphism. Let $H$ be a subgroup of $\Aut \big(\Conj(G)\big)$ generated by $H_1=f \big(\Aut(G) \big)$ and $H_2=g( \Sigma_{\Z(G)})$. Both $H_1,H_2$ are normal subgroups of $H$ (moreover, $H_2$ is normal in $\Aut \big(\Conj(G)\big)$) and $H_1\cap H_2= 1$. Hence, we have $H=H_1\times H_2 \cong \Aut(G)\times \Sigma_{\Z(G)}$.
$\blacksquare$    \end{proof}          

\begin{corollary} 
Let $G$ be a group. Then $\Z(G)$ embeds in $\Aut \big(\Conj (G)\big)$.
\end{corollary}

\begin{proof} 
Since $\Z(G)$ acts faithfully on itself by right multiplications, it follows that $\Z(G)$ is a subgroup of $\Sigma_{\Z(G)}$. If $G$ is abelian, then $\Conj(G)$ is  the trivial quandle and  $\Aut \big(\Conj(G)\big) \cong \Sigma_{G}$. If $G$ is not abelian, then by Proposition \ref{direct aut conj}, the group $\Sigma_{\Z(G)}$ is a subgroup of $\Aut \big(\Conj (G)\big)$.
$\blacksquare$    \end{proof}          

The next result proves an analogue of Proposition \ref{direct aut conj} for finite abelian groups.

\begin{proposition}\label{abbb} 
Let $G$ be a non-trivial finite abelian group. Then $\Aut(G)\times\Sigma_{G}$ embeds in $\Aut \big(\Conj(G)\big)$ if and only if $G \cong \mathbb{Z}_2$.
\end{proposition}

\begin{proof} 
Note that conjugation quandles of abelian groups are trivial quandles. Clearly, if $G\cong \mathbb{Z}_2$, then $\Aut \big(\Conj(G)\big) \cong  \Aut(G)\times\Sigma_{G}$. Conversely, suppose that $\Aut(G)\times\Sigma_{G}$ embeds in $\Aut \big(\Conj(G)\big)$. Since $G$ is an abelian group, we have $\Aut \big(\Conj(G)\big) =\Sigma_{G}$, which has order $|G|!$. But, ${\Aut(G)}\times\Sigma_{G}$ has order $|\Aut(G)| (|G|!)$. This implies that $|\Aut(G)|=1$, and hence $G \cong \mathbb{Z}_2$. 
$\blacksquare$    \end{proof}           
 
 The following result gives a description of all finite groups $G$ for which $\Aut \big(\Conj(G) \big) \cong \Aut(G)\times \Sigma_{\Z(G)}$ \cite[Theorem 1]{MR3948284}.

\begin{theorem} Let $G$ be a finite group. Then $\Aut\big(\Conj(G)\big) \cong \Aut(G)\times \Sigma_{\Z(G)}$ if and only if $\Z(G)$ is trivial or $G \cong \mathbb{Z}_2$.
\end{theorem}

\begin{proof}
Corollary \ref{equal} and Proposition \ref{abbb} imply that if $\Z(G)$ is trivial or $G \cong \mathbb{Z}_2$, then $\Aut \big(\Conj(G)\big) \cong \Aut(G)\times \Sigma_{\Z(G)}$.
\para

For the converse, suppose that $\Aut \big(\Conj(G)\big) \cong \Aut(G)\times \Sigma_{\Z(G)}$. If $\Z(G)$ is trivial, then the assertion follows from Corollary \ref{equal}. Suppose that $\Z(G)$ is non-trivial. Then, according to the proof of Proposition \ref{direct aut conj}, the group $\Aut(G)\times \Sigma_{\Z(G)}$ consists of automorphisms of the form
\begin{equation}\label{gfo}
x\mapsto\begin{cases}
\varphi(x)& ~\textrm{if}~x\notin \Z(G),\\
\sigma(x)&~\textrm{if}~ x\in \Z(G),
\end{cases}
\end{equation}
where $\varphi$ is an automorphism of $G$ and $\sigma$ is a permutation of $\Z(G)$. For an element $1 \neq a\in \Z(G)$, denote by $t_a$ the automorphism of $\Conj(G)$ of the form $t_a(x)= xa$. Then the map $t_a$ must have the form \eqref{gfo}. If there exist  elements $x,y \in G \setminus \Z(G)$ such that $xy \in G \setminus \Z(G)$, then \eqref{gfo} gives $$\varphi(x)\varphi(y)=\varphi(xy)=t_a(xy)=axy.$$ On the other hand,  from \eqref{gfo}, we have $\varphi(x)=ax$ and $\varphi(y)=ay$. This implies that $a=1$, which is a contradiction. Thus, we have proved that for all $x,y \in G \setminus \Z(G)$, the product $xy$ lie in $\Z(G)$. This means that $|G/\Z(G)|\leq 2$, and $G$ is abelian. Hence, by Proposition \ref{abbb}, we have $G \cong \mathbb{Z}_2$.
$\blacksquare$   
\end{proof}

By Proposition \ref{aut-conj}, the semi-direct product $\Z(G)\rtimes \Aut(G)$ embeds in $\Aut\big(\Conj (G)\big)$. The following result describes all finite groups $G$ satisfying $\Aut\big(\Conj (G)\big) \cong \Z(G)\rtimes \Aut(G)$ \cite[Theorem 2]{MR3948284}.

\begin{theorem}\label{solution}
Let $G$ be a finite group. Then $\Aut \big(\Conj (G)\big) \cong \Z(G)\rtimes \Aut(G)$ if and only if either $\Z(G)$ is trivial or $G \cong \mathbb{Z}_2$, $\mathbb{Z}_2^2$ or $\mathbb{Z}_3$.
\end{theorem}

\begin{proof} By Corollary \ref{equal}, if $\Z(G)$ is trivial, then $\Aut \big(\Conj (G)\big)=\Aut (G) \cong \Z(G)\rtimes \Aut (G)$. For the groups $\mathbb{Z}_2$, $\mathbb{Z}_2^2$ and $\mathbb{Z}_3$, we see that 
\begin{align}
\notag \Aut \big(\Conj(\mathbb{Z}_2)\big)&\cong \Sigma_2\cong\mathbb{Z}_2\cong\mathbb{Z}_2\rtimes \Aut(\mathbb{Z}_2),\\
\notag \Aut \big(\Conj(\mathbb{Z}_2^2)\big)&\cong \Sigma_4\cong\mathbb{Z}_2^2\rtimes\Sigma_3\cong\mathbb{Z}_2^2\rtimes\Aut(\mathbb{Z}_2^2),\\
\notag \Aut \big(\Conj(\mathbb{Z}_3)\big)&\cong \Sigma_3\cong\mathbb{Z}_3\rtimes\Sigma_2\cong\mathbb{Z}_3\rtimes\Aut(\mathbb{Z}_3).
\end{align}

Conversely, suppose that $G$ is finite and $\Aut \big(\Conj (G)\big) \cong \Z(G)\rtimes \Aut(G)$. For an element $a\in \Z(G)$, denote by $t_a$ the automorphism of the quandle $\Conj(G)$ of the form $t_a(x)=xa$. The group generated by the set $\{t_a \, \mid \, a\in \Z(G) \}$ is isomorphic to $\Z(G)$. In view of Proposition \ref{aut-conj}, the subgroup of $\Aut \big(\Conj(G)\big)$ which is isomorphic to $\Z(G)\rtimes  \Aut(G)$ consists of the automorphisms of the form $t_a\varphi$ for some $a\in \Z(G)$ and $\varphi\in\Aut(G)$. Thus, we need to prove that if every automorphism of $\Conj(G)$ has the form $t_a\varphi$ for some $a\in \Z(G)$ and $\varphi\in\Aut(G)$, then either  $\Z(G)$ is trivial or $G \cong \mathbb{Z}_2$, $\mathbb{Z}_2^2$ or $\mathbb{Z}_3$. 
\para

Suppose that $\Z(G)$ is non-trivial. For a permutation $\sigma$ of  $\Z(G)$, denote by $f$ the automorphism of $\Conj(G)$ of the form
$$f(x)=\begin{cases}
x& ~\textrm{if}~x\notin \Z(G),\\
\sigma(x)&~\textrm{if}~ x\in \Z(G).
\end{cases}$$
Suppose that $f=t_a\varphi$ for some $a\in \Z(G)$ and $\varphi\in \Aut(G)$. If there exist two elements $x,y\in G \setminus \Z(G)$ such that $xy \in G \setminus \Z(G)$, then we see that
\begin{eqnarray}
\notag x&=&f(x)=t_a\varphi(x)=a\varphi(x),\\
\notag y&=&f(y)=t_a\varphi(y)=a\varphi(y),\\
\notag xy&=&f(xy)=t_a\varphi(xy)=a\varphi(x)\varphi(y).
\end{eqnarray}
These equalities imply that $a=1$, and hence $f=\varphi$. Since $\Z(G)\neq 1$, we choose $\sigma$ such that $\sigma(1)\neq1$. In this case, $f$ cannot belong to $\Aut(G)$. Thus, for any two elements $x,y \in G \setminus \Z(G)$, their product $xy$ must lie in $\Z(G)$. This implies that $|G/\Z(G)|\leq 2$, and hence $G$ is abelian. Consequently, we have $|\Aut \big(\Conj(G)\big)|=|\Sigma_{G}|=|G|!$. Since every group automorphism of $G$ fixes the identity element, it follows that  $|\Aut(G)|\leq \big(|G|-1 \big)!$. Further, since $\Aut \big(\Conj (G)\big) \cong \Z(G)\rtimes \Aut(G)=G\rtimes \Aut(G)$, we get $|\Aut(G)|=\big(|G|-1 \big)!$. This means that every permutation of non-trivial elements of $G$ is a group automorphism of $G$. In particular, all non-trivial elements of $G$ have the same order. Suppose that $G$ has at least five distinct elements $\{1,x,y,xy,t\}$. The map $\varphi$ with $\varphi(1)=1$, $\varphi(x)=x$, $\varphi(y)=y$ and $\varphi(xy)=t$ can be extended to a group automorphism of $G$. This gives  $t=\varphi(xy)=\varphi(x)\varphi(y)=xy$, which is a contradiction. Hence,  $G$ has at most $4$ elements. Since all non-trivial elements of $G$ have the same order, it follows that $G \cong \mathbb{Z}_2$, $\mathbb{Z}_2^2$ or $\mathbb{Z}_3$.
$\blacksquare$   
\end{proof}          
\bigskip
\bigskip 

\section{Multi-transitive action of automorphism groups}\label{section6 doubly transitive}
Let $G$ be a group acting on a set $X$ from the left. For each $x \in X$, the subgroup $G_x= \{g \in G \, \mid \, g \cdot x=x \}$ is called the \index{stabilizer subgroup} {\it stabilizer subgroup} at $x$. For each $1 \leq k \leq |X|$, we say that $G$ acts \index{$k$-transitive group}{\it $k$-transitively} on $X$ if for each pair of $k$-tuples $(x_1,\dots,x_k)$ and $(y_1,\dots,y_k)$ of distinct elements of $X$, there exists $g \in G$ such that $g \cdot x_i = y_i$ for each $1 \leq i \leq k$. A 1-transitive action is simply called {\it transitive} and a 2-transitive action is also called {\it doubly transitive}.
\para

Let $X$ be a quandle and $\Inn(X)$ its inner automorphism group. For an integer $k\geq 1$, the quandle $X$ is called \index{$k$-transitive quandle} {\it $k$-transitive} if the group $\Inn(X)$ acts $k$-transitively on $X$. A 1-transitive quandle is simply  a connected quandle. A 2-transitive quandle is also called a \index{two-point homogeneous quandle} {\it two-point homogeneous quandle} in \cite{MR3127819}, where a classification of all two-point homogeneous quandles of prime order was given. A complete classification of finite two-point homogeneous quandles was later completed by Vendramin \cite[Theorem 3]{MR3685034} and Wada \cite{MR3379001}.
\para

McCarron  proved that if $X$ is a finite $k$-transitive quandle with $k\geq 2$ and with at least four elements, then $k = 2$ \cite[Proposition~5]{MR3206312}. Moreover, the dihedral quandle $\R_3$ of order three is the only $3$-transitive quandle. Thus, higher order transitivity does not exist in quandles with at least four elements. In view of the preceding remark, the following problem seems natural.

\begin{problem}
For an integer $k \geq 2$, classify all finite quandles $X$ for which $\Aut(X)$ acts $k$-transitively on $X$.
\end{problem}

We generalise a result of Ferman, Nowik, and Teicher \cite{MR2795262} about doubly transitive action of the automorphism group of a finite Alexander quandle of prime order on the underlying quandle. An alternate proof of the same result was given by Watanabe \cite{MR3319676}. Our approach is more general than \cite{MR2795262} and \cite{MR3319676}. 
\para 

Let $G$ be a finite abelian group, $\varphi \in \Aut(G)$ and 
$$\Aut\big(\Alex(G, \varphi)\big)_1=\big\{f \in  \Aut\big(\Alex(G, \varphi)\big) \, \mid \,  f(1)=1 \big\}$$
 the stabilizer subgroup of $\Aut\big(\Alex(G, \varphi)\big)$ at the identity element $1 \in G$. We first prove a stronger version of Theorem \ref{main-theorem} for finite abelian groups with fixed-point free automorphisms.

\begin{theorem}\label{fixed-point-free-theorem}
Let $G$ be a finite abelian group and $\varphi \in \Aut(G)$ a fixed-point free automorphism. Then the following  assertions hold:
\begin{enumerate}
\item $\C_{\Aut(G)}(\varphi)= \Aut\big(\Alex(G, \varphi)\big)_1$.
\item $ \Aut\big(\Alex(G, \varphi)\big) \cong G \rtimes \C_{\Aut(G)}(\varphi)$.
\item $ \Inn\big(\Alex(G, \varphi)\big) \cong G \rtimes \langle \varphi \rangle$.
\end{enumerate}
\end{theorem}

\begin{proof}
Let $f \in \Aut\big(\Alex(G, \varphi)\big)_1$. Then, for $a, b \in G$, we have
\begin{equation}\label{alex aut fix point free}
f\big(\varphi(a)\varphi(b)^{-1} b \big)=\varphi\big(f(a)\big)\varphi\big(f(b)\big)^{-1}f(b).
\end{equation}
Taking $a=1$ gives $f\big(\varphi(b)^{-1}b \big)=\varphi\big(f(b)\big)^{-1}f(b)$ and taking $b=1$ gives $f\big(\varphi(a)\big)=\varphi\big(f(a)\big)$. Using these identities in \eqref{alex aut fix point free}, we obtain
\begin{equation}\label{alex aut fix point free 2}
f\big(\varphi(a)\tilde{\varphi}(b)^{-1} \big)=f\big(\varphi(a)\varphi(b)^{-1}b \big)=f\big(\varphi(a)\big) f\big(\varphi(b)^{-1}b \big)= f\big(\varphi(a)\big) f\big(\tilde{\varphi}(b)^{-1} \big).
\end{equation}
Since $\varphi$ is fixed-point free, it follows from Theorem \ref{Bae-Choe-theorem} that $\tilde{\varphi} \in \Aut(G)$. Thus, \eqref{alex aut fix point free 2} gives $f(xy)=f(x)f(y)$ for all $x, y \in G$, and hence $f \in \C_{\Aut(G)}(\varphi)$. Conversely, if $f \in \C_{\Aut(G)}(\varphi)$, then $f(1)=1$ and $f(a*b)=f(a)*f(b)$. This proves assertion (1).
\para 
By Proposition \ref{joyce-embedding}, there is an embedding $\Phi:G \rtimes \C_{\Aut(G)}(\varphi)  \hookrightarrow \Aut\big(\Alex(G, \varphi)\big)$ given by $\Phi(a,f) = t_a  f$, where $t_a(x)=xa$. Let $f \in \Aut\big(\Alex(G, \varphi)\big)$. Define $h: G \to G$ by $h(a)=f(a)f(1)^{-1}$ for all $a \in G$. Then $h \in \Aut\big(\Alex(G, \varphi)\big)$ and $h(1)=1$. Hence, $h \in \C_{\Aut(G)}(\varphi)$ and $f= t_{f(1)} h \in \Phi\big(G \rtimes \C_{\Aut(G)}(\varphi) \big)$, which proves assertion (2).

Since $\varphi$ is fixed-point free, by Theorem \ref{Bae-Choe-theorem}, we have $G= \im(\tilde{\varphi})= \{ a^{-1}\varphi(a) \, \mid \, a \in G \}$. Thus, for each $a \in G$, $S_a=t_{\varphi(a)^{-1}a} \varphi \in \Phi(G \rtimes \langle \varphi \rangle)$. Conversely, since each $a \in G$ can be written as $a=\varphi(b)^{-1}b$ for some $b \in G$, we have $t_a \varphi=S_{b}$. This proves assertion (3).
$\blacksquare$    
\end{proof}          

The result below is well-known and follows easily.

\begin{lemma}\label{doubly-transitive-criteria}
Let $G$ be a group acting on a set $X$ and $x \in X$. Then the action is doubly transitive if and only if it is transitive and the stabilizer subgroup $G_x$ acts transitively on $X\setminus \{x\}$.
\end{lemma}

The following result is presumably well-known; however, we include a proof for the sake of completeness.

\begin{lemma}\label{transitive-criteria}
Let $G$ be a non-trivial finite $p$-group, where $p$ is a prime. Then $\Aut(G)$ acts transitively on $G \setminus \{1\}$ if and only if $G \cong (\mathbb{Z}_p)^n$ for some integer $n \geq 1$.
\end{lemma}

\begin{proof}
Suppose that $\Aut(G)$ acts transitively on $G \setminus \{1\}$. By Cauchy's Theorem for finite groups, there exists an element $a \in G$ of order $p$.  Given any $1 \ne b \in G$, there exists $\varphi \in \Aut(G)$ such that $\varphi(a)=b$. This implies that all non-trivial elements of $G$ have order $p$. Since $G$ is a non-trivial $p$-group, its center $\Z(G)$ is non-trivial. For each $1 \ne a \in G$ and $1 \ne b \in \Z(G)$, there exists $\varphi \in \Aut(G)$ such that $\varphi(b)=a$. Since  $\Z(G)$ is characteristic, it follows that  $a\in \Z(G)$. Hence, $G\cong  (\mathbb{Z}_p)^n$ for some integer $n \geq 1$. The converse is obvious, and hence omitted.
$\blacksquare$  
  \end{proof}          

The next result generalises \cite[Theorem 3.9]{MR2795262} to elementary abelian $p$-groups.

\begin{theorem}\label{FNT-generalisation}
Let $G=  (\mathbb{Z}_p)^n$ for some prime $p$ and integer $n \geq 1$. If $\varphi$ is multiplication by a non-trivial unit of $\mathbb{Z}_p$, then $\Aut\big(\Alex(G, \varphi)\big)$ acts doubly transitively on $\Alex(G, \varphi)$.
\end{theorem}

\begin{proof}
Since $\varphi$ is fixed-point free, by Theorem \ref{Bae-Choe-theorem}, the quandle $\Alex(G, \varphi)$ is connected. This implies that $\Inn\big(\Alex(G, \varphi)\big)$, and hence $\Aut\big(\Alex(G, \varphi)\big)$ acts transitively on $\Alex(G, \varphi)$. In view of Lemma \ref{doubly-transitive-criteria}, it suffices to prove that $\Aut\big(\Alex(G, \varphi)\big)_0$ acts transitively on $\Alex(G, \varphi) \setminus \{0\}$. By Theorem \ref{fixed-point-free-theorem},  $\Aut\big(\Alex(G, \varphi)\big)_0=\C_{\Aut(G)}(\varphi)$. But, $\C_{\Aut(G)}(\varphi)= \Aut(G)$, which by Lemma \ref{transitive-criteria} acts transitively on $\Alex(G, \varphi) \setminus \{0\}$.
$\blacksquare$    
\end{proof}          

\begin{remark}{\rm 
It is worth pointing out that $\Alex(G, \varphi)$ is not two-point homogeneous for $G =  (\mathbb{Z}_p)^n$ with $n \geq 2$. This follows from results of Vendramin \cite[Theorem 3]{MR3685034}  and Wada \cite[Corollary 4.5]{MR3379001}, where they classify all finite two-point homogeneous quandles. In fact, they proved that any such quandle is isomorphic to an Alexander quandle defined by primitive roots over a finite field.}
\end{remark}

Next, we define the center of a quandle. 

\begin{definition}
Let $X$ be a quandle. The subset $\Z(X)= \big\{x\in X \, \mid \, x*y=x~\text{for all}~y\in X \big\}$ is called the \index{center of a quandle} {\it center} of the quandle $X$.
\end{definition}

 If $X=\Conj(G)$, then the center $\Z(X)$ of the quandle $X$ coincides with the center $\Z(G)$ of the group $G$.

\begin{lemma}\label{center}
Let $X$ be a quandle with center $\Z(X)$ and $\varphi$ an automorphism of  $X$. Then $\varphi\big(\Z(X)\big)=\Z(X)$.
\end{lemma}
\begin{proof}
If $x\in \Z(X)$, then for all $y\in X$, we have $\varphi(x)*\varphi(y)=\varphi(x*y)=\varphi(x)$. Thus, we have $\varphi(x)\in \Z(X)$. Conversely, if $\varphi(x)\in \Z(X)$, then for all $y\in X$, we have $\varphi(x*y)=\varphi(x)*{\varphi(y)}=\varphi(x)$. Hence, $x*y=x$ for all $y\in X$, and therefore $x\in \Z(X)$.
$\blacksquare$    \end{proof}

\begin{theorem}\cite[Theorem 5]{MR3948284}\label{tran}
Let $X$ be a finite quandle. Then the following statements are equivalent:
\begin{enumerate}
\item $X$ is either the trivial quandle or $X \cong \R_3$.
\item $\Aut(X)=\Sigma_{X}$. 
\item $\Aut(X)$ acts $3$-transitively on $X$.
\end{enumerate}
\end{theorem}
\begin{proof}
The implications $(1)\Leftrightarrow(2)$ and $(2)\Rightarrow(3)$ are obvious. It remains to prove the implication $(3)\Rightarrow(1)$. 
\para 
Suppose that $\Z(X)\neq\varnothing$. Since $\Aut(X)$ acts $3$-transitively (and therefore $1$-transitively) on $X$, by Lemma \ref{center}, the quandle $X$ is trivial.
\para 
Suppose that $\Z(X)=\varnothing$. Then $|X|\geq3$ and for an element $x\in X$, there exists an element $y\in X$ such that $z:=x*y\neq x$. Since $x\neq y$, we have $z=x*y\neq y*y=y$, and hence $x,y,z$ are distinct elements. Let $t$ be an arbitrary element of $X$. Since $\Aut(X)$ acts $3$-transitively on $X$, if $t\notin\{x,y\}$, then there exists an automorphism $\varphi$ of $X$ such that $\varphi(x)=x$, $\varphi(y)=y$ and $\varphi(z)=t$. This gives $\varphi(z)=\varphi(x*y)=\varphi(x)*\varphi(y)=x*y=z$, and hence  $X=\{x,y,z\}$.
\para
Next, we identify the quandle structure on $X$. Since $\Aut(X)$ acts $3$-transitively on $X$, there exists an automorphism $\varphi$ such that $\varphi(x)=y$, $\varphi(y)=x$ and $\varphi(z)=z$. Thus, $y*x=\varphi(x)*{\varphi(y)}=\varphi(x*y)=\varphi(z)=z=x*y$.
Again, since $\Aut(X)$ acts $3$-transitively on $X$, there exists an automorphism $\psi$ such that $\psi(x)=z$, $\psi(y)=x$ and $\psi(z)=y$, and hence $y=\psi(z)=\psi(x*y)=\psi(x)*{\psi(y)}=z*x$. Repeating this argument for all possible triples constructed from the set $X=\{x,y,z\}$, we obtain the  relations
$$x*y=y*x=z,~x*z=z*x=y~\textrm{and}~y*z=z*y=x.$$
Hence, $X=\{x,y,z\}$ is isomorphic to the dihedral quandle $\R_3$.
$\blacksquare$    
\end{proof}          

\begin{corollary}
If $X$ is a finite quandle such that $\Aut(X)$ acts $3$-transitively on $X$, then it acts $k$-transitively on $X$ for each $k\geq3$.
\end{corollary}
\bigskip
\bigskip 


\section{Automorphisms of quandle extensions}\label{sec7}

In Chapter \ref{chapter homology YBE}, we will delve into the details of the homology and cohomology of racks and quandles. For this section, we revisit only the essential information. Let $X$ be a quandle and $S$ a set.  A map $\alpha : X\times X \to \Sigma_{S}$ is called \index{constant quandle 2-cocycle} {\it a constant quandle 2-cocycle} if it satisfies the following two conditions:
\begin{enumerate}
\item $\alpha(x*y,z)\, \alpha(x,y) = \alpha(x*z,y*z) \, \alpha(x,z)$ for all $x, y, z \in X$.
\item $\alpha(x,x) = \id_S$ for all $x \in X$.
\end{enumerate}

A 2-cocycle satisfying only condition (1) is called a \index{constant rack 2-cocycle} {\it constant rack 2-cocycle}. The following example gives a way to construct such 2-cocycles, which resembles the classification of \index{Yetter--Drinfel'd module}{Yetter--Drinfel'd modules} over group algebras \cite{MR1994219}.

\begin{example}\label{constant quandle cocycle examaple}
{\rm 
Let $X$ be a connected rack, $x_0\in X$ a fixed element, $G=\Inn(X)$ and  $H$ the subgroup of $G$ which fix $x_0$. Under the natural action of $G$ on $X$, the orbit map $g\mapsto g \cdot x_0$ induces a bijection $G/H\to X$. Fix a set-theoretic section $\tau:X\to G$ to the orbit map, that is, $\tau(x)\cdot x_0=x$ for all $x\in X$. Note that
$$\big(\tau(x* y)^{-1}S_y \tau(x) \big)\cdot x_0	=\tau(x* y)^{-1}S_y\cdot x	=\tau(x* y)^{-1}\cdot(x* y)=x_0$$
for all $x, y \in X$. This determines, for each $x,y\in X$, an element $h_{x,y}\in H$ such that $S_y \tau(x)=\tau(x* y)h_{x,y}$. Let $Z$ be a set and $\rho:H\to\Sigma_Z$ a group homomorphism. Then it is straightforward to see that the map $\beta:X\times X\to\Sigma_Z$ given by $\beta(x,y)=\rho(h_{x,y})$ is a constant rack 2-cocycle.}
\end{example}

Let $X$ be a quandle and $S$ a set. For a constant quandle 2-cocycle $\alpha:X\times X\to\Sigma_{S}$, consider the set $X\times S=\{(x,t)\, \mid \,x\in X,~t\in S\}$ with the binary operation $*$ given by
\begin{equation}\label{non-abelian-structure}
(x, t) * (y, s) = \big(x * y, ~\alpha(x,y)(t)\big)
\end{equation}
for $x, y \in X$ and $s, t \in S$. This turns the set $X\times S$ into a quandle called the \index{non-abelian extension}{\it non-abelian extension of $X$ by $S$ through $\alpha$}, and we denote it by $X \times_\alpha S$. There is a surjective quandle homomorphism  $X \times_\alpha S \to X$, namely, the projection onto the first coordinate. On the other hand, for each $x \in X$, the set $\{ (x,t) \, \mid \, t \in S \}$ is a trivial subquandle of $X \times_\alpha S$.  Such extensions were introduced by Andruskiewitsch and Gra\~{n}a in \cite{MR1994219}. We shall see a further generalisation of this construction in Subsection \ref{sec extension theory of quandles} of Chapter
\ref{chapter homology YBE}.

\begin{example}{\rm 
Let $X$ be a connected rack,  $Z$ a set, $\rho:H\to\Sigma_Z$ a group homomorphism and $\beta$ the constant rack 2-cocycle of Example \ref{constant quandle cocycle examaple}. Then the non-abelian extension  $X\times_{\beta}Z$ has the binary operation
$$(x,t) *(y,s)=\big(x * y,~ \rho(\tau(x * y)^{-1}S_y \tau(x))(t)\big)$$
for $x, y \in X$ and $s, t \in Z$. Even if $X$ is a quandle, $X\times_{\beta}Z$ is not a quandle in general, since it does not necessarily satisfy condition (2) in the definition of a constant quandle 2-cocycle. Observe that
$$\beta(x,x)=\rho \big(\tau(x * x)^{-1}S_x \tau(x) \big)
	=\rho \big(\tau(x)^{-1}S_x \tau(x) \big)=\rho(S_{\tau(x)^{-1}\cdot x})
	=\rho(S_{x_0})$$
for all $x \in X$. Thus, $X\times_{\beta}Z$ is a quandle if and only if $X$ is a quandle and $\rho(S_{x_0})=\id_Z$.}
\end{example}
\para

Given a quandle $X$ and a set $S$, we denote the set of all constant quandle 2-cocycles by $\mathcal{C}^2(X,S)$.  Two constant quandle 2-cocycles $\alpha, \beta\in \mathcal{C}^2(X,S)$ are said to be \index{cohomologous}{\it cohomologous} if there exists a map $\lambda:X \to \Sigma_{S}$ such that $$\alpha(x,y)=\lambda(x*y)^{-1} \beta(x,y)\lambda(x)$$ for all $x,y \in X$. The relation of being cohomologous is an equivalence relation on $\mathcal{C}^2(X,S)$. The equivalence class of a constant quandle 2-cocycle $\alpha$ is called \textit{the cohomology class of $\alpha$} and is denoted by $[\alpha]$. Let $\mathcal{H}^2(X,S)$ denote the set of cohomology classes of constant quandle 2-cocycles. We shall see a detailed exposition of quandle (co)homology in Part-III. The following result generalizes the result in \cite[Lemma 4.8]{MR1885217}, where it is formulated for Alexander quandles.

\begin{lemma}\label{cohomologous-isomorphic}
Let $X$ be a quandle and $S$ a set. If constant quandle 2-cocycles  $\alpha, \beta:X\times X\to \Sigma_{S}$ are cohomologous, then the quandles $X \times_\alpha S$ and $X \times_\beta S$ are isomorphic.
\end{lemma}

\begin{proof} 
We denote the quandle operations in $X\times_{\alpha}S$ and  $X\times_{\beta}S$ by $*$ and $\circ$, respectively. Since $\alpha$ and  $\beta$ are cohomologous, there exists a map $\lambda:X \to \Sigma_{S}$ such that $\alpha(x,y)=\lambda(x*y)^{-1} \beta(x,y)\lambda(x)$ for all $x,y \in X$. Define $f: X \times_\alpha S \to X \times_\beta S$ by $f \big((x,t)\big)=\big(x, \lambda(x)(t)\big)$. The map $f$ is obviously a bijection. Further, for $(x,t), (y,s) \in X \times_\alpha S$, we have
\begin{align}
\notag f \big((x,t)*(y,s) \big) &= f\big(x*y, \alpha(x,y)(t)\big) =  \big(x*y, \lambda(x*y) \alpha(x,y)(t)\big)\\
\notag&  =  \big(x*y, \beta(x,y)\lambda(x)(t) \big)  =  \big(x,\lambda(x)(t)\big)\circ \big(y,\lambda(y)(s)\big)\\
\notag&  =  f (x,t)\circ f (y,s),
\end{align}
and hence $f$ is an isomorphism.
$\blacksquare$    \end{proof}          

\begin{lemma}\label{aut-constant-action}
Let $X$ be a quandle, $S$ a set and $\alpha$ a constant quandle 2-cocycle. Then for $(\phi, \theta) \in \Aut(X)\times \Sigma_{S}$, the map $^{(\phi, \theta)}:\mathcal{C}^2(X,S)\to\mathcal{C}^2(X,S)$ given by $$^{(\phi,\theta)}\alpha(x,y)=\theta  \alpha \big(\phi^{-1}(x),~ \phi^{-1}(y) \big) \theta^{-1} $$
induces a left-action of $\Aut(X) \times \Sigma_{S}$ on $\mathcal{H}^2(X,S)$.
\end{lemma}

\begin{proof} It is easy to check that $^{(\phi,\theta)}\alpha$ is a constant quandle 2-cocycle and the map $^{(\phi, \theta)}:\mathcal{C}^2(X,S)\to\mathcal{C}^2(X,S)$ is a left-action of $\Aut(X)\times\Sigma_{S}$ on $\mathcal{C}^2(X,S)$. If $\alpha,\beta$ are cohomologous constant 2-cocycles, then there exists a map $\lambda:X \to \Sigma_{S}$ such that $\alpha(x,y)= \lambda(x*y)^{-1} \beta(x,y)\lambda(x)$ for all $x, y \in X$. Acting on this equality by the map $^{(\phi, \theta)}:\mathcal{C}^2(X,S)\to\mathcal{C}^2(X,S)$, we obtain
$$^{(\phi,\theta)}\alpha(x,y)= {\lambda^{\prime}(x*y)^{-1}} ~^{(\phi,\theta)}\beta(x,y)\, \lambda^{\prime}(x),$$ 
where the map $\lambda^{\prime}:Q \to\Sigma_{S}$ is given by $\lambda'(x)=\theta \lambda \big(\phi^{-1}(x) \big)\theta^{-1}$. Thus, $^{(\phi,\theta)}\alpha$ and $^{(\phi,\theta)}\beta$ are cohomologous constant quandle 2-cocycles. Consequently, $(\phi,\theta)$ maps the class $[\alpha]$ to the class $\left[^{(\phi,\theta)}\alpha\right]$, and gives a left-action of $\Aut(X) \times \Sigma_{S}$ on $\mathcal{H}^2(X,S)$.
$\blacksquare$    \end{proof}

The following result gives a relation between the groups $\Aut(X)\times\Sigma_{S}$ and $\Aut(X\times_{\alpha}S)$.

\begin{theorem}\label{exterrr}
Let $X$ be a quandle, $S$ a set and $\alpha:X\times X\to \Sigma_{S}$ a constant quandle 2-cocycle. Then there exists an embedding $\big(\Aut(X)\times \Sigma_{S} \big)_{\alpha}\hookrightarrow \Aut(X\times_{\alpha}S)$.
\end{theorem}

\begin{proof}
For $(\phi,\theta) \in \big(\Aut(X)\times\Sigma_{S} \big)_{\alpha}$, let $\gamma: X \times_\alpha S \to X \times_\alpha S$ be the map defined as $\gamma(x, t)= \big(\phi(x), \theta(t)\big)$ for all $x \in X$ and $t \in S$.  The map $\gamma$ is clearly a bijection. Since $^{(\phi, \theta)}\alpha=\alpha$, by definition of the action of $\Aut(X)\times\Sigma_{S}$, we have  $\theta  \alpha \big(\phi^{-1}(x),~ \phi^{-1}(y) \big) \theta^{-1}=\alpha(x,~y)$ for all $x, y \in X$. For $x, y \in X$ and  $t, s \in S$, we see that
\begin{align}
\notag\gamma \big((x,t)*(y,s) \big)& = \gamma\big((x*y, ~\alpha(x,y)(t))\big)=  \big(\phi(x)*\phi(y), ~\theta\alpha(x,y)(t)\big)\\
\notag  &= \big(\phi(x)*\phi(y), ~\alpha(\phi(x), \phi(y))\theta(t) \big)=  \big(\phi(x),~\theta(t)\big)*\big(\phi(y),~\theta(s)\big)\\
\notag  &=  \gamma (x,t)*\gamma (y,s).
\end{align}
Hence, $\gamma$ is an automorphism of $X\times_{\alpha}S$. This yields an injective map
$$\Psi: \big(\Aut(X) \times \Sigma_{S}\big)_{\alpha} \to \Aut(X \times_\alpha S)$$
given by $\Psi \big((\phi,\theta) \big)=\gamma$.   If $(\phi_1, \theta_1), (\phi_2, \theta_2) \in \big(\Aut(X) \times \Sigma_{S}\big)_{\alpha}$, then we see that
 $$\gamma_1 \gamma_2(x, t)= \gamma_1\big(\phi_2(x), ~\theta_2(t)\big)=\big(\phi_1\phi_2(x), ~\theta_1\theta_2(t)\big),$$
and hence $\Psi$ is a homomorphism.
$\blacksquare$    \end{proof}

\begin{remark}{\rm 
Let $(\phi, \theta)\in\Aut(X) \times \Sigma_{S}$ and $\gamma: X \times_\alpha S \to X \times_\alpha S$ be defined by $\gamma(x, t)= \big(\phi(x), \theta(t)\big)$. It is easy to see that $\gamma$ is an automorphism of $X\times_{\alpha}S$ if an only if $^{(\phi,\theta)}\alpha=\alpha$.}
\end{remark}

Let $A$ be an abelian group and $\psi: A \to \Sigma_{A}$ be defined by $\psi(a)=\psi_a$, where $\psi_a(b)=b+a$ for all $b \in A$. A map  $\mu : X \times X \to A$ is called a  {\it quandle  $2$-cocycle} if the map $\psi \,\mu : X \times X \to \Sigma_{A}$ is a constant quandle 2-cocycle.
Two quandle $2$-cocycles $\mu, \nu: X \times X \to A$ are called {\it cohomologous} if the corresponding constant quandle 2-cocycles $\psi\mu$ and $\psi\nu$ are cohomologous. The set of cohomology classes of quandle $2$-cocycles is the second cohomology group  $\Ho^2(X,A)$ of $X$ with coefficients in $A$.
\para 
For a given quandle $2$-cocycle $\mu$, by the formula \eqref{non-abelian-structure}, the quandle operation in $X\times_{\psi \,\mu} A$ is given by
$$(x, a) * (y, b) = \big(x * y, ~a + \mu(x, y) \big)$$
for $x, y \in X$ and $a, b \in A$. In this case, the quandle $X\times_{\psi  \,\mu} A$ is called an \index{abelian extension}{\it abelian extension} of $X$ by $A$. Many well-known quandles are abelian extensions of smaller quandles (see \cite{MR2008876} for details). As all results applicable to non-abelian extensions of quandles also hold for abelian extensions, Theorem \ref{exterrr}  leads to the following conclusion.

\begin{corollary}
Let $X$ be a quandle, $A$ an abelian group and $\mu : X \times X \to A$ a quandle $2$-cocycle. Then there is an embedding $\big(\Aut(X) \times \Aut(A) \big)_{\mu} \hookrightarrow  \Aut\big(X\times_{\psi  \,\mu} A\big)$.
\end{corollary}
\bigskip
\bigskip
 

\chapter{Residual finiteness of quandles}\label{chap residual finiteness}

\begin{quote}
This chapter investigates the property of residual finiteness within the framework of quandles. We begin by establishing foundational results concerning residually finite quandles and then focus on specific classes, including core, conjugation, and Alexander quandles, arising from residually finite groups. We further demonstrate that both free quandles and free products of residually finite quandles, under suitable conditions, preserve residual finiteness.  In the context of knot theory, we leverage these findings to establish that all link quandles exhibit residual finiteness. This conclusion, in turn, implies the solvability of the word problem for link quandles.
\end{quote}
\bigskip

\section{Properties of residually finite quandles}\label{basic properties of residually finite quandles}

Recall that, a group $G$ is said to be \index{residually finite quandle}{\it residually finite} if for each $g, h \in G$ with $g \neq h$, there exists a finite group $F$ and a homomorphism $\phi: G \to F$ such that $\phi(g) \neq \phi(h)$. 
\para
For example, finite groups, free groups and link groups \cite{MR0648524} are known to be residually finite. Residual finiteness of groups play a crucial role in combinatorial group theory. This notion can be analogously defined for quandles.

\begin{definition}
	A quandle $X$ is said to be {\it residually finite} if for all $x,y \in X$ with $x \neq y$, there exists a finite quandle $F$ and quandle homomorphism $\phi:X \rightarrow F$ such that $\phi(x) \neq \phi(y)$.
\end{definition}

In \cite{Malcev-1}, Mal$'$cev gave the definition of a residually finite algebra, and proved that for some algebras, residual finiteness implies that the word problem is solvable. The preceding definition is a particular case of Mal$'$cev's idea.
\para

Clearly, every finite quandle is residually finite, and every subquandle of a residually finite quandle is residually finite. We begin with some elementary observations.

\begin{proposition}\label{trivial-res-finite}
	Every trivial quandle is residually finite.
\end{proposition}

\begin{proof}
	Let $X$ be a trivial quandle. If $X$ is finite, then there is nothing to prove. Suppose that $X$ is infinite. Let $x, y \in X$ with $x \neq y$. Consider the trivial subquandle $\{x, y \}$ of $X$ and define $\phi: X \to \{x, y \}$ by $\phi(x)=x$ and $\phi(z)=y$ for all $z \neq x$. It is easy to see that $\phi$ is a quandle homomorphism with $\phi(x) \neq \phi(y)$, and hence $X$ is residually finite.
$\blacksquare$    \end{proof}          

Next, we investigate some closure properties of residually finite quandles. Let $\lbrace X_i\rbrace_{i \in I}$ be a family of quandles and $X= \prod_{i \in I}X_i$. Then $X$ is itself a quandle, called the {\it product quandle},  with binary operation given by  $$(x_i) \ast (y_i)= (x_i \ast y_i)$$ for $(x_i), (y_i) \in X$.  Further, for each $j \in I$, the projection map $$\pi_j: X \rightarrow X_j$$ given by $\pi_j \big((x_i)\big)=x_j$ is a quandle homomorphism.

\begin{proposition}\label{direct-prod-rf}
	Let $\{X_i\}_{i \in I}$ be a family of residually finite quandles. Then the product quandle $X = \prod_{i \in I}X_i$ is residually finite.
\end{proposition}

\begin{proof}
	Let $x=(x_i),  y=(y_i) \in X$ such that $x \neq y$. Then there exists an $i_0 \in I$ such that $x_{i_0} \neq y_{i_0} $. Since $X_{i_0}$ is residually finite, there exists a finite quandle $F$ and a homomorphism $ \phi: X_{i_0} \rightarrow F$ such that $ \phi(x_{i_0}) \neq \phi(y_{i_0})$. The homomorphism $ \phi ' := \phi ~ \pi_{i_0}$ satisfy $ \phi '(x) \neq \phi ' (y)$, and hence $X$ is a residually finite quandle.
$\blacksquare$    \end{proof}

\begin{proposition}
	The following statements are equivalent for a quandle $X$:
	\begin{enumerate}
		\item $X$ is residually finite.
		\item There exists a family $\{W_i\}_{i \in I}$ of finite quandles such that the quandle $X$ is isomorphic to a subquandle of the product quandle $\prod _{i \in I} W_i $.
	\end{enumerate}
\end{proposition}

\begin{proof}
	The implication $(2) \Rightarrow (1)$ follows from Proposition  \ref{direct-prod-rf} and the fact that a  subquandle of a residually finite quandle is residually finite. Conversely, suppose that $X$ is residually finite. For each pair $(x,y) \in X \times X$ such that $x \neq y$, there exists a finite quandle $W_{(x,y)}$ and a homomorphism $ \phi_{(x,y)} : X \rightarrow W_{(x,y)}$ such that $ \phi_{(x,y)}(x) \neq \phi_{(x,y)}(y)$. Consider the product quandle $$ W= \displaystyle\prod_{(x,y)\in X \times X,~ x \neq y} W_{(x,y)}.$$ Then the homomorphism $\psi : X \rightarrow W $ defined by $$ \psi= \displaystyle\prod_{(x,y)\in X \times X,~ x \neq y} \phi_{(x,y)}$$ is clearly injective, which proves the implication  $(1) \Rightarrow (2)$.
$\blacksquare$    \end{proof}          

\begin{definition}
A quandle $X$ is called \index{Hopfian quandle}{\it Hopfian} if every surjective quandle endomorphism of $X$ is injective.
\end{definition}

It is well-known that finitely generated residually finite groups are Hopfian  \cite{MR0003420}. We prove a similar result for quandles.
\begin{theorem}\label{hopfian-thm}
Every finitely generated residually finite quandle is Hopfian.
\end{theorem}

\begin{proof}
Let $X$ be a finitely generated residually finite quandle and $\phi:X \rightarrow X$ a surjective quandle homomorphism. Suppose that $\phi$ is not injective. Let $x_1, x_2 \in X$ such that $ x_1 \neq x_2$ and $\phi(x_1) = \phi(x_2)$. Since $X$ is residually finite, there exist a finite quandle $F$ and a quandle homomorphism $\tau:X \rightarrow F$ such that $\tau(x_1) \neq \tau(x_2)$.
\para

We claim that the maps $\tau \, \phi ^{n}:X \rightarrow F$ are distinct quandle homomorphisms for all $n \geq 0$. Let $0 \leq m <n$ be fixed integers. Since $$ \phi ^{m}:X \rightarrow X$$ is surjective, there exist $y_1, y_2 \in X$ such that $\phi ^{m}(y_1)=x_1$ and $\phi ^{m}(y_2)=x_2$. Thus, we have $$\tau \, \phi ^{m}(y_1) \neq \tau \, \phi ^{m}(y_2),$$ whereas $$\tau \, \phi ^{n}(y_1) = \tau \, \phi ^{n}(y_2),$$ which proves our claim. Thus, there are infinitely many quandle homomorphisms from $X$ to $F$. But, this is a contradiction, since $X$ is finitely generated and $F$ is finite. Hence, $\phi$ is an automorphism, and  $X$ is Hopfian.
$\blacksquare$    \end{proof}          

Analogous to groups, we define the \index{word problem} {\it word problem} for quandles as the problem of determining whether two given elements of a quandle are the same. The word problem is solvable for finitely presented residually finite groups  \cite[p. 55]{MR1357169}. Below is an analogous result for quandles.

\begin{theorem}\label{word-problem}
Every finitely presented residually finite quandle has a solvable word problem.
\end{theorem}

\begin{proof}
Let $X=\langle  S \, \mid \, R \rangle $ be a  finitely presented residually finite quandle, and  $w_1,w_2$ two words in the generators $S$. We describe two procedures which tell us whether or not $w_1 =w_2$ in $X$. The first procedure lists all the words  that we obtain by using  the relations of $X$ on the word $w_1$. If the word $w_2$ turns up at some stage, then $w_1=w_2$, and we are done.
\para
The second procedure lists all the finite quandles. Since $X$ is finitely generated, for each finite quandle $F$, the set $\Hom_{\mathcal{Q}}(X,F)$ of all quandle homomorphisms is finite. Now, for each homomorphism $\phi \in \Hom_{\mathcal{Q}}(X,F)$, we check whether or not $\phi(w_1)=\phi(w_2)$. Since $X$ is  residually finite, the above procedure must stop at some point. That is, there exists a finite quandle $F$ and  $\phi \in \Hom_{\mathcal{Q}}(X,F)$  such that  $\phi(w_1) \neq \phi(w_2)$ in $F$, and hence $w_1 \neq w_2$ in $X$.
$\blacksquare$    \end{proof}

\begin{remark}{\rm 
In \cite{MR3440592}, Belk and McGrail showed that the word problem for quandles is unsolvable in general by giving an example of a finitely presented quandle with unsolvable word problem. In view of Theorem \ref{word-problem}, such a quandle cannot be residually finite. Clearly, if a group $G$ has unsolvable word problem, then $\Conj(G)$ also has unsolvable word problem. But, if a group $G$ is finitely presented, then it is not clear whether $\Conj(G)$ is finitely presented.}
\end{remark}

\bigskip
\bigskip 


\section{Residual finiteness of quandles arising from groups}\label{sec-rf-other-quandles}

In this section, we examine the residual finiteness of conjugation, core, and Alexander quandles associated with residually finite groups. Additionally, we explore the residual finiteness of certain automorphism groups of residually finite quandles.

\begin{proposition}\label{conj-g-res-finite}
Let $G$ be a residually finite group. Then the following assertions hold:
\begin{enumerate}	
\item If $\alpha :G \rightarrow G$ is an inner automorphism, then $\Alex(G, \alpha)$ is a  residually finite quandle.
\item	$\Conj_n(G)$ and $\Core(G)$ are residually finite quandles for each $n$.
\end{enumerate}
\end{proposition}

\begin{proof}
For assertion (1), let $\alpha$ be the inner automorphism induced by $x \in G$. If $g_1, g_2 \in G$ such that $g_1 \neq g_2$, then there exists a finite group $F$ and a group homomorphism $\psi:G \rightarrow F$ such that $\psi(g_1) \neq \psi(g_2)$. Let  $\beta$ be the inner automorphism of $F$ induced by $\psi(x)$. It follows that $\psi$ viewed as a map $ \Alex(G, \alpha) \rightarrow \Alex(F, \beta)$  is a quandle homomorphism with $\psi(g_1) \neq \psi(g_2)$, and hence $\Alex(G, \alpha)$ is residually finite.
\para
Similarly, assertion (2) follows  from the fact that $\Conj_n$ and $\Core$ are functors from the category of groups to that of quandles.
$\blacksquare$    
\end{proof}          

It is well-known that the automorphism group of a finitely generated residually finite group is residually finite \cite[p. 414]{MR0207802}. An analogous result holds for inner automorphism groups of residually finite quandles \cite[Theorem 4.3]{MR3981139}.

\begin{theorem} If $X$ is a residually finite quandle, then $\Inn(X)$ is a residually finite group.
\end{theorem}

\begin{proof}
Let $S_{a_1}^{\epsilon_1} \cdots S_{a_m}^{\epsilon_m}\neq \id_X$ be an element of $\Inn(X)$, where $a_j \in X$ and $\epsilon_j \in \{1, -1\}$ for each $j$. Then there exists an element $x \in X$ such that $$S_{a_1}^{\epsilon_1}  \cdots S_{a_m}^{\epsilon_m}(x) \neq x$$ or equivalently
$$\big(\big((x \ast^{\epsilon_m} a_m) \ast^{\epsilon_{m-1}} a_{m-1} \big) \cdots \big) \ast ^{\epsilon_1} a_1 \neq x.$$ Since $X$ is residually finite, there exists a finite quandle $F$ and a surjective quandle homomorphism $\phi: X \rightarrow F$ such that
\begin{equation}\label{phi-big}
\phi\big(\big(\big((x \ast^{\epsilon_m} a_m) \ast^{\epsilon_{m-1}} a_{m-1} \big) \cdots \big) \ast ^{\epsilon_1} a_1\big) \neq \phi(x).
\end{equation}
Define a map $$\widetilde \phi: \big\{ S_{y}^{\pm 1} \,\mid\,y \in X \big\} \rightarrow \Inn(F)$$ by setting $$\widetilde \phi\big(S_{y}^{\pm 1}\big)=S_{\phi(y)}^{\pm 1}$$ 
for each $y \in X$. It is easy to see that if $S_{x_1}^{\mu_1}  \cdots S_{x_n}^{\mu_n}=\id_X$ is a relation in $\Inn(X)$, then $S_{\phi(x_1)}^{\mu_1}  \cdots S_{\phi(x_n)}^{\mu_n}=\id_F$, and hence $\widetilde \phi$ extends to a group homomorphism $\widetilde \phi: \Inn(X) \to \Inn(F)$. If $\widetilde \phi \big(S_{a_1}^{\epsilon_1} \cdots S_{a_m}^{\epsilon_m} \big) = 1$, then evaluating both the sides at $\phi(x)$ contradicts \eqref{phi-big}. Hence, $\Inn(X)$ is a residually finite group.
$\blacksquare$    \end{proof}

Next, we consider residual finiteness of automorphism groups of core and conjugation quandles.

\begin{proposition}
If $G$ is a finitely generated abelian group without $2$-torsion, then $\Aut\big(\Core(G)\big)$ is a residually finite group.
\end{proposition}

\begin{proof}
Since $G$ is a finitely generated abelian group, it is residually finite, and hence $\Aut(G)$ is also residually finite. Moreover, a semi-direct product of a finitely generated residually finite group by a residually finite group is residually finite. By Theorem \ref{main-theorem}, we have $\Aut\big(\Core(G)\big) \cong G \rtimes \Aut(G)$, and hence it is residually finite.
$\blacksquare$    \end{proof}

\begin{proposition}
If $G$ is a finitely generated residually finite group with trivial center, then $\Aut \big(\Conj(G)\big)$ is a residually finite group.
\end{proposition}

\begin{proof}
Since $G$ has trivial center, by Corollary \ref{aut conj trivial center},  we have $\Aut \big(\Conj(G)\big)=\Aut(G)$. Since $G$ is  residually finite, so it  $\Aut(G)$. $\blacksquare$    \end{proof}          
\bigskip 
\bigskip 

\section{Residual finiteness of free products of quandles}\label{residual finiteness of free quandles and free product}

Recall the construction of free quandles described in Section \ref{section free quandles} of Chapter \ref{chapter preliminaries on racks and quandles}.

\begin{theorem}\label{free-quandle-rf}
Every free quandle is residually finite.
\end{theorem}

\begin{proof}
	Let $FQ(S)$ be the free quandle on the set $S$. It is well-known that the free group $F(S)$ is residually finite \cite[Theorem 2.3.1]{MR2683112}. By Proposition \ref{conj-g-res-finite}, the quandle $\Conj\big(F(S)\big)$ is residually finite. Since $FQ(S)$ is a subquandle of $\Conj \big(F(S)\big)$, it follows that $FQ(S)$ is residually finite.
$\blacksquare$    \end{proof}

Theorems \ref{free-quandle-rf} and \ref{hopfian-thm} lead to the following result.

\begin{corollary}\label{fg-free-quandle-hopf}
Every finitely generated free quandle is Hopfian.
\end{corollary}

\begin{remark}{\rm 
Corollary \ref{fg-free-quandle-hopf} is not true for infinitely generated free quandles. Indeed, if $FQ_{\infty}$ is the free quandle freely generated by a countably infinite set $\{x_1, x_2, \ldots \}$, then we can define a homomorphism $\varphi : FQ_{\infty} \to FQ_{\infty}$ by setting
$$
\varphi(x_1) = x_1~\textrm{and}~\varphi(x_i) = x_{i-1}
$$
for $i\ge 2$. It is easy to see that $\varphi$ is an epimorphism which is not an automorphism since $\varphi(x_1) = \varphi(x_2)$.}
\end{remark}

The following is a well-known result for free groups \cite[p. 42]{MR0080089}.

\begin{theorem}\label{free-group-symm}
Let $F(S)$ be a free group on a set $S$ and $g\in F(S)$ a non-trivial element. Then there is a homomorphism $\rho: F(S) \to \Sigma_n$ for some $n$ such that $\rho(g) \neq 1$, where $\Sigma_n$ is the symmetric group on $n$ elements.
\end{theorem}
 
We apply the preceding theorem to establish an analogous result for free quandles. 

\begin{theorem}
Let $FQ(S)$ be the free quandle on a set $S$ and let $x,y \in FQ(S)$ be two distinct elements. Then there is a quandle homomorphism $\phi:FQ(S) \rightarrow \Conj(\Sigma_n)$ for some $n$ such that $\phi(x) \neq \phi(y)$.
\end{theorem}

\begin{proof}
Recall that the map $\Phi: FQ(S) \rightarrow \Conj \big(F(S)\big) $ in Theorem \ref{equivalence-two-models} is an injective quandle homomorphism. Let $(a_{1},{w_{1}}) \neq (a_{2},{w_{2}}) \in FQ(S)$. Then $g_{1} \neq g_{2} \in F(S)$, where $g_{1}=\Phi \big((a_{1},{w_{1}}) \big)$ and $g_{2}=\Phi \big((a_{2},{w_{2}})\big)$. Consequently, $g_{2}^{-1}g_{1}$ is a non-trivial element of $F(S)$. By Theorem \ref{free-group-symm}, there exists a group homomorphism $\rho:F(S) \rightarrow \Sigma_n$ for some $n$ such that $\rho(g_{1})\neq\rho(g_{2})$. We can view $\rho$ as a quandle homomorphism
$\Conj \big(F(S)\big) \rightarrow \Conj(\Sigma_n)$. Taking  $\phi := \rho\, \Phi : FQ(S) \rightarrow \Conj(\Sigma_n)$, we see that $\phi \big((a_{1},{w_{1}})\big) \neq \phi \big((a_{2},{w_{2}}) \big)$, and the proof is complete.
$\blacksquare$    \end{proof}

The following result is well-known in combinatorial group theory, first proved by Gruenberg \cite[Theorem 4.1]{MR0087652}. See also \cite{MR0484178, MR0463296}.

\begin{theorem}\label{free-product-of-residually-finite-groups}
A free product of residually finite groups is residually finite.
\end{theorem}

We prove an analogue of the preceding theorem for quandles provided their adjoint groups are residually finite. Recall that, every element of a quandle $X$ can be written in a left associated form $x_0 \ast^{\epsilon_1} x_1 \ast^{\epsilon_2} x_2 \ast ^{\epsilon_3} \cdots \ast ^{\epsilon_n} x_n$. Moreover, the expression $x_0 \ast^{\epsilon_1} x_1 \ast^{\epsilon_2} x_2 \ast ^{\epsilon_3} \cdots \ast ^{\epsilon_n} x_n$ is called a \index{reduced form}{\it reduced form} when $x_0 \neq x_1$ and if $x_i = x_{i+1}$, then $\epsilon_i = \epsilon_{i+1}.$ Note that the reduced form is not unique. For example, if $X = \{ t \} \star \R_3$ is the free product of the one element trivial quandle and the dihedral quandle $\R_3 = \{ 0, 1, 2 \}$, then
$$
(t * 1) * 2 = (t * 2) * (1 * 2) =  (t * 2) * 0.
$$

\begin{theorem} \cite[Theorem 4.3]{MR4075375}\label{finite-free-product-of-residually-finite-quandles}
	Let $X_1, \ldots,X_n$ be residually finite quandles. If each adjoint group $\Adj(X_i)$ is residually finite, then $X_1 \star \cdots \star X_n$ is a residually finite quandle.
\end{theorem}
\begin{proof}
It is enough to consider the case $n=2$. Set $X= X_1 \star X_2$. Let $x$ and $x'$ be two distinct elements of $X$, where
$$ x =a_0 \ast ^{\epsilon_1}a_1  \ast ^{\epsilon_2} a_2 \ast ^{\epsilon_3}\cdots \ast ^{\epsilon_n} a_n \quad \textrm{and} \quad 
x'=b_0 \ast ^{\epsilon_1'}b_1 \ast ^{\epsilon_2'} b_2  \ast ^{\epsilon_3'}\cdots \ast ^{\epsilon_m'} b_m $$
are their reduced expressions with $a_i, b_j  \in X_1 \sqcup X_2$. 
\begin{enumerate}
\item[Case 1:] $x, x' \in X_1$ or  $x, x' \in X_2$. Suppose that  $x, x' \in X_1$. Since $X_1$ is a residually finite quandle, there exist a finite quandle $F$ and a quandle homomorphism $\phi: X_1 \rightarrow F$ such that $\phi(x)\neq \phi(x')$. Let $a \in F$ be  some fixed element. Define a map $\tilde{\phi}: X \rightarrow {F}$ by setting
	$$
	\tilde{\phi}(q)= \left\{
	\begin{array}{ll}
	\phi(q) &  \textrm{if}~q  \in X_1,\\
	a~& \textrm{if}~ q \in X_2.
	\end{array} \right.
	$$
	Since $\tilde{\phi}$ preserve all the relations in $X$, it extends to a quandle homomorphism with $\tilde{\phi}(x) \neq \tilde{\phi}(x')$ in $F$.
\para
\item[Case 2:] $ x \in X_1$ and $x' \in X_2$. Let $Y= \{ a, b \}$ and $FQ(Y)$ the free quandle on $Y$. Define a map $\phi : X \rightarrow FQ(Y)$ by setting
	$$
	\phi(q)= \left\{
	\begin{array}{ll}
	a~&  \textrm{if}~q  \in X_1,\\
	b~& \textrm{if}~ q \in X_2.
	\end{array} \right.
	$$
	Since $\phi$ preserve all the relations in $X$, it extends to a quandle homomorphism with $\phi(x) \neq \phi(x')$ in $FQ(Y)$.
\para
\item[Case 3:] $ x \in X \setminus (X_1 \sqcup X_2)$ and $x ' \in X_1$. We can assume that either $a_0 \in X_1$, $a_1 \in X_2$ and $a_2, \ldots , a_n \in X_1 \sqcup X_2$ or $a_0 \in X_2$, $a_1 \in X_1$ and $a_2,\ldots , a_n \in X_1 \sqcup X_2$, that is,

	$$
	x= \left\{
	\begin{array}{ll}
	q_1 \ast^{\epsilon_1} q_2  \ast^{\epsilon_2} a_2 \ast ^{\epsilon_3} \cdots \ast^{\epsilon_n} a_n & ~~\textrm{where}~~q_1 \in X_1,~~q_2 \in X_2,\\
	\textrm{or} &\\
	q_2 \ast ^{\epsilon_1} q_1 \ast ^{\epsilon_2} a_2 \ast ^{\epsilon_3} \cdots \ast^{\epsilon_n} a_n & ~~\textrm{where}~~q_1 \in X_1,~~q_2 \in X_2.
	\end{array} \right.
	$$
	It follows from Proposition \ref{presentation-of-adjoint-group-of-free-product-of-quandles}, Theorem \ref{free-product-of-residually-finite-groups} and Proposition \ref{conj-g-res-finite} that $\Conj\big(\Adj(X)\big)$ is a residually finite quandle. Let	$i_X:X  \rightarrow \Conj\big(\Adj(X)\big)$	be the natural quandle homomorphism. Then, we have
	$$
	i_X(x_0 \ast ^{\epsilon_1}x_1  \ast ^{\epsilon_2} x_2 \ast ^{\epsilon_3}\cdots \ast ^{\epsilon_n} x_n )= (\textswab{a}_{x_1} ^{\epsilon_1}  \textswab{a}_{x_2} ^{\epsilon_2} \cdots  \textswab{a}_{x_n} ^{\epsilon_n})^{-1} \textswab{a}_{x_0} (\textswab{a}_{x_1} ^{\epsilon_1}  \textswab{a}_{x_2} ^{\epsilon_2} \cdots  \textswab{a}_{x_n} ^{\epsilon_n}).
	$$
	We claim that $i_X(x) \neq i_X(x')$. 
\para
\item[Case 3.1:] If $x =  q_1 \ast^{\epsilon_1} q_2  \ast^{\epsilon_2} a_2 \ast ^{\epsilon_3} \cdots \ast^{\epsilon_n} a_n, ~~\textrm{where}~~q_1 \in X_1,~~q_2\in X_2,~~ \textrm{and}~~a_2, \ldots, a_{n}   \in X_1 \sqcup X_2$, then
	\begin{align*}
	i_X(x) &= (\textswab{a}_{q_2} ^{\epsilon_1}  \textswab{a}_{a_2} ^{\epsilon_2} \cdots  \textswab{a}_{a_n} ^{\epsilon_n})^{-1} \textswab{a}_{q_1} (\textswab{a}_{q_2} ^{\epsilon_1}  \textswab{a}_{a_2} ^{\epsilon_2} \cdots  \textswab{a}_{a_n} ^{\epsilon_n})\\
	&= \textswab{a}_{a_n} ^{-\epsilon_n} \ldots \textswab{a}_{a_2} ^{-\epsilon_2} \textswab{a}_{q_2} ^{-\epsilon_1} \textswab{a}_{q_1} \textswab{a}_{q_2} ^{\epsilon_1}  \textswab{a}_{a_2} ^{\epsilon_2} \cdots  \textswab{a}_{a_n} ^{\epsilon_n}.
	\end{align*}
	Suppose that $i_X(x)=i_X(x')$. Then by the Remark \ref{image-under-eta} and the fact that elements of $\Adj(X_1)$ have no relations with elements of $\Adj(X_2)$ in the group $\Adj(X)$, it follows that either $\textswab{a}_{q_2} ^{\epsilon_1}  \textswab{a}_{a_2}^{\epsilon_2} \cdots  \textswab{a}_{a_n} ^{\epsilon_n}= 1$ in $\Adj(X)$ or $\textswab{a}_{q_2} ^{\epsilon_1}  \textswab{a}_{a_2} ^{\epsilon_2} \cdots  \textswab{a}_{a_n} ^{\epsilon_n}= \textswab{a}_{q_{i_1}}^{\epsilon _1} \textswab{a}_{q_{i_2}}^{\epsilon _2} \cdots \textswab{a}_{q_{i_k}}^{\epsilon _k} $, where $q_{i_1}, q_{i_2} , \ldots, q_{i_k}  \in X_1$ and $\epsilon_j = \pm 1$ for $1 \leq j \leq k$. Since $\Adj(X)$ has a left-action on the quandle $X$, this implies that in either situation $q_1 \cdot (\textswab{a}_{q_2} ^{\epsilon_1}  \textswab{a}_{a_2} ^{\epsilon_2} \ldots  \textswab{a}_{a_n} ^{\epsilon_n})$ belongs to $X_1$. Thus, $x=q_1 \ast ^{\epsilon_1} q_2 \ast ^{\epsilon_2} a_2 \ast ^{\epsilon_3}\cdots \ast ^{\epsilon_n} a_n  \in X_1$, which is a contradiction. Hence, we must have $i_X(x) \neq i_X(x')$.
\para
\item[Case 3.2:] If $x = q_2 \ast ^{\epsilon_1} q_1 \ast ^{\epsilon_2} a_2 \ast ^{\epsilon_3} \cdots \ast^{\epsilon_n} a_n,  ~~\textrm{where}~~q_1 \in X_1,~~q_2 \in X_2~~ \textrm{and}~~a_2, \ldots , a_{n} \in X_1 \sqcup X_2,$ then
	\begin{align*}
	i_X(x) &= (\textswab{a}_{q_1} ^{\epsilon_1}  \textswab{a}_{a_2} ^{\epsilon_2} \cdots  \textswab{a}_{a_n} ^{\epsilon_n})^{-1} \textswab{a}_{q_2} (\textswab{a}_{q_1} ^{\epsilon_1}  \textswab{a}_{a_2} ^{\epsilon_2} \cdots  \textswab{a}_{a_n} ^{\epsilon_n})\\
	&= \textswab{a}_{a_n} ^{-\epsilon_n} \ldots \textswab{a}_{a_2} ^{-\epsilon_2} \textswab{a}_{q_1} ^{-\epsilon_1} \textswab{a}_{q_2} \textswab{a}_{q_1} ^{\epsilon_1}  \textswab{a}_{a_2} ^{\epsilon_2} \cdots  \textswab{a}_{a_n} ^{\epsilon_n}.
	\end{align*}
	Clearly $i_X(x) \neq i_X(x')$, since they belong to different conjugacy classes in $\Adj(X)$.
\para
\item[Case 4:] $ x \in X \setminus (X_1 \sqcup X_2)$ and $x ' \in X_2$. This is similar to Case 3.
\para
\item[Case 5:] $ x, x' \in X \setminus (X_1 \sqcup X_2).$ This case can be reduced to one of the Cases (1--4) by repeated use of the second quandle axiom. More precisely, we can replace the element $x$ by $y$ and $x'$ by $y'$, where
	\begin{align*}
	y&= a_0 \ast ^{\epsilon_1} a_1 \ast ^{\epsilon_2} \cdots \ast^{\epsilon_n} a_n \ast^{-\epsilon_m'} b_m \ast^{-\epsilon' _{m-1}} b_{m-1} \ast^{-\epsilon' _{m-2}} \cdots \ast^{-\epsilon' _{1}} b_1,\\
	y'&= b_0.
	\end{align*}
\end{enumerate}
Hence, we conclude that $X= X_1 \star X_2$ is a residually finite quandle.
$\blacksquare$    \end{proof}          

We now extend  the preceding result to arbitrary family of quandles \cite[Theorem 4.4]{MR4075375}.

\begin{theorem}\label{arbitrary-free-product-of-residually-finite-quandles}
Let $\{X_i\}_{i \in I}$ be a family of residually finite quandles. If each $\Adj(X_i)$ is a residually finite group, then the free product $\star_{i \in I} X_i$ is a residually finite quandle.
\end{theorem}
\begin{proof}
Let $X= \star_{i \in I} X_i$ be the free product of residually finite quandles $X_i$. Let $x, x' \in X$ be two distinct elements such that
$$ x =a_0 \ast ^{\epsilon_1}a_1  \ast ^{\epsilon_2} a_2 \ast ^{\epsilon_3}\cdots \ast ^{\epsilon_n} a_n \quad \textrm{and} \quad 
x'=b_0 \ast ^{\varepsilon_1}b_1 \ast ^{\varepsilon_2} b_2  \ast ^{\varepsilon_3}\cdots \ast ^{\varepsilon_m} b_m.
$$
Consider the set $S = \{ a_i, b_j  \, \mid \,  1 \leq i \leq n~~\textrm{and}~~1 \leq j \leq m \}$. Then $S$ is a finite set contained in $ X_{i_1} \sqcup \cdots \sqcup X_{i_k}$ for some $i_1, \ldots, i_k \in I$. Let $a$ be a fixed element in $\sqcup_{i \in I} X_i \setminus ( X_{i_1} \sqcup \cdots \sqcup X_{i_k})$.
Define a map $$\phi: X \rightarrow X_{i_1} \star \cdots \star X_{i_k}$$ by setting
	
$$
\phi(q)= \left\{
\begin{array}{ll}
q &  \textrm{if}~~q  \in  X_{i_1} \sqcup \cdots \sqcup X_{i_k},\\
a~& \textrm{if}~ q \in \sqcup_{i \in I} X_i \setminus (X_{i_1} \sqcup \cdots \sqcup X_{i_k}).
\end{array} \right.
$$
Since $\phi$ preserves all the relations in $X$, it extends to a quandle homomorphism with $\phi(x) \neq \phi(x')$. Hence, it follows from Theorem \ref{finite-free-product-of-residually-finite-quandles} that $X$ is a residually finite quandle.
$\blacksquare$    \end{proof}          
 
We note that Theorem \ref{arbitrary-free-product-of-residually-finite-quandles} gives an alternate proof of Theorem \ref{free-quandle-rf}. Finally we discuss residual finiteness of adjoint groups of quandles. 

\begin{proposition}\label{ass-group-finite}
If $X$ is a finite quandle, then its adjoint group $\Adj(X)$ is a residually finite group.
\end{proposition}

\begin{proof}
Consider the natural group homomorphism  $\psi_{X}: \Adj(X) \rightarrow \Inn(X)$. Since $X$ is a finite quandle, the inner automorphism group $\Inn(X)$ is finite, and hence  $\Adj(X)/ \ker(\psi_{X})$ is finite. Moreover, $\ker(\psi_{X})$ is contained in the center $\Z \big(\Adj(X) \big)$ of $\Adj(X)$, and hence $\Adj(X)/ \Z \big(\Adj(X)\big)$ is finite. It is well-known that
if $G$ is a finitely generated group with infinitely generated center $\Z(G)$, then the quotient $G/\Z(G)$ is not finitely presented. Hence, $\Z \big(\Adj(X)\big)$ is a finitely generated abelian group, and therefore it is residually finite. It is known that if $N$ is a normal subgroup of finite index in a group $G$ and $N$ is residually finite, then $G$ is also residually finite \cite[Proposition 2.2.12]{MR2683112}. Hence,  $\Adj(X)$ is a residually finite group.
$\blacksquare$    \end{proof}          
\para

As a consequence of Proposition \ref{ass-group-finite} and Theorem \ref{arbitrary-free-product-of-residually-finite-quandles}, we obtain the following result.

\begin{corollary}
A free product of finite quandles is residually finite.
\end{corollary}
\bigskip 
\bigskip 


\section{Residual finiteness of link quandles}\label{sec-non-split}
In this section, we prove that the link quandle of any oriented link in $\mathbb{S}^3$ is residually finite. We recall the following definition from \cite{Malcev-1}.

\begin{definition}
A subgroup $H$ of a group $G$ is said to be \index{finitely separable subgroup}{\it finitely separable} if for any $g \in G \setminus H$, there exists a finite group $F$ and a group homomorphism $\phi: G \rightarrow F$ such that $\phi(g) \not\in \phi(H)$.
\end{definition}
\para

For example, if $G$ is a residually finite group and $H$ a finite subgroup of $G$, then $H$ is finitely separable in $G$. In analogy with groups, we introduce the following definition.

\begin{definition}
	A subquandle $Y$ of a quandle $X$ is said to be \index{finitely separable subquandle}{\it finitely separable} in $X$ if for each $x \in X\setminus Y$, there exists a finite quandle $F$ and a quandle homomorphism $\phi: X \to F$ such that $\phi(x) \not\in \phi(Y)$.
\end{definition}
 
The following result is of independent interest.
 
\begin{proposition}
Let $X$ be a residually finite quandle and $\alpha \in \Aut(X)$. If the fixed-point set $\Fix(\alpha)= \{x \in X \, \mid \,  \alpha(x)=x\}$ is non-empty, then it is a finitely separable subquandle of $X$.
\end{proposition}

\begin{proof}
Clearly $\Fix(\alpha)$ is a subquandle of $X$. Let $x_0 \in X \setminus \Fix(\alpha)$, that is, $\alpha(x_0) \neq x_0$. Since $X$ is residually finite, there exists a finite quandle $F$ and a quandle homomorphism $\phi: X \rightarrow F$ such that $\phi \big(\alpha(x_0) \big) \neq \phi(x_0)$. Define a map $\psi: X \rightarrow F\times F $ by $\psi(x)=\big(\phi(x), \, \phi\, \alpha(x)\big)$. Clearly $\psi$ is a quandle homomorphism with $\psi(x_0) \not\in \psi \big(\Fix(\alpha) \big)$, and hence $\Fix(\alpha)$ is finitely separable in $X$.
$\blacksquare$    \end{proof}          

The following result will be used to prove residual finiteness of link quandles \cite[Proposition 3.4]{MR4075375}.

\begin{proposition} \label{res-finite-quandle}
	Let $G$ be a group, $\{z_i \, \mid \, i \in I \}$ be a finite set of elements of $G$, and $\{H_i \, \mid \, i \in I \}$ subgroups of $G$ such that $H_i \le \C_G(z_i)$. If each $H_i$ is finitely separable in $G$, then the quandle $\bigsqcup_{i \in I}(G,H_i,z_i)$ is residually finite.
\end{proposition}

\begin{proof}
	Let $H_k a \neq	 H_j b$ be two elements of $\bigsqcup_{i \in I}(G,H_i,z_i)$.
\begin{enumerate}
\item[Case 1:] $k \neq j$. Let $F=\{a',b' \}$ be a two element trivial quandle. Define $$\phi: \bigsqcup_{i \in I}(G,H_i,z_i) \rightarrow F$$ by setting
	$$
	\phi(H_i x)= \left\{
	\begin{array}{ll}
	a' &  \textrm{if}~ i=k,\\
	b'~& \textrm{if}~ i\neq k.
	\end{array} \right.
	$$
	Then $\phi$ is a quandle homomorphism with $\phi(H_k a)\neq \phi(H_j b)$.
\para
\item[Case 2:] $k=j$. Since $H_k a \neq H_k b$, we have $a \neq hb$ for any $h \in H_k$. Further, since $H_k$ is finitely separable in $G$, there exist a finite group $F$ and a group homomorphism $\phi:G \rightarrow F$ such that $\phi(a) \neq \phi(h b)$ for each $h \in H_k$. For each $i \in I$, set $\overline{H}_i= \phi(H_i)$  and $\bar{z}_i=\phi(z_i)$. Then  $\overline{H}_i \le  \C_F( \bar{z}_i)$ and $\bigsqcup_{i \in I}(F, \overline{H}_i,\bar{z}_i)$ is a finite quandle. Further, the group homomorphism $\phi:G \rightarrow F$ induces a map
	
	 $$\bar{\phi}:\bigsqcup_{i \in I}(G,H_i,z_i) \rightarrow \bigsqcup_{i \in I}(F,\overline{H}_i,\bar{z}_i)$$ 	 
	 given by $$\bar{\phi}(H_i x)=\overline{H_i} \phi(x),$$ which is a quandle homomorphism. Also, $\bar{\phi}(H_k a) \neq \bar{\phi}(H_k b)$, otherwise $\phi(a)=\phi(hb)$ for some $h\in H_k$, which is a contradiction.  
\end{enumerate}	 
Hence, the quandle $\bigsqcup_{i \in I}(G,H_i,z_i)$ is residually finite.
$\blacksquare$    \end{proof}          

Let $M$ be a 3-manifold. A surface $F$ in $M$, which is not a 2-sphere, is said to be \index{incompressible subsurface}{\it incompressible} if the homomorphism $\pi_1(F, x) \to \pi_1(M, x)$ is injective for $x \in F$. The following result concerning  finitely separable subgroups of fundamental groups of irreducible 3-manifolds is due to Long and Niblo \cite[Theorem 1]{MR1109662}.

\begin{theorem}\label{Long-Niblo-thm}
	Suppose that $M$ is an orientable, irreducible compact 3-manifold and $X$ an incompressible connected subsurface of a component of $\partial(M)$. If $x\in X$ is a base point, then $\pi_1(X,x)$ is a finitely separable subgroup of $\pi_1(M,x)$.
\end{theorem}
\para

Joyce \cite{MR2628474, MR0638121} and Matveev \cite{MR0672410} independently introduced the link quandle $Q(L)$ of an oriented link $L$ in $\mathbb{S}^3$ (see Example \ref{construction link quandle}). We recall their topological construction here.
 \para
 
Let $L$ be an oriented link in $\mathbb{S}^3$ with component knots $K_1, \ldots, K_m$ and $V(L)=\bigsqcup_{i=1}^m V(K_i)$ be its tubular neighbourhood. Let $X(L)$ be the closure of $\mathbb{S}^3 \setminus V(L)$ in $\mathbb{S}^3$.  For each $i$, consider the inclusion maps $\kappa_i : \partial V(K_i) \rightarrow X(L)$ and induced group homomorphisms $$\kappa^*_i : \pi_1 \big(\partial V(K_i) \big) \rightarrow \pi_1\big(X(L)\big).$$ Since $\pi_1\big(\partial V(K_i)\big) \cong \mathbb{Z} \oplus \mathbb{Z}$, we may choose generators of this group corresponding to the \textit{meridian} $m_i$ and the \textit{longitude}  $\ell_i$ of the knot  $K_i$, which are well-defined up to homotopy in $\partial V(K_i)$.  The group $P_i=\kappa^*_i \big(\pi_1(\partial V(K_i)) \big)$ is called the $i$-th \textit{peripheral subgroup} of $\pi_1\big(X(L)\big)$. 
\para

Let $Q(L)$ be the set of homotopy classes of paths in the space $X(L)$, with initial point on $ \partial V(L)$ and a fixed endpoint $ x_0 \in X(L)$. We require that these conditions on the initial and the final points are preserved during the homotopy.  Define a binary operation on $ Q(L)$ by
\begin{eqnarray*}
	[a] * [b] = [ab^{-1}m_j b],
\end{eqnarray*}
where $a$ and $b$ are paths in $X(L)$ representing $[a]$ and $[b]$, respectively. It can be checked that the binary operation turns $Q(L)$ into a quandle, called the {\it link quandle} of the link $L$. The natural action of $\pi_1\big(X(L)\big)$ on $Q(L)$ given by $$[a] \cdot [g]=[a g]$$ yields the following homogeneous representation of $Q(L)$ in terms of $\pi_1 \big(X(L)\big)$ (see \cite[Section 4.9, Corollary 1]{MR2628474} or \cite[Theorem 3.6]{MR4075375}).

\begin{theorem}\label{coset representation of knot quandle}
Let $L$ be an oriented link in $\mathbb{S}^3$ with component knots $K_1, \ldots, K_m$. For each $i$, let $m_i$ and $\ell_i$ be the fixed meridian and the longitude of the component $K_i$, respectively. Then $Q(L) \cong \bigsqcup_{i=1}^m \big(\pi_1(X(L)),P_i,m_i \big)$, where $P_i=\langle m_i,\ell_i \rangle.$ 
\end{theorem}

Following is the main result of this section  \cite[Theorem 4.5]{MR4075375}.

\begin{theorem}\label{non-split-link-quandle-is-residually-finite}
The link quandle of any oriented link in $\mathbb{S}^3$ is residually finite.
\end{theorem}

\begin{proof} 
If $L$ is a trivial link, then its link quandle is a free quandle. By Theorem \ref{free-quandle-rf}, such a quandle is residually finite. Suppose that $L$ is a non-trivial and non-split link. By Theorem \ref{coset representation of knot quandle}, we have 
$$Q(L) \cong \bigsqcup_{i=1}^m \big(\pi_1(X(L)),P_i,m_i \big),$$
where $P_i=\langle m_i,\ell_i \rangle.$   Since $L$ is non-split, it follows from Theorem \ref{Long-Niblo-thm} that each $P_i$ is finitely separable in $\pi_1 \big(X(L)\big)$. Hence, by Proposition \ref{res-finite-quandle}, the link quandle $Q(L)$ is residually finite.
\para
Now, suppose that $L$ is a non-trivial and split link. Then the link quandle $Q(L)$ is a free product of link quandles of its non-split components. Using the fact that all link groups are residually finite \cite{MR0895623}, the result now follows from the preceding case and Theorem \ref{finite-free-product-of-residually-finite-quandles}.
$\blacksquare$    \end{proof}          

\begin{corollary}
The link quandle of any oriented link in $\mathbb{S}^3$ is Hopfian, has solvable word problem, and has residually finite inner automorphism group.
\end{corollary}

Since adjoint groups of finite quandles, free quandles and link quandles are residually finite, the following seems to be the case in general.
 
\begin{conjecture}\label{residually-finite-adjoint-group}
The adjoint group of a finitely presented residually finite quandle is a residually finite group.
\end{conjecture}
 
\begin{remark}
{\rm 
If Conjecture \ref{residually-finite-adjoint-group} is true, then by Theorem \ref{arbitrary-free-product-of-residually-finite-quandles}, a free product of finitely presented residually finite quandles would be residually finite, which would be an analogue of Theorem \ref{free-product-of-residually-finite-groups} for quandles.
}
\end{remark}

\begin{remark}
{\rm 
In \cite{DhanwaniSarafSingh}, building on results from Thurston's geometrisation program, Theorem \ref{non-split-link-quandle-is-residually-finite} was extended to all $n$-quandles of link quandles.}
\end{remark}


\chapter{Orderability of quandles}\label{chap orderability of quandles}

\begin{quote}
Orderability of an algebraic structure is expected to have profound implications on its inherent structure. In this chapter, we introduce the idea of orderability of quandles. We begin by exploring properties of linear orderings on quandles and provide many examples of orderable quandles. We also give some general constructions of orderable quandles, and observe that  free quandles  are orderable. As application to knot theory, we also delve into the examination of the orderability of link quandles.
\end{quote}
\bigskip

\section{Properties and examples of orderable quandles}
Existence of a linear order on a group is known to have profound implications on its structure. For instance, a left-orderable group cannot have torsion and a bi-orderable group cannot have even generalized torsion (product of conjugates of a non-trivial element being trivial). In this chapter, we investigate orderability of quandles with applications to link quandles.
\para

\begin{definition} 
A quandle $X$ is said to be \index{left-orderable quandle}{\it left-orderable} if there is a (strict) \index{linear order}linear order $<$ on $X$ such that $x<y$ implies $z*x<z*y$ for all $x,y,z\in X$. Similarly, a quandle $X$ is \index{right-orderable quandle}{\it right-orderable} if there is a linear order $<$ on $X$ such that $x< y$ implies $x*z< y*z$ for all $x,y,z\in X$.  Further, a quandle is \index{bi-orderable quandle}{\it bi-orderable} if it has a linear order with respect to which it is both left and right ordered.
\end{definition}

If $\T=\left\{x_1,x_2,\ldots\right\}$ is a trivial quandle, then it is clear that the linear order $x_1<x_2<\cdots$ is preserved under multiplication on the right, but is not preserved under multiplication on the left. Thus, a trivial quandle with more than one element is right-orderable but not left-orderable. 
This stands in contrast to the case of groups, where left-orderability and right-orderability are equivalent.
\para

The following result demonstrates that a wide variety of left- or right-orderable quandles can be constructed from bi-orderable groups (\cite[Proposition 7]{MR2320157} and (\cite[Proposition 2.6]{MR4330281}).

\begin{proposition}\label{conj right orderable}
The following  assertions hold for a bi-orderable group $G$:
\begin{enumerate}
\item  $\Conj_n(G)$ is a right-orderable quandle for each $n$.
\item $\Core(G)$ is a left-orderable quandle. 
\item If $\phi \in \Aut(G)$ is an order reversing automorphism, then $\Alex(G, \phi)$ is a left-orderable quandle.
\end{enumerate}
\end{proposition}

\begin{proof}
Let $G$ be bi-ordered with respect to the  linear order $<$ and $x, y, z \in G$ such that $x < y$. Then $$x*z=z^nxz^{-n} < z^nyz^{-n}=y*z$$ implies that $\Conj_n(G)$ is a right-orderable quandle. Similarly, 
$$z*x= xz^{-1}x <  yz^{-1}x <  yz^{-1}y= z*y$$ implies that $\Core(G)$ is a left-orderable quandle. This establishes (1) and (2).
\para
For assertion (3), bi-ordering of $G$ and $\phi$ being order reversing implies that $\phi(x)^{-1} < \phi(y)^{-1}$. This gives
$$z*x= \phi(zx^{-1})x=  \phi(z)\phi(x^{-1})x < \phi(z)\phi(x^{-1})y < \phi(z)\phi(y^{-1})y=z*y,$$
which proves that $\Alex(G, \phi)$ is left-orderable.
$\blacksquare$    \end{proof}

For example, free groups \cite{MR0031482},  fundamental group of any connected surface other than the projective plane or the Klein bottle \cite{MR2141698} and  pure braid groups  \cite{MR0975081} are bi-orderable.
\para

An immediate consequence of the preceding proposition is the following result.

\begin{corollary}
The following  assertions hold for any quandle $X$:
\begin{enumerate}
\item If $X$ is a subquandle of $\Conj_n(G)$ for some bi-orderable group $G$, then $X$ is right-orderable. 
\item If $X$ is a subquandle of $\Core(G)$ for some bi-orderable group $G$, then $X$ is left-orderable. 
\end{enumerate}
\end{corollary}

It is easy to see that a right- or left-orderable group must be infinite. But, this is not true for quandles since any finite trivial quandle is right-orderable. A quandle is said to be \index{semi-latin}{\it semi-latin} if left multiplication by each element is injective. Clearly, a latin quandle is semi-latin, but the converse does not hold in general. The following properties hold for quandles \cite[Proposition 3.7]{MR4450681}.

\begin{proposition}\label{orderable-implies-latin}
The following  assertions hold for any quandle $X$:
\begin{enumerate}
\item If $X$ is  right-orderable, then the $\langle S_y \rangle $-orbit of  $x$ is infinite for all $x, y \in X$ with $S_y(x) \ne x$.
\item If $X$ is  left-orderable, then it is semi-latin and the set $\{L_y^n(x)\}_{n \ge 0}$ is infinite for $x \ne y \in X$, where
$$
L_y^0(x) = x~\textrm{and}~L_y^{i}(x) =L_y \big(L_y^{i-1}(x) \big) ~\textrm{for}~i \ge 1.
$$
\end{enumerate}
\end{proposition}

\begin{proof}
If $x<S_y(x)$ and the $\langle S_y \rangle$-orbit of $x$ is finite, then right-orderability of $X$ implies that  $$x <  S_y(x) <  S_y^2(x) <  \cdots <  S_y^n(x)=x$$ for some integer $n$, which is a contradiction. Similarly, the assertion follows if $S_y(x) < x$.
\para

Suppose that there are elements $x, y, z \in X$ with $y\neq z$, say $y < z$, such that $x * y = x * z$. This is a contradiction to left-orderability of $X$, and hence $X$ must be semi-latin. Further, if $x \ne y$ are two elements of $X$ such that $x = L_y(x)$, then $x * x = y * x$, which contradicts the second quandle axiom. Hence, $x<L_y(x)$ or $L_y(x) < x$. Suppose that $x<L_y(x)$. Since $X$ is left-orderable, we have 
$$
x < L_y(x) < L_y^2(x) < \cdots < L_y^n(x)< \cdots,
$$
and hence $\{L_y^n(x)\}_{n \ge 0}$ is infinite. The case $L_y(x) < x$ is similar.
$\blacksquare$    \end{proof}          

The other sided orderability of quandles in Proposition \ref{conj right orderable} does not hold in general, as shown in \cite[Corollaries 3.8 and 3.9]{MR4450681}.

\begin{proposition}
The following  assertions hold for any non-trivial group $G$:
\begin{enumerate}
\item The quandle $\Conj_n(G)$ is not left-orderable for each $n$.
\item The quandle $\Core(G)$ is not right-orderable.
\item If $\phi\in \Aut(G)$ an involution, then the quandle $\Alex(G,\phi)$ is not right-orderable.
\end{enumerate}
\end{proposition}

\begin{proof}
Since $\Conj_n(G)$ is not semi-latin and $\Core(G)$ is involutory, the assertions follow from Proposition \ref{orderable-implies-latin}. If $\phi \in \Aut(G)$ is an involution, then $\Alex(G, \phi)$ is involutory and we obtain the third assertion.
$\blacksquare$    \end{proof}

\begin{proposition}\label{orderable quandle infinite}
Any non-trivial left or right-orderable quandle is infinite.
\end{proposition}

\begin{proof}
Let $X$ be a non-trivial quandle that is right-orderable. Then there exist elements $x \neq y$ in $X$ such that $S_y(x)\neq x$. It follows from Proposition \ref{orderable-implies-latin}(1) that the $\left\langle S_y\right\rangle$-orbit of $x$ is infinite, and hence $X$ must be infinite. On the other hand, if $X$ is left-orderable, then by Proposition \ref{orderable-implies-latin}(2), the set $\{L_y^n(x)\}_{n \ge 1}$ is infinite for any $x\neq y$ in $X$, and hence $X$ must be infinite.
$\blacksquare$    \end{proof}          

For free quandles, we have the following result \cite[Theorem 3.5]{MR4450681}.

\begin{theorem}\label{orderability-free-quandle}
Free quandles are right-orderable and semi-latin. In particular, link quandles of trivial links are right orderable.
\end{theorem}

\begin{proof}
Let $S$ be a non-empty set. It is well-known that the free group $F(X)$ on the set $S$ is bi-orderable  (see \cite[Theorem 2.1.9]{Deroin2014} and \cite{MR0031482}). It follows from  Proposition \ref{conj right orderable}(1) that $\Conj \big(F(S)\big)$ is right-orderable. By Proposition \ref{equivalence-two-models}, the free quandle $FQ(S)$ is a subquandle of $\Conj \big(F(S)\big)$, and hence  is right-orderable.
\para

Next, we prove that $FQ(S)$ is semi-latin. If $S$ is a singleton set, then $FQ(S)$ is the one-element trivial quandle and the assertion is evident. Suppose that $|S| > 1$ and there are elements $x, y, z \in FQ(S)$ such that $x \not= y$ and $z * x = z * y$. Since the free quandle $FQ(S)$ is a subquandle of $\Conj \big(F(S)\big)$ consisting of all conjugates of elements of $S$, we can write
$$
z = s_i^{z_0},~~x = s_j^{x_0}~\textrm{and}~y = s_k^{y_0}
$$
for some $s_i, s_j, s_k \in S$ and $z_0, x_0, y_0 \in F(S)$. Here, $a^b$ denotes the element $b^{-1} a b$. The identity $z * x = z * y$ gives
$$
s_i^{z_0 x_0^{-1} s_j x_0} = s_i^{z_0 y_0^{-1} s_k y_0}
$$
in $F(S)$, which is further equivalent to
$$
s_i^{z_0 x_0^{-1} s_j x_0 y_0^{-1} s_k^{-1}  y_0 z_0^{-1}} = s_i.
$$
Such an equality is possible in $F(S)$ if and only if 
\begin{equation} \label{e1}
z_0 x_0^{-1} s_j x_0 y_0^{-1} s_k^{-1}  y_0 z_0^{-1} = s_i^{\alpha}
\end{equation}
for some integer $\alpha$. Taking the quotient of $F(S)$ by its commutator subgroup $[F(S), F(S)]$, equation \ref{e1} gives
$$
\overline{s}_j \, \overline{s}_k^{-1} = \overline{s}_i^{\alpha},
$$
where $\overline{s}_j,  \overline{s}_k,  \overline{s}_i$ are the generators of the free abelian group $F(S)/[F(S), F(S)]$. This implies that $j = k$ and $\alpha = 0$. Thus, \eqref{e1} takes the form
$$
z_0 x_0^{-1} s_j x_0 y_0^{-1} s_j^{-1}  y_0 z_0^{-1} = 1,
$$
which further gives
$$
 s_j x_0 y_0^{-1} s_j^{-1}  y_0 x_0^{-1}= 1.
$$
This equality holds in $F(S)$ if and only if $y_0 x_0^{-1} = s_j^{\beta}$ for some integer $\beta$, that is, $y_0 = s_j^{\beta} x_0$.
Hence, $z, x, y$ have the form
$$
z = s_i^{z_0},~~ x = s_j^{x_0},~~ y= s_j^{x_0},
$$
which contradicts our assumption that $x \ne y$. Hence, $FQ(S)$ is a semi-latin quandle.
$\blacksquare$    \end{proof}          

For free involutory quandles, Propositions \ref{left associated form in quandles} and \ref{Lem:the canonical left associated form} take the following simpler form \cite[Lemma 4.4.2 and Theorem 4.4.4 ]{MR2634013}.

\begin{proposition}\label{free invol quandle element}
Any element in a free involutory quandle can be written in a left-associated product of the form
\begin{equation*}
x_n *x_{n-1}*\cdots *x_1\,,
\end{equation*}
where $x_1,x_2,\ldots,x_n\in S$ with $x_i\neq x_{i+1}$ for each $i$. Further, the multiplication of such products is given by 
\begin{equation*}
(x_n*x_{n-1}*\cdots*x_1)*(y_m*y_{m-1}*\cdots*y_1) = x_n*x_{n-1}*\cdots*x_1*y_1*\cdots*y_{m-1}*y_m*y_{m-1}*\cdots* y_1.
\end{equation*}
\end{proposition}

It also follows from Proposition \ref{orderable-implies-latin}(1) that a non-trivial involutory quandle is not right-orderable. We present the following result regarding the left-orderability of free involutory quandles \cite[Proposition 4.4]{MR4669143}.
 
\begin{theorem}\label{free invol left orderable} 
Free involutory quandles are left-orderable.
\end{theorem}

\begin{proof}
Let $F(S)$ be the free group on a non-empty set $S$. Recall that, $F(S)$ is the set of all conjugates of elements of $S$ in $F(S)$, equipped with the binary operation of conjugation. Given any involutory quandle $X$ and a set-theoretic map $\varphi: S \to X$, by the universal property of free quandles, $\varphi$ extends to a unique surjective quandle homomorphism $\tilde{\varphi}: FQ(S) \to X$. Since $X$ is involutory, $\tilde{\varphi}$  descends to a quandle homomorphism $\tilde{\varphi}: \big(FQ(S)\big)_2 \to X$, where $\big(FQ(S)\big)_2$ is the involutory quotient of $FQ(S)$. Hence, $\big(FQ(S)\big)_2$ is a model for the free involutory quandle on $S$. In fact, it follows that $\big(FQ(S)\big)_2$ is isomorphic to the subquandle of $\Conj(\hexstar_S \mathbb{Z}_2)$ consisting of all conjugates of generators in the free product $\hexstar_S \mathbb{Z}_2$.
\para

Let $\phi:\big(FQ(S)\big)_2 \to\Core \big(F(S)\big)$ be the map defined by 
\begin{equation*}
\phi(x_n*x_{n-1}*\cdots*x_1) = x_1^{d_1}x_2^{d_2}\cdots x_{n-1}^{d_{n-1}}x_n^{d_n}x_{n-1}^{d_{n-1}}\cdots x_2^{d_2}x_1^{d_1}\,,
\end{equation*} 
where $x_1,x_2,\ldots,x_n\in S$ with $x_i\neq x_{i+1}$ and $d_i=(-1)^{i+1}$ for $1\leq i\leq n$. By Proposition \eqref{free invol quandle element}, it can be checked that $\phi$ is a quandle homomorphism. Further, it turns out that $\phi$ is injective, and hence we have an embedding of $\big(FQ(S)\big)_2$ into $\Core \big(F(S)\big)$. Since free groups are well-known to be bi-orderable (see \cite[Theorem 2.1.9]{Deroin2014} and \cite{MR0031482}), it follows from Proposition \ref{conj right orderable}(2) that $\Core \big(F(S)\big)$ is a left-orderable quandle. Hence, the free involutory quandle $\big(FQ(S)\big)_2$ is left-orderable.
$\blacksquare$    
\end{proof}          

Theorems \ref{orderability-free-quandle} and \ref{free invol left orderable} suggest the following problem.

\begin{problem}
Determine whether free $k$-quandles are left-, right- or bi-orderable for $k \ge 3$.
\end{problem}

\begin{remark}{\rm 
In \cite{MR2069015}, Sikora introduced a natural topology on the set of left-orderings of a group and proved that for a countable group this space is compact, metrisable and totally disconnected. An analogous study was carried out in \cite{MR2320157} for the space of left-orderings on a groupoid. It was proved that such a space is compact and metrisable whenever the groupoid is countable.}
\end{remark}

The preceding remark together with Theorems \ref{orderability-free-quandle} and \ref{free invol left orderable} suggest the following problem.

\begin{problem}
Identify the space of right-orderings on free quandles and the space of left-orderings on free involutory quandles.
\end{problem}

Next, we present the following result for generalised Alexander quandles \cite[Proposition 4.5]{MR4669143}.

\begin{proposition}\label{biordered_alexander_quandle}
Let $<$ be a bi-ordering on a group $G$ and $\phi\in\Aut(G)$. Then the following statements are equivalent:
\begin{enumerate}
\item The linear order $<$ is a bi-ordering on $\Alex(G,\phi)$.
\item If $1$ is the identity element of $G$ and $1 <g$, then $1<\phi(g)<g$.
\end{enumerate}
\end{proposition}

\begin{proof}
For $(1) \Rightarrow (2)$, let $g\in G$ be such that $1<g$. Since $<$ is a bi-ordering on $\Alex(G,\phi)$, we have $1< g*1$ and $1<1*g$. This gives $1<\phi(g)$ and $1<\phi(g)^{-1}g$, and hence
$1<\phi(g)<g$.
\para
For $(2) \Rightarrow (1)$, let $x, y, z \in G$ be such that $x<y$. This gives $1<yx^{-1}$, and hence $1<\phi(y)\phi(x)^{-1}<yx^{-1}$. Since $<$ is a bi-ordering on $G$, it follows that $\phi(x) < \phi(y)$ and  $\phi(x)^{-1}x < \phi(y)^{-1} y$. Again using bi-ordering on $G$, we get  $\phi(x)\phi(z)^{-1}z<\phi(y)\phi(z)^{-1}z$ and $\phi(z)\phi(x)^{-1}x < \phi(z)\phi(y)^{-1}y$. By definition, this gives
$x*z<y*z$ and $z*x<z*y$, which is desired.
$\blacksquare$    \end{proof}          

We conclude by presenting a sufficient condition under which left-orderability fails in quandles.

\begin{proposition}\label{prop11}
Let $X$ be a quandle generated by a set $S$ such that the map $i_X: X \to \Conj \big(\Adj(X)\big)$ is injective. If there exist two distinct commuting elements in $\Adj(X)$ that are not inverses of each other and that are conjugates of elements from $i_X(S)^{\pm 1}$, then the quandle $X$ is not left-orderable.
\end{proposition}

\begin{proof}
Recall from Theorem \ref{presentation-of-adjoint-group} that the set $i_X(S)= \{ \textswab{a}_x \,\mid \,x \in S\}$ is a generating set for the adjoint group $\Adj(X)$. Let $i_X(S)^{-1}$ denote the set of inverses of elements in $i_X(S)$. Let $\textswab{a}_{a},\textswab{a}_{b}\in\Adj(X)$ with ${\textswab{a}_{a}}^{\,\pm1}\neq\textswab{a}_{b}$ be two commuting elements that are conjugates of elements from $i_X(S)^{\pm1}$. Then we can write
\begin{align*}
\textswab{a}_{a}&=\textswab{a}_{x_1}^{-d_1} \textswab{a}_{x_2}^{-d_2}\cdots \textswab{a}_{x_{m-1}}^{-d_{m-1}}\textswab{a}_{x_m}^{d_m}\:\!\textswab{a}_{x_{m-1}}^{d_{m-1}}\cdots \textswab{a}_{x_2}^{d_2}\textswab{a}_{x_1}^{d_1},\\
\textswab{a}_{b}&=\textswab{a}_{y_1}^{-\epsilon_1} \textswab{a}_{y_2}^{-\epsilon_2}\cdots \textswab{a}_{y_{n-1}}^{-\epsilon_{n-1}}\textswab{a}_{y_n}^{\epsilon_n}\:\! \textswab{a}_{y_{n-1}}^{\epsilon_{n-1}}\cdots \textswab{a}_{y_2}^{\epsilon_2} \textswab{a}_{y_1}^{\epsilon_1}\,,
\end{align*}
where $\textswab{a}_{x_i},\textswab{a}_{y_i}\in i_X(X)$ and $d_i,\epsilon_i\in\{-1,1\}$ for all $i$. For each $i$, there exist $x_i,y_i\in X$ such that $\textswab{a}_{x_i}=i_X(x_i)$ and $\textswab{a}_{y_i}=i_X(y_i)$. We see that
\begin{align*}
\textswab{a}_{a}^{d_m}&=i_X(x_1)^{-d_1} i_X(x_2)^{-d_2} \cdots i_X(x_{m-1})^{-d_{m-1}}i_X(x_m)\:\!i_X(x_{m-1})^{d_{m-1}}\cdots i_X(x_2)^{d_2}i_X(x_1)^{d_1}\\
&=i_X(x_m)*^{d_{m-1}}i_X(x_{m-1})*^{d_{m-2}}\cdots*^{d_1}i_X(x_1),~\textrm{by quandle operation in}~ \Conj \big(\Adj(X)\big)\\
&=i_X\!\left(x_m*^{d_{m-1}}x_{m-1}*^{d_{m-2}}\cdots*^{d_1}x_1\right),~\textrm{since $i_X$ is a quandle homomorphism}\\
&=i_X(a),
\end{align*}
and similarly
\begin{equation*}
\textswab{a}_{b}^{\epsilon_n}=i_X(b),
\end{equation*}
where $a=x_m*^{d_{m-1}}x_{m-1}*^{d_{m-2}}\cdots*^{d_1}x_1$ and $b=y_n*^{\epsilon_{n-1}}y_{n-1}*^{\epsilon_{n-2}}\cdots*^{\epsilon_1}y_1$. 
\para

Suppose on the contrary that the quandle $X$ is left-ordered with respect to a linear order $<$. Since $\textswab{a}_{a}^{\pm1}\neq\textswab{a}_{b}$, we get $\textswab{a}_{a}^{d_m}\neq \textswab{a}_{b}^{\epsilon_n}$, and hence $i_X(a)\neq i_X(b)$. This implies that $a\neq b$. Suppose that $a\, <\, b$. Then, we have $a=a*a\, <\, a*b$. Since $\textswab{a}_{b}^{-1}\textswab{a}_{a} \textswab{a}_{b}=\textswab{a}_{a}$, we have $\textswab{a}_{b}^{-1} \textswab{a}_{a}^{d_m}\textswab{a}_{b}=\textswab{a}_{a}^{d_m}$, and hence
\begin{equation*}
i_X(a*^{\epsilon_n}b)=i_X(a)*^{\epsilon_n}i_X(b)=i_X(b)^{-\epsilon_n}i_X(a)i_X(b)^{\epsilon_n}=\textswab{a}_{b}^{-1} \textswab{a}_{a}^{d_m}\textswab{a}_{b}=\textswab{a}_{a}^{d_m}=i_X(a).
\end{equation*}
The map $i_X$ being a monomorphism gives $a*^{\epsilon_n}b=a$, and hence $a*b=a$. Since $a\, <\, a*b$, we get a contradiction.
$\blacksquare$    \end{proof}

\begin{corollary}
Let $K$ be a prime knot such that its knot quandle $Q(K)$ is generated by a set $X$. If there exist two distinct commuting elements in the knot group $G(K)$ that are not inverses of each other and that are conjugates of elements from $i_{Q(K)}(X)^{\pm1}$, then $Q(K)$ is not left-orderable.
\end{corollary}

\begin{proof}
If $K$ is a prime knot, then by Theorem \ref{prime-quandles-injective-in-adjoint-group}, the map $i_{Q(K)}: Q\left(K\right)\to\Conj \big(G(K)\big)$ is a monomorphism of quandles. The result now follows from Proposition \ref{prop11}.
$\blacksquare$    
\end{proof} 
\bigskip
\bigskip 


\section{Properties of linear orderings on quandles}\label{properties-orderable-quandle}

In this section, we examine basic properties of linear orderings on quandles. Recall that a quandle essentially has two binary operations $*$ and $*^{-1}$. Thus, comprehending the behaviour of a linear order concerning both of these binary operations becomes imperative.
\para

Let $<$ be a linear order on a quandle $X$. We introduce an additional symbol $>$ by writing $y>x$ if and only if $x<y$.

\begin{definition}
Let $<$ be a linear order on a quandle $X$ and $\mathcal{O}$ be the set $\left\{=,<,>\right\}$. For a quadruple $( \diamond_1, \diamond_2, \diamond_3, \diamond_4)\in\mathcal{O}^4$, the order $<$ is said to be of {\it type} $( \diamond_1, \diamond_2, \diamond_3, \diamond_4)$ if the following  assertions hold for all $x,y,z\in X$ with $x<y$:
\begin{enumerate}[(1)]
\item $x*z\, \diamond_1\, y*z$.
\item $x*^{-1}z\, \diamond_2\, y*^{-1}z$.
\item $z*x\, \diamond_3\, z*y$.
\item $z*^{-1}x\, \diamond_4\, z*^{-1}y$.
\end{enumerate}
We say the order $<$ is of type $(\underline{\;\;}, \diamond_2,\underline{\;\;},\underline{\;\;})$ if the second condition is true, it is of the type $( \diamond_1,\underline{\;\;}, \diamond_3,\underline{\;\;})$ if the first and third conditions are true, it is of the type $( \diamond_1, \diamond_2,\underline{\;\;}, \diamond_4)$ if the first, second and fourth conditions are true, etc.
\end{definition}

The second quandle axiom implies that if $x< y$, then 
\begin{equation}\label{eq17}
x*z\neq y*z\qquad\textrm{and}\qquad x*^{-1}z\neq y*^{-1}z.
\end{equation}
\para

\begin{lemma} \cite[Lemma 3.2]{MR4330281} \label{lem2}
Let $<$ be a linear order on a quandle $X$ and $ \diamond\in\{<,>\}$. Then the order $<$ is of the type $( \diamond,\underline{\;\;},\underline{\;\;},\underline{\;\;})$ if and only if it is of the type $(\underline{\;\;}, \diamond,\underline{\;\;},\underline{\;\;})$.
\end{lemma}

\begin{proof}
Let $ \diamond\in\{<,>\}$. Define $ \diamond^{-1}$ to be $>$ if $ \diamond$ is $<$ and define it as $<$ if $ \diamond$ is $>$. Furthermore, we write $ \diamond^1$ as $ \diamond$. By \eqref{eq17}, we note that $x*z\, \diamond^d\,y*z$ and $x*^{-1}z\, \diamond^e\,y*^{-1}z$ for some $d,e\in\{-1,1\}$ whenever $x,y,z\in X$ and $x<y$.
\para
For the forward implication, suppose on the contrary that $x*^{-1}z\, \diamond^{-1}\,y*^{-1}z$ for some $x,y,z\in X$ with $x<y$. Since the order $<$ is of the type $( \diamond,\underline{\;\;},\underline{\;\;},\underline{\;\;})$, this implies that $\left(x*^{-1}z\right)*z\, \diamond^{-1}\,\left(y*^{-1}z\right)*z$ if $ \diamond$ is $<$ and $\left(x*^{-1}z\right)*z\, \diamond\,\left(y*^{-1}z\right)*z$ if $ \diamond$ is $>$. In other words, $\left(x*^{-1}z\right)*z>\left(y*^{-1}z\right)*z$, that is, $x>y$, which is a contradiction.
\para
For the backward implication, suppose on the contrary that $x*z\, \diamond^{-1}\,y*z$ for some $x,y,z\in X$ with $x<y$. This implies that $\left(x*z\right)*^{-1}z\, \diamond^{-1}\,\left(y*z\right)*^{-1}z$ if $ \diamond$ is $<$ and $\left(x*z\right)*^{-1}z\, \diamond\,\left(y*z\right)*^{-1}z$ if $ \diamond$ is $>$, since the order $<$ is of the type $(\underline{\;\;}, \diamond,\underline{\;\;},\underline{\;\;})$. This gives $\left(x*z\right)*^{-1}z>\left(y*z\right)*^{-1}z$, that is, $x>y$, which is a contradiction.
$\blacksquare$    \end{proof}          

\begin{lemma}\cite[Lemma 3.3]{MR4330281}\label{lem3}
Let $<$ be a linear order on a quandle $X$. Then the following statements are equivalent:
\begin{enumerate}[(1)]
\item The quandle $X$ is trivial.
\item The order $<$ is of the type $(\underline{\;\;},\underline{\;\;},=,\underline{\;\;})$.
\item The order $<$ is of the type $(\underline{\;\;},\underline{\;\;},\underline{\;\;},=)$.
\end{enumerate}
\end{lemma}

\begin{proof}
It is clear that (1) $\Rightarrow$ (2) and (1) $\Rightarrow$ (3).
\para
For (2) $\Rightarrow$ (1), let $x, y \in X$. If $x=y$, then by the idempotency axiom, $x*y=x$. If $x<y$ or $y<x$, then by (2), $x*y=x*x=x$. This proves that $x*y=x$ for all $x,y\in X$.
\para
For (3) $\Rightarrow$ (1), let  $x, y \in X$. If $x=y$, then by the idempotency axiom, $x*y=x$. If $x<y$ or $y<x$, then by (3), $x*^{-1}y=x*^{-1}x=x$. This implies that $(x*^{-1}y)*y=x*y$. Hence, by cancellation, we get $x*y=x$. This proves that $x*y=x$ for all $x,y\in X$.
$\blacksquare$    \end{proof}          

\begin{theorem}\cite[Theorem 3.4]{MR4330281}\label{th2}
Let $<$ be a linear order on a quandle $X$ of type $( \diamond_1, \diamond_2, \diamond_3, \diamond_4)$ for some $( \diamond_1, \diamond_2, \diamond_3, \diamond_4)\in \mathcal{O}^4$. Then the following assertions hold:
\begin{enumerate}[(1)]
\item $ \diamond_1, \diamond_2\in\left\{<,>\right\}$.
\item $ \diamond_1$ and $ \diamond_2$ are the same.
\item The quandle $X$ is trivial $\Leftrightarrow$ $ \diamond_3$ is the equality \textquoteleft\,$=$\textquoteright\, $\Leftrightarrow$ $ \diamond_4$ is the equality \textquoteleft\,$=$\textquoteright\,.
\item The quadruple $( \diamond_1, \diamond_2, \diamond_3, \diamond_4)$ is one of the following:\\
$(<,<,=,=)$,\quad $(<,<,<,>)$,\quad $(<,<,>,<)$\quad{or}\quad $(>,>,<,<)$.
\end{enumerate}
\end{theorem}

\begin{proof}
Assertion (1) follows from \eqref{eq17}, assertion (2) follows from Lemma \ref{lem2}, and assertion (3) follows from Lemma \ref{lem3}. 
\para

If $ \diamond_3$ or $ \diamond_4$ is the equality \textquoteleft\,$=$\textquoteright, then by (3), the quandle $X$ is trivial. In this case, $( \diamond_1, \diamond_2, \diamond_3, \diamond_4)$ must be $(<,<,=,=)$. If $ \diamond_3, \diamond_4\in\left\{<,>\right\}$, then by (1) and (2), $( \diamond_1, \diamond_2, \diamond_3, \diamond_4)$ must be one of the following quadruples:
\begin{align*}
\textrm{(a)}\;\,&(<,<,<,<),&\textrm{(b)}\;\,&(<,<,<,>),&\textrm{(c)}\;\,&(<,<,>,<),&\textrm{(d)}\;\,&(<,<,>,>),\\
\textrm{(e)}\;\,&(>,>,<,<),&\textrm{(f)}\;\,&(>,>,<,>),&\textrm{(g)}\;\,&(>,>,>,<),&\textrm{(h)}\;\,&(>,>,>,>).
\end{align*}
To prove assertion (4), we have to rule out the cases (a), (d), (f), (g) and (h). Let us consider the case (g) first. Suppose that $( \diamond_1, \diamond_2, \diamond_3, \diamond_4)=(>,>,>,<)$. Let $x,y,z\in X$ and $x<y$. Then we have
\begin{align}
\nonumber &z*x>z*y&&(\textrm{since}\, \diamond_3\,\textrm{is}\,>),\\
\nonumber  \Rightarrow\quad&(z*x)*^{-1}x<(z*y)*^{-1}x&&(\textrm{since}\, \diamond_2\,\textrm{is}\,>),\\
\Rightarrow\quad&z<z*y*^{-1}x&&(\textrm{by right cancellation}).\label{eq1 right cancel}
\end{align}
Furthermore, we have
\begin{align*}
\nonumber  &(z*y)*^{-1}x<(z*y)*^{-1}y&&(\textrm{since}\, \diamond_4\,\textrm{is}\,<),\\
\Rightarrow\quad&z*y*^{-1}x<z&&(\textrm{by right cancellation}).
\end{align*}
This is a contradiction to \eqref{eq1 right cancel}, and hence the case (g) does not arise. 
\para

Next, we consider the case (h). Suppose that $( \diamond_1, \diamond_2, \diamond_3, \diamond_4)=(>,>,>,>)$. Let $x,y,z\in X$ and $x<y$. Then we have 
\begin{align}
\nonumber&z*x>z*y&&(\textrm{since}\, \diamond_3\,\textrm{is}\,>),\\
\nonumber \Rightarrow\quad&x*^{-1}(z*x)<x*^{-1}(z*y)&&(\textrm{since}\, \diamond_4\,\textrm{is}\,>),\\
\Rightarrow\quad&x*^{-1}z*x<x*^{-1}y*^{-1}z*y.&& \label{eq2 order}
\end{align}
Furthermore, we have
\begin{align}
\nonumber &\left(x*^{-1}z\right)*x>\left(x*^{-1}z\right)*y&&(\textrm{since}\, \diamond_3\,\textrm{is}\,>),\\
\Rightarrow\quad&x*^{-1}z*x>x*^{-1}z*y. &&\label{eq3}
\end{align}
Combining \eqref{eq2 order} with \eqref{eq3}, we get
\begin{align}
\nonumber &x*^{-1}z*y<x*^{-1}y*^{-1}z*y,&&\\
\nonumber\Rightarrow\quad&\left(x*^{-1}z*y\right)*^{-1}y>\left(x*^{-1}y*^{-1}z*y\right)*^{-1}y&&(\textrm{since}\, \diamond_2\,\textrm{is}\,>),\\
\nonumber\Rightarrow\quad&x*^{-1}z>x*^{-1}y*^{-1}z&&(\textrm{by right cancellation}),\\
\nonumber\Rightarrow \quad&\left(x*^{-1}z\right)*z<\left(x*^{-1}y*^{-1}z\right)*z&&(\textrm{since}\, \diamond_1\,\textrm{is}\,>),\\
\Rightarrow\quad&x<x*^{-1}y&&(\textrm{by right cancellation}).\label{eq4}
\end{align}
We also have
\begin{align*}
&x*^{-1}x>x*^{-1}y&&(\textrm{since}\, \diamond_4\,\textrm{is}\,>),\\
\Rightarrow\quad&x>x*^{-1}y. &&
\end{align*}
This is a contradiction to \eqref{eq4}, and hence the case (h) does not arise. The cases (a), (d) and (f) can be ruled out similarly.
$\blacksquare$    \end{proof}          

\begin{corollary}
Let $<$ be a linear order on a quandle $X$. Then $X$ is trivial if and only if the linear order $<$ is of the type $(<,<,=,=)$.
\end{corollary}

We note that all the four possibilities for the quadruple $( \diamond_1, \diamond_2, \diamond_3, \diamond_4)$ in Theorem \ref{th2}(4) can be realized, as shown by the following example.

\begin{example}\label{ex1}
{\rm
Consider the group $(\mathbb{R},+)$. For a non-zero $u\in\mathbb{R}$, let $\phi_u$ be the automorphism of $\mathbb{R}$ given by $\phi_u(x)=ux$. Then for the Alexander quandle $\Alex(\mathbb{R},\phi_u)$, the quandle operation $*$ and its inverse operation $*^{-1}$ are given by $x*y=ux+(1-u)y$ and $x*^{-1} y=u^{-1}x+\left(1-u^{-1}\right)y$, respectively. With the usual linear order $<$ on $\mathbb{R}$, it is easy to check the following assertions:
\begin{enumerate}[(i)]
\item If $0<u<1$, then $<$ is a bi-ordering on $\Alex(\mathbb{R},\phi_u)$.
\item If $u<1$, then $<$ is a left-ordering on $\Alex(\mathbb{R},\phi_u)$.
\item If $0<u$, then $<$ is a right-ordering on $\Alex(\mathbb{R},\phi_u)$.
\end{enumerate}
Further, the following properties of the ordering $<$ on $\Alex(\mathbb{R},\phi_u)$ can be checked easily.
\begin{enumerate}[(i)]
\item If $u=1$, then the order $<$ is of the type $(<,<,=,=)$.
\item If $0<u<1$, then the order $<$ is of the type $(<,<,<,>)$.
\item If $1<u$, then the order $<$ is of the type $(<,<,>,<)$.
\item If $u<0$, then the order $<$ is of the type $(>,>,<,<)$.
\end{enumerate}}
\end{example}

\begin{remark}{\rm 
For $0<u<1$, the automorphism $\phi_u$ satisfy statement (2) of Proposition \ref{biordered_alexander_quandle}.}
\end{remark}

\begin{remark}{\rm 
It was asked in  \cite[Question 3.6]{MR4450681} whether there exists an infinite non-commutative bi-orderable quandle. We see that for $u\in(0,1)\setminus\left\{\frac{1}{2}\right\}$, the quandle $\Alex(\mathbb{R},\phi_u)$ with the usual order $<$ on $\mathbb{R}$ is an infinite non-commutative bi-orderable quandle, thereby answering the question in an affirmative.}
\end{remark}

Example \ref{ex1} confirms the general result below \cite[Proposition 3.8]{MR4330281}.

\begin{proposition}\label{prop3}
Let $<$ be a linear order on a quandle $X$. Then the order $<$ is a bi-ordering on $X$ if and only if it is of the type $(<,<,<,>)$.
\end{proposition}

\begin{proof}
It is trivial that if the ordering $<$ is of the type $(<,<,<,>)$, then it is a bi-ordering on $X$. Conversely, suppose that $<$ is a bi-ordering on $X$. Then we know that the ordering $<$ is of the type $(<,\underline{\;\;},<,\underline{\;\;})$. Hence, by Lemma \ref{lem2}, the ordering $<$ is of the type $(<,<,<,\underline{\;\;})$. Suppose on the contrary that $<$ is not of the type $(<,<,<,>)$. Then $z*^{-1}x<z*^{-1}y$\, for some $x,y,z\in X$ with $x<y$. This gives
\begin{align}
\nonumber &\left(z*^{-1}y\right)*x<\left(z*^{-1}y\right)*y&&(\textrm{since}\,<\,\textrm{is a left ordering on}\;X),\\
\Rightarrow\quad&z*^{-1}y*x<z&&(\text{by right cancellation})\label{eq16}.
\end{align}
Also, we have
\begin{align*}
&\left(z*^{-1}x\right)*x<\left(z*^{-1}y\right)*x&&(\textrm{since}\,<\,\textrm{is a right ordering on}\;X),\\
\Rightarrow\quad&z<z*^{-1}y*x&&(\text{by right cancellation}).
\end{align*}
But, this contradicts \eqref{eq16}, and hence  $<$ is not of the type $(<,<,<,>)$.
$\blacksquare$   
 \end{proof}          
\bigskip
\bigskip

\section{Constructions of orderable quandles}\label{construction-order-quandles-automorphisms}

An \textit{action} of a quandle $X$ on a quandle $Y$ is a quandle homomorphism $$\phi : X \to \Conj\big(\Aut(Y)\big),$$ where $\Aut(Y)$ is the group of quandle automorphisms of $Y$. In other words, if $\phi_{x}$ denotes the image of $x$ under $\phi$, then $\phi_{x*y}=\phi_{y}\phi_{x}\phi_{y}^{-1}$. Viewing any set $Y$ as a trivial quandle and noting that $\Aut(Y)=\Sigma_Y$, the symmetric group on $Y$,  we obtain the definition of an action of a quandle $X$ on a set $Y$.

\begin{example}{\rm 
Some basic examples of quandle actions are:
\begin{enumerate}
\item If $X$ is a quandle, then the map $\phi: X \to \Conj \big(\Aut(X)\big)$ given by $q \mapsto S_q$ is a quandle homomorphism. Thus, every quandle acts on itself by inner automorphisms.
\item Let $G$ be a group acting on a set $X$. That is, there is a group homomorphism $\phi: G \to \Sigma_X$. Viewing both $G$ and $\Sigma_X$ as conjugation quandles and observing that a group homomorphism is also a quandle homomorphism between corresponding conjugation quandles, it follows that the quandle $\Conj(G)$ acts on the set $X$.
\end{enumerate}}
\end{example}

We begin with the following result \cite[Theorem 4.2]{MR4330281}.

\begin{theorem}\label{th5}
If a semi-latin quandle is right-orderable, then it acts faithfully on a linearly ordered set by order-preserving bijections. Conversely, if a quandle acts faithfully on a well-ordered set by order-preserving bijections, then it is right-orderable.
\end{theorem}

\begin{proof}
Let $X$ be a semi-latin quandle that is right-ordered with respect to a linear order $<$. Then $\phi: X \to \Conj(\Sigma_X)$ given by $\phi(q)=S_q$ gives an action of $X$ on the ordered set $X$. Further, if $q, x, y \in X$ such that $x<y$, then right-orderability of $X$ implies that
$$\phi(q)(x)=S_q(x)=x*q < y*q=S_q(y)=\phi(q)(y).$$
Further, if $p, q \in X$ such that $\phi(p)=\phi(q)$, then $X$ being semi-latin implies that $p=q$. Hence, $X$ acts faithfully on itself by order-preserving bijections.
\para

Conversely, suppose that $\phi: X \to \Conj(\Sigma_S)$ is a faithful action of a quandle $X$ on a well-ordered set $S$ by order-preserving bijections. Let $<$ be the well-order on $S$. We use the order $<$ to define an order on the quandle $X$ as follows. For $p, q \in X$ with $p \neq q$, consider the set $A_{p, q}=\{x \in S \, \mid \,\phi(p)(x) \neq \phi(q)(x) \}$. Faithfulness of the action implies that $\phi(p) \neq \phi(q)$, and hence $A_{p, q}$ is a non-empty subset of $S$. Since $<$ is a well-ordering on $S$, the set $A_{p, q}$ must admit the smallest element, say $x_0$, with respect to $<$. We define $p\prec q$ if $\phi(p)(x_0) < \phi(q)(x_0)$ and $q\prec p$ if $\phi(q)(x_0) < \phi(p)(x_0)$. 
\para

It is enough to check transitivity of $\prec$ on $X$. Let $p \prec q$ and $ q \prec r$. Let $A_{p, q} = \{ x \in S  \,\mid\, \phi(p)(x) \neq \phi(q)(x) \}$, $A_{q, r} = \{ x \in S  \, \mid \, \phi(q)(x) \neq \phi(r)(x) \}$ and $A_{p, r} = \{ x \in S \, \mid \, \phi(p)(x) \neq \phi(r)(x) \}$. Since $p \prec q$ and $ q \prec r$, it follows that $A_{p, r}$ is non-empty. Let $x_0$, $y_0$ and $z_0$ be the smallest elements of the sets $A_{p, q}$, $A_{q, r}$ and $A_{p, r}$, respectively. Then we have the following cases:
\begin{enumerate}[(i)]
\item If $x _0 < y_0$, then $\phi(p)(x_0) < \phi(q)(x_0)=\phi(r)(x_0)$, which implies that $z_0 \leq x_0$. If $z_0 < x_0$, then $\phi(p)(z_0)=\phi(q)(z_0) \neq \phi(r)(z_0)$, which contradicts the fact that $y_0$ is the smallest element of $A_{q, r} $. Hence, $x_0=z_0$ and $ p \prec r$.
\item If $x_0 = y_0$, then $\phi(p)(x_0) < \phi(q)(x_0)< \phi(r) (x_0)$, which gives $z_0 \leq x_0$. If $z_0 < x_0$, then $\phi(p)(z_0)=\phi(q)(z_0) \neq \phi(r)(z_0)$, which is a contradiction to the fact that $y_0$ is the smallest element of $A_{q, r}$. Hence, $z_0=x_0$ and $p \prec r$.
\item If $x_0 > y_0$, then $\phi(p)(y_0)=\phi(q)(y_0) < \phi(r)(y_0)$, which further gives $z_0 \leq y_0$. If $z_0 < y_0$, then $\phi(p)(z_0)=\phi(q)(z_0) \neq \phi(r)(z_0)$, which is again a contradiction to the fact that $y_0$ is the smallest element of $A_{q, r}$. Hence, $z_0=y_0$ and $ p \prec r$.
\end{enumerate}
Now, suppose that $p, q, r \in X$ such that $p \prec q$. Let $A_{p, q}=\{x \in S \, \mid \,\phi(p)(x) \neq \phi(q)(x) \}$ and $A_{p*r, q*r}=\{x \in S \, \mid \,\phi(p*r)(x) \neq \phi(q*r)(x) \}$. Since both $A_{p, q}$ and $A_{p*r, q*r}$ are non-empty, they admit smallest elements, say $x_0$ and $y_0$, respectively. Note that the bijection $\phi(r)$ maps $A_{p, q}$ onto $A_{p*r, q*r}$. Since $\phi(r)$ is order-preserving with respect to $<$, we have $\phi(r)(x_0)=y_0$. Since $p \prec q$, we have $\phi(p)(x_0) < \phi(q)(x_0)$, which implies that $\phi(p)\phi(r)^{-1}(y_0) < \phi(q)\phi(r)^{-1}(y_0)$. Since $\phi$ is a quandle homomorphism and $\phi(r)$ is order-preserving, this gives 
\begin{eqnarray*}
\phi(p*r)(y_0)&=& \phi(p)*\phi(r)(y_0)\\
&=& \phi(r)\phi(p)\phi(r)^{-1}(y_0) \\
&<& \phi(r)\phi(q)\phi(r)^{-1}(y_0)\\
&=& \phi(q)*\phi(r)(y_0)\\
&=& \phi(q*r)(y_0),
\end{eqnarray*}
and hence $p*r \prec q*r$. Thus, $X$ is a right orderable quandle, and the proof is complete.
$\blacksquare$    \end{proof}          

Next, we give three constructions of orderable quandles. 

\begin{proposition}\label{prop6}
Let $(X_1,*)$ and $(X_2,\circ)$ be right-orderable quandles and $\sigma: X_1 \to \Conj\big(\Aut(X_2) \big)$ and $\tau: X_2 \to \Conj\big(\Aut(X_1) \big)$ be order-preserving quandle actions. Suppose that
\begin{enumerate}
\item $\tau(z)(x)* y=\tau\big(\sigma(y)(z)\big)(x* y)$ for $x, y \in X_1$ and $z \in X_2$,
\item $\sigma(z)(x)\circ y=\sigma\big(\tau(y)(z)\big)(x\circ y)$ for $x, y \in X_2$ and $z \in X_1$.
\end{enumerate}
Then $X=X_1 \sqcup X_2$ with the binary operation
$$
x\star y=\begin{cases}
x*y& \textrm{if}~~~x, y \in X_1, \\
x\circ y &\textrm{if}~~~ x, y \in X_2, \\
{\tau(y)}(x) &\textrm{if}~~~x \in X_1, y \in X_2, \\
{\sigma(y)}(x) &\textrm{if}~~~x \in X_2, y \in X_1,
\end{cases} 
$$
is a right-orderable quandle.
\end{proposition}

\begin{proof}
It follows from Proposition \ref{new2} that $X$ is a quandle. Let $<_1$ and $<_2$ be the right-orders on $X_1$ and $X_2$, respectively. Define an order $<$ on $X$ by setting $x< y$ if $x, y \in X_1$ and $x<_1 y$ or $x, y \in X_2$ and $x<_2 y$ or $x\in X_1$ and $y \in X_2$. A direct check shows that $<$ is indeed a linear order on $X$. We claim that $<$ turns $X$ into a right orderable quandle. Let $x, y, z \in X$ such that $x<y$. We have the following cases:
\begin{enumerate}[(i)]
\item $x, y, z \in X_1$ or $x, y,z \in X_2$: In this case, since $X_1$ and $X_2$ are right-orderable, we get $x\star z < y\star z$.
\item $x, y \in X_1$ and $z \in X_2$: In this case, since $\tau(z)$ is order preserving, we have $x \star z=\tau(z)(x) <_1 \tau(z)(y)=y \star z$, and hence $x \star z < y \star z$.
\item $x, y \in X_2$ and $z \in X_1$: In this case, $\sigma(z)$ being order preserving implies that $x \star z=\sigma(z)(x) <_2 \sigma(z)(y)=y \star z$, and hence $x \star z < y \star z$.
\item $x, z \in X_1$ and $y \in X_2$: In this case, $x \star z= x*z \in X_1$ and $y \star z= \sigma(z)(y) \in X_2$, and hence $x \star z < y \star z$.
\item $x \in X_1$ and $y, z \in X_2$: In this case, $x \star z= \tau(z)(x) \in X_1$ and $y \star z= y \circ z \in X_2$, and hence $x \star z < y \star z$.
\end{enumerate}
Thus, $X$ is a right-orderable quandle.
$\blacksquare$    \end{proof}          

If $\sigma: X_1 \to \id_{X_2}$ and $\tau: X_2 \to \id_{X_1}$ are the trivial actions, then conditions (1) and (2) of Proposition \ref{prop6} always hold. Thus, the disjoint union of two right-orderable quandles is right-orderable. It is clear that the disjoint union of two left-orderable quandles is not left-orderable. For example, any trivial quandle with more than one element is not left-orderable.
\para

Let $\{X_i, *_i\}_{i \in \Lambda}$ be a family of quandles and $X= \prod_{i \in \Lambda} X_i$ their cartesian product. Then $X$ is a quandle with $(x_i)\star (y_i)= (x*_i y_i)$ and called the {\it product quandle}. The following observation is rather immediate.

\begin{proposition}\label{prop1}
The product of right-orderable quandles is a right-orderable quandle. Similarly, the product of left-orderable (bi-orderable) quandles is a left-orderable (bi-orderable) quandle.
\end{proposition}

\begin{proof}
Let $\{X_i, *_i\}_{i \in \Lambda}$ be a family of right-orderable quandles. Let $<_i$ be the right-order on $X_i$ for $i \in \Lambda$ and $X$ their product quandle. By axiom of choice, we can take a well-ordering $<$ on the indexing set $\Lambda$. Let $(x_i), (y_i) \in X$ such that $(x_i) \neq (y_i)$. Then there exists the least index $\ell \in \Lambda$ such that $x_\ell \neq y_\ell$. We define $(x_i) \prec (y_i)$ if $x_\ell <_\ell y_\ell$ and $(y_i) \prec (x_i)$ if $y_\ell <_\ell x_\ell$. It is easy to check that $\prec$ is a linear order on $X$.
\para

Let $(x_i), (y_i), (z_i) \in X$ such that $(x_i) \prec (y_i)$. Then $x_\ell <_\ell y_\ell$, where $\ell$ is the least index such that $x_\ell \neq y_\ell$. The second quandle axiom in $X$ implies that $(x_i *_i z_i)= (x_i)\star(z_i) \neq (y_i) \star (z_i)=(y_i *_i z_i)$. It turns out that $\ell$ is also the least index for which $x_\ell *_\ell z_\ell \neq y_\ell *_\ell z_\ell$. Since $x_\ell <_\ell y_\ell$ and $X_\ell$ is right orderable, it follows that $x_\ell *_\ell z_\ell <_\ell y_\ell *_\ell z_\ell$. By definition of $\prec$, we have $ (x_i)\star(z_i) \prec (y_i) \star (z_i)$. Thus, $X$ is a right-orderable quandle. The second assertion follows analogously.
$\blacksquare$    \end{proof}          
\medskip

The construction of free racks and free quandles was generalised in  \cite{MR4129183}, which we discuss next. Let $G$ be a group and $A$ a subset of $G$. Then the set $A \times G$ becomes rack under the following binary operation
\begin{equation*}
(a,u)*(b,v)= (a, vbv^{-1} u)
\end{equation*}
for $a,b\in A$ and $u,v\in G$. The rack defined as above is known as a $(G,A)$-rack and is denoted by $R(G,A)$. 
\para

Let $Q(G,A)$ be the quotient of the set $A \times G$ by the equivalence relation $(a, vu) \sim (a,u)$ if and only if $v \in \C_G(a)=\{x \in G \,\mid \, xa=ax\}$. Let $[(a,u)] \in Q(G,A)$ denote the equivalence class of $(a,u) \in A \times G$. The set $Q(G,A)$ becomes quandle under the binary operation
\begin{equation*}
[(a,u)]*[(b,v)]= [(a,vbv^{-1}u)]
\end{equation*}
for $a,b\in A$ and $u,v\in G$, and this quandle is known as a $(G,A)$-quandle. For simplicity of notation, we write $(a,u)$ instead of $[(a,u)]$ throughout this section. There is a natural rack homomorphism $\epsilon: R(G,A) \to \Conj(G)$ defined as $\epsilon(a,u)=u^{-1}a u$. Moreover, this map induces a quandle homomorphism $\overline{\epsilon}: Q(G, A) \to \Conj(G)$ defined as $\overline{\epsilon}(a,u)=u^{-1} a u$.
\para

We prove the following result \cite[Theorem 5.1]{MR4330281}.

\begin{theorem}\label{th3}
Let $G$ be a group and $A$ a subset of $G$. Then the following assertions hold:
\begin{enumerate}
\item If $G$ is right-orderable, then the rack $R(G,A)$ is right-orderable.
\item If $G$ is bi-orderable, then the quandle $Q(G,A)$ is right-orderable.
\end{enumerate}
\end{theorem}

\begin{proof}
Let $<$ be a right ordering on $G$. We define a linear order $<'$ on $R(G,A)$ as follows. Let $(a,u)$ and $(b,v)$ be two distinct elements of $R(G,A)$.
\begin{enumerate}[(i)]
\item If $a \neq b$, define $(a,u)<'(b,v)$ if $a<b$ and $(b,v)<'(a,u)$ if $b<a$.
\item If $a=b$, define $(a,u)<'(a,v)$ if $u<v$ and $(a,v)<'(a,u)$ if $v<u$.
\end{enumerate}
Let $(a,u), (b,v), (c,w) \in R(G,A)$ such that $(a,u)<'(b,v)$. If $a \neq b$, then $a<b$, and hence $(a,u)*(c,w)<'(b,v)*(c,w)$. If $a=b$, then $u<v$. Since $G$ is right-ordered with respect to $<$, it follows that $uw^{-1}cw<vw^{-1}cw$, and hence $(a,u)*(c,w)<'(a,v)*(c,w)$. This shows that $R(G,A)$ is a right-orderable rack, which is assertion (1).
\medskip

Let $<$ be a bi-ordering on $G$. Define a linear order $<'$ on $Q(G,A)$ as follows. Let $(a,u)$ and $(b,v)$ be two distinct elements of $Q(G,A)$. 
\begin{enumerate}[(i)]
\item If $a \neq b$, then define $(a,u)<'(b,v)$ if $a<b$ and $(b,v)<'(a,u)$ if $b<a$.
\item If $a=b$, then we define the order using the image of $(a,u)$ and $(a,v)$ under the map $\overline{\epsilon}: Q(G,A) \to G$. Note that, if $(a,u) \neq (a,v)$ in $Q(G,A)$, then $\overline{\epsilon}(a,u)\neq \overline{\epsilon}(a,v)$. For, if $u^{-1}au=v^{-1}av$, then $vu^{-1}a=a vu^{-1}$; this implies that $vu^{-1} \in C_G(a)$ and hence $(a,u) =(a,v)$. Now, define $(a,u)<'(a,v)$ if $u^{-1} a u<v^{-1} a v$ and $(a,v)<'(a,u)$ if $v^{-1}av<u^{-1}au$.
\end{enumerate}
We claim that $Q(G,A)$ is right-ordered with respect to $<'$. Let $(a,u), (b,v), (c,w) \in R(G,A)$ such that $(a,u)<'(b,v)$. If $a \neq b$, then $a<b$, and hence $(a,u)*(c,w)<'(b,v)*(c,w)$. If $a=b$, then $u^{-1}au<v^{-1} a v$. Since $G$ is bi-ordered with respect to $<$, we have $w^{-1}c^{-1}wu^{-1}auw^{-1}cw<w^{-1}c^{-1}wv^{-1}avw^{-1}cw$, and hence $(a,u)*(c,w)<'(a,v)*(c,w)$. This shows that $Q(G,A)$ is right-orderable, which is assertion (2).
$\blacksquare$    
\end{proof}          

\begin{corollary}
Free racks are right-orderable.
\end{corollary}

\begin{proof}
Recall from Section \ref{section free quandles} of Chapter \ref{chapter preliminaries on racks and quandles} that the free rack $FR(A)$ on a set $A$ is the rack $R \big(F(A),A \big)$, where $F(A)$ is the free group on $A$. Since $F(A)$ is right-orderable, the result follows from Theorem \ref{th3}(1). $\blacksquare$    
\end{proof}

\begin{remark}
{\rm If $A$ is a set of representatives of conjugacy classes of a group $G$, then $Q(G, A)\cong \Conj(G)$, and we recover Proposition \ref{conj right orderable}(1). Similarly, if $G$ is the free group on a set $A$, then  $Q(G, A)$ is the free quandle on $A$, and we retrieve Theorem \ref{orderability-free-quandle}.}
\end{remark}
\bigskip
\bigskip


\section{Orderability of link quandles}\label{orderability-knot-quandle}
Problem 3.16 in \cite{MR4450681} asked whether link quandles are orderable. In this section, we investigate orderability of link quandles and provide solution to this problem in some cases. We begin with the following result.

\begin{proposition}\label{prop4}
Let $X$ be a quandle such that the natural map $i_X:X\to \Conj \big(\Adj(X) \big)$ is injective. If $\Adj(X)$ is a bi-orderable group, then $X$ is a right-orderable quandle.
\end{proposition} 

\begin{proof}
Since $\Adj(X)$ is a bi-orderable group, by Proposition \ref{conj right orderable}(1), $\Conj \big(\Adj(X) \big)$ is a right-orderable quandle. Since $i_X$ is injective, it follows that $X$ is right-orderable.
$\blacksquare$    \end{proof}           

\begin{corollary}
If $X$ is a commutative, latin or simple quandle such that $\Adj(X)$ is a bi-orderable group, then $X$ is right-orderable.
\end{corollary}

\begin{proof}
It is not difficult to see that the map $i_X$ is injective for a commutative, latin or simple quandle.
$\blacksquare$    \end{proof}          

\begin{corollary}\label{cor1}
If the knot group of a prime knot is bi-orderable, then its knot quandle is right-orderable.
\end{corollary}

\begin{proof}
Let $K$ be a prime knot such that its knot group $G(K)$ is bi-orderable. Since $K$ is prime, by Theorem \ref{prime-quandles-injective-in-adjoint-group}, the map $i_X: Q\left(K\right)\to \Conj \big(G(K)\big)$ is injective. Hence, by Proposition \ref{prop4}, the knot quandle $Q\left(K\right)$ is right-orderable.
$\blacksquare$    \end{proof}          

\begin{corollary}\label{cor4}
If all the roots of the Alexander polynomial of a fibered prime knot are real and positive, then its knot quandle is right-orderable.
\end{corollary}

\begin{proof}
Let $K$ be a fibered prime knot all the roots of whose Alexander polynomial are real and positive. Then, by \cite[Theorem 1.1]{MR1990838}, the knot group $G(K)$ is bi-orderable. The result now follows from Corollary \ref{cor1}.
$\blacksquare$    \end{proof}          

As a special case, it follows that the knot quandle of the figure eight knot is right-orderable.
\medskip

Next, we focus on orderability of link quandles of torus links. Recall that, two links $L_1$ and $L_2$ are called {\it weakly equivalent} if $L_1$ is ambient isotopic to either $L_2$ or the reverse of the mirror image of $L_2$. It is known that link quandles of weakly equivalent links are isomorphic \cite[Theorem 5.2 and Corollary 5.3]{MR1194995}. For any $m,n\geq1$, since the torus link $T(m,n)$ is invertible, it is weakly equivalent to its reverse, mirror image and the reverse of its mirror image, and hence their link quandles are isomorphic to that of $T(m,n)$. Recall that a torus link $T(m,n)$ is a knot (a one component link) if and only if $\gcd(m,n)=1$.
\para

We begin with the following result \cite[Proposition 7.1]{MR4330281}.

\begin{proposition}\label{prop5}
The link quandle of the torus link $T(m,n)$ is generated by $a_1,\ldots,a_m$ and has the following defining relations
\begin{equation*}
a_i=a_{n+i}*a_n*a_{n-1}*\cdots*a_1\quad\textrm{for}\;\;1 \le i \le m,
\end{equation*}
where $a_{mj+k}=a_k$ for $j\in\mathbb{Z}$ and $1 \le k \le m$.
\end{proposition}

\begin{proof}
Since a torus link $T(m,n)$ is the closure of the braid $\tau(m,n)=\left(\sigma_1\sigma_2\cdots\sigma_{m-1}\right)^n$, with reference to Figure \ref{toricbraid2}, it is enough to prove that
\begin{equation}\label{eq5}
c_i=a_{n+i}*a_n*a_{n-1}*\cdots*a_1\quad\textrm{for}\;\;1 \le i \le m.
\end{equation}

\begin{figure}[hbtp]
\centering
\includegraphics[height=2.8cm]{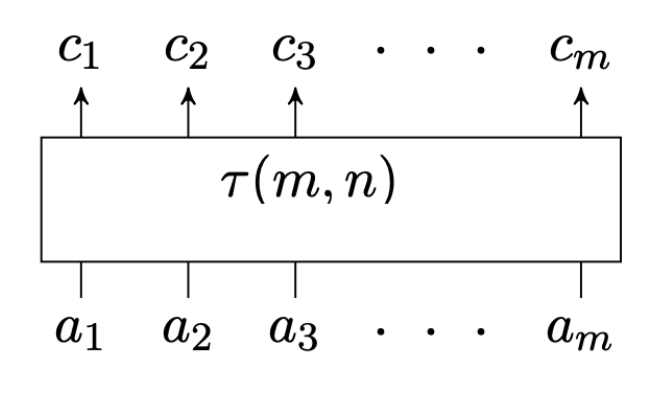}
\caption{Toric braid $\tau(m, n)$.}
\label{toricbraid2}
\end{figure}

We prove \eqref{eq5} by induction on $n$. Observing Figure \ref{toricbraid3}, we see that $c_i=a_{i+1}*a_1$ for $1 \le i \le m$.
\begin{figure}[hbtp]
\centering
\includegraphics[height=3.6cm]{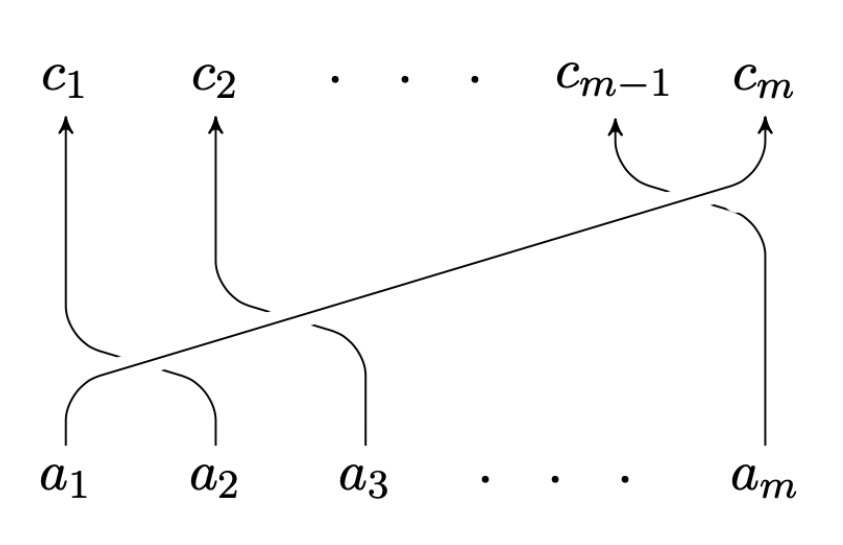}
\caption{Toric braid $\tau(m, 1)$.}
\label{toricbraid3}
\end{figure}
 Thus,  \eqref{eq5} hold for $n=1$. Suppose that \eqref{eq5} hold for a positive integer $n-1$. Since $\tau(m,n)=\tau(m,n-1)\tau(m,1)$ (see Figure \ref{toricbraid}), we have 
\begin{equation}\label{eq6}
c_i=b_{i+1}*b_1 \quad\textrm{for}\;\;1 \le i \le m
\end{equation}
where $b_{m+1}=b_1$. By induction hypothesis, we have
\begin{equation}\label{eq7}
b_{i+1}=a_{n+i}*a_{n-1}*a_{n-2}*\cdots*a_1\quad\textrm{for}\;\;1 \le i \le m.
\end{equation}

\begin{figure}[hbtp]
\centering
\includegraphics[height=5.6cm]{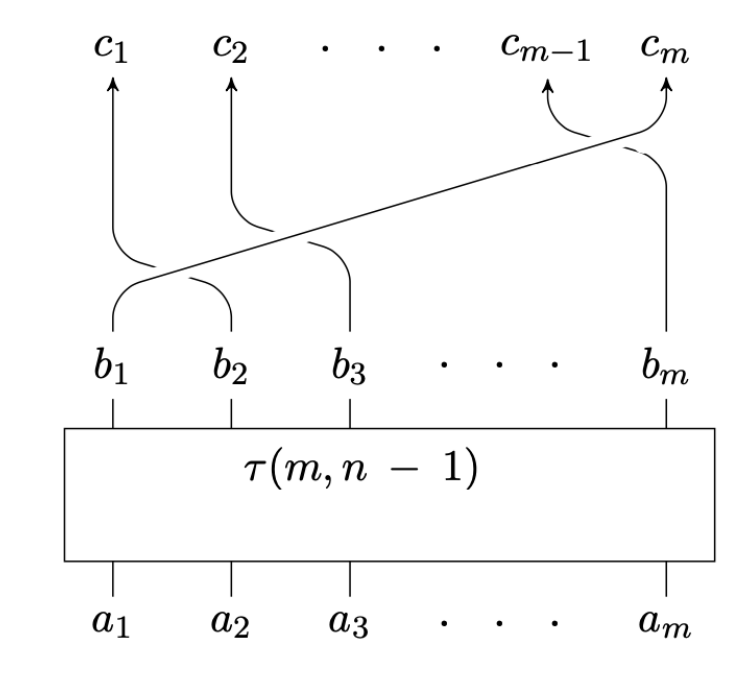}
\caption{Toric braid $\tau(m, n)$ seen as  $\tau(m, n-1)$ $\tau(m, 1)$.}
\label{toricbraid}
\end{figure}

\noindent Using \eqref{eq7} in \eqref{eq6}, we get
\begin{align*}
c_i&=\left(a_{n+i}*a_{n-1}*a_{n-2}*\cdots*a_1\right)*\left(a_n*a_{n-1}*a_{n-2}*\cdots*a_1\right)\\
&=a_{n+i}*a_{n-1}*a_{n-2}*\cdots*a_1*^{-1}a_1*^{-1}a_2*^{-1}\cdots*^{-1}a_{n-1}*a_n*a_{n-1}*\cdots*a_1\\
&=a_{n+i}*a_n*a_{n-1}*\cdots*a_1\quad\;\textrm{for}\;\;1 \le i \le m.
\end{align*}
This proves that \eqref{eq5} hold for all $n$. $\blacksquare$    
\end{proof}          

If $<$ is a right-ordering on a quandle $X$ and $x,y,z_1,z_2,\ldots,z_n\in X$ with $x\, \diamond\, y$ for $ \diamond\in\{<,>\}$, then we have
\begin{equation}\label{eq14}
x*z_1*z_2*\cdots*z_n\, \diamond\, y*z_1*z_2*\cdots*z_n~\textrm{and}~
x*^{-1}z_1*^{-1}z_2*^{-1}\cdots*^{-1}z_n\, \diamond\, y*^{-1}z_1*^{-1}z_2*^{-1}\cdots*^{-1}z_n.
\end{equation}

We now present the main result of this section \cite[Theorem 7.2]{MR4330281}.

\begin{theorem}\label{th4}
Let $m,n\geq2$ be integers such that one is not a multiple of the other. Then the link quandle of the torus link $T(m,n)$ is not right-orderable.
\end{theorem}

\begin{proof}
Note that the torus links $T(m,n)$ and $T(n,m)$ are ambient isotopic. Thus, we can assume that $m<n$ by switching $m$ and $n$, if required. If $d=\gcd(m,n)$, then $d<m$. By Proposition \ref{prop5}, the link quandle $Q \big(T(m,n) \big)$ is generated by $a_1,\ldots,a_m$ and has relations
\begin{equation}\label{eq10}
a_i=a_{n+i}*a_n*a_{n-1}*\cdots*a_1\quad\textrm{for}\;\;1 \le i \le m,
\end{equation}
where $a_{mj+k}=a_k$ for $j\in\mathbb{Z}$ and $1 \le k \le m$. Using \eqref{eq10}, we obtain
\begin{equation}\label{eq11}
a_i=a_{n+i}*a_n*a_{n-1}*\cdots*a_1\quad\textrm{for all}\;\;i\in\mathbb{Z},
\end{equation}
where $a_{mj+k}=a_k$ for $j\in\mathbb{Z}$ and $1 \le k \le m$. We can rewrite \eqref{eq11} as
\begin{equation}\label{eq12}
a_{i-n}=a_i*a_n*a_{n-1}*\cdots*a_1\quad\textrm{for all}\;\;i\in\mathbb{Z}.
\end{equation}
\noindent Also, \eqref{eq11} can be written as 
\begin{equation}\label{eq13}
a_{n+i}=a_i*^{-1}a_1*^{-1}a_2*^{-1}\cdots*^{-1}a_n\quad\textrm{for all}\;\;i\in\mathbb{Z}.
\end{equation}
	
Suppose on the contrary that the quandle $Q \big(T(m,n)\big)$ is right-ordered with respect to a linear order $<$. By the proof of Proposition \ref{prop5} (see Figures \ref{toricbraid2}, \ref{toricbraid3} and \ref{toricbraid}), the generators $a_1,\ldots,a_m$ of $Q \big(T(m,n)\big)$ correspond to some of the arcs in the standard diagram of the closed toric braid representing $T(m,n)$.  Note that $i_X(a_1),\ldots,i_X(a_m)$ are the meridional elements that generate the link group $G\big(T(m,n)\big)$, where $i_X:Q\big(T(m,n)\big)\to G\big(T(m,n)\big)$ is the natural map. According to \cite[Corollary 1.5]{MR0889376}, the elements $i_X(a_1),\ldots,i_X(a_m)$ must be pairwise distinct in $G \big(T(m,n)\big)$, and hence so are the elements $a_1,\ldots,a_m$ in $Q\big(T(m,n)\big)$. In particular, we have $a_1\neq a_{d+1}$, and hence $a_1\, \diamond\, a_{d+1}$ for some $ \diamond\in\{<,>\}$. A repeated application of \eqref{eq14} together with \eqref{eq12} and \eqref{eq13} yields
\begin{equation}\label{eq15}
a_{nk+1}\, \diamond\, a_{nk+d+1}\quad\textrm{for all}\;\;k\in\mathbb{Z}.
\end{equation}
Let $l$ be an integer. Since $\gcd(m,n)=d$, we have $dl=mj+nk$ for some integers $j$ and $k$. This implies that $nk+1\equiv dl+1\;(\!\!\!\!\mod m)$ and $nk+d+1\equiv dl+d+1\;(\!\!\!\!\mod m)$. By \eqref{eq15}, we have $a_{dl+1}\, \diamond\, a_{dl+d+1}$. Thus, $a_{dl+1}\, \diamond\, a_{dl+d+1}$ for any integer $l$. Using this repeatedly, we get $a_1\, \diamond\, a_{d+1}\, \diamond\, a_{2d+1}\, \diamond\,\cdots\, \diamond\, a_{cd+1}\, \diamond\, a_1$, where $c=\frac{m}{d}-1$. This implies that $a_1<a_1$ or $a_1>a_1$, a contradiction.
$\blacksquare$    \end{proof}          

As a consequence of Theorem \ref{th4}, we retrieve the following result of Perron and Rolfsen \cite[Proposition 3.2]{MR1990838}.

\begin{corollary}\label{cor3}
The knot group of a non-trivial torus knot is not bi-orderable.
\end{corollary}

\begin{proof}
Let $K$ be a non-trivial torus knot. By Theorem \ref{th4}, the knot quandle of $K$ is not right-orderable. Hence, by Corollary \ref{cor1}, the knot group of $K$ is not bi-orderable.
$\blacksquare$    \end{proof}          

We conclude with the following result.

\begin{corollary}
The knot quandle of the trefoil knot is neither left nor right orderable.
\end{corollary}

\begin{proof}
Note that the trefoil knot is the torus knot $T(2,3)$. By Theorem \ref{th4}, the knot quandle $Q\big(T(2,3)\big)$ is not right-orderable. We claim that it is not left-orderable as well. The quandle $Q\big(T(2,3)\big)$ is generated by $a_1, a_2$ and has  relations 
\begin{align}
a_1&=a_2*a_1*a_2,\label{eq8}\\
a_2&=a_1*a_2*a_1.\label{eq9}
\end{align}
Assume contrary that $Q\big(T(2,3)\big)$ is left-ordered with respect to a linear order $<$. Since $Q\big(T(2,3)\big)$ is non-trivial, we must have $a_1\neq a_2$. Hence, $a_1\, \diamond\, a_2$ for some $ \diamond\in\{<,>\}$. We see that
\begin{align*}
&\quad a_1\, \diamond\, a_2&&\\
\Rightarrow&\quad a_1*a_1\, \diamond\, a_1*a_2&&\textrm{(since $<$ is left-ordering)},\\
\Rightarrow&\quad a_1\, \diamond\, a_1*a_2&&\textrm{(by idempotency)},\\
\Rightarrow&\quad a_2*a_1\, \diamond\, a_2*(a_1*a_2)&&\textrm{(since $<$ is left-ordering)},\\
\Rightarrow&\quad a_2*a_1\, \diamond\, a_2*a_1*a_2,&&\\
\Rightarrow&\quad a_2*a_1\, \diamond\, a_1&&\textrm{(by \eqref{eq8})},\\
\Rightarrow&\quad a_1*(a_2*a_1)\, \diamond\, a_1*a_1&&\textrm{(since $<$ is left-ordering)},\\
\Rightarrow&\quad a_1*a_2*a_1\, \diamond\, a_1,&&\\
\Rightarrow&\quad a_2\, \diamond\, a_1&&\textrm{(by \eqref{eq9})}.
\end{align*}
This gives a contradiction, since we cannot have $a_1\, \diamond\, a_2$ and $a_2\, \diamond\, a_1$ together.
$\blacksquare$   
 \end{proof} 

\begin{remark}
{\rm 
Further work on the orderability of knot quandles includes a non-left-orderability criterion for involutory quandles of non-split links \cite{arXiv.2310.05735}, used to show that those of non-trivial alternating links are not left-orderable. Biorderability of knot quandles for prime knots up to eight crossings has also been studied \cite{arXiv:2505.19573}.
}
\end{remark}


\chapter{Quandle rings and algebras}\label{chaper quandle rings}

\begin{quote}
This chapter explores the construction and properties of quandle rings and algebras, drawing analogies with the classical theory of group rings and algebras. We examine fundamental aspects of quandle rings and how their algebraic structure reflects the properties of the underlying quandles. Key topics include power-associativity, commutator width, and the behavior of associated graded rings. The chapter also includes a discussion of the isomorphism problem for quandle rings, highlighting structural challenges and open questions in the field.
\end{quote}
\bigskip

\section{Properties of quandle rings and algebras}\label{section1 of chap quandle rings}
Let $(X, *)$ be a quandle and $\mathbb{k}$ an associative ring with unity {\bf 1}.  Let $e_x$ be a unique symbol corresponding to each $x \in X$. Let $\mathbb{k}[X]$ be the set of all formal expressions of the form $\sum_{x \in X }  \alpha_x e_x$, where $\alpha_x \in \mathbb{k}$ such that all but finitely many $\alpha_x=0$.  The set $\mathbb{k}[X]$ has a free $\mathbb{k}$-module structure with basis $\{e_x \, \mid \, x \in X \}$ and admits a product given by 
 $$ \Big( \sum_{x \in X }  \alpha_x e_x \Big) \Big( \sum_{ y \in X }  \beta_y e_y \Big)
 =   \sum_{x, y \in X } \alpha_x \beta_y e_{x * y},$$
where $x, y \in X$ and $\alpha_x, \beta_y \in \mathbb{k}$. This turns $\mathbb{k}[X]$ into a ring (rather a $\mathbb{k}$-algebra) called the \index{quandle ring} {\it quandle ring} or the \index{quandle algebra} {\it quandle algebra} of $X$ with coefficients in $\mathbb{k}$. Even though the coefficient ring $\mathbb{k}$ is associative, the quandle ring $\mathbb{k}[X]$ is non-associative when $X$ is a non-trivial quandle. The quandle $X$ can be identified as a subset of $\mathbb{k}[X]$ via the natural map $x \mapsto {\bf 1}e_x=e_x$.  We will see that, unlike groups, the quandle ring structure of a trivial quandle is interesting in its own.

\begin{remark}{\rm 
 If $X$ is a rack, then $\mathbb{k}[X]$ is referred to as the \index{rack ring}{\it rack ring} or the \index{rack algebra} {\it rack algebra} of $X$ with coefficients in $\mathbb{k}$. Most of the concepts introduced in this chapter also apply to racks.}
\end{remark}

\begin{remark}{\rm 
Note that a quandle with a (left) identity element has only one element. Thus, $\mathbb{k}[X]$ is a non-associative ring without unity, unless $X$ has only one element.}
\end{remark}

Analogous to group rings, we define the \index{augmentation map} {\it augmentation map}
$$\varepsilon: \mathbb{k}[X] \to \mathbb{k}$$
 by setting $$\varepsilon \Big(\sum_{x \in X} \alpha_x e_x \Big)= \sum_{x \in X} \alpha_x.$$
Clearly, $\varepsilon$ is a surjective ring homomorphism. Let $\Delta_{\mathbb{k}}(X)$ be the kernel of the homomorphism $\varepsilon$. Then $\Delta_{\mathbb{k}}(X)$ is a two-sided ideal of $\mathbb{k}[X]$ and referred to as the {\it augmentation ideal} of $\mathbb{k}[X]$. Further, we have
$$\mathbb{k}[X]/\Delta_{\mathbb{k}}(X) \cong \mathbb{k}$$ as rings. In the case $\mathbb{k} = \mathbb{Z}$, we denote the augmentation ideal simply by $\Delta(X)$. The following result is straightforward.

\begin{proposition}
Let $X$ be a rack and $\mathbb{k}$ an associative ring with unity. Then $\{e_x- e_y\,\mid\,x, y \in X \}$ is  a generating set for $\Delta_{\mathbb{k}}(X)$ as a $\mathbb{k}$-module. Further, if $x_0 \in X$ is a fixed element, then the set $\big\{e_x-e_{x_0} \,\mid\, x \in X \setminus \{ x_0\} \big\}$ is a basis for $\Delta_{\mathbb{k}}(X)$ as an $\mathbb{k}$-module.
\end{proposition}

The following is a noteworthy observation.

\begin{proposition}
Let $X$ be a quandle and $\mathbb{k}$ an associative ring with unity. Then $e_x e_y+ e_y e_x \equiv e_x+e_y~ \mod \Delta_{\mathbb{k}}^2(X)$ for all $x, y \in X$.
\end{proposition}

\begin{proof}
It follows from the first quandle axiom that $e_x^2=e_x$ in $\mathbb{k}[X]$ for all $x \in X$. We see that
\begin{eqnarray*}
e_x & = & (e_x-e_y+e_y)^2\\
& = & (e_x-e_y)^2+e_y^2+(e_x-e_y) e_y+e_y(e_x-e_y)\\
& \equiv & e_x e_y+e_y e_x-e_y \mod  \Delta_{\mathbb{k}}^2(X),
\end{eqnarray*}
and hence $e_x e_y+e_y e_x \equiv e_x+e_y \mod \Delta_{\mathbb{k}}^2(X)$ for all $x, y \in X$. $\blacksquare$
\end{proof}

Clearly, if $Y$ is a subrack of a rack $X$, then $\Delta_{\mathbb{k}}(Y) \subseteq \Delta_{\mathbb{k}}(X)$. It is natural to look for conditions under which $\Delta_{\mathbb{k}}(Y)$ is a two-sided ideal of $\mathbb{k}[X]$. For trivial racks, we have the following result.

\begin{proposition}
Let $X$ be a trivial rack, $Y$ a subrack of $X$ and $\mathbb{k}$ an associative ring with unity. Then $\Delta_{\mathbb{k}}(Y)$ is a two-sided  ideal of $\mathbb{k}[X]$.
\end{proposition}

\begin{proof}
Note that $\Delta_{\mathbb{k}}(Y)$ is generated as an $\mathbb{k}$-module by the set $\{e_y-e_z \,\mid \,y,z \in Y \}$. Then, for any $\sum_{x \in X}\alpha_x e_x \in \mathbb{k}[X]$, we have
$$\big(\sum_{x \in X} \alpha_x e_x\big)(e_y -e_z)=\sum_{x \in X} \alpha_x (e_x e_y - e_x  e_z)=\sum_{x \in X}\alpha_x (e_x - e_x)= 0 \in \Delta_{\mathbb{k}}(Y),$$
and
$$(e_y -e_z) \big(\sum_{x \in X}\alpha_x e_x\big)=\sum_{x \in X}\alpha_x (e_y  e_x -e_z  e_x)=\sum_{x \in X}\alpha_x (e_y -e_z) \in \Delta_{\mathbb{k}}(Y).$$
Hence, $\Delta_{\mathbb{k}}(Y)$ is a two-sided ideal of $\mathbb{k}[X]$. $\blacksquare$
\end{proof}

The next result characterises trivial quandles in terms of their augmentation ideals \cite[Theorem 3.5]{MR3977818}.

\begin{theorem}\label{deltasqzero}
Let $X$ be a quandle and $\mathbb{k}$ an associative ring with unity. Then $X$ is trivial if and only if $\Delta_{\mathbb{k}}^2(X)=\{0\}$.
\end{theorem}

\begin{proof}
Suppose that $\Delta_{\mathbb{k}}^2(X)=\{0\}$. Let $x_0 \in X$ be a fixed element. Observe that $\Delta_{\mathbb{k}}^2(X)$ is generated as an $\mathbb{k}$-module by the set $\big\{(e_x-e_{x_0})(e_y-e_{x_0}) \, \mid \, x, y \in X\setminus \{x_0 \} \big\}$. It follows that $(e_x-e_{x_0})(e_y-e_{x_0})=0$ for all $x, y \in X\setminus \{x_0 \}$. In particular, $(e_x-e_{x_0})(e_x-e_{x_0})=0$ for all $x \in X\setminus \{x_0 \}$, which yields
$$e_x-e_{x_0} e_x - e_x  e_{x_0}+e_{x_0}=0$$
for all $x \in X\setminus \{x_0 \}$. Since it is an expression in $\mathbb{k}[X]$, the terms must cancel off with each other. Suppose that $e_{x_0} e_x=e_x$ and $e_x e_{x_0}=e_{x_0}$ for some $x \in X\setminus \{x_0 \}$. Also, we have $e_x  e_x=e_x$ by the first quandle axiom. Thus, by second quandle axiom, we have $e_x=e_{x_0}$, which is a contradiction. Hence, we must have $e_x  e_{x_0}=e_x$ and $e_{x_0} e_x=e_{x_0}$ for all $x \in X\setminus \{x_0 \}$. This means $x_0$ acts trivially on all elements of $X$. Since $x_0$ was arbitrary, it follows that $X$ is trivial.
\para

Conversely, suppose that $X$ is a trivial quandle. Let $y, z, y', z' \in X$. Then
$$(e_y-e_z)(e_{y'}-e_{z'})=e_y  e_{y'}-e_{z}  e_{y'}-e_y  e_{z'}+e_{z} e_{z'}=e_y-e_z-e_y+e_z=0,$$
and hence $\Delta_{\mathbb{k}}^2(X)=\{0\}$. $\blacksquare$
\end{proof}

\begin{corollary}
A group $G$ is abelian if and only if $\Delta_{\mathbb{k}}^2\big(\Conj(G)\big)=\{0\}$.
\end{corollary}

\begin{remark}{\rm 
Clearly, if $X$ is a trivial rack, then $\Delta_{\mathbb{k}}^2(X)=\{0\}$. However, the converse is not true for racks that are not quandles. For example, take $X= \{0, 1, \dots, n-1\}$ with the rack structure given by $$i*j=i+1 \mod n.$$ Then $X$ is not a trivial rack. On the other hand, for $0 \le i, j \le n-1$, we have
$$(e_{i}-e_{1})(e_{j}-e_{1})=e_{i}e_{j} - e_{i}e_{1} - e_{1}e_{j} + e_{1}e_{1}=e_{i+1}-e_{i+1}-e_{2}+e_{2}=0,$$ and hence  $\Delta_{\mathbb{k}}^2(X)=\{0\}$.}
\end{remark}

 \begin{remark}{\rm 
 If $\phi:X \to Y$ is a homomorphism of quandles, then it extends linearly to a homomorphism of corresponding quandle rings $\mathbb{k}[X]\to \mathbb{k}[Y]$. Throughout the text, we denote this induced homomorphism by $\hat{\phi}: \mathbb{k}[X]\to \mathbb{k}[Y]$.}
 \end{remark}
 \bigskip
 \bigskip 
 

\section{Subquandles and ideals of quandle rings}\label{section4 quandle rings}
In this section, we explore the relationships between subquandles of a given quandle and the ideals of its associated quandle ring. We first consider subquandles associated to ideals.
\para
 Let $X$ be a quandle and $\mathbb{k}$ an associative ring with unity.  For each $x_0 \in X$ and each two-sided ideal $I$ of $\mathbb{k}[X]$, we define
$$X_{I, x_0}=\{ x \in X \, \mid \, e_x-e_{x_0} \in I \}.$$

Note that, if $I= \Delta_{\mathbb{k}}(X)$, then $X_{I, x_0}=X$, and if $I = \{0\}$, then $X_{I,x_0} = \{x_0\}$. In general, we have the following result \cite[Theorem 4.1]{MR3977818}. 

\begin{theorem}\label{ideal-to-subquandle}
Let $X$ be a finite quandle and $\mathbb{k}$ an associative ring with unity. Then for each $x_0 \in X$ and a two-sided ideal $I$  of $\mathbb{k}[X]$, the set $X_{I, x_0}$ is a subquandle of $X$. Further, there is a subset $\{ x_1, \ldots, x_m \}$ of $X$  such that $X$ is the disjoint union
$$
X = X_{I,x_1} \sqcup \cdots \sqcup X_{I,x_m}.
$$
\end{theorem}

\begin{proof}
Clearly, we have $x_0 \in X_{I, x_0}$. If $x, y \in X_{I, x_0}$, then 
$$e_{x*y}- e_{x_0}=e_{x}e_{y}- e_{x_0}= e_x e_y- e_{x_0} e_y+e_{x_0} e_y -e_{x_0} e_{x_0}=(e_x-e_{x_0})e_y+e_{x_0}(e_y -e_{x_0}) \in I,$$
and hence $x* y \in X_{I, x_0}$. Further, the inner automorphism $S_y$ restricts to a map $X_{I, x_0}  \to X_{I, x_0}$. Since $S_y$ is injective and $X$ is finite, it follows that  $S_y\big(X_{I, x_0} \big)= X_{I, x_0}$. Hence, if $x, y \in X_{I, x_0}$, then there exists a unique element $z \in X_{I, x_0} $ such that $x=z*y$, which proves that $X_{I, x_0}$ is a subquandle of $X$.
\para

For the second assertion, it is sufficient to prove that if $x_0, y_0 \in X$, then the subquandles $X_{I, x_0}$ and $X_{I, y_0}$ intersect if and only if they are equal. Let $z \in X_{I, x_0} \cap X_{I, y_0}$. Then we have $e_{x_0} - e_{y_0} \in I$, which further implies that $x_0 \in X_{I, y_0}$. Now, if $x \in X_{I, x_0}$, then $e_x-e_{y_0}=e_x-e_{x_0}+e_{x_0}-e_{y_0} \in I$, which implies that $X_{I, x_0} \subseteq X_{I, y_0}$. By interchanging roles of $x_0$ and $y_0$, we get $X_{I, x_0} = X_{I, y_0}$. $\blacksquare$
\end{proof}

\begin{remark}{\rm 
If $X$ is an involutory quandle (not necessarily finite), then $X_{I, x_0}$ is always a subquandle of $X$, being closed under the quandle multiplication.}
\end{remark}

Given a two-sided ideal $I$ of $\mathbb{k}[X]$ and elements $x_0, y_0 \in X$, it is natural to ask whether there is any relation between the subquandles $X_{I, x_0}$ and $X_{I, y_0}$. We answer this question in the following result \cite[Theorem 4.3]{MR3977818}.

\begin{theorem}\label{id-subq}
Let $X$ be a finite involutory quandle and $\mathbb{k}$ an associative ring with unity.  If $I$ is a two-sided ideal of $\mathbb{k}[X]$ and $x_0, y_0 \in X$ are in the same orbit under action of $\Inn(X)$, then $X_{I, x_0} \cong X_{I, y_0}$.
\end{theorem}

\begin{proof}
We first claim that if $x_0, y_0 \in X$, then the map $f_{y_0}: X_{I, x_0} \to X_{I, x_0*y_0}$ given by $f_{y_0}(x)=x*y_0$ is an injective quandle homomorphism. For, if $x \in X_{I, x_0}$, then $e_{x}-e_{x_0} \in I$. Consequently, $e_{x} e_{y_0}-e_{x_0} e_{y_0}= (e_{x}-e_{x_0})e_{y_0} \in I$, which further implies $f_{y_0}(x)=x*y_0 \in X_{I, x_0*y_0}$. The claim now follows by observing that $f_{y_0}$ is simply the restriction of the inner automorphism $S_{y_0}$ on the subquandle $X_{I, x_0}$.
\para
Now, suppose that $x_0, y_0 \in X$ such that there exists $f \in \Inn(X)$ with $f(x_0)=y_0$. Since $X$ is involutory, we can write $f= S_{x_k}  \cdots  S_{x_1}$ for some $x_i \in X$. Then we have $$y_0= S_{x_k}  \cdots   S_{x_1}(x_0)= \big(\cdots ((x_0* x_1) * x_2) \cdots \big) * x_k.$$
It follows from the preceding claim that there is a sequence of embeddings of subquandles of $X$
$$ X_{I, x_0} \hookrightarrow X_{I, x_0* x_1} \hookrightarrow X_{I, (x_0* x_1)* x_2} \hookrightarrow \cdots \hookrightarrow  X_{I, y_0}.$$
Thus, we obtain an embedding of $X_{I, x_0}$ into $X_{I, y_0}$. Writing $x_0= S_{x_1}  \cdots   S_{x_k}(y_0)$, we obtain an embedding of $X_{I, y_0}$ into $X_{I, x_0}$. Since $X$ is finite, it follows that $X_{I, x_0} \cong X_{I, y_0}$. $\blacksquare$
\end{proof}

As an immediate consequence, we have the following result.

\begin{corollary}
Let $X$ be a finite connected involutory quandle and $\mathbb{k}$ an associative ring with unity. Then $X_{I, x_0} \cong X_{I, y_0}$ for any ideal $I$ of $\mathbb{k}[X]$ and $x_0, y_0 \in X$.
\end{corollary}

For example,  dihedral quandles of odd order are connected  and involutory.

\begin{remark}{\rm 
It is worth noting that the results of the preceding discussion hold for racks as well. It would be interesting to explore whether Theorem \ref{id-subq} holds if $X$ is not involutory.}
\end{remark}
\bigskip

Next, we proceed in the reverse direction of associating an ideal of $\mathbb{k}[X]$ to a subquandle of $X$. Let $f:X \to Z$ be a quandle homomorphism. Consider the equivalence relation $\sim$ on $X$ given by $x_1\sim x_2$ if $f(x_1)=f(x_2)$. Let $X/_\sim$ be the set of equivalence classes, where equivalence class of an element $x$ is denoted by $$X_x:= \{x' \in X \, \mid \, f(x')=f(x) \}.$$

\begin{proposition}
$X_x$ is a subquandle of $X$ for each $x \in X$.
\end{proposition}

\begin{proof}
Clearly, $X_x$ is non-empty since $x \in X_x$. Let $x_1, x_2 \in X_x$. Then $f(x_1 * x_2)=f(x_1) *f(x_2)=f(x)*f(x)=f(x)$, and hence $x_1*x_2 \in X_x$. Further, if $x_3 \in X$ is the unique element such that $x_3*x_1=x_2$, then applying $f$ yields $f(x_3) *f(x)=f(x)$. This together with the first quandle axiom imply that $f(x_3)=f(x)$, and hence $x_3 \in X_x$. $\blacksquare$
\end{proof}

The following is an analogue of the first isomorphism theorem for quandles \cite[Theorem 4.9]{MR3977818}.

\begin{theorem}\label{quotient-quandle}
Let $f:X \to Z$ be a quandle homomorphism. Then the binary operation given by $X_{x_1} \circ X_{x_2}= X_{x_1*x_2}$ gives a quandle structure on $X/_\sim$. Further, if $f$ is surjective, then $X/_\sim \cong Z$ as quandles.
\end{theorem}

\begin{proof}
The binary operation $X_{x_1} \circ X_{x_2}= X_{x_1*x_2}$ is clearly well-defined. We only need to check the second quandle axiom. Let $X_{x_1}, X_{x_2} \in X/_\sim$.  If $x_3 \in X$ is the unique element such that $x_3*x_1=x_2$, then $f(x_3)*f(x_1)=f(x_2)$ and $X_{x_3} \circ X_{x_1}=X_{x_2}$. Suppose that there exists another element $X_{x_3'} \in X/_\sim$ such that $X_{x_3'} \circ X_{x_1}=X_{x_2}$, then  $f(x_3')*f(x_1)=f(x_2)$. By the second quandle axiom, we must have $f(x_3)=f(x_3')$, and hence $X_{x_3} =X_{x_3'}$.
\para

Suppose that $f:X \to Z$ is surjective. By definition of the equivalence relation, there is a well-defined bijective map $\bar{f}: X/_\sim \to Z$ given by $\bar{f}(X_x)=f(x)$. If $X_{x_1}, X_{x_2} \in X/_\sim$, then
$$\bar{f}(X_{x_1}\circ X_{x_2})=\bar{f}(X_{x_1* x_2})=f(x_1 * x_2)=f(x_1)*f(x_2)=\bar{f}(X_{x_1})*\bar{f}(X_{x_2}),$$
and hence $\bar{f}$ is an isomorphism of quandles. $\blacksquare$
\end{proof}

We know that a subgroup of a group is normal if and only if it is the kernel of some group homomorphism. In a similar way, we say that a subquandle $Y$ of a quandle $X$ is {\it normal} if $Y=X_{x_0}$ for some $x_0 \in X$ and some quandle homomorphism $f:X \to Z$. In this case, we say that $Y$ is {\it normal based at} $x_0$.
\para

A {\it pointed quandle}, denoted $(X, x_0)$, is a quandle $X$ together with a fixed base point $x_0$. Let $f:(X, x_0) \to (Z, z_0)$ be a homomorphism of pointed quandles, and $Y=X_{x_0}$ a normal subquandle based at $x_0$. In this case, we consider $Y$ as the base point of $X/_\sim$, and denote $X/_\sim$ by $X/Y$. Then the natural map
$$\pi: (X, x_0) \to (X/Y, X_{x_0})$$ 
given by $\pi(x) = X_x$ is a surjective homomorphism of pointed quandles.  This further induces a surjective ring homomorphism $$\hat{\pi}: \mathbb{k}[X] \to \mathbb{k}[X/Y]$$ with $\ker(\hat{\pi})$ being a two-sided ideal of $\mathbb{k}[X]$.

\para
Let $(X, x_0)$ be a pointed quandle, $\mathcal{I}$ the set of two-sided ideals of $\mathbb{k}[X]$ and $\mathcal{S}$ the set of normal subquandles of $X$ based $x_0$. Then there exist maps $\Phi: \mathcal{I} \to \mathcal{S}$ given by
$$\Phi(I)=X_{I,x_0}$$
and $\Psi: \mathcal{S} \to \mathcal{I}$ given by
$$\Psi(Y)=\ker(\hat{\pi}).$$

With this set up, we have the following result \cite[Theorem 4.10]{MR3977818}.

\begin{theorem}\label{dictionary}
Let $(X, x_0)$ be a pointed quandle and $\mathbb{k}$ an associative ring. Then $\Phi \Psi= \id_{\mathcal{S}}$ and $\Psi \Phi \neq \id_{\mathcal{I}}$.
\end{theorem}

\begin{proof}
Let $Y=X_{x_0}$ be a normal subquandle of $X$, that is,  $Y=\{x \in X \, \mid \, f(x)=f(x_0) \}$. Then we have
$$ 
\Phi \Psi(Y)  =  \Phi \big( \ker(\hat{\pi}) \big)= \big\{x \in X  \mid e_{x}-e_{x_0} \in  \ker(\hat{\pi}) \big\}= \big\{x \in X  \mid  X_x=X_{x_0}\big\}=
\big\{x \in X \mid  f(x)=f(x_0)\big\}=Y.
$$
Obviously, $\Psi \Phi \neq \id_{\mathcal{I}}$, since $\Psi \Phi \big(\mathbb{k}[X] \big)=\Psi (X)= \Delta_{\mathbb{k}}(X) \neq \mathbb{k}[X]$. $\blacksquare$
\end{proof}
\bigskip

We define the \index{center of quandle ring}{\it center} $\Z \big(\mathbb{k}[X] \big)$ of the quandle ring $\mathbb{k}[X]$ as
$$
\Z \big(\mathbb{k}[X] \big) = \big\{ u \in \mathbb{k}[X]  \,\mid \, u v = v u~\mbox{for all}~v \in \mathbb{k}[X] \big\}.
$$
Note that $\Z \big(\mathbb{k}[X] \big)$ is clearly an additive subgroup of the quandle ring $\mathbb{k}[X]$. Since $\mathbb{k}[X]$, in general, is non-associative, it follows that $\Z \big(\mathbb{k}[X] \big)$ need not be a subring of  $\mathbb{k}[X]$. Obviously, if $X$ is a trivial quandle, then $\Z \big(\mathbb{k}[X] \big)$ is a subring of $\mathbb{k}[X]$. However, even in this case, $\Z \big(\mathbb{k}[X] \big)$ need not be an ideal of $\mathbb{k}[X]$.
\para 
Recall that, a quandle is called {\it latin} if the left multiplication by each element is a permutation of the quandle. For example, odd order dihedral quandles are latin. Note that a latin quandle is always connected.

\begin{proposition}
Let $X = \{ x_1, \ldots, x_n \}$ be a finite latin quandle and $\mathbb{k}$ an associative ring. Then the following  assertions hold:
\begin{enumerate}
\item The element $w = e_{x_1} + \cdots + e_{x_n}$ lies in the center $\Z \big(\mathbb{k}[X] \big)$.
\item If $\mathbb{k}$ is a field with characteristic not dividing $n$, then $\frac{1}{n}w$ is an idempotent.
\end{enumerate}
\end{proposition}

\begin{proof}
For each $u \in \mathbb{k}[X]$, let $\hat{S}: \mathbb{k}[X] \to \mathbb{k}[X]$ be the map given by
$\hat{S}_u (w)= wu$ for all $w \in  \mathbb{k}[X]$. In particular, $\hat{S}_{e_x}$ is an automorphism of $\mathbb{k}[X]$ for each $x \in X$ since it is induced by the inner automorphism $S_x:X \to X$. For each $e_{x_i}$, we have
$$
w  e_{x_i} = \hat{S}_{e_{x_i}}(e_{x_1})+ \cdots + \hat{S}_{e_{x_i}}(e_{x_n}) = w,
$$
since $\hat{S}_{e_{x_i}}$ is an automorphism of $X$ and acts as a permutation. On the other hand, since the quandle $X$ is latin, we have $e_{x_i}  w = w$, which proves (1). If the characteristic of $\mathbb{k}$ does not divide $n$, then a direct computation proves assertion (2). 
\end{proof}
\bigskip
\bigskip


\section{Associated graded rings of quandle rings}\label{section6}
Let $X$ be a quandle (respectively, rack) and $\mathbb{k}$ an associative ring with unity. Consider the direct sum
$$
\mathcal{X}_{\mathbb{k}}(X) = \bigoplus_{i\geq 0} \Delta^i_{\mathbb{k}}(X) / \Delta^{i+1}_{\mathbb{k}}(X)
$$
of $\mathbb{k}$-modules $\Delta^i_{\mathbb{k}}(X) / \Delta^{i+1}_{\mathbb{k}}(X)$. We regard $\mathcal{X}_{\mathbb{k}}(X)$ as a graded $\mathbb{k}$-module with the convention that the elements of $\Delta^i_{\mathbb{k}}(X) / \Delta^{i+1}_{\mathbb{k}}(X)$ are homogeneous of degree $i$. Define multiplication in $\mathcal{X}_{\mathbb{k}}(X)$ as follows. For $u_i \in \Delta^i_{\mathbb{k}}(X)$, let $\underline{u}_i = u_i + \Delta^{i+1}_{\mathbb{k}}(X)$. Then, for $\underline{u}_i \in \Delta^i_{\mathbb{k}}(X) / \Delta^{i+1}_{\mathbb{k}}(X)$, $\underline{u}_j \in \Delta^j_{\mathbb{k}}(X) / \Delta^{j+1}_{\mathbb{k}}(X)$, we define 
$$\underline{u}_i ~ \underline{u}_j = \underline{u_i u_j}.$$ The product of two arbitrary elements of $\mathcal{X}_{\mathbb{k}}(X)$ is defined by extending the above product linearly. In other words, if $\underline{u} = \sum_i \underline{u}_i$ and $\underline{v} = \sum_j \underline{v}_j$ are decompositions of $\underline{u}, \underline{v} \in \mathcal{X}_{\mathbb{k}}(X)$ into homogeneous components, then
$$
\underline{u}~ \underline{v} = \sum_{i,j} \underline{u}_i ~ \underline{v}_j.
$$
With this multiplication, $\mathcal{X}_{\mathbb{k}}(X)$ becomes a graded ring, called the \index{associated graded ring}{\it associated graded ring} of $\mathbb{k}[X]$. We have a complete description of associated graded rings of trivial quandles.

\begin{proposition}
If $\T$ is a trivial quandle and $x_0 \in \T$, then $\mathcal{X}_{\mathbb{k}}(\T) = \mathbb{k} x_0\oplus \Delta_{\mathbb{k}}(\T)$.
\end{proposition}

\begin{proof}
It follows from the definition of $\mathcal{X}_{\mathbb{k}}(\T)$ and Theorem \ref{deltasqzero}. $\blacksquare$
\end{proof}

In order to understand $\mathcal{X}_{\mathbb{k}}(X)$,  it is essential to examine the quotients $\Delta^i_{\mathbb{k}}(X) / \Delta^{i+1}_{\mathbb{k}}(X)$. In what follows,  we compute  powers of the integral augmentation ideals of dihedral quandles. Note that $\R_1$ and $\R_2$ are trivial quandles. We first consider dihedral quandles of odd orders \cite[Theorem 6.2]{MR3915329}. 

\begin{theorem} 
Let $n>1$ be an odd integer. Then $\Delta^k(\R_n) / \Delta^{k+1}(\R_n) \cong \mathbb{Z}_n$ for all $k \geq 1$. 
\end{theorem}

\begin{proof} 
We write $\R_n= \{0, 1, 2, \ldots, n-1\}$ and prove the result by induction on $k$. Throughout the proof, all the computations are modulo $\Delta^{2}(\R_n)$. Set $f_i:= e_{i}-e_{0}$ for $1 \le i \le n-1 $. First, we claim that $\Delta(\R_n)/\Delta^2(\R_n)$ is generated as an abelian group by the coset of $f_1$ and that $\Delta(\R_n) / \Delta^2(\R_n) \cong \mathbb{Z}_n$. Recall that $\Delta(\R_n)$ is generated by $\{ f_1, f_2, \ldots, f_{n-1} \}$ and $\Delta^2(\R_n)$ is generated by  $\{f_i  f_j \, \mid \, 1 \le i, j \le n-1\}$ as abelian groups. We see that
$$f_{2i}+f_{n-2i}= e_{2i}- e_0 +e_{n-2i}- e_0= -(e_{2i}-e_0)(e_{i}-e_0)=-f_{2i}f_{i},$$
for each $1\leq i\leq \frac{n-1}{2}$. This implies that $f_{2i}=-f_{n-2i}$ for each $1\leq i\leq \frac{n-1}{2}$. In particular, for  $i=\frac{n-1}{2}$, we have $f_1=-f_{n-1}$. Note that, $f_1 f_1=0$ implies $f_2=2f_1$, whereas $f_{n-1} f_1=0$ implies $f_3=3f_1$.  Continuing in this manner, we get $f_{n-1}=(n-1)f_1$. Since $f_1=-f_{n-1}$, we obtain $n f_1=0$. Hence, $\Delta(\R_n)/\Delta^2(\R_n)$ is generated as an abelian group by the coset of $f_1$, and $\Delta(\R_n) / \Delta^2(\R_n) \cong \mathbb{Z}_n$. 
\para
Now, assume that $\Delta^k(\R_n) / \Delta^{k+1}(\R_n) \cong \mathbb{Z}_n$ for some $k \ge 1$.  Let the coset of $\alpha \in  \Delta^k(\R_n)$ generate $\Delta^k(\R_n) / \Delta^{k+1}(\R_n)$. Obrserve that the cosets of the elements $f_1\alpha, \ldots, f_{n-1}\alpha$~ generate $\Delta^{k+1}(\R_n) / \Delta^{k+2}(\R_n)$. Since $f_k=k f_1 \mod \Delta^2(\R_n)$ for each $1\leq k\leq n-1$, it follows that $f_k \alpha=k f_1 \alpha \mod \Delta^{k+2}(\R_n)$ for each $1\leq k\leq n-1$. Thus, the coset of $f_1\alpha$ generates $\Delta^{k+1}(\R_n) / \Delta^{k+2}(\R_n)$. Since $n f_1=0 \mod \Delta^2(\R_n)$, it follows that $n f_1 \alpha =0 \mod \Delta^{k+2}(\R_n)$, and hence
$$\Delta^{k+1}(\R_n) / \Delta^{k+2}(\R_n) \cong \mathbb{Z}_n,$$ 
which is desired. $\blacksquare$
\end{proof}

Next, we calculate powers of the augmentation ideal for $\R_4$. Setting $f_i: = e_i - e_0$ for $1 \le i \le 3$, we see that $\Delta(\R_4)$ is generated as an abelian groups by $\{f_1, f_2, f_3 \}$.

\begin{proposition}\label{ass-graded-2}
The following statements hold for $\R_4$:
\begin{enumerate}
\item $\Delta^2(\R_4)$ is generated as an abelian group by $\{f_1 - f_2 - f_3,~ 2 f_2\}$
and $\Delta(\R_4) / \Delta^2(\R_4) \cong \mathbb{Z} \oplus \mathbb{Z}_2$.
\item If $k > 2$, then $\Delta^k(\R_4)$ is generated as an abelian group by $\{2^{k-2}(f_1 - f_2 - f_3), ~2^{k-1} f_2\}$
and $\Delta^{k-1}(\R_4) / \Delta^k(\R_4) \cong \mathbb{Z}_2 \oplus \mathbb{Z}_2$.
\end{enumerate}
\end{proposition}

\begin{proof}
The abelian group $\Delta^2(\R_4)$ is generated by the products $f_i f_j$. The following table shows these products $f_i f_j$, where $f_i$ is taken from the leftmost column and $f_j$ is taken from the topmost row.
\begin{center}
\begin{tabular}{|c||c|c|c|}
  \hline
 $\cdot$ & $f_1$ & $f_2$ & $f_3$ \\
   \hline
     \hline
$f_1$ & $f_1 - f_2 -f_3$ & $0$ & $f_1 - f_2 -f_3$ \\
  \hline
$f_2$  & $-2 f_2$ & $0$ & $-2 f_2$ \\
  \hline
$f_3$ & $-f_1 - f_2+f_3$ & $0$ & $-f_1 - f_2+f_3$ \\
  \hline
\end{tabular}\\
\end{center}
It follows that $\Delta^2(\R_4)$ is generated by the set
$$
\big\{f_1 - f_2 -f_3,~~ -2 f_2,~~ -f_1 - f_2+f_3 \big\}
$$
and $\{f_1 - f_2 - f_3, 2 f_2\}$ forms a basis of $\Delta^2(\R_4)$. The isomorphism $\Delta(\R_4) / \Delta^2(\R_4)\cong \mathbb{Z} \oplus \mathbb{Z}_2$ follows from properties of abelian groups, which proves assertion (1).
\para
Next, we consider assertion (2). In order to find a basis of $\Delta^3(\R_4)$, we multiply the basis of $\Delta^2(\R_4)$ by the elements $f_1, f_2, f_3$. Using the preceding multiplication table, we get
\begin{align*}
&(f_1 - f_2 - f_3) f_1 = 2 (f_1 + f_2 - f_3), & & f_1 (f_1 - f_2 - f_3) = 0, & &(f_1 - f_2 - f_3) f_2 = f_2 (f_1 - f_2 - f_3) = 0,\\
&(f_1 - f_2 - f_3) f_3 = 2 (f_1 + f_2 - f_3), & &f_3 (f_1 - f_2 - f_3) = 0,& & 2 f_1  f_2 = 0, \\
&2 f_2 f_1 = - 4 f_2, & &  2 f_2 f_2 = 0, & & 2 f_2 f_3 = - 4 f_2,\\
& 2 f_3 f_2 = 0. & & &
\end{align*}
Hence, the set $\{2 (f_1 - f_2 - f_3), 4 f_2\}$ forms a basis of $\Delta^3(\R_4)$. We obtain the general formulas using induction on $k$. Finally, the  assertion $\Delta^{k-1}(\R_4) / \Delta^k(\R_4) \cong \mathbb{Z}_2 \oplus \mathbb{Z}_2$ follows from properties of abelian groups. $\blacksquare$
\end{proof}

Let $\R_n$ be the dihedral quandle of order $n$ and $f_i: = e_i - e_0$ for $1 \le i \le n-1$. A direct check gives
\begin{equation}\label{neranga2}
f_i  f_j= f_{2j-i}-f_{n-i}-f_{2j}
\end{equation}
for all $1 \le i, j \le n-1$. Next, we have the following result \cite[Theorem 6.3]{MR3915329}.

\begin{theorem}\label{neranga3}
If $n=2k$ for some integer $k \ge 2$, then $\Delta(\R_n) / \Delta^2(\R_n)\cong \mathbb{Z}\oplus \mathbb{Z}_{k}$. 
\end{theorem}

\begin{proof}
Throughout the proof, all the computations are modulo $\Delta^{2}(\R_n)$.  We claim that the following assertions hold: 
\begin{enumerate}
\item  $\Delta(\R_n) / \Delta^{2}(\R_n)$ is generated as an abelian group by the cosets of $f_1$ and $f_2$.
\item   The abelian subgroup generated by the coset of $f_2$ is isomorphic to $\mathbb{Z}_k$. 
 \item  The abelian subgroup generated by the coset of $f_1$ is isomorphic to $\mathbb{Z}$. 
 \end{enumerate}

Using \eqref{neranga2}, we  have
$f_{n-l} f_{1}=-f_{2}-f_{l}+f_{l+2}$. Since $f_{n-l} f_{1}=0$, this gives 
\begin{equation}\label{relation among fis}
f_{l+2}=f_l+f_2.
\end{equation}
We  claim that
\begin{equation}\label{Star relations}
f_l = \left\{
\begin{array}{ll}
\frac{l}{2}\,f_2 &  \textnormal{if}\,\,l\,\,\textnormal{is even},\\[0.3cm]
\lfloor \frac{l}{2} \rfloor\,f_2+f_1 &  \textnormal{if}\,\,l\,\,\textnormal{is odd}.
\end{array}
\right.
\end{equation}
 for each $1 \le l \le  n-1$, where $\lfloor - \rfloor$ is the floor function. The relations are obvious for $l=1, 2$.  We proceed by induction on $l$. First, assume that $l \ge 4$ is even. Since $0=f_{n-2}  f_1=-2f_2+f_4$, we obtain $f_4=2f_2$, which is the case $l=4$. Suppose that \eqref{Star relations} holds for $l \ge 4$ even. Using the induction hypothesis and \eqref{relation among fis}, we get
$$f_{l+2}=f_l+f_2=\frac{l}{2}\,f_2+f_2=\Big(\frac{l+2}{2}\Big)\,f_2,$$
which is desired. Next, we assume that $l \ge 3$ is odd. The division algorithm gives $l=\lfloor \frac{l}{2} \rfloor\,2+1$, and hence
we can write $f_{l-1}=\frac{(l-1)}{2} f_2=\lfloor \frac{l}{2} \rfloor\,f_2$.  It follows from \eqref{neranga2} that 
\begin{equation}\label{second relations amongs fis}
f_l=f_{l-1}+f_1.
\end{equation}
Since $0=f_{n-1}  f_1=-f_1-f_2+f_3$, we obtain $f_3=f_2+f_1$, which is the case $l=3$. Suppose that \eqref{Star relations} holds for $l \ge 3$ odd. Using induction hypothesis, \eqref{relation among fis} and \eqref{second relations amongs fis}, we get
$$f_{l+2}=f_{l-1}+f_1+f_2=\Big(\frac{l-1}{2}\Big)\,f_2+f_1+f_2=\Big(\frac{l+1}{2}\Big)\,f_2+f_1=\lfloor \frac{l+2}{2} \rfloor\,f_2+f_1.$$
This proves \eqref{Star relations}, from which assertion (1) follows immediately.
\para

Next, if $i$ is even, then \eqref{neranga2} and \eqref{Star relations} give
$$
-f_i  f_j=-f_{2j-i}+f_{n-i}+f_{2j}=-\Big(\frac{2j-i}{2}\Big)\,f_2+\Big(\frac{n-i}{2}\Big)\,f_2+\Big(\frac{2j}{2}\Big)\,f_2 =\frac{n}{2}\,f_2 =kf_2.
$$
Similarly, if $i$ is odd, then $-f_i f_j=kf_2$. Since $f_i f_j=0$, it follows that $kf_2=0$. Further, we  have $0 \ne f_4 = 2f_2, \ldots, 0 \neq f_{2k-2} =(k-1)f_2$,
it follows that the abelian subgroup generated by the coset of $f_2$ is isomorphic to $\mathbb{Z}_k$, which is assertion (2).
\para

The relation \eqref{Star relations} implies that $f_1 =f_3-f_2$. Thus, if $mf_1 =0$, then $mf_2 =mf_3$, which implies that $m=0$. Hence, the abelian group generated by the coset of $f_1$ is torsion-free, which is assertion (3). Thus, $\Delta(\R_n) / \Delta^2(\R_n)\cong \mathbb{Z}\oplus \mathbb{Z}_{k},$  which is desired. $\blacksquare$
\end{proof}
\bigskip
\bigskip


\section{Power-associativity of quandle rings}\label{section7}

Following Albert \cite{MR0026044}, we have the following general form of associativity for rings.

\begin{definition}
A ring $A$ is called \index{power-associative ring}{\it power-associative}  if every element of $A$ generates an associative subring of $A$.
\end{definition}

Every associative ring is power-associative. In particular, if $\T$ is a trivial quandle and $\mathbb{k}$ an associative ring, then $\mathbb{k}[\T]$ is power-associative.  It is well-known \cite{MR0026044} that a ring $A$ of characteristic zero is power-associative if and only if 
$$(x  x) x=x (x  x) \quad \textit{and} \quad (x  x) (x  x)=\big((x  x) x\big)  x$$
for all $x \in A$. The following result shows that  quandle rings of non-trivial quandles are not power-associative \cite[Theorem 3.4]{MR3915329}.

\begin{theorem}\label{power ass quandle algebra}
Let $X$ be a non-trivial quandle and $\mathbb{k}$ an associative ring such that $\charac(\mathbb{k})\neq 2, 3$.  Then the quandle ring $\mathbb{k}[X]$ is not power-associative.
\end{theorem}

\begin{proof}
We claim that the following assertions hold:  
\begin{enumerate}
\item  If there exist $x, y\in X$ with $x\neq y$ such that $x* y = y* x$, then $\mathbb{k}[X]$ is not power-associative.
\item  Suppose that for all $x, y\in X$ with $x\neq y$, we have  $x* y \neq y* x$. If $\mathbb{k}[X]$ is power-associative, then $x* y= x$.
\end{enumerate}

Suppose that there exist elements $x,y \in X$ with $x\neq y$ such that $x* y = y* x$. Assume that $\mathbb{k}[X]$ is power-associative. If we take $u=a e_{x}+be_{y}$, then we must have $(u u) (u u)=\big((u u) u\big) u$ for all choices of $a, b \in \mathbb{k}$. Direct computations give
$$u u=a^2 e_{x}+2ab e_{x* y}+ b^2 e_{y},$$
\begin{equation}\label{NK1}
(u u) (u u)=a^4 e_{x} +4 a^3be_{x* (x* y)}+4ab^3 e_{y*(x* y)}+6a^2b^2 e_{x* y}  +b^4 e_y
\end{equation}
and 
\begin{equation}\label{NK2}
\begin{split}
\big((u u) u \big) u&=a^4 e_{x} +b^4 e_{y} + a^3b (e_{x* y} + e_{(x* y)* x} + 2e_{((x* y)* x)* x} )\cr
&+ ab^3 (e_{y* x} + e_{(y* x)* y} + 2e_{((x* y)* y)* y})\cr
&+ a^2b^2 (e_{(x* y)* y} + 2e_{((x* y)* x)* y} + 2e_{((x* y)* y)* x} + e_{(y* x)* x)}).
\end{split}
\end{equation}
We now analyse the identity $(u u) (u u)=\big((u u) u \big) u$. Since $\charac(\mathbb{k})\neq 2, 3$, it is easy to check (by plugging in different values for $a$ and $b$) that the coefficients (which are the elements of $\mathbb{k}[X]$) of $a^3b$, $ab^3$ and $a^2b^2$ on both the sides of the relation $(u u) (u u)= \big((u u) u \big) u$ must be pairwise equal. In particular, by comparing the coefficients of $a^2b^2$ in \eqref{NK1} and \eqref{NK2}, we get $x* y =(x* y)* y$, which implies $x = x* y$.  Since $x* y = y* x$, it follows that $x=y$, which contradicts our assumption that $x\neq y$. Thus, $\mathbb{k}[X]$ is not power-associative, which is assertion (1). 
\para 
Let $x, y\in X$ with $x\neq y$ such that $x* y \neq y* x$. Suppose that $\mathbb{k}[X]$ is power-associative.  Let $u=a e_{x}+ b e_{y}$ such that $ab \ne 0$. Then, direct computations give 
$$u u=a^2 e_{x}+ab e_{x* y}+ ab e_{y* x}+b^2 e_{y},$$
\begin{equation}\label{NK3}
(u u) u = a^3 e_{x}+ b^3 e_{y} + a^2b (e_{x* y} + e_{(x* y)* x} + e_{(y* x)* x} ) 
+ ab^2( e_{y* x} + e_{(x* y)* y} + e_{(y* x)* y}) 
\end{equation}
and 
\begin{equation}\label{NK4}
u (u u) = a^3 e_{x}+ b^3 e_{y} + a^2b (e_{y* x} + e_{x* (x* y)} + e_{x* (y* x)})
+ ab^2 (e_{x* y} + e_{y* (x* y)} + e_{y* (y* x)}).
\end{equation}
Due to power-associativity, we must have $(u u) u = u (u u)$. Since $\charac(\mathbb{k})\neq 2$, by equating \eqref{NK3} and \eqref{NK4}, and letting $a=b=1$ and $a=-1,b=1$, we get 
$$e_{x* y} + e_{(x* y)* x} + e_{(y* x)* x} = e_{y* x} + e_{x* (x* y)} + e_{x* (y* x)}.$$
Since $(x* y)* x=x* (y* x)$ in $X$, the preceding equation becomes 
$$e_{x* y} + e_{(y* x)* x} = e_{y* x} + e_{x* (x* y)}.$$
Since $x* y \neq y* x$, we get $x* y=x* (x* y)$, which implies that $x* y = x$, which is assertion (2).  But, this implies that $X$ is trivial, a contradiction.
$\blacksquare$
\end{proof}		
\para

Next, we examine the connections between quandle algebras and other algebras. Group algebras are associative, allowing the use of tools from associative algebra in their study. In contrast, quandle algebras of non-trivial quandles are non-associative. Various types of non-associative algebras, such as alternative algebras, Jordan algebras, and Lie algebras, have been extensively studied. This raises the natural question of whether quandle algebras fall into any of these well-known classes.
\para

We recall definitions of some well-known algebras \cite{MR0026044, MR0668355}. Let $A$ be an algebra over a commutative and associative ring $\mathbb{k}$ with unity. Then, we say that $A$ is:
\begin{enumerate}
\item  an  \index{alternative algebra}{\it alternative algebra} if
$$
a^2 b = a (a b) ~\mbox{and}~ a b^2 = (a b) b~\mbox{for all}~a, b \in A.
$$
\item a \index{Jordan algebra} {\it Jordan algebra} if 
$$
ab = ba ~\mbox{and}~(a^2 b) a = a^2 (b a)~\mbox{for all}~a, b \in A.
$$
\item a \index{Lie algebra} {\it Lie algebra} if
$$
a^2 = 0 ~\mbox{and}~ (a b) c + (b c) a + (c a) b = 0~\mbox{for all}~a, b, c \in A.
$$
\item an \index{elastic algebra} {\it elastic algebra} if 
$$
(x y) x = x (y x)~\textrm{for all}~ x, y \in A.
$$
\item  a {\it power-associative algebra} if every element of $A$ generates an associative subalgebra of $A$.
\end{enumerate}
\para

For example, any commutative or associative algebra is elastic. In particular quandle algebras of trivial quandles are elastic being associative. Also, by Theorem \ref{power ass quandle algebra},  quandle algebras of non-trivial quandles over rings of characteristic other than 2 and 3 are not power-associative.

\begin{proposition}\label{elasticity-quandle-ring}
Let $X$ be a non-trivial quandle. If $\mathbb{k}$ is a ring of characteristic other than 2 and 3, then $\mathbb{k}[X]$ cannot be a Jordan algebra, an alternative algebra, or an elastic algebra.
\end{proposition}

\begin{proof}
By Theorem \ref{power ass quandle algebra}, the quandle algebra $\mathbb{k}(X)$ is not power-associative. On the other hand,  any Jordan algebra is power-associative \cite[Chapter 2]{MR0668355}. Also, by Artin's theorem any alternative algebra is also power-associative  \cite[Chapter 2]{MR0668355}. Further, if $\mathbb{k}[X]$ is elastic, then $(xx)x=x(xx)$ and $(xx)(xx)= \big((xx)x\big)x$ for all $x \in \mathbb{k}[X]$, which implies that $\mathbb{k}[X]$ is power-associative, a contradiction. $\blacksquare$
\end{proof}

\begin{remark}{\rm 
If $X$ is a quandle, then the third quandle axiom implies that $(x y) x = x (y x)$ for all $x, y \in X$. Thus, all quandles satisfy the elasticity condition which, by Proposition \ref{elasticity-quandle-ring}, is in contrast to their quandle algebras.}
\end{remark}

There are two natural constructions on any algebra $A = \langle A; +, \cdot \rangle$. Define an algebra $A^{(-)} = \langle A; +, \circ \rangle$ with multiplication
$$
x \circ y = x y - y x.
$$
If $A$ is an associative algebra, then $A^{(-)}$ is a Lie algebra. If the ring $\mathbb{k}$ contains $\frac{1}{2}$, then we define $A^{(+)} = \langle A; +, \odot \rangle$ with multiplication given by
$$
x \odot y = \frac{1}{2} (x y + y x).
$$
If $A$ is an associative algebra, then $A^{(+)}$ is a Jordan algebra. We conclude with the following result.

\begin{theorem}
Let $\T = \{ x_1, x_2,x_3,  \ldots \}$ be a trivial quandle with more than one element, $A=\mathbb{k}[\T]$, $L = A^{(-)}$ the corresponding  Lie algebra and $J = A^{(+)}$  the corresponding Jordan algebra. If $L^k$ is the subalgebra of $L$ generated by products of $k$ elements of $L$, then the following  assertions hold:
\begin{enumerate}
\item $L^2 = L^3$ and this algebra has a basis
$$
\{e_{x_1} - e_{x_2}, e_{x_2} - e_{x_3}, \ldots \}.
$$
In particular, if $\T = \T_n$ contains $n$ elements, then $L$ has rank $n-1$.

\item $(L^2)^2 = \{0\}$, that is, $L$ is a metabelian algebra.

\item $J^2 = J$.
\end{enumerate}
\end{theorem}

\begin{proof} The algebra $L^2$ is generated by the products $e_{x_i} \circ e_{x_j} = e_{x_i} - e_{x_j}$, where $i < j$. Setting $f_i = e_{x_i} - e_{x_{i+1}}$, we claim that any product $e_{x_i} \circ e_{x_j}$ is a linear combination of $f_i$. Indeed, if $j = i+1$, then this product is $f_i$. If $j-i > 1$, then
$$
e_{x_i} - e_{x_j} = f_i + f_{i+1} + \cdots + f_{j-1}.
$$
A direct check shows that the elements $f_1, f_2,f_3, \ldots$ are linearly independent and
$$
f_i \circ e_{x_j} =  f_i \quad \textrm{and} \quad e_{x_j} \circ e_i = -f_i.
$$
Hence, we have $L^2 \subseteq L^3$, and assertion (1) follows.
\para

The algebra $(L^2)^2$ is generated by the products $f_i \circ f_j$. Direct calculation gives
$$
f_i \circ f_j = f_i f_j - f_j f_i = (e_{x_i} - e_{x_{i+1}}) (e_{x_j} - e_{x_{j+1}}) - (e_{x_j} - e_{x_{j+1}})  (e_{x_i} - e_{x_{i+1}}) = 0,
$$
which is assertion (2).
\para

For assertion (3), note that
$$
e_{x_i} \odot e_{x_j} = \frac{1}{2} (e_{x_i} e_{x_j} + e_{x_j} e_{x_i}) = \frac{1}{2} (e_{x_i} + e_{x_j}).
$$
It follows that $J^2$ contains the elements $e_{x_i} = e_{x_i} \odot e_{x_i}$, and hence $J^2 = J$. $\blacksquare$
\end{proof}
\bigskip
\bigskip


\section{Commutator width of quandle algebras}\label{commutator-width}
Let $X$ be a quandle and $\mathbb{k}$ an associative ring with unity. Define the \textit{commutator} of elements $u, v \in \mathbb{k}[X]$ as the element $$[u, v]=uv-vu.$$ Then the \index{commutator subalgebra}\textit{commutator subalgebra} $\mathbb{k}[X]'$ of $\mathbb{k}[X]$ is the $\mathbb{k}$-algebra generated by the set of all commutators in $\mathbb{k}[X]$. 

\begin{lemma}\label{commutator-augmentation}
Let $X$ be a quandle and $\mathbb{k}$ an associative ring with unity. Then the following assertions hold:
\begin{enumerate}
\item  If $X$ is a commutative quandle, then $\mathbb{k}[X]' =\{0\}$. 
\item If $\mathbb{k}$ is commutative, the $\mathbb{k}[X]' \le \Delta_{\mathbb{k}}(X)$.
\end{enumerate}
\end{lemma}

\begin{proof}
Assertion (1) is obvious. Since $\mathbb{k}$ is commutative and $\varepsilon \big([u, v] \big)=0$ for each $[u, v] \in \mathbb{k}[X]'$, assertion (2) follows.
\end{proof}

The equality in Lemma \ref{commutator-augmentation}(2) does not hold in general. For example, the dihedral quandle $\R_3$ is commutative, and hence $\mathbb{Z}[\R_3]'=\{0\}$. On the other hand, we have $\Delta(\R_3) =\langle e_1-e_0,  e_2-e_0 \rangle \neq \{0\}$. 
\para

We define the \index{commutator length} \textit{commutator length} $\cl(u)$ of an element $u \in \mathbb{k}[X]'$ as 
$$\cl(u)=\min \big\{n \,\mid\, u= \sum_{i=1}^n \alpha_i [u_i,v_i],~\textrm{where}~\alpha_i \in R ~\textrm{and}~u_i, v_i \in \mathbb{k}[X]\big\}.$$
The \textit{commutator width} $\cw \big(\mathbb{k}[X] \big)$ of $\mathbb{k}[X]$ is defined as 
$$\cw \big(\mathbb{k}[X] \big)=\sup \big\{\cl(u) \, \mid \, u \in \mathbb{k}[X]'\big\}.$$

It follows from the definition of commutator width that a quandle $X$ is commutative if and only if  $\cw \big(\mathbb{k}[X] \big)=0$. Consequently, we have $\cw \big(\mathbb{k}[\R_3] \big)=0$.
\para

The computation of commutator width has been studied in the context of free Lie rings \cite{MR1803583}, free metabelian Lie algebras \cite{MR3328893}, and absolutely free as well as free solvable Lie rings of finite rank \cite{MR3601334}. In the remainder of this section, we compute the commutator width for several quandle rings.
\para

We say that a quandle $X$ is {\it strongly non-commutative} if for every pair of distinct elements $x, y \in X$ there exist elements $a, b \in X$ such that $a*b = x$ and $b*a = y$. Obviously, every strongly non-commutative quandle is non-commutative.

\begin{theorem}\cite[Theorem 7.2]{MR4450681} \label{commutator-width-general}
Let $X$ be a strongly non-commutative quandle or a non-commutative quandle admitting a 2-transitive action by $\Aut(X)$. Let $\mathbb{k}$ be a commutative and associative ring with unity. Then the following  assertions hold:
\begin{enumerate}
\item $\mathbb{k}[X]'=\Delta_{\mathbb{k}}(X)$.
\item If $X$ has order $n$, then $1 \le \cw \big(\mathbb{k}[X] \big) \le n-1$.
\end{enumerate}
\end{theorem}

\begin{proof}
In view of Lemma \ref{commutator-augmentation}(2), we only need to show that $\Delta_{\mathbb{k}}(X) \le \mathbb{k}[X]'$. Since $\Delta_{\mathbb{k}}(X)$ is generated as an $\mathbb{k}$-module by elements of the form $e_x-e_y$, where $x, y \in X$ are distinct elements, it suffices to show that each such element is a commutator. Suppose first that $X$ is strongly non-commutative. Let $x, y \in X$ be two distinct elements. Then, there exist $a, b \in X$ such that $e_x-e_y = e_ae_b - e_be_a  = [e_a, e_b] \in \mathbb{k}[X]'$. Now, suppose that $X$ is non-commutative. Then, there exist elements $c, d \in X$ such that $c*d \neq d*c$. By 2-transitivity of $\Aut(X)$ action on $X$, there exists $\phi \in \Aut(X)$ such that $\phi(c*d)=x$ and  $\phi(d*c)=y$. Thus, we have $e_x-e_y=\hat{\phi}(e_c e_d)-\hat{\phi}(e_d e_c)=[\hat{\phi}(e_c), \hat{\phi}(e_d)] \in  \mathbb{k}[X]'$, where $\hat{\phi}$ is the induced automorphism of $\mathbb{k}[X]$. This proves assertion (1).
\para
Let $X=\{x_0, x_1, \ldots, x_{n-1} \}$. By assertion (1), we have $\mathbb{k}[X]'=\Delta_{\mathbb{k}}(X)$. Hence, any $u \in \mathbb{k}[X]'$ has the form $u= \sum_{i=1}^{n-1} (e_{x_i}-e_{x_0})$, where each $(e_{x_i}-e_{x_0})$ is a commutator as shown in the proof of (1). Since $X$ is non-commutative, we obtain $1 \le \cw \big(\mathbb{k}[X] \big) \le n-1$. $\blacksquare$
\end{proof}

\begin{corollary}
Let $G$ be an elementary abelian $p$-group with $p>3$ and $\phi \in \Aut(G)$ act as multiplication by a non-trivial unit of $\mathbb{Z}_p$. Then $1 \le \cw\big(\mathbb{k}[\Alex(G,\phi)]\big) \le |G|-1$.
\end{corollary}

\begin{proof}
The quandle $\Alex(G,\phi)$ is non-commutative for $p>3$. Further, by Theorem \ref{FNT-generalisation}, it admits a 2-transitive action by $\Aut \big(\Alex(G,\phi) \big)$.   $\blacksquare$
\end{proof}

We remark that a complete description of finite 2-transitive quandles has been given by Bonatto \cite{MR4075376} by extending results of Vendramin \cite{MR3685034}. The following result computes commutator width of some finite quandles \cite[Theorem 7.4]{MR4450681}.

\begin{theorem} 
Let $\mathbb{k}$ be a commutative and associative ring with unity. Then the following assertions hold:
\begin{enumerate}
\item If $\T$ is a trivial quandle, then $\cw \big(\mathbb{k}[\T] \big)=1$.
\item $\cw \big(\mathbb{k}[\R_4] \big) =1$.
\item $\cw \big(\mathbb{k}[\J_3] \big)=1$.
\end{enumerate}
\end{theorem}

\begin{proof}
By Theorem \ref{commutator-width-general}, we have $\mathbb{k}[\T]'=\Delta_{\mathbb{k}}(\T)$, and hence any element of $\mathbb{k}[\T]'$ has the form $u= \sum_{i=1}^n \alpha_i (e_{x_i}-e_{x_0})$. Taking $v= (1- \sum_{i=1}^n \alpha_i)e_{x_0} +  \sum_{i=1}^n \alpha_i e_{x_i}$ and $w=e_{x_0}$, we see that
\begin{eqnarray*}
[v, w] &=& vw-wv\\
&=& \varepsilon(w) v -\varepsilon(v) w\\
&=& (1- \sum_{i=1}^n \alpha_i)e_{x_0} +  \sum_{i=1}^n \alpha_i e_{x_i} -e_{x_0}\\
&=&  \sum_{i=1}^n \alpha_i (e_{x_i}-e_{x_0})=u,
\end{eqnarray*}
and hence $\cw \big(\mathbb{k}[\T] \big)=1$, which proves (1).
\para

By Lemma \ref{commutator-augmentation}(2), we have $\mathbb{k}[\R_4]' \le \Delta_{\mathbb{k}}(\R_4)$. Note that $\Delta_{\mathbb{k}}(\R_4)$ is generated by $\{f_1, f_2, f_3\}$, where $f_i=e_i-e_0$ for $i=1,2,3$. A direct check shows that
$$
f_1 = [e_3, e_2],~~f_2=[e_2, e_0], ~~f_3 = -[e_2, e_1],
$$ 
and hence $\mathbb{k}[\R_4]' = \Delta_{\mathbb{k}}(\R_4)$. 
If $w=\gamma_1 f_1+ \gamma_2 f_2 + \gamma_3 f_3 \in \Delta_{\mathbb{k}}(\R_4)$ for some $\gamma_i \in \mathbb{k}$, then we can write it as
$$
w=[e_2, \gamma_2 e_0 - \gamma_3 e_1 - \gamma_1 e_3].
$$ 
Thus, every element of $\mathbb{k}[\R_4]'$ is a commutator, and hence $\cw \big(\mathbb{k}[\R_4] \big) = 1$, which proves (2).
\para

Let $\J_3=\{1,2,3\}$ be the quandle with multiplication table given by
$$
\J_3=\begin{tabular}{|c||c|c|c|}
    \hline
$*$ & 1 & 2 & 3  \\
  \hline \hline
1 & 1 & 1 & 1  \\
\hline
2 & 3 & 2 & 2  \\
\hline
3 & 2 & 3 & 3  \\
  \hline
\end{tabular}.
$$
The augmentation ideal $\Delta_{\mathbb{k}}(\J_3)$ is generated by $f_1=e_3-e_2$ and $f_2=e_1-e_2$. Further, multiplication rules in $\J_3$ show that $f_1=e_3 e_2-e_2 e_3=[e_3, e_2]$ and $f_2=e_1 e_3-e_3 e_1= [e_1, e_3]$, and hence $\mathbb{k}[\J_3]'=\Delta_{\mathbb{k}}(\J_3)$.  Let $w=\gamma_1f_1+ \gamma_2f_2 \in \Delta_{\mathbb{k}}(\J_3)$, where $\gamma_i \in \mathbb{k}$. A direct check gives
$$f_1 e_2=f_1,~~f_2 e_2=f_2,~~e_2 f_1=0,~~e_2 f_2=f_1.$$
Now, taking $u=e_2+ (\gamma_1+\gamma_2)f_1 + \gamma_2 f_2$ and $v=e_2$, we see that
$$w=uv-vu=[u, v].$$
Hence, $\cw \big(\mathbb{k}[\J_3] \big)=1$, which establishes assertion (3). $\blacksquare$
\end{proof}

\begin{problem}
Determine commutator width of quandle algebras of dihedral quandles, free quandles and other interesting classes of quandles.
\end{problem}
\bigskip
\bigskip


\section{Quandle algebras as representation spaces}
Let $X=X_1 \sqcup\cdots\sqcup X_k$ be the decomposition of a finite quandle $X$ by its connected components. Each $x \in X$ gives rise to  maps $L_x: X\rightarrow X$ and $S_x: X\rightarrow X$ given by $L_x(y)=x*y$ and $S_x(y)=y*x$, respectively. Thus, we have a map  $$\psi: X\rightarrow\mathcal{M}(X)$$
given by $x \mapsto L_x$. Similarly, we have a map
$$\varphi: X\rightarrow \Sigma_{n}$$
given by $x \mapsto S_x$. Note that this map further restricts to  maps $\varphi_i: X\rightarrow \Sigma_{n_i}$, where $n_i=|X_i|$ and $n=\sum\limits_{i=1}^{k} n_i$. Let $H(X)$ denote the semigroup generated by $\psi(X)$. Note that $H(X)$ is not necessarily a group. On the other hand, we know that the semigroup generated by  $\varphi(X)$ is the group $\Inn(X)$ of inner automorphisms of $X$. For each $i$, let $G(X_i)$ denote the group generated by $\varphi_i (X)$, which is a subgroup of $\Sigma_{n_i}$.
\para

Let $X$ be a quandle and $\mathbb{k}$ a field. Let $V$ denote the underlying $\mathbb{k}$-vector space of the quandle algebra $\mathbb{k}[X]$. To study right ideals of the algebra $\mathbb{k}[X]$, we view $V$ as a representation space of $\Inn(X)$. Clearly, we have 
$$V=\underset{i=1}{\overset{k}{\oplus}}V^i,$$
 where $V^i= \mathbb{k}[X_i]$. Thus, it is enough to  understand  decomposition of  $V^i$ into irreducible representations of $G(X_i)$. Since $X_i$ is finite, $V^i$ contains a one-dimensional subspace $V_{\textrm{triv}}^i:=\mathbb{k}v_{\textrm{triv}}$ which is invariant under the action of $G(X_i)$, where $v_{\textrm{triv}}=\sum\limits_{x\in X_i} e_x$.

\begin{definition}
Let $X=X_1 \sqcup\cdots\sqcup X_k$ be a finite quandle. Then $X$ is said to be:
\begin{enumerate}
\item  \textit{left $2$-transitive} if the semigroup  $H(X)$ acts $2$-transitively on $X$.
\item  \textit{right $2$-transitive} if the group $\Inn(X)$ acts $2$-transitively on $X$. 
\item  \textit{right orbit $2$-transitive} if  the group $G(X_i)$ acts $2$-transitively on $X_i$ for each $i$.
\end{enumerate}
\end{definition}

\para
Let $X$ be a finite quandle of order $n$ and  $V_{\textrm{st}}$ the subspace of $V$ orthogonal to $V_{\textrm{triv}}$.  The following result gives the list of subgroups $G\leqslant  S_n$ for which the representation $V_{\textrm{st}}$ is irreducible. 
\para

\begin{theorem} \label{TransClassif}
Let $X$ be a finite quandle of order $n$ and $\mathbb{k}$ a field. Then the following statements hold:
\begin{enumerate}
\item If $\charac(\mathbb{k})=0$,  then $V_{\textrm{st}}$ is an irreducible representation of the subgroup $G < S_n$ if and only if  $G$ is $2$-transitive and $A_n\not\leq G$. 
 \item If $\charac(\mathbb{k})=p>3$,  then $V_{\textrm{st}}$ is an irreducible representation of the subgroup $G < S_n$ if and only if  $G$ is $2$-transitive and $A_n\not\leq G$ except when:
 \begin{enumerate}
 \item[(i)] $G \leqslant  A\Gamma L(m, q)$ and $p$ divides $q$;
 \item[(ii)] $G \leqslant  P\Gamma L(m, q)$, $m\geq 3$ and $p$ divides $q$;
 \item[(iii)] $G \leqslant  P\Gamma U(3, q)$ and $p$ divides $q+1$;
 \item[(iv)] $G \leqslant  Sz(q)$ and $p$ divides $q+1+m$, where $m^2=2q$;
 \item[(v)] $G \leqslant  Re(q)$ and $p$ divides $(q+1)(q+1+m)$, where $m^2=3q$.
\end{enumerate} 
\end{enumerate}
\end{theorem}

\begin{proof}
Assertion (1) is stated in \cite[Theorem 1(b)] {MR0913206}, while assertion (2) is established in \cite[Main Theorem (ii)]{MR1807270}. The description of groups $(i)-(v)$ can be found in \cite[Chapter 7]{MR1409812}. In short, the groups $(i)-(v)$ are the affine general semilinear groups, projective general semilinear groups, projective semilinear unitary groups, Suzuki groups and Ree groups, respectively. $\blacksquare$
\end{proof}

If $A$ is a finite semigroup and $e \in A$ is an idempotent, then $eAe$ is a monoid with identity $e$. The group of units of $eAe$ is called the \textit{maximal subgroup} of $A$ at $e$.

\begin{corollary}\label{RightModDecomp}
Let $X$ be a finite quandle of order $n$ such that $A_n \not\le \Inn(X)$ and $\mathbb{k}$ a field. Then the following statements hold:
\begin{enumerate}
\item  If $X=X_1 \sqcup \cdots \sqcup X_k$ is a right orbit 2-transitive quandle and $\charac(\mathbb{k})=0$, then  $\mathbb{k}[X]\cong \underset{i=1}{\overset{k}{\bigoplus}}(V^i_{\textrm{st}}\oplus V^i_{\textrm{triv}})$, where the summands on the right hand side are simple right ideals.

\item  If $X=X_1 \sqcup \cdots \sqcup X_k$ is a right orbit 2-transitive quandle such that none of the groups $G(X_i)$ is in the list of groups in Theorem~\ref{TransClassif}(2) and  $\charac(\mathbb{k})=p>3$, then $\mathbb{k}[X]\cong \underset{i=1}{\overset{k}{\bigoplus}}(V^i_{\textrm{st}}\oplus V^i_{\textrm{triv}})$, where the summands on the right hand side are simple right ideals.

\item If $X$ is such that  $H(X)$ contains a maximal subgroup with a $2$-transitive action on $X$ and $\charac(\mathbb{k})=0$, then $\mathbb{k}[X]\cong V_{\textrm{st}}\oplus V_{\textrm{triv}}$, where the summands on the right hand side are simple left ideals.

\item If $X$ is such that  $H(X)$ contains a maximal subgroup with a $2$-transitive action on $X$, which is not in the list of groups in  Theorem~\ref{TransClassif}(2) and $\charac(\mathbb{k})=p>3$, then $\mathbb{k}[X]\cong V_{\textrm{st}}\oplus V_{\textrm{trivial}}$, where the summands on the right hand side are simple left ideals.
\end{enumerate}
\end{corollary}

\begin{definition}
A finite quandle $X$ is of \index{cyclic type quandle} \textit{right cyclic type} if for each $x \in X$, the inner automorphism $S_x$ acts on $X \setminus \{x\}$ as a cycle of length $|X|-1$.
\para
 Similarly, $X$ is said to be of \textit{left cyclic type} if for each $x \in X$, the semigroup element $L_x$ acts on $X \setminus \{x\}$ as a cycle of length $|X| - 1$.
\end{definition}
\para

The following result relates  2-transitivity of finite quandles to cyclicity \cite[Corollary 4]{MR3685034}.

\begin{theorem}
Every finite right  $2$-transitive quandle is of right cyclic type.
\end{theorem}

\begin{remark}
{\rm 
If $X$ is a quandle of left cyclic type, then it is left $2$-transitive.  
For, let $(a, b),(c, d) \in X \times X$. Since $X$ is of left cyclic type, there exist $i, j$ such that $L_b^i(a)= c$ and   $L_c^j(b) = d$. Consequently, $L_c^j \,L_b^i$ maps the pair $(a, b)$ to the pair $(c, d)$.
}
\end{remark}

The preceding observations allow to strengthen Corollary~\ref{RightModDecomp}.

\begin{corollary}
Let $X$ be a finite quandle of order $n$ such that $A_n \not\le \Inn(X)$ and $\mathbb{k}$ a field. Then the following statements hold:
\begin{enumerate}
\item  If $X=X_1 \sqcup\cdots\sqcup X_k$ is such that each subquandle $X_i$ is of right cyclic type, then $\mathbb{k}[X]\cong \underset{i=1}{\overset{k}{\bigoplus}}(V^i_{st}\oplus V^i_{triv})$, where the summands on the right hand side are simple right ideals.
\item If $X$ is of left cyclic type, then $\mathbb{k}[X]\cong V_{st}\oplus V_{triv}$, where the summands on the right hand side are simple left ideals.
\end{enumerate}
\end{corollary}
\para

Table \ref{2-transitive quandles} gives the number of right as well as left $2$-transitive quandles of order up to eight.

\begin{table}[H]
\begin{center}
\begin{tabular}{ |c||c||c||c| } 
 \hline
Order & Number of quandles & Number of right  &  Number of left  \\ 
&&$2$-transitive quandles&$2$-transitive quandles \\
  \hline  \hline
 $3$ & $3$ & $3$ & $2$ \\ 
  \hline
 $4$ & $7$ & $6$  & $3$  \\ 
 \hline
$5$ &  $22$ & $16$  & $7$ \\ 
 \hline
$6$ &  $73$ & $42$  & $14$ \\ 
 \hline
$7$ &  $298$ & $151$   & $39$\\ 
 \hline
$8$ &  $1581$ & $656$  & $105$ \\ 
 \hline
\end{tabular}
\end{center}
\caption{$2$-transitive quandles.}\label{2-transitive quandles}
\end{table}

For instance, the dihedral quandle $\R_n$ with odd $n\ge 5$ is not right $2$-transitive.

\begin{remark}{\rm 
Let $X= \Conj(G)$ be the conjugation quandle of a group $G$. Then the problem of decomposition of $\mathbb{k}[X]$ into simple right ideals is equivalent to decomposing $\mathbb{k}[G]$ into simple right ideals with respect to the conjugation action of $G$. This has been studied in detail in \cite{MR0289672, MR0232868}.}
\end{remark}
\bigskip
\bigskip

\section{Isomorphism problem for quandle rings}\label{iso}	
This section is motivated by the classical isomorphism problem for group rings.  We begin with the following definition.
	
\begin{definition}
Let $X$ be a finite quandle of order $n$. We say that $X$ has \index{partition type} \textit{partition type} $\lambda = (\lambda_1, \ldots, \lambda_n)$ if $X$ has $\lambda_i$ orbits of cardinality $i$ for each $1 \le i \le n$.
\end{definition}	

\begin{example}
The partition type of the dihedral quandle $\R_{2n}$ is $\lambda=(0,0, \ldots,0,2,0,\ldots,0)$, where $\lambda$ is a $2n$-tuple with 2 at $n$-th coordinate and 0 elsewhere. 
\end{example}	

We first derive a necessary condition for quandle rings of two quandles to be isomorphic \cite[Theorem 5.3]{MR3915329}.

\begin{theorem} \label{ISO}
Let $\mathbb{k}$ be an associative ring with unity such that $\charac(\mathbb{k})\neq 2,3$. Let  $X$ and $Y$ be right orbit $2$-transitive quandles
such that $\mathbb{k}[X] \cong \mathbb{k}[Y]$ and $A_n \not\le \Inn(X)$. Then $X$ and $Y$ are of the same partition type. 
\end{theorem}

\begin{proof}
Let $s$ denote the number of orbits in $X$, and $d$ the number of orbits in $Y$. Let $\lambda$ and $\mu$ be partition types of $X$ and $Y$, respectively.   Corollary~\ref{RightModDecomp}(1) implies that $\mathbb{k}[X]\cong \underset{i=1}{\overset{s}{\bigoplus}}(V^i_{\textrm{st}}\oplus V^i_{\textrm{triv}})$ and $\mathbb{k}[Y]\cong \underset{j=1}{\overset{d}{\bigoplus}}(W^j_{\textrm{st}}\oplus W^j_{\textrm{triv}})$, where the summands are simple right ideals. The isomorphism $\varphi: \mathbb{k}[X] \rightarrow \mathbb{k}[Y]$ induces another decomposition $\mathbb{k}[Y]=\overset{s}{\underset{j=1}{\oplus}}\varphi(V^i_{\textrm{st}})\oplus \varphi(V^i_{\textrm{triv}})$, where each summand is a simple right ideal in $\mathbb{k}[Y]$. The Krull-Schmidt theorem asserts that the decomposition of $\mathbb{k}[Y]$ is unique up to a permutation of summands, from which we conclude that $\lambda=\mu$. $\blacksquare$
\end{proof}

The \index{isomorphism problem}{\it isomorphism problem} for group rings asks whether an isomorphism $\mathbb{k}[G]\cong \mathbb{k}[H]$ between the group rings of two groups $G$ and $H$  necessarily implies that $G \cong H$. It is known that nilpotent groups and extensions of abelian groups by $p$-groups are determined by their integral group rings. A counterexample to this problem was provided by Hertweck \cite{MR1847590}, nearly sixty years after the question was first posed by Higman \cite{MR0002137}.
\para
 In analogy with group rings, we can formulate an isomorphism problem for quandle (respectively, rack) rings: Given quandles (respectively, racks) $X$ and $Y$, does an isomorphism $\mathbb{k}[X]\cong \mathbb{k}[Y]$ imply $X \cong Y$? We show that, in general, the answer is negative by presenting examples of non-isomorphic quandles whose quandle rings are isomorphic.

\begin{example}{\rm 
\label{CounterEx1}
Let $\mathbb{k}$ be a field with $\charac(\mathbb{k})=3$. Let $X$ and $Y$ be quandles with multiplication tables given by

		$$X=\begin{array}{|c||c c c c|} 
		\hline
		* \ &  x_1 & x_2 & x_3 & x_4  \\ \hline \hline
		\ x_1 &  x_1 & x_1 & x_2 & x_2  \\
		\ x_2 &  x_2 & x_2 & x_1 & x_1  \\
		\ x_3 &  x_3 & x_3 & x_3 & x_3  \\ 
		\ x_4 &  x_4 & x_4 & x_4 & x_4  \\ \hline
		\end{array}~~~~~~
	\quad \textit{and}	\quad 
	 Y=\begin{array}{|c||c c c c|}
		\hline
		\ * &  y_1 & y_2 & y_3 & y_4  \\ \hline \hline
		\ y_1 &  y_1 & y_1 & y_2 & y_1  \\ 
		\ y_2 &  y_2 & y_2 & y_1 & y_2  \\ 
		\ y_3 &  y_3 & y_3 & y_3 & y_3  \\
		\ y_4 &  y_4 & y_4 & y_4 & y_4  \\ \hline
		\end{array}.$$
The tables show that $X$ admits two elements whose corresponding inner automorphisms are non-trivial, whereas $Y$ has unique such element. Thus, we have $X \not\cong Y$. Clearly, the linear map represented by 
$$\varphi=\left(\begin{array}{cccc}
1 & 0 & 0 & 1\\
0 & 1 & 0 & 1\\
0 & 0 & 1 & 1\\
0 & 0 & 0 & 1\\
\end{array}\right)$$
is an isomorphism $\varphi:\mathbb{k}[X] \longrightarrow\mathbb{k}[Y]$ of underlying $\mathbb{k}$-vector spaces. To see that $\phi$ preserves ring multiplication, it suffices to check that $\phi(e_{x_i}  e_{x_j})=\phi(e_{x_i})  \phi(e_{x_j})$ for all $i, j$. Note that $\phi(e_{x_i})=e_{y_i}$ for $i=1,2,3$ and the subquandles $\{x_1,x_2,x_3\}$ and $\{y_1,y_2,y_3\}$, of respectively $X$ and $Y$, are isomorphic. Thus, we just need to check that $\phi(e_{x_4}  e_{x_i})=\phi(e_{x_4})  \phi(e_{x_i})$ and $\phi(e_{x_i}  e_{x_4})=\phi(e_{x_i})  \phi(e_{x_4})$ for all $i$. We see that
$\phi(e_{x_1}e_{x_4})=\phi(e_{x_2})=e_{y_2}=3e_{y_1}+e_{y_2}= e_{y_1}(e_{y_1}+e_{y_2}+e_{y_3}+e_{y_4})= \phi(e_{x_1}) \phi(e_{x_4}),$ since  $\charac(\mathbb{k})=3$.  In a manner, we obtain $\phi(e_{x_2}e_{x_4})=\phi(e_{x_1})=e_{y_1}=3e_{y_2}+e_{y_1}= \phi(e_{x_2}) \phi(e_{x_4})$ and
$\phi(e_{x_3}e_{x_4})=\phi(e_{x_3})=e_{y_3}=4e_{y_3}=  \phi(e_{x_1}) \phi(e_{x_4})$. For the remaining case, we have $\phi(e_{x_4}e_{x_j})=\phi(e_{x_4})=e_{y_1}+e_{y_2}+e_{y_3}+e_{y_4}=\phi(e_{x_4})\phi(e_{x_j})$, and hence $\varphi$ is an isomorphism of quandle rings.
} 
\end{example}
\para 

\begin{example}{\rm 
\label{CounterEx2}
Let $\mathbb{k}$ be a field with $\charac(\mathbb{k})=0$. Let $X$ and $Y$ be quandles with multiplication tables given by
		$$X=\begin{array}{|c||c c c c c c c|} 
		\hline
* \ &  x_1 & x_2 & x_3 & x_4 & x_5 & x_6 &  x_7 \\ \hline \hline
		\ x_1 &  x_1 & x_1 & x_1 & x_1 & x_2 & x_2 & x_1\\
		\ x_2 &  x_2 & x_2 & x_2 & x_2 & x_1 & x_1 & x_2 \\
		\ x_3 &  x_3 & x_3 & x_3 & x_3 & x_3 & x_4 & x_3 \\ 
		\ x_4 &  x_4 & x_4 & x_4 & x_4 & x_4 & x_3 & x_4 \\
		\ x_5 &  x_5 & x_5 & x_5 & x_5 & x_5 & x_5 & x_5 \\
		\ x_6 &  x_6 & x_6 & x_6 & x_6 & x_6 & x_6 & x_6 \\
		\ x_7 &  x_7 & x_7 & x_7 & x_7 & x_7 & x_7 & x_7\\ \hline
		\end{array}~~~~~~
		\quad  \textit{and}	\quad 
		 Y=\begin{array}{|c||c c c c c c c|} 
			\hline
* \ &  y_1 & y_2 & y_3 & y_4 & y_5 & y_6 &  y_7  \\ \hline \hline
		\ y_1 &  y_1 & y_1 & y_1 & y_1 & y_2 & y_1 & y_1\\
		\ y_2 &  y_2 & y_2 & y_2 & y_2 & y_1 & y_2 & y_2 \\
		\ y_3 &  y_3 & y_3 & y_3 & y_3 & y_3 & y_4 & y_3 \\ 
		\ y_4 &  y_4 & y_4 & y_4 & y_4 & y_4 & y_3 & y_4 \\
		\ y_5 &  y_5 & y_5 & y_5 & y_5 & y_5 & y_5 & y_5  \\
		\ y_6 &  y_6 & y_6 & y_6 & y_6 & y_6 & y_6 & y_6 \\
		\ y_7 &  y_7 & y_7 & y_7 & y_7 & y_7 & y_7 & y_7  \\ \hline
		\end{array}.$$
The tables shows that $X$ admits an element (namely $x_6$) whose corresponding inner automorphism is a product of two disjoint transpositions. On the other hand, $Y$ has no such element, and hence $X \not\cong Y$. We see that the linear map represented by 
\[\varphi=\left(\begin{array}{ccccccc}
1 & 0 & 0 & 0 & 0 & 0 & 0\\
0 & 1 & 0 & 0 & 0 & 0 & 0\\
0 & 0 & 1 & 0 & 0 & 0 & 0\\
0 & 0 & 0 & 1 & 0 & 0 & 0\\
0 & 0 & 0 & 0 & 1 & 1 & 0\\
0 & 0 & 0 & 0 & 0 & 1 & 0\\
0 & 0 & 0 & 0 & 0 & -1 & 1\\
\end{array}\right)\]
is an isomorphism $\varphi:\mathbb{k}[X] \longrightarrow\mathbb{k}[Y]$ of underlying $\mathbb{k}$-vector spaces. As in the preceding example, it only suffices to check that $\phi(e_{x_6}  e_{x_i})=\phi(e_{x_6})  \phi(e_{x_i})$ and $\phi(e_{x_i}  e_{x_6})=\phi(e_{x_i})  \phi(e_{x_6})$  for all $i$. These can be verified by straightforward computation.  For example,  $\phi(e_{x_1}e_{x_6})=\phi(e_{x_2})=e_{y_2}=e_{y_1}(e_{y_5}+e_{y_6}-e_{y_7})=\phi(e_{x_1})  \phi(e_{x_6})$ and $\phi(e_{x_6}e_{x_6})= e_{y_5}+e_{y_6}-e_{y_7}=(e_{y_5}+e_{y_6}-e_{y_7})(e_{y_5}+e_{y_6}-e_{y_7}) =\phi(e_{x_6})  \phi(e_{x_6}).$  
}
\end{example}
  \para
  
\begin{problem}\label{ques1}
Classify quandles for which the isomorphism problem has an affirmative solution over some ring.
\end{problem}
\bigskip
\bigskip

\chapter{Zero-divisors and idempotents in quandle rings} \label{chap zero-divisors and idempotents}

\begin{quote}
In this chapter, we investigate zero-divisors and idempotents in quandle rings. Since each element of the quandle can be viewed as an idempotent in its quandle ring, it turns out idempotents in quandle rings are natural objects of interest. We examine an intimate relation between quandle coverings and idempotents in quandle rings. We show that integral quandle rings of free quandles admit only trivial idempotents, thus providing an infinite family of quandles with this property. We also propose some conjectures pertaining structure of idempotents in quandle rings.
\end{quote}
\bigskip

\section{Zero-divisors in quandle rings}\label{orderability-zero-divisors}
Recall that, a non-zero element $u$ of a ring is called a \index{zero-divisor} {\it zero-divisor} if there exists a non-zero element $v$ such that either $uv=0$ or $vu=0$. Every non-zero nilpotent element of an associative ring is a zero-divisor. Determining whether group rings of torsion-free groups over fields have zero-divisors is a classical and still open problem in the theory of group rings. In this section, we investigate the analogous problem for quandle rings.
\para
Let $\mathbb{k}$ be an integral domain with unity. It is easy to see that if $\T$ is a trivial quandle with more than one element, then $\mathbb{k}[\T]$ contains zero-divisors. If $G$ is a group with an element $g$  of finite order, say $n>1$, then the element $$\tilde{g}:= 1+ g+ \cdots+ g^{n-1}$$ of the integral group ring $\mathbb{Z}[G]$ satisfies
$\tilde{g}(1-g)=(1-g)\tilde{g}=0$, and hence $\mathbb{Z}[G]$ has a  zero-divisor. By analogy, the next proposition gives a condition under which a quandle ring fails to be an integral domain.
	
\begin{proposition}
Let $X$ be a quandle admitting a finite component $Y$ such that $1 < |Y|<\infty$. Then $\mathbb{k}[X]$ contains zero-divisors.
\end{proposition}

\begin{proof}
Indeed, it follows from the second quandle axiom that  $$\Big(\sum_{z\in Y} e_{z} \Big) (e_{x}-e_{y})=0,$$ where $x$ and $y$ are any two distinct elements of $X$. $\blacksquare$
\end{proof}

We first formulate some sufficient conditions under which a quandle ring contains zero-divisors. We say that a quandle $X$ containing more than one element is \index{inert} {\it inert} if there is a finite subset $A = \{ a_1, \ldots, a_n \}$ of $X$ and two distinct elements $x, y \in X$ such that $A x = A y$, where $A z = \{ a_1 *z, \ldots, a_n *z \}$.

\begin{proposition}\label{two different elements}
The following  assertions hold:
\begin{enumerate}
\item  If $X$ is a quandle containing a trivial subquandle with more than one element, then $\mathbb{k}[X]$ contains zero-divisors.
\item If $X$ is an inert quandle, then $\mathbb{k}[X]$ contains zero-divisors. In particular, if $X$ contains a finite subquandle with more than one element, then $\mathbb{k}[X]$ contains zero-divisors.
\item If $X$ is not semi-latin, then $\mathbb{k}[X]$ contains zero-divisors.
\end{enumerate}
\end{proposition}

\begin{proof}
Let $T = \{ x, y  \}$ be a trivial subquandle in $X$. Taking $u = e_{x} - e_{y} \in \mathbb{k}[X]$ gives
$$
u^2 = (e_{x} - e_{y}) (e_{x} - e_{y}) = 0,
$$
which proves assertion (1).
\para

For assertion  (2), let $x$ and $y$ be two distinct elements in $X$ and  $A = \{ a_1, \ldots, a_m \}$ such that $A x = A y$. Then we have
$$
(e_{a_1} + \cdots + e_{a_m}) (e_{x} - e_{y}) = 0.
$$
If $X$ contains a finite subquandle $A$ with more than one element, then we can take $x$ and $y$ to be two distinct elements of $A$.
\para
For assertion (3), suppose that for some $x\in X$ there exist distinct $y, z \in X$ such that $x*y = x*z$. Then we have
$e_{x} (e_{y}-e_{z}) = 0.$ $\blacksquare$
\end{proof}

The following is a natural problem to explore.

\begin{problem}\label{Elhamdadi-problem}
Determine non-trivial quandles whose quandle rings are free of zero-divisors.
\end{problem}

The following definition is motivated from group theory (see, for example, \cite[Chapter 13]{MR0470211}).

\begin{definition}
A quandle $X$ is said to be a {\it up-quandle} (unique product quandle) if given any two non-empty finite subsets $A$ and $B$ of $X$, there exists an element $x \in X$ that has a unique representation of the form $x = a *b$ for some $a \in A$ and $b\in B$. 
\medskip

A quandle $X$ is said to be a {\it tup-quandle} (two unique product quandle) if given any two non-empty finite subsets $A$ and $B$ of $X$ with $|A| + |B| > 2$, there exist elements $x, y \in X$ that have unique representations of the form $x = a *b$ and $y = c *d$, where $a, c \in A$ and $b, d \in B$. 
\end{definition}

It is clear that every tup-quandle is a up-quandle. The following observation is an analogue of the corresponding result for groups \cite[Chapter 13, Lemma 1.9]{MR0470211}.

\begin{proposition} \cite[Proposition 3.3]{MR4450681} \label{up-zero}
If $X$ is a up-quandle and $\mathbb{k}$ an integral domain with unity, then $\mathbb{k}[X]$ has no zero-divisors.
\end{proposition}

\begin{proof}
Let $u$ and $v$ be non-zero elements of $\mathbb{k}[X]$. We write
$$
u = \sum_{i=1}^n \alpha_i e_{x_i} \quad \textrm{and} \quad v = \sum_{j=1}^m \beta_j e_{y_j},
$$
where $\alpha_i, \beta_j$ are non-zero elements of $\mathbb{k}$. Take $A = \{ x_1, \ldots, x_n \}$ and $B = \{ y_1, \ldots, y_m \}$. Then
$$
u v = \sum_{i, j} \alpha_i \beta_j e_{x_i}  e_{y_j},
$$
where each $\alpha_i \beta_j  \not= 0$ since $\mathbb{k}$ has no zero-divisors. Since $X$ is a up-quandle, there exists a uniquely represented element in the product $A B$, say $z = x_i *y_j$. It follows that the non-zero summand $\alpha_i \beta_j e_{x_i} e_{y_j}$ cannot be cancelled by any other term in the product $uv$. Hence, $u v \not= 0$ and $\mathbb{k}[X]$ has no zero-divisors. $\blacksquare$
\end{proof}

We now consider orderable quandles to give explicit examples of up-quandles. It is known that the  group ring of a right-orderable group has no zero-divisors \cite[Chapter 13]{MR0470211}. On the other hand, a trivial quandle  with more than one element is right-orderable and its quandle ring always has zero-divisors. However, for semi-latin quandles we have the following result, which is a quandle analogue of \cite[Chapter 13, Lemma 1.7]{MR0470211}, and also contributes to Problem \ref{Elhamdadi-problem}.

\begin{proposition} \cite[Proposition 3.10]{MR4450681}\label{orderable-tup}
Let $X$ be a semi-latin quandle which is right- or left-orderable. Then $X$ is a tup-quandle.
\para
In fact, if $A$ and $B$ are non-empty  finite subsets of $X$, then there exist $b', b'' \in B$ such that the product $a_{max} *b'$ and $a_{min} * b''$ are uniquely represented in $A B$, where $a_{max}$ denotes the largest element in $A$ and $a_{min}$ the smallest.
\end{proposition}

\begin{proof} Suppose that $X$ is a semi-latin and right-orderable quandle. Let
$$
A = \{a_{min} = a_1 < a_2 < \cdots < a_n = a_{max} \}, ~~\textrm{where}~~n \geq 2,
$$ 
and 
$$
B = \{ b_1 < b_2 < \cdots < b_m \}
$$ 
be two non-empty finite subsets of $X$. We write the elements of the product $A B$ in the tabular form

$$
\begin{array}{ccccccc}
a_1 *b_1 & < & a_2 *b_1 & < & \cdots & < & a_n *b_1, \\
a_1 *b_2 & < & a_2 *b_2 & < & \cdots & < & a_n *b_2, \\
\vdots & \vdots & \vdots & \vdots & \vdots & \vdots & \vdots \\
a_1 *b_m & < & a_2 *b_m & < & \cdots & < & a_n *b_m, \\
\end{array}
$$
where the inequalities in the rows follow from the right-ordering of $X$. Since $X$ is semi-latin, it follows that all  the entries in each column are distinct.
\para
Let $b_i \in B$ be the element  such that $a_1 *b_i$ is the smallest element in the first column.  We claim that $a_{min} *b_i = a_1 *b_i$ is uniquely represented in $A B$. It suffices to prove that $a_1 *b_i < a_k *b_l$ for any pair $(k, l) \not= (1, i)$. If $k=1$, then $l \neq i$ and the inequality $a_1 *b_i < a_1 *b_l$ follows from the choice of $b_i$. If $k >1$, then $a_1 *b_i \leq a_1 *b_l$ and the inequalities in the $l$-th row  imply that $a_1 *b_i < a_k *b_l$.
\para
Let $b_j \in B$ be the element such that $a_n *b_j$ is the largest element in the last column. We claim that $a_k *b_l < a_n *b_j$ for each $(k, l) \not= (n, j)$. If $k = n$, then the inequality follows from the choice of $b_j$. If $k < n$, then  inequalities in the $l$-th row gives
$$
a_k *b_l < a_n *b_l \leq a_n *b_j.
$$
Hence, the product $a_{max} *b_j = a_n *b_j$ is  uniquely represented in $A B$. The case when $X$ is left-orderable is similar. $\blacksquare$
\end{proof}

Proposition \ref{up-zero} and Proposition \ref{orderable-tup}  together yield the following result.

\begin{theorem} \label{orderable-latin}
Let  $X$ be a semi-latin quandle that is right- or left-orderable. Then $\mathbb{k}[X]$ has no zero-divisors.
\end{theorem}

Theorem \ref{orderability-free-quandle} and Theorem \ref{orderable-latin}  leads to the following result.

\begin{corollary}\label{free-quandle-zero-divisors}
Quandle rings of free quandles have no  zero-divisors.
\end{corollary}

As a consequence of Proposition \ref{conj right orderable}, Proposition \ref{orderable-implies-latin}(2) and Theorem \ref{orderable-latin}, we have the following results.

\begin{corollary}
Let $G$ be a bi-orderable group. Then the following assertions hold:
\begin{enumerate}
\item $\mathbb{k}[\Core(G)]$ has no  zero-divisors.
\item If $\phi \in \Aut(G)$ is an order reversing automorphism, then $\mathbb{k}[\Alex(G, \phi)]$ has no  zero-divisors.

\end{enumerate}
\end{corollary}
\para

Proposition \ref{two different elements} and Theorem \ref{orderable-latin} suggest the following analogue of \index{Kaplansky's zero-divisor conjecture}Kaplansky's zero-divisor conjecture for quandles.

\begin{conjecture}\label{kaplansky-analogue}
Let $\mathbb{k}$ be an integral domain with unity and $X$ a non-inert semi-latin quandle.  Then the quandle ring $\mathbb{k}[X]$ has no zero-divisors.
\end{conjecture}
\bigskip
\bigskip

\section{Idempotents in quandle rings}\label{sec-idempotents}
Units in group rings play a fundamental role in the structure theory of group rings. In contrast, it turns out that idempotents are the most natural objects in quandle rings since each quandle element is, by definition, an idempotent of the quandle ring. In general, the computation of idempotents is an important problem in ring theory.  It is well-known that integral group rings do not have non-trivial idempotents (see \cite[p. 123]{MR0349646} or \cite[p. 38]{MR0470211}). In contrary, we shall see that integral quandle rings of many non-trivial quandles admit non-trivial idempotents.
\para

Let $X$ be a quandle and $\mathbb{k}$ an integral domain with unity.  A non-zero element $u \in \mathbb{k}[X]$ is called an \index{idempotent} {\it idempotent} if $u^2=u$. Let $\mathcal{I}\big(\mathbb{k}[X] \big)$ denote the set of all idempotents of $\mathbb{k}[X]$. It is clear that the basis elements $\{e_x \, \mid \, x \in X\}$ are idempotents of $\mathbb{k}[X]$, and we refer them as \index{trivial idempotent} {\it trivial idempotents}.  A non-trivial idempotent  is an element of $\mathbb{k}[X]$ that is not of the form $e_x$ for any $x \in X$. 
\para

Clearly, if $Y$ is a subquandle of $X$, then $\mathcal{I} \big(\mathbb{k}[Y] \big) \subseteq \mathcal{I} \big(\mathbb{k}[X] \big)$. The construction of a quandle ring is functorial. Hence, a quandle homomorphism $\phi: X \rightarrow Z$ induces a ring homomorphism $$\hat{\phi}: \mathbb{k}[X] \rightarrow \mathbb{k}[Z]$$ which maps $ \mathcal{I} \big(\mathbb{k}[X] \big)$ into $\mathcal{I} \big(\mathbb{k}[Z] \big)$. Since the augmentation map $\varepsilon: \mathbb{k}[X] \rightarrow \mathbb{k}$ is a ring homomorphism, it follows that $\varepsilon(u)=0$ or $\varepsilon(u)=1$ for each idempotent $u$ of $\mathbb{k}[X]$.
\para

It is a well-known result of Swan \cite[p. 571]{MR138688} that if $G$ is a finite group, then the group ring $\mathbb{k}[G]$ has a non-trivial idempotent if and only if some prime divisor of $|G|$ is invertible in $\mathbb{k}$. Although, we do not have Lagrange's theorem for finite quandles, a partial one way analogue of this result does hold for finite quandles.
 
\begin{proposition}\label{prop swan analogue}
Let $X$ be a quandle and $\mathbb{k}$ an integral domain with unity. If $X$ admits a subquandle $Y$ with more than one element such that $|Y|$ is invertible in $\mathbb{k}$, then $\mathbb{k}[X]$ has a non-trivial idempotent.
\end{proposition}

\begin{proof}
Since the subquandle $Y$ has more than one element, a direct check shows that the element $u= \frac{1}{|Y|}\sum_{y \in Y} e_y$ is a non-trivial idempotent of $\mathbb{k}[X]$. $\blacksquare$
\end{proof}
 
\begin{remark}{\rm 
The converse of Proposition \ref{prop swan analogue} does not hold in general. For example, consider the quandle 
$$
\J_3=\begin{tabular}{|c||c|c|c|}
    \hline
$*$ & 1 & 2 & 3  \\
  \hline \hline
1 & 1 & 1 & 1  \\
\hline
2 & 3 & 2 & 2  \\
\hline
3 & 2 & 3 & 3  \\
  \hline
\end{tabular}
$$
Then, the quandle ring $\mathbb{Z}[\J_3]$ has non-trivial idempotents of the form $\alpha e_2+ (1-\alpha) e_3$ for $\alpha \in \mathbb{Z}$, but $J_3$ has no subquandle $Y$ with more than one element such that $|Y|$ is invertible in $\mathbb{Z}$.}
\end{remark}
  
\begin{proposition}\label{trivial subquandle idempotent}
Let $X$ be a quandle and $\mathbb{k}$ an integral domain with unity. If $X$ contains a  trivial subquandle $Y$ with more than one element, then $\mathbb{k}[X]$ has non-trivial idempotents.
\end{proposition}

\begin{proof}
Consider the element $u= \sum_{i=1}^n \alpha_i e_{y_i}$, where $n \ge2$,  $y_i \in Y$ and $\alpha_i \in \mathbb{k}$ with  $\sum_{i=1}^n \alpha_i=1$. A direct check shows that $u^2=u$, and hence $u$ is a non-trivial idempotent of $\mathbb{k}[X]$. $\blacksquare$
\end{proof}

Recall that a quandle is \index{faithful}{\it faithful} if the natural map $S: X \to \Inn(X)$ given by $S(x)=S_x$ is injective.

\begin{proposition}\label{fixed point idempotent}
Let $X$ be a faithful quandle such that $S_z$ has more than one fixed-point for some $z \in X$. Then $\mathbb{k}[X]$ has non-trivial idempotents.
\end{proposition}

\begin{proof}
Since $S_z$ has a non-trivial fixed-point, we have $x*z=x$ for some $x \ne z$. Since $S_zS_x=S_{x * z} S_z$, it follows that $S_x$ and $S_z$ commute.  Thus, the identity $S_x S_z=S_{z * x} S_x$ implies that $S_{z * x}=S_z$.  Since $X$ is faithful, we obtain $z * x=z$, and hence $\{x, z \}$ is a trivial subquandle of $X$. The result now follows from Proposition \ref{trivial subquandle idempotent}. $\blacksquare$
\end{proof}

\begin{proposition}
If $G$ is a non-trivial group and $\mathbb{k}$ an integral domain with unity, then $\mathbb{k}[\Conj(G)]$ has non-trivial idempotents.
\end{proposition}

\begin{proof}
Note that, for each non-identity element $x \in G$ and distinct integers $i, j$, the set $ \{x^i, x^j \}$ forms a trivial subquandle of $\Conj(G)$. The result now follows from Proposition \ref{trivial subquandle idempotent}. $\blacksquare$
\end{proof}

As an application to quandle rings of link quandles, we have the following result.

\begin{proposition}\label{Hopf link idempotent}
Let $L$ be a link containing the Hopf link as its sublink and $\mathbb{k}$ an integral domain with unity. Then $\mathbb{k}[Q(L)]$ has non-trivial idempotents.
\end{proposition}

\begin{proof}
Let $H$ be the Hopf link. It follows from the construction of the link quandle (see Example \ref{construction link quandle}) that $Q(L)$ contains $Q(H)$ as a subquandle, which is a trivial quandle with two elements. The result now follows from Proposition \ref{trivial subquandle idempotent}. $\blacksquare$
\end{proof}

It is tempting to speculate that connected quandles have only trivial idempotents. We give two examples showing that this is not true in general.

\begin{example}\label{involutive and faithful example1}
{\rm 
Let $X= \{1, 2, 3, 4, 5, 6 \}$ be the connected quandle of order 6 (see \cite{rig}) with multiplication table given by
$$X=\begin{tabular}{|c||c|c|c|c|c|c|}
    \hline
$*$ & 1 & 2 & 3 & 4 &5 &6 \\
  \hline \hline
1 & 1 & 1 & 5 & 6 & 3 & 4 \\
\hline
2 & 2 & 2 & 6 & 5 & 4 & 3 \\
\hline
3 & 5 & 6 & 3  & 3 & 1 & 2\\
\hline
4 & 6 & 5 & 4 & 4 & 2 & 1 \\
\hline
5 & 3 & 4 & 1  & 2 & 5 & 5\\
\hline
6 & 4 & 3 & 2  & 1 & 6 & 6\\
  \hline
\end{tabular}.
$$
Note that $X$ has trivial subquandles $\{1,2\}$, $\{3, 4\}$ and $\{5, 6\}$. Hence, by Proposition \ref{trivial subquandle idempotent}, the elements
$\alpha e_1+(1-\alpha) e_2$, $\beta e_3 +(1-\beta)e_4$ and $\gamma e_5+(1-\gamma)e_6$  are non-trivial idempotents of $\mathbb{Z}[X]$ for each $\alpha, \beta, \gamma \in \mathbb{Z}$.}
\end{example}

\begin{example}\label{involutive and faithful example2}
{\rm 
Consider the connected quandle $X=\{1, 2, \ldots,12\}$ of order $12$ (see \cite{rig}) with multiplication table given by
$$X=\begin{tabular}{|c||c|c|c|c|c|c|c|c|c|c|c|c|}
    \hline
$*$ & 1 & 2 & 3 & 4 & 5 & 6 & 7 & 8 & 9 & 10 & 11 & 12 \\
  \hline \hline
1 & 1 & 1 & 1 & 1 & 9 & 10 & 11 & 12 & 7 & 8 & 5 & 6\\
\hline
2 & 2 & 2 & 2 & 2 & 10 & 9 & 12 & 11 & 8 & 7 & 6 & 5\\
\hline
3 & 3 & 3 & 3  & 3 & 11 & 12 & 9 & 10 & 5 & 6 & 7 & 8\\
\hline
4 & 4 & 4 & 4 & 4 & 12 & 11 & 10 & 9 & 6 & 5 & 8 & 7\\
\hline
5 & 11 & 12 & 9  & 10 & 5 & 5 & 5 & 5 & 1 & 2 & 3 & 4\\
\hline
6 & 12 & 11 & 10  & 9 & 6 & 6 & 6 & 6 & 2 & 1 & 4 & 3\\
  \hline
7 & 9 & 10 & 11  & 12 & 7 & 7 & 7 & 7 & 3 & 4 & 1 & 2\\
\hline
8 & 10 & 9 & 12  & 11 & 8 & 8 & 8 & 8 & 4 & 3 & 2 & 1\\
\hline
9 & 5 & 6 & 7  & 8 & 3 & 4 & 1 & 2 & 9 & 9 & 9 & 9\\
\hline
10 & 6 & 5 & 8  & 7 & 4 & 3 & 2 & 1 & 10 & 10 & 10 & 10\\
\hline
11 & 7 & 8 & 5  & 6 & 1 & 2 & 3 & 4 & 11 & 11 & 11 & 11\\
\hline
12 & 8 & 7 & 6  & 5 & 2 & 1 & 4 & 3 & 12 & 12 & 12 & 12\\
\hline
\end{tabular}.
$$

We see that $X$ has trivial subquandles $\{1,2,3,4\}$, $\{5,6,7,8\}$ and $\{9,10,11,12\}$. By Proposition \ref{trivial subquandle idempotent}, the elements $\alpha e_1+\beta e_2+\gamma e_3 +(1-\alpha-\beta-\gamma)e_4$, $\alpha e_5+\beta e_6+\gamma e_7 +(1-\alpha-\beta-\gamma)e_8$ and $\alpha e_9+\beta e_{10}+\gamma e_{11} +(1-\alpha-\beta-\gamma)e_{12}$ are non-trivial idempotents of $\mathbb{Z}[X]$ for each $\alpha, \beta, \gamma \in \mathbb{Z}$.}
\end{example}

A computer assisted check for quandles of order less than seven suggests the following conjecture \cite[Conjecture 3.1]{MR4565221}.

\begin{conjecture}\label{Main conjecture}
The integral quandle ring of a semi-latin quandle has only trivial idempotents. In particular, the integral quandle ring of a finite latin quandle has only trivial idempotents.
\end{conjecture}

It follows from Theorem \ref{orderability-free-quandle} that free quandles are semi-latin. We shall prove in Theorem \ref{idempotents in free products} that free quandles satisfy Conjecture \ref{Main conjecture}. 
\para

For a quandle $X$, we denote by $\Aut_{\textrm{algebra}} \big(\mathbb{k}[X] \big)$ the group of $\mathbb{k}$-algebra automorphisms of $\mathbb{k}[X]$, that is, ring automorphisms of $\mathbb{k}[X]$ that are $\mathbb{k}$-linear. 

\begin{proposition}\label{auto quandle ring free products}
Let $X$ be a quandle and $\mathbb{k}$ an integral domain with unity. If $\mathbb{k}[X]$ admits only trivial idempotents, then 
$\Aut_{\textrm{algebra}}\big(\mathbb{k}[X]\big) \cong \Aut_{\textrm{quandle}}(X).$
\end{proposition}

\begin{proof}
Clearly, each automorphism of $X$ induces an automorphism of $\mathbb{k}[X]$. Conversely, if $\phi \in \Aut_{\textrm{algebra}}\big(\mathbb{k}[X]\big)$, then $\phi$ is a bijection of the set $\mathcal{I}\big(\mathbb{k}[X]\big)$ of all idempotents. Since $\mathbb{k}[X]$ has only trivial idempotents, we have $X \cong \mathcal{I}\big(\mathbb{k}[X]\big)$ via the map $x \mapsto e_x$, and hence $\phi$ can be viewed as an automorphism of $X$, proving the assertion. $\blacksquare$
\end{proof}

\para
In what follows, we give a complete description of idempotents for some elementary quandles.

\begin{proposition} \label{trivial-quandle-idempotents}
If $\T$ is a trivial quandle and $\mathbb{k}$ an integral domain with unity, then $\mathcal{I}(\mathbb{k}[\T]) =  e_{x_0} + \Delta_{\mathbb{k}} (\T)$, where  $x_0 \in \T$ is a fixed element.
\end{proposition}

\begin{proof}
Since $\T$ is trivial, by Theorem \ref{deltasqzero}, $\Delta_{\mathbb{k}}^2 (\T) = 0$. It follows that non-zero idempotents do not lie in $\Delta_{\mathbb{k}} (\T)$. Hence, a non-zero idempotent has the form
 $u = e_{x_0} + \delta$, where $\delta \in \Delta_{\mathbb{k}} (\T)$ and $x_0 \in \T$ some fixed element. Indeed, 
$$
u^2 = e_{x_0}^2 + e_{x_0} \delta + \delta e_{x_0} + \delta^2= e_{x_0}+ \delta=u
$$
since $e_{x_0}^2 = e_{x_0}$, $e_{x_0} \delta = \delta^2 = 0$ and $\delta e_{x_0} = \delta$. $\blacksquare$
\end{proof}

\begin{proposition}\cite[Proposition 4.2]{MR4450681} \label{joyce-quandle}
Let $\J_3$ be the 3-element quandle given by
$$
\J_3=\begin{tabular}{|c||c|c|c|}
    \hline
$*$ & 1 & 2 & 3  \\
  \hline \hline
1 & 1 & 1 & 1  \\
\hline
2 & 3 & 2 & 2  \\
\hline
3 & 2 & 3 & 3  \\
  \hline
\end{tabular}.
$$
Then $\mathcal{I} \big(\mathbb{Z}[\J_3] \big) = \big\{ (1 - \beta) e_2 +  \beta e_3,~ \alpha e_2 + \alpha e_3 +  (1 - 2 \alpha) e_1~ \mid~ \beta, \alpha \in \mathbb{Z} \big\}$.
\end{proposition}

\begin{proof}
The quandle $\J_3$ is a disjoint union of two trivial subquandles $\{1\}$ and $\{ 2, 3 \} $. If $w = \alpha e_2 + \beta  e_3 + \gamma e_1 \in \mathbb{Z}[\J_3]$, then
$$
w^2 = (\alpha^2 + \alpha \beta + \beta \gamma) e_2 + (\alpha \beta + \beta^2 + \alpha \gamma) e_3 + (\alpha \gamma + \beta \gamma + \gamma^2 )e_1.
$$
Thus, $w$ is an idempotent if and only if the system of equations
\begin{eqnarray*}
\alpha &=& \alpha^2 + \alpha \beta + \beta \gamma,\\
\beta &=& \beta^2 + \alpha \beta +  \alpha \gamma,\\
\gamma &=& \gamma^2 + \alpha \gamma + \beta \gamma,
\end{eqnarray*}
is simultaneously  solvable over integers. Suppose that $\gamma = 0$. Then we obtain the equations
$$\alpha = \alpha^2 + \alpha \beta \quad \textrm{and} \quad \beta = \beta^2 + \alpha \beta.$$
If $\alpha = 0$, then $\beta = 0, 1$. The first case gives $w=0$, while the second case gives $w = e_3$. If $\alpha \not= 0$, then $\alpha = 1 - \beta$ and we have idempotents
$$
w = (1 - \beta) e_2 +  \beta e_3,~~\beta \in \mathbb{Z}.
$$
These are idempotents of the quandle ring $\mathbb{Z}[\{ 2, 3 \}]$.
\para

Now, suppose that $\gamma \not= 0$. Then the third equation of the system gives $\gamma = 1 - \alpha - \beta$. Substituting this expression in the first and the second equations gives
$$
(\alpha - \beta) (\alpha + \beta) = (\alpha - \beta).
$$
If $\alpha - \beta \not= 0$, then we obtain the same idempotent as in the previous case. If $\alpha - \beta = 0$, then we have idempotents
$$
w = \alpha e_2 + \alpha e_3 +  (1 - 2 \alpha) e_1,~~\alpha \in \mathbb{Z}.
$$
Hence, we have proved that $\mathcal{I} \big(\mathbb{Z}[\J_3] \big) = \big\{ (1 - \beta) e_2 +  \beta e_3, \alpha e_2 + \alpha e_3 +  (1 - 2 \alpha) e_1 \, \mid \, \beta, \alpha \in \mathbb{Z} \big\}.
$ $\blacksquare$
\end{proof}

\begin{remark}{\rm 
 If a quandle $X = X_1 \cup \cdots \cup X_n$ is a union of subquandles and $\mathbb{k}$ an integral domain with unity, then
\begin{equation}\label{idempotents-union}
  \mathcal{I} \big(\mathbb{k}[X_1] \big) \cup \cdots \cup \mathcal{I} \big(\mathbb{k}[X_n] \big) \subseteq \mathcal{I} \big(\mathbb{k}[X] \big).
\end{equation}
The inclusion is, in general, not an equality. In fact, the preceding example shows that
$$
 \mathcal{I} \big(\mathbb{Z}[\{ 2, 3 \}] \big) \cup \mathcal{I} \big(\mathbb{Z}[\{1\}] \big) \neq \mathcal{I} \big(\mathbb{Z}[\J_3] \big).
$$}
\end{remark}
 \medskip

Let  $\R_n = \{0, 1, \ldots, n-1 \}$ be the dihedral quandle of order $n$, where $i*j=2j-i \mod  n$. We examine idempotents in $\mathbb{Z}[\R_n]$ for $1 \le n\le 4$. Note that $\R_1$ and $\R_2$ are trivial quandles.

\begin{proposition}\cite[Proposition 4.3]{MR4450681} \label{idempotents-r3}
$ \mathcal{I} \big(\mathbb{Z}[\R_3] \big) = \{ e_0, e_1, e_2 \}. $
\end{proposition}

\begin{proof}
Let $z = \alpha_0 e_0 + \alpha_1 e_1 + \alpha_2 e_2 \in \mathbb{Z}[\R_3]$ be an idempotent. Then we have $\varepsilon(z) = 0, 1$.
\medskip

 Case 1: $\varepsilon(z) = 0$, that is, $\alpha_0 = - \alpha_1 - \alpha_2$. Then $z = \alpha_1 f_1 + \alpha_2 f_2$, where $f_i = e_i - e_0$. The elements $f_1$ and $f_2$ generate $\Delta (\R_3)$ and have the following multiplication table.
\begin{center}
\begin{tabular}{|c||c|c|}
  \hline
 $\cdot$ & $f_1$ & $f_2$  \\
   \hline
     \hline
$f_1$ & $f_1 - 2 f_2$ & $-f_1 - f_2$  \\
  \hline
$f_2$  & $-f_1 - f_2$ & $-2 f_1 + f_2$ \\
  \hline
\end{tabular}\\
\end{center}
Thus, we have
$$
z^2 = (\alpha_1^2 - 2 \alpha_1 \alpha_2 - 2 \alpha_2^2) f_1 + (\alpha_2^2 - 2 \alpha_1 \alpha_2 - 2 \alpha_1^2) f_2.
$$
The equality $z^2 = z$ leads to the equations
$$\alpha_1 = \alpha_1^2 - 2 \alpha_1 \alpha_2 - 2 \alpha_2^2 \quad \textrm{and} \quad \alpha_2 = \alpha_2^2 - 2 \alpha_1 \alpha_2 - 2 \alpha_1^2.$$
Subtracting the second equation from the first yields
$$
\alpha_1 - \alpha_2 = 3 (\alpha_1 - \alpha_2) (\alpha_1 + \alpha_2).
$$
It is not difficult to see that in this case the only solution is $\alpha_1=\alpha_2=\alpha_3=0$.
\para 

Case 2: $\varepsilon(z) = 1$. In this case
 $\alpha_0 = 1 - \alpha_1 - \alpha_2$. Then $z = e_0 + \alpha_1 f_1 + \alpha_2 f_2$ and we get
$$
z^2 = e_0^2+ (2 \alpha_2 + \alpha_1^2 - 2 \alpha_1 \alpha_2 - 2 \alpha_2^2) f_1 + (2 \alpha_1 + \alpha_2^2 - 2 \alpha_1 \alpha_2 - 2 \alpha_1^2) f_2.
$$
 From $z^2 = z$, we obtain the equations
$$\alpha_1 - 2 \alpha_2 = \alpha_1^2 - 2 \alpha_1 \alpha_2 - 2 \alpha_2^2 \quad \textrm{and} \quad \alpha_2 - 2 \alpha_1 = \alpha_2^2 - 2 \alpha_1 \alpha_2 - 2 \alpha_1^2.$$
Subtracting the second from the first gives
$$
\alpha_1 - \alpha_2 =  (\alpha_1 - \alpha_2) (\alpha_1 + \alpha_2).
$$
If $\alpha_1 = \alpha_2$, then the system is equivalent to the equation $\alpha_1 = 3 \alpha_1^2$, whose only solution is $\alpha_1=\alpha_2=0$. Thus, $\alpha_0=1$, and hence $z=e_0$ in this case. If $\alpha_1 \not = \alpha_2$, then the system has solutions $(\alpha_1, \alpha_2) = (1,0)$ or $(0,1)$. In this case $\alpha_0 = 0$, and hence $z=e_1$ or $e_2$. $\blacksquare$
\end{proof}
\para

Note that $\R_4$ is disconnected and is a disjoint union of trivial subquandles $\{0, 2 \} $ and $\{1, 3 \}$. Thus, we have
$$
  \mathcal{I} \big(\mathbb{Z}[\{e_0, e_2 \}] \big) \cup \mathcal{I} \big(\mathbb{Z}[\{e_1, e_3 \}] \big) \subseteq \mathcal{I} \big(\mathbb{Z}[\R_4] \big).
$$
In fact, we have an equality in this case \cite[Proposition 4.4]{MR4450681}.

\begin{proposition}\label{R4-idempotents}
$
\mathcal{I} \big(\mathbb{Z}[\R_4] \big) = \Big\{ \alpha e_0 +  (1 - \alpha) e_2, ~\beta e_1 + (1-\beta) e_3  \, \mid \,  \alpha, \beta \in \mathbb{Z} \Big\}.
$
\end{proposition}

\begin{proof}
If $z \in \mathbb{Z}[\R_4]$ is an idempotent, then $\varepsilon (z) = 0, 1$. We consider the following two cases:
\para 

Case 1: $\varepsilon (z) = 0$, that is, $z \in \Delta (\R_4)$. In this case, we can write
$$
z = \alpha f_1 + \beta f_2 + \gamma f_3 ~\mbox{for some}~ \alpha, \beta, \gamma \in \mathbb{Z},
$$
where $f_i = e_i - e_0$ for $i = 1, 2, 3$. Using the multiplication table

\begin{center}
\begin{tabular}{|c||c|c|c|}
  \hline
 $\cdot$ & $f_1$ & $f_2$ & $f_3$ \\
   \hline
     \hline
$f_1$ & $f_1 - f_2 -f_3$ & $0$ & $f_1 - f_2 -f_3$ \\
  \hline
$f_2$  & $-2 f_2$ & $0$ & $-2 f_2$ \\
  \hline
$f_3$ & $-f_1 - f_2+f_3$ & $0$ & $-f_1 - f_2+f_3$ \\
  \hline
\end{tabular}\\
\end{center}
we obtain
$$
z^2 = (\alpha^2 - \gamma^2) f_1 - (\alpha^2 + 2 \alpha \beta + 2 \alpha \gamma + 2 \beta \gamma + \gamma^2)  f_2 + (-\alpha^2  + \gamma^2) f_3.
$$
Thus, $z^2=z$ if and only if
\begin{eqnarray*}
 \alpha &=& \alpha^2 - \gamma^2,\\
\beta &=& - (\alpha^2 + 2 \alpha \beta + 2 \alpha \gamma + 2 \beta \gamma + \gamma^2),\\
\gamma &=& -\alpha^2  + \gamma^2.
\end{eqnarray*}
Adding the first and the third equations gives $\alpha + \gamma = 0$. Then it follows from the system of equations that $\alpha = \beta = \gamma = 0$. Thus, $\mathbb{Z} [\R_4]$ does not have non-trivial idempotents with augmentation value 0.
\para

Case 2: $\varepsilon (z) = 1$, that is, $z = e_0 + \delta$, where $\delta \in \Delta (\R_4)$ and
$$
\delta = \alpha f_1 + \beta f_2 + \gamma f_3 ~\mbox{for some}~ \alpha, \beta, \gamma \in \mathbb{Z}.
$$
We have $z^2 = e_0 + \delta e_0 + e_0 \delta + \delta^2$. Since $\Delta (\R_4)$ is a two-sided ideal, we have $\delta e_0, e_0 \delta \in \Delta (\R_4)$. Using the formulas
$$
f_1 e_0 = f_3,~~f_2 e_0 = f_2,~~f_3 e_0 = f_1,~~e_0 f_1 = f_2,~~ e_0 f_2 = 0,~~ e_0 f_3 = f_2,
$$
we obtain
$$
\delta e_0 = \alpha f_3 + \beta f_2 + \gamma f_1 \quad \textrm{and} \quad e_0 \delta  = \alpha f_2 +  \gamma f_2.
$$
Using the expression for $\delta^2$ from Case 1 gives
$$
z^2 = e_0 + (\gamma + \alpha^2 - \gamma^2) f_1 + (\beta + \alpha + \gamma - \alpha^2 - \gamma^2 - 2 \alpha \gamma - 2 \alpha \beta - 2 \beta \gamma) f_2 + (\alpha - \alpha^2 + \gamma^2) f_3.
$$
Now, $z^2=z$ if and only if the system of equations
\begin{eqnarray*}
\alpha &=& \gamma + \alpha^2 - \gamma^2,\\
0 &=&  \alpha + \gamma - \alpha^2 - \gamma^2 - 2 \alpha \beta - 2 \alpha \gamma - 2 \beta \gamma,\\
\gamma &=& \alpha - \alpha^2  + \gamma^2,
\end{eqnarray*}
has integral solutions. The first equation has the form
$$
(\alpha - \gamma) = (\alpha - \gamma) (\alpha + \gamma).
$$
Suppose that $\alpha = \gamma$, then the second equation has the form $0 = \alpha (1 - 2 \alpha - 2 \beta)$. If $\alpha = 0$, then for arbitrary $\beta$, we have the idempotent $z = e_0 + \beta (e_2 - e_0)$. If $\alpha \not= 0$, then the second equation does not have solutions. Now, suppose that $\alpha \neq \gamma$, then $\gamma = 1 - \alpha$ and the second equation gives $\beta = 0$. Hence, for arbitrary $\alpha$, we have the idempotent $z = e_3 + \alpha(e_1 - e_3).$ $\blacksquare$
\end{proof}

The union construction for two quandles has a twisted version when the quandles act on each other by automorphisms (see Proposition \ref{new2}). We consider a simple case of this construction when both the quandles are trivial. Note that the automorphism group of a trivial quandle is the symmetric group on the underlying set. Let $X, Y$ be trivial quandles, $f \in \Aut(X)$ and $g \in \Aut(Y)$. For $x \in X$ and $y \in Y$,  setting 
$$x*y = f(x) \quad  \textrm{and} \quad y*x = g(y)$$
defines a quandle structure on the disjoint union $X \sqcup Y$, and we denote this quandle by $X \sqcup_{f, g}Y$. The following result gives the precise description of idempotents for $X \sqcup_{f, g}Y$ \cite[Proposition 6.3]{MR4565221}.

\begin{proposition}\label{idempotens in trivial unions}
Let $X$ and $Y$ be trivial quandles of orders $n$ and $m$, respectively. Let $\mathbb{k}$ be an integral domain such that characteristic of $\mathbb{k}$ is coprime to both $n$ and $m$. Let $f \in Aut(X)$ and $g \in Aut(Y)$ be automorphisms acting transitively on $X$ and $Y$, respectively. Then 
\begin{small}
$$\mathcal{I}\big(\mathbb{k}[X \sqcup_{f, g} Y]\big)= \mathcal{I}\big(\mathbb{k}[X]\big) \sqcup~ \mathcal{I}\big(\mathbb{k}[Y]\big) \sqcup \Big\{\alpha \big(\sum_{x \in X} e_x \big) + \beta \big(\sum_{y \in Y} e_y \big) ~\bigl\vert~ \alpha, \beta \in \mathbb{k}~~\textrm{~such that~}~~\alpha n+ \beta m=1 \Big\}.$$
\end{small}
\end{proposition}

\begin{proof}
Note that any $u \in \mathbb{k}[X \sqcup_{f, g} Y]$ can be written uniquely as $u=v+w$, where $v= \sum_{x \in X} \alpha_x e_x \in \mathbb{k}[X]$, $w= \sum_{y \in Y} \beta_y e_y \in \mathbb{k}[Y]$ and $\alpha_x, \beta_y \in \mathbb{k}$. If $u=u^2$, then
$$v+w=v^2+w^2+vw+wv= \varepsilon(v)v + \varepsilon(w)w + \varepsilon(w) \sum_{x \in X} \alpha_x e_{f(x)}+ \varepsilon(v)\sum_{y \in Y} \beta_y e_{g(y)},$$
and consequently
$$v= \varepsilon(v)v + \varepsilon(w) \sum_{x \in X} \alpha_x e_{f(x)} \quad \textrm{and} \quad w= \varepsilon(w)w+ \varepsilon(v)\sum_{y \in Y} \beta_y e_{g(y)}.$$
Comparing coefficients give
\begin{equation}\label{eq union 1}
\alpha_x= \varepsilon(v) \alpha_x+ \varepsilon(w) \alpha_{f^{-1}(x)}
\end{equation}
and
\begin{equation}\label{eq union 2}
\beta_y = \varepsilon(w) \beta_y + \varepsilon(v) \beta_{g^{-1}(y)}
\end{equation}
for all $x\in X$ and $y \in Y$. Adding \eqref{eq union 1} for all $x \in X$ gives $\varepsilon(v)= \varepsilon(v) \varepsilon(u)$. Similarly, adding 
\eqref{eq union 2} for all $y \in Y$ gives $\varepsilon(w)= \varepsilon(w) \varepsilon(u)$. If $\varepsilon(u)=0$, then $\varepsilon(v)=\varepsilon(w)=0$, and hence $u=0$, a contradiction. So, we can assume that $\varepsilon(u)=1$, and hence at least one of $\varepsilon(v)$ or $\varepsilon(w)$ is non-zero. If $\varepsilon(v) \neq 0$, then \eqref{eq union 2} gives $\beta_y = \beta_{g^{-1}(y)}$ for all $y \in Y$. Since $g$ acts transitively on $Y$, it follows that $\beta_y= \beta$ (say) for all $y \in Y$. If $\beta=0$, then $w=0$. In this case, $u=\sum_{x \in X} \alpha_x e_x$, where $\sum_{x \in X} \alpha_x =1$, and hence $u \in  \mathcal{I}\big(\mathbb{k}[X]\big)$. If $\beta \neq 0$, then $\varepsilon(w) = m\beta \neq 0$, and \eqref{eq union 1} gives $\alpha_x= \alpha_{f^{-1}(x)}$ for all $x \in X$. Since $f$ also acts transitively on $X$, it follows that $\alpha_x= \alpha$ (say) for all $x \in X$. Thus, we have $$u= \alpha \big(\sum_{x \in X} e_x \big) + \beta \big(\sum_{y \in Y} e_y \big),$$ where $n \alpha +m \beta=1$. Similarly, if $\varepsilon(w) \neq 0$ and $\alpha = 0$, then we get $v=0$. In this case, $u=\sum_{y \in Y} \beta_y e_y$, where $\sum_{y \in Y} \beta_y =1$, and hence $u \in \mathcal{I}\big(\mathbb{k}[Y]\big)$. This completes the proof.  $\blacksquare$
\end{proof}
\bigskip
\bigskip

\section{Quandle coverings and idempotents}\label{coverings and idempotents}
In this section, we use the idea of a quandle covering for giving a precise description of idempotents in quandle rings of quandles of finite type. The notion of a quandle covering is attributed to the work of Eisermann \cite{MR1954330, MR3205568}. 

\begin{definition}
    A quandle homomorphism $p: X  \to Y$ is called a \index{quandle covering}{\it quandle covering} if it is surjective and $S_x=S_{x'}$ whenever $p(x) = p(x')$ for any $x, x' \in X$. Clearly, an isomorphism of quandles is a quandle covering, called a \index{trivial covering}{\it trivial covering}. 

\end{definition}

\begin{example}{\rm 
Some basic examples of quandle coverings are as follows:
\begin{enumerate}
\item A surjective group homomorphism $p: G \to H$ yields a quandle covering $\Conj(G) \to \Conj(H)$ if and only if $\ker(p)$ is a central subgroup of $G$.
\item A surjective group homomorphism $p: G\to H$ yields a quandle covering $\Core(G) \to \Core(H)$ if and only if $\ker(p)$ is a central subgroup of  $G$ of exponent two.
\item Let $X$ be a quandle and  $F$ a non-empty set viewed as a trivial quandle. Consider $X \times F$ with the product quandle structure $(x,s) * (y,t) = (x *y,s)$. Then the projection $p: X \times F \to X$ given by $(x,s) \to x$ is a quandle covering.
\item Let $X$ be a quandle and $A$ an abelian group. Recall that, a quandle 2-cocycle is a map $\alpha: X \times X \to A$ satisfying
\begin{equation*} \label{group-coefficient-cocycle-condition}
\alpha_{x, y}~\alpha_{x*y, z}= \alpha_{x, z}~\alpha_{x* z, y*z}
\end{equation*}
and 
\begin{equation*}\label{normalised-cocycle-condition1}
\alpha_{x, x}=1
\end{equation*}
for $x, y, z \in X$. Given a 2-cocyle  $\alpha$, the set $X \times A$ turns into a quandle with the binary operation
\begin{equation*}\label{dynamical-quandle-operation}
(x, s)* (y,t)= \big( x* y, ~s~\alpha(x, y) \big),
\end{equation*} 
for $x, y \in X$ and $s, t \in A$. The quandle so obtained is called an  {\it extension} of $X$ by $A$ through $\alpha$, and is denoted by $X \times_{\alpha} A$. A direct check shows that the projection $p:X \times_{\alpha} A \to X$ given by $p(x, s)=x$ is a quandle covering.
\end{enumerate}}
\end{example}

The following lemma summarises some basic properties of quandle coverings, which we shall use without stating explicitly.

\begin{lemma}\label{quandle covering properties}
Let $p: X  \to Y$ be a quandle covering. Then the following  assertions hold:
\begin{enumerate}
\item Each fibre $p^{-1}(y)$ is a trivial subquandle of $X$. 
\item Each inner automorphism of $X$ permutes the fibres.
\item The fibres over any two elements of the same connected component of $Y$ are isomorphic.
\end{enumerate}
\end{lemma}

\begin{proof}
If $p: X  \to Y$ is a quandle homomorphism, then each fibre $p^{-1}(y)$ is a subquandle of $X$. Since $p$ is a covering, $S_x=S_{x'}$ whenever $x, x' \in p^{-1}(y)$. This gives $x*x'=S_{x'}(x)=S_x(x)=x$ and $x'*x=S_{x}(x')=S_{x'}(x')=x'$, which proves assertion (1).
\para
For assertion (2), it is enough to check that if $x_1, x_2 \in p^{-1}(y)$, then $S_x(x_1)$ and $S_x(x_2)$ are in the same fibre. Indeed, $p(S_x(x_1))=p(x_1*x)=
y*p(x)=p(x_2*x)=p(S_x(x_2))$, and we are done.
\para
Let $y, y'$ be elements of the same connected  component of $Y$. Then there exists elements $y_1, y_2, \ldots, y_n \in Y$ and $\mu_1, \mu_2, \ldots, \mu_n \in \{1, -1 \}$ such that $y'=y *^{\mu_1}y_1*^{\mu_2}y_2 \cdots *^{\mu_n}y_n$. Here the parentheses are left normalised. For each $i$,  choose one element $x_i \in p^{-1}(y_i)$. If $x \in p^{-1}(y)$, then we see that 
$$p(x*^{\mu_1} x_1*^{\mu_2} x_2 \cdots *^{\mu_n} x_n)= y*^{\mu_1}y_1*^{\mu_2}y_2 \cdots *^{\mu_n}y_n= y'.$$
Thus, the inner automorphism $S_{x_n}^{\mu_n}S_{x_{n-1}}^{\mu_{n-1}} \cdots S_{x_1}^{\mu_1}$ maps the fibre $p^{-1}(y)$ bijectively onto  $p^{-1}(y')$, which proves assertion (3). $\blacksquare$
\end{proof}

\begin{proposition}\label{covering with non-trivial idempotents}
Let $p:X \to Y$ be a non-trivial quandle covering and $\mathbb{k}$ an integral domain with unity. Then $\mathbb{k}[X]$ has non-trivial idempotents.
\end{proposition}

\begin{proof}
Since $p$ is a non-trivial covering, there is at least one connected component of $Y$ such that $|p^{-1}(y)| \ge 2$ for all elements $y$ of that connected component. By assertion (1) of Lemma \ref{quandle covering properties}, $p^{-1}(y)$ is a trivial subquandle of $X$. The result now follows from Proposition \ref{trivial subquandle idempotent}. $\blacksquare$
\end{proof}

A \index{long knot}{\it long knot} $L$ is the image of a smooth embedding $\ell: \mathbb{R} \hookrightarrow \mathbb{R}^3$ such that $\ell(t) = (t, 0, 0)$ for all $t$ outside some compact interval. We consider long knots only up to isotopy with compact support. The closure of a long knot is a usual knot defined in the obvious way.

\begin{proposition}\label{idempotents long knots}
Let $L$ be a non-trivial long knot in $\mathbb{R}^3$ and $\mathbb{k}$ an integral domain with unity. Then $\mathbb{k}[Q(L)]$ has non-trivial idempotents.
 \end{proposition}
 
\begin{proof}
Let $L$ be a long knot and $K$  its corresponding closed knot defined in the obvious way. Let $Q(L)$ and $Q(K)$ be knot quandles of $L$ and $K$, respectively. Note that $Q(K)$ is obtained from $Q(L)$ by adjoining one extra relation corresponding to the first and the last arc of $L$. By \cite[Theorem 35]{MR1954330}, the natural projection $p: Q(L) \to Q(K)$ is a non-trivial quandle covering, and the result follows from Proposition \ref{covering with non-trivial idempotents}. $\blacksquare$
 \end{proof}
 
Let $p:X \to Y$ be a quandle covering, and $\mathcal{F}(Y)$ the set of all finite subsets of $Y$. For each $y \in Y$, let $\mathcal{F}(p^{-1}(y))$ be the set of all finite subsets of $p^{-1}(y)$, and denote a typical element of this set by $I_y$. Given elements $x,y$ in a quandle $X$ of finite type, we set
$$[e_x]_y:= e_x+ e_{x*y}+ e_{x*y*y}+ \cdots +e_{x\underbrace{*y*y*\cdots *y}_{(n_y-1)-\textrm{times}}}=e_x+ e_{S_y(x)}+ e_{S_y^2(x)}+ \cdots +e_{S_y^{n_y-1}(x)},$$
the sum of the basis elements in the $S_y$-orbit of $e_x$. We have the following result \cite[Theorem 4.5]{MR4565221}.

\begin{theorem}\label{quandle covering theorem}
Let $X$ be a quandle of finite type, $p:X \to Y$ a non-trivial quandle covering and $\mathbb{k}$ an integral domain with unity. If $\mathbb{k}[Y]$ has only trivial idempotents, then the set of idempotents of $\mathbb{k}[X]$ is
\begin{small}
\begin{eqnarray}\label{quandle covering theorem idempotents}
\nonumber \mathcal{I} \big(\mathbb{k}[X]\big) &=& \Big\{ \sum_{y \in J} \Big( \sum_{x \in I_y,~~\sum \alpha_x =0}  \alpha_x~ [e_x]_{x_0} \Big) + \Big(\sum_{x' \in I_{y_0},~~\sum \alpha_{x'} =1} \alpha_{x'} ~e_{x'}\Big) ~\bigl\vert \\
 && J \in \mathcal{F}(Y),~~I_y \in \mathcal{F}(p^{-1}(y)),~~I_{y_0} \in \mathcal{F}(p^{-1}(y_0)),~~ x_0 \in  I_{y_0}, ~~y_0 \in Y, ~~\alpha_x, \alpha_{x'}  \in \mathbb{k} \Big\}.
\end{eqnarray}
\end{small}
\end{theorem}
 
\begin{proof}
Since $p$ is a quandle covering, we have $S_x=S_{x'}$ for any $x, x ' \in p^{-1}(y)$. Hence, the automorphisms of the quandle ring $\mathbb{k}[X]$ induced by $S_x$ and $S_{x'}$ are identical for any $x, x ' \in p^{-1}(y)$. This together with direct computations give
\begin{eqnarray}\label{even dihedral eq1}
&&\Big(\sum_{x \in J} \beta_{x}~ e_{x}\Big)\Big(\sum_{x' \in I_y,~~\sum \alpha_{x'} =1} \alpha_{x'} ~e_{x'}\Big)\\
\nonumber &=& \sum_{x' \in I_y,~~\sum \alpha_{x'} =1}\alpha_{x'} \Big( \sum_{x \in J} \beta_{x}~ e_{x} \Big) e_{x'}\\
\nonumber &=&\sum_{x' \in I_y,~~\sum \alpha_{x'} =1} \alpha_{x'} \Big( \sum_{x \in J} \beta_{x}~ e_{x} \Big)e_{x_0},\quad \textrm{for any fixed}~x_0\in I_y\\
\nonumber &=&\sum_{x' \in I_y,~~\sum \alpha_{x'} =1} \alpha_{x'} \Big( \sum_{x \in J} \beta_{x} ~e_{x*x_0} \Big)\\
\nonumber  &=&   \sum_{x \in J} \beta_{x}~ e_{x*x_0},\quad \textrm{since}~\sum_{x' \in I_y,~~\sum \alpha_{x'} =1} \alpha_{x'} =1,
\end{eqnarray}
and
\begin{eqnarray}\label{even dihedral eq2}
&&  \Big(\sum_{x \in J} \beta_{x} ~e_{x}\Big) \Big( \sum_{x' \in I_y,~~\sum \alpha_{x'} =0}  \alpha_{x'} ~e_{x'} \Big)\\
\nonumber &=&  \sum_{x' \in I_y,~~\sum \alpha_{x'} =0}  \alpha_{x'} \Big(\sum_{x \in J} \beta_{x} ~e_{x}\Big) e_{x'}\\
\nonumber &=&  \sum_{x' \in I_y,~~\sum \alpha_{x'} =0}  \alpha_{x'}\Big(\sum_{x \in J} \beta_{x}~ e_{x*x_0}\Big), \quad \textrm{for any fixed}~x_0 \in I_y\\
\nonumber &=& \Big(\sum_{x' \in I_y,~~\sum \alpha_{x'} =0}  \alpha_{x'} \Big) \Big(\sum_{x \in K} \beta_{x} ~e_{x*x_0}\Big)\\
\nonumber &=& 0,\quad \textrm{since}~ \sum_{x' \in I_y,~~\sum \alpha_{x'} =0}  \alpha_{x'}=0,
\end{eqnarray}
where $J \in \mathcal{F}(X)$, $I_y \in \mathcal{F}(p^{-1}(y))$,  $y \in Y$ and $\beta_x, \alpha_{x'}  \in \mathbb{k}$. Let $u=v+w$, where 
\begin{eqnarray*}
v &=&\sum_{y \in J} \Big( \sum_{x \in I_y,~~\sum \alpha_x =0}  \alpha_x ~[e_x]_{x_0}\Big),\\
w &=& \sum_{x' \in I_{y_0},~~\sum \alpha_{x'} =1} \alpha_{x'} ~e_{x'},\\
\end{eqnarray*}
$J \in \mathcal{F}(Y)$  and $x_0 \in  I_{y_0}$ a fixed element. Equations \eqref{even dihedral eq1} and \eqref{even dihedral eq2} imply that $w^2=w$, $wv=0$ and  $v^2=0$. By definition of the element $[e_x]_{x_0}$, it follows that $\big([e_x]_{x_0}\big)e_{x_0}=[e_x]_{x_0}$. Consequently, $vw=v$, and hence $u^2=u$. 
\para
For the converse, let $u$ be a non-zero idempotent of $\mathbb{k}[X]$. Since $X$ is the disjoint union of fibres of $p$, we can write $u$ uniquely in the form
$$u=\sum_{y \in J} \Big( \sum_{x \in I_y}  \alpha_x ~e_x \Big)$$
for some $J \in \mathcal{F}(Y)$ and $I_y \in \mathcal{F}(p^{-1}(y))$ for each $y \in J$. If $\hat{p}: \mathbb{k}[X] \to \mathbb{k}[Y]$ is the induced homomorphism of rings, then $\hat{p}(u)$ is an idempotent of $\mathbb{k}[Y]$. It follows from the decomposition of $u$ that
$$\hat{p}(u)=  \sum_{y \in J} \Big( \sum_{x \in I_y}  \alpha_x \Big) e_y.$$
Since $\mathbb{k}[Y]$ has only trivial idempotents, it follows that either $\hat{p}(u)=0$ or precisely one of the coefficients of $\hat{p}(u)$ is 1 and all other coefficients are 0. If $\hat{p}(u)=0$, then $\sum_{x \in I_y}  \alpha_x =0$ for each $y \in J$. Writing 
$$u=\sum_{y \in J} \Big( \sum_{x \in I_y, ~~\sum \alpha_{x} =0}  \alpha_x ~e_x \Big),$$
 it follows from \eqref{even dihedral eq2} that $u=u^2=0$, which is a contradiction as $u \ne 0$. Hence, there exists $y_0 \in J$ such that $\sum_{x' \in I_{y_0}}  \alpha_{x'}=1$ and $\sum_{x \in I_y}  \alpha_x=0$ for all $y \neq y_0$. We can write $u= v+w$, where 
$$v= \sum_{y \in J, ~y \neq y_0} \Big( \sum_{x \in I_y,~~\sum \alpha_x =0}  \alpha_x~ e_x \Big)$$ and $$w= \sum_{x' \in I_{y_0},~~\sum \alpha_{x'} =1} \alpha_{x'}~ e_{x'}.$$ Again, equations \eqref{even dihedral eq1} and \eqref{even dihedral eq2} imply that $w^2=w$, $wv=0$ and  $v^2=0$. Thus, we have
$$u=u^2=v^2+w^2+vw+wv= w+vw,$$
and consequently $v=vw$. This implies that 
$$\sum_{y \in J, ~y \neq y_0} \Big( \sum_{x \in I_y,~~\sum \alpha_x =0}  \alpha_x~ e_x \Big)= \sum_{y \in J, ~y \neq y_0} \Big( \sum_{x \in I_y,~~\sum \alpha_{x} =0}  \alpha_{x}~ e_{x} \Big) e_{x_0}$$
for some fixed $x_0 \in I_{y_0}$. Hence, it follows that $v$ has the form
$$v=\sum_{y \in J', ~y \neq y_0} \Big( \sum_{x \in I_y,~~\sum \alpha_x =0}  \alpha_x ~[e_x]_{x_0} \Big),$$
where $J' \subseteq J$ is a set of representatives of orbits of the action of $S_y$ on $J$. This completes the proof of the theorem. $\blacksquare$
 \end{proof}
 
 \begin{corollary}
If $X$ is a trivial quandle and $\mathbb{k}$ an integral domain with unity, then 
$$\mathcal{I}\big(\mathbb{k}[X]\big) =\Big\{ \sum_{x \in J} \alpha_x ~e_x  \, \mid \,    J \in \mathcal{F}(X)~\textrm{such that}~ \sum_{x \in J} \alpha_x=1 \Big\}.$$
 \end{corollary}
 \begin{proof}
If $\{z \}$ is a one element quandle, then the constant map $c: X \to \{z \}$ is a quandle covering. The proof now follows from Theorem \ref{quandle covering theorem}. $\blacksquare$
 \end{proof}

\begin{corollary}
Let $p:X \to Y$ be a non-trivial quandle covering and $\mathbb{k}$ an integral domain with unity. If $\mathbb{k}[Y]$ has only trivial idempotents, then every idempotent of $\mathbb{k}[X]$ has augmentation value 1.
\end{corollary}

\begin{proof}
The assertion follows from the proof of the converse part of Theorem \ref{quandle covering theorem}. Note that we do not need our quandles to be of finite type. $\blacksquare$
\end{proof}

The preceding result together with computer-assisted computations for quandles of order less than six suggests the following conjecture \cite[Conjecture 6.2]{MR4579329}.

\begin{conjecture}
Any non-zero idempotent of the integral quandle ring of a quandle has augmentation value 1.
\end{conjecture}

By Proposition \ref{idempotents-r3},  the integral quandle ring of $\R_3$ has only trivial idempotents. A computer assisted check shows that the same assertion holds for the integral quandle ring of $\R_5$ as well \cite[Section 6]{MR4579329}. As an application of the preceding theorem, we characterise idempotents in quandle rings of certain dihedral quandles of even order under the assumption of Conjecture \ref{Main conjecture}.

\begin{corollary}\label{even dihedral corollary}
Let $n=2m+1$ be an odd integer with $m \ge 1$ and $\mathbb{k}$ an integral domain with unity. Assume that $\mathbb{k}[\R_{n}]$ has only trivial idempotents. Then the set of idempotents of $\mathbb{k}[\R_{2n}]$ is given by
\begin{small}
$$\mathcal{I}\big(\mathbb{k}[\R_{2n}]\big) =\Big\{ \big(\beta e_j + (1-\beta)e_{n+j}\big)+ \sum_{i=0}^{m} \alpha_i \big(e_i-e_{n+i} + e_{2j-i}- e_{n+2j-i}\big) ~\bigl\vert ~ 0 \le j \le n-1 ~\textrm{and~} ~~ \alpha_i, \beta \in \mathbb{k} \Big\}.$$
\end{small}
 \end{corollary}
 
 \begin{proof}
Note that the natural map $p:\R_{2n} \to \R_{n}$ given by the reduction modulo $n$ is a two-fold non-trivial quandle covering. Further, for each $i \in \R_n$,  we have $p^{-1}(i)=\{i, n+i\}$. The result now follows from Theorem \ref{quandle covering theorem}. $\blacksquare$
 \end{proof}

\begin{proposition}\label{zero-divisors for coverings}
Let $p:X \to Y$ be a non-trivial quandle covering and $\mathbb{k}$ an integral domain with unity. Then $\mathbb{k}[X]$ has right zero-divisors.
 \end{proposition}
 
\begin{proof}
Let $J \in \mathcal{F}(X)$, $y \in Y$ and $I_y \in \mathcal{F}(p^{-1}(y))$ such that $|I_y| \ge 2$. Then for any $\sum_{x \in I_y,~~\sum \alpha_x =0}  \alpha_x e_x$ and $\sum_{x' \in J} \beta_{x'} e_{x'}$, it follows from \eqref{even dihedral eq2} that
$$\Big(\sum_{x' \in J} \beta_{x'} ~e_{x'}\Big) \Big( \sum_{x \in I_y,~~\sum \alpha_x =0}  \alpha_x ~e_x \Big)=0$$
and hence  $\sum_{x \in I_y,~~\sum \alpha_x =0}  \alpha_x e_x$ is a right zero-divisor of $k[X]$. $\blacksquare$
 \end{proof}

\begin{proposition}
Let $X$ be an involutory quandle and $\mathbb{k}$ an integral domain with unity. Let $A$ be a non-trivial abelian group and  $\alpha: X \times X \to A$ a quandle 2-cocycle.
\begin{enumerate}
\item  If $\alpha$ satisfy $\alpha_{x*y, y}=\alpha_{x, y}^{-1}$ for all $x, y \in X$, then the extension $X \times_{\alpha} A$ is involutory. 
\item If $\mathbb{k}[X]$ has only trivial idempotents, then $\mathbb{k}[X \times_{\alpha} A]$ has non-trivial idempotents.
\end{enumerate}
\end{proposition}

\begin{proof}
A direct check shows that the condition $\alpha_{x*y, y}=\alpha_{x, y}^{-1}$ is equivalent to the extension $X \times_{\alpha} A $ being involutory. Since the map $p:X \times_{\alpha} A \to X$ is a non-trivial quandle covering, the result follows from Theorem \ref{quandle covering theorem}. In fact, Theorem \ref{quandle covering theorem} gives the precise set of idempotents.  $\blacksquare$
\end{proof}

It turns out that the set of all idempotents of a quandle ring has the structure of a quandle in some cases \cite[Proposition 4.11]{MR4565221}.

\begin{proposition}\label{idempotents forming quandle covering}
Let $X$ be a quandle of finite type, $p : X \to Y$ a quandle covering and  $\mathbb{k}$ an integral domain with unity. Then the set 
\begin{eqnarray}
\nonumber  \mathcal{I} \big(\mathbb{k}[X]\big)  &=& \Big\{ \sum_{y \in J} \Big( \sum_{x \in I_y,~~\sum \alpha_x =0}  \alpha_x ~[e_x]_{x_0} \Big) + \Big(\sum_{x' \in I_{y_0},~~\sum \alpha_{x'} =1} \alpha_{x'}~ e_{x'}\Big) ~\bigl\vert \\
\nonumber && J \in \mathcal{F}(Y),~~I_y \in \mathcal{F}(p^{-1}(y)),~~I_{y_0} \in \mathcal{F}(p^{-1}(y_0)),~~ x_0 \in  I_{y_0}, ~~y_0 \in Y,~~\alpha_x, \alpha_{x'}  \in \mathbb{k}\Big\}
\end{eqnarray}
of idempotents of $\mathbb{k}[X]$ is a quandle with respect to the ring multiplication.
\end{proposition}

\begin{proof}
Consider the elements 
\begin{eqnarray*}
u &=& \sum_{y \in J_1} \Big( \sum_{x \in I_y,~~\sum \alpha_x =0}  \alpha_x ~ [e_x]_{x_1}  \Big) + \Big(\sum_{x' \in I_{y_1},~~\sum \alpha_{x'} =1} \alpha_{x'}~ e_{x'}\Big),\\
v &=&\sum_{y \in J_2} \Big( \sum_{x \in I_y,~~\sum \beta_x =0}  \beta_x ~[e_x]_{x_2}  \Big) + \Big(\sum_{x' \in I_{y_2},~~\sum \beta_{x'} =1} \beta_{x'}~ e_{x'}\Big)
\end{eqnarray*}
in $ \mathcal{I} \big(\mathbb{k}[X]\big) $, where $J_i \in \mathcal{F}(Y)$ and $I_{y} \in \mathcal{F}(p^{-1}(y))$, $y_i \in Y$ and $x_i \in I_{y_i}$. Then we have
\begin{eqnarray*}
& uv&\\
 &=& \Big(\sum_{y \in J_1} \Big( \sum_{x \in I_y,~~\sum \alpha_x =0}  \alpha_x ~ [e_x]_{x_1} \Big) + \Big(\sum_{x' \in I_{y_1},~~\sum \alpha_{x'} =1} \alpha_{x'}~ e_{x'}\Big) \Big) \\
&& \Big( \sum_{y \in J_2} \Big( \sum_{x \in I_y,~~\sum \beta_x =0}  \beta_x ~[e_x]_{x_2}  \Big) + \Big(\sum_{x' \in I_{y_2},~~\sum \beta_{x'} =1} \beta_{x'}~ e_{x'}\Big)\Big) \\
&=& \Big(\sum_{y \in J_1} \Big( \sum_{x \in I_y,~~\sum \alpha_x =0}  \alpha_x ~[e_x]_{x_1} \Big) + \Big(\sum_{x' \in I_{y_1},~~\sum \alpha_{x'} =1} \alpha_{x'} ~e_{x'}\Big) \Big) \Big(\sum_{x' \in I_{y_2},~~\sum \beta_{x'} =1} \beta_{x'}~ e_{x'}\Big),\\
&& \quad \textrm{by}~ \eqref{even dihedral eq2} \\
&=& \Big(\sum_{y \in J_1} \Big( \sum_{x \in I_y,~~\sum \alpha_x =0}  \alpha_x ~[e_x]_{x_1} \Big) + \Big(\sum_{x' \in I_{y_1},~~\sum \alpha_{x'} =1} \alpha_{x'} ~e_{x'}\Big) \Big) e_{x_2}, \quad \textrm{by} ~\eqref{even dihedral eq1}\\
&=& \sum_{y \in J_1} \Big( \sum_{x \in I_y,~~\sum \alpha_x =0}  \alpha_x ~\big([e_x]_{x_1}\big)e_{x_2} \Big) + \Big(\sum_{x' \in I_{y_1},~~\sum \alpha_{x'} =1} \alpha_{x'} ~ e_{x'*x_2}\Big),\\
&=& \sum_{y \in J_1} \Big( \sum_{x*x_2 \in I_{y*y_2},~~\sum \alpha_x =0}  \alpha_x~[e_{x*x_2}]_{x_1*x_2} \Big) + \Big(\sum_{x'*x_2 \in I_{y_1*y_2},~~\sum \alpha_{x'} =1} \alpha_{x'} ~e_{x'*x_2}\Big),
\end{eqnarray*}
since $\big([e_x]_{x_1}\big)e_{x_2}=[e_{x*x_2}]_{x_1*x_2}$ due to right-distributivity. Note here that $x_1*x_2 \in I_{y_1*y_2}$. Thus, we have proved that $uv \in  \mathcal{I} \big(\mathbb{k}[X]\big) $. The preceding computation also shows that the right multiplication by $v$ is precisely the right multiplication by 
$e_{x_2}$ for any fixed $x_2 \in I_{y_2}$. In other words, the right multiplication by $v$ is the ring automorphism $\hat{S}_{x_2}$ of $\mathbb{k}[X]$. This proves that the set $ \mathcal{I} \big(\mathbb{k}[X]\big) $ is a quandle. $\blacksquare$
\end{proof}
 
An immediate consequence of Proposition \ref{idempotents forming quandle covering} is the following result.
 
 \begin{corollary}\label{idempotens form quandle}
 Let $X$ be a quandle of finite type, $p : X \to Y$ a quandle covering and  $\mathbb{k}$ an integral domain with unity. Suppose that $\mathbb{k}[Y]$ has only trivial idempotents. Then the following  assertions hold:
\begin{enumerate}
\item The set of all idempotents of $\mathbb{k}[X]$ is  a quandle with respect to the ring multiplication.
\item The right multiplication by each idempotent of $\mathbb{k}[X]$ is a ring automorphism induced by some trivial idempotent of $\mathbb{k}[X]$.
\end{enumerate}
 \end{corollary}
\bigskip
\bigskip

\section{Idempotents of free products}\label{sec idempotents in free products}
Let us recall the construction of free product of quandles. Let $X_i = \langle S_i \, \mid \, R_i \rangle$ be a collection of $n \ge 2$ quandles given in terms of presentations. Then their {\it free product} $X_1   \star \cdots \star X_n$ is the quandle defined by the presentation
$$
X_1  \star \cdots \star X_n = \langle S_1 \sqcup \cdots \sqcup S_n  \, \mid \, R_1 \sqcup \cdots \sqcup R_n\rangle.
$$
For example, the free quandle $FQ_n$ of rank $n$ can be seen as 
$$
FQ_n = \langle x_1\rangle  \star \cdots \star\langle x_n\rangle,
$$
the free product of $n$ copies of trivial one element quandles $\langle x_i\rangle$.
\para

Recall that, each element of a quandle $X$ has a canonical left-associated expression of the form $x_0*^{\epsilon_1}x_1*^{\epsilon_2}\cdots*^{\epsilon_n}x_n$, where $x_0\neq x_1$, and if $x_i=x_{i+1}$ for any $1 \le i \le n-1$, then $\epsilon_i= \epsilon_{i+1}$.
\para
Lack of associativity in quandles makes it hard to have a normal form for elements in free products of quandles. We overcome this difficulty by defining a length for elements in free products.  Let $X= X_1 \star \cdots \star X_n$ be the free product of $n \ge 2$ quandles. Given an element $w \in X$, we define the length $\ell(w)$ of $w$ as 
\begin{eqnarray*}
\ell(w) &=& \min \Big\{r \, \mid \, w ~\textrm{can be written as a canonical left associated product of} ~r\\
&& \quad \quad  \textrm{elements from}~ X_1 \sqcup \cdots \sqcup X_n \Big\}.
\end{eqnarray*}
Note that each $w \in X$ has a reduced left associated expression attaining the length $\ell(w)$. This can be done by gathering together all the leftmost alphabets in a left associated expression of $w$ that lie in the same component quandle $X_i$, and rename it as a single element of $X_i$. This shows that $\ell(w)=1$ if and only if $w \in X_i$ for some $i$. Equivalently,  $\ell(w) \ge 2$ if and only if $w \in X\setminus (\sqcup_{s=1}^n X_s)$. 
\para
For example, if $x_1, x_2 \in X_i$ and $y_1, y_2 \in X_j$ for $i \ne j$, then $\ell(x_1*x_2)=1$, $\ell(x_1*x_2*^{-1}x_1)=1$, $\ell(x_1*y_1)=2$, $\ell(x_1*y_1*y_2)=3$ and $\ell(x_1*y_1*y_2*x_2)=4$.
\para
Note that, if $X= X_1 \star \cdots \star X_n$, then every $u \in \mathbb{k}[X]$ can be written uniquely in the form 
\begin{equation}\label{elements in rings of free products}
u= u_1+u_2+ \cdots + u_n + v,
\end{equation} 
where each $u_i \in \mathbb{k}[X_i]$,~ $v = \sum_{k=1}^m\gamma_k e_{w_k}$ with each $\ell(w_k) \ge 2$ and $\gamma_k \in \mathbb{k}$.

\begin{proposition}
Let $X= X_1 \star \cdots \star X_n$ be the free product of $n$ quandles and $\mathbb{k}$ an integral domain with unity. If each $\mathbb{k}[X_i]$ has only trivial idempotents, then any idempotent $u$ of $\mathbb{k}[X]$ can be written uniquely as
$$u=\alpha_1 e_{x_1} + \alpha_2 e_{x_2} + \cdots +\alpha_n e_{x_n} + v,$$
where $x_i \in X_i$, $v=\sum_{k=1}^m\gamma_k e_{w_k}$ with $\ell(w_k) \ge 2$ and $\alpha_i, \gamma_k \in \mathbb{k}$ for all $i$ and $k$.
\end{proposition}

\begin{proof}
For each $i$, fix an element $z_i \in X_i$. Then the maps $p_i: X \to X_i$ defined by setting 
$$
p_i(x)= \left\{
\begin{array}{ll}
x \quad \textrm{if}~~ x \in X_i, &  \\
z_i \quad \textrm{if}~~ x \in X_j  ~~\textrm{for}~~j \ne i.
\end{array} \right.
$$
The universal property of free products implies that each $p_i$ is a quandle homomorphism. Let $u= u_1+u_2+ \cdots + u_n + v$ be an idempotent of $\mathbb{k}[X]$, where each $u_i \in \mathbb{k}[X_i]$,~ $v = \sum_{k=1}^m\gamma_k e_{w_k}$ with $\ell(w_k) \ge 2$ and $\gamma_k \in \mathbb{k}$. Then $\hat{p_i}(u)$ is an idempotent in $\mathbb{k}[X_i]$ for each $i$. Since each $\mathbb{k}[X_i]$ has only trivial idempotents and
$$\hat{p_i}(u)= u_i+ \sum_{j \ne i, ~j=1}^n \epsilon(u_j) e_{z_i} + \hat{p_i}(v),$$
it follows that $u_i=\alpha_i e_{x_i}$ for some $x_i \in X_i$ and $\alpha_i \in \mathbb{k}$. Note that if $ \epsilon(u_j) \neq 0$ for any $j \ne i$, then $x_i=z_i$. Thus, $u=\alpha_1 e_{x_1} + \alpha_2 e_{x_2} + \cdots +\alpha_n e_{x_n} + v$, and we are done. $\blacksquare$
\end{proof}

\begin{lemma}\label{free product lemma}
Let $X= X_1 \star \cdots \star X_n$ be the free product of $n$ quandles with $n \ge 2$ and $\mathbb{k}$ an integral domain with unity. Let $u \in \mathbb{k}[X]$ be an idempotent and
$u=u_1+u_2+ \cdots + u_n + v$ be its unique decomposition as in \eqref{elements in rings of free products}. Suppose that 
\begin{enumerate}
\item $\ell(w_k*w_l) \ge 2$ for any $k$ and $l$,
\item $\ell(w_k*x ) \ge 2$ for any $x \in \sqcup_{s=1}^n X_s$ and any $k$.
\end{enumerate}
Then $u_i$ is an idempotent of $\mathbb{k}[X_i]$ for each $i$.
\end{lemma}

\begin{proof}
Since $u=u^2$, we have
\begin{equation}\label{elements in rings of free products2}
u_1+u_2+ \cdots + u_n + v= u_1^2+ u_2^2 + \cdots + u_n^2 + v^2 + \sum_{i \neq j, ~i, j=1}^n u_i u_j +  \sum_{i=1}^n u_i v +  \sum_{j=1}^n v u_j.
\end{equation} 
If $v=0$, then \eqref{elements in rings of free products2} takes the form
\begin{equation}\label{elements in rings of free products3}
u_1+u_2+ \cdots + u_n= u_1^2+ u_2^2 + \cdots + u_n^2 + \sum_{i \neq j, ~i, j=1}^n u_i u_j.
\end{equation}
For $1 \le i \ne j \le n$, each basis element of $\mathbb{k}[X]$ appearing in a product $u_i u_j$ corresponds to a quandle element from $X \setminus (\sqcup_{s=1}^n X_s)$. For each $1 \le i \le n$, gathering all the summands on the right hand side of \eqref{elements in rings of free products3} corresponding to elements from the quandle $X_i$ implies that $u_i=u_i^2$, which is desired.
\para
Now, suppose that $v \neq 0$. For each $1 \le k, l\le m$, the condition $\ell(w_k*w_l) \ge 2$ implies that the basis element of $\mathbb{k}[X]$ corresponding to the quandle element $w_k*w_l$ does not appear as a summand for any $u_j$. Further, each basis element appearing in a product $u_i v$ corresponds to a quandle element of the form $x*w_k$ for some $x \in X_i$ and some $1 \le k \le m$. But, we have $\ell(x*w_k) \ge 2$ for such elements. Lastly, the condition $\ell(w_k*x ) \ge 2$ for any $x \in \sqcup_{s=1}^n X_s$ also implies that the basis element of $\mathbb{k}[X]$ corresponding to the quandle element $w_k*x$ does not appear as a summand for any $u_j$. For each $1 \le i \le n$, gathering together all the summands on the right hand side of \eqref{elements in rings of free products2} corresponding to elements from the quandle $X_i$ imply that $u_i=u_i^2$, which is desired. $\blacksquare$
\end{proof}

We are now ready to prove the main result of this section \cite[Theorem 5.3]{MR4565221}.

\begin{theorem}\label{idempotents in free products}
Let $FQ_n$ be the free quandle of rank $n \ge 1$. Then $\mathbb{Z}[FQ_n]$ has only trivial idempotents. The same assertion holds for the free quandle of countably infinite rank. 
\end{theorem}

\begin{proof}
By Theorem  \ref{IvanovKadantsevKuznetsov theorem}, every subquandle of a free quandle is free. Let $FQ_2= \langle x\rangle  \star \langle y \rangle$ be the free quandle of rank two. Note that $FQ_2$ can be identified with the subquandle of the conjugation quandle of the free group $F(x, y)$ on the set $\{x, y\}$. This together with the fact that the set $\{x, yx y^{-1}, y^2x y^{-2}, \ldots,  y^{n-1}x y^{-n+1} \}$ generates a free subgroup of $F(x, y)$ of rank $n$ imply that 
$$FQ_n \cong \langle x\rangle  \star \langle x*y \rangle \star \langle x*y*y \rangle \star \cdots \star \langle x*\underbrace{y*y* \cdots *y}_{(n-1)~\textrm{times}} \rangle,$$
and hence $FQ_n$ embeds as a subquandle of $FQ_2$ for each $n \ge 3$. Similarly, one can identify the free quandle of countably infinite rank as 
$$FQ_{\infty} \cong \langle x\rangle  \star \langle x*y \rangle \star \langle x*y*y \rangle \star \cdots \star \langle x*\underbrace{y*y* \cdots *y}_{(n-1)~\textrm{times}} \rangle \star \cdots ,$$
and hence it also embeds as a subquandle in $FQ_2$.   Thus, it suffices to prove that $\mathbb{Z}[FQ_2]$ has only trivial idempotents.
\para
Let  $u= \alpha e_x + \beta e_y +v$  be an idempotent of $\mathbb{Z}[FQ_2]$, where $v=\sum_{k=1}^m \gamma_k e_{w_k}$ with $\ell(w_k) \ge 2$ and $\alpha, \beta, \gamma_k \in \mathbb{Z}$. If $v=0$, then  Lemma \ref{free product lemma} implies that $\alpha e_x=  \alpha^2 e_x$ and $ \beta e_y=  \beta^2 e_y$. Hence, either $u=e_x$ or $u=e_y$, and $u$ is a trivial idempotent.
\para
Now, suppose that $v \neq 0$. Note that the first two alphabets (from left) in the reduced left associated expression of each $w_k$ are distinct. We claim that $\gamma_k=1$ for each $k$. This will be achieved by transforming the idempotent $u$ into a new idempotent such that conditions of Lemma \ref{free product lemma} are satisfied.  Fix a $k$ such that $1 \le k \le m$ and write
$$w_k= x_0 *^{\epsilon_1}x_1 *^{\epsilon_2} x_2 *^{\epsilon_3} \cdots  *^{\epsilon_r} x_r,$$
in its reduced left associated expression, where $x_i \in \{ x, y \}$ and $\epsilon_i \in \mathbb{Z}$ for each $i$. Since the expression is reduced, without loss of generality, we can assume that $x_0=x$ and $x_1=y$. Consider the inner automorphism 
$$\phi=S_{x_0}S_{x_0}S_{x_1}^{-\epsilon_1}S_{x_2}^{-\epsilon_2} \cdots S_{x_{r-1}}^{-\epsilon_{r-1}}S_{x_r}^{-\epsilon_r} $$
of $FQ_2$. We analyse the effect of $\phi$ on each summand of $u$. First note that $\phi(w_k)= x_0=x$. Consider any fixed $w_i$ for $i \neq k$ and write $w_i=y_0 *^{\mu_1}y_1 *^{\mu_2} y_2 *^{\mu_3} \cdots  *^{\mu_s} y_s$ in its reduced left associated expression, where $y_t \in \{ x, y \}$ and $\mu_t \in \mathbb{Z}$ for each $t$. We have
$$\phi(w_i)= y_0 *^{\mu_1}y_1 *^{\mu_2} y_2 *^{\mu_3} \cdots  *^{\mu_s} y_s *^{-\epsilon_r}  x_r *^{-\epsilon_{r-1}}  x_{r-1}*^{-\epsilon_{r-2}} \cdots *^{-\epsilon_1}x_1 * x_0 * x_0.$$
Considering the cases $s=r$, $s>r$ and $s<r$, and using the fact that the set of alphabets is $\{x, y \}$, we obtain $\ell(\phi(w_i)) \ge 3$. This clearly implies that $\phi(w_i)*x, \phi(w_i)*y \not\in \{x, y \}$ for any $i \ne k$. Now, consider another $w_j$ for $j \neq k$ and $j \neq i$ and write $w_j=z_0 *^{\nu_1}z_1 *^{\nu_2} z_2 *^{\nu_3} \cdots  *^{\nu_l} z_l$ in its reduced left associated expression, where $z_t \in \{ x, y \}$ and $\nu_t \in \mathbb{Z}$ for each $t$. Then Lemma 
\ref{Lem:the canonical left associated form} gives
\begin{eqnarray*}
&&\phi(w_i) *\phi(w_j)\\
&=&\phi(w_i*w_j)\\
&=&\big((y_0 *^{\mu_1}y_1 *^{\mu_2}  \cdots  *^{\mu_s} y_s)(z_0 *^{\nu_1}z_1 *^{\nu_2}  \cdots  *^{\nu_l} z_l) \big)\\
&&*^{-\epsilon_r}  x_r *^{-\epsilon_{r-1}}  x_{r-1}*^{-\epsilon_{r-2}} \cdots *^{-\epsilon_1}x_1 * x_0 * x_0 \\
&=& y_0 *^{\mu_1}y_1 *^{\mu_2}  \cdots  *^{\mu_s} y_s *^{-\nu_l} z_l *^{-\nu_{l-1}} z_{l-1} *^{-\nu_{l-2}} \cdots *^{-\nu_1}z_1 *z_0 *^{\nu_1}z_1 *^{\nu_2}  \cdots  *^{\nu_l} z_l\\
&&*^{-\epsilon_r}  x_r *^{-\epsilon_{r-1}}  x_{r-1}*^{-\epsilon_{r-2}} \cdots *^{-\epsilon_1}x_1 * x_0 * x_0.
\end{eqnarray*}
As before, by comparing $\ell(w_i*w_j)$ and $r$, we obtain $\ell \big(\phi(w_i) *\phi(w_j)\big) \ge 3$. If $\alpha$ and $\beta$ are non-zero, then $\ell(\phi(x)), \ell(\phi(y)) \ge 3$ for the same reason. Thus, the only summand of the idempotent $\phi(u)= \alpha e_{\phi(x)} + \beta e_{\phi(y)} + \sum_{k=1}^m \gamma_k e_{\phi(w_k)}$ that corresponds to an element from $\{x, y \}$ is $\phi(w_k)$, and all the summands  corresponding to $\phi(w_i)$ for $i \ne k$ satisfy the conditions of  Lemma \ref{free product lemma}. Thus, we obtain $\gamma_k e_{\phi(w_k)} =(\gamma_k e_{\phi(w_k)})^2$, and hence $\gamma_k=1$, which proves the claim. On plugging this information back to $u$, we can write $u= \alpha e_x + \beta e_y + \sum_{k=1}^m e_{w_k}$. Since $u$ is an idempotent, we have
\begin{eqnarray*}
\alpha e_x + \beta e_y + \sum_{k=1}^m e_{w_k} &=& \alpha^2 e_x + \beta^2 e_y + \sum_{k, ~l=1}^m e_{w_k *w_l} + \alpha \beta e_{x*y} +\alpha \beta e_{y*x}\\
&&+  \alpha \sum_{k=1}^m e_{x* w_k} + \alpha \sum_{k=1}^m e_{w_k* x} + \beta \sum_{k=1}^m e_{y* w_k} + \beta \sum_{k=1}^m e_{w_k* y}.
\end{eqnarray*}
Comparing coefficients of $e_x$ gives
$$\alpha= \alpha^2, \quad \alpha= \alpha^2 + \sum_{w_k *w_l=x} 1, \quad \alpha= \alpha^2 + \beta \quad \textrm{or} \quad \alpha= \alpha^2 + \sum_{w_k *w_l=x} 1 + \beta.$$
Similarly, comparing coefficients of $e_y$ gives
$$\beta = \beta ^2, \quad \beta = \beta ^2 + \sum_{w_k *w_l=y} 1, \quad \beta = \beta ^2 + \alpha \quad \textrm{or} \quad \beta = \beta ^2 + \sum_{w_k *w_l=y} 1 + \alpha.$$
Using the fact that the coefficients are from $\mathbb{Z}$, a direct check shows that the only possible cases are
\begin{eqnarray*}
\alpha= \alpha^2 &\textrm{and}& \beta=\beta ^2,\\
\alpha= \alpha^2 + \beta &\textrm{and}& \beta=\beta ^2,\\
\alpha= \alpha^2 &\textrm{and}& \beta=\beta ^2+ \alpha,\\
\alpha= \alpha^2 + \beta&\textrm{and}& \beta=\beta ^2+ \alpha.
\end{eqnarray*}
This together with the fact that $\varepsilon(u)=\alpha+ \beta +m$ shows that $\alpha=\beta=0$ and $m=1$. Hence, $u=e_w$ for some $w \in FQ_2$, and the proof is complete. $\blacksquare$
\end{proof}

Since the link quandle of a trivial link with $n$ components is the free quandle of rank $n$, we obtain the following result.

\begin{corollary}\label{idempotents in trivial link quandles}
If $L$ is a trivial link, then $\mathbb{Z}[Q(L)]$ has only trivial idempotents.
\end{corollary}

Let $\WB_n$ be the \index{welded braid group} {\it welded braid group} on $n$-strands. See \cite{MR3689901} for a nice survey of these groups.
Proposition \ref{auto quandle ring free products}, Theorem \ref{idempotents in free products} and \cite{MR1410467} yield the following result.

\begin{corollary}
$\Aut_{\textrm{Algebra}}\big(\mathbb{Z}[FQ_n]\big) \cong \Aut_{\mathcal{Q}}(FQ_n)\cong \WB_n$ for each $n \ge 1$.
\end{corollary}

We conjecture that the following statement holds.

\begin{conjecture}
Let $\{X_i\}_{i \in I}$ be a family of quandles such that each $\mathbb{Z}[X_i]$ has only trivial idempotents. Then $\mathbb{Z}[\star_{i \in I} X_i]$ admit only trivial idempotents. 
\end{conjecture}

\begin{remark}
{\rm In \cite{MR4579329}, understanding of idempotents in quandle rings has been applied to construct proper enhancements of the classical quandle coloring invariant of links in the 3-space.}
\end{remark}


\chapter*{Part III}
\begin{center}
{\huge{Homology and cohomology of solutions}}
\end{center}

\chapter{Trunks and rack spaces}\label{chapter trunks and rack space}

\begin{quote}
This chapter builds on Fenn, Rourke, and Sanderson's concept of the rack space, a topological analogue of a group's classifying space. It introduces key concepts such as trunks, cubical sets, and square sets.  The construction of the rack space of a rack is then presented, from which the quandle space of a quandle is derived. The chapter also explores homotopy types of selected rack spaces and their second homotopy groups, offering insights into the topological aspects of racks and quandles.
\end{quote}
\bigskip

\section{Trunks}\label{trun}

Through an exploration of racks, particularly the rack space, Fenn, Rourke, and Sanderson \cite{MR1364012} were prompted to explore the notion of a trunk. In broad terms, a trunk can be considered as an entity akin to a category. We begin the section by presenting the subsequent definition.

\begin{definition}
A \index{trunk}{trunk} $T$ consists of a class of objects such that there is a set $\Hom_T(A,B)$ of morphisms for each ordered pair $(A,B)$ of
objects, and $T$ admits certain preferred squares
	\begin{equation}\label{preferred}
		\begin{tikzcd}[row sep=huge]
			A \arrow[r, "f"] \arrow[d, "k"]  & B  \arrow[d, "g"]  \\
			C \arrow[r, "h"]   & D 
		\end{tikzcd}
	\end{equation}	
of morphisms (a notion analogous to that of composition in a category). 	
\end{definition}

Note that, a trunk can be thought of as a directed graph, together with a collection of oriented squares called preferred squares.  There is a natural notion of a \index{trunk map}{trunk map} (or a functor). 

\begin{definition}
 Given trunks $T_1$ and $T_2$, a \index{trunk map}{trunk map} $F: T_1 \rightarrow T_2$ is a map that assigns to each  object $S$ of $T_1$ an object $F(S)$ of $T_2$, and to each morphism $f:A \rightarrow B$ of $T_1$ a morphism $F(f): F(A) \rightarrow F(B)$ of $T_2$  such that $F$ preserves preferred squares. 
 \end{definition}
 
 In other words, when the functor $F$ is covariant, the preferred square~(\ref{preferred}) gives the preferred square
\begin{equation*}
	\begin{tikzcd}[row sep=huge]
		F(A) \arrow[r, "F(f)"] \arrow[d, "F(k)"]  & F(B)  \arrow[d, "F(g)"]  \\
		F(C) \arrow[r, "F(h)"]   & F(D), 
	\end{tikzcd}
\end{equation*}	
and when the functor $F$ is contravariant, the preferred square~(\ref{preferred}) yields the preferred square

\begin{equation*}
	\begin{tikzcd}[row sep=huge]
		F(A)   & F(B)  \arrow[l, "F(f)"]  \\
		F(C) \arrow[u, "F(k)"]   & F(D).   \arrow[l, "F(h)"]   \arrow[u, "F(g)"] 
	\end{tikzcd}
\end{equation*}

\begin{example}\label{trunk examples}
Let us consider some examples of trunks.
\begin{enumerate}
\item For each category $\mathcal{C}$,  there is a trunk $\mathcal{T}(\mathcal{C})$ which has the same objects and morphisms as $\mathcal{C}$, and whose preferred squares are the commutative diagrams in $\mathcal{C}$. 	
\item  Each rack $X$ defines a trunk $\mathcal{T}(X)$, called the \index{rack trunk}{rack trunk}, which has a single vertex $\bullet$ and has $X$ as its set of edges.  The preferred squares are of the  form	
\begin{equation*}
			\begin{tikzcd}[row sep=huge]
				\bullet \arrow[r, "x"] \arrow[d, "y"]  & 	\bullet \arrow[d, "y"]  \\
\bullet \arrow[r, "x*y"]   & 	\bullet,
\end{tikzcd}
\end{equation*}	
one for each ordered pair $(x,y)$ of elements of $X$. 
\item For each rack $X$, we define another trunk $\mathcal{T}'(X)$ as follows. For each element $x \in X$,  there is an object of $\mathcal{T}'(X)$ (also denoted by $x$); and for each ordered pair $(x,y)$ of elements of $X$, there is a morphism $\textbf{a}_{x,y}: x \rightarrow x*y$ and a morphism $\textbf{b}_{x,y}: y \rightarrow x*y$ such that the squares
		\begin{equation*}
			\begin{tikzcd}[row sep=huge]
				x \arrow[r, "\textbf{a}_{x,y}"] \arrow[d, "\textbf{a}_{x,z}"]  & 	x*y \arrow[d, "\textbf{a}_{x*y,z}"]  \\
				x*z \arrow[r, "\textbf{a}_{x*z,y*z}"]   & 	(x*y)*z=(x*z)*(y*z)
			\end{tikzcd}
		\end{equation*}	
		and 
		\begin{equation*}
			\begin{tikzcd}[row sep=huge]
				y \arrow[r, "\textbf{b}_{x,y}"] \arrow[d, "\textbf{a}_{x,z}"]  & 	x*y \arrow[d, "\textbf{a}_{x*y,z}"]  \\
				y*z \arrow[r, "\textbf{b}_{x*z,y*z}"]   & 	(x*y)*z=(x*z)*(y*z)
			\end{tikzcd}
		\end{equation*}	
are preferred squares for all $x,y,z \in X$.
\item The naturally oriented one-dimensional skeleton of the $n$-cube $I^n=[0,1]^n$ in $\mathbb{R}^n$ forms a trunk, where the preferred squares are the two-dimensional faces of $I^n$.
\end{enumerate}	
\end{example}

\begin{remark}\label{remark on trunk map}
{\rm 
Let $\mathcal{A}$ be the category of abelian groups. For each rack $X$, we define a trunk map $A: \mathcal{T}'(X) \to \mathcal{A}$ that assigns an abelian group $A_x$ to each $x \in X$ and assign  group homomorphisms $\eta_{x,y} \colon A_x \rightarrow A_{x*y}$ and $\xi_{x,y} \colon A_y \rightarrow A_{x*y}$ for each ordered pair $(x, y)$ of elements of $X$ such that 
$$	\eta_{x*y,z}\eta_{x,y}=\eta_{x*z,y*z} \eta_{x,z} \quad \textrm{and} \quad 
\eta_{x*y,z}\xi_{x,y} = \xi_{x*z,y*z} \eta_{x,z}
$$
for all $x,y,z \in X$.  Later, this trunk map will be used to define the notion of a quandle module in Chapter \ref{chapter homology YBE}.
}		
\end{remark}

\begin{definition}\label{corner}
A \index{corner trunk}{corner trunk} is a trunk that satisfies the following two corner axioms.
\begin{enumerate}
\item Given morphisms $f:A \rightarrow B$ and $k:A \rightarrow C$, there exist unique morphisms $g:B \rightarrow D$ and $h:C \rightarrow D$ such that the square 
		\begin{equation*}
			\begin{tikzcd}[row sep=huge]
				A \arrow[r, "f"] \arrow[d, "k"]  & B  \arrow[d, "g"]  \\
				C \arrow[r, "h"]   & D 
			\end{tikzcd}
		\end{equation*}	
is preferred.  
\item In the following diagram, if the center, south, and west squares are preferred, then it is possible to complete the diagram so that the east and north squares are also preferred.
		\begin{equation*}
			\begin{tikzcd}[row sep=huge]
				E \arrow{rrr} &&& F\\
				&A \arrow{r} \arrow{ul}& B \arrow{ur} &\\    
				&C  \arrow[r, "f"] \arrow[dl, "g"]  \arrow[u, "h"] &  D \arrow{dr}  \arrow{u}  & \\
				H   \arrow{rrr} \arrow{uuu}&&&G \arrow{uuu}
			\end{tikzcd}
		\end{equation*}	
In other words, assuming condition (1),  the morphisms $f$, $g$ and $h$ determine the entire diagram of preferred squares. 
\end{enumerate}
\end{definition}
\medskip

Next, we consider cubes and square sets. Let $I=[0,1]$ and $I^n=[0,1]^n$ be the $n$-cube in $\mathbb{R}^n$. A vertex of $I^n$ is a point of the form $(\epsilon_1, \ldots, \epsilon_n)$, where each $\epsilon_i=0$ or $1$.  A $p$-face of $I^n$ is a subset given by setting $(n-p)$ coordinates equal to either $0$ or $1$. For example, a $1$-face (also called an edge) is of the form 
$$\{ \epsilon_1 \} \times \cdots \times  \{ \epsilon_{j-1} \} 
\times I_j  \times \{ \epsilon_{j+1} \}  \times \cdots \times   \{ \epsilon_{n}\}$$
for some $j$ with $1 \leq j \leq n$,
where each $\epsilon_k = 0$ or $1$ and $I_j=I$. Similarly, a $2$-face has the form 
$$  \{ \epsilon_1 \} \times \cdots \times  \{ \epsilon_{i-1} \} 
\times I_i  \times \{ \epsilon_{i+1} \}  \times \cdots \times
\{ \epsilon_{j-1} \}  \times I_j  \times \{ \epsilon_{j+1} \}
\times \cdots \times  \{ \epsilon_{n}\}  $$
for some $i$ and $j$ with $1 \leq i<  j \leq n$, where each $\epsilon_k = 0$ or $1$ and $I_i=I_j=I$.  Faces can be given canonical orientations.  For simplicity of notation, we will think of a vertex of $I^n$ as a subset of $\{1, \ldots, n\}$ as follows: A subset $A \subseteq \{1, \ldots, n\}$ corresponds to a vertex with $\epsilon _i=1$ if and only if $i \in A$.  Then an edge has the form $A \rightarrow A\cup \{i\}$ with $i \notin A$.  
\para 

The $n$-cube $I^n$ can be regarded as a trunk (called the $n$-cube trunk) with the 2-faces
\begin{equation*}
	\begin{tikzcd}[row sep=huge]
		A \arrow[r] \arrow[d]  & 	A \cup \{i\} \arrow[d]  \\
		A \cup \{j\} \arrow[r]   & 	A \cup \{i\} \cup \{j\}
	\end{tikzcd}
\end{equation*}	
as preferred squares, where $1 \le i<j \le n$. In the following diagram, we see that if the center, south and west squares are preferred, then the diagram can be completed such that the east and north squares are also preferred.
\begin{equation*}
	\begin{tikzcd}[row sep=huge]
		B \cup \{i,k \}\arrow{rrr} &&&  B \cup \{i,j,k\}\\
		&B \cup \{k\} \arrow{r} \arrow{ul}& B  \cup \{j,k\} \arrow{ur} &\\    
		&B  \arrow[r, "f"] \arrow[dl, "g"]  \arrow[u, "h"] &  B \cup \{j\} \arrow{dr}  \arrow{u}  & \\
		B \cup \{i\}   \arrow{rrr} \arrow{uuu}&&&  B \cup \{i,j\} \arrow{uuu}.
	\end{tikzcd}
\end{equation*}	

Next, we provide the definition of the square category, which will be utilized subsequently in conjunction with trunks.

\begin{definition}
	The \index{square category}{square category} (~$\square$-category) is the category whose objects are the $n$-cubes $I^n$ for $n\ge 0$ and whose morphisms are the face maps.
\end{definition}

The face maps $\delta_i^\epsilon: I^{n-1} \rightarrow I^n$ corresponding to the $(n-1)$-faces of $I^n$ are given by $$ \delta_i^\epsilon(x_1, \ldots, x_{n-1})= (x_1, \ldots, x_{i-1},\epsilon, x_i, \ldots x_{n-1}),$$ where $\epsilon \in \{0,1\}$. It is not difficult to see that 

\begin{eqnarray}
	\delta_i^\epsilon \delta_{j-1}^\omega = \delta_{j}^\omega \delta_i^\epsilon
\end{eqnarray}
 for $1 \leq i<j\leq n$ and $\epsilon, \omega \in \{0,1\}$. Any face map $\lambda$ can be uniquely written as a composition $\delta_{i_k}^{\epsilon_k} \ldots \delta_{i_1}^{\epsilon_1}$, where each $\epsilon_j \in \{0,1\}$ and we assume that either $i_1<\cdots <i_k$ or $i_k \leq \cdots \leq i_1$.

\begin{definition}
A $\square$-set is a functor $C: \square^{\textrm{op}} \rightarrow \mathcal{S}$, where  $\square^{\textrm{op}}$ is the opposite category of the $\square$-category and $\mathcal{S}$ is the category of sets.  A $\square$-map between two $\square$-sets is a natural transformation of functors.  
\end{definition}

\begin{definition}
The \index{nerve}{nerve} of a trunk $T$ is the $\square$-set $NT$ satisfying the following conditions:
\begin{enumerate}
\item The set $NT_n$ is defined to be the set $\Hom(I^n,T)$ of trunk-maps from the cube trunk $I^n$ to the trunk $T$.  
\item The  map $\lambda^*$ is given by $\lambda^*(f)= f \circ \lambda$, where $\lambda \colon I^p \rightarrow I^n$ is a face map and $f \colon I^n \rightarrow T$ is a trunk map.
\end{enumerate}
\end{definition}

For an $n$-cube $I^n$ and a face map $\lambda$, we denote $C(I^n)$ by $C_n$ and  $C(\lambda)$ by $\lambda^*$.  We also denote $C(\delta_i^\epsilon)$ by  $\partial_i^\epsilon$. The geometric realization $|C|$ of a $\square$-set $C$ is defined to be the quotient space of the disjoint union $\bigsqcup_{n \geq 0} \; \big(C_n \times I^n \big)$ by the equivalence relation generated by
$$\big(\lambda^*(x),t \big) \sim \big(x, \lambda(t) \big),$$ 
where $x \in C_n$ and $t \in I^n$.  Note that $|C|$ is a CW-complex with one $n$-cell for each element of $C_n$ and  each $n$-cell has a canonical characteristic map from the $n$-cube.

\begin{remark}
{\rm A trunk has a natural cubical nerve, analogous to the simplicial nerve of a category. The classifying space of the trunk is the geometric realization of this nerve.   The rack space is the realization of a semi-cubical complex, and the natural structure is cubical, rather than simplicial as in the case of a classifying space. 
}    
\end{remark}
\bigskip


\section{Rack and quandle spaces}\label{Rspace}
The rack space $BX$ for a rack $X$ is the analogue of the classifying space $BG$ for a group $G$.  It defines a functor from the category of racks to the category of topological spaces.  Following \cite[Section 3]{MR1364012}, we have the following definition.

\begin{definition}
Let $X$ be a rack and $\mathcal{T}(X)$ the corner trunk associated to $X$ (Example \ref{trunk examples}(2)).  Then the geometric realization of the nerve $N\mathcal{T}(X)$ of $\mathcal{T}(X)$ is called the \index{rack space}{rack space} for $X$ and is denoted by $BX$.
\end{definition}

If $X$ is a rack, then we have $N\mathcal{T}(X)_n = X^n$,  the $n$-fold cartesian product of $X$ with itself. Further, the face maps are given by
$$\partial_i^0(x_1, \ldots, x_n)= (x_1, \ldots, x_{i-1}, x_{i+1}, \ldots, x_n)$$
and
$$\partial_i^1(x_1, \ldots, x_n)= (x_1*x_i, \ldots, x_{i-1}*x_i, x_{i+1}, \ldots, x_n)$$
for all $1 \le i \le n$. Then the rack space $BX$ consists of the following:
	
\begin{enumerate}
\item  One $0$-cell $\bullet$.  
\item One $1$-cell $\bullet  \xrightarrow{x} \bullet $ for each element $x\in X$.
\item	One $2$-cell 
\begin{equation*}
\begin{tikzcd}[row sep=huge]
\bullet \arrow[r, "x"] \arrow[d, "y"]  & 	\bullet \arrow[d, "y"]  \\
\bullet \arrow[r, "x*y"]   & 	\bullet
\end{tikzcd}
\end{equation*}	
for each ordered pair $(x,y)$ of elements of $X$.		
\item One $3$-cell 
\begin{equation*}
\begin{tikzcd}[row sep=huge]
\bullet \arrow[rrr, "(x*y)*z=(x*z)*(y*z)"] &&& \bullet\\
&\bullet \arrow[r, "x*z"] \arrow[ul, "y*z"] & \bullet \arrow[ur, "y*z"] &\\    
&\bullet  \arrow[r, "x"] \arrow[dl, "y"]  \arrow[u, "z"] &  \bullet \arrow[dr, "y"]  \arrow[u, "z"]  & \\
\bullet  \arrow[rrr, "x*y"] \arrow[uuu, "z"] &&&\bullet \arrow[uuu, "z"]
\end{tikzcd}
\end{equation*}	
for each ordered triple $(x,y,z)$ of elements of $X$.
\item	One $n$-cell for each ordered $n$-tuple $(x_1, \ldots, x_n)$ of elements of $X$ for any integer $n> 3$.
\end{enumerate}

By construction, the rack space $BX$ is the quotient of the disjoint union $\bigsqcup_{n\ge 0} (X^n \times I^n)$ by the equivalence relation generated by
$$(x_1,\ldots, x_n,t_1, \ldots,t_{i-1},0,t_{i+1},\dots,t_n ) \sim (x_1,\ldots,x_{i-1},x_{i+1}, \dots,x_n,t_1, \ldots,t_{i-1},t_{i+1},\dots,t_n ) $$ and $$(x_1,\ldots, x_n,t_1, \ldots,t_{i-1},1,t_{i+1},\dots,t_n ) \sim (x_1*x_i,\ldots,x_{i-1}*x_i,x_{i+1}, \dots,x_n,t_1, \ldots,t_{i-1},t_{i+1},\dots,t_n ), $$ for $x_j \in X$ and $t_j \in [0,1]$.  If $X$ is a quandle, then the quandle space $\widehat{BX}$ for $X$ was defined by Nosaka in \cite{MR2786669} by modifying the construction of the rack space $BX$ for the underlying rack of $X$.

\begin{definition}
Let $X$ be a quandle. The quandle space $\widehat{BX}$ for $X$ is the space obtained from the rack space $BX$ by attaching $3$-cells which bound $2$-cells labeled $(x,x)$ for all $x \in X$.
\end{definition}

Since rack and quandle spaces are $CW$-complexes, all their homology and cohomology groups can be taken to be cellular (see \cite{MR1867354} for cellular homology).  Next, we attempt to understand homotopy-theoretic aspects of rack and quandle spaces. Before that, we present a connection between their fundamental groups and adjoint groups of corresponding racks and quandles \cite[Proposition 5.1]{MR2255194}.

\begin{proposition} \label{BXisAdj}
If $X$ is a quandle, then $\pi_1(BX) \cong \pi_1(\widehat{BX}) \cong \Adj(X)$.
\end{proposition}

\begin{proof}
Note that the space $BX$ has a single vertex and its edges are in bijection with elements of $X$. Thus, by the van-Kampen theorem,  elements of $X$ generate $\pi_1(BX)$. Further, the relations are given by the preferred squares, that is, $xy=y(x*y)$ for all $x,y \in X$. This is also the presentation of the adjoint group of the rack $X$. Since $BX$ is path-connected, it follows that $\widehat{BX}$ is also path-connected. Thus, we have $\pi_1(BX) \cong \pi_1(\widehat{BX}) \cong \Adj(X)$. $\blacksquare$
\end{proof}

It is well-known that, if $Z$ is a path-connected topological space, then its fundamental group $\pi_1(Z)$ acts canonically on its higher homotopy groups $\pi_n(Z)$ for each $n \ge 1$, which is simply the conjugation action when $n=1$ \cite[Chapter 4, p. 342]{MR1867354}.

\begin{proposition} \label{pi1 action on pin}\cite[Proposition 5.2]{MR2255194}
If $X$ is a rack, then the canonical action of $\pi_1(BX)$ on $\pi_n(BX)$ is trivial for each $n \ge 2$.
\end{proposition}

\begin{proof}
	Let $\alpha \in \pi_1(BX)$ be a generator corresponding to an element $x \in X$.  Let $\beta \in \pi_n(BX)$ be represented by a labelled diagram $D$ in $\mathbb{R}^n$.  Then the element $\alpha \cdot \beta$ is represented by the diagram which comprises a framed standard sphere in $\mathbb{R}^n$ labelled by $x$ and containing $D$ in its interior.  Now, pull the sphere under $D$ to one side (without changing any labels on $D$) and then eliminate it.  This diagram cobordism proves that the two diagrams are equivalent.
	$\blacksquare$
\end{proof}

The following result provides information on the rack space for a trivial rack \cite[Theorem 5.12]{MR2255194}.

\begin{theorem}
Let $\T_n$ be the trivial rack on $n$ elements.  Then the rack space $B \T_n$ has the homotopy type of the loop space of $\vee _n\mathbb{S}^2$, that is, $B\T_n \simeq \Omega(\vee_n  \mathbb{S}^2)$.
\end{theorem}

Recall that, two framed links $L_1, L_2$ in $\mathbb{R}^n$ are cobordant if there is a framed link $L$ properly embedded in $\mathbb{R}^n  \times I$ which meets $\mathbb{R} \times \{0, 1\}$ in $L_1, L_2$, respectively. If $X$ is a rack, then there is an analogous notion of cobordism of links with representation in $X$ (that is, a homomorphism of the fundamental rack to $X$) extending the given representations on the ends. The subsequent theorem delineates the homotopy type of the rack space of a free rack \cite[Theorem 5.13]{MR2255194}.

\begin{theorem}\label{homotopytype}
Let $FR_n$ denote the free rack of rank $n$.  Then the rack space $BFR_n$ has the homotopy type of $\vee_n \mathbb{S}^1$.
\end{theorem}

\begin{proof}
It follows from  \cite[Theorem 3.11]{MR2255194} that $\pi_n(BFR_n) \cong \mathcal{F}(\mathbb{R}^{n+1}, FR_n)$, where $\mathcal{F}(\mathbb{R}^{n+1},FR_n)$ is 
the set of cobordism classes of framed links in $\mathbb{R}^{n+1}$ with representation in $FR_n$.  Observe that $FR_n$ is the fundamental rack of $D_n$, which is the link comprising $n$ framed points in the $2$-disk $D^2$.  We see from \cite[Theorem 5.4]{MR2255194} that if $f$ represents an element of $\pi_n(BFR_n)$ ($n \geq 2$), then $f$ pulls back from a transverse map of $\mathbb{R}^{n+1}$ to $D_n$.  Since $D^2$ is contractible, this map is null homotopic. Applying relative transversality, we obtain a null cobordism of $f$.  By Proposition~\ref{BXisAdj}, $\pi_1(BFR_n) \cong F_n$,  the free group of rank $n$.  Hence, the result follows. $\blacksquare$
\end{proof}	

As an immediate consequence, we obtain the following result.

\begin{proposition}\label{homology free rack space}
Let $FR_n$ denote the free rack of rank $n$ and let $A$ be an abelian group.  Then
$$
\Ho_m^{\rm Cell}(BFR_n; A) \cong
\left\{
\begin{array}{ll}
A  &    \textrm{if}~ m=0,\\
A^n   &    \textrm{if}~ m=1,\\
1  & \textrm{otherwise,}
\end{array}
\right.
$$
where $A^n$ is the direct product of $n$ copies of $A$.
\end{proposition}
\bigskip
\bigskip


\section{Homotopy groups of rack and quandle spaces}\label{secondhomotopy}
Generally, determining the higher homotopy groups and other algebraic topological properties of rack and quandle spaces poses a non-trivial challenge. To present some results in this regard, we review some basic facts about spectral sequences. We shall also need it in Section \ref{section homology permutation racks} of Chapter \ref{chapter homology computations}. Loosely speaking, a spectral sequence is a tool for computing homology (or cohomology) from filtrations of a chain complex.  In general,  given a sequence of subcomplexes $\{F_pC\}_{p \in \mathbb{Z}}$ of a chain complex $C$ with $F_{p-1}C\subset F_pC$, it is reasonable to try to obtain information about the homology $H_*(C)$ of the chain complex $C$ in terms of the homology groups $H_*(F_p C/F_{p-1}C)$ of the subquotients. 
Often when computing homology, one encounters double complexes.  It is a family of doubly graded abelian groups of the form $(M_{i,j},d^h,d^v)$, where $d^h_{p,q}:M_{p,q} \rightarrow M_{p-1,q}$ and $d^v_{p,q}:M_{p,q} \rightarrow M_{p,q-1}$ such that 
$$ d^h_{p,q-1} ~ d^v_{p,q} + d^v_{p-1,q} ~  d^h_{p,q}=0$$
 is satisfied. In other words,  the squares in the diagram
	
	$$\begin{CD}
		@VVd^vV	@VVd^v V          @VVd^v V   @VVd^vV\\	
		\ldots@>d^h>>	M_{p,q}        @>d^h>> M_{p-1,q} @>d^h>> \ldots\\
		@VVd^vV	@VVd^v V          @VVd^v V  @VVd^vV \\
		\ldots @>d^h>>	M_{p,q-1}@>d^h>> M_{p-1,q-1} @>d^h>> \ldots\\
		@VVd^vV	@VVd^v V          @VVd^v V   @VVd^vV\\	
	\end{CD}$$
anti-commute.
	
\begin{definition}
		Let $M=M_{*,*}$ be a \index{doubly graded complex}{doubly graded complex}. Then its total complex $\Tot(M)$ is the complex with the $n$-th term
$$ \Tot(M)_n= \bigoplus_{ p +q=n}M_{p,q}$$
and the differential $D_n \Tot(M)_n \rightarrow \Tot(M)_{n-1}$ given by 
$$ D_n=\sum_{p+q=n}d^h_{p,q} + d^v_{p,q}.$$
\end{definition}
	
It is a routine computation based on the anti-commutativity of the preceding diagram to check that $D_{n-1} ~ D_n=0$, and hence the total complex $\Tot(M)$ is a chain complex.  We usually assume that $M_{p,q}=0$ whenever $p$ or $q$ is negative, thus obtaining what is usually called a {\it first quandrant double complex}.  Each column of the double complex is itself a chain complex with homology
$$ E^1_{p,q} :=\Ho_p(M_{p,*}, d^v).
$$	
The horizontal differential induces a map 
$$	d^1:=(d^h)_*: E^1_{p,q} \rightarrow E^1_{p-1,q},
$$
and hence we can take the homology of these horizontal complexes to get
$$	E^2_{p,q}:=\Ho_p(E^1_{*,q}, d)=\Ho_p^h\big(H_q^v(M_{*,*})\big).
$$
There is a map  $$	d^2:E^2_{p,q} \rightarrow E^2_{p-2,q+1} $$
obtained by letting $[x] \in E^2_{p,q} $ be the homology class of the cycle $x \in E^1_{p,q} $ (that is, $d^1(x)=0)$.  This element $x$ is itself the class of some vertical cycle $\tilde{x} \in M_{p,q}$  (that is, $d^v(\tilde{x})=0)$.  The cycle condition $d^1(x)=(d^h)_*(x)=0$ is equivalent to $d^h(\tilde{x})=d^v(y)$ for some $y \in M_{p-1,q+1}$.  From all these relations, we obtain 
$$ d^v \big(d^h(y) \big)=-d^h \big(d^v(y) \big)=-d^h \big(d^h(\tilde{x}) \big)=0.$$
 Hence, $d^h(y)$ is a vertical cycle in $M_{p-2,q+1}$ and determines a homology class $z=[d^h(y)]
	\in E^1_{p-2,q+1}$. To compute $d^1(z)$, we see that
$$	d^1(z)=(d^h)_*(z)=d^h_* \big(d^h(y) \big)=0. $$
Thus, $z$ is a cycle for $d^1$ and determines $d^2[x]:=[z] \in E^2_{p-1,q+1}$.  We can now give the definition of a spectral sequence.
	
\begin{definition}
A \index{spectral sequence}{spectral sequence} $E=\{E^r,d^r\}$ is a sequence of $\mathbb{Z}$-bigraded modules, each with diffrential $d^r:E^r_{p,q}\longrightarrow E^r_{p-r,q+r-1},$ of bidegree 	$(-r,r-1)$ such that there are isomorphisms $E^{r+1}\cong \Ho(E^r,d^r)$. 
	\end{definition}
	
The reader may refer to \cite{MR1793722} or \cite[Chapter 7]{MR0672956} for more on spectral sequences.  The next result shows that the second homotopy group $\pi_2(BX)$ of the rack space $BX$ of a rack $X$ is well-behaved \cite[Theorem 3.1]{MR2786669}.
	
\begin{theorem}\label{finiterack}
If $X$ is a finite quandle, then $\pi_2(BX)$ and $\pi_2 \big(\widehat{BX} \big)$ are finitely generated.
\end{theorem}
	
\begin{proof} 
Since $\widehat{BX}$ is obtained from $BX$ by attaching 3-cells, the inclusion $BX  \hookrightarrow \widehat{BX}$ induces an epimorphism $\pi_2(BX) \to  \pi_2 \big(\widehat{BX} \big)$ of homotopy groups. Thus, it suffices to show that $\pi_2(BX)$ is finitely generated. Note that $BX$ is a path-connected $CW$-complex. In view of Proposition \ref{pi1 action on pin}, the canonical action of $\pi_1(BX)$ on $\pi_2(BX)$ is trivial. Hence, by \cite[Lemma $8^{bis}.27$]{MR1793722} or \cite[Exercise II.5.1]{MR0672956}, we obtain the following exact sequence 
		\begin{eqnarray}\label{exactseq}
			\Ho_3^{\rm Grp} \big(\pi_1(BX); \mathbb{ Z} \big) \rightarrow \pi_2(BX) \rightarrow \Ho_2^{\rm Cell}(BX; \mathbb{ Z}).
		\end{eqnarray}
It suffices to prove that $\Ho_3^{\rm Grp} \big(\pi_1(BX); \mathbb{ Z} \big)$ and $\Ho_2^{\rm Cell}(BX; \mathbb{ Z})$ are finitely generated.	Since $BX$ has finitely many $2$-cells, it follows that $\Ho_2^{\rm Cell}(BX; \mathbb{ Z})$ is finitely generated.  Let $\phi: \Adj(X) \rightarrow \Inn(X)$ be the map given by $\textswab{a}_x \mapsto S_x$. Since $X$ is finite, it follows that both $\Adj(X)$ and $\Inn(X)$ are finitely generated. Consequently, we deduce that $\ker(\phi)$ is finitely generated (see also \cite[Proposition 10]{MR1769723}).  Now, the short exact sequence 
$$1 \rightarrow \ker(\phi) \rightarrow \Adj(X) \rightarrow \Inn(X) \rightarrow 1$$ of groups gives the Hochschild-Serre spectral sequence  $\{E^r,d^r\}$, whose second term is $$E_2^{p,q} \cong  \Ho_p^{\rm Grp} \big(\Inn(X); \Ho_q^{\rm Grp} (\ker(\phi); \mathbb{Z})\big)$$ and which converges to the group homology $\Ho_*^{\rm Grp} \big(\Adj(X); \mathbb{ Z} \big)$ (see \cite{MR0052438} or \cite[Theorem 6.3]{MR0672956}).  Since $\Inn(X)$ is finitely generated, the term $E_2^{p,q}$ is finitely generated, and hence $\Ho_*^{\rm Grp} \big(\Adj(X); \mathbb{ Z} \big)$ is finitely generated.  By Proposition \ref{BXisAdj}, we have $\pi_1(BX) \cong \Adj(X)$, and hence $\Ho_3^{\rm Grp} \big(\pi_1(BX); \mathbb{ Z} \big)$ is  finitely generated.  Now, it follows from the exact sequence~\eqref{exactseq}  that $\pi_2(BX)$ is finitely generated.
		$\blacksquare$
	\end{proof}

The second homotopy group for the case of dihedral quandles is given by the following result \cite[Proposition 4.1]{MR2786669}.

\begin{proposition}
Let $\R_p$ be the dihedral quandle of order $p$, where $p$ is an odd prime.  Then $\pi_2 \big(\widehat{B\R_p} \big)\cong \mathbb{ Z}_p$.
\end{proposition}

The following result identifies the second homotopy group for a specific Alexander quandle \cite[Proposition 4.5]{MR2786669}.

\begin{proposition}
Let $X=\mathbb{ Z}_2[T]/(T^2+T+1)$ equipped with the quandle operation $$x*y=Tx+(1-T)y.$$  Then $\pi_2(BX) \cong \mathbb{ Z}_8 \oplus \mathbb{ Z}_2$.  
\end{proposition}

\begin{remark}
{\rm Table \ref{table of homotopy groups} gives  second homotopy groups of some quandle spaces \cite[Table 1]{MR3054333}. We refer the reader to \cite[Appendic C.2]{MR3729413} for more details.  Recall that, the quandle $\mathbb{Z}_n[T]/(T-\omega)$ is the ring $\mathbb{Z}_n$ equipped with the quandle operation $$x*y=\omega x + (1-\omega)y,$$ where $\omega \neq \pm 1$ is an invertible element of $\mathbb{Z}_n$. Note that, when $\omega = -1$,  then $x * y =2y-x$, and we obtain the dihedral quandle $\R_n$.
	
\begin{table}[H]
\begin{center}
\begin{tabular}{|c||c|c|c|}
 \hline $\ \ \ \ \ \ \ \ \ \text{Quandle} \;X \ \ \ \ \ \ \ \ \ $ &  $\pi_2 \big(\widehat{BX} \big)$ \\
  \hline \hline
$\R_3 $&  $\mathbb{Z}_3$ \\ \hline
$\R_5 $&  $\mathbb{Z}_5$ \\ \hline
$\mathbb{Z}_5[T]/(T-\omega)$&  $1$ \\ \hline
$ \R_7 $&  $\mathbb{Z}_7$ \\ \hline
$\mathbb{Z}_7[T]/(T-\omega)$ &  $1$ \\ \hline
$\mathbb{Z}_2[T]/(T^3+T^2+1)$& $\mathbb{Z}_2$ \\\hline
$\mathbb{Z}_2[T]/(T^3+T+1)$& $\mathbb{Z}_2$ \\\hline
$ \R_9 $&  $\mathbb{Z}_9$ \\ \hline
$\mathbb{Z}_3[T]/(T^2+1)$&  $(\mathbb{Z}_3)^3$ \\ \hline
$\mathbb{Z}_3[T]/(T^2+T-1)$&  1 \\ \hline
$\mathbb{Z}_3[T]/(T^2-T-1)$& 1 \\
\hline 
\end{tabular}
\caption{Second homotopy groups of some quandle spaces.} \label{table of homotopy groups}
\end{center}
\end{table}
}
\end{remark}


\chapter{Homology and cohomology of solutions to the Yang--Baxter equation}\label{chapter homology YBE}

\begin{quote}
Fenn, Rourke, and Sanderson introduced a homology theory for racks as the homology of the rack space, as outlined in Chapter \ref{chapter trunks and rack space}. This initial development inspired the formulation of a (co)homology theory for racks and quandles by Carter, Jelsovsky, Kamada, Langford, and Saito. Later, a module theoretic approach to quandle cohomology was developed by Andruskiewitsch and Gra\~{n}a.  Subsequently,  Carter, Elhamdadi, and Saito extended the homology theory to encompass general solutions to the Yang--Baxter equation. These (co)homology theories have since undergone further expansion by Lebed and Vendramin, which encompasses and unifies all previously known theories in this domain. This chapter summarizes these advancements, includes selected computations, and explores their applications in knot theory.
\end{quote}
\bigskip

\section{Homology and cohomology of racks and quandles}\label{Homology} 
The subsequent construction of the chain complex leads to a homology theory for racks and quandles. This construction was introduced in \cite{MR1990571} as a variant of the homology of the rack space for a rack.
\para 
Let $(X,*)$ be a quandle. For each $n \ge 0$, let  $C_n^{\rm R}(X)$ be the free  abelian group generated by $n$-tuples $(x_1, \ldots, x_n)$ of elements of $X$. For each $n \geq 2$, define  $\partial_{n}: C_{n}^{\rm R}(X) \to C_{n-1}^{\rm R}(X)$ by 
\begin{eqnarray}
\partial_{n}(x_1, x_2, \dots, x_n) &=& \sum_{i=2}^{n} (-1)^{i} \big( (x_1, x_2, \dots, x_{i-1}, x_{i+1},\dots, x_n) \\
\nonumber && - (x_1 \ast x_i, x_2 \ast x_i, \dots, x_{i-1}\ast x_i, x_{i+1}, \dots, x_n) \big)
\end{eqnarray}
and set $\partial_n=0$ for $n \leq 1$. A simple verification demonstrates that $\partial_{n} ~ \partial_{n+1}=0$, and hence  
$C_\ast^{\rm R}(X)= \{C_n^{\rm R}(X), \partial_n \}$ is a chain complex.
\para 

For $n \geq 2$, let $C_n^{\rm D}(X)$ be the subgroup of $C_n^{\rm R}(X)$ generated
by $n$-tuples $(x_1, \dots, x_n)$ with $x_{i}=x_{i+1}$ for some $i \in \{1, \dots,n-1\}$, and otherwise set $C_n^{\rm D}(X)=0$. If $X$ is a quandle, then it follows that $\partial_n(C_n^{\rm D}(X)) \subseteq C_{n-1}^{\rm D}(X)$. Hence,
$C_\ast^{\rm D}(X) = \{ C_n^{\rm D}(X), \partial_n \}$ becomes a sub-complex of
$C_\ast^{\rm	R}(X)$.  Let $$C_n^{\rm Q}(X) := C_n^{\rm R}(X)/ C_n^{\rm D}(X)$$ and  $C_\ast^{\rm Q}(X) = \{ C_n^{\rm Q}(X), \overline{\partial}_n \}$ be the quotient complex, where $\overline{\partial}_n$ is the induced homomorphism, which we shall denote by $\partial_n$ again.
\para 

Let ${\rm W}={\rm D}$, ${\rm R}$, ${\rm Q}$. For each $n \ge 0$ and an abelian group $A$, define
$$ C_n^{\rm W}(X;A) = C_n^{\rm W}(X) \otimes_{\mathbb{Z}} A \quad \textrm{and} \quad \partial'_n :=\partial_n \otimes {\rm id}: C_n^{\rm W}(X;A) \to C_{n-1}^{\rm W}(X;A).$$
Similarly, we define
$$C^n_{\rm W}(X;A) = \Hom_{\mathcal{G}} \big(C_n^{\rm	W}(X), A \big) \quad \textrm{and} \quad \partial^n := \Hom(\partial_{n+1}, G): C^n_{\rm W}(X;A) \to C^{n+1}_{\rm W}(X;A).
$$
For convenience, we shall denote $\partial'_n$ by $\partial_n$.
\para 

Once we have the chain and the cochain complexes, we can define homology groups
$$\Ho_n^{\rm W}(X; A) = \Ho_{n} \big(C_\ast^{\rm W}(X;A) \big)$$
and cohomology groups
$$\Ho^n_{\rm W}(X; A)  = \Ho^{n}\big(C^\ast_{\rm W}(X;A)\big).$$
For ${\rm W}={\rm D}$, ${\rm R}$, ${\rm Q}$ and $n \ge 0$, the group $\Ho_n^{\rm W}(X; A)$ is called the $n$-th \index{degenerate homology}{degenerate}, \index{rack homology}{rack} and \index{quandle homology}{quandle homology} group of $X$ with coefficients in the abelian group $A$, respectively. Similarly, the group $\Ho^n_{\rm W}(X; A)$ is called the $n$-th \index{degenerate cohomology}{degenerate}, \index{rack cohomology}{rack} and \index{quandle cohomology}{quandle cohomology} group of $X$ with coefficients in the abelian group $A$, respectively.   We make use of the convention that whenever the coefficients of the (co)homology groups are not specified, we will mean $\mathbb{Z}$ coefficients. We denote by $\beta_n^{\rm W}(X)$ the $n$-th \index{Betti number}{Betti number} of $X$ determined by the homology group $\Ho_n^{W}(X)$.
\para 

It is natural to expect that the homology of racks and that of rack spaces are related. It follows from the construction of the rack space, the definition of cellular homology of a CW-complex and the definition of the rack homology of a rack,  that the cellular homology of the rack space is the same as the rack homology of the rack  \cite[Proposition 5.11]{MR3729413}.

\begin{proposition}\label{cellular and rack same}
Let $X$ be a rack.  Then the cellular chain complex of the rack space $BX$ and the rack chain complex of $X$ are isomorphic as chain complexes. As a consequence, 
$$\Ho_n^{\rm R}(X; A) \cong \Ho_n^{\rm Cell}(BX; A)\quad \textrm{and} \quad \Ho^n_{\rm R}(X; A) \cong \Ho^n_{\rm Cell}(BX; A)$$ for each $n \ge 0$ and for each abelian group $A$.
\end{proposition}

The simillar result holds for quandles and quandle spaces \cite[Proposition 5.11]{MR3729413}.

\begin{proposition}\label{cellular and quandle same}  
Let $X$ be a quandle.  Then the cellular chain complex of the quandle space $\widehat{BX}$  and the quandle chain complex of $X$ are isomorphic as chain complexes. As a consequence, 
		$$\Ho_n^{\rm Q}(X; A) \cong \Ho_n^{\rm Cell} \big(\widehat{BX}; A \big)\quad \textrm{and} \quad \Ho^n_{\rm Q}(X; A) \cong \Ho^n_{\rm Cell}\big(\widehat{BX}; A\big)$$ for each $n \ge 0$ and for each abelian group $A$.
	\end{proposition}

In view of Proposition \ref{cellular and rack same} and Proposition \ref{homology free rack space}, we have the following result for free racks.

\begin{proposition}
Let $FR_n$ be the free rack of rank $n$ and $A$ an abelian group.  Then
$$
\Ho_m^{\rm R}(FR_n; A) \cong
\left\{
\begin{array}{ll}
A  &    \textrm{if}~~ m=0,\\
A^n   &    \textrm{if}~~ m=1,\\
1  & \textrm{otherwise.}
\end{array}
\right.
$$
\end{proposition}

See \cite[Theorem 3.1]{MR3196064} and \cite[Remark 2.6]{arXiv:2011.04524} for independent proofs of the preceding proposition. Further, a similar result holds for free quandles \cite[Theorem 4.1]{MR3196064}.

\begin{proposition}
Let $FQ_n$ be the free quandle of rank $n$ and $A$ an abelian group.  Then
$$
\Ho_m^{\rm Q}(FQ_n; A) \cong
\left\{
\begin{array}{ll}
A  &    \textrm{if}~~ m=0,\\
A^n   &    \textrm{if}~ ~m=1,\\
1  & \textrm{otherwise.}
\end{array}
\right.
$$
\end{proposition}

\begin{proposition}\cite[Proposition 3.3]{MR1812049}
If $X$ is a quandle, then there is a long exact sequence 
 \begin{eqnarray}\label{Basic Homology Long Exact Sequence}
		\cdots \stackrel{\updelta_n}{\longrightarrow} \Ho_n^{\rm D}(X;A) \stackrel{i_n}{\longrightarrow} \Ho_n^{\rm R}(X;A)
		\stackrel{j_n}{\longrightarrow} \Ho_n^{\rm Q}(X;A)
		\stackrel{\updelta_{n-1}}{\longrightarrow} \Ho_{n-1}^{\rm D}(X;A) \longrightarrow \cdots
	\end{eqnarray}
of quandle homology groups which is natural with respect to homomorphisms induced by quandle homomorphisms. Here, each $\updelta_n$ is a connecting homomorphism.
\end{proposition}

\begin{proof}
By standard homological algebra, the short exact sequence of chain complexes    
$$1 \longrightarrow C_*^{\rm D}(X)\otimes_{\mathbb{Z}} A \stackrel{i\otimes \id}{\longrightarrow} C_n^{\rm R}(X)\otimes_{\mathbb{Z}}  A \stackrel{j\otimes \id}{\longrightarrow} C_n^{\rm Q}(X) \otimes_{\mathbb{Z}}  A \longrightarrow 1$$ 
induces the desired long exact sequence of homology groups. $\blacksquare$
\end{proof}

We have the following split exact sequences \cite[Proposition 3.4]{MR1812049}.  

\begin{proposition}\label{Universal Coefficient Theorem}
	There exist split exact sequences
	\begin{eqnarray*}
		1 \longrightarrow \Ho_n^{\rm W}(X) \otimes A \longrightarrow \Ho_n^{\rm W}(X; A) \longrightarrow
		{\rm Tor} \big(\Ho_{n-1}^{\rm W}(X), A \big) \longrightarrow 1 \\
		1 \longrightarrow {\rm Ext} \big(\Ho_{n-1}^{\rm W}(X), A \big) \longrightarrow \Ho^n_{\rm W}(X; A) \longrightarrow
		{\rm Hom} \big(\Ho_n^{\rm W}(X), A \big) \longrightarrow 1,
	\end{eqnarray*}
	where ${\rm W} = {\rm R}$ if $X$ is
	a rack, or one of ${\rm D}$, ${\rm R}$, ${\rm Q}$ if $X$ is a quandle.
\end{proposition}

\begin{proof}
Since $C_*^W(X)$ is a chain complex of free abelian groups, the result follows from the \index{Universal Coefficient Theorem}Universal Coefficient Theorems \cite[Theorems 53.1 and 55.1]{MR0755006}.  $\blacksquare$
\end{proof}

Let us consider some elementary examples of computations \cite[Examples 3.5 and 3.6]{MR1812049}.

\begin{example}\label{dihedral quandle of order three homology}
{\rm 
Let $\R_3$ be the dihedral quandle of three elements. Direct calculations using the definition gives
\begin{eqnarray*}
\Ho_1^{\rm Q}(\R_3) = \mathbb{Z} \quad \textrm{and} \quad  \quad \Ho_2^{\rm Q}(\R_3) = 1. 
\end{eqnarray*}
Thus, we have
\begin{eqnarray*}
\Ho_2^{\rm Q}(\R_3;A) = 1 \quad \textrm{and} \quad  \quad \Ho^2_{\rm Q}(\R_3;A) = 1 
\end{eqnarray*}
for any coefficient group $A$.
\para

Let $\R_4$ be the dihedral quandle of four elements. Direct calculations give
\begin{eqnarray*}
\Ho_1^{\rm Q}(\R_4) = \mathbb{Z}^2 \quad \textrm{and} \quad  \Ho_2^{\rm Q}(\R_4) = \mathbb{Z}^2 \oplus (\mathbb{Z}_2)^2. 
\end{eqnarray*}
Hence, we obtain
\begin{eqnarray*}
\Ho_2^{\rm Q}(\R_4;\mathbb{Z}_2) = (\mathbb{Z}_2)^4 \quad \textrm{and} \quad  \quad \Ho^2_{\rm Q}(\R_4;\mathbb{Z}_2) = (\mathbb{Z}_2)^4
\end{eqnarray*}
and
\begin{eqnarray*}
\Ho_2^{\rm Q}(\R_4;\mathbb{Z}_m) = (\mathbb{Z}_m)^2 \quad \textrm{and} \quad  \quad \Ho^2_{\rm Q}(\R_4;\mathbb{Z}_m) = (\mathbb{Z}_m)^2
\end{eqnarray*}
for any odd positive integer $m$.
}
\end{example}

\begin{example}\label{homology of trivial quandle}
{\rm 
Let $\T_m$ be the trivial quandle with $m < \infty$ elements. Since each $\partial_n :C_n^{\rm R}(\T_m) \to C_{n-1}^{\rm R}(\T_m)$ is the zero map, the connecting homomorphisms $\updelta_\ast$ in the long exact homology sequence with $A=\mathbb{Z}$ are
zero maps. Hence, the long exact sequence is decomposed into the short exact sequences
\begin{eqnarray*}
1 \to \Ho_n^{\rm D}(\T_m) \to \Ho_n^{\rm R}(\T_m) \to \Ho_n^{\rm Q}(\T_m) \to 1,
\end{eqnarray*}
which are further identified with the short exact sequences 
\begin{eqnarray*}
1 \to C_n^{\rm D}(\T_m) \to C_n^{\rm R}(\T_m) \to C_n^{\rm Q}(\T_m) \to 1.
\end{eqnarray*}
In particular, we have
\begin{eqnarray*}
\beta_n^{\rm D}(\T_m) &=& {\rm rank} \big(C_n^{\rm D}(\T_m) \big),\\
\beta_n^{\rm R}(\T_m) &=& m^n,\\ 
\beta_n^{\rm Q}(\T_m) &=& {\rm rank} \big(C_n^{\rm Q}(\T_m) \big).
\label{eqn:BettiTrivial}
\end{eqnarray*}
}
\end{example}
\para 

Let $X$ be a quandle and $\mathcal{O}(X)$ the set of orbits of $X$. For simplicity, we assume that $|\mathcal{O}(X)| = m < \infty$. Let $\pi: X \to \mathcal{O}(X) \cong \T_m$ be the natural projection from $X$ onto its set of orbits, which is identified with the trivial quandle $\T_m$. The naturality of the
 long exact homology  sequence gives the following  commutative diagram
\begin{eqnarray} \label{eqn:BasicTrivializingDiagram}
\begin{array}{ccccccccc}
\cdots & \stackrel{\updelta_n}{\longrightarrow}
& \Ho_n^{\rm D}(X) & \stackrel{i_n}{\longrightarrow}
& \Ho_n^{\rm R}(X) & \stackrel{j_n}{\longrightarrow}
& \Ho_n^{\rm Q}(X) & \stackrel{\updelta_{n-1}}{\longrightarrow} & \cdots
\\ {} & {} & \downarrow {\pi_\ast} & {} & \downarrow {\pi_\ast} & {} & \downarrow {\pi_\ast} & {} & {} \\
1 & \longrightarrow
& \Ho_n^{\rm D}(\T_m) & \stackrel{i_n}{\longrightarrow}
& \Ho_n^{\rm R}(\T_m) & \stackrel{j_n}{\longrightarrow}
& \Ho_n^{\rm Q}(\T_m) & \longrightarrow 1. & 
\end{array}
\end{eqnarray}

\begin{proposition}\cite[Proposition 3.8]{MR1812049}.
Let $X$ be a quandle. Then the following  assertions hold:
\begin{enumerate}
\item $\Ho_1^{\rm D}(X)= 1$.
\item $\Ho_1^{\rm R}(X)\cong \Ho_1^{\rm Q}(X) \cong \mathbb{Z}^{|\mathcal{O}(X)|}$.
\end{enumerate}
\end{proposition}

\begin{proof}
By definition, $\Ho_1^{\rm D}(X)= \Ho_0^{\rm D}(X)=1$ for any quandle $X$. Thus, it follows from the long exact homology sequence that $\Ho_1^{\rm R}(X) \cong \Ho_1^{\rm Q}(X)$ for any quandle $X$.  The group $Z_1^{\rm R}(X)$ of 1-cycles is freely generated by elements of $X$ and the  group $B_1^{\rm R}(X)$ of 1-boundaries is generated by the images $\partial_2((x,y))= x - x \ast y$ for all $x,y \in X$. Hence, $\Ho_1^{\rm R}(X)$ is generated by $\{ [x_\omega] \, \mid \, \omega \in \mathcal{O}(X) \}$, where $[x_\omega]$ is the homology class of $x_\omega$, which is a representative of the orbit $\omega$ of $X$. Let $m =|\mathcal{O}(X)|$ and $\T_m$ the trivial quandle with $m$ elements. By Example \ref{homology of trivial quandle}, $\beta_1^{\rm R}(\T_m) = m$, and hence $\Ho_1^{\rm R}(\T_m)\cong \mathbb{Z}^m$ is the free abelian group generated by $\T_m \cong \mathcal{O}(X)$. Considering the commutative diagram \eqref{eqn:BasicTrivializingDiagram} for $n=1$ and noting that $\pi_\ast: \Ho_1^{\rm R}(X) \to \Ho_1^{\rm R}(\T_m)$ maps $[x_\omega]$ to $[\omega]$,  we see that $\pi_\ast$ an isomorphism.
$\blacksquare$
\end{proof}

A result of Litherland and Nelson \cite{MR1952425} shows that the long exact sequence \eqref{Basic Homology Long Exact Sequence} decomposes into short exact sequences of quandle homology groups. We need the following set-up to prove this result. 
\para 
For $x = (x_1, \ldots, x_n)\in X^n$ and $y \in X$, we set
$x \diamond y := (x_1, \ldots, x_n, y)\in X^{n+1}$ and $ x*y: =
(x_1*y, \ldots, x_n*y) \in X^n$. Then, for $c \in C^{\rm R}_n(X)$, we define
$c \diamond y \in C^{\rm R}_{n+1}(X)$ and $c*y\in C^{\rm R}_n(X)$ by linearity. More precisely, if $c= \sum_{i=1}^k \alpha_i (x_{i_1}, x_{i_2}, \ldots,  x_{i_n})$, then 
$$c \diamond y= \sum_{i=1}^k \alpha_i (x_{i_1}, x_{i_2}, \ldots,  x_{i_n}, y) \quad \textrm{and} \quad c * y= \sum_{i=1}^k \alpha_i (x_{i_1}*y, x_{i_2}*y, \ldots,  x_{i_n}*y).$$ 
Also, note that $$\partial_{n+1}(c \diamond y) = \partial_n(c) \diamond y +(-1)^{n+1}(c - c*y).$$ Next we define homomorphisms ${\alpha_n}:{C^{\rm R}_n(X)} \rightarrow {C^{\rm R}_n(X)}$ by induction on $n$. We take $\alpha_1$ to be the identity map. For $n \geq 1$, $x = (x_1, \ldots, x_n)\in X^n$ and $y\in X$, we define $\alpha_{n+1}$ by the recursive formula
$$
 \alpha_{n+1}( x \diamond y) = \alpha_n( x)\diamond y -\alpha_n( x)\diamond x_n.
$$
We also define homomorphisms ${\beta_n}:{C^{\rm R}_n(X)} \rightarrow {C^{\rm R}_{n+1}(X)}$
by $\beta_n( x) = \alpha_n( x)\diamond x_n$.  Then, for any $c \in
C^{\rm R}_n(X)$ and $y\in X$, we have the relation
\begin{equation}\label{equation in litherland2}
 \alpha_{n+1}(c \diamond y) = \alpha_n(c)\diamond y - \beta_n(c).
\end{equation}

\begin{lemma}\label{litherland nelson lemma} 
The homomorphisms ${\alpha_n}:{C^{\rm R}_n(X)} \rightarrow {C^{\rm R}_n(X)}$
form a chain map ${\alpha}:{C^{\rm R}_*(X)} \rightarrow {C^{\rm R}_*(X)}$.
\end{lemma}

\begin{proof} 
Induction on $n$ shows that $\alpha_n( x*y) = \alpha_n(x)*y$ for $ x \in X^n$ and $y\in X$. We prove by induction on $n$ that $\partial_n \alpha_n = \alpha_{n-1}\partial_n$ for $n \geq 2$. For $n = 2$, we have $\alpha_2(x,y) = (x,y)-(x,x)$ and hence $\partial_2 \alpha_2(x,y) = \partial_2(x,y) =
\alpha_1\partial_2(x,y)$.  Suppose that the result is true for some
$n \geq 2$, and let $ x= (x_1, \ldots, x_n) \in X^n$ and $y\in X$. We compute
\begin{eqnarray*}
\partial_{n+1} \alpha_{n+1}( x \diamond y)
&=&\partial_{n+1} \bigl(\alpha_n( x) \diamond y\bigr) -
\partial_{n+1} \bigl(\alpha_n( x) \diamond  x_n\bigr)\\
&=&\partial_n \alpha_n( x) \diamond y + (-1)^{n+1}\bigl(\alpha_n( x) -
\alpha_n( x)*y\bigr) \\
&&- \partial_n  \alpha_n( x) \diamond  x_n
- (-1)^{n+1}\bigl(\alpha_n( x)
- \alpha_n( x)*{x_n}\bigr) \\
&=&\alpha_{n-1} \partial_n( x) \diamond y -
\alpha_{n-1}\partial_n ( x) \diamond x_n
+ (-1)^n\alpha_n( x*y) -
(-1)^n\alpha_n( x*{x_n})
\end{eqnarray*}
and
\begin{eqnarray*}
\alpha_n\partial_{n+1} ( x \diamond y) &=& \alpha_n\bigl(\partial_n ( x) \diamond y
+ (-1)^{n+1}( x -  x * y)\bigr)\\
&=& \alpha_{n-1}\partial_n ( x)\diamond y - \beta_{n-1}\partial_n ( x)
- (-1)^n\alpha_n( x) + (-1)^n\alpha_n( x*y).
\end{eqnarray*}
Hence, $\partial_{n+1} \alpha_{n+1}( x \diamond y) = \alpha_n\partial_{n+1} ( x \diamond y)$ if and only if
\begin{equation}\label{equation in litherland}
\alpha_{n-1} \partial_n(x) \diamond x_n + (-1)^n\alpha_n( x*{x_n})
= \beta_{n-1}\partial_n ( x) + (-1)^n\alpha_n( x).    
\end{equation}

We shall need the components of the boundary maps $\partial_n$ of the chain complex. For $1 \leq i \leq n$, $\epsilon = 0$ or 1, let $ {\partial ^\epsilon_i}: {X^n} \rightarrow {X^{n-1}}$ be the maps given by
\begin{eqnarray*}
 \partial^0_i(x_1,\ldots,x_n) & = & (x_1,\ldots,x_{i-1},x_{i+1},\ldots,x_n)\\
 \hbox{and}\qquad
 \partial^1_i(x_1,\ldots,x_n) & = & (x_1*{x_i},\ldots,x_{i-1}*{x_i},x_{i+1},\ldots,x_n).
\end{eqnarray*}
Now, for $1 \leq i < n$ and $\epsilon = 0$ or 1, the element  $\partial_i^\epsilon( x)$ of $X^{n-1}$ has last entry $x_n$. Thus, $\alpha_{n-1}\partial_i^\epsilon( x)\diamond x_n
=\beta_{n-1}\partial_i^\epsilon( x)$. Further, using \eqref{equation in litherland2}, we obtain
\begin{eqnarray*}
\alpha_{n-1}\partial_n^0( x)\diamond x_n-\beta_{n-1}\partial_n^0( x)
&=& \alpha_n\bigl(\partial_n^0( x)\diamond x_n\bigr) =\alpha_n( x)\\
\hbox{and}\qquad
\alpha_{n-1}\partial_n^1( x)\diamond x_n -\beta_{n-1}\partial_n^1( x)
&=&\alpha_n\bigl(\partial_n^1( x)\diamond x_n\bigr)
=\alpha_n( x*{x_n}).
\end{eqnarray*}
It follows that
$$\alpha_{n-1}\partial_n( x)\diamond x_n -\beta_{n-1}\partial_n( x)
 =(-1)^n\bigl(\alpha_n( x)-\alpha_n( x*{x_n})\bigr),
$$
which establishes \eqref{equation in litherland}, and hence we obtain the lemma.
$\blacksquare$
\end{proof}

\begin{theorem}\cite[Theorem 4]{MR1952425}
If $X$ is a quandle, then each connecting homomorphism $\updelta_\ast: \Ho_n^{\rm Q}(X) \to \Ho_{n-1}^{\rm D}(X)$ is zero. Consequently, the long exact sequence \eqref{Basic Homology Long Exact Sequence} decomposes into short exact sequences
	\begin{eqnarray}
		1 \longrightarrow \Ho_n^{\rm D}(X)  \stackrel{i_n}{\longrightarrow} \Ho_n^{\rm R}(X)
		\stackrel{j_n}{\longrightarrow}
		\Ho_n^{\rm Q}(X) \longrightarrow 1
  \end{eqnarray}
  of homology groups for each $n$.
\end{theorem}

\begin{proof} 
By Lemma \ref{litherland nelson lemma}, the maps $\gamma_n: C^{\rm R}_n(X) \to C^{\rm R}_n(X)$ given by $\gamma_n(c)= c - \alpha_n(c)$ gives a chain map $\gamma_*: C^{\rm R}_*(X) \to C^{\rm R}_*(X)$. We show that each $\gamma_n$ is a projection onto the subcomplex $C^{\rm D}_n(X)$. Thus, we need to prove the following assertions:
\begin{enumerate}
\item If $c \in C^{\rm D}_n(X)$, then $\alpha_n(c) = 0$.
\item If $c \in C^{\rm R}_n(X)$, then $c-\alpha_n(c) \in C^{\rm D}_n(X)$.
\end{enumerate}
For $n=1$, $C^{\rm D}_n(X)=0$, and hence assertion (1) is true in this case. Let $n \geq 1$, $x = (x_1, \ldots, x_n)\in X^n$ and $y\in X$. Suppose that $ x \diamond y \in C^{\rm D}_{n+1}(X)$.  Then either $ x \in C^{\rm D}_n(X)$ or $x_n = y$. If $ x \in C^{\rm D}_n(X)$, then  using the formula $ \alpha_{n+1}( x \diamond y) = \alpha_n( x)\diamond y -\alpha_n( x)\diamond x_n$ and the induction hypothesis, we obtain $\alpha_{n+1}( x \diamond y) = 0$. And, if $x_n = y$, then $\alpha_{n+1}( x \diamond y)=\alpha_{n+1}( x \diamond x_n) =\alpha_n(x) \diamond x_n-\alpha_n(x) \diamond x_n =0$.  Since $\alpha_n$ is linear, if $c \in C^{\rm D}_n(X)$ is a general element, then $\alpha_n(c) = 0$.  This proves assertion (1).  
\par

Assertion (2) is clear for $n=1$. Let  $n \geq 1$, $ x \in X^n$ and $y\in X$. Then we have
$$ x \diamond y - \alpha_{n+1}( x \diamond y) = \bigl(x -\alpha_n( x)\bigr)\diamond y -\bigl( x -\alpha_n( x)\bigr)\diamond x_n + x \diamond x_n.$$
Since $ x \diamond x_n \in C^{\rm D}_{n+1}(X)$ and $x -\alpha_n( x) \in C^{\rm D}_{n}(X)$ by induction hypothesis, it follows that  $ x \diamond y -\alpha_{n+1}( x \diamond y)$, and assertion (2) follows. It follows that $\gamma_* \,i_*: C^{\rm D}_*(X) \to C^{\rm D}_*(X)$ is the identity chain map, and hence the induced maps $i_n: \Ho_n^{\rm D}(X) \rightarrow \Ho_n^{\rm R}(X)$ of homology groups are injective. Exactness of the long exact homology sequence implies that the connecting homomorphisms $\updelta_n:\Ho_n^{\rm Q}(X) \rightarrow \Ho^{\rm D}_{n-1}(X)$ are zero, which proves the theorem. $\blacksquare$
\end{proof}

\begin{remark}
The preceding theorem was conjectured in \cite[p. 142]{MR1812049} and verified therein for $n=3$ \cite[Proposition 3.9]{MR1812049}. Further, the theorem was proved for $n=4$ in \cite{MR1825963}.
\end{remark}

In \cite[Section 4]{MR1948837}, Etingof and Gra\~na extended and sharpened results of \cite{MR1812049, MR1952425, MR1960136} on Betti numbers of finite quandles.  More precisely, they proved that if $X$ is a finite quandle, then the lower bounds for the Betti numbers of rack, quandle and degenerate cohomologies of $X$ given in \cite{MR1812049} are in fact equalities. 

\begin{theorem}
Let $X$ be a finite quandle with $m$ orbits under the action of $\Adj(X)$ on $X$. Then the Betti numbers of the rack, quandle and degenerate cohomologies of $X$ are given by
$$\beta_n^{\rm D}(X) = a_n, \quad  \beta_n^{\rm R}(X) = m^n \quad \textrm{and} \quad  \beta_n^{\rm Q}(X) = b_n,$$
where $a_n+b_n=|X|^n$ and $b_n= |X| \big(|X|-1 \big)^{n-1}$ for $n \ge 1$.
\end{theorem}

\begin{remark}
{\rm
It is well-known that the bounded cohomology of groups is intimately related to word metrics on groups, and is a key tool in the geometry of manifolds and rigidity theory. Building on this insight, Kedra \cite{MR4779104} introduced the concept of bounded cohomology for racks and quandles, linking it to a natural metric derived from their internal symmetries. This work appears to have initiated the incorporation of geometric perspectives into the study of racks and quandles.}
\end{remark}
\bigskip
\bigskip


\section{Module theoretic approach to quandle cohomology}
In this section, we consider a general cohomology theory for racks and quandles developed by Andruskiewitsch and Gra\~{n}a \cite{MR1994219}.  Braided vector spaces arising from groups play a fundamental role in the classification of finite-dimensional complex pointed Hopf algebras. The most important class of braided vector spaces arising from groups is the class of braided vector spaces $(\mathbb{C}X, r_\phi )$, where $X$ is a rack and $\phi$ is a 2-cocycle on
$X$ with coefficients in $\mathbb{C}$. This lead Andruskiewitsch and Gra\~{n}a to study cohomology of racks and quandles, thereby developing a general cohomology theory that encompasses the previously defined cohomology theories.  
\para 

\subsection{Cohomology with coefficients in a quandle module}
Let $X$ be a rack or a quandle and $\Omega (X)$ the free ${\mathbb{Z}}$-algebra  generated by the symbols $\eta_{x,y},$ $\xi_{x,y}$   for $x,y \in X$ with the condition that $\eta_{x,y}$ is invertible for every $x,y \in X$.  Then we have the following definition \cite[Definition 2.22]{MR1994219}. 

\begin{definition}\label{QuandAlg}
The  \index{rack algebra}{\it rack algebra} $\mathbb{Z}\{X\}$ of a rack $X$ is defined to be the  quotient  of $\Omega(X)$ by the relations
\begin{eqnarray}
  \label{rack alg rel 1}\eta_{x*y,z}\eta_{x,y} & =&  \eta_{x*z,y*z}\eta_{x,z},\\
 \label{rack alg rel 2}\eta_{x*y,z}\xi_{x,y} & =&  \xi_{x*z,y*z}\eta_{y,z},\\
 \label{rack alg rel 3}\xi_{x*y,z}& =& \eta_{x*z,y*z} \xi_{x,z} + \xi_{x*z,y*z}\xi_{y,z}.
\end{eqnarray}
The  \index{quandle algebra}{\it quandle algebra} $\mathbb{Z}(X)$ of a quandle $X$ is defined to be the  quotient  of $\Omega(X)$ by the relations \eqref{rack alg rel 1}--\eqref{rack alg rel 3} together with the relation
\begin{eqnarray}
\label{quandle alg rel 1} \xi_{x,x} + \eta_{x,x} &=& 1.
\end{eqnarray}
\end{definition}

\begin{remark}{\rm 
Notice the abuse of terminology. The quandle algebra defined above is not the same as the quandle algebra defined in Chapter \ref{chaper quandle rings}.}
\end{remark}

\begin{definition}
Let $X$ be a rack. A  representation of the rack algebra $\mathbb{Z}\{X\}$ is an abelian group $A$ together a family $\{\eta_{x,y} \in \Aut(A) \mid x, y \in X \}$ of automorphisms  and a family $\{ \xi_{x,y}\in \End(A) \mid x, y \in X$ of endomorphisms  such that relations \eqref{rack alg rel 1}--\eqref{rack alg rel 3} hold. 
\para
If $X$ is a quandle, then a representation of the quandle algebra $\mathbb{Z}(X)$ is an abelian group $A$ together a collection of automorphisms $\{\eta_{x,y} \in \Aut(A) \mid x, y \in X \}$ and a collection of endomorphisms  $\{ \xi_{x,y}\in \End(A) \mid x, y \in X$ such that relations \eqref{rack alg rel 1}--\eqref{quandle alg rel 1} hold. 
\end{definition}

In other words, a representation of the quandle algebra $\mathbb{Z}(X)$ means that there is a $\mathbb{Z}$-algebra homomorphism  $$\mathbb{Z}(X) \rightarrow \End(A).$$ We denote the  image of a generator by the same symbol. Given a representation $A$ of  $\mathbb{Z}(X)$, we say that $A$ is a $\mathbb{Z}(X)$-module or \index{quandle module}{\it a quandle module}.  The action of $\mathbb{Z}(X)$ on $A$ is written as a left-action  and denoted by $$(\rho, g) \mapsto \rho g  =  \rho (g)$$ for  $g\in A$ and  $\rho \in \End(A)$, where $\rho$ is the image of the element  $\rho \in \mathbb{Z}(X)$ under the homomorphism $\mathbb{Z}(X) \rightarrow \End(A)$

\begin{example}
{\rm Let $X$ be a quandle.  Let $\Lambda = \mathbb{Z}[T,T^{-1}]$ denote the ring of Laurent polynomials over the ring of integers and $A$ a $\Lambda$-module. Then $A$ is a left $\mathbb{Z}(X)$-module by setting 
$$\eta_{x,y} (a)=Ta \quad \textrm{and}\quad \xi_{x,y} (b) = (1-T) (b)$$
for $x,y \in X$ and $a, b \in A$. Similarly, any left $\Adj(X)$-module $A$ can be viewed as a left $\mathbb{Z}(X)$-module by setting 
\begin{equation}\label{module over rack algebra via Adj}
\eta_{x,y} (a)=\textswab{a}_y a \quad \textrm{and}\quad \xi_{x,y} (b) = b- \textswab{a}_{x*y}b
\end{equation}
for $x, y \in X$ and $a, b \in A$.}
\end{example}

We now present the construction of the generalized quandle homology theory due to Andruskiewitsch and Gra\~{n}a \cite{MR1994219}. Let $(X,*)$ be a quandle and $C_n(X)$  the free  right  $\mathbb{Z}(X)$-module  with basis $X^n$, where $X^0=\{ x_0 \}$ for some fixed element $x_0 \in X$.  Let $\partial_n: C_{n+1}(X) \rightarrow C_n(X)$ be the $\mathbb{Z} ( X )$-linear map defined on the basis elements by 
\begin{eqnarray} \label{Andruskiewitsch grana cohomology}
&& \partial_n (x_1, \ldots, x_{n+1} )\\
\nonumber  &=& {\displaystyle  \sum_{i=2}^{n+1} (-1)^i \eta_{ [x_1, \ldots, x_{i-1}, x_{i+1}, \ldots, x_{n+1} ], [x_i, \ldots, x_{n+1}] }
		(x_1, \ldots, x_{i-1}, x_{i+1}, \ldots, x_{n+1} ) } \\
 \nonumber & & \\
\nonumber & & - {\displaystyle  \sum_{i=2}^{n+1} (-1)^i 
		(x_1*x_i, \ldots, x_{i-1}*x_i, x_{i+1},  \ldots, x_{n+1} ) }  \\
\nonumber 	& & \\
\nonumber & & +~	\xi_{[x_1, x_3, \ldots, x_{n+1}], [x_2, x_3, \ldots, x_{n+1}]} (x_2, \ldots, x_{n+1} ),
\end{eqnarray}
for $n \ge 1$ and
\begin{equation}\label{Andruskiewitsch grana cohomology at 0}
 \partial_0(x) = - \xi_{x *^{-1} x_0, x_0}(x_0),
\end{equation}
where
$$	[x_1, x_2, \ldots, x_n] = ( ( \cdots ( x_1 * x_2 ) * 	x_3  ) * \cdots ) * x_n$$
is the left-normalised product. We then have the following result.

\begin{lemma}
$\{C_n(X), \partial_n\}$ is a chain complex.
\end{lemma}

\begin{proof}
We write $\partial_n=  \sum_{i=1}^{n+1} (-1)^i  \partial_n^i$, where
	\begin{eqnarray*} 
		\partial_n^i (x_1, \ldots, x_{n+1} ) &=&
			\eta_{ [x_1, \ldots, x_{i-1}, x_{i+1}, \ldots, x_{n+1} ], [x_i, \ldots, x_{n+1}] }
			(x_1, \ldots, x_{i-1}, x_{i+1}, \ldots, x_{n+1} )  \\
		& & - (x_1*x_i, \ldots, x_{i-1}*x_i, x_{i+1},  \ldots, x_{n+1} ) 
	\end{eqnarray*}
for $2 \le i \leq n+1$ and
	\begin{eqnarray*} 
	\partial_n^{1} (x_1, \ldots, x_{n+1} ) &=& 	\xi_{[x_1, x_3, \ldots, x_{n+1}], [x_2, x_3, \ldots, x_{n+1}]}
	(x_2, \ldots, x_{n+1} ).
\end{eqnarray*}
A direct calculation verifies that
\begin{equation}\label{module theoretic boundary}
\partial_{n-1}^i \partial_n^j= \partial_{n-1}^{j-1} \partial_n^i
\end{equation}
for $1 \le i < j \le n+1$. We write 
\begin{equation}\label{module theoretic boundary2}
\partial_{n-1}\partial_n= \sum_{i=1}^n\sum_{j=1}^{n+1} (-1)^{i+j}\partial^i_{n-1} \partial^j_n= \sum_{1 \le i < j \le n+1}(-1)^{i+j}\partial^i_{n-1} \partial^j_n + \sum_{1 \le j \le i \le n} (-1)^{i+j}\partial^i_{n-1} \partial^j_n.
\end{equation}
Then, using \eqref{module theoretic boundary} in \eqref{module theoretic boundary2}, we obtain
$$\partial_{n-1}\partial_n=  \sum_{1 \le i < j \le n+1}(-1)^{i+j} \partial_{n-1}^{j-1} \partial_n^i + \sum_{1 \le j \le i \le n} (-1)^{i+j}\partial^i_{n-1} \partial^j_n=0,$$
which completes the proof. $\blacksquare$ 
\end{proof}

Let $X$ be a quandle and  $A$ a right $\mathbb{Z}(X)$-module. Take $C^n(X; A)= \Hom_{\mathbb{Z}(X)} \big(C_n(X), A \big)$ and let $\partial^n$ be the coboundary map induced from $\partial_n$. This gives a cochain complex $\{C^n(X; A), \partial^n \}$ leading to cohomology groups with coefficients in $A$, which we denote by $\Ho^n_{\rm Q}(X; A)$. If $A$ is a left $\mathbb{Z}(X)$-module, then we define $C_n(X; A)= C_n(X) \otimes_{ \mathbb{Z}(X)} A$, which leads to homology groups $\Ho_n^{\rm Q}(X; A)$ with coefficients in $A$.
\para

\begin{remark}
{\rm  When $X$ is a rack, then replacing the quandle algebra by the rack algebra $\mathbb{Z}\{X\}$ in the preceding construction, we obtain the generalized rack homology  $\Ho_n^{\rm R}(X; A)$  and  the generalized rack cohomology $\Ho^n_{\rm R}(X; A)$ of $X$ with coefficients in $A$.}
\end{remark}

\begin{remark}
{\rm If $X$ is a rack, then the  $2$-cocycle  condition for a $2$-cochain $\kappa$  in this cohomology theory 
is written as 
\begin{equation}\label{2-cocycle condition for gen cohomo}
\eta_{x*y, z} (\kappa_{x,y}) + \kappa_{x*y, z} = \eta_{x*z, y*z} (\kappa_{x,z})  +  \xi_{x*z, y*z} (\kappa_{y,z})
+ \kappa_{x*z,y*z}
\end{equation}
for any $x,y,z \in X$. Here, $\eta_{x, y}$ and $\xi_{x, y}$ corresponds to the right $\mathbb{Z}(X)$-module $A$. Further, 2-cocycles $\kappa$ and $\kappa'$ are cohomologous if there exists a map $f : X \to A$ such that 
\begin{equation}\label{2-coboundary condition for gen cohomo}
\kappa'_{x,y}- \kappa_{x,y}= \eta_{x,y} \big(f(x) \big)-f(x*y) +\xi_{x, y} \big(f(y) \big)
\end{equation}
for any $x,y,z \in X$. We refer to equation \eqref{2-cocycle condition for gen cohomo} as the {\it generalized $2$-cocycle condition}. If $X$ is a quandle, then $\kappa $ further satisfies the condition $\kappa_{x,x}=0$ for all $x \in X$.} 
\end{remark}

\begin{remark}
{\rm The (co)homology theory defined in this section recovers the  twisted quandle (co)homology theory introduced in \cite{MR1885217} by taking $\Lambda=\mathbb Z [T, T^{-1}]$ as coefficients, where $$\eta_{x,y} (a)=Ta \quad \textrm{and} \quad \xi_{x,y} (b) = (1-T) (b)$$ for all $x,y \in X$. In this case, the boundary map is explicitly given by 
\begin{eqnarray*}
	\partial^T_{n}(x_1, x_2, \dots, x_n)  \nonumber & = &
	\sum_{i=1}^{n} (-1)^{i} \Big( T (x_1, x_2, \dots, x_{i-1}, x_{i+1},\dots, x_n) \\
	& &
	- (x_1 \ast x_i, x_2 \ast x_i, \dots, x_{i-1}\ast x_i, x_{i+1}, \dots, x_n) \Big)
\end{eqnarray*}
for $n \geq 2$ and $\partial^T_n=0$ for $n \leq 1$.}
\end{remark}
\bigskip
\bigskip

\subsection{Relation between quandle cohomology and group cohomology}

In \cite{MR1948837}, a group-theoretical interpretation of the second cohomology of racks is given, which we discuss next.  Let $X$ be a connected quandle and  $x_0\in X$ a fixed element. Consider the unique surjective group homomorphism  $$\delta: \Adj(X) \to \mathbb{Z}$$ given by $\delta(\textswab{a}_{x})=1$ for each $x \in X$. By Proposition \ref{adjoint group as kernel}, we have
$$\Adj(X) \cong \ker (\delta) \rtimes \langle \textswab{a}_{x_0}\rangle \cong [\Adj(X),\Adj(X)] \rtimes \mathbb{Z}.$$ 
Since $X$ is connected, for $x, y \in X$, there exists $g \in \Adj(X)$ such that $g \cdot x=y$. Let $h=\textswab{a}_y^{-\delta(g)}g$. Then we see that $\delta(h)= \delta(\textswab{a}_y^{-\delta(g)})=0$, and hence $h \in \ker (\delta)=[\Adj(X),\Adj(X)]$. Also, we have $$h \cdot x= (\textswab{a}_y^{-\delta(g)}g) \cdot x= \textswab{a}_y^{-\delta(g)}\cdot y=y.$$ Thus, the action of $\Adj(X)$ on $X$ restricts to a transitive action of $[\Adj(X),\Adj(X)]$ on $X$. 
\para 
Let $N_0$ denote the stabilizer of $x_0$ under the action of $[\Adj(X),\Adj(X)]$ on $X$. Let $A$ be an abelian group and $\Map(X,A)$ the left $\Adj(X)$-module of maps $X \rightarrow A$, where 
\begin{equation}\label{adj action on set of maps}
(\textswab{a}_x \cdot f)(y)=f(y*^{-1} x)
\end{equation}
for $x,y\in X$ and $f\in  \Map(X,A)$.  The proof of the following result is a direct check \cite[Proposition 5.1]{MR1948837}.

\begin{proposition}\label{shift} 
Let $X$ be a rack and $A$ a trivial $\Adj(X)$-module. Then there is a natural isomorphism of cochain complexes
 $J: C^*(X;A)\to C^*\big(X; \Map(X,A)\big)$, where for each $n \ge 1$, $J_n: C^n(X;A)\to C^{n-1} \big(X; \Map(X,A)\big)$ is given by
$$(J_n\phi)(x_1,\ldots,x_{n-1})(x)=\phi(x_1,\ldots,x_{n-1}, x).$$ In particular, $J$ induces an isomorphism $\Ho^n_{\rm R}(X;A)\cong \Ho^{n-1}_{\rm R} \big(X; \Map(X,A)\big)$ of cohomology groups. 
\end{proposition}

Note that Proposition \ref{shift} fails when the action of $\Adj(X)$ on $A$ is non-trivial.  Using the preceding result, Etingof and Gra\~na proved the following result \cite[Proposition 5.3 and Corollary 5.4]{MR1948837}. 

\begin{theorem}\label{relation rack cohomo and group cohomo}
Let $X$ be a rack and $A$ an abelian group. Then the following  assertions hold:
\begin{enumerate}
\item $\Ho^1_{\rm Grp}\big(\Adj(X); A\big) \cong \Ho^1_{\rm R}(X; A)$.
\item If $A$ is a trivial $\Adj(X)$-module, then 
$$\Ho^1_{\rm Grp} \big(\Adj(X); \Map(X,A) \big)\cong \Ho^2_{\rm R}(X;A)$$
where $\Map(X,A)$ is the left $\Adj(X)$-module of maps from $X$ to $A$.
\end{enumerate}
Here, $\Ho^*_{\rm Grp}$ is the group cohomology of $\Adj(X)$ with coefficients in the left $\Adj(X)$-module $A$ for (1) and with coefficients in the left $\Adj(X)$-module $\Map(X,A)$ for (2). And, $\Ho^*_{\rm R}$ is the rack cohomology of $X$ (in the sense of \eqref{Andruskiewitsch grana cohomology}) with coefficients in a left $\mathbb{Z}(X)$-module $A$ (in the sense of  \eqref{module over rack algebra via Adj}) with respect to the given $\Adj(X)$-module structure on $A$.
\end{theorem}

\begin{proof}
We first prove assertion (1). Let $C^* \big(\Adj(X);A \big)$ be the standard cochain complex that gives the group cohomology of the group $\Adj(X)$ with coefficient in a left $\Adj(X)$-module $A$. Let  $i_X: X\rightarrow \Adj(X)$, given by $i_X(x)= \textswab{a}_x$, be the natural map. Let $i_X^*: C^1 \big(\Adj(X);A \big)\rightarrow  C^1(X;A)$ be the homomorphism induced by $i_X$, that is, $i_X^*(f)= f\,i_X$ for $f \in C^1 \big(\Adj(X);A \big)$.
\para 
We claim that $i_X^*$ restricts to a bijection from $Z^1 \big(\Adj(X);A \big)$ onto $Z^1(X;A)$. For $ g \in \Adj(X)$ and $a \in A$, let $g \cdot a$ denote the given left  group action of  $\Adj(X)$ on $A$.  Let $f \in Z^1 \big(\Adj(X);A \big)$ be a group $1$-cocycle. Then, by \eqref{group coboundary 2}, we have $f(gh)=f(g)+g \cdot f(h)$ for all $g,h \in \Adj(X)$. For $x,y \in X$, we see that
\begin{eqnarray*}
i_X^*(f) (x*y) + \textswab{a}_{x*y} \cdot i_X^*(f) (y)&=& f(\textswab{a}_{x*y}) + \textswab{a}_{x*y} \cdot f(\textswab{a}_y)\\
&=& f(\textswab{a}_{x*y} \textswab{a}_{y}), \quad \textrm{since $f$ is a group 1-cocycle}\\
&=& f(\textswab{a}_{y} \textswab{a}_{x}), \quad \textrm{since $\textswab{a}_{x*y} \textswab{a}_{y}=\textswab{a}_{y} \textswab{a}_{x}$ in $\Adj(X)$}\\
&=& f(\textswab{a}_{y})+  \textswab{a}_{y} \cdot f(\textswab{a}_{x}), \quad \textrm{since $f$ is a group 1-cocycle}\\
&=& i_X^*(f)(y)+  \textswab{a}_{y} \cdot i_X^*(f)(x),
\end{eqnarray*}
and hence $i_X^*(f)$ is a rack 1-cocycle. 
\para 
Next, suppose that $i_X^*(f)=i_X^*(h)$ for $f, h \in Z^1 \big(\Adj(X);A \big)$. By definition, we have $f(\textswab{a}_x)=h(\textswab{a}_x)$ for all $x \in X$. Recall that, a map $\gamma: \Adj(X) \rightarrow  A$ is a group $1$-cocycle if and only if the map $\hat{\gamma}:\Adj(X) \rightarrow \Adj(X) \ltimes A$ given by $\hat{\gamma}(g)=  (g,\gamma(g))$ is a group homomorphism. This gives $\hat{f}(\textswab{a}_x)= \big(\textswab{a}_x, f(\textswab{a}_x) \big)= \big(\textswab{a}_x, h(\textswab{a}_x)\big)=\hat{h}(\textswab{a}_x)$. Since $\hat{f}$ and $\hat{h}$ are group homomorphisms, it follows that $\hat{f}=\hat{h}$, and hence $f=g$, which shows that $i_X^*:Z^1 \big(\Adj(X);A \big) \to Z^1(X;A)$ is injective.
\para 

Let $f\in Z^1(X;A)$ be a rack $1$-cocycle, that is,
\begin{equation}\label{rack cocycle with andruskiewich action}
f(x*y)= \textswab{a}_{y} \cdot f(x) + f(y)-\textswab{a}_{x*y} \cdot f(y),
\end{equation}
for $x, y \in X$. The map $f$ induces a map $\tilde{f}: X\rightarrow  \Adj(X) \ltimes A$ given by $\tilde{f}(x)=\big(\textswab{a}_x,\, f(x) \big)$. We show that $\tilde{f}$ further extends to a group homomorphism $\tilde{f}:\Adj(X) \rightarrow \Adj(X) \ltimes A$. Recall that, $\Adj(X)$ has a presentation
$$\Adj(X)= \big\langle \textswab{a}_{x},~ x\in X \mid \textswab{a}_{x*y} \textswab{a}_y=\textswab{a}_y \textswab{a}_x ~\textrm{for all}~x, y \in X \big\rangle.$$ For $x, y \in X$, we compute
\begin{eqnarray*}
\tilde{f}(\textswab{a}_{x*y}) \tilde{f}(\textswab{a}_{y}) &=& \big(\textswab{a}_{x*y},~f(x*y) \big) \big(\textswab{a}_{y},~f(y) \big)\\
&=& \big(\textswab{a}_{x*y}\textswab{a}_{y},~f(x*y) + \textswab{a}_{x*y} \cdot f(y)\big)\\
&=&  \big(\textswab{a}_{y}\textswab{a}_{x},~f(y) + \textswab{a}_{y} \cdot f(x)\big)
\quad \textrm{by}\quad \eqref{rack cocycle with andruskiewich action}\\
&=& (\textswab{a}_{y}, ~f(y)) (\textswab{a}_{x}, ~f(x))\\
&=& \tilde{f}(\textswab{a}_{y}) \tilde{f}(\textswab{a}_{x}),
\end{eqnarray*}
and hence $\tilde{f}$ is a group homomorphism. Thus, it follows that the map $\pi \tilde{f}:\Adj(X) \rightarrow  A$ is a group 1-cocycle, where $\pi: \Adj(X) \ltimes A \to A$ is the projection onto the second factor. Observing that 
$$i_X^*(\pi \tilde{f})(x)= \pi \tilde{f} i_X(x)= \pi \tilde{f} (\textswab{a}_x)= \pi \big(\textswab{a}_x, f(x) \big)= f(x)$$ for all $x \in X$, it follows that $i_X^*:Z^1 \big(\Adj(X);A \big) \to Z^1(X;A)$  is surjective, and hence bijective.
\para 

Let $\lambda \in B^1 \big(\Adj(X);A \big)$ be a group 1-coboundary. We claim that $i_X^*(\lambda)$ is a rack 1-coboundary.  By definition, there exist $a \in A$ such that $\lambda(g)=g \cdot a -a$ for all $g \in \Adj(X)$. Recall from the construction of rack homology \eqref{Andruskiewitsch grana cohomology at 0} that, $\partial_0:C_1(X) \to C_0(X)$ is given by $$\partial_0(x)=-\xi_{x*^{-1}x_0 ,x_0}(x_0)= \textswab{a}_x \cdot x_0- x_0$$ for all $x \in X$, where $C_0(X)$ is the free $\mathbb{Z}(X)$-module over $\{x_0\}$ for some fixed $x_0 \in X$. Define $\alpha \in C^0(X; A)=\Hom_{\mathcal{G}} \big(C_0(X), A \big)$ by $\alpha(x_0)=a$. Then we see that
$$\partial^0(\alpha)(x)=\alpha(\partial_0(x))=\alpha(\textswab{a}_x \cdot x_0- x_0)=  \textswab{a}_x \cdot \alpha(x_0)- \alpha(x_0)=\textswab{a}_x \cdot a-a= \lambda i_X(x)= i_X^*(\lambda)(x)$$
for all $x \in X$. Hence, $i_X^*(\lambda)$ is a rack  1-coboundary, and we conclude that the map $i_X^*: C^1 \big(\Adj(X);A \big)\rightarrow  C^1(X;A)$ induces a group homomorphism $ \Ho^1_{\rm Grp}\big(\Adj(X);A \big) \rightarrow  \Ho^1_{\rm R}(X;A)$. Since $i_X^*:Z^1 \big(\Adj(X);A \big) \to Z^1(X;A)$ is already shown to be bijective, the proof of assertion (1) is complete. Finally, Propositions \ref{shift} and assertion (1) imply assertion (2).  $\blacksquare$
\end{proof}

In \cite{MR3671570}, Garc\'{i}a Iglesias and Vendramin established that, for a finite quandle $X$, the depicted diagram commutes, and its columns are exact.
\begin{tiny}
$$
\xymatrix{
	& 1\ar[d] & 1\ar[d]\\
	& \Ho^1_{\rm Grp} \big(\langle \textswab{a}_{x_0} \rangle ; \Map(X,A)^{[\Adj(X), \Adj(X)]} \big)\ar[r]^{\cong}\ar[d]_{\inf} & A\ar[d]^{\iota}\\
	\Ho^2_{\rm Q}(X;A)\ar[r]^{\cong\hspace*{1cm}} & \Ho^1_{\rm Grp} \big( \Adj(X);\Map(X,A)\big)
	\ar[d]_{\mbox{res}} & A\times \Ho^1_{\rm Grp}\big([\Adj(X), \Adj(X)];\Map(X,A)\big)\ar[d]^{\pi}\\
	& \Ho^1_{\rm Grp}\big([\Adj(X), \Adj(X)];\Map(X,A)\big)^{\langle \textswab{a}_{x_0} \rangle}\ar@{=}[r]\ar[d] & \Ho^1_{\rm Grp}\big([\Adj(X), \Adj(X)];\Map(X,A)\big)\ar[r]^{\hspace*{1cm}\cong}\ar[d] & \Hom_{\mathcal{G}} (N_0,A)\\
	& 1 & 1
}
$$
\end{tiny}
Here, $\mbox{inf}$ and $\mbox{res}$ denote the inflation and restriction maps in group cohomology, respectively; and $\iota$ and $\pi$ denote the canonical inclusion and projection maps, respectively. Recall that, $N_0$ is the stabilizer of $x_0$ under the action of $[\Adj(X), \Adj(X]$ on $X$.  Let the map $$\Ho^1_{\rm Grp}\big( \Adj(X); \Map(X,A)\big)\to \Hom_{\mathcal{G}} (N_0,A)$$ deduced from the preceding diagram be denoted by $f\mapsto f_0$. Using Proposition~\ref{shift} and the preceding diagram, Garc\'{i}a Iglesias and Vendramin gave an explicit description of $2$-cocycles of finite connected quandles with coefficients in an abelian group \cite[Theorem 1.1]{MR3671570}.

\begin{theorem}
Let $X$ be a finite connected quandle, $x_0 \in X$ a fixed element and $A$ an abelian group with trivial $\Adj(X)$-action. Then the map $\phi \mapsto \big(\phi(x_0,x_0),(J\phi)_0 \big)$ induces an isomorphism
\begin{equation}
	\label{eqn:map}
	\Ho^2_{\rm Q}(X;A)\cong  A\times\Hom_{\mathcal{G}}  (N_0,A).
\end{equation}
\end{theorem}

In particular, the preceding result shows that non-constant 2-cocycles on $X$ are controlled by a finite group. A new proof of the following result of Eisermann \cite[Theorem 1.12]{MR3205568} is also deduced in \cite{MR3671570}.

\begin{theorem}
Let $X$ be a finite connected quandle and $x_0 \in X$ a fixed element. Then
$$\Ho_2^{\rm Q}(X; \mathbb{Z})\cong \big([\Adj(X), \Adj(X] \cap \C_{\Adj(X)}(\textswab{a}_{x_0}) \big)_{\rm ab} \cong (N_0)_{\rm ab}$$
where $N_0$ is the stabilizer of  $x_0$ under the action of $[\Adj(X), \Adj(X]$ on $X$.
\end{theorem}

We conclude this subsection with a beautiful homological characterisation of the unknot due to Eisermann\cite[Theorem 1]{MR1954330}.

\begin{theorem}
Let $Q(K)$ be the knot quandle of a knot $K$. Then $\Ho_2^{\rm Q} \big(Q(K); \mathbb{Z} \big)=1$ if and only if $K$ is trivial.
\end{theorem}
\bigskip
\bigskip


\subsection{Extension theory of quandles}\label{sec extension theory of quandles} 
An extension theory of quandles using quandle cocycles was initially introduced in \cite{MR1885217, MR1973510}. A far reaching generalization of quandle extensions using dynamical 2-cocycles, which encompassed the previous constructions, was given in \cite{MR1994219}. The following result presents a further extension of this framework, as proved in \cite[Proposition 3.1]{MR4282648}.

\begin{proposition}\label{set-cocycle}
Let $X$ and $S$ be two sets, $\alpha: X \times X \to \Map(S \times S, S)$ and $\beta: S \times S \to \Map(X \times X, X)$ two maps. Then the set $X \times S$ with the binary operation 
\begin{equation}\label{genralised-quandle-operation}
(x, s)* (y,t)= \big( \beta_{s, t}(x, y), ~\alpha_{x, y}(s, t) \big)
\end{equation}
forms a quandle if and only if the following conditions hold:
\begin{enumerate}
\item $\beta_{s, s}(x, x)=x$ and $\alpha_{x, x}(s, s)=s$ for all $x \in X$, $s \in S$;
\item  for each $(y, t) \in X \times S$, the map $(x, s) \mapsto  \big( \beta_{s, t}(x, y),~ \alpha_{x, y}(s, t) \big)$ is a bijection;
\item $\beta_{\alpha_{x, y}(s, t), u}\Big(\beta_{s, t}(x, y),~ z \Big)= \beta_{\alpha_{x, z}(s, u), \alpha_{y, z}(t, u)}\Big(\beta_{s, u}(x, z), ~\beta_{t, u}(y, z) \Big)$ 
\item[] and
\item[] $\alpha_{\beta_{s, t}(x, y), z}\Big(\alpha_{x, y}(s, t), ~u \Big)= \alpha_{\beta_{s, u}(x, z), \beta_{t, u}(y, z)}\Big(\alpha_{x, z}(s, u), ~\alpha_{y, z}(t, u) \Big)$.
\end{enumerate}
\end{proposition}

\begin{proof}
The idempotency axiom holds if and only if $\beta_{s, s}(x, x)=x$ and $\alpha_{x, x}(s, s)=s$ for all $x \in X$, $s \in S$. The invertibility of the right multiplication axiom is equivalent to condition (2). Finally, if $(x, s), (y, t), (z, u) \in X \times S$, then 
\begin{eqnarray*}
&& \big((x, s)*(y, t) \big)* (z, u) \\
&=& \big( \beta_{s, t}(x, y), ~\alpha_{x, y}(s, t) \big) *(z, u)\\
&=& \Big(\beta_{\alpha_{x, y}(s, t), u}\big(\beta_{s, t}(x, y), z \big),~ \alpha_{\beta_{s, t}(x, y), z}\big(\alpha_{x, y}(s, t), u \big) \Big)
\end{eqnarray*}
and
\begin{eqnarray*}
&& \big((x, s)*(z, u) \big)* \big((y, t)*(z, u) \big)\\
&=& \big( \beta_{s, u}(x, z), ~\alpha_{x, z}(s, u) \big) *\big( \beta_{t, u}(y, z), ~\alpha_{y, z}(t, u) \big)\\
&=& \Big(\beta_{\alpha_{x, z}(s, u), ~\alpha_{y, z}(t, u)}\big(\beta_{s, u}(x, z), ~\beta_{t, u}(y, z) \big), ~\alpha_{\beta_{s, u}(x, z), ~\beta_{t, u}(y, z)}\big(\alpha_{x, z}(s, u), ~\alpha_{y, z}(t, u) \big)\Big).
\end{eqnarray*}
Thus, the right distributivity axiom of the quandle is equivalent to condition (3).
$\blacksquare$    \end{proof}          

We denote the quandle obtained in Proposition \ref{set-cocycle} by $X~ _{\alpha}\times_{\beta} S$. 

\begin{corollary}
Let $X$ and $S$ be two sets, $\alpha: X \times X \to \Map(S \times S, S)$ and $\beta: S \to \Map(X \times X, X)$ two maps. Then the set $X \times S$ with the binary operation 
$$(x, s)* (y,t)= \big( \beta_{t}(x, y), \alpha_{x, y}(s, t) \big)$$ forms a quandle if and only if the following conditions hold:
\begin{enumerate}
\item $\beta_{s}(x, x)=x$ and $\alpha_{x, x}(s, s)=s$ for all $x \in X$, $s \in S$;
\item $\beta_{s}(-, y): X \to X$ and  $\alpha_{x, y}(-, t): S \to S$ are bijections for all $x, y \in X$, $s \in S$;
\item $\beta_{u}\big(\beta_{t}(x, y), ~z \big)= \beta_{\alpha_{y, z}(t, u)}\big(\beta_{u}(x, z), ~\beta_{u}(y, z) \big)$ 
\item[] and
\item[] $\alpha_{\beta_{t}(x, y), z}\big(\alpha_{x, y}(s, t), ~u \big)= \alpha_{\beta_{u}(x, z), ~\beta_{u}(y, z)}\big(\alpha_{x, z}(s, u), ~\alpha_{y, z}(t, u) \big)$.
\end{enumerate}
\end{corollary}
\medskip

A special case of Proposition \ref{set-cocycle} was considered by Andruskiewitsch and Gra\~{n}a \cite[Section 2.1]{MR1994219} when $X$ is a quandle and $\beta_{s, t}(x, y)=x*y$ for all $x, y \in X$ and $s, t \in S$. In this case, the map $\alpha: X \times X \to \Map(S \times S, S)$ satisfy
\begin{equation}\label{dynamical-cocycle-condition1a}
\alpha_{x, x}(s, s)=s,
\end{equation}
\begin{equation}\label{dynamical-cocycle-condition2a}
\alpha_{x, y}(-, t): S \to S~\textrm{is a bijection}
\end{equation}
and the  {\it cocycle condition}
\begin{equation}\label{dynamical-cocycle-condition3a}
\alpha_{x*y, z}\big(\alpha_{x, y}(s, t), ~u \big)= \alpha_{x* z, y*z}\big(\alpha_{x, z}(s, u),~ \alpha_{y, z}(t, u) \big)
\end{equation}
for all $x, y \in X$ and $s, t \in S$. Such an $\alpha$ is called a {\it dynamical 2-cocycle}. Further, in this case, the quandle operation \eqref{genralised-quandle-operation}  on $X \times S$ becomes 
\begin{equation}\label{dynamical-quandle-operation1a}
(x, s)* (y,t)= \big( x* y, ~\alpha_{x, y}(s, t) \big).
\end{equation} 
The quandle so obtained is called a  \index{dynamical extension}{\it dynamical extension} of $X$ by $S$ through $\alpha$, and is denoted by $X \times_{\alpha} S$.

\begin{corollary}\label{isomorphic-s-quandles}
Let $X$ be a connected quandle. Then the following  assertions hold:

\begin{enumerate}
\item For each $x \in X$, the set $S$ forms a quandle be defining
$$s *_x t :=\alpha_{x,x}(s,t)$$ for all $s, t \in S$.
\item $(S, *_x) \cong (S, *_y)$ for all $x, y \in X$.
\end{enumerate}
\end{corollary}

\begin{proof}
A direct verification shows that $(S, *_x)$ is a quandle for each $x \in X$. Let $x, z \in X$ and $u \in S$. Define a map $\phi: (S, *_x) \to (S, *_{x*z})$ by 
$$\phi(s)= \alpha_{x, z}(s, u)$$ for all $s \in S$.
Then, for $s, t \in S$, we have
\begin{eqnarray*}
\phi(s *_x t) &=& \alpha_{x, z}(s *_x t, u)\\
&=& \alpha_{x, z}\big( \alpha_{x, x}(s, t), u\big)\\
&=&\alpha_{x* z, x*z}\big(\alpha_{x, z}(s, u),~ \alpha_{x, z}(t, u) \big),~\textrm{by \eqref{dynamical-cocycle-condition3}}\\
&=&\phi(s) *_{x*z} \phi(t),
\end{eqnarray*}
and $\phi$ is an isomorphism of quandles. Since $X$ is connected, it follows that $(S, *_x) \cong (S, *_y)$ for all $x, y \in X$.
$\blacksquare$    
\end{proof}

\begin{example} \label{example of dynamical cocycle eta xi kappa}
	{\rm
		Let $A$ be  a $\mathbb{Z}(X)$-module for the quandle $X$,
		and $\kappa$ a generalized $2$-cocycle. 
		For $a,b \in A$, let 
		$$\alpha_{x,y}(a,b) =  \eta_{x,y}(a)  +  \xi_{x,y}(b)
		+ \kappa_{x,y} .$$ 
		Then it can be verified directly that $\alpha$ is a dynamical cocycle. 
		In particular, even with $\kappa=0$, 
		a $\mathbb{Z}(X)$-module structure on the abelian 
		group $A$ defines a quandle structure
		$A\times_\alpha X$. }
\end{example}

If $A$ is an abelian group, then a normalized quandle 2-cocycle is a map $\alpha: X \times X \to A$ satisfying
$$
\alpha_{x, y}~\alpha_{x*y, z}= \alpha_{x, z}~\alpha_{x* z, y*z} \quad \textrm{and} \quad \alpha_{x, x}=1
$$
for all $x, y, z \in X$. A normalized quandle 2-cocycle $\alpha: X \times X \to A$ gives rise to a dynamical 2-cocycle $\alpha': X \times X \to \Map(A \times A, A)$ given as $$\alpha'_{x, y}(s, t)=s \,\alpha_{x, y}.$$ In this case, the quandle $X \times_\alpha A$ is called the {\it abelian extension} of $X$ by $A$ through $\alpha$. Such extensions appeared first in \cite{MR1973510}.

\begin{proposition}\cite[Corollary 2.5]{MR1994219}\label{AGlemma}   
Let $p: Y \rightarrow X$ be a surjective quandle homomorphism  such that the cardinality of $p^{-1}(x)$ is constant for all $x \in X$. Then $Y$ is isomorphic to a dynamical  extension  $X \times_{\alpha} S$ of $X$ by $S$ through some dynamical 2-cocycle $\alpha$, where $S$ is a set with $|S|=|p^{-1}(x)|$.
\end{proposition}

\begin{proof}
Since the cardinality of $p^{-1}(x)$ is constant for each $x \in X$, we can choose a set $S$ and bijections $g_x: p^{-1}(x) \to S$ for each $x \in X$. Define $\alpha: X \times X \to \Map(S \times S, S)$ by 
$$\alpha_{x_1, x_2}(s_1,s_2)=g_{x_1*x_2} \big(g_{x_1}^{-1}(s_1) * g_{x_2}^{-1}(s_2)\big)$$ for $x_1, x_2 \in X$ and $s_1, s_2 \in S$. Clearly, $\alpha$ satisfy \eqref{dynamical-cocycle-condition1a} and \eqref{dynamical-cocycle-condition2a}. Further, for  $x_1, x_2, x_3 \in X$ and $s_1, s_2, s_3 \in S$, we have
\begin{eqnarray*}
\alpha_{x_1*x_2, x_3} \big(\alpha_{x_1,x_2}(s_1, s_2), s_3 \big)&=& \alpha_{x_1*x_2, x_3} \big(g_{x_1*x_2}(g_{x_1}^{-1}(s_1) * g_{x_2}^{-1}(s_2)), s_3 \big)\\
&=& g_{(x_1*x_2)*x_3}\big(g_{x_1*x_2}^{-1} (g_{x_1*x_2}(g_{x_1}^{-1}(s_1) * g_{x_2}^{-1}(s_2)))* g_{x_3}^{-1}(s_3) \big)\\
&=& g_{(x_1*x_3)*(x_2*x_3)}\big((g_{x_1}^{-1}(s_1) * g_{x_2}^{-1}(s_2))* g_{x_3}^{-1}(s_3) \big)\\
&=&g_{(x_1*x_3)(x_2*x_3)} \big((g_{x_1}^{-1}(s_1) * g_{x_3}^{-1}(s_3))* (g_{x_2}^{-1}(s_2)*g_{x_3}^{-1}(s_3))\big)\\
&=&\alpha_{x_1*x_3, x_2*x_3} \big(\alpha_{x_1,x_3}(s_1, s_3), \alpha_{x_2,x_3}(s_2, s_3) \big).
\end{eqnarray*}
Hence, $\alpha$ satisfies \eqref{dynamical-cocycle-condition3a}, and therefore it is a dynamical 2-cocycle.
\para
We can write $Y= \sqcup_{x \in X}~ p^{-1}(x)$. Define $\varphi: Y \to X \times_\alpha S$ by 
$$\varphi(y)=\big(p(y), ~g_{p(y)}(y)\big)$$
for $y \in Y$. For $y_1, y_2\in Y$, we see that 
\begin{eqnarray*}
\varphi(y_1*y_2) &=& \big( p(y_1*y_2), ~g_{p(y_1*y_2)}(y_1*y_2)\big)\\
&=& \big( p(y_1)*p(y_2), ~g_{p(y_1)*p(y_2)}(g_{p(y_1)}^{-1}(g_{p(y_1)}(y_1))*g_{p(y_2)}^{-1}(g_{p(y_2)}(y_2)))\big)\\
&=& \big( p(y_1)*p(y_2), ~\alpha_{p(y_1), p(y_2)}(g_{p(y_1)}(y_1), g_{p(y_2)}(y_2))\big)\\
&=& \big( p(y_1), g_{p(y_1)}(y_1) \big) * \big( p(y_2), g_{p(y_2)}(y_2) \big)\\
&=& \varphi(y_1)*\varphi(y_2),
\end{eqnarray*}
and $\varphi$ is a quandle homomorphism. A direct check shows that the map $\psi: X \times_\alpha S \to Y$ given by $\psi(x, s)= g_x^{-1}(s)$ is the inverse to $\varphi$, and hence it is an isomorphism of quandles.
$\blacksquare$
\end{proof}

\begin{remark}
{\rm
Most of the outcomes discussed in this section can also be applied to racks. We encourage readers to verify these assertions on their own as an exercise.
}    
\end{remark}
\bigskip
\bigskip

\subsection{From group extensions to dynamical extensions of quandles} \label{extsec}  
In this subsection, we derive dynamical extensions of quandles from extensions of groups \cite{MR2166720}. We shall utilise basic extension theory of groups, which is summarised in Section \ref{obstruction to SLB arising from RBG} of Chapter \ref{chap RBG and YBE}. 
\para

Let $1 \rightarrow N \stackrel{i}{\rightarrow} E  \stackrel{\pi}{\rightarrow} H \rightarrow 1$ be a short exact sequence of groups, where $N$ is abelian. Let $s:H \to E$ be a normalised set-theoretic section. Then there is a group 2-cocycle $\theta: H \times H \to N$ given by
$$  s(x) s(y)  s(x y)^{-1}= i \big(\theta(x, y)\big)
$$
for  all $x,  y \in H$. We denote the conjugation action of $H$ on $N$ by $(x, a) \mapsto x \cdot a$ for all $x \in H$ and $a \in A$. Recall that the set $N \times H$  is a group with the binary operation $$(a, x)(b, y)=\big(a + x \cdot b+ \theta(x, y), xy\big)$$
for all $x, y \in H$ and $a, b \in N$, and we denote this group by $N \times_\theta H$. A direct check shows that the map $\psi: E \to N \times_\theta H$ given by $\psi \big(i(a)s(x)\big)=(a, x)$ is an isomorphism of groups. With the preceding set-up, we have the following result.

\begin{proposition}   \label{qmoduleprop}
Let $1 \rightarrow N \stackrel{i}{\rightarrow} E  \stackrel{\pi}{\rightarrow} H \rightarrow 1$ be a short exact sequence of groups, where $N$ is abelian. Then the following  assertions hold:
\begin{enumerate}
\item $\Conj(E)$ is a dynamical extension of $\Conj(H)$  through a dynamical 2-cocycle $\alpha: H \times H \rightarrow \Map(N \times N, N)$.
\item The dynamical 2-cocycle $\alpha$ is given by 
$$\alpha_{x,y} (a,b ) =\eta_{x,y} (a) + \xi_{x,y}(b) + \kappa_{x,y}$$ 
for all $x, y \in H$ and $a, b \in N$, where 
$$\eta_{x,y} (a)=y  \cdot a, \quad \xi_{x,y}(b)= b-(x*y) \cdot b, \quad \kappa_{x,y}=\theta(y,x) - (yx) \cdot \theta (y^{-1}, y) + \theta(yx, y^{-1})$$
and $\theta: H \times H \to N$ is a group 2-cocycle associated with a set-theoretic section $s:H \to E$.
\end{enumerate}
\end{proposition}

\begin{proof} 
The surjective group homomorphism $\pi: E \to H$ gives a surjective quandle homomorphism $\pi: \Conj(E) \to \Conj(H)$ such that $|\pi^{-1}(x)|=|\ker(\pi)|=|N|$ for all $x \in H$. Hence, by Lemma~\ref{AGlemma}, $\Conj(E)$ is  a dynamical extension of $\Conj(H)$ by $N$ through a dynamical 2-cocycle $\alpha: H \times H \rightarrow \Map(N \times N, N)$. We use the isomorphism $E \cong N \times_\theta H$ to derive the expression of $\alpha$.  Clearly, $(0,1)$ is the identity element in $N \times_\theta H$ and 
	$$ (b, y)^{-1} = \big(-y^{-1} \cdot b - \theta(y^{-1}, y), \,y^{-1} \big) =  \big(-y^{-1} \cdot b - y^{-1} \cdot \theta (y, y^{-1}),\, y^{-1} \big)$$ 
 for any $(b,y)\in N \times_\theta H$. Expanding the quandle operation in $\Conj (N \times_\theta H)$ gives
\begin{eqnarray*} 
(a, x)* (b, y) &=& (b, y) (a, x) (b, y)^{-1} \\
&= & \big( b + y \cdot a - (yxy^{-1})\cdot b + \theta(y,x)- (yx) \cdot \theta (y^{-1}, y) + \theta(yx, y^{-1}),~ y x y^{-1} \big)\\
&= & \big( y \cdot a +b - (x*y) \cdot b + \theta(y,x)- (yx) \cdot \theta (y^{-1}, y) + \theta(yx, y^{-1}),~ x * y \big).
\end{eqnarray*}
Taking
$$\eta_{x,y} (a)=y \cdot a,  \quad \xi_{x,y}(b)= b-(x*y) \cdot b,\quad 
\kappa_{x,y}=\theta(y,x) - (yx) \cdot \theta (y^{-1}, y) + \theta(yx, y^{-1})$$
and
$$ \alpha_{x,y}(a,b)=\eta_{x,y}(a) + \xi_{x,y}(b) + \kappa_{x,y},$$
we can write 
$$(a, x)* (b, y)= \big(\alpha_{x,y}(a,b), \, x * y \big).$$
By Example \ref{example of dynamical cocycle eta xi kappa}, $\alpha$ is a dynamical 2-cocycle, which proves the proposition. $\blacksquare$
\end{proof}

\begin{example} \label{wreathex} {\rm 
The wreath product of groups gives  specific examples. Let $G$ be any abelian group and $N=G^n$ for some integer $n \ge 2$. Then the  symmetric group  $H=\Sigma_n$ acts on $N$  by  permutation of factors, that is,
$$\sigma(x_1,\ldots, x_n)= \big(x_{\sigma^{-1} (1)}, \ldots, x_{\sigma^{-1} (n)} \big)$$  for $\sigma \in \Sigma_n$ and $(x_1,\ldots, x_n) \in N$. In this case, the semi-direct product $E=N \rtimes H= G \wr \Sigma_n$ is  called
a wreath product. Since the extension $1 \to N \to E \to H \to 1$ splits, we can choose a set-theoretic section for which $\theta=0$, and hence $\kappa=0$.
} 
\end{example} 

Recall that, in case of groups, the group $2$-cocycle $\theta$ is an obstruction to the section $s:H \to E$ being  a group homomorphism. There is a similar  interpretation for quandle $2$-cocycles. 

\begin{proposition}\label{sectlemma}
Let $1\rightarrow N \stackrel{i}{\rightarrow} E  \stackrel{\pi}{\rightarrow} H \rightarrow 1$ be a short exact sequence of groups, where $N$ is abelian. Then the following  assertions hold:
\begin{enumerate}
\item The map $\kappa$ in Proposition~\ref{qmoduleprop} is a quandle $2$-cocycle.
\item The map $\kappa$ satisfy $s(x)*s(y)=i(\kappa_{x,y}) s(x*y)$. 
\item The map $\kappa$ satisfy  $\kappa_{x,y} = \theta(y, x) - \theta(yxy^{-1}, y )$.
\item If $\theta$ is a group 2-coboundary, then $\kappa$ is a quandle 2-coboundary.
\end{enumerate}
\end{proposition} 

\begin{proof} 
Since $\alpha$ is a dynamical 2-cocycle, it satisfies 
$$\alpha_{x*y, z}\big(\alpha_{x, y}(a, b), ~c \big)= \alpha_{x* z, y*z}\big(\alpha_{x, z}(a, c),~ \alpha_{y, z}(b, c) \big)$$
for $x, y, z \in H$ and $a, b,c \in N$. Substituting $\alpha_{x,y} (a,b ) =\eta_{x,y} (a) + \xi_{x,y}(b) + \kappa_{x,y}$ and using the relations \eqref{rack alg rel 1}-\eqref{rack alg rel 3}, we obtain
$$ \eta_{x*y,z}(\kappa_{x,y}) + \kappa_{x*y,z} = \eta_{x*y,y*z}(\kappa_{x,z}) + \xi_{x*z,y*z}(\kappa_{y,z}) + \kappa_{x*z,y*z}$$
for all $x, y, z \in H$. This shows that $\kappa$ is a quandle $2$-cocycle (see \eqref{2-cocycle condition for gen cohomo}), which is assertion (1). 
\para
  
Using the isomorphism $E \cong N \times_\theta H$, we can write $s(x)=(0, x)$ and $i(a)=(a, 1)$. Then we have
$$s(x)*s(y)=(0,x)*(0,y)= \big(\alpha_{x,y}(0,0), x*y \big)=(\kappa_{x,y}, x*y)= i(\kappa_{x,y}) s(x*y),$$
which is assertion (2).
\para 
By assertion (2), we have $s(y) s(x) s(y)^{-1} = i(\kappa_{x,y}) s(yxy^{-1} )$. This gives
$$ i \big(\theta(y,x)\big) s(yx) =i (\kappa_{x,y}) s(yxy^{-1} ) s(y) =i(\kappa_{x,y})  i \big(\theta(yxy^{-1}, y )\big) s(yx)= i \big(\kappa_{x,y} +\theta(yxy^{-1}, y )\big)s(yx).
$$
Since $i$ is injective, we obtain the formula in assertion (3).
\para 

Suppose that $\theta$ is a group 2-coboundary. That is, there is a $1$-cochain $\gamma:H \to N$ such that 
$$ \theta(x,y) =(\partial^1_{\rm Group} \gamma) (x, y) =\gamma(xy)-\gamma(x)-x \cdot \gamma(y),$$
	where $\partial^1_{\rm group}$ denotes the group 1-coboundary operator. By assertion (3), we have
	$$\kappa_{x,y} = \gamma(x*y) -y \cdot \gamma(x) - \gamma(y) + (x*y) \cdot \gamma(y)=\gamma(x*y) -\eta_{x, y} \big(\gamma(x)\big) - \xi_{x, y} \big(\gamma(y)\big) = \big(\partial^1_{\rm Quandle} (-\gamma)\big)(x,y),$$
where  $\partial^1_{\rm quandle}$ denotes the quandle 1-coboundary operator. Hence, $\kappa$ is a quandle 2-coboundary, which is assertion (4). $\blacksquare$
\end{proof}

\begin{remark}
{\rm 
Let   $E \cong N \times_{\theta} H $ be as in the preceding result. Let $X$ be a subquandle of $\Conj(H)$
and $\tilde{X}=\pi^{-1}(X)$. Then $\tilde{X}$ is a subquandle of $\Conj(E)$  and $\pi$ induces the quandle homomorphism $\pi:\tilde{X} \rightarrow X$, which is a dynamical extension.}
\end{remark}

Let $X$ be a quandle and $A$ a right $\mathbb Z(X)$-module.  Consider the dynamical extension $A \times_{\alpha} X$ of $X$ by $A$ through the dynamical 2-cocycle  $\alpha_{x,y}=\eta_{x,y}+\xi_{x,y}$. Let  $f \in C^1(X;A)=\Hom_{\mathbb Z(X)} \big(C_1(X),A \big)$ be a quandle $1$-cochain. Define $\hat{f}: X \rightarrow A \times_{\alpha} X$ by 
$\hat{f}(x)=\big(f(x), x \big)$. Since $\pi \, \hat{f}={\rm id}_X$, the map $\hat{f}$ is a set-theoretic section to the projection  $\pi: A \times_{\alpha} X \rightarrow X$.

\begin{proposition}\label{onecocylemma}
The quandle 1-cochain $f$ is a 1-cocycle if and only if the section $\hat{f}$ is a quandle homomorphism.
\end{proposition} 

\begin{proof}
It follows from the definition that $f:X \to A$ is a quandle $1$-cocycle if $$f(x*y)=\eta_{x,y} \big(f(x)\big)+\xi_{x,y} \big(f(y)\big)$$
for all $x, y \in X$. On the other hand, we see that
$$ \hat{f} (x*y) = \big(f(x*y), ~x*y \big)$$
and
$$\hat{f}(x) * \hat{f}(y) = \big(f(x), ~x \big) *\big(f(y), ~y \big) =\big(\alpha_{x,y}(f(x),f(y)), ~x*y \big) = \big(\eta_{x,y}f(x)+\xi_{x,y}f(y), ~x*y \big)$$
for all $x, y \in X$. Thus, $f$ is a quandle 1-cocycle if and only if $\hat{f}$ is a quandle homomorphism.
	$\blacksquare$
\end{proof}
\bigskip
\bigskip

\section{Homology and cohomology of general solutions via colorings}  \label{homsec colorings}

In this section, we present a homology and cohomology theory for general solutions to the Yang--Baxter equation via colorings, and give applications to knot theory.

\subsection{Yang--Baxter coloring, homology and cohomology}
We begin by delving into the geometric construction outlined in \cite{MR2128041}. This construction enables us to elucidate low-dimensional differentials and understand their origins through the associated diagrammatic representations.
\para 

Let us describe the $n$-dimensional cubes, inspired by the preferred squares approach to rack spaces introduced by Fenn, Rourke, and Sanderson  \cite{MR1257904,MR2255194, MR1364012}. Let $I = [0, 1]$ be the unit interval. For a positive integer $n$, let ${\mathcal I}_n$  be the $n$-dimensional cube $I^n$ regarded as a  \index{cubical complex}CW (cubical) complex. Denote the $k$-skeleton  of ${\mathcal I}_n$ by ${\mathcal I}_n^{(k)}$. For a positive integer $k$, every $k$-dimensional face of ${\mathcal I}_n$ is another $k$-dimensional cube. We give the orientation  to each $k$-face to be the one  defined from the order of the coordinate axes.   Every $k$-face $\sigma$ is 
regarded as having this orientation fixed, and the same $k$-face with the opposite orientation is denoted by $-\sigma$.  In particular, note that every $2$-face can be written as 
$$  \{ \epsilon_1 \} \times \cdots \times  \{ \epsilon_{i-1} \} 
\times I_i  \times \{ \epsilon_{i+1} \}  \times \cdots \times
\{ \epsilon_{j-1} \}  \times I_j  \times \{ \epsilon_{j+1} \}
\times \cdots \times  \{ \epsilon_{n}\}  $$
for some $i$ and $j$ such that $1 \leq i<  j \leq n$, where $\epsilon_k \in \{0,1 \}$ and $I_i$, $I_j$ denote the $i$-th and the $j$-th factors of a copy of $I$, respectively.  We abbreviate parentheses for the sake of  simplicity.
When the integer $0$ or $1$ is placed at the $i$-th factor,  we denote it by $0_i$ or $1_i$, respectively.

\begin{definition} 
Let $(X, r)$ be a solution to the Yang--Baxter equation such that $r(x, y)=\big(\sigma_x(y), \tau_y(x)\big)$ for all $x, y \in X$. Let $E({\mathcal I}_n) $ denote the set of edges ($1$-faces) of ${\mathcal I}_n$ with each edge oriented as above. A  \index{Yang--Baxter coloring}(Yang--Baxter) coloring of ${\mathcal I}_n$ by $(X, r)$ is a map $L: E({\mathcal I}_n) \rightarrow X$  such that if 
$$ L(  \epsilon_1  \times \cdots \times  I_i  \times \cdots \times  0_j  \times  \cdots \times   \epsilon_{n}  )  =  x \quad  \textrm{and}\quad   L(  \epsilon_1 \times \cdots \times   1_i   \times \cdots \times  I_j \times  \cdots \times   \epsilon_{n}  )= y, $$
then
$$ L(  \epsilon_1  \times \cdots \times   0_i   \times \cdots \times  I_j  \times  \cdots \times   \epsilon_{n}  ) = \sigma_x(y) \quad  \textrm{and}\quad L(  \epsilon_1  \times \cdots \times  I_i  \times \cdots \times  1_j  \times  \cdots \times   \epsilon_{n}  ) =  \tau_y(x).$$ 
\end{definition}

\begin{figure}[h]
	\begin{center}
		\mbox{
			\epsfxsize=2.3in
			\epsfbox{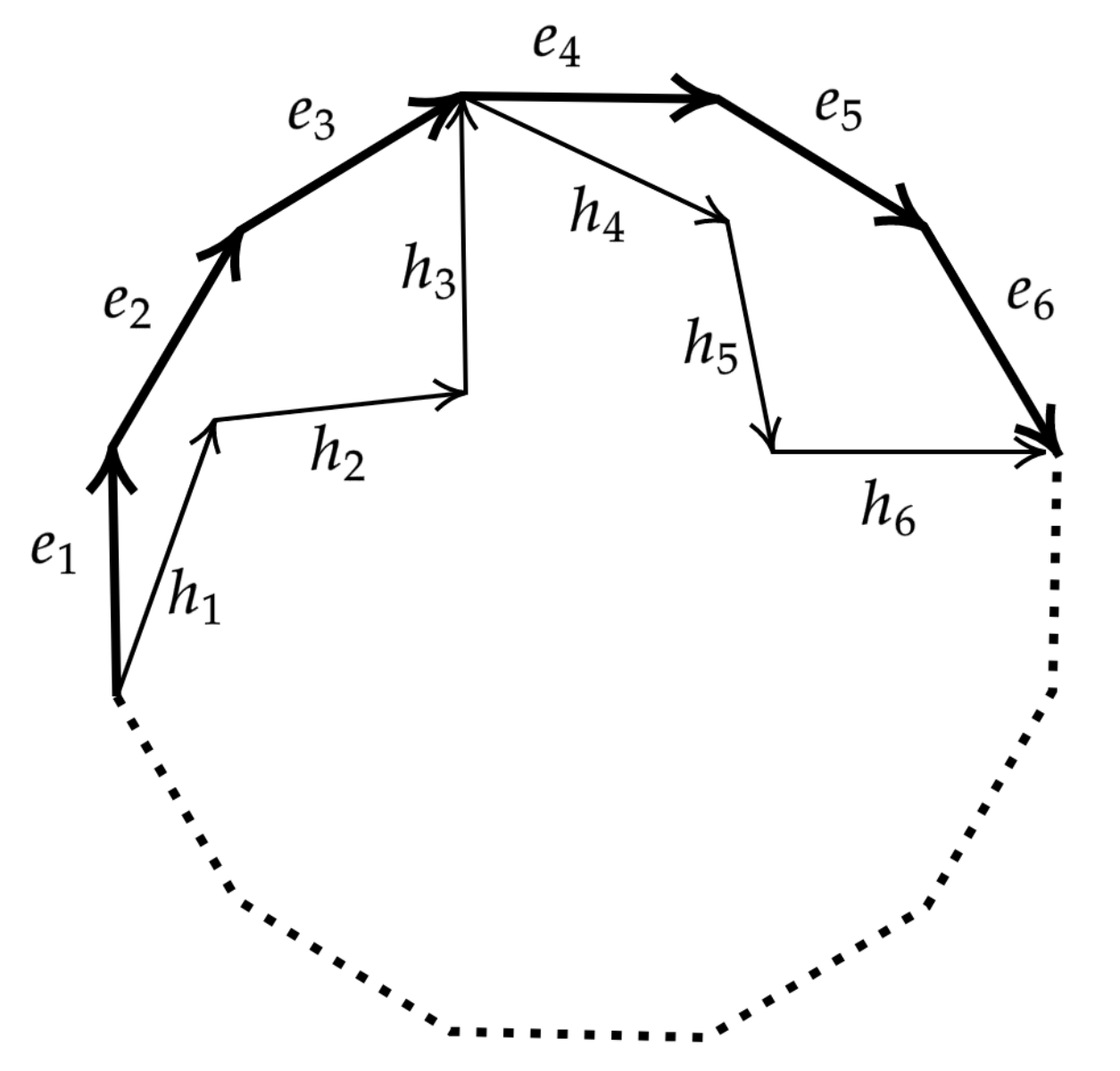} 
		}
	\end{center}
	\caption{ The initial path of ${\mathcal I}_6$. }
	\label{initialpath} 
\end{figure}

To see how many Yang--Baxter colorings  by $(X, r)$ does ${\mathcal I}_n$ admits,  we specify the initial path in  ${\mathcal I}_n$. The sequence $( e_1, \ldots, e_{n})$ of edges of  ${\mathcal I}_n$, where 
\begin{eqnarray*}
	e_1 & = & I_1 \times  0_2 \times \cdots \times  0_n, \\
	e_2 &=&  1_1  \times I_2 \times  0_3   \times \cdots \times  0_n, \\
	\vdots &  \vdots & \quad  \quad  \quad  \quad  \quad  \vdots \\
	e_n  &=&   1_1  \times \cdots \times  1_{n-1}   \times I_n ,
\end{eqnarray*}
is called the \index{initial path of $n$-cube}{\it initial path} of   ${\mathcal I}_n$. Note that the orientations of the edges $e_i$  of the initial path are consistent in the sense that the terminal point of $e_i$ is the initial point of $e_{i+1}$ for all $i=1, \ldots, n-1$.
\para 
For a sub-cube 
$$C=   \epsilon_1   \times \cdots \times 
I_{j_1}  \times \cdots \times 
I_{j_2}  \times \cdots \times 
I_{j_k}  \times \cdots \times   \epsilon_n  $$
of dimension $k$,  identify $C$ with   ${\mathcal I}_k$ by  the  obvious map sending the $h$-th factor of $I$ in $I^k$  to $I_{j_h}$ in $C$. The orientations are preserved by this identification.  Then the sequence of edges of $C$ corresponding to the initial path of  ${\mathcal I}_k$ under this identification is called  the {\it initial path} of $C$. Figure~\ref{initialpath} depicts a  projection of   ${\mathcal I}_6$ and its initial path is labeled by $e_1, \ldots, e_6$  along  the top edges. The following result assures the existence of a Yang--Baxter coloring \cite[Lemma 2.3]{MR2128041}.

\begin{lemma} \label{yblemma}
Let $(X,r)$ be a solution to the Yang--Baxter equation and $(e_1, \ldots, e_n)$
the initial path of ${\mathcal I}_n$. For any $n$-tuple $(x_1, \ldots, x_n)$ of elements of $X$,  there exists a unique 	Yang--Baxter coloring $L$ of   ${\mathcal I}_n$ by $(X,r)$ 	such that $L(e_i)=x_i$ for all $i=1, \ldots, n$. 
\end{lemma}
\para

Let $(X,r)$ be a solution to the Yang--Baxter equation. Let $C_n^{\rm YB}(X)$ be the free 
abelian group generated by $n$-tuples $(x_1, \dots, x_n)$ of elements of $X$. Consider a Yang--Baxter coloring $L$ of  ${\mathcal I}_n$  with $L(e_i)=x_i$ for all $i=1, \ldots, n$, which exists uniquely by Lemma~\ref{yblemma}. Consider any $k$-face (subcube)  ${\mathcal J}$ of ${\mathcal I}_n$.  Let $(f_1, \ldots, f_k)$ be the initial path of  ${\mathcal J}$. Then there is a unique $k$-tuple $(y_1, \ldots, y_k)$ of elements of $X$ such that $L(f_j)=y_j$ for all $j=1, \ldots, k$. We abbreviate this as $L( {\mathcal J})=(y_1, \ldots, y_k)$.
\para 

Let $\partial_n^C$ denote the $n$-dimensional  boundary map in  \index{ cubical homology theory}cubical homology theory \cite{MR0617135, MR0052766, MR0045386}. Hence, $\partial_n^C ({\mathcal I}_n) = \sum_{i=1}^{2n} \epsilon_i {\mathcal J}_i$, where  ${\mathcal J}_i$ is  an $(n-1)$-face of ${\mathcal I}_n$ and
$\epsilon_i=\pm 1$ depending on whether the orientation of $ {\mathcal J}_i$ matches the induced orientation on $ {\mathcal J}_i$. For the induced orientation, we take the convention that the inward pointing normal to an $(n-1)$ face appears last in a sequence of vectors that specifies an orientation, and the orientation of the $(n-1)$-face is chosen so that  this sequence agrees with the orientation of the $n$-cube.  For instance, the $(n-1)$-face $I_1 \times \cdots \times I_{n-1} \times \{0\}$ has  a  compatible orientation.
\para

We define a homomorphism
$$\partial_{n}: C_{n}^{\rm YB}(X) \to C_{n-1}^{\rm YB}(X)$$ by setting
$$\partial_n ( (x_1, \ldots, x_n)) = \sum_{i=1}^{2n} \epsilon_i L(C_i),$$
where each $C_i$ is the $(n-1)$-dimensional sub-cube of $\mathcal{I}_n$. Since $\partial_{n-1}^C \, \partial_{n}^C =0 $,  we have $\partial_{n-1} \, \partial_{n}=0$. Thus, we obtain a chain complex  $\{C_n^{\rm YB}(X),  \partial_n\}$ and denote its homology groups by $\Ho_*^{\rm YB}(X)$. If $A$ is any abelian group, as usual, we can define homology groups with coefficients in $A$, cochain groups and cohomology groups with coefficients in $A$. We denote these groups by  $\Ho_n^{\rm YB}(X;A)$, $C^n_{\rm YB}(X;A)$ and  $\Ho^n_{\rm YB}(X;A)$, respectively. The homology and cohomology theories defined above are called the {\it homology and cohomology theories of solution to the Yang--Baxter equation} \cite{MR2128041}. We also refer to cycles and cocycles in these homology and cohomology theories as \index{Yang--Baxter cycles}{\it Yang--Baxter cycles} and \index{Yang--Baxter cocycles}{\it Yang--Baxter cocycles}, respectively.

\begin{remark}
{\rm The cochain complex in the preceding paragraph can also be regarded as the diagonal part of Eisermann's Yang--Baxter cohomology, which controls deformations of the solution \cite{MR3240919, MR2153117}.}
\end{remark}

\begin{remark}
{\rm
In  \cite{MR3381331}, Przytycki  defined another homology theory for pre-Yang--Baxter operators, having a nice graphic visualization. Later, in \cite{MR3835755}, Przytycki and Wang showed that this homology theory is equivalent to the one due to Carter, Elhamdadi, and Saito \cite{MR2128041} discussed above.
}
\end{remark}

We present explicit formulas for the boundary homomorphisms in lower dimensions in the following examples.

\begin{figure}[hbtp]
	\centering
	\includegraphics[height=4.5cm]{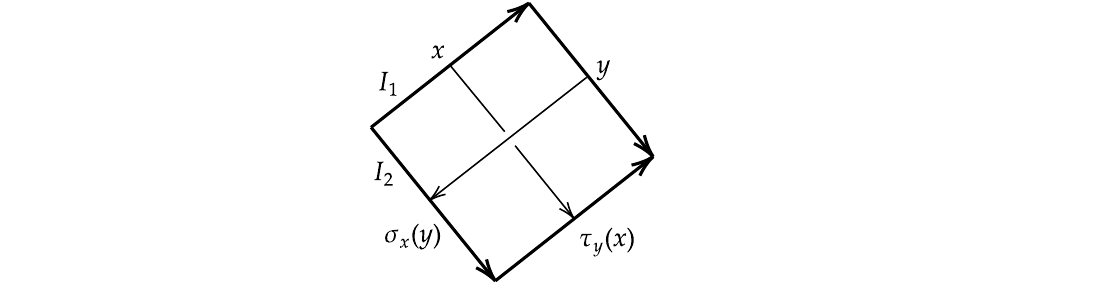}
	\caption{$2$-dimensional boundary homomorphism.}
	\label{2dboundary} 
\end{figure}

\begin{example} {\rm 
Let $(X,r)$ be a solution to the Yang--Baxter equation. Figure \ref{2dboundary} depicts the $2$-dimensional cube ${\mathcal I}_2$. The top two edges form the initial path, and the edges are colored by $x,y \in X$. The bottom edges are colored by $\sigma_x(y)$ and $\tau_y(x)$, so that these assignments indeed defines a Yang--Baxter coloring of the cube ${\mathcal I}_2$. Hence, the boundary homomorphism in this case is given by 
\begin{equation}\label{dim two boundary homo}
\partial_2 (x,y) = (x) + (y) - \sigma_x(y) - \tau_y(x).
\end{equation}
In Figure \ref{2dboundary}, a correspondence between a square and a positive crossing point used for  classical knot diagrams is depicted. We use this correspondence in the following examples.
} \end{example}

\begin{figure}[hbtp]
	\centering
	\includegraphics[height=12cm]{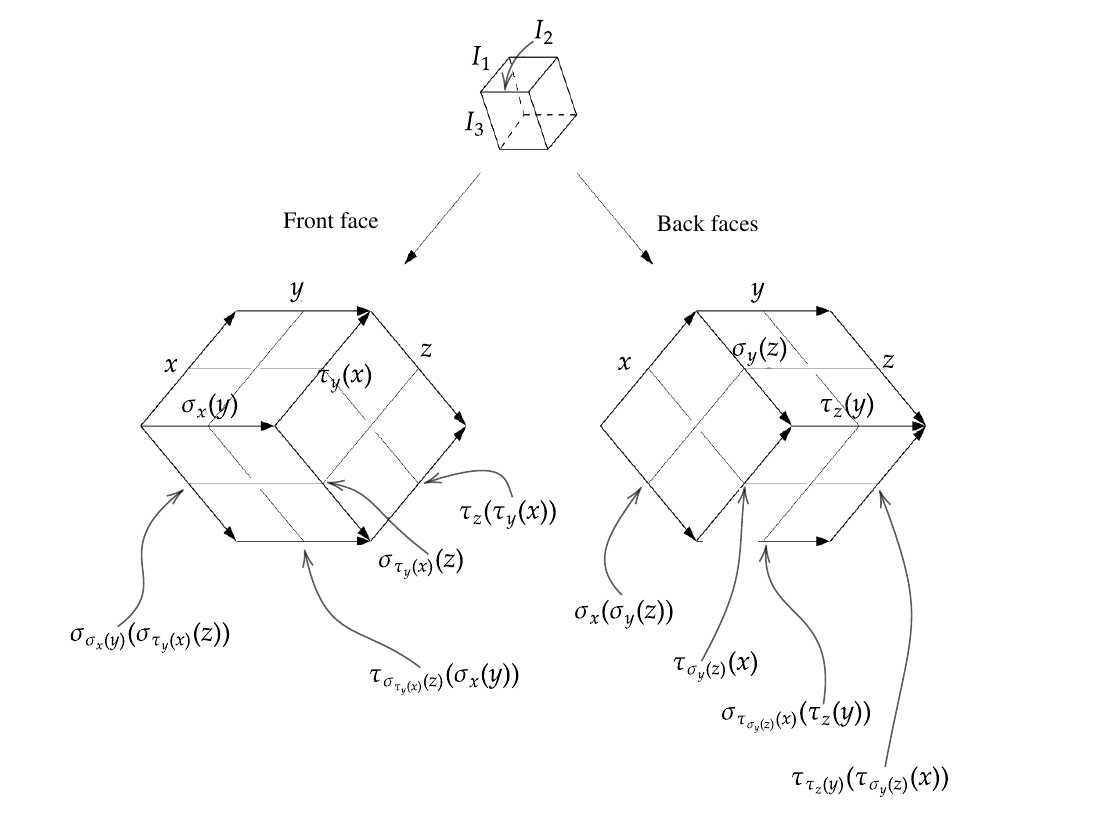}
	\caption{$3$-dimensional boundary homomorphism.}
	\label{3dboundary} 
\end{figure}

\begin{example} {\rm 
In Figure~\ref{3dboundary}, the $3$-dimensional cube ${\mathcal I}_3$ is depicted at the top of the figure. The top three edges form the initial path, to which the elements $x,y,z$ are assigned. The cube, as depicted in Figure \ref{3dboundary}, consists of  three front faces and three back faces from the reader's perspective. In the top figure of the cube, the three back faces are located behind the three front  faces, and depicted by dotted lines.  In the bottom of Figure \ref{3dboundary}, the front faces (left)  and the back faces (right) are depicted separately. The front faces determine, via the condition of a Yang--Baxter coloring, the colors assigned to all the edges of the front faces as depicted. The same is true for the edges of the back faces. There are six edges that are common for both front and back faces. Three edges of the initial path (labeled by $x,y,z$) at the top, and three edges at the bottom. The expressions of the colors assigned to these bottom three edges obtained from front faces and those obtained from back faces are different, as depicted in the figure. These are, of course, the same elements in $X$ for each edge, 
 since the map $r$ is a solution to the Yang--Baxter equation. Indeed, these  elements are the same is equivalent to the condition that $r$ is a solution to the Yang--Baxter equation (see Proposition \ref{component conditions for a solution}). The relation between cubes and crossings is again depicted in the bottom figure, and it is seen that the Yang--Baxter equation corresponds to the Reidemeister move $R3$ from knot theory. The boundary homomorphism in this case is given by 
$$\partial_3(x,y,z)   = (x,y)+\big(\tau_y(x), z \big)+ \big(\sigma_x(y), \sigma_{\tau_y(x)}(z) \big) - (y,z) - \big(x, \sigma_y(z) \big) - \big(\tau_{\sigma_y(z)}(x), \tau_z(y)\big). 
$$
} 
\end{example}
\bigskip
\bigskip

\subsection{Extensions of solutions and obstruction cocycles}
Low dimensional Yang--Baxter cocycles have interpretations as obstruction cocycles in extensions of solutions to the Yang--Baxter equation. This is similar to cohomology theories for groups and other algebraic structures, and also generalizes the quandle cohomology case as well.
\para 

Let $1 \rightarrow N \stackrel{i}{\rightarrow}  G \stackrel{p}{\rightarrow}  A \rightarrow 1$  be an exact sequence of abelian groups  with a normalised set-theoretic section $s: A \rightarrow G$. Let $n$ be a positive integer and  $f \in C^n_{\rm YB}(X;A)= \Hom_{\mathcal{G}} \big(C_{n-1}^{\rm YB}(X), A\big)$ an $n$-cocycle. Then $\partial^n f (x_1, \ldots, x_n)$ has $2n$ terms, each of which  has $(n-1)$ arguments. Let us write  $$\partial^n f (x_1, \ldots, x_n)=T_1 + \cdots + T_{2n}.$$ Consider the expression $h=s(T_1) + \cdots + s(T_{2n})$. Then $p(h)=0 \in A$, since $p \, s = \id_A$, $p$ is a homomorphism and $f$ is a cocycle. Thus, there is a unique  element, say $\psi(x_1, \ldots, x_n) \in N$,  such that  $i \, \psi(x_1, \ldots, x_n) = h$.  For example, for $n=2$, by \eqref{dim two boundary homo}, we obtain 
\begin{equation}\label{defining2} 
	sf(x) + sf(y)= i \psi (x, y) +  sf \big(\sigma_x(y)\big) + sf \big( \tau_y(x)\big).
\end{equation}

\begin{proposition} 
The map $\psi$ is an $(n+1)$-cocycle of $X$ with coefficients in $N$.
\end{proposition}
\begin{proof}
The $2n$ terms $T_1, \ldots, T_{2n}$ are in one-to-one correspondence with $(n-1)$-faces of ${\mathcal I}_n$, whose initial path is labeled by $x_1, \ldots, x_n$ via $L$ (by Lemma \ref{yblemma}). We assign $i \psi(x_1, \ldots, x_n)$ to this cube ${\mathcal I}_n$. In ${\mathcal I}_{n+1}$, whose initial path  is labeled by $(x_1, \ldots, x_{n+1})$, the $n$-faces ${\mathcal I}_n$ are in one-to-one correspondence with the assigned $i\psi$. Then $\partial^{n+1} (\psi) =\partial^{n+1} \partial^n(f)=0$ since $\partial^{n+1} \partial^n=0$.
$\blacksquare$
\end{proof}

We call such a cocycle $\psi$ an \index{obstruction cocycle}{\it obstruction $(n+1)$-cocycle}. Explicit computations can be carried out using this correspondence. For example, when $n=2$, we have
\begin{eqnarray*}
 &&  \underline{sf(x) + sf(y)} + sf(z) \\
&=&  i \psi (x, y) +  sf \big(\sigma_x(y)\big) + \underline{sf \big(\tau_y(x) \big) + sf(z)} \\
	&=&  i \psi (x, y) +  i \psi \big(\tau_y(x), z \big)  + \underline{ sf \big(\sigma_x(y) \big) + sf \big( \sigma_{\tau_y(x)}(z) \big)}  + sf \big( \tau_z(\tau_y(x)) \big) \\
	&=&  i \psi (x, y) +  i \psi \big(\tau_y(x), z \big) 	+  i \psi \big(\sigma_x(y), \sigma_{\tau_y(x)}(z) \big) \\ 
	& & 	+  sf \big(\sigma_{\sigma_x(y)}(\sigma_{\tau_y(x)}(z)) \big)	+ sf \big( \tau_{ \sigma_{\tau_y(x)}(z)}(\sigma_x(y)) \big)  +  sf \big( \tau_z( \tau_y(x)) \big)
\end{eqnarray*} 
and
\begin{eqnarray*} 
 && sf(x) +  \underline{sf(y)+ sf(z)} \\
  &=& i \psi (y,z) + \underline{sf(x) + sf \big(\sigma_y(z) \big)} + sf \big(\tau_z(y)\big) \\
	&=&  i \psi (y,z) +  i \psi \big(x, \sigma_y(z) \big) + sf \big( \sigma_x(\sigma_y(z)) \big) +  \underline{sf \big( \tau_{\sigma_y(z)}(x) \big)+sf \big( \tau_z(y)\big) } \\
	&=&  i \psi (y,z) +  i \psi \big(x, \sigma_y(z) \big) + i \psi \big( \tau_{ \sigma_y(z)}(x),  \tau_z(y)\big) \\
 & & + sf \big( \sigma_x(\sigma_y(z))\big) + sf \big( \sigma_{\tau_{ \sigma_y(z)}(x)}(\tau_z(y))\big) + sf \big( \tau_{\tau_z(y)}(\tau_{\sigma_y(z)}(x))\big),
\end{eqnarray*} 
where the defining relation \eqref{defining2} of $\psi$ is applied to the underlined terms. We recover the formula for the second coboundary homomorphism this way.  This computation  is  directly visualized from Figure \ref{3dboundary}.
\para 

\begin{proposition} \label{extprop}
Let $A$ be a ring, $(X,r)$ a solution to the Yang--Baxter equation with $r(x, y)= \big(\sigma_x(y), \tau_y(x)\big)$, and $\psi_1, \psi_2 \in  C^2_{\rm YB}(X;A)$.
Let $V=A \times X$ and $s: V \times V \rightarrow  V \times V$ be defined by 
	$$s \big( (a_1, x_1), (a_2, x_2) \big) = \big( (a_2 + \psi_1(x_1, x_2) , \sigma_{x_1}(x_2)),
	(a_1 + \psi_2(x_1, x_2),  \tau_{x_2}(x_1) ) \big)$$
	for  $(a_i, x_i) \in V$. If  $(V, s)$ is a solution to the Yang--Baxter equation, then 
	$\psi_1 + \psi_2 \in  Z^2_{\rm YB}(X;A)$.
\end{proposition}
\begin{proof} 
We compute
\begin{eqnarray}\label{sum of cochains solution 1}
&& (s \times \id )( \id \times s)(s \times \id) \big(( a_1, x_1), (a_2, x_2), (a_3, x_3) \big) \\
\nonumber &=& \Big(\Big(a_3 + \psi_1( \tau_{x_3}(x_1), x_3) + \psi_1(\sigma_{x_1}(x_2), ~\sigma_{\tau_{x_2}(x_1)}(x_3)), ~\sigma_{\sigma_{x_1}(x_2)}(\sigma_{\tau_{x_2}(x_1)}(x_3) ) \Big),~\\
\nonumber 	& & \Big( a_2 + \psi_1(x_1, x_2) + \psi_2( \sigma_{x_1}(x_2), ~\sigma_{\tau_{x_2}(x_1)}(x_3)), ~ 
 \tau_{\sigma_{\tau_{x_2}(x_1)}(x_3)}(\sigma_{x_1}(x_2)) \Big), ~\\
	\nonumber & & \Big( a_1 +  \psi_2(x_1, x_2) + \psi_2( \tau_{x_2}( x_1) , x_3), ~\tau_{x_3}(\tau_{x_2}(x_1))  \Big)  \Big)
\end{eqnarray} 
and 
\begin{eqnarray}\label{sum of cochains solution 2}
&& (\id \times s)(s \times \id)( \id \times s)\big(( a_1, x_1), (a_2, x_2), (a_3, x_3) \big) \\
\nonumber &=& \Big(  \Big(a_3 + \psi_1( x_2, x_3) + \psi_1( x_1, \sigma_{x_2}(x_3) ),~  \sigma_{x_1}(\sigma_{x_2}(x_3))\Big), \\
\nonumber & & \Big( a_2 +  \psi_2( x_2, x_3) + \psi_1(\tau_{\sigma_{x_2}(x_3)}(x_1),   \tau_{x_3}(x_2) ), ~\sigma_{\tau_{\sigma_{x_2}(x_3)}(x_1)}(\tau_{x_3}(x_2)) \Big), \\
\nonumber  & &  \Big( a_1 + \psi_2(x_1, \sigma_{x_2}(x_3) ) +   \psi_2(\tau_{ \sigma_{x_2}(x_3)}(x_1), \tau_{x_3}(x_2)), ~ \tau_{ \tau_{x_3}(x_2)}(\tau_{\sigma_{x_2}(x_3)}(x_1)) \Big) \Big).
\end{eqnarray} 
Since $(V, s)$ is a solution to the Yang--Baxter equation, equating \eqref{sum of cochains solution 1} and \eqref{sum of cochains solution 1} gives
\begin{eqnarray*}
& & \Big( \psi_1( \tau_{x_3}(x_1), x_3) + \psi_1(\sigma_{x_1}(x_2),~\sigma_{\tau_{x_2}(x_1)}(x_3)), ~\sigma_{\sigma_{x_1}(x_2)}(\sigma_{\tau_{x_2}(x_1)}(x_3) ) \Big)\\
&=&  \Big( \psi_1( x_2, x_3) + \psi_1( x_1, \sigma_{x_2}(x_3) ),~  \sigma_{x_1}(\sigma_{x_2}(x_3))\Big),\\
& &\\
 & & \Big( \psi_1(x_1, x_2) + \psi_2( \sigma_{x_1}(x_2), ~\sigma_{\tau_{x_2}(x_1)}(x_3)), ~  \tau_{\sigma_{\tau_{x_2}(x_1)}(x_3)}(\sigma_{x_1}(x_2)) \Big) \\
 &=& \Big(  \psi_2( x_2, x_3) + \psi_1(\tau_{\sigma_{x_2}(x_3)}(x_1),   \tau_{x_3}(x_2) ), ~\sigma_{\tau_{\sigma_{x_2}(x_3)}(x_1)}(\tau_{x_3}(x_2)) \Big),\\
 & &\\
& & \Big( \psi_2(x_1, x_2) + \psi_2( \tau_{x_2}( x_1) , x_3), ~\tau_{x_3}(\tau_{x_2}(x_1))  \Big) \\
 &=& \Big( \psi_2(x_1, \sigma_{x_2}(x_3) ) +   \psi_2(\tau_{ \sigma_{x_2}(x_3)}(x_1), \tau_{x_3}(x_2)), ~ \tau_{ \tau_{x_3}(x_2)}(\tau_{\sigma_{x_2}(x_3)}(x_1)) \Big).
\end{eqnarray*} 
Adding these equalities give  the $2$-cocycle  condition for $\psi_1+ \psi_2$, and the proof is complete. $\blacksquare$
\end{proof}

Note that, although $\psi_1 + \psi_2$ is a cocycle, each of $\psi_1$ and $\psi_2$ can not be a cocycle. Taking $\psi_1=\psi_2$ in Proposition~\ref{extprop}, we obtain the following corollary.

\begin{corollary} \label{extcor}
Let $A$ be a ring in which $2$ is invertible and $V=A \times X$. Let $s: V \times V \rightarrow  V \times V $ be defined by 
	$$s\big( (a_1, x_1), (a_2, x_2) \big) = \big( (a_2 + \psi(x_1, x_2) , \sigma_{x_1}(x_2)),
	(a_1 + \psi(x_1, x_2),  \tau_{x_2}(x_1) ) \big)$$
	for  $(a_i, x_i) \in V$. 	If  $(V, s)$ is a solution to the Yang--Baxter equation, then 
	$\psi \in  Z^2_{\rm YB}(X;A)$.
\end{corollary}

\begin{proof}
In Proposition~\ref{extprop}, setting $\psi_1=\psi_2$ gives twice the cocycle condition at the end. The result follows since  $2$ is invertible in $A$.
$\blacksquare$
\end{proof}

\begin{example} \label{omegaextex} {\rm 
Let $X=\mathbb Z_q [s^{\pm 1}, t^{\pm 1} ] / (1-s)(1-t)$. Then, by Example~\ref{YBEexamples}, the matrix  $\displaystyle r= \left[ \begin{array}{cc} 1-s & s \\ t & 1-t \end{array} \right]$  gives a solution $(X, r)$ to the Yang--Baxter equation. Let $A=X$ and define $\psi_i: X \times X \rightarrow A$ by  $\psi_i(x, y)=u_i(y-x)$ for some $u_i \in X$, where $i=1,2$. Let $V=A \times X$ and $s: V \times V \rightarrow V \times V$ be as defined in Proposition~\ref{extprop}. That is, 
$$s\big( (a_1, x_1), (a_2, x_2) \big) = \big( (a_2 + \psi_1(x_1, x_2) , \sigma_{x_1}(x_2) ), ~(a_1 + \psi_2(x_1, x_2),  \tau_{x_2}(x_1) )\big)
$$
for any $(a_i, x_i) \in V$. Then $s$ is represented by the matrix  
$$\displaystyle  s= \left[ \begin{array}{cccc}
0 & -u_1 & 1 & u_1 \\
0 & 1-s & 0 & s \\
1 & -u_2 & 0 & u_2 \\
0 & t & 0 & 1-t 
\end{array} \right]. $$ 
Taking $\displaystyle Y= \left[ \begin{array}{cc} 1 & u_1 \\ 0 & s \end{array} \right]$
and $\displaystyle Z= \left[ \begin{array}{cc} 1 & -u_2 \\ 0 & t \end{array} \right]$, we see that  
$s$ can be written  as $\displaystyle s= \left[ \begin{array}{cc} I-Y & Y \\ Z & I-Z \end{array} \right]$, where $Y,Z$ invertible and $I$ is the identity matrix. If $u_2(1-s)+u_1(1-t)=0$, then $YZ=ZY$ and $(I-Y)(I-Z)=0$. Thus, by Example~\ref{YBEexamples}, $(V,s)$ is a solution to the Yang--Baxter equation and is an extension of $(X,r)$ by $(\psi_1, \psi_2)$.
} \end{example}

\begin{example} \label{extex} {\rm  
 Let $\mathbb{k}$ be a ring in which $2$ is invertible. Let $X=\mathbb{k}^2$ and $\displaystyle r= \left[ \begin{array}{cc} I-Y & Y \\ Z & I-Z \end{array} \right]$, where $\displaystyle Y=\left[ \begin{array}{cc} 1 & s \\ 0 & 1 \end{array} \right]$, $\displaystyle Z=\left[ \begin{array}{cc} 1 & t \\ 0 & 1 \end{array} \right]$ and $s, t \in \mathbb{k}$. Then, by Example~\ref{YBEexamples}, $(X, r)$ is a solution to the Yang--Baxter equation. Let $A=X$ and $V=A \times X$.  Let $\psi : X \times X \rightarrow A$ be defined by 
$$\displaystyle \psi \left( \left[ \begin{array}{cc} x_1 \\ x_2 \end{array} \right] ,
\left[ \begin{array}{cc} y_1 \\ y_2 \end{array} \right] \right)
=  \left[ \begin{array}{cc} w (y_2 - x_2)  \\ 0 \end{array} \right] , $$
where $w \in \mathbb{k}$. Take $s : V \times V \rightarrow  V \times V$ as in Proposition~\ref{extprop}, that is, 
$$s\big( (\vec{a}, \vec{x}), (\vec{b}, \vec{y}) \big) = \big( (\vec{b} + \psi(\vec{x}, \vec{y}) , \sigma_{\vec{x}}( \vec{y}) ), ~(\vec{a }+ \psi(\vec{x}, \vec{y}),  \tau_{\vec{y}}(\vec{x}) ) \big)$$ 
for $(\vec{a}, \vec{x}), (\vec{b},\vec{y})  \in V$.  The map $s$ above can be written as
$$ \left[ \begin{array}{cccc} O & -W & I & W \\
O & I-Y & O & Y \\ I & -W &  O & W \\ O & Z & O & I-Z \end{array} \right] , $$
where $O$ denotes the zero matrix and  $\displaystyle W =  \left[ \begin{array}{cc} 0 & w \\ 0 & 0 \end{array} \right]$. We can rewrite the matrix of $s$ as  $\displaystyle \left[ \begin{array}{cc} I - Y' & Y' \\ Z' & I-Z'  \end{array} \right]$, where $\displaystyle Y'=  \left[ \begin{array}{cc}I & W \\ O & Y  \end{array} \right]$ and  $\displaystyle Z'=  \left[ \begin{array}{cc}I & -W \\ O & Z  \end{array} \right]$. We see that $(I-Y')(I-Z')=O$.
Thus, by Example~\ref{YBEexamples}, $(V,s)$ is a solution to the Yang--Baxter equation and is an extension of $(X,r)$ by  $\psi$. Further, Corollary~\ref{extcor} implies that $\psi \in Z^2_{\rm YB}(X;A)$.
} \end{example}

\begin{proposition} 
Let $\mathbb{k}$ be a ring such that $2$ is invertible in $\mathbb{k}$. Let $X=\mathbb{k}^2$ and $(X,r)$ a solution to the Yang--Baxter equation, where $\displaystyle r= \left[ \begin{array}{cc} I-Y & Y \\ Y & I-Y \end{array} \right]$ and $\displaystyle Y=\left[ \begin{array}{cc} 1 & t \\ 0 & 1 \end{array} \right]$. Then $\Ho^2_{\rm YB} (X;A) \neq 0$. 
\end{proposition}

\begin{proof}
Let $\psi : X \times X \rightarrow A$ be defined by 
$$\displaystyle \psi \left( \left[ \begin{array}{cc} x_1 \\ x_2 \end{array} \right] ,
\left[ \begin{array}{cc} y_1 \\ y_2 \end{array} \right] \right)
=  \left[ \begin{array}{cc} w (y_2 - x_2)  \\ 0 \end{array} \right],$$
where $0 \neq w \in \mathbb{k}$. By Example \ref{extex}, we have $\psi \in Z^2_{\rm YB}(X;A)$. Note that 
$\displaystyle \psi \left( \left[ \begin{array}{cc} 1 \\ 0 \end{array} \right] ,
\left[ \begin{array}{cc} 0 \\ 1 \end{array} \right] \right)
=  \left[ \begin{array}{cc} w   \\ 0 \end{array} \right] $
and 
$\displaystyle \psi \left( \left[ \begin{array}{cc} 0 \\ 1 \end{array} \right] ,
\left[ \begin{array}{cc} 1 \\ 0 \end{array} \right] \right)
=  \left[ \begin{array}{cc} - w   \\ 0 \end{array} \right] $. Since $w \neq 0$ and $2$ is invertible in $\mathbb{k}$, the two vectors are distinct. Suppose that $\psi=\partial^1(f)$ for some  $f \in C^1_{\rm YB}(X;A)$.  Recall that, for all $\vec{u}, \vec{v} \in \mathbb{k}^2 $, we have 
$$\partial^1(f) (\vec{u}, \vec{v})= f(\vec{u})+f( \vec{v})-f \big(\sigma_{\vec{u} } (\vec{v})\big)-f \big(\tau_{\vec{v} } (\vec{u})\big).$$
 Also, notice that $\sigma_{\vec{u} } (\vec{v}) =(I-Y) \vec{u}+Y \vec{v}$ and $\tau_{\vec{v} } (\vec{u}) = (I-Y) \vec{v}+Y \vec{u}$.  Then we see that 
\begin{eqnarray*}
	\displaystyle \psi \left( \left[ \begin{array}{cc} 1 \\ 0 \end{array} \right] ,
	\left[ \begin{array}{cc} 0 \\ 1 \end{array} \right] \right)
	&=&  \partial^1(f)  \left( \left[ \begin{array}{cc} 1 \\ 0 \end{array} \right] ,
	\left[ \begin{array}{cc} 0 \\ 1 \end{array} \right] \right)\\	
	&= &f  \left( \left[ \begin{array}{cc}  1 \\ 0 \end{array} \right] \right)+  f  \left( \left[ \begin{array}{cc} 0 \\ 1 \end{array} \right] \right)	-  f  \left( \left[ \begin{array}{cc} t \\ 1 \end{array} \right] \right)
	- f  \left( \left[ \begin{array}{cc}1 -t 
		\\ 0  \end{array} \right] \right)	\\ 
 &=&  \partial^1(f)  \left( \left[ \begin{array}{cc} 0 \\ 1 \end{array} \right] ,
	\left[ \begin{array}{cc} 1 \\ 0 \end{array} \right] \right)\\
	&=&\displaystyle \psi \left( \left[ \begin{array}{cc} 0 \\ 1 \end{array} \right] ,
	\left[ \begin{array}{cc} 1 \\ 0 \end{array} \right] \right),
\end{eqnarray*}
which is a contradiction.  $\blacksquare$
\end{proof}
\bigskip
\bigskip


\section{Homology of non-degenerate  solutions}
We present a (co)homology theory for left non-degenerate solutions to the Yang--Baxter equation based on the work \cite{MR3558231} of Lebed and Vendramin, which is a generalization of the (co)homology theory given in \cite{MR2128041}.  The (co)homology theory introduced in \cite{MR2128041} is a generalization of the quandle (co)homology theory defined in \cite{MR1990571}, which itself is a generalization of the (co)homology theory of Fenn, Rourke, and Sanderson \cite{MR1257904, MR2063665,MR2255194, MR1364012}.
\bigskip

\subsection{Pre-cubical structure}

Lebed and Vendramin \cite{MR3558231} studied left non-degenerate solutions to the Yang Baxter equation with the focus on two classes of algebraic structures, namely, racks and cycle sets.  They consider two cohomology theories of such solutions which extend previous works and place existing (co)homology theories into a broader perspective.  They showed that for a certain type of left non-degenerate  solutions, including quandles and non-degenerate cycle sets, the (co)homologies split into the degenerate and the normalized parts. Furthermore, they gave an interpretation of the second cohomology group of left non-degenerate solutions in terms of group cohomology generalizing a similar result for racks established by Etingof and Gra\~na \cite{MR1948837}. We now review the constructions from \cite{MR3558231}.

\begin{definition}
Let $\mathcal{C}$ be a category. A \index{pre-cubical structure}{pre-cubical structure} in $\mathcal{C}$ is a sequence of objects $\{C_k\}_{k \geqslant 0}$, together with two sequences of morphisms $d^+_i, d^-_i : C_k \to C_{k-1}$ (called boundaries) for each $k \ge 1$ and $1 \le i \le k$, such that the following  conditions hold:
\begin{eqnarray}\label{Equation pre Cubical}
d^\varepsilon_i d^\zeta_j = d^\zeta_{j-1}d^\varepsilon_i && ~\text{if}~ i < j \:\text{and } \varepsilon,\zeta \in \{+,-\}.
\end{eqnarray} 
A pre-cubical structure is called a \index{weak skew cubical structure}{weak skew cubical structure} if it also includes morphisms $s_i : C_k \to C_{k+1}$ (called degeneracies) for each $k \ge 1$ and $1 \le i \le k$, such that the following conditions hold:
\begin{eqnarray}
\nonumber d^\varepsilon_i s_j = s_{j-1}d^\varepsilon_i & & ~\text{if}~ i < j \:\text{and } \varepsilon \in \{+,-\},\\
\nonumber  d^\varepsilon_i s_j = s_{j}d^\varepsilon_{i-1} & & ~\text{if}~i > j+1 \:\text{and } \varepsilon \in \{+,-\},\\
d^\varepsilon_i s_i = d^\varepsilon_{i+1}s_{i} && ~\text{if}~ i \:\text{and } \varepsilon \in \{+,-\}.\label{Equation Weak Cubical 3}
\end{eqnarray}
A \index{skew cubical structure}{skew cubical structure} (respectively, \index{semi-strong skew cubical structure}{semi-strong skew cubical structure}) additionally satisfies  the conditions
	\begin{eqnarray*}
		s_i s_j = s_{j+1}s_i & &~\text{if}~ i \le j
	\end{eqnarray*}
	and the enhanced version
	\begin{eqnarray}
		d^\varepsilon_i s_i = d^\varepsilon_{i+1}s_{i} = \id \label{Equation Weak Cubical 3'}
		\end{eqnarray}
	of condition~\eqref{Equation Weak Cubical 3} for all $\varepsilon \in \{+,-\}$ (respectively, for $\varepsilon = +$ only).
\end{definition}

With the preceding definition,  we have the following result \cite[Theorem 2.2]{MR3558231}.

\begin{theorem}
Let $\{C_k,\,d^+_i,\, d^-_i \}$ be a pre-cubical structure in the category of modules over a commutative ring $\mathbb{k}$ with unity.
\begin{enumerate}
\item The $\mathbb{k}$-modules $C_k$ endowed with the alternating sum maps 
$$
\partial^{(\alpha, \beta)}_k = \alpha\sum\nolimits_{i=1}^k (-1)^{i-1}d^+_i + \beta \sum\nolimits_{i=1}^k (-1)^{i-1} d^-_i
$$
form a chain complex $\{C_k,\partial^{(\alpha, \beta)}_k \}$ for each $\alpha, \beta \in \mathbb{k}$. 
\item If the boundaries $(d^+_i,\, d^-_i)$ can be completed with degeneracies~$s_i$, then the modules $C^{D}_k= \sum_{i=1}^{k-1} \im (s_i)$ form a sub-complex of $\{C_k,\partial^{(\alpha, \beta)}_k \}$ for each $\alpha, \beta \in \mathbb{k}$.
\item If $\{C_k,\, d^+_i,\, d^-_i,\, s_i \}$ is moreover semi-strong skew cubical, then there are $\mathbb{k}$-module decompositions
\begin{equation} \label{Equation splitting general formula}
C_k = C^{D}_k\oplus C^{N}_k,
\end{equation}
where $C^{N}_k = \im (\eta_k)$ and  $\eta_k = ( \id- s_1d^+_2)( \id- s_2d^+_3)\cdots( \id- s_{k-1}d^+_k).$ It yields a chain complex splitting for $\{C_k,\partial^{(\alpha, \beta)}_k \}$ for each $\alpha, \beta \in \mathbb{k}$.
	\end{enumerate}
\end{theorem}

\begin{proof}
Assertions (1) and (2) can be verified by a direct computation. We give a complete proof of assertion (3). We fix $\alpha, \beta \in \mathbb{k}$ and put $\partial_k = \partial^{(\alpha, \beta)}_k$. We first show that $\eta_k$ form an endomorphism of the chain complex $\{C_k,\partial_k\}$. It suffices to verify the relations 
	\begin{align}\label{Equation eta complex map}
		&\sum\nolimits_{i=1}^k (-1)^{i-1}d^\varepsilon_i \eta_k = \eta_{k-1}\sum\nolimits_{i=1}^k (-1)^{i-1}d^\varepsilon_i \quad \textrm{for} \quad \varepsilon \in \{+,-\}.
	\end{align}
We put $p_i = \id - s_{i}d^+_{i+1}$ and rewrite $\eta_k$ as $p_1 \cdots p_{k-1}$. The weak skew cubical conditions imply that
\begin{eqnarray}
\nonumber d^\varepsilon_i p_j &=& p_{j-1}d^\varepsilon_i \quad \textrm{if} \quad i < j,\\
\label{Equation weak cubical p} d^\varepsilon_i p_j &=& p_{j}d^\varepsilon_{i} \quad \textrm{if} \quad i>j + 1,\\
\label{Equation weak cubical 3p} (d^\varepsilon_i-d^\varepsilon_{i+1})p_{i} &= & d^\varepsilon_i-d^\varepsilon_{i+1}.
\end{eqnarray}
Further, semi-strong skew cubical conditions yield
	\begin{align}
		p_id^\varepsilon_{i+1}p_{i+1} &= d^\varepsilon_{i+1}p_{i+1}.\label{Equation weak cubical 4p}
	\end{align}
Since
	\begin{align*}
		s_{i}d^+_{i+1}d^\varepsilon_{i+1}s_{i+1}d^+_{i+2} &\overset{\eqref{Equation pre Cubical}}{=}s_{i}d^\varepsilon_{i+1}d^+_{i+2}s_{i+1}d^+_{i+2} \overset{\eqref{Equation Weak Cubical 3'}}{=} s_{i}d^\varepsilon_{i+1}d^+_{i+2} \overset{\eqref{Equation pre Cubical}}{=} s_{i}d^+_{i+1}d^\varepsilon_{i+1},
	\end{align*}
we have 
$$
s_{i}d^+_{i+1}d^\varepsilon_{i+1}( \id-s_{i+1}d^+_{i+2}) = 0.
$$
Using \eqref{Equation weak cubical p}, we rewrite the right-hand side of~\eqref{Equation eta complex map} as
$$ p_1 \cdots p_{k-2}\sum\nolimits_{i=1}^k (-1)^{i-1}d^\varepsilon_i = \sum\nolimits_{i=1}^{k-1} (-1)^{i-1} p_1 \cdots p_{i-1} d^\varepsilon_i p_{i+1} \cdots p_{k-1} + (-1)^{k-1}p_1 \cdots p_{k-2} d^\varepsilon_k.
$$
In order to conclude, we obtain the identical expression for the left-hand side of~\eqref{Equation eta complex map} by repeatedly using the following for $s < k-1$:
\begin{eqnarray*}
&& \sum\nolimits_{i=s}^k  (-1)^{i-1} d^\varepsilon_i p_s \cdots p_{k-1} \\
&\overset{\eqref{Equation weak cubical p},\eqref{Equation weak cubical 3p}}{=}&\sum\nolimits_{i=s+2}^k (-1)^{i-1} p_s d^\varepsilon_i p_{s+1} \cdots p_{k-1}  + \big((-1)^{s-1} d^\varepsilon_s + (-1)^{s} d^\varepsilon_{s+1}\big) p_{s+1} \cdots p_{k-1}\\ 
&\overset{\eqref{Equation weak cubical 4p}}{=} & (-1)^{s-1} d^\varepsilon_s p_{s+1} \cdots p_{k-1} + p_s \big(\sum\nolimits_{i=s+1}^k (-1)^{i-1} d^\varepsilon_i \big) p_{s+1} \cdots p_{k-1}.
\end{eqnarray*}
	
Thus, both $C^{D}_k = \sum_{i=1}^{k-1} \im(s_i)$ and $C^{N}_k = \im(\eta_k)$ give sub-complexes of $\{C_k,\partial^{(\alpha, \beta)}_k\}$. It remains to establish the $\mathbb{k}$-module decomposition $C_k= C^{D}_k\oplus C^{N}_k$. The definition of $\eta_k$ directly imply that $\im ( \id- \eta_k) \subseteq \sum_{i=1}^{k-1} \im(s_i)$. 
Further, the semi-strong skew cubical conditions imply that
\begin{eqnarray*}
s_i p_j &=& p_{j+1}s_i \quad \textrm{if} \quad i \le j,\\
s_i p_j &=& p_{j}s_{i} \quad \textrm{if} \quad i > j+1,\\
p_is_i &=& 0,\\
p_is_{i+1} &=& s_{i+1} - s_i,\\
(s_i-s_{i+1})p_{i} &=& s_i-s_{i+1}.
\end{eqnarray*}
These relations further imply that the map~$\eta_k$ vanishes on $\im(s_i)$ for all $1 \le i \le k-1$. Consider the sequence $$C^{D}_k \overset{\iota_k}{\hookrightarrow} C_k \overset{\eta_k}{\to} C_k$$ of $\mathbb{k}$-modules and homomorphisms, where~$\iota_k$ is the inclusion map. Note that $\eta_k \iota_k = 0$ and $\im ( \id- \eta_k) \subseteq \im(\iota_k)$. The condition $\im ( \id- \eta_k) \subseteq \im(\iota_k)$ implies that $C_k=\im(\iota_k) + \im (\eta_k)$. The relation $\eta_k \iota_k = 0$  means that $\im (\iota_k) \subseteq \ker(\eta_k)$. Thus, we have $\im(\id -\eta_k) \subseteq  \ker(\eta_k)$, which translates to $\eta_k^2 = \eta_k$. But this yields $\ker (\eta_k) \cap \im (\eta_k)=0$, and hence $\im (\iota_k) \cap \im (\eta_k)=0$. Thus, we obtain the direct sum decomposition $C_k=\im (\iota_k) \oplus \im(\eta_k)$. $\blacksquare$
\end{proof}

The decomposition~\eqref{Equation splitting general formula} induces a decomposition at the level of homology groups.  The ${D}$- and ${N}$- parts of the complexes and homology groups above are called \emph{degenerate} and \emph{normalized}, respectively.
\bigskip
\bigskip


\subsection{Braided homology}\label{S:BrHom}  
As an application of the notion of a pre-cubical structure, we present an extension of the (co)homology theory for solutions developed in \cite{MR3081627}, while also generalising the (co)homology theory initiated in \cite{MR2128041}.
\para 

\begin{figure}[!ht]
 \begin{center}
\includegraphics[height=3.3cm]{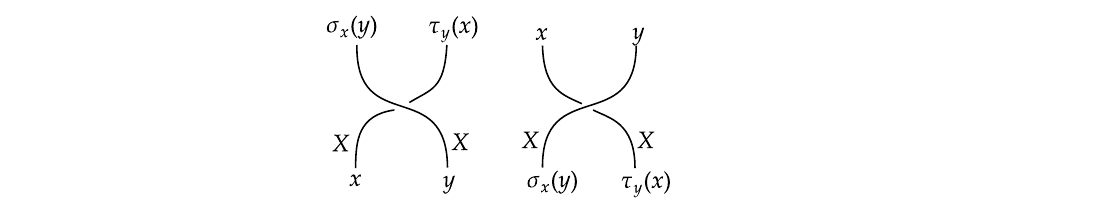}
\caption{Labelling by a solution at a crossing and its opposite.}
\label{graphic calculus rule} 
\end{center}
\end{figure}

\begin{figure}[!ht]
 \begin{center}
\includegraphics[height=6cm]{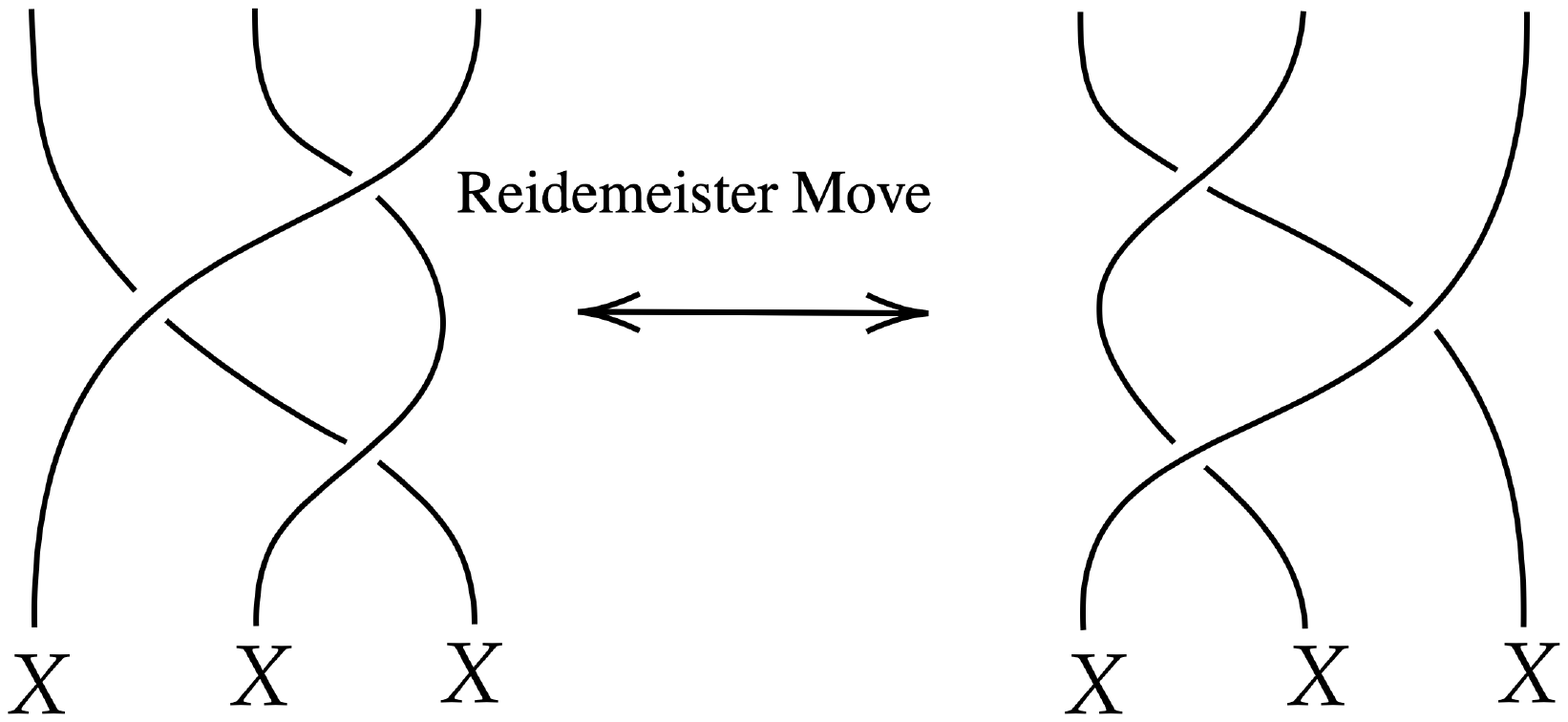}
\caption{Reidemeister move $R3$ representing the Yang--Baxter equation.}
\label{graphical YBE} 
\end{center}
\end{figure}

Before proceeding further, we recall some graphical calculus which has been used by many researchers and which will make the constructions intuitive \cite{MR3558231}. A braided diagram represents maps between sets, a set being associated to each strand. The horizontal glueing corresponds to Cartesian product, vertical glueing to composition (which should be read from bottom to top), straight vertical lines to identity maps, crossings to the braiding $r$, and opposite crossings to its inverse (whenever it exists), as shown in Figure \ref{graphic calculus rule}. With these conventions, the Yang--Baxter equation for the braiding $r$ becomes the diagram in Figure \ref{graphical YBE}, which is precisely the Reidemeister move $R3$ from knot theory.

\begin{figure}[!ht]
 \begin{center}
\includegraphics[height=3cm]{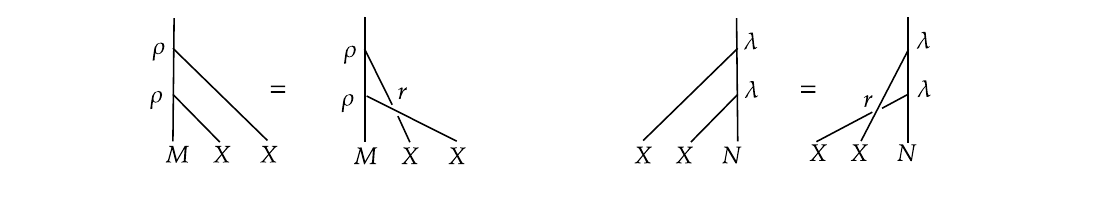}
\caption{Right and left braided modules.}
\label{left right modules YBE} 
\end{center}
\end{figure}

\begin{definition}\label{braided module over a solution}
Let $(X, r)$ be a solution to the Yang--Baxter equation. A right (braided) module over $(X,r)$ is a pair $(M, \rho)$, where $M$ is a set and $\rho : M \times X \to M$ is a map  written as $\rho(m,x) = m \cdot x$ such that the following  condition
$$(m \cdot x) \cdot y = \big(m \cdot  \sigma_x(y) \big) \cdot \tau_y(x)$$ 
 hold for all $m \in M$ and $x,y \in X$. \\
 Similarly, a left (braided) module over $(X,r)$ is a pair $(N, \lambda)$, where  $\lambda : X \times N \to N$ is a map  written as $\lambda(x, n) = x \cdot n$ such that 
$$ x \cdot (y \cdot n)  = \sigma_x(y) \cdot \big(\tau_y(x) \cdot n \big)$$ 
 hold for all $n \in N$ and $x,y \in X$. See Figure \ref{left right modules YBE} for a graphical interpretation of braided modules.
\end{definition}

\begin{example} \label{Example Trivial}
{\rm Let us consider some examples.
\begin{enumerate}
\item  The singleton set $I = \{\ast\}$ with the unique map $X \to I$ yields an example of a right (as well as left) $(X,r)$-module for any solution $(X,r)$. It is referred to as the \textit{trivial} right (left) $(X,r)$-module.
\item If $(X,r)$ is a solution, then $X$ is a right module over itself with the action $\rho(x, y) =\tau_y(x)$  (see \eqref{condition for solution 3}). Similarly, $X$ it is a left module over itself with the action $\lambda(x, y) = \sigma_x(y)$ (see \eqref{condition for solution 1}). These modules are called \index{adjoint module} adjoint modules. More generally, any power $X^n$ is a right and a left module over $T(X) = \sqcup_{i \ge 0} X^i$ with the module structure adjoint to the extension of the braiding $(X,r)$ to $T(X)$.
\end{enumerate}}
\end{example}

Let $(M,\rho)$ be a right module  and $(N,\lambda)$ a left module  over a solution $(X,r)$ to the Yang--Baxter equation. Let  $r_i, \rho_0, \lambda_n: M \times X^n \times N \to M \times X^{(n-1)}\times N$ be maps given by
\begin{eqnarray*}
r_i &=& \id_{M} \times \id_{X}^{(i-1)} \times r \times \id_{X}^{ (n-i-1)} \times \id_N,\\
\rho_0 &=& \rho \times \id_{X}^{(n-1)} \times \id_N,\\
\qquad \lambda_n &=& \id_{M} \times \id_{X}^{ (n-1)} \times \lambda.
\end{eqnarray*}

The  proof of the following result is a computational check and left as an exercise \cite[Theorem 3.5]{MR3558231}.

\begin{figure}[!ht]
 \begin{center}
\includegraphics[height=10cm]{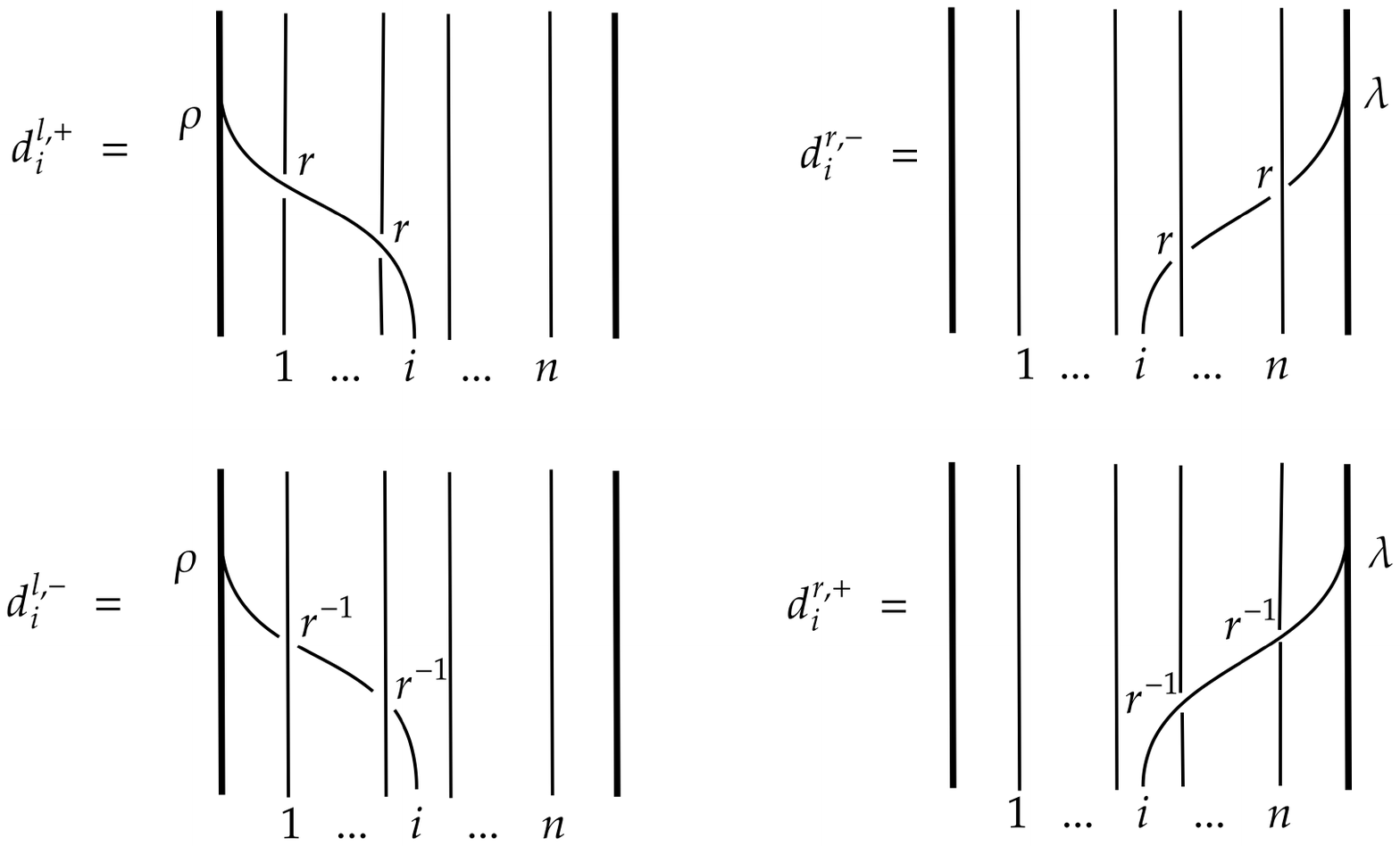}
\caption{Differentials in braided homology.}
\label{differentials of Theorem BrHom} 
\end{center}
\end{figure}

\begin{theorem}\label{Theorem module sol pre cubical}
Let $(M,\rho)$ be a right module  and $(N,\lambda)$ a left module  over a solution $(X,r)$ to the Yang--Baxter equation. Define $C_n = M \times X^n \times N$. Then the following assertions hold:
\begin{enumerate} 
\item  The  maps $d_i^{l,+}, d_i^{r,-}: C_n \to C_{n-1}$ given by 
$$
d_i^{l,+} = \rho_0 \circ r_1 \circ \cdots \circ r_{i-1}
\quad \textrm{and} \quad 
 d_i^{r,-} = \lambda_n \circ r_{n-1} \circ \cdots \circ r_{i}
$$
give a pre-cubical structure $\{C_n, d_i^{l,+}, d_i^{r,-} \}$.
\item Suppose that the braiding $r$ is invertible and define
$$ d_i^{l,-} = \rho_0 \circ r^{-1}_1 \circ \cdots \circ r^{-1}_{i-1}
\quad \textrm{and} \quad 
d_i^{r,+} = \lambda_n \circ r^{-1}_{n-1} \circ \cdots \circ r^{-1}_{i}.
$$
Then for any choice of $\varepsilon,\zeta \in \{l,r\}$, the triplet $\{C_n, d_i^{\varepsilon,+},\, d_i^{\zeta,-} \}$ with $n \ge 1$ and $1 \le i \le n$, forms a pre-cubical structure. See Figure \ref{differentials of Theorem BrHom} for a graphical interpretation of the differentials.
\end{enumerate}
\end{theorem}

\begin{remark}
{\rm 
The homology theory arising from Theorem \ref{Theorem module sol pre cubical} is called the \index{braided homology}{braided homology}. By choosing trivial modules (see Example \ref{Example Trivial}) as coefficients with  $\alpha = 1$ and $\beta = -1$, one gets the chain complex and hence homology of \cite{MR2128041} described in Section \ref{homsec colorings} of this chapter. See also \cite[Theorem 3.1]{MR3835755} for a direct chain  homotopy equivalence between the two complexes. We recommend the reader to see \cite{MR3558231} for elegant graphical and intuitive interpretations of the differentials defined in this section.}
\end{remark}
\bigskip


\subsection{Birack homology}\label{S:BirackHom} 
The origins of a homology theory for biracks lie in the work \cite{MR1257904} of Fenn, Rourke and Sanderson, and was developed in the language of augmented biracks in \cite{MR3266530}. It is worth noting that, in \cite{MR2906433, MR3381331}, Przytycki studied homology theory of distributive magmas which specialises to homology of racks and quandles.

\begin{definition}
 A solution to the Yang--Baxter equation is called a \index{birack}{birack} if it is bijective and non-degenerate.
\end{definition}

In other words, a solution $(X, r)$ with $r(x, y)= \big(\sigma_x(y), \tau_y(x)\big)$ is a birack if the braiding $r:X \times X \to X \times X$ and the maps $\sigma_x, \tau_x:X \to X$ are bijective  for all $x \in X$.

\begin{remark}
{\rm It is worth noting that in several other works, for instance \cite{MR4406425, MR2100870, MR4125863}, the term birack is used for an algebraic structure with two binary operations satisfying certain mixing axioms. It turns out that the birack defined above gives rise to such an algebraic structure and the converse also holds. Special types of biracks gives rise to what are called biquandles, which have been used extensively for constructing invariants of virtual knots and links \cite{MR3055555, MR3666513, MR2100870, MR2994593, MR3868206}.
}
\end{remark}

We have the following result \cite[Theorem 4.2]{MR3558231}.

\begin{theorem}\label{Theorem birack pre cubical}
Let $(X,r)$ be a birack and $C_n=X^n$. Then the following assertions hold:
\begin{enumerate}
\item The maps $d_i, d_i':C_n\to C_{n-1}$ given by
\begin{eqnarray*}
&& d_i : (x_1, \ldots, x_n) \mapsto \big(\tau_{x_i}^{-1}(x_1), \ldots, \tau_{x_i}^{-1}(x_{i-1}), \sigma_{\tau_{x_{i+1}}^{-1}(x_i)}(x_{i+1}), \ldots, \sigma_{\tau_{x_n}^{-1}(x_i)}(x_n)\big),\\
&& d_i'  : (x_1, \ldots, x_n) \mapsto (x_1, \ldots, x_{i-1}, x_{i+1}, \ldots, x_n), 
\end{eqnarray*}
form a pre-cubical structure $\{C_n, d_i, d_i' \}$.
\item The assertion remains true when the maps $d_i$ are replaced with
\begin{eqnarray*}
&& d^\star_i : (x_1, \ldots, x_n) \mapsto \big( \sigma_{\tau_{x_1}^{-1}(x_i)}(x_1),\ldots, \sigma_{\tau_{x_{i-1}}^{-1}(x_i)}(x_{i-1}), \tau_{x_i}^{-1}(x_{i+1}), \ldots, \tau_{x_i}^{-1}(x_n)\big).
\end{eqnarray*}
\end{enumerate}
\end{theorem}

\begin{proof}
The proof of assertion (1) follows from the fact that $d_i\; d'_j=d'_{j-1} \; d_i$, which is straightforward from the definitions of $d_i$ and $d'_i$. Similarly, definitions of $d_i^\star$ and $d'_i$ yield assertion (2). $\blacksquare$
\end{proof}

The homology theory arising from Theorem \ref{Theorem birack pre cubical} is called the \index{birack homology}{birack homology}. It turns out that the braided homology and the birack homology are intimately related. Let us set some terminology to state the precise result. 

\begin{figure}[!ht]
 \begin{center}
\includegraphics[height=5.6cm]{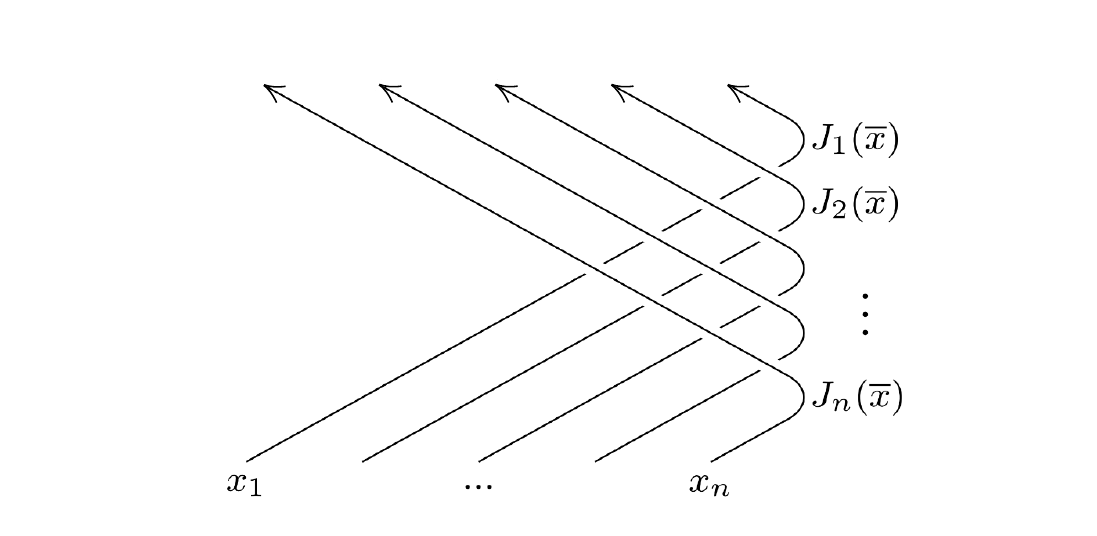}
\caption{The $n$-guitar map.}
\label{guitar map} 
\end{center}
\end{figure}

\begin{definition}
Let  $(X, r)$ be a solution to the Yang--Baxter equation and $n \ge 1$. The $n$-guitar map for $(X, r)$ is the map $J: X^n\to  X^n$ defined by $$J(\bar{x})= \big(J_1(\bar{x}), \ldots, J_n(\bar{x}) \big),$$ where $J_i(\bar{x})=\tau_{x_n} \tau_{x_{n-1}}\cdots \tau_{x_{i+1}}(x_i)$ and $\bar{x}=(x_1, \ldots, x_n)\in X^n$.  See Figure \ref{guitar map} for a graphical interpretation of the $n$-guitar map.
\end{definition}

Recall that, a shelf is a non-empty set with a right-distributive binary operation.

\begin{lemma}\label{LNDS guita map lemma}
Let $(X, r)$ be a left non-degenerate solution to the Yang--Baxter equation. Then the following assertions  hold:
\begin{enumerate}
\item The binary operation $x*_r y= \tau_y \big(\sigma_{\tau_x^{-1}(y)}(x)\big)$ defines a shelf structure on $X$.
\item The shelf $(X, *_r)$ is a rack if and only if $r$ is bijective.
\item The $n$-guitar map $J$ for $(X, r)$ is bijective for all $n \ge 1$.
\end{enumerate}
\end{lemma}

\begin{proof}
Assertion (1) is immediate from the left non-degeneracy of the solution $(X, r)$. Similarly, the bijectivity of $J$ is equivalent to the left non-degeneracy of $r$, which gives assertion (3). For assertion (2), suppose that $r$ is bijective. This implies that, for a given $y \in X$, the map $x \mapsto \sigma_{\tau_x^{-1}(y)}(x)$ is bijective. Since the map $x \mapsto \tau_y(x)$ is bijective by left non-degeneracy, it follows that the map $x \mapsto x *_r y$ is bijective, and hence $(X, *_r)$ is a rack. Conversely, the bijectivity of $x \mapsto x *_r y$ and of $x \mapsto \tau_y(x)$ yields the bijectivity of $x \mapsto \sigma_{\tau_x^{-1}(y)}(x)$. Then $r^{-1}$ is defined by $r^{-1}(z, y) = \big(\tau_x^{-1}(y), x \big)$, where $x$ is the unique element of $X$ satisfying  $\sigma_{\tau_x^{-1}(y)}(x)=z$. $\blacksquare$
\end{proof}

Let $(X, r)$ be a left non-degenerate solution. In view of Lemma \ref{LNDS guita map lemma}(3), let $J^{-1}= (J_1^{-1}, \ldots, J_n^{-1})$ denote  the inverse of the $n$-guitar map $J$. Define maps $\chi_i: X^n \to X$ by
$$\chi_i(x_1, \ldots, x_n)=J_1^{-1}\big((x_i *_r x_{i-1})*_r \cdots *_r x_1, x_1, x_2, \ldots, x_{i-1}, x_{i+1}, \ldots, x_n \big),$$

The next result extends the birack homology to the more general left non-degenerate solutions and also add coefficients to the involved chain complexes. Further, it establishes an equivalence between this extended birack homology and the braided homology using the guitar map \cite[Theorem 7.1]{MR3558231}. 

\begin{theorem}\label{equality of braided and birack homology}
Let $(X, r)$ be a left non-degenerate solution to the Yang--Baxter equation. Let $(M,\rho)$ be a right module  and $(N,\lambda)$ a left module  over $(X,r)$. Consider the sets $C_n = M \times X^n \times N$. Then the following assertions hold:
\begin{enumerate} 
\item  The  maps $d_i, d_i': C_n \to C_{n-1}$ defined by 
\begin{eqnarray*}
&& d_i: (a, x_1, \ldots, x_n, b) \mapsto  \big(a, \tau_{x_i}^{-1}(x_1), \ldots, \tau_{x_i}^{-1}(x_{i-1}), \sigma_{\tau_{x_{i+1}}^{-1}(x_i)}(x_{i+1}), \ldots, \sigma_{\tau_{x_n}^{-1}(x_i)}(x_n),  \lambda(x_i, b) \big),\\
&& d_i'  : (x_1, \ldots, x_n) \mapsto \big(\rho(a, \chi_i(\bar{x})),  x_1, \ldots, x_{i-1}, x_{i+1}, \ldots, x_n, b \big) 
\end{eqnarray*}
for $\bar{x}=(x_1, \ldots, x_n) \in X^n$, $a \in M$ and $b \in N$, give a pre-cubical structure $\{C_n, d_i, d_i' \}$.
\item The extended guitar map $J := \id_M \times J \times  \id_N$ yields an isomorphism between the pre-cubical structure $\{C_n, d_i^{r,-}, d_i^{l,+}\}$ of Theorem  \ref{Theorem module sol pre cubical} and the structure $\{C_n, d_i, d_i' \}$ of assertion (1) above.
\end{enumerate}
\end{theorem}

\begin{remark}
{\rm Note that when the solution is a birack and the coefficients $M$ an $N$ are trivial (see Example \ref{Example Trivial}), the pre-cubical structure from Theorem \ref{equality of braided and birack homology} specializes to that from Theorem \ref{Theorem birack pre cubical}. Thus, the braided homology and the birack homology are identical for biracks.
}
\end{remark}

\begin{remark}\label{star 2-cocycle}
{\rm
Let $(X, r)$ be a  left non-degenerate solution to the Yang--Baxter equation. Changing the pre-cubical structure $(d_i, d'_i)$ to $(d^{\star}_i, d'_i)$ in Theorem~\ref{Theorem birack pre cubical}, one obtains  a (co)homology theory for which a map $f: X \times X \to A$ is a $2$-cocycle if 
\begin{eqnarray}\label{left non-dgen cocycle condition}
  f(x,z)+f \big(\tau_x^{-1}(y),\tau_x^{-1}(z) \big)=f(y,z)+f \big(\sigma_{\tau_x^{-1}(y) } (x),\tau_y^{-1}(z) \big)
\end{eqnarray}
for all $x, y, z \in X$.}
\end{remark}
\medskip

\begin{lemma}	\label{L:FunAsMod}
Let $(X,r)$ be a left non-degenerate solution to the Yang--Baxter equation and~$A$ an abelian group. Then the map
	$$\Map(X,A) \times X \to \Map(X,A)$$
given by $(\gamma, x) \mapsto \gamma \cdot x$, where  $(\gamma \cdot x)(y)= \gamma \big(\tau_x^{-1} (y)\big)$ for all $x, y \in X$,
extends to a right $G(X,r)$-module structure on the abelian group $\Map(X,A)$ of maps from $X$ to $A$.
\end{lemma}

\begin{proof}
Since $(X, r)$ is left non-degenerate, each $\tau_x$ is a bijection. Hence, by Proposition \ref{p122}, there is a left-action of the group $G(X,r)$ on the set $X$ obtained by 
extending the map $X \to \Sigma_X$, given by $x \mapsto  \tau_x^{-1}$, to a group homomorphism $G(X,r) \to \Sigma_X$. This left-action induces a right $G(X, r)$-module structure on the abelian group $\Map(X,A)$ by setting $(\gamma, u) \mapsto \gamma \cdot u$, where  
$$(\gamma \cdot u)(x)= \gamma(u \cdot x)=\gamma \big(\tau_u^{-1} (x) \big) \quad \textrm{and} \quad \tau_u^{-1}= \tau_{x_k}^{-\epsilon_k} \cdots \tau_{x_1}^{-\epsilon_1}
$$ for $x \in X$, $\gamma \in \Map (X, A)$ and $u= x_1^{\epsilon_1}\cdots x_k^{\epsilon_k} \in G(X, r)$. $\blacksquare$
\end{proof}

Compare the action of Lemma \ref{L:FunAsMod} with the one given by \eqref{adj action on set of maps} for racks. The group cohomology with these choices turns out to be useful in studying the cohomology of our left non-degenerate solution. We now present the following result \cite[Theorem 11.2]{MR3558231}, which generalizes Theorem \ref{relation rack cohomo and group cohomo}(2). 

\begin{theorem}\label{T:GroupCohom}
Let $(X,\, r)$ be a left non-degenerate solution to the Yang--Baxter equation and~$A$ an abelian group. Then there is an isomorphism
$$\Ho^2_{\rm LND}(X;A)\cong \Ho^1_{\rm Grp} \big(G(X, r);\Map(X,A) \big)$$
of groups, where on the left hand side $\Ho^2_{\rm LND}$ stands for the cohomology of $(X,\, r)$ with coefficients in $A$ (as in Remark \ref{star 2-cocycle}), and on the right hand side $\Ho^1_{\rm Grp}$ stands for the group cohomology of $G(X, r)$ with coefficients in the right module $\Map(X,A)$ (as in Lemma \ref{L:FunAsMod}).
\end{theorem}
	
\begin{proof}
Consider the map 
$$ \Phi: Z^1\big(G(X, r);\Map(X,A) \big) \longrightarrow Z^2(X;A)$$
given by
$$\Phi(\theta)(x,y)=\theta(x)(y)$$
for $\theta \in Z^1\big(G(X, r);\Map(X,A) \big)$ and $x, y \in X$.
Consider another map
$$\Psi: Z^2(X;A) \to  Z^1\big(G(X, r); \Map(X,A) \big)$$
given by
$$f \mapsto \Psi(f),$$
where $\Psi(f):G(X, r) \to \Map(X,A)$ has the form
$$\Psi(f)(u)(x)= f(u,x)$$
for $f \in Z^2(X;A)$, $u \in G(X, r)$ and $x \in X$.  We claim that the following assertions hold:
\begin{enumerate}
\item $\Phi $ is well-defined.
\item $\Psi$ is well-defined.
\item $\Phi$ maps $1$-coboundaries to $2$-coboundaries.
\item $\Psi$ maps $2$-coboundaries to $1$-coboundaries.
\end{enumerate}

In view of the preceding claims, we see that $\Phi$ and~$\Psi$ are inverses of each other.  Let us prove the claims one-by-one.
\begin{enumerate}
\item Let $\theta \in Z^1\big(G(X, r); \Map(X,A) \big)$ be a group 1-cocycle. Then, by \eqref{group coboundary 1}, $\theta$ satisfies the relation 
			$$\theta(g_1g_2) =  \theta(g_2) + \theta(g_1)\cdot g_2 = \theta(g_2) + \theta(g_1)\tau_{g_2}^{-1}$$
			for $g_1,g_2 \in G(X, r)$.	For $f := \Phi(\theta)$, we need to check the 2-cocycle condition \eqref{left non-dgen cocycle condition} given by
			$$f(x,z)+ f \big(\tau_x^{-1}(y),\tau_x^{-1}(z) \big)=f(y,z)+f \big(\sigma_{\tau_x^{-1}(y) } (x),\tau_y^{-1}(z) \big)$$
			for all $x,y,z \in X$, which can be rewritten as
			$$\theta(x)(z)+\theta \big(\tau_x^{-1}(y) \big) \big(\tau_x^{-1}(z) \big)=\theta(y)(z)+ \theta \big(\sigma_{\tau_x^{-1}(y)}(x) \big) \big(\tau_y^{-1}(z) \big).$$
			Recalling the definition of the $G(X, r)$-action on $\Map(X,A)$, the preceding equation transforms into
			$$\theta(x)(z)+\big(\theta(\tau_x^{-1}(y)) \cdot x \big)(z)=\theta(y)(z)+ \big(\theta(\sigma_{\tau_x^{-1}(y)}(x)) \cdot y \big)(z).$$
   This gives
   $$\theta(x)+\theta \big(\tau_x^{-1}(y) \big) \cdot x=\theta(y)+\theta \big(\sigma_{\tau_x^{-1}(y)}(x) \big) \cdot y$$
   for all $x, y \in X$. By the $1$-cocycle condition for~$\theta$, the preceding equation is the same as
$$\theta \big((\tau_x^{-1}(y)) x \big)= \theta \big(( \sigma_{\tau_x^{-1}(y) } (x)) y \big).$$
But this relation follows from the relation $\big(\tau_x^{-1}(y) \big) x = \big( \sigma_{\tau_x^{-1}(y) } (x)\big) y$, which always hold in the structure group $G(X, r)$.
			\item Recall that, a map $\theta: G(X, r)\rightarrow  \Map(X,A)$ is a group $1$-cocycle if and only if the map
			$G(X, r)\to G(X, r)\ltimes \Map(X,A)$ given by $g\mapsto (g,\theta(g))$ is a group homomorphism. Here, $G(X, r)\ltimes \Map(X,A)$ is the set  $G(X, r)\times \Map(X,A)$ endowed with the group operation 
   $$(g,\gamma)(g',\gamma')=(gg',~\gamma \cdot g' +\gamma')=(gg',~\gamma \tau^{-1}_{g'} +\gamma')$$ 
in the semi-direct product.   We will now show that, for each $f \in Z^2(X;A)$, the map $\hat{f} : X \to G(X, r)\ltimes  \Map(X,A)$, given by $\hat{f}(x)=\big(x,f(x,-)\big)$, extends to a unique group homomorphism $G(X, r) \to G(X, r)\ltimes \Map(X,A)$. For this, it suffices to check that
$$\hat{f} \big(\tau_x^{-1}(y)\big) \hat{f}(x) = \hat{f} \big( \sigma_{\tau_x^{-1}(y) } (x) \big) \hat{f}(y)$$
for all $x,y \in X$. The preceding equation reads as
			$$\big(\tau_x^{-1}(y),~ f(\tau_x^{-1}(y) ,-)\big)\big(x,f(x,-)\big)=\big( \sigma_{\tau_x^{-1}(y) } (x), ~f( \sigma_{\tau_x^{-1}(y) } (x) ,-)\big)\big(y,f(y,-)\big),$$
			which further simplifies to	
$$ \big((\tau_x^{-1}(y))x, ~f(\tau_x^{-1}(y), 	\tau_x^{-1}(-)) + f(x,-)\big)= \big(( \sigma_{\tau_x^{-1}(y) } (x))y,  ~f(\sigma_{\tau_x^{-1}(y) } (x) , 	\tau_{y}^{-1}(-) )+f(y,-)\big).$$
Since the relation $ \big(\tau_x^{-1}(y) \big) x = \big(\sigma_{\tau_x^{-1}(y)}(x) \big) y$ always holds in $G(X, r)$, it remains to show that
			$$ f \big(\tau_x^{-1}(y), 	\tau_x^{-1}(-) \big) + f(x,-)= f \big(\sigma_{\tau_x^{-1}(y) } (x) , 	\tau_{y}^{-1}(-) \big)+f(y,-),$$
			but this is precisely the definition of a quandle $2$-cocycle.
\item A group $1$-coboundary in $C^1\big(G(X, r); \Map(X,A) \big)$ is a map $\lambda: G(X, r) \to \Map(X,A)$ of the form 			
   $$\lambda(g)= \gamma - \gamma \cdot g= \gamma - \sigma_{\tau_{g}^{-1}(\gamma) } (g)$$
    for some $\gamma \in \Map(X,A).$  Its image $\Phi(\lambda)$ is the map given by
			$$\Phi(\lambda)(x, y)= \lambda(x)(y)= (\gamma - \gamma \cdot x)(y)= \gamma(y) - \gamma \big( \tau_{x}^{-1}(y)\big) = - \partial^1 \gamma(x,y),$$
			yielding $\Phi(\lambda) = \partial^1 (-\gamma)$.
\item A $2$-coboundary in $C^2(X;A)$ is a map of the form $\partial^1(\gamma)$ for some $\gamma \in \Map(X,A)$. Since~$\Phi$ and~$\Psi$ are  inverses of each other, the computation above implies the relation $\Psi \big(\partial^1 (\gamma)\big) = \partial_{{\rm group}}^0 (-\gamma)$.
		\end{enumerate}
  This completes the proof of the theorem. $\blacksquare$
\end{proof}
\bigskip
\bigskip


\section{Knot invariants  from generalized quandle $2$-cocycles} \label{cocyinvsec}

Following \cite{MR2166720}, we give some applications of generalized 2-cocycles with coefficients in quandle modules to define quandle cocycle invariants for classical knots and links.
\para 

\subsection{Boltzmann weights and quandle cocycle invariants} Let $X$ be a quandle and $\mathbb{Z}(X)$ its quandle algebra with generators 
$$ \big\{ \eta_{x,y}^{\pm 1}  \, \mid \,  x, y \in X \big\} \quad \textrm{and} \quad \big\{ \xi_{x,y}  \, \mid \,  x, y \in X \big\}.$$
 Then there is a surjective $\mathbb{Z}$-algebra homomorphism from $\mathbb{Z}(X)$ onto the group algebra $\mathbb{Z}[\Adj(X)]$ given by  $\eta_{x,y } \mapsto \textswab{a}_y$ and $\xi_{x,y} \mapsto 1-\textswab{a}_{x*y}$. Thus, any left $\mathbb{Z}[\Adj(X)]$-module has the structure of a $\mathbb{Z}(X)$-module. For the rest of this subsection, we assume that $A$ is an abelian group with  this  $\mathbb Z (X)$-module structure. Let $\kappa_{x,y}$ be a generalized quandle $2$-cocycle of $X$ with coefficients in the abelian group $A$. By \eqref{2-cocycle condition for gen cohomo}, in this setting, the generalized $2$-cocycle condition is written as 
\begin{equation}\label{gen quandle cocycle condition for application}
\textswab{a}_z \kappa_{x,y} + \kappa_{x*y, z} + \textswab{a}_{(x*y)*z} \kappa_{y,z}=  \kappa_{y,z} + \textswab{a}_{y*z}  \kappa_{x, z} + \kappa_{x*z,y*z}.    
\end{equation}
We will define a cocycle invariant for knots and links using this $2$-cocycle. 
\para 

Let $D(K)$  be a given oriented classical knot or link diagram. A vector perpendicular to an arc in the diagram is called a normal vector. We choose the normal so that the ordered pair (tangent, normal) agrees with the orientation of the plane. This normal is called the {\it orientation normal}. There are four regions near each crossing, divided by the arcs of the diagram. The unique  region into which  both normals (to over- and under-arcs) point is called the {\it target region}.   Let $\gamma$ be an arc from the region at infinity of the plane to the target region of a given crossing $r$, that intersect  $D(K)$ in  finitely many points transversely and missing the crossing points. Let $a_1, a_2, \ldots, a_k$ be the ordered set of arcs of $D(K)$ that intersect $\gamma$ as we move from the region at infinity to the crossing $r$. Let ${\mathcal C}$ be a coloring of $D(K)$ by a fixed finite quandle $X$.
\para

We have the following definition \cite[Definition 6.1]{MR2166720}.

\begin{definition}
The \index{Boltzmann weight}Boltzmann weight $B( {\mathcal C}, r, \gamma)$ for the crossing $r$, for a coloring ${\mathcal C}$, with respect to $\gamma$, is defined by 
		$$ B( {\mathcal C},r, \gamma)=
		\pm  \big(\textswab{a}_{{\mathcal C}(a_1)}^{\epsilon (a_1) } \textswab{a}_{{\mathcal C}( a_2)}^{\epsilon (a_2)}
		\cdots \textswab{a}_{{\mathcal C}( a_k)}^{\epsilon (a_k) } \big)  \, \kappa_{x, y} \in A, $$
where $x,y$ are the colors at the crossing $r$ ($x$ is assigned on the under-arc from which the normal of the over-arc
points, and $y$ is assigned to the over-arc), and the element 
		$$ \textswab{a}_{{\mathcal C}(a_1)}^{\epsilon (a_1) } \textswab{a}_{{\mathcal C}( a_2)}^{\epsilon (a_2) }
		\cdots \textswab{a}_{{\mathcal C}( a_k)}^{\epsilon (a_k) }  \in \Adj(X)$$ 
acts on $A$ via  the quandle module structure. The sign $\pm$ in front is determined by whether $r$ is a positive crossing or a negative crossing. For each $j$, the exponent $\epsilon (a_j)$ is $1$ if the arc $\gamma$ crosses the arc $a_j$ against its normal, and is $-1$ otherwise.
\end{definition}

\begin{figure}[!ht]
 \begin{center}
\includegraphics[height=6cm]{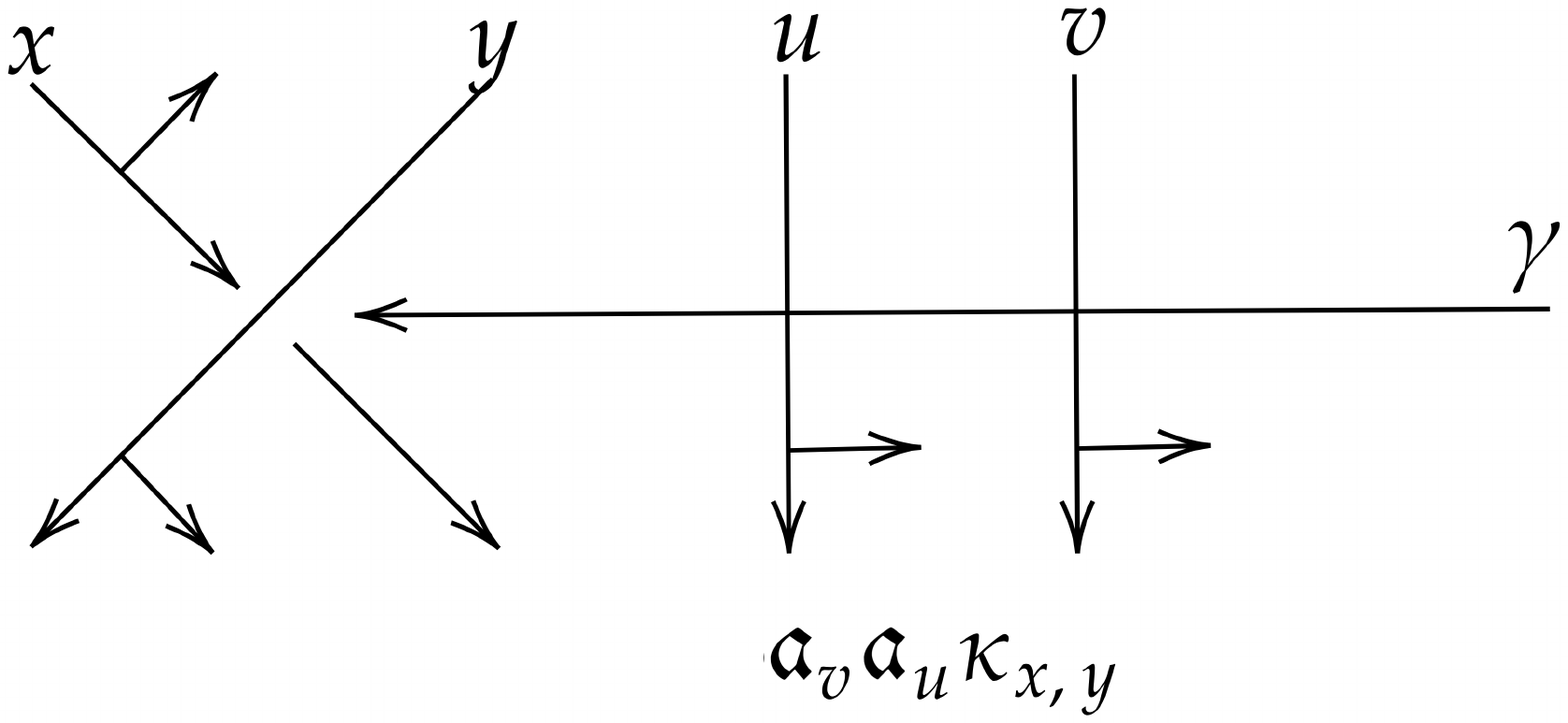}
\caption{Weights on arcs.}
\label{weightR}
\end{center}
\end{figure}

The set-up is depicted in Figure \ref{weightR},  when the arc $\gamma$ intersects two arcs with colors $u$ and $v$ in this order.

\begin{lemma} 
The Boltzmann weight $B( {\mathcal C}, r, \gamma)$ does not depend on the choice of the arc $\gamma$.
\end{lemma} 

\begin{proof}
It is sufficient to check the changes in the coefficient $ \textswab{a}_{{\mathcal C}(a_1)}^{\epsilon (a_1) } \textswab{a}_{{\mathcal C}( a_2)}^{\epsilon (a_2) }	\cdots \textswab{a}_{{\mathcal C}( a_k)}^{\epsilon (a_k)}  \in \Adj(X)$ when the arc is homotoped. When a pair of intersection points with the knot diagram $D(K)$ is canceled or introduced by a path that  zig-zags, then their colors are inverses of each other, and are adjacent in the above sequence, so that the product does not change. 
It remains  to check what happens when an arc is homotoped through 	each crossing, and the effect of the incident is depicted in Figure \ref{pathcross}. 

\begin{figure}[!ht]
 \begin{center}
\includegraphics[height=8cm]{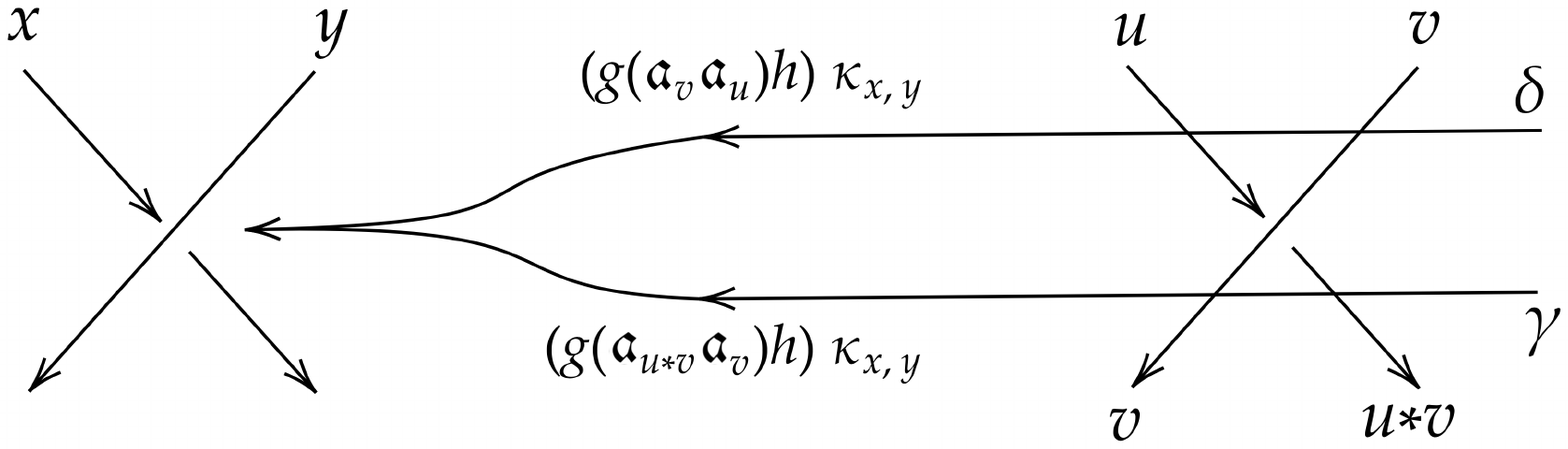}
\caption{When an arc passes a crossing.}
\label{pathcross} 
\end{center}
\end{figure}

Consider the two arcs $\delta$ and $\gamma$ as depicted in  Figure \ref{pathcross}. The given crossing $r$ is the left crossing with colors $x$ and $y$.  The arc $\gamma$ is obtained from $\delta$ by homotoping $\delta$ through a crossing with colors $u$ and $v$ as depicted. One sees that $B( {\mathcal C},r, \delta)= \big(g (\textswab{a}_v \textswab{a}_u)h \big) \kappa_{x,y}$,	where  $g$ and $h$ are sequences of colors  that $\delta$ intersects before and after, respectively, it  intersects  $u$ and $v$. For the arc $\gamma$, we have   $B( {\mathcal C},r, \gamma)= \big(g (\textswab{a}_{u*v} \textswab{a}_v)h \big) \kappa_{x,y}$, which  agrees with $B( {\mathcal C},r, \delta)$ due to the defining relation $\textswab{a}_{u*v}=\textswab{a}_v\textswab{a}_u \textswab{a}_v^{-1}$ in $\Adj(X)$. 	Alternative orientations, and signs of crossings follow similarly. $\blacksquare$
\end{proof}

Henceforth, we denote the Boltzmann weight by  $B( {\mathcal C}, r)$. 

\begin{definition}\label{def quandle gen cocycle invariant}
The family $\Phi_{\kappa} \big(D(K)\big) =  \big\{ \sum_{r}  B( {\mathcal C}, r) \big\}_{{\mathcal C}}$
is called the \index{quandle cocycle invariant}quandle cocycle invariant with respect to the generalized 2-cocycle $\kappa$. 
\end{definition}

The invariant agrees with the quandle cocycle invariant $\Phi_{\phi}\big(D(K)\big)$ defined  in \cite{MR1990571}  when the quandle module structure  on the coefficient group $A$ is trivial,  and with $\Phi_{\phi} \big(D(K)\big)$ defined  in \cite{MR1885217}, modulo the Alexander numbering convention in  the Boltzmann weight in \cite{MR1885217}, when the coefficient  group $A$ is  a $\Z [t, t^{-1}]$-module with the quandle module  structure is given by $\eta_{x,y}(a)=Ta$ and $\xi_{x,y}(b)=b-Tb$ for $x, y \in X$ and $a, b \in A$. Note that, in \cite{MR1885217, MR1990571}, the state-sum form is used, instead  of families.
\para 

Note that, in Definition \ref{def quandle gen cocycle invariant}, the quadruple $\big(D(K), X, A, \kappa \big)$ are 
chosen and fixed. 

\begin{figure}[!ht]
 \begin{center}
\includegraphics[height=11cm]{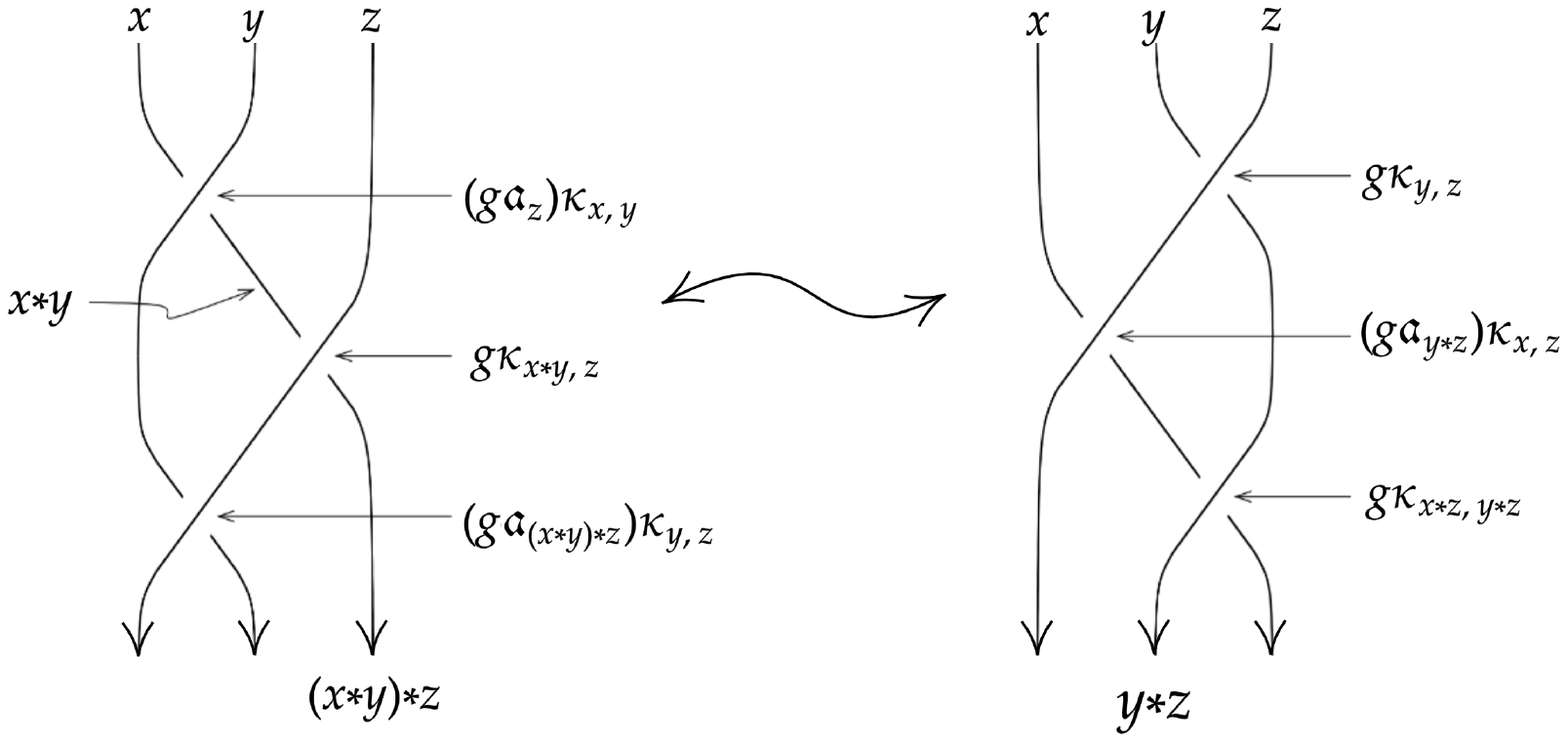}
\caption{Reidemeister move $R3$ and the $2$-cocycle condition.}
\label{weight3R} 
\end{center}
\end{figure}

\begin{theorem} \cite[Theorem 6.4]{MR2166720}
The family $\Phi_{\kappa} \big(D(K)\big)$ does not depend on the choice of the diagram $D(K)$ of the knot or the link $K$. 
\end{theorem}
\begin{proof} 
The proof is a routine check of Reidemeister moves, and it is straightforward for moves $R1$ and $R2$.  The move $R3$  is depicted, for one of the orientation choices, in Figure\ref{weight3R}, where $g$ denotes the sequence of colors of arcs that appear before the arc $\gamma$ intersects the crossings in consideration. Recall that $\kappa$ is a generalised 2-cocycle, and hence satisfy 
$$\textswab{a}_z \kappa_{x,y} + \kappa_{x*y, z} + \textswab{a}_{(x*y)*z} \kappa_{y,z}=  \kappa_{y,z} + \textswab{a}_{y*z}  \kappa_{x, z} + \kappa_{x*z,y*z}  $$
for all $x, y, z \in X$. It follows that the contribution of the Boltzmann weights from the three crossings involved does not change before  and after the Reidemeister move $R3$. The other orientation possibilities can be checked similarly. $\blacksquare$
\end{proof} 

The preceding theorem gives a well-defined invariant of knots and links taking values in $\mathbb{Z}[A]$, and we denote it by $\Phi_{\kappa}(K)$.

\begin{proposition} \label{coblemma}
The quandle cocycle invariant $\Phi_{\kappa}(K)$ depends only on the cohomology class of  the generalized quandle 2-cocycle $\kappa$. In particular, if $\kappa$ is a 2-coboundary, then $\Phi_{\kappa}(K)$ is trivial.
\end{proposition}

\begin{proof}
The argument is similar to the one given in \cite[Proposition 7.7]{MR1885217}. Let $\kappa$ be a generalized quandle 2-cocycle such that $\kappa=\partial^1(\lambda)$ for some $\lambda \in C^1(X;A)$. Then, by \eqref{Andruskiewitsch grana cohomology}, we have 
$$\kappa_{x,y}=\eta_{x, y} \big(\lambda(x)\big) -\lambda(x*y)+\xi_{x, y} \big(\lambda(y) \big) = \textswab{a}_y\lambda(x)  -\lambda(x*y) +\lambda(y)- \textswab{a}_{x*y} \lambda(y)$$ for all $x, y \in X$.  Consider a diagram of a knot or a link $K$ and remove a small neighbourhood of each of its crossings. Let $\gamma_1, \ldots, \gamma_m$ be the resulting arcs.
The end points of these arcs are located near the crossings, and depicted by bold dots (see Figure \ref{cocycle invariant from coboundary}). Assign each term of the right-hand side of the expression of $\kappa_{x,y}$ to these end points as  depicted in Figure \ref{cocycle invariant from coboundary}. The situation at two adjacent crossings is depicted in Figure \ref{cocycle invariant from coboundary 2}. Note that  the argument in $\lambda$  coincides with the color (a quandle element) of the arc. Then it can be seen that the terms assigned to  the two end points of each arc are the same, but with opposite signs.  Hence, the contribution to the state-sum for any coloring is $1$,
and the state-sum is a positive integer (which is the number of  colorings).	 Note that the two middle terms contribute zero and  the same works for two adjacent terms.  Hence, the cocycle invariant $\Phi_\kappa(K)$ is trivial. $\blacksquare$

\begin{figure}[!ht]
 \begin{center}
\includegraphics[height=4.5cm]{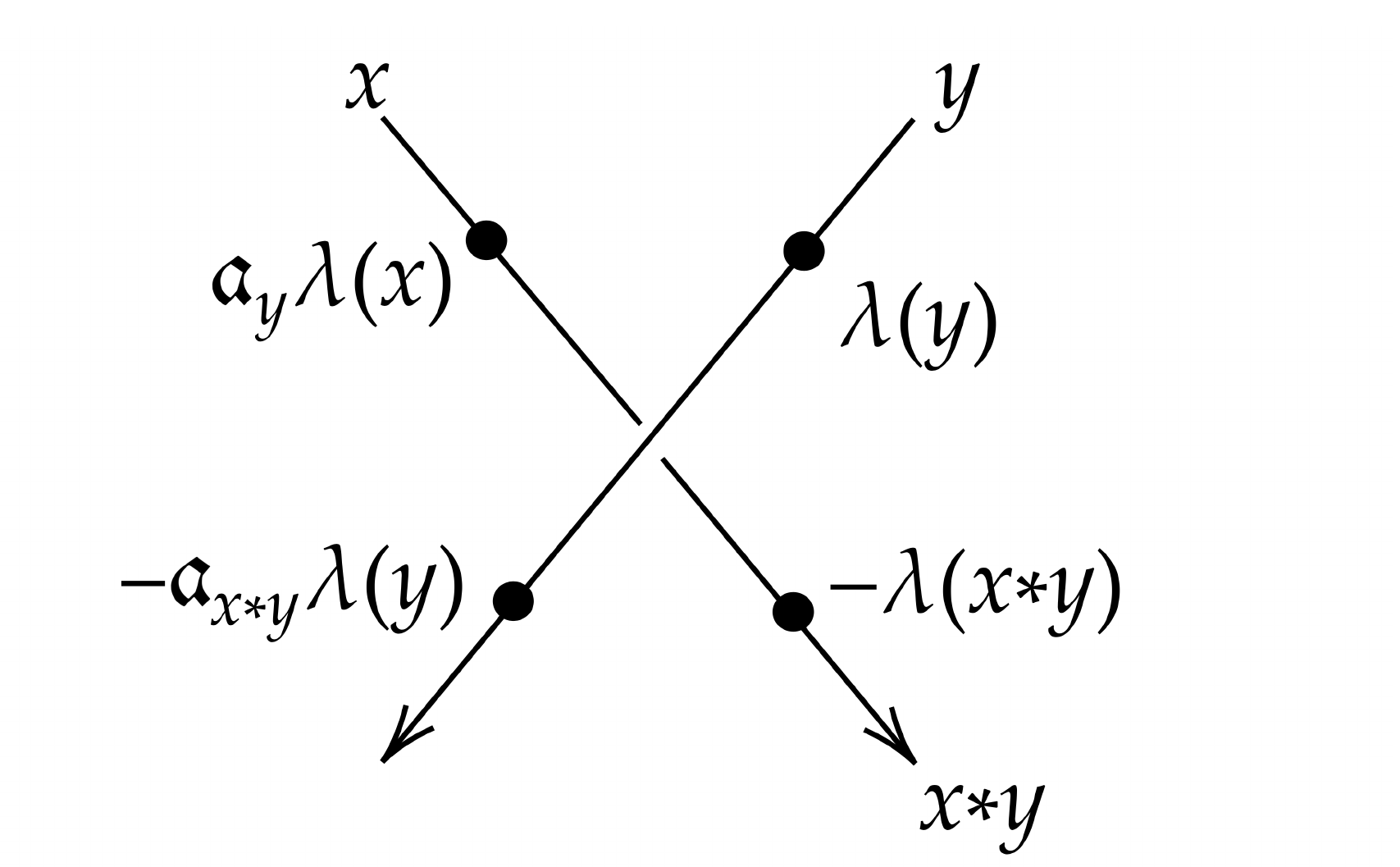}
\caption{Trivial invariant from a coboundary.}
\label{cocycle invariant from coboundary} 
\end{center}
\end{figure}
\begin{figure}[!ht]
 \begin{center}
\includegraphics[height=8cm]{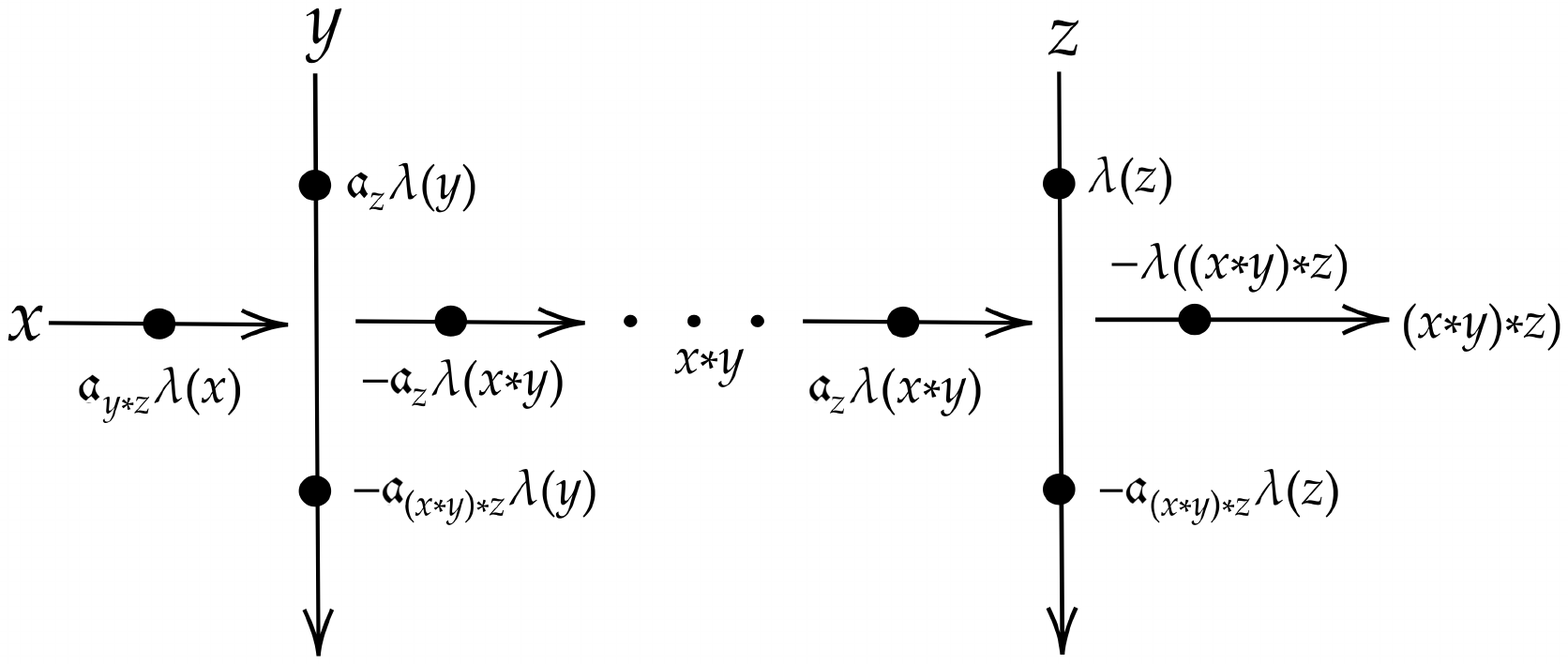}
\caption{Trivial invariant from a coboundary.}
\label{cocycle invariant from coboundary 2} 
\end{center}
\end{figure}
\end{proof}
\bigskip 


\subsection{Braids and quandle cocycle invariants} 
Recall that, the braid group $\B_n$ on $n \ge 2$ strands is the group generated by $(n-1)$ elements $\sigma_1,\ldots,\sigma_{n-1}$ with the following defining relations:
\begin{eqnarray*}
 \sigma_i\sigma_{i+1}\sigma_i & =& \sigma_{i+1}\sigma_i\sigma_{i+1} \quad \textrm{for}\quad 1 \le i \le n-2,\\
\sigma_i\sigma_j &=& \sigma_j\sigma_i \quad \textrm{for}\quad |i-j|\geq 2. 
\end{eqnarray*}

\begin{figure}[!ht]
 \begin{center}
\includegraphics[height=4.5cm]{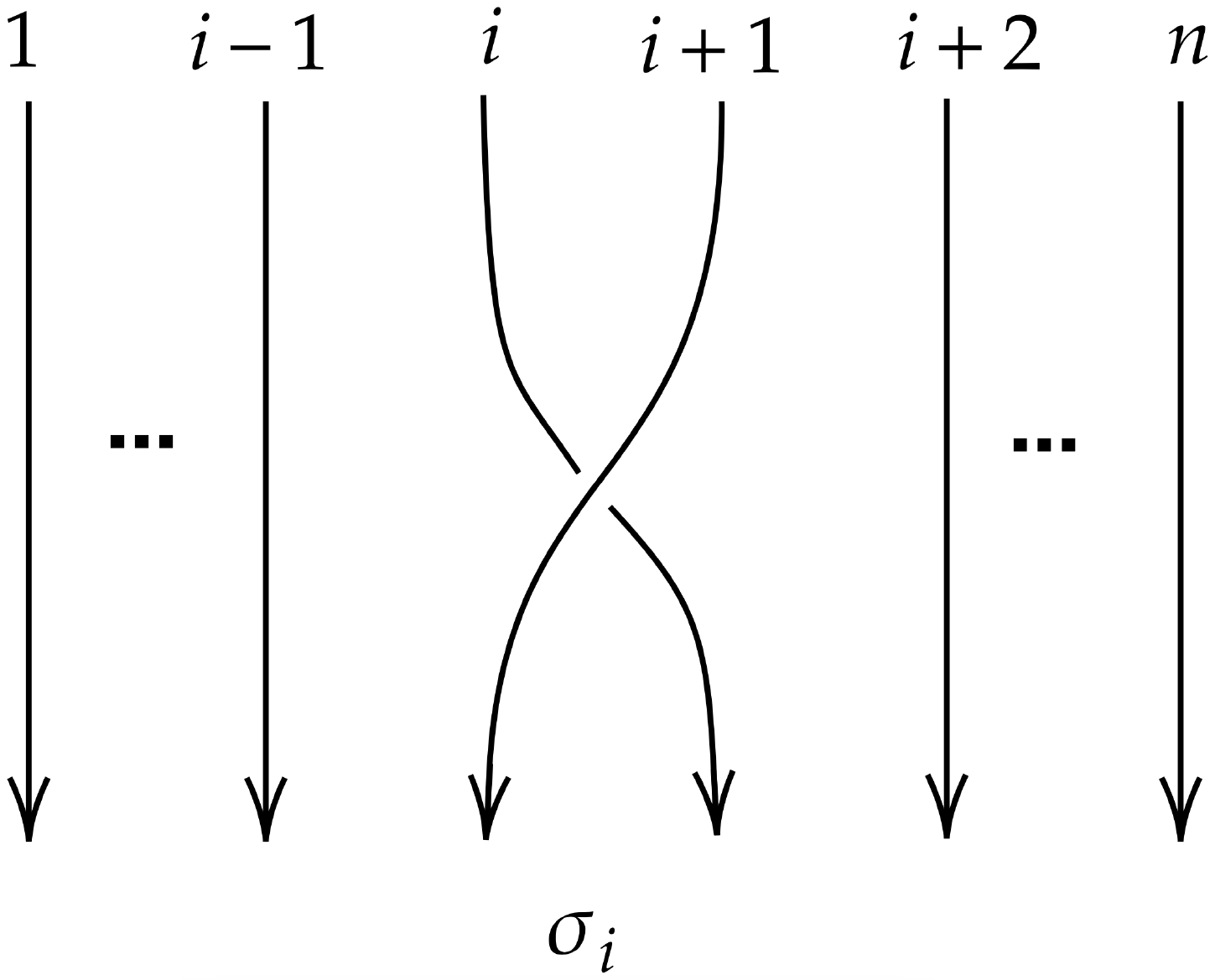}
\caption{A positive generator of the braid group.}
\label{positive generator of braid} 
\end{center}
\end{figure}

\begin{figure}[!ht]
 \begin{center}
\includegraphics[height=5cm]{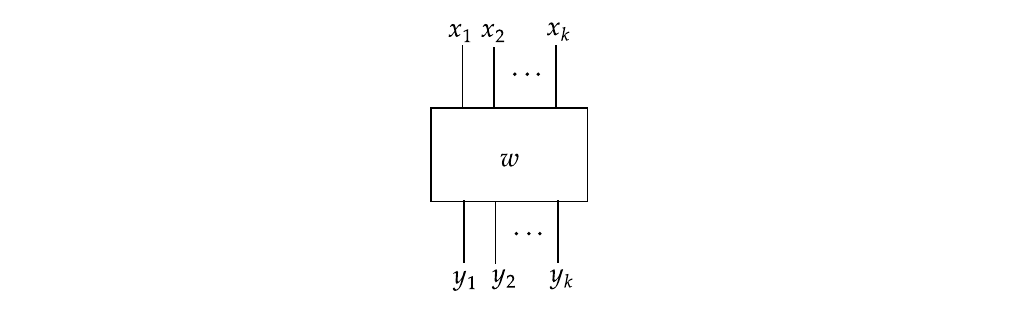}
\caption{Quandle coloring of a braid.}
\label{quandle coloring braid} 
\end{center}
\end{figure}

A braid word $w$ on $k$-strands is  a product of standard  generators $\sigma_1, \ldots, \sigma_{k-1}$  of the braid group $B_k$ and their inverses. Geometrically, 
a positive generator of the braid group is represented in Figure \ref{positive generator of braid} and the word $w$ is represented by a diagram in a rectangular box with  $k$ end points  at the top and $k$ end points at the bottom, where the strands  go downwards monotonically (see Figure \ref{quandle coloring braid}).  Each generator or its inverse is represented by a crossing in the diagram, and we use the same letter $w$ for a choice of such a diagram.  Let $\hat{w}$  denote the closure of  the geometric diagram of $w$.  Quandle colorings of  $w$  are  defined in  exactly the same manner as in the case of knots. However, the quandle elements used for coloring the top and the bottom of the diagram of  $w$ do not necessarily coincide.  The  quandle elements used at the top and the bottom of the diagram of  $w$ must coincide when we consider a coloring of the link  $\hat{w}$.

Let $X$ be a quandle and $\gamma_1, \ldots, \gamma_k$ be the top arcs of $w$. For a given vector $\vec{x}=(x_1, \ldots, x_k) \in X^k$, we assign the elements $ x_1, \ldots, x_k$ to the arcs $\gamma_1, \ldots, \gamma_k$ as their colors, respectively. Then from the definition,  a coloring ${\mathcal C}$ of $w$ by $X$ is uniquely determined
such that  
${\mathcal C}(\gamma_i)=x_i$ for $i=1, \ldots, k$.
Such a coloring ${\mathcal C} $ is called the coloring induced from $\vec{x}$. Let $\delta_1, \ldots, \delta_k$ be the arcs at the bottom.
Let 
$\vec{y}=(y_1, \ldots, y_k)= \big({\mathcal C}(\delta_1), \ldots, {\mathcal C}(\delta_k) \big)\in X^k$ 
be the colors assigned to the bottom arcs, that are uniquely determined from 
$\vec{x}$. We denote this situation by a left-action 
$$\vec{y} =  w \cdot \vec{x}.$$ The colors $\vec{x}$ and $\vec{y}$ are 
called {\it top} and {\it bottom colors}
or {\it color vectors}, respectively.
See Figure \ref{quandle coloring braid} for a diagrammatic illustration.

Let $X$ be a quandle and $G$ be a quandle module over $X$.  Let $$\alpha= \eta + \xi + \kappa$$ be a  dynamical 2-cocycle which acts on $G^2$ by
$$\alpha_{x,y}(a,b) =\eta_{x,y} (a )  +  \xi_{x,y}(b) + \kappa_{x,y}$$
for any $(x,y) \in X^2$ and  $(a,b) \in G^2$. Let 
$\tilde{X}=G \times_{\alpha} X$ be the associated dynamical extension. If $\vec{r} =\big( (a_1, x_1), \ldots, (a_k, x_k) \big)$ and $\vec{s} =\big( (b_1, y_1), \ldots, (b_k, y_k) \big)$
are the top and the bottom colors of a braid $w$  in $\B_k$ by elements of $\tilde{X}$,
respectively, then we write this situation as 
$$\vec{b} =  M(w, \vec{x}) \cdot \vec{a},$$
where $\vec{a}=(a_1, \ldots, a_k)$ and  $\vec{b}=(b_1, \ldots, b_k)$ are vectors of elements from $G^k$.  Hence, $M(w, \vec{x})$ represents a map 
$$M(w, \vec{x}): G^k \rightarrow G^k.$$

\begin{lemma} \label{indeplemma}
If $w,w'$ are two braids representing the same element in the braid group $\B_k$, then $M(w, \vec{x})=M(w', \vec{x}): G^k \rightarrow G^k$ for each vector $\vec{x} \in X^k$.
\end{lemma}

\begin{proof}
The invariance under the braid relations can be checked using the definition. In particular,  we have
\begin{eqnarray*} 
	\lefteqn{  M \big( \sigma_1 \sigma_2 \sigma_1 , (x,y,z) \big)(a,b,c)  } \\
	&= & \big(c,\  \eta_{y,z}( b) +   \xi_{y,z}(c) + \kappa_{y,z} , \
	\eta_{x*y, z} \eta_{x,y}(a) +  \eta_{x*y, z} \xi_{x,y}(b)
	+  \xi_{x*y, z}(c) +  \eta_{x*y, z}( \kappa_{x,y}) + \kappa_{x*y, z} \big), \\
	\lefteqn{ M \big( \sigma_2 \sigma_1 \sigma_2 , (x,y,z) \big) (a,b,c)  } \\
	&= & \big(c, \   \eta_{y,z}(b) +   \xi_{y,z}(c) + \kappa_{y,z}, \,\eta_{x*z, y*z} \eta_{x,z}(a) +  \xi_{x*z, y*z} \eta_{y,z}(b)+  ( \eta_{x*z, y*z} \xi_{x,z} +  \xi_{x*z, y*z} \xi_{y,z} )(c) \\ & &  
	+ \eta_{x*z, y*z}( \kappa_{x,z}) + \xi_{x*z, y*z} ( \kappa_{y,z}) +
	\kappa_{x*z, y*z} \big),
\end{eqnarray*}
and the equality follows from the quandle module conditions and  the generalised $2$-cocycle condition. $\blacksquare$
\end{proof}

We note that $M(-, \vec{x}): \B_k \rightarrow \Map(G^k, G^k)$ is not a representation of the braid group since it depends on the coloring of the braid $w$ by elements of $X$. However, when the coloring by $X$ is trivial, that is, $x_1=x_2=\cdots = x_k$ and $\kappa=0$, then
it is a representation of the braid group.  In that case, the map $M(-, \vec{x}):  \B_k \rightarrow \Map(G^k, G^k)$ is called a \index{colored representation}{\it colored representation}. 
\para 

\begin{figure}[!ht]
 \begin{center}
\includegraphics[height=5cm]{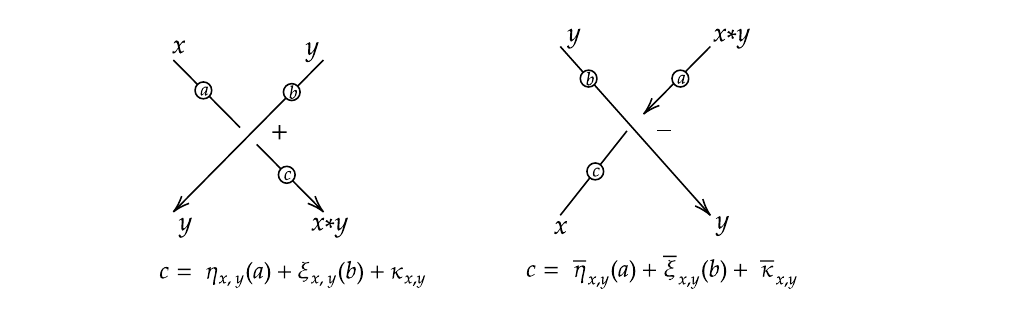}
\caption{Quandle coloring of a braid generator.}
\label{quandle coloring braid generator} 
\end{center}
\end{figure}

Let $X$ be a quandle and $G$ a quandle module over $X$. In Figure \ref{quandle coloring braid generator}, the left figure represents the standard braid generator $\sigma_j$, and the right represents its inverse $\sigma_j^{-1}$.  Further, the colors by $X$ are assigned to arcs.  Elements of $G$ are put in small circles on arcs.  We  imagine these circles sliding up through a crossing, at which  the dynamical cocycle $\alpha$ acts and changes the circled elements when a circle goes under a crossing.  Going over a crossing does not change the circled element. 

\begin{figure}[!ht]
 \begin{center}
\includegraphics[height=8.5cm]{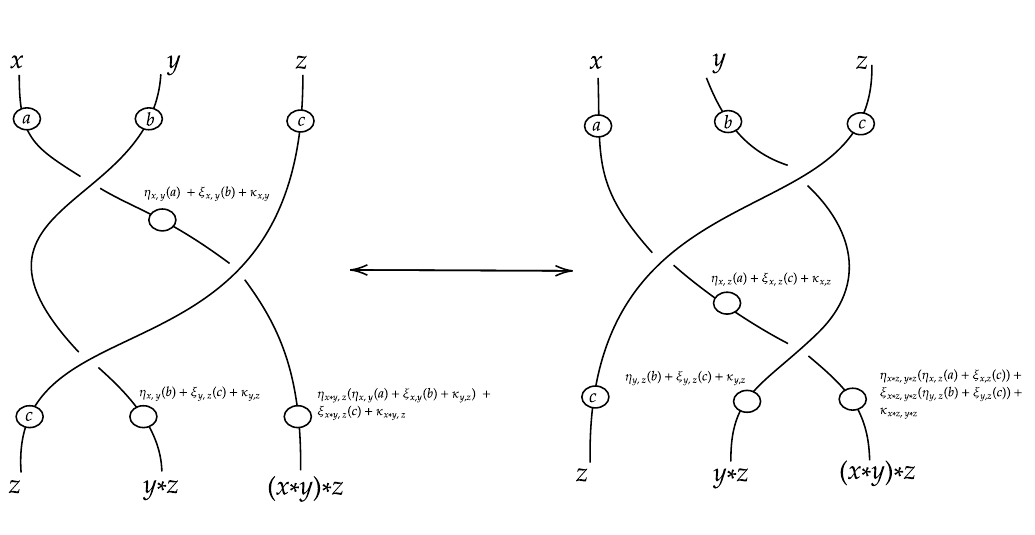}
\caption{Quandle coloring and Reidemeister move $R3$.}
\label{quandle coloring R3 move} 
\end{center}
\end{figure}

Next, we have the following result \cite[Theorem 4.2]{MR2166720}.

\begin{theorem}
Let $X$ be a quandle and $G$ a quandle module over $X$. Let $\alpha= \eta + \xi $ be a dynamical 2-cocycle which acts on $G \times G$ by
$$ \alpha_{x,y}(a,b) =\eta_{x,y} (a )  +  \xi_{x,y}(b)$$
for any $(x,y) \in X^2$ and  $(a,b) \in G \times G$.  Let $L$ be a link represented as the closure $\hat{w}$ of a braid  word $w$ on $k$-strands  and  ${ \mbox{Col}_X(L) }$ be the set of colorings of $L$ by the quandle $X$. For $\mathcal{C} \in  \mbox{Col}_X(L)$, 	let $\vec{x}$ be the color vector of top strings of $w$ that is the restriction of $\mathcal{C}$.  Then the family  
$$\tilde{\Phi}(X, \alpha\ ; L)= \big\{G^k/ \im (M(w,\vec{x})-I)  \big\}_{{\mathcal C} \in  {\rm Col}_X(L)}$$ of isomorphism classes of modules presented by the maps $\{M(w,\vec{x})-I \}$, where $I$ denotes the identity map,  is independent of the choice of the braid word $w$, and thus defines an invariant of the link $L$.
\end{theorem}

\begin{proof}
By Lemma~\ref{indeplemma}, the colored representation $M(w) :=M(w, \vec{x})$ does not depend on the choice of a braid word. We use Markov's theorem to prove the assertion.
First note that the set of colorings remains unchanged by a  stabilization, in the sense that  a coloring $(x_1, \ldots, x_k)$ on the top strings  of a braid word $w$ extends uniquely to a coloring $(x_1, \ldots, x_k, x_{k})$ of the stabilization $w\sigma_{k}^{\pm 1}$ of $w$. There is a bijection of colorings between conjugate words as well.
Hence, it is sufficient to prove that,  for a given coloring, the isomorphism classes of modules defined in the statement of the theorem remains unchanged by 
conjugation and stabilization for an induced coloring. The  invariance under conjugation is seen from  the fact that  conjugation by a braid word 
induces a  conjugation  by a matrix, and the module  is isomorphic under conjugate presentation matrices. Thus, it remains to investigate the stabilization. 
\para 

We represent maps $G^k \to G^k$ by $k \times k$ matrices, whose each entry itself represents a map of $G$.   The braid generator $\sigma_k$ in the stabilization  $w\sigma_{k}^{\pm 1}$ 
is represented by   the  matrix 
$$M(\sigma_k)=I_{k-1} 
\oplus 
\left[ \begin{array}{cc} O & I\\ W & I - W \end{array} \right],$$ where $I$ denotes the identity map of $G$ and $I_k$ denotes the identity map on $G^k$.   This is because the $k$-th and $(k+1)$-th strands receive  the same color. The block matrix  $\left[ \begin{array}{cc} O & I \\ W &  I -W \end{array} \right]$ represents the map 
$$(a,b) \mapsto \big(b, \,\alpha_{x,x}(a,b) \big)= \big(b, \,\eta_{x,x}(a) + \xi_{x,x}(b) \big),$$
where $x=x_k$. Here, $W$ corresponds to the action by $\eta$, and hence it is an isomorphism of $G$.  The fact that $\xi_{x,x}$
corresponds to the matrix $I - W$  follows from the condition that $\eta_{x,x}+\xi_{x,x}$ is the identity map in a quandle algebra. We express $M(w)$ as  a $(k+1) \times (k+1)$ matrix of  maps, that is,
$$M(w)= \left[ \begin{array}{cccc} M_{11} & \cdots &  M_{1k} &O\\ 
	\vdots & & \vdots &  \vdots\\
	M_{k1} &  \cdots & M_{kk} & O
	\\ O &  \cdots & O & I \end{array} \right],$$
where $w$ is regarded as  a braid word on $(k+1)$-strands after the stabilization, $M_{ij}$ ($1 \leq i, j \leq k$) are maps of $G$ and $I$ denotes the identity map on $G$. 
Then $M(w \sigma_k)$ is represented by 
the matrix $$ M(w) M(\sigma_k)=
\left[ \begin{array}{ccccc} 
	M_{11} &  \cdots  & M_{1\ (k-1)} &O & M_{1 k} \\ 
	\vdots &  & \vdots & \vdots & \vdots \\
	M_{k1} &  \cdots & M_{k\ (k-1)} & O & M_{k k} \\
	O& \cdots  & O& W & I - W 
\end{array} \right] .  $$
Hence, we have 
$$ M(w) M(\sigma_k) - I_{k} = 
\left[ \begin{array}{ccccc} 
	M_{11} - I  &  \cdots  & M_{1\ (k-1)} &O & M_{1 k} \\ 
	\vdots &  & \vdots & \vdots & \vdots \\
	M_{(k-1)\ 1} & \cdots &  M_{(k-1)\ (k-1)} - I & O &  M_{(k-1)\ k} \\
	M_{k1} &  \cdots & M_{k\ (k-1)} & -I  & M_{k k} \\
	O& \cdots  & O& W &  - W 
\end{array} \right] ,
$$
which can be column reduced to 
$$
\left[ \begin{array}{ccccc} 
	M_{11} - I  &  \cdots  & M_{1\ (k-1)} &O & M_{1 k} \\ 
	\vdots &  & \vdots & \vdots & \vdots \\
	M_{(k-1)\ 1} & \cdots &  M_{(k-1)\ (k-1)} - I & O &  M_{(k-1)\ k} \\
	M_{k1} &  \cdots & M_{k\ (k-1)} & -I  & M_{k k} -I \\
	O& \cdots  & O& W &  O 
\end{array} \right]. 
$$
This can further be row reduced to 
$$
\left[ \begin{array}{ccccc} 
	M_{11} - I  &  \cdots  & M_{1\ (k-1)} &O & M_{1 k} \\ 
	\vdots &  & \vdots & \vdots & \vdots \\
	M_{(k-1)\ 1} & \cdots &  M_{(k-1)\ (k-1)} - I & O &  M_{(k-1)\ k} \\
	M_{k1} &  \cdots & M_{k\ (k-1)} & O  & M_{k k} -I \\
	O& \cdots  & O& W &  O 
\end{array} \right], 
$$
which represents a module that is isomorphic to the one represented by $M(w)- I_k$,  since $W$ is an isomorphism.  The invariance under stabilization by $\sigma_k^{-1}$ follows similarly, and the proof is complete.  $\blacksquare$
\end{proof}

The preceding theorem leads to  the following definition.

\begin{definition} 
The family of modules	$\tilde{\Phi}(X, \alpha \ ; L) =\big\{ G^k/\im(M(w,\vec{x})-I) \big\}_{{\mathcal C} \in {\rm Col}_X(L) } $
is called the \index{quandle module invariant}{quandle module invariant}.
 \end{definition}

\begin{definition}
Let $X$ be a quandle and $A$ a $\mathbb{Z}(X)$-module with  
$$\eta_{x,y }(a)= \textswab{a}_y a \quad  \textrm{and} \quad  \xi_{x,y}(b)= b-(\textswab{a}_{x*y})b$$ for $x, y \in X$ and $a, b \in A$. Let $\vec{x}=(x_1, \ldots, x_k)\in X^k$ for some positive integer $k$. A \index{weighted sum}{weighted sum} on $A^k$ with respect to $\vec{x}$ is
defined by
$$ \mbox{WS}_{\vec{x}}(\vec{a})= \sum_{i=1}^k u_i a_i = (\textswab{a}_{x_k} \cdots \textswab{a}_{x_2}  ) a_1 +   (\textswab{a}_{x_k} \cdots \textswab{a}_{x_3}) a_2 + \cdots +  \textswab{a}_{x_k} a_{k-1} +  a_k,$$
where $\vec{a} =(a_1, \ldots, a_k) \in A^k$ and $u_i=\textswab{a}_{x_k} \cdots \textswab{a}_{x_{i+1}}$ for $i=1, \ldots, k-1$ and $u_k=1$.
\end{definition}

Let $\alpha^0_{x,y}= \eta_{x,y} + \xi_{x,y} $ be a quandle 2-cocycle with coefficient group $A$ and  $\kappa=0$. Then any  color vector $\vec{a}=(a_1, \ldots, a_k) \in A^k$ at the bottom strings of $w$ uniquely determines a color vector  $\vec{b}=(b_1, \ldots, b_k) \in A^k$ at the top strings of $w$ with respect to $\alpha$, that is, $\vec{b}=M(w, \vec{x}) \cdot \vec{a}$. 

\begin{lemma} \label{wslemma}
Let  $\alpha^0_{x,y}= \eta_{x,y} + \xi_{x,y} $ be a quandle $2$-cocycle  and $\vec{b}=M(w, \vec{x}) \cdot \vec{a}$, as above. Then  $\mbox{WS}_{\vec{x}}(\vec{a})=\mbox{WS}_{\vec{y}}(\vec{b})$.
\end{lemma}

\begin{proof} 
It is sufficient to prove the statement for a braid group generator and its inverse. Since the inverse case is similar,  we consider 
the case when $w=\sigma_i$ for some $i$ with $1 \leq i < k$.  The only difference in the weighted sum is the $i$-th and $(i+1)$-st terms. 
Let $u=\textswab{a}_{x_k} \cdots \textswab{a}_{x_{i+2}}$.  For the bottom colors, the $i$-th and $(i+1)$-st terms are 
$$ \mbox{WS}_{\vec{x}}(\vec{a})=
\cdots + u \textswab{a}_{x_{i+1}}  a_i + u   a_{i+1} + \cdots .$$
On the other hand, we have
$$
\mbox{WS}_{\vec{y}}(\vec{b}) = 
\cdots + u (\textswab{a}_{x_i * x_{i+1}}) a_{i+1} 
+ u 
\big( \textswab{a}_{x_{i+1}} a_i + (1-  \textswab{a}_{x_i * x_{i+1}}) a_{i+1} \big)\ + \cdots.$$
We only need to check that  $\textswab{a}_{x_{i+1}}  a_i + a_{i+1} = \textswab{a}_{y_{i+1}}  b_i + b_{i+1}$, but this comes from the fact that $y_i=x_{i+1}$, $y_{i+1}=x_i*x_{i+1}$, $b_i=a_{i+1}$ and $ b_{i+1}=\eta_{x_1,x_{i+1}}a_i +\xi_{x_i,x_{i+1}}a_{i+1}$. This completes the proof of the lemma. $\blacksquare$
\end{proof}

\begin{theorem} \cite[Theorem 6.9]{MR2166720}\label{wsthm}
Let $w$ be a braid on $k$-strands with crossings $r_1, \ldots, r_h$ and $\vec{x} \in X^k$ a bottom color vector. Let $A$ be a $\mathbb{Z}(X)$-module,  $\alpha_{x,y}= \eta_{x,y} + \xi_{x,y} + \kappa_{x,y}$ be a quandle $2$-cocycle with coefficients in $A$ (with possibly non-zero $\kappa$) and $\vec{b}=M(w, \vec{x}) \cdot \vec{a}$. Then we have 
	$$ \mbox{WS}_{\vec{y}}(\vec{b}) - \mbox{WS}_{\vec{x}}(\vec{a})
	= \sum_{\ell=1}^{h}  B({\mathcal C},r_\ell). $$
\end{theorem} 

\begin{proof} 
Let $M^0(w, \vec{x})$ denote the map corresponding to  the $2$-cocycle  $\alpha^0_{x,y}= \eta_{x,y} + \xi_{x,y}$, which is obtained from $\alpha$ by setting $\kappa=0$. The theorem follows by induction once we prove it for  a braid generator and its inverse, and we only consider the case
$w=\sigma_i$. Using Figure \ref{braid coloring}, we compute
$$
\mbox{WS}_{\vec{y}}(\vec{b})=\mbox{WS}_{\vec{y}} \big( M(\sigma_i, \vec{x} )\cdot  \vec{a} \big)=	\mbox{WS}_{\vec{y}} \big( M^0(\sigma_i, \vec{x} )\cdot \vec{a} \big)	+ B({\mathcal C}, \sigma_i)= \mbox{WS}_{\vec{x}}( \vec{a} ) - B({\mathcal C}, \sigma_i),
$$
where the first equality follows from the definitions, and the second equality follows from Lemma~\ref{wslemma}.
\begin{figure}[!ht]
 \begin{center}
\includegraphics[height=4.5cm]{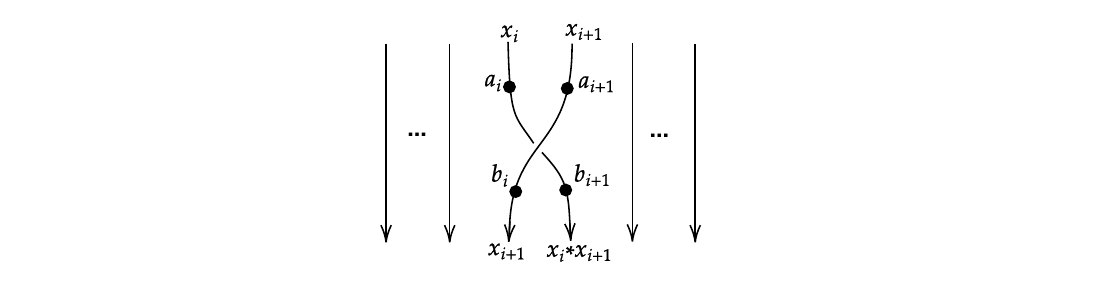}
\caption{Coloring of the braid generator $\sigma_i$.}
\label{braid coloring} 
\end{center}
\end{figure}
Note that the braid generator represents a negative crossing in the definition of the quandle module invariant, and a positive crossing in the cocycle invariant, so that there is a negative sign  for the weight $B({\mathcal C}, \sigma_i)$.  $\blacksquare$
\end{proof}

\begin{theorem} \cite[Theorem 6.10]{MR2166720}\label{colorextthm}
If a  coloring ${\mathcal C}$ of a link $L$ by $X$  lifts to a coloring of $L$ by the extension $E= A \times_{\alpha} X$, then the coloring ${\mathcal C}$ contributes  	a trivial term (that is, an integer in $\mathbb{Z}[A]$)
to  the generalized cocycle invariant $\Phi_{\kappa} (L) $.
\end{theorem}

\begin{proof}
A given coloring agrees  on the bottom and top strings of a braid, so that  $\vec{x}=\vec{y}$.  If the given coloring lifts to $E= A \times_{\alpha} X$, then we have  $ \vec{b}= M(\sigma_i, \vec{x} )\cdot  \vec{a}=  \vec{a}$, and in particular,  $\mbox{WS}_{\vec{y}}( \vec{b})= \mbox{WS}_{\vec{x}}( \vec{a} )$. The theorem now follows from Theorem~\ref{wsthm}.
$\blacksquare$
\end{proof}

\begin{remark} {\rm 
Let $1 \rightarrow A \rightarrow E \rightarrow H \rightarrow 1$ be a split short exact  sequence of groups, where
$A$ is  abelian. Let $K$ be a knot, $G(K)=\pi_1(\mathbb{S}^3\setminus K)$ its knot group and $\rho :G(K) \rightarrow H$ a group homomorphism. Since there is an action of $H$ on $A$ via conjugation in $E$, it induces an action of $G(K)$ on $A$ given via the homomorphism $\rho$.  In \cite{MR1331749}, Livingston examined the situation when there is a lift $\tilde{\rho}: G(K) \rightarrow E$  of $\rho$, thereby generalizing Perko's theorem \cite{MR0377852} that any homomorphism $\rho: G(K) \rightarrow \Sigma_3$ lifts to a homomorphism $\tilde{\rho}: G(K) \rightarrow \Sigma_4$. In this case, the corresponding extension is 
$$1 \to \mathbb Z_2 \times \mathbb Z_2 \to \Sigma_4 \to \Sigma_3 \to 1,$$
 where $\Sigma_3$ acts on $\mathbb Z_2 \times \mathbb Z_2 \cong \{ (1), (12)(34), (13)(24), (14)(23) \}$  via the conjugation and the obvious section $s:\Sigma_3 \rightarrow \Sigma_4$. It follows from \cite[Proposition 2.3]{MR0672956} that  there is a one-to-one correspondence between $A$-conjugacy classes of lifts of $\rho$  and the cohomology group $\Ho^1_{\rm Grp} \big(G(K); \{A\}\big)$, where the  coefficients $\{A\}$ are twisted by the action of $G(K)$  on $A$. Since a quandle coloring of $K$ by $\Conj(H)$, that is, a quandle homomorphism from the knot quandle $Q(K)$ to $\Conj(H)$ induces a group homomorphism $\rho: G(K) \to H$ as above, we have the following  cohomological characterization of lifting of colorings.} 
\end{remark}

\begin{corollary}
Let $K$ be a knot and $1 \rightarrow A \rightarrow E \rightarrow H \rightarrow 1$ a split short exact  sequence of groups, where $A$ is abelian. If a coloring $\mathcal{C}$ of $K$ by $\Conj(H)$ lifts to $\Conj(E)$, then it  contributes a constant term to the cocycle invariant of Theorem~\ref{colorextthm}.
\end{corollary}

We now present a computational method by viewing a knot as a closed braid. We take a positive crossing as a positive generator of the braid group (see Figure \ref{positive generator of braid}). Let $w=w_1 \cdots w_h$ be a braid word,  where $w_s=\sigma_{j(s)}^{\epsilon(s)}$ is a standard generator of the braid group or its inverse for each $s=1, \ldots, h$. The braids are oriented downwards and the normals point  to the right. Then the target region is to the right  of each crossing. 

\begin{figure}[!ht]
 \begin{center}
\includegraphics[height=5.5cm]{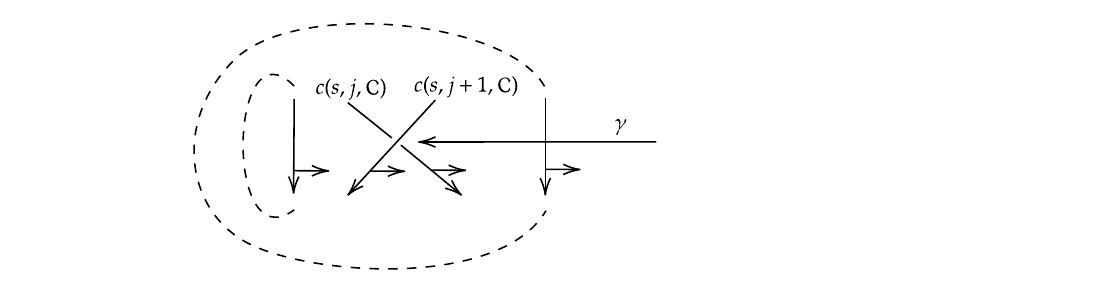}
\caption{Colors and arcs for closed braids.}
\label{braidconvention} 
\end{center}
\end{figure}

Let $X$ be a quandle and ${\mathcal C}$ a coloring of the closure $\hat{w}$ of $w$  by $X$, where the closure $\hat{w}$ is taken to the left as depicted in Figure \ref{braidconvention}. Let $c(s, j,  {\mathcal C})$ be the color of the $j$-th string from the left, immediately above the $s$-th crossing $w_s$ for the coloring ${\mathcal C}$.  Take the right-most region as the region at infinity. Now, take the arc  $\gamma$ to the target region as an arc that goes 
horizontally from the right to the left of the crossing as in Figure \ref{braidconvention}.
Then the Boltzmann weight at the  $s$-th crossing is 
given by 
$$ B(  {\mathcal C}, w_s)
= \big(\textswab{a}_{c(s, n, {\mathcal C})}  \cdots   \textswab{a}_{c(s, j+2, {\mathcal C})} \big) \,
\kappa_{c(s, j,  {\mathcal C}),  c(s, j+1,  {\mathcal C})}
$$ 
for a positive crossing,  and the $2$-cocycle evaluation is replaced by  $\kappa_{c(s+1, j,  {\mathcal C}),  c(s, j,  {\mathcal C})}$ for a negative crossing, and with a negative sign in front. Note that the expression in front of $\kappa$ represents  the action of this group element on the coefficient group, so that 
in the case of wreath product, for example, this can be written in terms of matrices. Also, if $j=n-1$, then the  expression in front is understood to be
empty.  This formula can be used to evaluate the invariant for knots  in closed braid form, and can be implemented  in a computer.

\begin{example} \label{r3classicalex} {\rm   
The preceding formula for the dihedral quandle $\R_3$ was implemented in {\it Maple} and confirmed also by {\it Mathematica}, thus giving the following results.
\para
 We write $\R_3=\{ 1,2,3 \}$ and identify it as a subquandle of $\Conj(\Sigma_3)$ via the map $1 \mapsto (2\ 3)$,  $2 \mapsto (1\ 3)$ and $3 \mapsto (1\ 2)$. Consider the wreath product
$\mathbb{Z}_q \wr \Sigma_3=(\mathbb{Z}_q)^3 \rtimes \Sigma_3$. For $q=0$, any $2$-cocycle gives trivial invariant (each coloring contributes the zero vector and the family of vectors is a family of zero vectors) for all $3$-colorable knots in the knot table up to $9$ crossings.  Let $h(i,j)=\big( f_1(i,j),  f_2(i,j),  f_3(i,j)\big)$ denote a vector valued 2-cochain. Then,  for $q=3$, setting
\begin{eqnarray*}
f_1(3, 1)\ =\ f_3(3,1)\ =\  f_1(2,3)\ =\  f_1(2,1)\ =\  f_3(2,1) &= & 1, \\
f_2(3, 1)\  =\ f_2(3,2)\ =\ f_2(1,3)\ =\  f_1(1,2)\ =\   f_2(1,2) &= & 2
\end{eqnarray*} 
and taking all the other values as zeros, defines a 2-cocycle. For positive integers $m$ and $n$, let $m_n$ denote the $n$-th knot in the knot table \cite{RolfsenTable} with $m$ crossings. Let $\{ \sqcup_\ell (a,b,c) \}$ represent the family consisting of $\ell$ copies of  $(a,b,c)$. Then we obtain the following results:
\begin{enumerate}
\item $\Phi_{\kappa}(K)=\big\{ \sqcup_9 (0,0,0) \big\}$ for  $K=  6_1,  8_{10}, 8_{11}, 8_{20}, 9_1, 9_6, 9_{23}, 9_{24}$.
\item $\Phi_{\kappa}(K)=\big\{ \sqcup_3 (0,0,0), \,\sqcup_6 (1,1,1)  \big\}$ for  $K=3_1, 7_4, 7_7, 9_{10}, 9_{38}$.
\item $\Phi_{\kappa}(K)=\big\{ \sqcup_3(0,0,0), \, \sqcup_6 (2,2,2) \big\}$ for $K=8_5, 8_{15}, 8_{19}, 8_{21}, 9_2, 9_4, 9_{11}, 9_{15}, 9_{16}, 9_{17}, 9_{28}, 9_{29}, 9_{34}, 9_{40}$.
\item $\Phi_{\kappa}(K)=\big\{  \sqcup_9 (0,0,0), \, \sqcup_{18} (1,1,1) \big\} $ for $K= 9_{35}, 9_{47}, 9_{48}$.
\item $\Phi_{\kappa}(K)=\big\{  \sqcup_3 (0,0,0), \, \sqcup_{12} (1,1,1), \, \sqcup_{12} (2,2,2) \big\} $ for $K=8_{18}$.
\item $\Phi_{\kappa}(K)=\big\{  \sqcup_{15} (0,0,0), \, \sqcup_6(1,1,1), \, \sqcup_6 (2,2,2) \big\} $ for $K=9_{37}, 9_{46}$.
\end{enumerate}
In the original cohomology theory with trivial action on the coefficient group, we have $\Ho^2_{\rm Q}(\R_3; A)=1$ for any coefficient group $A$ (see Example \ref{dihedral quandle of order three homology}). Hence, $\R_3$ gives rise to a trivial cocycle invariant. The example shows that with the wreath product action, $\R_3$ has non-trivial second cohomology group  and gives rise to a non-trivial cocycle invariant. Hence, the cohomology theory with non-trivial actions on coefficients is more general and strictly stronger than the original cohomology theory.} 
\end{example}
\bigskip
\bigskip


\chapter{Computations of homology and cohomology of solutions}\label{chapter homology computations}

\begin{quote}
The computation of the homology and cohomology of racks, quandles, and other solutions to the Yang--Baxter equation is typically challenging. Although  the number of computational techniques for these structures is relatively limited compared to homology theories for other algebraic structures, there is a significant body of computational results available in the literature. This chapter surveys key findings related to the computational aspects of homology and cohomology for solutions to the Yang--Baxter equation.
\end{quote}
\bigskip

\section{Low dimensional homology of racks and quandles}
A thorough examination of the existing literature on homology and cohomology computations would necessitate a dedicated text focused exclusively on the subject. For a more in-depth exploration of the topic, we recommend referring to the following references: \cite{MR4038329, MR4065999, MR3613433, MR3342681, MR2854230, MR4602407, MR4310563, MR3196064, MR4491426, MR3226796, MR2988494, MR2657689, MR3898509, MR2521709, MR2494367, MR2673749, MR2777023, MR3042590, MR4400200, MR3571388, MR3937311}. We use the convention that, whenever the coefficients of (co)homology groups are not specified, they are assumed to be $\mathbb{Z}$ coefficients. 
\para 

The free parts of the rack, degenerate, and quandle homology groups of finite racks and quandles are well understood (\cite[Corollary 4.3]{MR1948837} and \cite[Theorem 2]{MR1952425}).

\begin{theorem}  \label{Theorem 1.4}
Let $\mathcal{O}$ be the set of orbits of a finite rack or a quandle $X$. Then the following assertions hold:
	\begin{enumerate}
		\item $\rank \big( \Ho_n^{\rm R}(X) \big)=|\mathcal{O}|^{n}$.
		\item $\rank \big(\Ho_n^{\rm Q}(X) \big)=|\mathcal{O}| \big(|\mathcal{O}|-1 \big)^{n-1}$.
		\item $\rank \big( \Ho_n^{\rm D}(X) \big)=|\mathcal{O}|^{n}-|\mathcal{O}| \big(|\mathcal{O}|-1 \big)^{n-1}$.
	\end{enumerate}
\end{theorem}

If $|\mathcal{O}|=1$, that is, $X$ is a connected rack, then $\rank \big(\Ho_n^{\rm R}(X) \big)=1$. For latin quandles, Niebrzydowski and Przytycki provided a natural description of the kernel of the second boundary map  \cite[Lemma 2.1]{MR2777023}.

\begin{lemma}
Let $X$ be a latin quandle and $a_0 \in X$ a fixed element. Let $\langle(a_0,a) \rangle$ be the subgroup of $C^{\rm Q}_2(X)$ generated by elements of the form $(a_0,a)$ for any $a\in X$.  Then the quotient map $C^{\rm Q}_2(X) \rightarrow	C^{\rm Q}_2(X)/ \langle (a_0,a)\rangle$ is an isomorphism when restricted to $\ker( \partial_2)$.
\end{lemma}

As a consequence, the following result is obtained \cite[Corollary 2.2]{MR2777023}.

\begin{corollary} 
If $X$ is a latin quandle and $a_0 \in X$ a fixed element, then $$\Ho^{\rm Q}_2(X) \cong \mathbb{Z}^{X\times X} / \big\langle (a,a),(a_0,b), \im(\partial_3) \big\rangle.$$
\end{corollary}

Recall that, if $G$ is a group, then $\Core(G)$ is an involutory quandle with the binary operation $g*h=hg^{-1}h$. If $G$ is abelian, then $\Core(G)=\T(G)$, the Takasaki quandle of $G$. Note that the Takasaki quandle of a group of odd order is latin. In order to state the next result, we recall the definition of the exterior square of a module.

\begin{definition}
	Let $\mathbb{k}$ be a commutative ring and let $A$ be a $\mathbb{k}$-module.  Then the exterior square $\wedge^2 A$ is defined as the quotient of $A \otimes_{\mathbb{k}} A$ by the submodule generated by the set $\{ x \otimes x \; \mid \; x \in A \}$.
\end{definition}

The following result holds \cite[Theorem 1.1]{MR2777023}. 

\begin{theorem}\label{exterior}
If $G$ is an abelian group of odd order, then
	$$\Ho_2^{\rm Q} \big(T(G) \big) \cong  \wedge^2 G.$$ 
In particular, if	$G=\mathbb{Z}_k  ^n$, where $k$ is odd, then
	$$\Ho_2^{\rm Q} \big(T(\mathbb{Z}_k^n) \big) \cong \mathbb{Z}_k^{\frac{n(n-1)}{2}}.$$
\end{theorem}

Next, we consider computations of second quandle homology of Core quandles of some groups \cite[Example 4.2]{MR2777023}.

\begin{example}
There are five non-isomorphic groups of order $27$. The second quandle homology of their core quandles are as follows:
\begin{enumerate}
\item If $G=\mathbb{Z}_{27}$, then  $\Ho^{\rm Q}_2 \big(T(G) \big) \cong 1$.
\item If $G=\mathbb{Z}_{3}\oplus \mathbb{Z}_9$, then  $\Ho^{\rm Q} _2 \big(T(G) \big) \cong \mathbb{Z}_3$.
\item If  $G=\mathbb{Z}_{3}^3$, then $\Ho^{\rm Q} _2 \big(T(G)\big) \cong \mathbb{Z}_3^3$.
\item If $G=\langle s,t\ \mid \  s^9=t^3=1,~st=ts^4 \rangle $, then $\Ho^{\rm Q} _2 \big(\Core(G) \big) \cong\mathbb{Z}_3$.
\item If $G=\langle x,y,z\ \mid \  x^3 = y^3 = z^3 = 1, ~yz = zyx,~ xy = yx, ~xz = zx \rangle $, then $\Ho^{\rm Q}_2 \big(\Core(G) \big) \cong \mathbb{Z}_3^3$.	
\end{enumerate}
\end{example}

Mochizuki proved the following result \cite[Theorem 3.1(ii)]{MR1960136}.

\begin{theorem}\label{Mochizuki theorem on third homology}
If $p$ is an odd prime and $\R_p$ is the dihedral quandle of order $p$, then $\Ho_3^{\rm Q}(\R_p;\mathbb{ Z}_p ) \cong \mathbb{ Z}_p$.
\end{theorem}

Note that the preceding result is equivalent to $ \Ho^{\rm R}_3(\R_p) \cong\mathbb{ Z} \oplus \mathbb{ Z}_p$, which was originally conjectured by Fenn, Rourke, and Sanderson  \cite[Conjecture 5.12 ]{MR2065029}.  Furthermore, they proposed  the following general conjecture  \cite[Conjecture 16]{MR2494367}.

\begin{conjecture}
If $p$ is an odd prime and  $\R_p$ is the dihedral quandle of order $p$, then the torsion subgroup of $\Ho^{\rm R}_n(\R_p)$ is annihilated by $p$ for each $n>1$.
\end{conjecture}

The following results are proved using some homology operations on homology of quandles \cite[p. 1537]{MR2673749}.

\begin{proposition}
Let $\R_n$ be the dihedral quandle of order $n$. Then the following assertions hold:
\begin{enumerate}
\item $\Ho_3^{\rm Q}(\R_8) \cong \mathbb{Z}^2 \oplus \mathbb{Z}_8^2$.
\item $\Ho_4^{\rm Q}(\R_8) \cong \mathbb{Z}^2 \oplus \mathbb{Z}_2^4 \oplus \mathbb{Z}_4^4 \oplus \mathbb{Z}_8^2$.
\item $\Ho^{\rm Q}_3(\R_{16}) \cong \mathbb{Z}^2 \oplus \mathbb{Z}_{16}^2$.
\end{enumerate}
\end{proposition}

In general, the following statement was conjectured \cite[Conjecture 14]{MR2673749}.

\begin{conjecture}\label{Conjecture 3.1}
If $k \neq 2^t$ for any $t>1$, then $k$ annihilates the torsion in $\Ho_n^{\rm R}(\R_{2k})$. If $k=2^t$ for some $t>1$, then $2k$ is the smallest integer annihilating torsion in $ \Ho_n^{\rm R}(\R_{2k})$. 
\end{conjecture}

It is a classical result that the reduced homology of a finite group $G$ is annihilated by its order \cite[Corollary 10.2]{MR0672956}. Przytycki and Yang
proved the following analogous result for latin quandles \cite[Theorem 2.1]{MR3346517}, which was conjectured by Niebrzydowski and Przytycki  \cite[Conjecture 42]{MR2673749}.

\begin{theorem} \label{1}
Let $X$ be a finite latin quandle. Then the torsion subgroup of $\Ho_n^{\rm R}(X)$ is annihilated by $|X|$.
\end{theorem}

Note that every Latin quandle is connected; however, a connected quandle need not be Latin. It is observed in \cite{MR3346517} that the previous result does not hold for connected quandles that are not Latin.

\begin{example}\label{connected but not latin homology}
	Let $X=\{1,2,3,4,5,6\}$ be the quandle given by the following multiplication table. 
$$
\begin{tabular}{|c||c|c|c|c|c|c|}
\hline
$\ast$ & 1 & 2 & 3 & 4 & 5 & 6 \\
\hline \hline
1 & 1 & 1 & 6 & 5 & 3 & 4 \\ \hline
2 & 2 & 2 & 5 & 6 & 4 & 3 \\ \hline
3 & 5 & 6 & 3 & 3 & 2 & 1 \\ \hline
4 & 6 & 5 & 4 & 4 & 1 & 2 \\ \hline
5 & 4 & 3 & 1 & 2 & 5 & 5 \\ \hline
6 & 3 & 4 & 2 & 1 & 6 & 6	\\ \hline
\end{tabular}
$$	
Then $X$ is isomorphic to the subquandle of $\Conj(\Sigma_4)$ consisting of all 4-cycles. It has been proved in \cite[Example 7.5]{MR1825963} that $\Ho^{\rm Q}_3(X; \mathbb{Z}) \cong \mathbb{ Z}_{24}$.
\end{example}
\bigskip
\bigskip


\section{Homology of almost latin quandles}

In \cite{MR3571388}, the concept of a latin quandle was extended to that of an $m$-almost latin quandle, employing the terminology `quasigroup' in place of `latin' \cite[Definition 2.1]{MR3571388}.

\begin{definition}
Let $m\geq 1$ be an integer.  A quandle $X$ is said to be \index{$m$-almost latin quandle}{\it $m$-almost latin} if it satisfies the following conditions:
\begin{enumerate}
\item The stabilizer $\Stab(a) =\{ x \in X | \; a*x=a \}$ has order $m$  for each $a \in X$.
\item The equation $a*x=b$ has a unique solution for each $b \in X \setminus \Stab(a)$.
	\end{enumerate}
\end{definition}

If $m=1$, then $X$ is a latin quandle. Unlike latin quandles, $m$-almost latin quandles are not connected in general.

\begin{example}
{\rm The following are examples of $m$-almost latin quandles:
\begin{enumerate}
\item Every finite trivial quandle $X$ is an $|X|$-almost latin quandle.
\item The quandle $X$ of Example \ref{connected but not latin homology} is a $2$-almost latin quandle.
\item The subquandle of $\Conj(\Sigma_n)$ consisting of all $2$-cycles is a $\big(\binom{n-2}{2}+1\big)$-almost latin quandle if $n \geq 2$.
\item Let $(X,*)$ be a finite latin quandle and $\T_n$ the trivial quandle of order $n$. Consider the dynamical quandle 2-cocycle $\alpha : \T_n \times \T_n \to \Map(X\times X, X)$ defined by $\alpha_{a,b}(s,t)= s*t$ if $a=b$ and $\alpha_{a,b}(s,t)= s$ if $a\neq b$. Then the dynamical extension $X\times_{\alpha}\T_{n}$ of $\T_n$ by $X$ through $\alpha$ is a $\big((n-1)|X|+1\big)$-almost latin, but not connected for $n \geq 2$.  \cite[Example 2.5]{MR3571388}
\end{enumerate} }
\end{example}

A special case of Example (4) above is the quandle $\R_{3}\times_{\alpha} \T_{2}$ (see the table below for its multiplication table).
$$
\begin{tabular}{|c||c|c|c|c|c|c|}
\hline $\ast$ & 1 & 2 & 3 & 4 & 5 & 6 \\
\hline \hline
		1 & 1 & 3 & 2 & 1 & 1 & 1 \\ \hline
		2 & 3 & 2 & 1 & 2 & 2 & 2 \\ \hline
		3 & 2 & 1 & 3 & 3 & 3 & 3 \\ \hline
		4 & 4 & 4 & 4 & 4 & 6 & 5 \\ \hline
		5 & 5 & 5 & 5 & 6 & 5 & 4 \\ \hline
		6 & 6 & 6 & 6 & 5 & 4 & 6 \\ \hline
\end{tabular}
$$

The following homology computations are known for this quandle \cite[Example 2.5]{MR3571388}.

\begin{proposition}
Let $\R_{3}\times_{\alpha} \T_{2}$ be the quandle as above. Then the following assertions hold: 
\begin{enumerate}
\item $ \Ho_2^{\rm R}(\R_{3}\times_{\alpha} \T_{2}) \cong \mathbb{Z}^{4},\quad$ $\Ho_3^{\rm R}(\R_{3}\times_{\alpha} \T_{2}) \cong\mathbb{Z}^{8}\oplus\mathbb{Z}_{3}^{2}\quad$ and $\quad \Ho_4^{\rm R}(\R_{3}\times_{\alpha} \T_{2}) \cong\mathbb{Z}^{16}\oplus\mathbb{Z}_{3}^{8}$.
\item  $\Ho_2^{\rm Q}(\R_{3}\times_{\alpha} \T_{2})  \cong\mathbb{Z}^{2}, \quad$ $\Ho_3^{\rm Q}(\R_{3}\times_{\alpha} \T_{2}) \cong\mathbb{Z}^{2}\oplus\mathbb{Z}_{3}^{2}\quad$ and $\quad\Ho_4^{\rm Q}(\R_{3}\times_{\alpha} \T_{2}) \cong\mathbb{Z}^{2}\oplus\mathbb{Z}_{3}^{6}$.
\end{enumerate}
\end{proposition} 
\para

Using homotopy of certain chain maps, the following result was established in \cite[Theorem 3.2]{MR3571388}.

\begin{theorem}  
Let $X$ be a finite $m$-almost latin quandle. Suppose that for each $x \in X$, the stabilizer $\Stab(x)$ is a trivial subquandle of $X$. Then the torsion subgroup of $\Ho_n^{\rm R}(X)$ is annihilated by $m\big(\lcm(|X|,|X|-m) \big)$.
\end{theorem}

Let $X$ be a quandle and $\Inn(X)$  be its inner automorphism group.  Then the map $S:X \rightarrow \Conj \big(\Inn(X) \big)$ is a quandle homomorphism.  In \cite{MR4407086}, the sequence of homomorphisms $$\rho_n: C_n^{\rm R}(X) \rightarrow C_n^{\rm R} \big(\Conj(\Inn(X)) \big)$$ defined by $\rho_n(x_1, \ldots, x_n)= (S_{x_1}, \ldots , S_{x_n})$ is considered, which constitute a chain map. This map is then used to establish the following result \cite[Theorem 2.1]{MR4407086}.

\begin{theorem}
Let $X$ be a finite connected quandle.  Then, for each $n >1$,  the image of the torsion subgroup of $\Ho_n^{\rm R}(X)$ under the induced homomorphism $\rho_n:  \Ho_n^{\rm R}(X) \rightarrow \Ho_n^{\rm R} \big(\Conj(\Inn(X))\big)$ is annihilated by $| \Inn(X) |$.
\end{theorem}   

We can define the reduced homology of quandles in a manner similar to the reduced homology of groups (see \cite[p. 211]{MR0672956}).

\begin{definition}
The reduced homology of a quandle $X$ is obtained from the augmented chain complex 
$$ \cdots \rightarrow C_3(X) \rightarrow C_2(X) \rightarrow C_1(X) \xrightarrow{\partial_1} \mathbb{Z},$$ where $\partial_1(x)=1$ for all $x \in X$.
\end{definition}

With the preceding definition, we have the following result \cite[Corollary 2.3]{MR4407086}. 

\begin{corollary}
Let $X$ be a finite connected quandle.  For each $n> 1$, if the induced homomorphism $\rho_n :  \Ho_n^{\rm R}(X) \rightarrow \Ho_n^{\rm R} \big(\Conj(\Inn(X))\big)$ is injective, then both the reduced quandle homology $\tilde{\Ho}^{\rm Q}_n(X)$ and the torsion subgroup of $\Ho_n^{\rm R}(X)$ are annihilated by $| \Inn(X) |$.
\end{corollary}

\begin{remark}
{\rm 
Table \ref{table of homology of conjugation quandles} demonstrates computations of integral quandle homology for some connected quandles that are not Latin. The table is borrowed from \cite[Example 2.2]{MR4407086} and the quandles are taken from \cite{rig}, where the notation $Q(n,i)$ stands for the $i$-th connected quandle of order $n$ in the RIG package. Here, $\Sl(2,\mathbb{F}_3)$ denotes the special linear group of $2 \times 2$ matrices over the field of three elements and $A_n$ denotes the alternating group on $n$ letters.}
\end{remark}
	\begin{table}[H]
		\begin{center}
			\begin{tabular}{|c||c||c||c||c||c|} \hline $\ \ \ \ \ \  \text{Quandle} \;X \ \ \ \ \ \  $& $|X|$  & $\Inn(X)$& $\Ho^{\rm Q}_2(X)$ & $\Ho^{\rm Q}_2 \big(\Conj(\Inn(X))\big)$ & $\Ho^{\rm Q}_3(X)$ \\ \hline \hline
				$Q(6,1) $ & 6 &$ \Sigma_4$ & $\mathbb{Z}_2$ & $\mathbb{Z}^{20} \oplus \mathbb{Z}_2 \oplus \mathbb{Z}_2 \oplus \mathbb{Z}_6 $& $\mathbb{Z}_6$  \\ \hline
				$Q(6,2) $& 6& $ \Sigma_4$ & $\mathbb{Z}_4$ & $\mathbb{Z}^{20} \oplus \mathbb{Z}_2 \oplus \mathbb{Z}_2 \oplus \mathbb{Z}_6 $& $\mathbb{Z}_{24}$   \\ \hline
				$ Q(8,1) $& 8& $ \Sl(2,\mathbb{F}_3)$ & $1$ & $\mathbb{Z}^{42} \oplus \mathbb{Z}_2 \oplus \mathbb{Z}_2  \oplus \mathbb{Z}_2 \oplus \mathbb{Z}_2 \oplus \mathbb{Z}_4 $& $\mathbb{Z}_8$  \\ \hline
				$Q(10,1)$& 10& $ \Sigma_5$ & $\mathbb{Z}_2$ & $?$ & $\mathbb{Z}_2$ \\ \hline
				$Q(12,1)$ & 12 & $ \Sigma_4$ & $\mathbb{Z}_2$ & $\mathbb{Z}^{20} \oplus \mathbb{Z}_2 \oplus \mathbb{Z}_2 \oplus \mathbb{Z}_6 $& $\mathbb{Z}_{24}$  \\\hline
				$Q(12,2)$ & 12& $ A_4 \rtimes \mathbb{Z}_4 $ & $\mathbb{Z}_2$&? &$\mathbb{Z}_{24}$  \\\hline
				$Q(12,3)$ & 12 &$A_5$ & $\mathbb{Z}_{10}$ & ?& $\mathbb{Z}_2 \oplus \mathbb{Z}_{60} $  \\ \hline
\end{tabular}
\caption{Integral quandle homology of small connected quandles that are not latin.}\label{table of homology of conjugation quandles}
\end{center}
\end{table}	
\bigskip
\bigskip


\section{Homology of Alexander quandles}

Recall that, an Alexander quandle is a $\mathbb{Z}[t, t^{-1}]$-module $X$ with the quandle operation $x* y = tx + (1- t)y$ for all $x, y \in X$.  Since $X$ is a $\mathbb{Z}[t, t^{-1}]$-module, it is a $\mathbb{Z}$-module, and hence we can consider its tensor product $X \otimes_{\mathbb{Z}} X$ and the exterior product $X\wedge_{\mathbb{Z}} X$. In \cite{MR4065999}, a connection was established between the second quandle homology of quandles and the Schur multipliers of suitable associated groups. Further, a description of the second quandle homology of some Alexander quandles was also given \cite[Theorem 4.1]{MR4065999}\label{prob1}. 

\begin{theorem} 
Let $X$ be an Alexander quandle such that $1-t$ invertible. Then the map $\mathbb{Z} ^{ X \times X}  \to X \otimes_{\mathbb{Z}} X$, given by $(x,y) \mapsto x \otimes (1-t)y$, gives an isomorphism 
\begin{equation*}\label{abc2} \Ho_2^{\rm Q}(X;\mathbb{Z}) \cong
X\otimes_{\mathbb{Z}} X / \langle x \otimes y - y \otimes t x \mid x,y \in X \rangle.
\end{equation*}	
\end{theorem}

Note that, the map $ X \otimes_{\mathbb{Z}} X \rightarrow X \otimes_{\mathbb{Z}} X$ given by  $x \otimes y \mapsto x \otimes (1-t)y$ induces an isomorphism
$$
X\otimes_{\mathbb{Z}} X / \langle x \otimes y - y \otimes t x \mid x,y \in X \rangle \longrightarrow X\wedge_{\mathbb{Z}} X /
  \langle x \wedge y - tx \wedge ty \mid x,y \in X \rangle.
$$
In fact, its inverse is induced by the map $ X \otimes_{\mathbb{Z}} X \rightarrow X \otimes_{\mathbb{Z}} X$ given by $x \otimes y \mapsto x \otimes (1-t)^{-1}y$.  This leads to the following result \cite[Corollary 4.3] {MR4065999}.

\begin{corollary} \label{prob7}
Let $X$ be an Alexander quandle such that $1-t$ invertible. Then there is an isomorphism
	\begin{equation*}\label{abc} \Ho_2^{\rm Q}(X;\mathbb{Z}) \cong X\wedge_{\mathbb{Z}} X / \langle x \wedge y - tx \wedge ty  \mid x,y \in X \rangle. 
	\end{equation*}
\end{corollary}

The preceding result is used to compute the second homology of the Alexander quandle $\mathbb{F}_p [t]/(1+t+ \dots +t^{n-1})$, where $p$ is a prime  \cite[Proposition 4.6] {MR4065999}.

\begin{proposition}\label{prob5} 
Let $p$ be a prime and $n > 2$ an integer such that $\gcd(n, p)=1$.	Let $X$ be the Alexander quandle $\mathbb{F}_p [t]/(1+t+ \dots +t^{n-1})$. Then there is an isomorphism
	\begin{equation*}\label{abc34} \Ho_2^{\rm Q}(X;\mathbb{Z}) \cong (\mathbb{Z}_p)^{ \lfloor \frac{n-1}{2}\rfloor } .
	\end{equation*}
\end{proposition}

The next result gives computations of some higher homology groups \cite[Example 13(ii)]{MR2673749}.

\begin{proposition}
Let $X = \mathbb{Z}_2[t]/(1+t+t^2)$ be the four element Alexander quandle. Then the following assertions hold:
$$ \Ho^{\rm Q}_2(X) \cong \mathbb{Z}_2, \;\; \quad \quad \quad 	\Ho^{\rm Q}_3(X) \cong \mathbb{Z}_2\oplus \mathbb{Z}_4, \quad \Ho^{\rm Q}_4(X) \cong \mathbb{Z}_2^{2}\oplus \mathbb{Z}_4,$$
$$\Ho^{\rm Q}_5(X) \cong \mathbb{Z}_2^{5}\oplus \mathbb{Z}_4, \quad \Ho^{\rm Q}_6(X) \cong \mathbb{Z}_2^{9}\oplus \mathbb{Z}_4^{2}, \quad 	\Ho^{\rm Q}_7(X) \cong \mathbb{Z}_2^{17}\oplus \mathbb{Z}_4^{3}.
$$
\end{proposition}

The preceding computations lead to the following conjecture \cite[Example 13(ii)]{MR2673749}.

\begin{conjecture}
The torsion subgroup of $\Ho^{\rm Q}_n \big(\mathbb{Z}_2[t]/(1+t+t^2) \big)$ is isomorphic  to $\mathbb{Z}_4^{f_n} \oplus \mathbb{Z}_2^{f'_n}$, where 
$\{f_n\}_{n \ge 1}$  is the sequence of delayed Fibonacci numbers and 
$$f'_n = \log_2 \big(| \Tor \Ho_n^{\rm Q}(\mathbb{Z}_2[t]/(1+t+t^2))| \big) -2f_n.$$ 
\end{conjecture}

It is well-known \cite{MR1848966} that any  connected quandle of order a prime $p$ is isomorphic to one of the Alexander quandles defined on $X=\mathbb{Z}_p$ with the quandle operation given by $$x*y=\omega x+(1-\omega)y,$$ where $\omega \neq 0,1$.  Nosaka computed the integral quandle homology of Alexander quandles of prime orders \cite[Theorem 2.1]{MR3042590}.

\begin{theorem}\label{Nosaka result on alexander quandles} 
Let $X=\mathbb{ Z}_p$ be the Alexander quandle of order $p$ with quandle operation $x*y=\omega x+(1-\omega)y,$ where $\omega \neq 0,1$.  Then
	$$\Ho^{\rm Q}_1(X) \cong \mathbb{ Z} \oplus \mathbb{ Z}_p^{b_1} \quad  \textrm{and} \quad \Ho^{\rm Q}_n(X) \cong  \mathbb{ Z}_p^{b_n}\quad  \textrm{for each} \quad n \geq 2,$$  where $b_{n+2e}=b_n+b_{n+1} + b_{n+2}$,  $b_1=b_2=\cdots =b_{2e-2}=0$,  $b_{2e-1}=b_{2e}=1$ and $e>0$ is the minimal integer satisfying $\omega^e=1$.  
\end{theorem} 

As an interesting consequence of Theorem \ref{Nosaka result on alexander quandles}, we have the following result, which was originally conjectured by Niebrzydowski and Przytycki  \cite[Conjecture 5]{MR2494367}.

\begin{corollary}
If $p$ is an odd prime and $\R_p$ is the dihedral quandle of order $p$. Then the torsion subgroup of $\Ho^{\rm Q}_n(\R_p)$ is isomorphic to $ \mathbb{Z}_p^{f_n}$, where $\{f_n\}_{n \ge 1}$ is the sequence of delayed Fibonacci numbers given by $$f_1=f_2=0, ~f_3=1 \quad \textrm{and} \quad f_n = f_{n-1} + f_{n-3}~\textrm{for}~ n \ge 4.$$
\end{corollary}

As another consequence of Theorem \ref{Nosaka result on alexander quandles}, we have  the following result \cite[Corollary 2.3]{MR3042590}.

\begin{corollary} 
Let $X=\mathbb{ Z}_p$ be the Alexander quandle of order $p$ with quandle operation $x*y=\omega x+(1-\omega)y,$ where $\omega \neq -1,0,1$. Then $\Ho^{\rm Q}_n(X)$ is trivial for $n=2,3,4$.
\end{corollary}
\para

In \cite{MR2988494},  Kabaya presented a construction of quandle cocycles derived from group cocycles. For an integer $p \geq 3$, he obtained quandle cocycles of the dihedral quandle $\R_p$ from group cocycles of the cyclic group $\mathbb{Z}_p$. In particular, it was shown  that a group 3-cocycle of $\mathbb{Z}_p$ gives rise to a non-trivial quandle 3-cocycle of $\R_p$.  By Theorem \ref{Mochizuki theorem on third homology}, for an odd prime $p$, the cohomology group $\Ho_{\rm Q}^3(\R_p; \mathbb{ Z}_p)$ is one-dimensional. Hence, this quandle $3$-cocycle is a constant multiple of the Mochizuki 3-cocycle up to a coboundary. Below we review Kabaya's construction.
\para

 For a group $G$ and a $G$-module $A$, we denote the homology and cohomology of $G$ with coefficients in $A$ by $\Ho_n^{\rm{Grp}} (G;A)$ and $\Ho^n_{\rm{Grp}}(G;A)$, respectively. Let $G$ be the cyclic group $\mathbb{Z}_p$, where $p >2$ is an integer not necessarily a prime. The first cohomology group $\Ho^1_{\rm{Grp}}(G; \mathbb{Z}_p) = \Hom_{\mathcal{G}} (\mathbb{Z}_p, \mathbb{Z}_p)$ is generated by the group 1-cocycle $b_1$ given by
$$
b_1(x) = x.
$$
The connecting homomorphism $\delta \colon \Ho^1_{\rm{Grp}}(G; \mathbb{Z}_p) \to \Ho^2_{\rm{Grp}}(G; \mathbb{Z})$ of the long exact sequence of cohomology groups corresponding to the short exact sequence $$1 \to \mathbb{Z} \to \mathbb{Z} \to \mathbb{Z}_p \to 1$$ maps $b_1$ to a generator of $\Ho^2_{\rm{Grp}}(G; \mathbb{Z})$. Further,  the homomorphism $\Ho^2_{\rm{Grp}}(G; \mathbb{Z}) \to  \Ho^2_{\rm{Grp}}(G; \mathbb{Z}_p)$ induced by the natural projection $\mathbb{Z} \to \mathbb{Z}_p $,  maps it to a generator $b_2$ of $\Ho^2_{\rm{Grp}}(G; \mathbb{Z}_p)$. More explicitly, we have
\begin{equation*}
	\label{eq:2_cocycle_of_cyclic_group}
	b_2(x,y) = \frac{1}{p} \big(\overline{x}+\overline{y} - \overline{x+y} \big)
	=
	\left\{ 
	\begin{array}{ll}
		1 & \textrm{if $\overline{x}+\overline{y} \geq p$,} \\  
		0 & \textrm{otherwise,}
	\end{array} \right. 
\end{equation*}
where $\overline{x}$ is an integer $0 \leq \overline{x} < p$ such that $\overline{x} \equiv x \mod p$. When $p$ is an odd prime, it is known that  any element of $\Ho^m_{\rm{Grp}}(G; \mathbb{Z}_p)$ can be presented by a cup product of $b_1$'s and $b_2$'s \cite[Proposition 3.5.5]{MR1110581}. 
\para

For an integer $p > 2$, let $f$ be a normalized group $k$-cocycle representing an element of $\Ho^k_{\rm{Grp}}(G; \mathbb{Z}_p)$. Viewing $\R_p$ as $\mathbb{Z}_p$, 
we obtain a map $f \colon (\R_p)^{k+1} \to \mathbb{Z}_p$ satisfying the following conditions:
\begin{enumerate}
	\item[(1)] $\displaystyle\sum_{i=0}^{k+1} (-1)^i f(x_0, \dots, \widehat{x_i}, \dots ,x_{k+1}) = 0$.
	\item[(2)] $f(x_0, \dots , x_{k} ) = 0$ if $x_i=x_{i+1}$ for some $i$.
\end{enumerate}
If $f$ also satisfies the condition 
\begin{enumerate}
	\item[(3)] $f(x_0 * y, \dots, x_k * y) =  f(x_0, \dots, x_k)$ for any $y \in \R_p$,
\end{enumerate}
then $f$ gives rise to a quandle cohomology class in $\Ho_{\rm Q}^k(\R_p; \mathbb{Z}_p)$.  Define $\tilde{f}: (\R_p)^{k+1} \to \mathbb{Z}_p$ by
$$	\tilde{f}(x_0, \dots, x_k) = f(x_0, \dots, x_k) + f(-x_0,\dots ,-x_k).
$$
Then the map $\tilde{f}$ satisfies conditions (1), (2) and (3), and hence is a quandle $k$-cocycle. We use this construction to give an explicit presentation of the quandle $3$-cocycle arising from $b_1 b_2 \in \Ho^3_{\rm{Grp}}(G; \mathbb{Z}_p)$. We define
$$d(x,y)= b_2(x,y) - b_2(-x,-y).  $$
By the preceding construction, $d$ is a quandle 2-cocycle. By definition, we have
$$	d(-x,-y) = -d(x,y)$$
and 
$$
	d(x,y) =  \left\{ 
	\begin{array}{ll}
		1  & ~\textrm{if}~~ \overline{x} + \overline{y} > p, \\  
		-1 & ~\textrm{if}~~ \overline{x} + \overline{y} < p,~x \neq 0~ \textrm{and}~y \neq 0, \\
		0  & ~\textrm{otherwise}. \\
	\end{array} \right. 
$$

The preceding discussion leads to the following result \cite[Proposition 7.1]{MR2988494}.

\begin{proposition}
The quandle $3$-cocycle arising from $b_1 b_2 \in \Ho_{\rm{Grp}}^3(G; \mathbb{Z}_p)$ has the form
	\begin{equation}
		\label{eq:the_cocycle}
		(x,y,z) \mapsto 2z \big(d(y-x,z-y) + d(y-x,y-z) \big)
	\end{equation}
for all $x,y,z \in \R_p$. Further, it gives a non-trivial cohomology class in  $\Ho^3_{\rm Q}(\R_p; \mathbb{Z}_p)$.

\end{proposition}

Since 2 is divisible in $\mathbb{Z}_p$ when $p$ is odd, we have the following result \cite[Corollary 7.3]{MR2988494}.

\begin{corollary}
If $p > 2$ is an odd integer, then the map $\R_p \times \R_p \times \R_p \to \mathbb{Z}_p$, defined by
	\begin{equation*}
		(x,y,z) \mapsto z \big(d(y-x,z-y) + d(y-x,y-z) \big),
	\end{equation*}
gives a non-trivial cohomology class in  $\Ho^3_{\rm Q}(\R_p; \mathbb{Z}_p)$.
\end{corollary}

We conclude this section with Table \ref{table of homology of quandles}, which lists second and third integral quandle homology groups of some Alexander quandles \cite[Table 1]{MR3054333}.

\begin{table}[H]
\begin{center}
\begin{tabular}{|c||c||c|c|} \hline $\ \ \ \ \ \ \ \ \ \text{Quandle} \;X \ \ \ \ \ \ \ \ \ $ & $\Ho_2^{\rm Q}(X)$ & $\Ho_3^{\rm Q}(X)$ \\ \hline \hline
				$\R_3 $& 1& $\mathbb{Z}_3$ \\ \hline
				$\R_5 $& 1& $\mathbb{Z}_5$  \\ \hline
				$\mathbb{Z}_5[T]/(T-\omega)$& 1 & $1$ \\ \hline
				$ \R_7 $& 1& $\mathbb{Z}_7$  \\ \hline
				$\mathbb{Z}_7[T]/(T-\omega)$& 1 & $1$  \\ \hline
				$\mathbb{Z}_2[T]/(T^3+T^2+1)$& 1 & $\mathbb{Z}_2$ \\\hline
				$\mathbb{Z}_2[T]/(T^3+T+1)$&1 & $\mathbb{Z}_2$ \\\hline
				$ \R_9 $& 1& $\mathbb{Z}_9$  \\ \hline
				$\mathbb{Z}_3[T]/(T^2+1)$& $\mathbb{Z}_3$ & $(\mathbb{Z}_3)^3$ \\ \hline
				$\mathbb{Z}_3[T]/(T^2+T-1)$& 1 & 1\\ \hline
				$\mathbb{Z}_3[T]/(T^2-T-1)$&1 & 1 \\
				\hline
\end{tabular}
\caption{Second and third integral quandle homology groups of some Alexander quandles.}\label{table of homology of quandles}
\end{center}
\end{table}	
\bigskip
\bigskip


\section{Homology of permutation racks}\label{section homology permutation racks}

In \cite{MR3937311}, Szymik provided an interpretation of rack homology using tools from Quillen's homotopical algebra.  Employing this approach, in  \cite{arXiv:2011.04524}, Lebed and Szymik computed the entire homology of permutation racks with $\mathbb{Z}$-coefficients. They constructed a spectral sequence that reduced the problem to dealing with certain equivariant homology.  

\begin{definition}
A permutation rack is a non-empty set $X$ with a binary operation defined by $x*y=f(x)$, where $f$ is a fixed permutation of $X$.
\end{definition}

We denote such a rack by $(X,f)$.  The orbit of an element $x \in X$ is given by $\{f^n(x) \mid n \in \mathbb{Z} \}$.  A permutation rack $(X,f)$ can be thought of as a set $X$ with an action of the infinite cyclic group $\mathbb{Z}$, where the integer $n$ acts by $f^n$.  Further, a permutation rack $(X,f)$ is called \index{free permutation rack}{\it free} if the corresponding action is free. Let $\mathcal{O}$ be the set of orbits of the permutation rack $(X, f)$. Let  $\mathbb{Z}^{\mathcal{O}}$ denote the free abelian group on the set $\mathcal{O}$ and let
 $$\mathbb{Z}^{\mathcal{O}}_0=\Big\{ \sum_{x_i \in \mathcal{O}}n_i x_i\;|\; \sum n_i=0 \Big\}$$ be the free abelian subgroup of $\mathbb{Z}^{\mathcal{O}}$ consisting of elements whose coefficients add up to zero.  Then the following result is known \cite[Theorem B]{arXiv:2011.04524}.

 \begin{theorem}
If $(X,f)$ be a free permutation rack with set of orbits $\mathcal{O}$, then
$$	\Ho_n^{\rm R}(X) \cong \underbrace{(\mathbb{Z}^{\mathcal{O}}_0 \otimes \cdots \otimes \mathbb{Z}^{\mathcal{O}}_0)}_{(n-1)~{\rm times}}\otimes  \mathbb{Z}^{\mathcal{O}}.$$
In particular, the homology is a free abelian group in each degree. If $|\mathcal{O}| < \infty$, then
the $n$-th Betti number of $(X, f)$ is $|\mathcal{O}|\big(|\mathcal{O}| -1\big)^{n-1}$.
 \end{theorem}

A $G$-set is a non-empty set $X$ equipped with an action of a discrete group $G$. It is well-known that there exists a contractible space $EG$ on which $G$ acts freely, and whose quotient $EG/G$ is the classifying space of the group $G$.  We recall a well-known construction in equivariant topology due to Borel \cite{MR0116341}. Consider an equivalence relation on $EG \times X$ generated by $(g \cdot x,y) \sim (x,g^{-1} \cdot y)$ for any $x \in EG$, $y \in X$ and $g \in G$.  The quotient space $$X//G= (EG \times X)/\sim$$ is called the \index{Borel space}{\it Borel space} of the $G$-set $X$.  The singular homology $\Ho_n(X//G)$ of the space $X//G$ is called the \index{$G$-equivariant homology}{\it $G$-equivariant homology} of $X$ and is denoted by $\Ho_n^G(X)$.
\para 

  Now, let $(X,f)$ be a permutation rack (viewed as a $\mathbb{Z}$-set) and let  $X//f$ be the Borel space of the $\mathbb{Z}$-set $X$.  Then the equivariant  homology $ \Ho_n^{\mathbb{Z}}(X)$ of the $\mathbb{ Z}$-set $X$ is given by the following result \cite[Proposition 3.2]{arXiv:2011.04524}.
  
  \begin{proposition}
  	Let $(X,f)$ be a permutation rack.  Let $\mathcal{O}$ be the set of all orbits of $(X,f)$ and $\mathcal{O}_{\rm{finite}}$ the subset of all finite orbits.  Then the equivariant homology of $(X,f)$ is given by
$$
  	\Ho_n^{\mathbb{Z}}(X)\cong \begin{cases}
  		\mathbb{ Z}^\mathcal{O} & \text{if}  \;n=0,\\
  		\mathbb{ Z}^{\mathcal{O} _{\rm{finite}}} & \text{if} \; n=1,\\
1 & \text{otherwise}.
  	\end{cases}
  $$
  \end{proposition}

The next result shows the existence of a spectral sequence converging to the rack homology of a permutation rack \cite[Theorem A]{arXiv:2011.04524}.

\begin{theorem}\label{main theorem Lebed and Szymik}
Let $(X,f)$ be a permutation rack.   Then there is a homological spectral sequence whose second page $E^2_{*, *}$ is given by
$$	  E^2_{p,q} \cong \tilde{\Ho}_p(X//f)^{\otimes (q-1)} \otimes \Ho_p(X//f)$$
and  which converges to the rack homology $\Ho^{\rm R}_{p+q}(X)$,  where $\tilde{\Ho}_p(X//f)$ stands for the reduced singular homology.  Furthermore, this spectral sequence always degenerates from its $E^2_{*, *}$ page onwards.
\end{theorem} 
 
In general, the spectral sequence in Theorem \ref{main theorem Lebed and Szymik} is concentrated in the region of the first quadrant where $p\le q$, which implies the vanishing of the differentials in low dimensions. This yields the following main result \cite[Theorem C]{arXiv:2011.04524}.

\begin{theorem}
Let $(X,f)$ be a permutation rack.  Let $\mathcal{O}$ be the set of all orbits of $(X,f)$ and $\mathcal{O}_{\rm{finite}}$ the subset of all finite orbits. Then the following assertions hold:
\begin{enumerate}
\item $\Ho^{\rm R}_0(X; \mathbb{Z}) \cong \mathbb{Z}$.
\item $\Ho^{\rm R}_1(X; \mathbb{Z}) \cong  \mathbb{Z}^\mathcal{O} \cong \mathbb{Z}^{|\mathcal{O}|}$.
\item  $\Ho^{\rm R}_2(X; \mathbb{Z}) \cong  \big(\mathbb{Z}^{\mathcal{O}}_0 \otimes \mathbb{Z}^\mathcal{O} \big) \oplus \mathbb{Z}^{\mathcal{O}_{\rm finite}} \cong \mathbb{Z}^{|\mathcal{O}| \big(|\mathcal{O}|-1 \big) + |\mathcal{O}_{\rm finite}|}$.
\end{enumerate}
\end{theorem} 

The following result can be deduced by establishing the degeneracy of the spectral sequence \cite[Corollary 6.2]{arXiv:2011.04524}.
 
\begin{corollary}
Let $(X,f)$ be a permutation rack. Then all its homology groups $\Ho^{\rm R}_n(X; \mathbb{Z})$ are free abelian.
\end{corollary}


\chapter{Automorphism groups, second cohomology and extensions of quandles}\label{chap cohomology of quandles}

\begin{quote}
This chapter introduces a precise four-term exact sequence that links quandle 1-cocycles, second quandle cohomology, and a particular group of automorphisms arising from an abelian extension of quandles. We examine quandle extensions that originate from group extensions. Additionally, we develop natural group homomorphisms connecting the second cohomology of a group to that of its associated core and conjugation quandles, revealing deeper structural relationships among these objects.
\end{quote}
\bigskip

\section{Action of automorphism groups on dynamical 2-cocycles}\label{Action of automorphisms on dynamical 2-cocycles}

In \cite{MR1994219}, Andruskiewitsch and Gra\~{n}a considered a special case of Proposition \ref{set-cocycle} when $X$ is a quandle and $\beta_{s, t}(x, y)=x*y$ for all $x, y \in X$ and $s, t \in S$. In this case, the map $\alpha$ satisfies the
conditions
\begin{equation}\label{dynamical-cocycle-condition1}
\alpha_{x, x}(s, s)=s,
\end{equation}
\begin{equation}\label{dynamical-cocycle-condition2}
\alpha_{x, y}(-, t): S \to S~\textrm{is a bijection}
\end{equation}
and the  {\it cocycle condition}
\begin{equation}\label{dynamical-cocycle-condition3}
\alpha_{x*y, z}\big(\alpha_{x, y}(s, t), ~u \big)= \alpha_{x* z, y*z}\big(\alpha_{x, z}(s, u),\, \alpha_{y, z}(t, u) \big)
\end{equation}
for all $x, y \in X$ and $s, t \in S$. Such an $\alpha$ is referred to as a {\it dynamical 2-cocycle}. Further, in this case, the quandle operation \eqref{genralised-quandle-operation}  on $X \times S$ becomes 
\begin{equation}\label{dynamical-quandle-operation1}
(x, s)* (y,t)= \big( x* y, \,\alpha_{x, y}(s, t) \big).
\end{equation} 
The quandle so obtained is called an  {\it extension} of $X$ by $S$ through $\alpha$, and is denoted by $X \times_{\alpha} S$. It is easy to see that, for each $x \in X$, the set $S$ forms a quandle with the binary operation $$s *_x t :=\alpha_{x,x}(s,t)$$ for all $s, t \in S$.
\para

Let $X$ be a quandle, $\Sigma_S$ the symmetric group on a set $S$ and $\mathcal{Z}^2_{\rm Q}(X; S)$ the set of all dynamical 2-cocycles. For $(\phi, \theta) \in \Aut(X) \times \Sigma_S$ and a dynamical 2-cocycle $\alpha \in \mathcal{Z}^2_{\rm Q}(X; S)$, define $$^{(\phi, \theta)}\alpha:X \times X \to \Map(S \times S, S)$$ by setting
\begin{equation}\label{auto-action}
^{(\phi, \theta)}\alpha_{x, y}(s, t):= \theta \Big(\alpha_{\phi^{-1}(x), \phi^{-1}(y)}\big(\theta^{-1}(s), ~\theta^{-1}(t)\big)  \Big)
\end{equation}
for $x, y \in X$ and $s, t \in S$.

\begin{lemma}\label{aut action on cocycles}
The map $\big(\Aut(X) \times \Sigma_S\big) \times \mathcal{Z}^2_{\rm Q}(X; S) \to \mathcal{Z}^2_{\rm Q}(X; S)$ given by
$$\big(( \phi, \theta), \alpha\big) \mapsto~ ^{(\phi, \theta)}\alpha$$ defines a left-action of the group $\Aut(X) \times \Sigma_S$ on the set $\mathcal{Z}^2_{\rm Q}(X; S)$.
\end{lemma}
\begin{proof} 
Let $(\phi, \theta) \in \Aut(X) \times \Sigma_S$ and $\alpha \in \mathcal{Z}^2_{\rm Q}(X; S)$ a dynamical 2-cocycle. We first prove that $^{(\phi, \theta)}\alpha$ is a dynamical 2-cocycle. For $x, y \in X$ and $s, t \in S$, we see that
\begin{eqnarray*}
&  & ^{(\phi, \theta)}\alpha_{x*y, z}\Big(^{(\phi, \theta)}\alpha_{x, y}(s, t), u \Big) \\
& = & \theta \Big(\alpha_{\phi^{-1}(x)*\phi^{-1}(y), ~\phi^{-1}(z)}\Big(\theta^{-1}\big(^{(\phi, \theta)}\alpha_{x, y}(s, t)\big),~ \theta^{-1}(u)\Big)  \Big)\\
& = & \theta \Big(\alpha_{\phi^{-1}(x)*\phi^{-1}(y), ~\phi^{-1}(z)}\Big(\alpha_{\phi^{-1}(x), ~\phi^{-1}(y)}\big(\theta^{-1}(s),~ \theta^{-1}(t)\big),~ \theta^{-1}(u)\Big)  \Big)\\
& = & \theta \Big( \alpha_{\phi^{-1}(x)* \phi^{-1}(z),~ \phi^{-1}(y)*\phi^{-1}(z)}\Big(\alpha_{\phi^{-1}(x),~ \phi^{-1}(z)} \big(\theta^{-1}(s), ~\theta^{-1}(u)\big), ~\alpha_{\phi^{-1}(y), ~\phi^{-1}(z)}\big(\theta^{-1}(t), ~\theta^{-1}(u)\big) \Big) \Big)\\
& & \textrm{by the cocycle condition \eqref{dynamical-cocycle-condition3} for}~\alpha \\
& = & \theta \Big( \alpha_{\phi^{-1}(x*z), \phi^{-1}(y*z)}\Big(\alpha_{\phi^{-1}(x), ~\phi^{-1}(z)}\big(\theta^{-1}(s),~ \theta^{-1}(u)\big),~ \alpha_{\phi^{-1}(y), ~\phi^{-1}(z)}\big(\theta^{-1}(t),~ \theta^{-1}(u)\big) \Big) \Big)\\
& = & ^{(\phi, \theta)}\alpha_{x* z, y*z}\Big(^{(\phi, \theta)}\alpha_{x, z}(s, u),~~ ^{(\phi, \theta)}\alpha_{y, z}(t, u) \Big),
\end{eqnarray*}
and hence $^{(\phi, \theta)}\alpha$ satisfies \eqref{dynamical-cocycle-condition3}. A direct check shows that $ ^{(\phi, \theta)}\alpha$ satisfy \eqref{dynamical-cocycle-condition1} and \eqref{dynamical-cocycle-condition2}, which shows that $ ^{(\phi, \theta)}\alpha \in \mathcal{Z}^2_{\rm Q}(X; S)$.
\para

Given $(\phi_1, \theta_1), (\phi_2, \theta_2) \in \Aut(X) \times \Sigma_S$, we see that
\begin{eqnarray*}
^{(\phi_1, \theta_1)(\phi_2, \theta_2)}\alpha_{x, y}(s, t) & = & ^{(\phi_1\phi_2, ~\theta_1\theta_2)}\alpha_{x, y}(s, t) \\
& = & \theta_1\theta_2 \Big(\alpha_{\phi_2^{-1}\phi_1^{-1}(x), ~\phi_2^{-1}\phi_1^{-1}(y)}\big(\theta_2^{-1}\theta_1^{-1}(s), ~\theta_2^{-1}\theta_1^{-1}(t)\big)  \Big)\\
& = & \theta_1 \Big( ^{(\phi_2, \theta_2)}\alpha_{\phi_1^{-1}(x), ~\phi_1^{-1}(y)}\big(\theta_1^{-1}(s),~ \theta_1^{-1}(t)\big) \Big)\\
& = & ^{(\phi_1, \theta_1)}\big(^{(\phi_2, \theta_2)}\alpha\big)_{x, y}(s, t),
\end{eqnarray*}
and hence the group $\Aut(X) \times \Sigma_S$ acts on $\mathcal{Z}^2_{\rm Q}(X; S)$ from the left. 
$\blacksquare$    \end{proof}          

Two dynamical 2-cocycles $\alpha, \beta \in \mathcal{Z}^2_{\rm Q}(X; S)$ are said to be {\it cohomologous } if there exists a map $\lambda:X \to \Sigma_S$ such that
\begin{equation}\label{cohomologous}
\beta_{x, y}(s, t)=\lambda_{x*y} \Big(\alpha_{x, y}\big(\lambda_x^{-1}(s),~ \lambda_y^{-1}(t) \big)\Big)
\end{equation}
for all $x, y \in X$ and $s, t \in S$. Clearly, being cohomologous is an equivalence relation on $\mathcal{Z}^2_{\rm Q}(X; S)$, and let $\mathcal{H}^2_{\rm Q}(X; S)$ denote the set of cohomology classes of dynamical 2-cocycles. For $\alpha \in \mathcal{Z}^2_{\rm Q}(X; S)$, we denote its equivalence class  by $[\alpha]  \in \mathcal{H}^2_{\rm Q}(X; S)$.

\begin{proposition} \cite[Proposition 4.2]{MR4282648}
The map $\big(\Aut(X) \times \Sigma_S\big) \times \mathcal{H}^2_{\rm Q}(X; S) \to \mathcal{H}^2_{\rm Q}(X; S)$ given by
$$\big(( \phi, \theta), [\alpha] \big) \mapsto~ [^{(\phi, \theta)}\alpha]$$ defines a left-action of the group $\Aut(X) \times \Sigma_S$ on the set $\mathcal{H}^2_{\rm Q}(X; S)$.
\end{proposition}

\begin{proof}
It suffices to prove that cohomologous cocycles are mapped onto cohomologous cocycles under the action in Lemma \ref{aut action on cocycles}. Let $(\phi, \theta) \in \Aut(X) \times \Sigma_S$,  $\alpha, \beta \in \mathcal{Z}^2_{\rm Q}(X; S)$ and $\lambda:X \to \Sigma_S$ such that equation \eqref{cohomologous} is satisfied. Define a map $\lambda':X \to \Sigma_S$ by setting $\lambda'_x= \theta \lambda_{\phi^{-1}(x)} \theta^{-1}$ for $x \in X$. Then, for $x, y \in X$ and $s, t \in S$, we see that
\begin{eqnarray*}
&  & ^{(\phi, \theta)}\beta_{x, y}(s, t) \\
& = & \theta \Big(\beta_{\phi^{-1}(x),~ \phi^{-1}(y)}\big(\theta^{-1}(s), ~~\theta^{-1}(t)\big)  \Big)\\ 
& = & \theta \Big(\lambda_{\phi^{-1}(x)* \phi^{-1}(y)} \Big(\alpha_{\phi^{-1}(x), ~\phi^{-1}(y)} \Big( \lambda_{\phi^{-1}(x)}^{-1}\big(\theta^{-1}(s) \big), ~~\lambda_{\phi^{-1}(y)}^{-1}\big(\theta^{-1}(t) \big) \Big) \Big) \Big),~ \textrm{by using}~ \eqref{cohomologous}\\
& = & \theta \Big(\lambda_{\phi^{-1}(x*y)} \Big(\alpha_{\phi^{-1}(x),~ \phi^{-1}(y)} \Big( \lambda_{\phi^{-1}(x)}^{-1}\theta^{-1}(s), ~~\lambda_{\phi^{-1}(y)}^{-1}\theta^{-1}(t) \Big) \Big) \Big)\\
& = & \Big( \theta \lambda_{\phi^{-1}(x*y)} \theta^{-1} \Big) \theta \Big(\alpha_{\phi^{-1}(x),~ \phi^{-1}(y)} \Big( \theta^{-1} (\theta \lambda_{\phi^{-1}(x)}\theta^{-1})^{-1}(s), ~~\theta^{-1}(\theta \lambda_{\phi^{-1}(y)}\theta^{-1})^{-1}(t)  \Big) \Big)\\  
& = & \lambda'_{x*y} ~ \theta \Big(\alpha_{\phi^{-1}(x), ~\phi^{-1}(y)} \Big( \theta^{-1} {\lambda'_x}^{-1}(s), ~~\theta^{-1}{\lambda'}_y^{-1}(t)  \Big) \Big)\\  
& = & \lambda'_{x*y}  ~ \Big(^{(\phi, \theta)}\alpha_{x, y} \Big({\lambda'_x}^{-1}(s), ~~{\lambda'_y}^{-1}(t)  \Big) \Big),  
\end{eqnarray*}
and hence $[^{(\phi, \theta)}\alpha]=[^{(\phi, \theta)}\beta]$.
$\blacksquare$    \end{proof}          

\begin{remark}\label{dynamical to ordinary cocycle}
{\rm 
Let $X$ be a quandle, $S=A$ a group (not necessarily abelian), and $\alpha: X \times X \to A$ a map. Then the map $\alpha': X \times X \to \Map(A \times A, A)$ given by $$\alpha'_{x, y}(s, t)= s ~\alpha_{x, y},$$ where $x, y \in X$ and $s, t \in A$, satisfy \eqref{dynamical-cocycle-condition1} and \eqref{dynamical-cocycle-condition3} if and only if $\alpha$ satisfy
\begin{equation}\label{group-coefficient-cocycle-condition1}
\alpha_{x, y}~\alpha_{x*y, z}= \alpha_{x, z}~\alpha_{x* z, y*z}
\end{equation}
and 
\begin{equation}\label{normalised-cocycle-condition}
\alpha_{x, x}=1
\end{equation}
for $x, y, z \in X$. If $A$ is abelian, then \eqref{group-coefficient-cocycle-condition1} is the cocycle condition defining the usual second quandle cohomology $\Ho^2_{\rm Q}(X; A)$ (see Section \ref{Homology} of Chapter  \ref{chapter homology YBE}).}
\end{remark}
\bigskip
\bigskip


\section{Wells-type exact sequence for dynamical extensions of quandles}\label{Wells-type exact sequence via dynamical cohomology}
Let $X$ be a quandle, $S$ a set, $\alpha\in \mathcal{Z}^2_{\rm Q}(X; S)$ a dynamical 2-cocycle and $E=X \times_\alpha S$ the extension of $X$ by $S$  through $\alpha$.  Fixing a base point $x_0 \in X$, we define
$$\Aut^{x_0}(X)= \big\{\phi \in \Aut(X) \, \mid \,  \phi(x_0)=x_0 \big\},$$ 
and 
\begin{small}
$$\Aut^{x_0}_S(E)=\Big\{ \psi \in \Aut(E) \, \mid \, \psi(x, s)=\big(\phi(x),~ \tau(x, s)\big)~\textrm{for some}~\phi \in \Aut^{x_0}(X)~\textrm{and}~\tau \in \Map(X \times S, S) \Big\}.$$
\end{small}

For $x \in X$ and $s \in S$, we denote $\tau(x, s)$ by $\tau_x(s)$.

\begin{lemma}
The following  assertions hold:
\begin{enumerate}
\item If $\psi \in \Aut(E)$ such that $\psi(x, s)=\big(\phi(x), \tau_x(s)\big)$, then $\tau_x \in \Sigma_S$.
\item $\Aut^{x_0}_S(E)$ is a subgroup of $\Aut(E)$.
\end{enumerate}
\end{lemma}

\begin{proof}
Assertion (1) follows from the bijectivity of $\psi$. For assertion (2), let $\psi \in \Aut^{x_0}_S(E)$ such that $\psi(x, s)=\big(\phi(x), \tau_x(s)\big)$ for all $x \in X$ and $s \in S$. Define $\varphi(x,s) = \big(\phi^{-1}(x), ~\tau_{\phi^{-1}(x)}^{-1}(s)\big)$. Since $\phi$ and $\tau_x$ are bijections for all $x \in X$, it follows that $\varphi$ is also a bijection. Further, since $\psi$ is a quandle homomorphism, for $x, y \in X$ and $s, t \in S$, we have
\begin{equation*}
\tau_{x*y}\big(\alpha_{x, y}(s, t)\big)=\alpha_{\phi(x), \phi(y)}\big(\tau_x(s),~ \tau_y(t)\big),
\end{equation*}
which is equivalent to
\begin{equation}\label{useful-in-subgroup}
\alpha_{\phi^{-1}(x), \phi^{-1}(y)} \big(\tau_{\phi^{-1}(x)}^{-1}(s),~ \tau_{\phi^{-1}(y)}^{-1}(t) \big)=\tau_{\phi^{-1}(x*y)}^{-1} \big(\alpha_{x, y}(s,t) \big).
\end{equation}

Now, we compute
\begin{eqnarray*}
\varphi \big( (x, s)*(y, t) \big) &=& \varphi \big(x*y,~ \alpha_{x, y}(s, t) \big)\\
 &=&  \big(\phi^{-1}(x*y), ~\tau_{\phi^{-1}(x*y)}^{-1}\big(\alpha_{x, y}(s, t) \big)\big)\\
  &=&  \big(\phi^{-1}(x)*\phi^{-1}(y), ~\alpha_{\phi^{-1}(x), \phi^{-1}(y)} \big(\tau_{\phi^{-1}(x)}^{-1}(s),~ \tau_{\phi^{-1}(y)}^{-1}(t)\big)\big),~\textrm{by}~\eqref{useful-in-subgroup}\\
  &=&    \big(\phi^{-1}(x), ~\tau_{\phi^{-1}(x)}^{-1}(s)\big)* \big(\phi^{-1}(y), ~\tau_{\phi^{-1}(y)}^{-1}(t)\big)\\
  &=&  \varphi(x,s) * \varphi(y,t).
\end{eqnarray*}
Thus, $\varphi \in \Aut^{x_0}_S(E)$ and a direct check shows that it is the inverse of $\psi$. Further, given $\psi(x, s)=\big(\phi(x), ~\tau_x(s)\big)$ and $\psi'(x, s)=\big(\phi'(x),~ \tau'_x(s)\big)$, we see that
\begin{equation}\label{product-psi}
\psi \psi'(x, s)=\psi \big( \phi'(x), ~\tau'_x(s)\big)=\big(\phi\big( \phi'(x)\big), ~\tau_{\phi'(x)}\big(\tau'_x(s)\big) \big)=\big(\phi \phi'(x), ~\tau_{\phi'(x)}\tau'_x(s) \big),
\end{equation}
and hence $\Aut^{x_0}_S(E)$ is a subgroup of $\Aut(E)$.
$\blacksquare$    \end{proof}          
\medskip

Restricting the action of $\Aut(X) \times \Sigma_S$ to that of its subgroup $\Aut^{x_0}(X) \times \Sigma_S$, let $$\Omega_{[\alpha]}: \Aut^{x_0}(X) \times \Sigma_S \to \mathcal{H}^2_{\rm Q}(X; S)$$ be the orbit map given by 
\begin{equation}\label{wells-map}
\Omega_{[\alpha]}(\phi, \theta)= ~^{(\phi, \theta)}[\alpha].
\end{equation}
For ease of notation, we denote  $\Omega_{[\alpha]}$ by $\Omega$. Further, let $$\Phi: \Aut^{x_0}_S(E) \to \Aut^{x_0}(X) \times \Sigma_S$$ be the restriction map given by 
\begin{equation}\label{restriction-map}
\Phi(\psi)= ~(\phi, \theta),
\end{equation} 
where $\psi(x, s)=\big(\phi(x),~ \tau_x(s)\big)$ and $\theta(s):=\tau_{x_0}(s)$ for all $x \in X$ and $s \in S$.  Define a map $\lambda:X \to \Sigma_S$ by setting $$\lambda_x:= \theta\tau_{\phi^{-1}(x)}^{-1}.$$ 
Note that $\lambda_{x_0}=\id_S$ and $\tau_x$ can be uniquely written as 
 \begin{equation}\label{expression-tau}
 \tau_x=(\tau_x \theta^{-1}) \theta=\lambda_{\phi(x)}^{-1}\theta
\end{equation}
for all $x \in X$. With this set-up, we prove the following result.

\begin{proposition}\label{image of phi equals kernel of Omega}
The map $\Phi$ is a group homomorphism and $\im(\Phi)=\Omega^{-1}\big([\alpha]\big)$.
\end{proposition}

\begin{proof}
If $\psi, \psi' \in \Aut^{x_0}_S(E)$, then by \eqref{product-psi}, we obtain 
$$\psi \psi'(x, s)=\big(\phi \phi'(x), ~\tau_{\phi'(x)}\tau'_x(s) \big),$$
where $\tau_{\phi'(x_0)}\tau'_{x_0}(s)=\tau_{x_0}\tau'_{x_0}(s)=\theta \theta'(s)$. This shows that $\Phi(\psi \psi' )=\Phi(\psi)\Phi(\psi' )$, that is, $\Phi$ is a homomorphism. Note that $$\Omega^{-1}\big([\alpha]\big)= \big( \Aut^{x_0}(X) \times \Sigma_S \big)_{[\alpha]},$$ the stabilizer of the action of $\Aut^{x_0}(X) \times \Sigma_S$ at $[\alpha]$. Let $\psi \in \Aut^{x_0}_S(E)$ such that $\Phi (\psi)= (\phi, \theta)$. Since $\psi$ is a quandle homomorphism, for $x, y \in X$ and $s, t \in S$, we have
\begin{equation}\label{ext-homo-condition1}
\tau_{x*y}\big(\alpha_{x, y}(s, t)\big)=\alpha_{\phi(x), \phi(y)}\big(\tau_x(s),~ \tau_y(t)\big).
\end{equation}
Further, since $\phi \in \Aut^{x_0}(X)$ and $\theta \in \Sigma_S$, we can replace $x, y$ by $\phi^{-1}(x),\phi^{-1}(y)$ and $s, t$ by $\theta^{-1}(s),\theta^{-1}(t)$, respectively,  in \eqref{ext-homo-condition1}. This gives
\begin{equation}\label{ext-homo-condition2}
\tau_{\phi^{-1}(x*y)}\Big(\alpha_{\phi^{-1}(x), \phi^{-1}(y)}\big(\theta^{-1}(s), ~\theta^{-1}(t)\big)\Big)=\alpha_{x, y}\Big(\tau_{\phi^{-1}(x)}\big(\theta^{-1}(s)\big), ~\tau_{\phi^{-1}(y)}\big(\theta^{-1}(t)\big)\Big).
\end{equation}
Define a map $\lambda: X \to \Sigma_S$ by setting $\lambda_x:=\theta \tau^{-1}_{\phi^{-1}(x)}$. Then, equation \eqref{ext-homo-condition2} takes the form
\begin{equation}\label{ext-homo-condition3}
\lambda_{x*y}^{-1} ~\theta \Big(\alpha_{\phi^{-1}(x), ~\phi^{-1}(y)}\big(\theta^{-1}(s), ~\theta^{-1}(t)\big)\Big)=\alpha_{x, y}\big(\lambda_x^{-1}(s), ~\lambda_y^{-1}(t) \big),
\end{equation}
which is further equivalent to
\begin{equation}\label{ext-homo-condition4}
^{(\phi, \theta)}\alpha_{x, y}(s, t)=\theta \Big(\alpha_{\phi^{-1}(x), ~\phi^{-1}(y)}\big(\theta^{-1}(s), ~\theta^{-1}(t)\big)\Big)= \lambda_{x*y} \Big(\alpha_{x, y}\big(\lambda_x^{-1}(s), ~\lambda_y^{-1}(t) \big) \Big).
\end{equation}
This shows that $\im(\Phi)\subseteq \big( \Aut^{x_0}(X) \times \Sigma_S \big)_{[\alpha]}$.
\para
For the converse, let $(\phi, \theta) \in \big( \Aut^{x_0}(X) \times \Sigma_S \big)_{[\alpha]}$. Then, there exists a map $\lambda: X \to \Sigma_S$ such that \eqref{ext-homo-condition4} holds. Define $\psi: E \to E$ by setting
$$\psi(x, s)= \big(\phi(x),~ \lambda_{\phi(x)}^{-1} \theta (s)\big).$$
Since $\phi$ and $\theta$ are bijective, so is $\psi$. By replacing $x, y$ by $\phi(x),\phi(y)$ and $s, t$ by $\theta(s),\theta(t)$, respectively, in \eqref{ext-homo-condition4}, we obtain
\begin{equation}\label{ext-homo-condition5}
\lambda_{\phi(x)*\phi(y)}^{-1}~\theta \big(\alpha_{x, y}(s, ~t)\big)=  \alpha_{\phi(x), \phi(y)}\big(\lambda_{\phi(x)}^{-1}\theta(s), ~\lambda_{\phi(y)}^{-1}\theta(t) \big).
\end{equation}
Now, we compute
\begin{eqnarray*}
\psi \big((x, s)*(y, t) \big) & = & \psi \big(x*y, ~\alpha_{x, y}(s, t) \big) \\
&  =&  \Big(\phi(x)*\phi(y), ~\lambda_{\phi(x)*\phi(y)}^{-1} \theta \big(\alpha_{x, y}(s, t)\big) \Big) \\
&  =&  \Big(\phi(x)*\phi(y), ~\alpha_{\phi(x), \phi(y)}\big(\lambda_{\phi(x)}^{-1}\theta(s), ~\lambda_{\phi(y)}^{-1}\theta(t) \big) \Big),~ \textrm{by}~\eqref{ext-homo-condition5}\\
&  =&  \big(\phi(x), ~\lambda_{\phi(x)}^{-1}\theta(s) \big)* \big(\phi(y), ~\lambda_{\phi(y)}^{-1}\theta(t) \big)\\ 
&  =& \psi (x, s)*\psi (y, t).
\end{eqnarray*}
Thus, $\psi \in\Aut^{x_0}_S(E)$ and $\Phi(\psi)=(\phi, \theta)$, which completes the proof of the proposition.
$\blacksquare$    \end{proof}          

\begin{proposition}\label{kernel Phi description}
$\ker(\Phi) \cong \Big\{\lambda: X \to \Sigma_S \, \mid \, \lambda_{x_0}=\id_S~\textrm{and}~\alpha_{x, y}(s, t)=\lambda_{x*y} \big(\alpha_{x, y}\big(\lambda_x^{-1}(s),\lambda_y^{-1}(t)\big)\big) \Big\}$.
\end{proposition}
\begin{proof}
Note that $\psi \in \ker(\Phi)$ if and only if $\Phi(\psi)=(\id_X, \id_S)$. In view of \eqref{expression-tau}, we have $$\psi(x, s)=\big(x, \lambda_x^{-1}(s)\big),$$ and $\lambda_{x_0}=\id_S$. Further, by \eqref{ext-homo-condition1}, $\psi$ is a quandle homomorphism if and only if $$\alpha_{x, y}(s, t)=\lambda_{x*y} \big(\alpha_{x, y}\big(\lambda_x^{-1}(s),\lambda_y^{-1}(t)\big)\big)$$ for all $x, y \in X$ and $s,t \in S$. Thus, $\ker(\Phi)$ is the desired set via the map $\psi \mapsto \lambda$. Here, in view of \eqref{product-psi}, the right hand side is a group with operation $(\lambda\lambda')_x :=\lambda_x \lambda'_x$ for all $x \in X$.
$\blacksquare$    \end{proof}          

Combining propositions \ref{image of phi equals kernel of Omega} and \ref{kernel Phi description}, we obtain the following exact sequence \cite[Theorem 5.4]{MR4282648}, which we refer to  as a Wells-type exact sequence following \cite{MR0272898}.

\begin{theorem}
Let $X$ be a quandle, $x_0 \in X$, $S$ a set, $\alpha\in \mathcal{Z}^2_{\rm Q}(X; S)$ and $E=X \times_\alpha S$ the extension of $X$ by $S$ using $\alpha$. Then there exists an exact sequence
\begin{equation}\label{well-sequence}
1 \longrightarrow \ker(\Phi) \longrightarrow \Aut^{x_0}_S(E) \stackrel{\Phi}{\longrightarrow} \Aut^{x_0}(X) \times \Sigma_S \stackrel{\Omega}{\longrightarrow} \mathcal{H}^2_{\rm Q}(X; S),
\end{equation}
where exactness at $ \Aut^{x_0}(X) \times \Sigma_S$ means that $\im(\Phi)=\Omega^{-1}\big([\alpha]\big)$.
\end{theorem}
\medskip

As a consequence, we obtain the following short exact sequence.

\begin{corollary}
Let $X$ be a quandle, $x_0 \in X$, $S$ a set, $\alpha\in \mathcal{Z}^2_{\rm Q}(X; S)$ and $E=X \times_\alpha S$ the extension of $X$ by $S$ using $\alpha$. Then there exists a short exact sequence of groups
\begin{equation}\label{der-aut-short-exact}
1 \longrightarrow \ker(\Phi) \longrightarrow \Aut^{x_0}_S(E) \stackrel{\Phi}{\longrightarrow} \big(\Aut^{x_0}(X) \times \Sigma_S\big)_{[\alpha]} \longrightarrow 1.
\end{equation}
\end{corollary}
\medskip

It is worthwhile to investigate the conditions under which the short exact sequence \eqref{der-aut-short-exact} splits.

\begin{proposition}
Let $X$ and $S$ be quandles and $x_0 \in X$. If $\alpha \in \mathcal{Z}^2_{\rm Q}(X; S)$ is the 2-cocycle given by $\alpha_{x, y}(s, t)= s*t$ for all $x, y \in X$ and $s, t \in S$, then the short exact  sequence 
\begin{equation}\label{der-aut-short-exact2}
1 \longrightarrow \ker(\Phi) \longrightarrow \Aut^{x_0}_S(E) \stackrel{\Phi}{\longrightarrow} \big(\Aut^{x_0}(X) \times \Aut(S)\big)_{[\alpha]} \longrightarrow 1
\end{equation}
 splits.
\end{proposition}

\begin{proof}
For the given 2-cocycle $\alpha$, the quandle operation \eqref{genralised-quandle-operation} on $E=X \times_\alpha S$ becomes $$(x, s)* (y,t)= (x* y, s* t ),$$ which is just the product quandle structure. Further, in this case, we have
$$\ker(\Phi) \cong \Big\{\lambda: X \to \Aut(S) \, \mid \, \lambda_{x_0}=\id_S~\textrm{and}~\lambda_{x*y}(s*t)=\lambda_x(s) * \lambda_y(t)\Big\}$$
and
$$\Aut^{x_0}_S(E)=\Big\{ \psi \in \Aut(E) \, \mid \, \psi(x, s)=\big(\phi(x),~ \tau_x(s)\big)~\textrm{for some}~\phi \in \Aut^{x_0}(X)~\textrm{and}~\tau: X \to \Aut(S) \Big\}.$$
Now, for any $(\phi, \theta) \in \Aut^{x_0}(X) \times \Aut(S)$, we see that
$$^{(\phi, \theta)}\alpha_{x, y}(s, t)= \theta \Big(\alpha_{\phi^{-1}(x), \phi^{-1}(y)}\big(\theta^{-1}(s), ~\theta^{-1}(t)\big)  \Big)=\theta \big(\theta^{-1}(s)*\theta^{-1}(t)\big)= s*t=\alpha_{x, y}(s, t)$$
for $x, y \in X$ and $s, t \in S$. Thus, we have $\big(\Aut^{x_0}(X) \times \Aut(S) \big)_{[\alpha]}=\Aut^{x_0}(X) \times \Aut(S)$. Define 
$$\xi:  \Aut^{x_0}(X) \times \Aut(S) \to \Aut^{x_0}_S(E)$$ by setting
$\xi(\phi, \theta)=\psi$, where $\psi(x, s)= \big(\phi(x), \theta(s)\big)$ for $x \in X$ and $s \in S$. Then $\xi$ is a group homomorphism with $\Phi \, \xi$ being the identity map, and sequence \eqref{der-aut-short-exact2} splits.
$\blacksquare$    \end{proof}          
\bigskip
\bigskip


\section{Wells-type exact sequence for quandle cohomology}\label{abelian-quandle-cohomology-section}
Upon specializing to the standard quandle cohomology with coefficients in an abelian group and applying certain appropriate modifications, the exact sequence \eqref{well-sequence} assumes a simpler and more usable form, which we derive in this section. First, we briefly review the quandle cohomology theory presented in Chapter \ref{chapter homology YBE}.
\para

Let $X$ be a quandle and $A$ an abelian group. Let $C_n^{\rm R}(X)$ be the free abelian group generated by $n$-tuples of elements of $X$. Define a homomorphism
$\partial_{n}: C_n^{\rm R}(X) \to C_{n-1}^{\rm R}(X)$ by 

\begin{eqnarray*}
\partial_{n}(x_1, x_2, \dots, x_n)  & = & \sum_{i=2}^{n} (-1)^{i} \big( (x_1, x_2, \dots, x_{i-1}, x_{i+1},\dots,  x_n)\\
& & - (x_1 \ast x_i, x_2 \ast x_i, \dots, x_{i-1}\ast x_i, x_{i+1},x_{i+2},  \dots, x_n) \big)
\end{eqnarray*}
for $n \geq 2$ and $\partial_n=0$ for $n \leq 1$. Then $\{C_n^{\rm R}(X), \partial_n \}$ forms a chain complex. For $n \geq 2$, let $C_n^{\rm D}(X)$ be the free abelian subgroup of $C_n^{\rm R}(X)$ generated by degenerate $n$-tuples $(x_1, \dots, x_n)$ with $x_{i}=x_{i+1}$ for some $1 \le i \le n-1$. Set $C_n^{\rm D}(X)=0$   for $n \leq 1$. It can be checked that if $X$ is a quandle, then $\partial_n\big(C_n^{\rm D}(X)\big) \subseteq C_{n-1}^{\rm R}(X)$. Setting $C_n^{\rm Q}(X) = C_n^{\rm R}(X)/ C_n^{\rm D}(X)$, we see that $\{C_n^{\rm Q}(X), \partial_n \}$ forms a chain complex, where by abuse of notation, $\partial_n$ also denotes the induced homomorphism.  Define $C^n_{\rm Q}(X;A) := \Hom \big(C_n^{\rm Q}(X), A\big)$ and $\delta^n :C^n_{\rm Q}(X;A) \to C^{n+1}_{\rm Q}(X;A)$ given by
$$\delta^n(f)= (-1)^n (f \, \partial_{n+1}).$$ 
This turns $\{C^n_{\rm Q}(X;A), \delta^n \}$ into a cochain complex, and $\Ho^n_{\rm Q}(X;A)$ is the $n$-th quandle cohomology group of  $X$ with coefficients in $A$. We denote the group of $n$-th cocycles and the group of $n$-th coboundaries by  $Z^n_{\rm Q}(X;A)$ and $B^n_{\rm Q}(X;A)$, respectively.
\medskip

Following Remark \ref{dynamical to ordinary cocycle}, we see that $Z^2_{\rm Q}(X;A)$ can be identified with the set of maps $\alpha: X \times X \to A$ satisfying the 2-cocycle conditions \eqref{group-coefficient-cocycle-condition1} and \eqref{normalised-cocycle-condition}. In other words, 
$$\alpha_{x, y}~\alpha_{x*y, z}= \alpha_{x, z}~\alpha_{x* z, y*z}$$
and $$\alpha_{x, x}=1$$
for $x, y, z \in X$.  Further, two 2-cocycles $\alpha, \beta \in Z^2_{\rm Q}(X;A)$ are cohomologous if there exists a 1-cochain $\lambda \in C^1(X;A)$ such that $\beta ~\delta^1 (\lambda)= \alpha $, that is, $$\beta_{x, y}=\lambda_{x*y} \lambda_x^{-1}\alpha_{x, y}$$ for $x, y \in X$. 
\medskip

For $(\phi, \theta) \in \Aut(X) \times \Aut(A)$  and $\alpha \in Z^2_{\rm Q}(X;A)$, we define
\begin{equation}\label{abelian-action}
^{(\phi, \theta)}\alpha_{x, y}:= \theta \big(\alpha_{\phi^{-1}(x), \phi^{-1}(y)}  \big)
\end{equation}
for $x,y \in X$. We prove the following result \cite[Proposition 6.1]{MR4282648}.

\begin{proposition}
The group $\Aut(X) \times \Aut(A)$ acts by automorphisms on $\Ho^2_{\rm Q}(X; A)$ as
$$^{(\phi, \theta)}[\alpha] := [^{(\phi, \theta)}\alpha]$$ 
for $(\phi, \theta) \in \Aut(X) \times \Aut(A)$ and $[\alpha]\in \Ho^2_{\rm Q}(X; A)$.
\end{proposition}

\begin{proof}
Equation  \eqref{abelian-action} defines an action of $\Aut(X) \times \Aut(A)$ on $Z^2_{\rm Q}(X;A)$ follows along the lines of proof of Lemma \ref{aut action on cocycles}. Let $(\phi, \theta) \in \Aut(X) \times \Aut(A)$  and $\alpha, \beta \in Z^2_{\rm Q}(X;A)$. For $x, y \in X$, we see that
\begin{eqnarray*}
^{(\phi, \theta)}(\alpha \beta)_{x, y} &= & \theta \big({(\alpha\beta)}_{\phi^{-1}(x), \phi^{-1}(y)}  \big)\\
&= & \theta \big(\alpha_{\phi^{-1}(x), \phi^{-1}(y)}~\beta_{\phi^{-1}(x), \phi^{-1}(y)}  \big)\\
&= & \theta \big(\alpha_{\phi^{-1}(x), \phi^{-1}(y)}\big) ~ \theta\big(\beta_{\phi^{-1}(x), \phi^{-1}(y)}  \big), ~\textrm{since}~\theta \in \Aut(A)\\
&= &  ^{(\phi, \theta)}\alpha_{x, y} ~^{(\phi, \theta)}\beta_{x, y}\\
&= &  (^{(\phi, \theta)}\alpha ~^{(\phi, \theta)}\beta)_{x, y},
\end{eqnarray*}
and hence the action is by automorphisms. Let $\alpha,\beta \in Z^2_{\rm Q}(X;A)$ be two cohomologous quandle 2-cocycles. Then there exists a map $\lambda:X \to A$ such that $$\beta_{x, y}=\alpha_{x, y}~\lambda_x~ \lambda_{x*y}^{-1}$$
for all $x, y \in X$. If we set $\lambda'_x:=\theta (\lambda_{\phi^{-1}(x)})$, then
\begin{eqnarray*}
^{(\phi, \theta)}\beta_{x, y} &=& \theta \big(\beta_{\phi^{-1}(x), \phi^{-1}(y)}  \big)\\
 &=& \theta \big(\alpha_{\phi^{-1}(x), \phi^{-1}(y)} \lambda_{\phi^{-1}(x)} \lambda_{\phi^{-1}(x*y)}^{-1} \big)\\
 &=& ^{(\phi, \theta)}\alpha _{x, y} ~\theta \big(\lambda_{\phi^{-1}(x)} \big) \theta \big(\lambda_{\phi^{-1}(x*y)}^{-1} \big)\\
 &=& ^{(\phi, \theta)}\alpha _{x, y} ~ \lambda'_x ~{\lambda'_{x*y}}^{-1}
\end{eqnarray*}
for all $x, y \in X$. Thus, the group $\Aut(X) \times \Aut(A)$ acts by automorphisms on  $\Ho^2_{\rm Q}(X; A)$.
$\blacksquare$    
\end{proof}

Let $X$ be a finite quandle and $A$ a finite abelian group. Let $\alpha: X\times X \to A$ be a quandle 2-cocycle and $\big(\Aut(X) \times \Aut(A)\big)_{[\alpha]}$ the stabilizer of the action of  $\Aut(X) \times \Aut(A)$ at $[\alpha]$. Then the cardinality of the orbit of $[\alpha]$ equals
$$\frac{|\Aut(X) \times \Aut(A)|}{|\big(\Aut(X) \times \Aut(A)\big)_{[\alpha]}|}.$$
This yields the following result.

\begin{corollary}
If $X$ is a finite quandle, $A$ a finite abelian group and $\alpha: X\times X \to A$ a quandle 2-cocycle, then
$$|\Ho^2_{\rm Q}(X; A)| \ge \frac{|\Aut(X) \times \Aut(A)|}{|\big(\Aut(X) \times \Aut(A)\big)_{[\alpha]}|}.$$
\end{corollary}
\medskip

Let $X$ be a quandle and $\alpha\in Z^2_{\rm Q}(X; A)$ a fixed quandle 2-cocycle. Note that the group $\Ho^2_{\rm Q}(X; A)$ acts freely and transitively on itself by left multiplication. For any $(\phi, \theta) \in \Aut(X) \times \Aut(A)$, since $[\alpha]$ and $^{(\phi, \theta)}[\alpha]$ are elements of $\Ho^2_{\rm Q}(X; A)$, there exists a unique element $\Theta_{[\alpha]} (\phi, \theta) \in \Ho^2_{\rm Q}(X; A)$ such that
\begin{equation}\label{abelian-wells-map}
^{\Theta_{[\alpha]}  (\phi, \theta)} \big(^{(\phi, \theta)}[\alpha] \big)= [\alpha].
\end{equation}
This defines a map $$\Theta_{[\alpha]} : \Aut(X) \times \Aut(A) \to \Ho^2_{\rm Q}(X; A).$$ 
For convenience of notation, we denote $\Theta_{[\alpha]}$ by $\Theta$. By Remark \ref{dynamical to ordinary cocycle}, the map $\alpha': X \times X \to \Map(A\times A, A)$ given by $\alpha'_{x, y}(s, t)=s~\alpha_{x, y}$, where $x, y \in X$ and $s, t \in A$,  is a dynamical 2-cocycle. Thus, by Proposition \ref{set-cocycle}, the binary operation
\begin{equation}\label{abelian-extension-operation}
(x, s)*(y, t):=(x*y,~s~\alpha_{x, y})
\end{equation}
defines a quandle $E:= X \times_\alpha A$, called the {\it abelian extension} of $X$ by $A$  through $\alpha$. Next, we define
\begin{small}
$$\Aut_A(E)=\Big\{ \psi \in \Aut(E) \, \mid \, \psi(x, s)=\big(\phi(x),~ \lambda_x\theta(s)\big)~\textrm{for some}~(\phi, \theta) \in \Aut(X)\times \Aut(A)~\textrm{and map}~\lambda:X \to A \Big\}.$$
\end{small}

With the preceding set-up, we have the following result.

\begin{proposition}
Let $X$ be a quandle, $A$ an abelian group and $\alpha\in Z^2_{\rm Q}(X; A)$ a quandle 2-cocycle. If $E$ is the abelian extension of $X$ by $A$ through $\alpha$,  then $\Aut_A(E)$ is a subgroup of $\Aut(E)$.
\end{proposition}

\begin{proof}
Let $\psi \in\Aut_A(E)$ such that $\psi(x, s)=\big(\phi(x),~ \lambda_x\theta(s)\big)$ for $x \in X$ and $s \in A$. Define $\varphi(x,s) := \big(\phi^{-1}(x), ~\theta^{-1} \big(\lambda_{\phi^{-1}(x)}^{-1}s\big)\big)$. Since $\phi$ and $\theta$ are bijections, it follows that $\varphi$ is also a bijection. Further, since $\psi$ is a quandle homomorphism, for $x, y \in X$ and $s, t \in S$, we have
\begin{equation}\label{abelian-homo-condition}
\lambda_{x*y} ~\theta (\alpha_{x, y})=\lambda_x~ \alpha_{\phi(x), \phi(y)},
\end{equation}
which is equivalent to
\begin{equation}\label{useful-in-abelian-subgroup}
\theta^{-1} \big(\lambda_{\phi^{-1}(x*y)}^{-1} ~\alpha_{x, y}\big)=\theta^{-1} \big(\lambda_{\phi^{-1}(x)}^{-1}\big)  ~\alpha_{\phi^{-1}(x), \phi^{-1}(y)}.
\end{equation}

Next, we see that
\begin{eqnarray*}
\varphi \big( (x, s)*(y, t) \big) &=& \varphi \big(x*y,~ s ~\alpha_{x, y} \big)\\
 &=&  \big(\phi^{-1}(x*y), ~\theta^{-1} \big(\lambda_{\phi^{-1}(x*y)}^{-1} s ~\alpha_{x, y}\big)\big)\\
  &=&    \big(\phi^{-1}(x)*\phi^{-1}(y), ~\theta^{-1} \big(\lambda_{\phi^{-1}(x)}^{-1}s\big)  ~\alpha_{\phi^{-1}(x), \phi^{-1}(y)}  \big),\\
  & &~\textrm{due to}~ \eqref{useful-in-abelian-subgroup}~\textrm{and the fact that $A$ is abelian}\\
  &=&    \big(\phi^{-1}(x), ~\theta^{-1} \big(\lambda_{\phi^{-1}(x)}^{-1}s\big)\big)* \big(\phi^{-1}(y), ~\theta^{-1} \big(\lambda_{\phi^{-1}(y)}^{-1}t\big) \big)\\
  &=&  \varphi(x,s) * \varphi(y,t).
\end{eqnarray*}
Thus, $\varphi \in \Aut_A(E)$ and a direct check shows that $\varphi$ is the inverse of $\psi$. Further, given $\psi(x, s)=\big(\phi(x), ~\lambda_x \theta(s)\big)$ and $\psi'(x, s)=\big(\phi'(x),~ \lambda'_x \theta'(s)\big)$, we see that
\begin{equation}\label{product-psi-abelian}
\psi \psi'(x, s)=\psi \big( \phi'(x), ~\lambda'_x \theta'(s) \big)=\big(\phi\big( \phi'(x)\big), ~\lambda_{\phi'(x)} \theta \big(\lambda'_x \theta'(s) \big) \big)=\big(\phi \phi'(x), ~\lambda_{\phi'(x)} ~\theta (\lambda'_x) ~\theta\theta'(s) \big),
\end{equation}
and hence $\Aut_A(E)$ is a subgroup of $\Aut(E)$.
$\blacksquare$    \end{proof}          
\medskip

The definition of $\Aut_A(E)$ lead to the natural restriction map $$\Psi: \Aut_A(E) \to \Aut(X) \times \Aut(A)$$ given by $\Psi(\psi)=(\phi, \theta)$.

\begin{proposition}
The map $\Psi$ is a group homomorphism and $\im(\Psi)=\Theta^{-1}(0)$.
\end{proposition}

\begin{proof}
It is clear from \eqref{product-psi-abelian} that $\Psi$ is a group homomorphism. Observe that $$\Theta^{-1}(0)=\big(\Aut(X) \times \Aut(A) \big)_{[\alpha]}.$$  If $(\phi, \theta)\in \big(\Aut(X) \times \Aut(A) \big)_{[\alpha]}$, then there exists a map $\lambda:X \to A$ such that
\begin{equation}\label{showing-psi-kernel}
\theta \big(\alpha_{\phi^{-1}(x), \phi^{-1}(y)} \big)=\alpha_{x, y}~\lambda_x~ \lambda_{x*y}^{-1}
\end{equation}
for all $x, y \in X$. Define $\psi: E \to E$ by setting $\psi(x, s)= \big(\phi(x), \lambda_{\phi(x)} \theta(s)\big)$. Clearly, $\psi$ is bijective and is a quandle homomorphism due to \eqref{showing-psi-kernel}. Further, $\Psi(\psi)=(\phi, \theta)$, and hence $\big(\Aut(X) \times \Aut(A) \big)_{[\alpha]} \subseteq \im(\Psi)$.
\par
For the converse, if $(\phi, \theta)\in \im(\Psi)$, then there exists a map $\lambda: X \to A$ such that \eqref{showing-psi-kernel} holds. But this is equivalent to $^{(\phi, \theta)}[\alpha]=[\alpha]$, which is desired.
 $\blacksquare$    
 \end{proof}          

The next result identifies kernel of $\Psi$ and the group of quandle 1-cocycles.

\begin{proposition}
$\ker(\Psi) \cong Z^1_{\rm Q}(X;A) =\{\lambda: X \to A \, \mid \, \lambda_x=\lambda_{x*y}~\textrm{for all}~x, y \in X\}$.
\end{proposition}

\begin{proof}
Note that $\psi \in \ker(\Psi)$ if and only if $\Psi(\psi)=(\id_X, \id_S)$. This gives $\psi(x, s)=\big(x, \lambda_x~ s\big)$ for all $x \in X$ and $s \in A$. Further, by \eqref{abelian-homo-condition}, $\psi$ is a quandle homomorphism if and only if $\lambda_x=\lambda_{x*y}$. Hence, $\ker(\Psi)$ is the desired group $Z^1_{\rm Q}(X;A)$ of 1-cocycles identified via the map $\psi \mapsto \lambda$.
$\blacksquare$    \end{proof}          
\medskip

The preceding results lead to the following theorem \cite[Theorem 6.6]{MR4282648}.

\begin{theorem}\label{abelian extension main theorem}
Let $X$ be a quandle, $A$ an abelian group, $\alpha\in Z^2_{\rm Q}(X; A)$ and $E=X \times_\alpha A$ the abelian extension of $X$ by $A$ through $\alpha$. Then there exists an exact sequence
\begin{equation}\label{abelian-well-sequence}
1 \longrightarrow Z^1_{\rm Q}(X;A) \longrightarrow \Aut_A(E) \stackrel{\Psi}{\longrightarrow} \Aut(X) \times \Aut(A) \stackrel{\Theta}{\longrightarrow} \Ho^2_{\rm Q}(X; A),
\end{equation}
where exactness at $ \Aut(X) \times \Aut(A)$ means that $\im(\Psi)=\Theta^{-1}(0)$.
\end{theorem}

Theorem \ref{abelian extension main theorem} can be reformulated as follows.

\begin{corollary}
Let $X$ be a quandle, $A$ an abelian group, $\alpha\in Z^2_{\rm Q}(X; A)$ and $E=X \times_\alpha A$ the abelian extension of $X$ by $A$ through $\alpha$. Then there exists a short exact sequence
\begin{equation}\label{abelian-well-sequence reformulattion}
1 \longrightarrow Z^1_{\rm Q}(X;A) \longrightarrow \Aut_A(E) \stackrel{\Psi}{\longrightarrow} \big(\Aut(X) \times \Aut(A)\big)_{[\alpha]} \longrightarrow 1.
\end{equation}
\end{corollary}

Restricting the action of $\Aut(X) \times \Aut(A)$ on $\Ho^2_{\rm Q}(X; A)$ to that of its subgroups $\Aut(X)$ and $\Aut(A)$ gives the following result.

\begin{corollary}
Every automorphism in $\Aut(X)_{[\alpha]}$ and $\Aut(A)_{[\alpha]}$ can be extended to an automorphism in $\Aut_A(E)$.
\end{corollary}
\bigskip

Next, we investigate the map $\Theta$ in detail. Let $X$ be a quandle and $A$ an abelian group. Since the group $\Aut(X) \times \Aut(A)$ acts on the group $\Ho^2_{\rm Q}(X; A)$, we have the semi-direct product group $\Ho^2_{\rm Q}(X; A) \rtimes \big(\Aut(X) \times \Aut(A) \big)$. Further, note that the group $\Ho^2_{\rm Q}(X; A)$ acts on itself by left multiplication.

\begin{proposition}
For $(\phi, \theta) \in  \Aut(X) \times \Aut(A)$ and $[\alpha], [\beta]\in  \Ho^2_{\rm Q}(X; A)$ setting
$$^{[\alpha](\phi, \theta)} [\beta]=~^{[\alpha]}{(^{(\phi, \theta)} [\beta])}$$
defines an action of $~\Ho^2_{\rm Q}(X; A) \rtimes \big(\Aut(X) \times \Aut(A) \big)$ on $\Ho^2_{\rm Q}(X; A)$.
\end{proposition}

\begin{proof}
Let $ (\phi_1, \theta_1), (\phi_2, \theta_2) \in \Aut(X) \times \Aut(A)$ and $[\alpha_1], [\alpha_2], [\beta] \in \Ho^2_{\rm Q}(X; A)$. Then, we see that
\begin{eqnarray*}
^{\big([\alpha_1](\phi_1, \theta_1)\big)\big([\alpha_2](\phi_2, \theta_2)\big)}[\beta] & = & ^{\big([\alpha_1]~^{(\phi_1, \theta_1)}[\alpha_2] \big) \big((\phi_1, \theta_1)(\phi_2, \theta_2)\big)}[\beta]\\
& = & ^{\big([\alpha_1]~^{(\phi_1, \theta_1)}[\alpha_2]\big)}{\big(^{\big((\phi_1, \theta_1)(\phi_2, \theta_2)\big)}[\beta] \big)}\\
& = & ^{[\alpha_1]}{\big(^{^{(\phi_1, \theta_1)}{[\alpha_2]}}{\big(^{(\phi_1, \theta_1)}{\big(^{(\phi_2, \theta_2)}[\beta]\big)} \big)} \big)}\\
& = & ^{[\alpha_1]}{\big({^{(\phi_1, \theta_1)}{[\alpha_2]}}{\big(^{(\phi_1, \theta_1)}{\big(^{(\phi_2, \theta_2)}[\beta]\big)} \big)} \big)},\\
& &  \textrm{since}~^{(\phi_1, \theta_1}{[\alpha_2]} \in  \Ho^2_{\rm Q}(X; A),~\textrm{which acts on itself by left translation}\\
& = & ^{[\alpha_1]}{\big({^{(\phi_1, \theta_1)}{\big([\alpha_2] \big(^{(\phi_2, \theta_2)}[\beta]\big)}} \big)\big)},\\
&  & \textrm{since}~\Aut(X) \times \Aut(A)~\textrm{acts by automorphisms on}~\Ho^2_{\rm Q}(X; A)\\
& = & ^{\big([\alpha_1](\phi_1, \theta_1)\big)}{\big(^{\big([\alpha_2](\phi_2, \theta_2)\big)}[\beta] \big)},\\
& & ~ \textrm{since}~\Ho^2_{\rm Q}(X; A)~\textrm{acts on itself by left translation}.
\end{eqnarray*}
Hence, the group $\Ho^2_{\rm Q}(X; A) \rtimes \big(\Aut(X) \times \Aut(A) \big)$ acts on $\Ho^2_{\rm Q}(X; A)$.
$\blacksquare$    \end{proof}          
\medskip

To understand the map $\Theta$, we recall some group cohomology (see Section \ref{obstruction to SLB arising from RBG} of Chapter \ref{chap RBG and YBE}). Let $A$ be an abelian group and $G$ a group acting on $A$ from left by automorphisms.  The first cohomology of $G$ with coefficients in $A$ is defined as $\Ho^1_{\rm Grp}(G; A)=Z^1(G; A)/B^1(G; A)$, where
$$Z^1(G; A) =  \big\{f:G  \to A \, \mid \,  f(xy)= f(x)~{^{x}}f(y)~ \textrm{for all}\ x,y\in G\big\}$$
is the group of  1-cocycles and
$$B^1(G;A) = \big\{f:G\to A \, \mid \, \textrm{there exists}~a \in A~\textrm{such that}~ f(x)=({^x}a)a^{-1}~\textrm{for all}~x\in G \big\}$$
is the group of 1-coboundaries.
\medskip

Further, a complement of a subgroup $H$ in a group $G$ is another subgroup $K$ of $G$ such that $G=HK$ and $H\cap K=1$.  The following result relating 1-cocycles and complements is well-known \cite[11.1.2]{MR1357169}.

\begin{lemma}\label{1-cocycle-complement}
Let $G$ be a group acting on an abelian group $H$ by automorphisms. Then the map $f \mapsto \{f(g)g \, \mid \, g \in G \}$ gives a bijection from the set $\Z^1(G; H)$ of 1-cocycles to the set $\big\{K \, \mid \, G=HK~\textrm{and}~H \cap K=1 \big\}$ of complements of $H$ in $G$. 
\end{lemma}

The following result gives a description of the map $\Theta$ \cite[Theorem 7.3]{MR4282648}.

\begin{theorem}
The map $\Theta: \Aut(X) \times \Aut(A) \to \Ho^2_{\rm Q}(X; A)$ is a group-theoretical 1-cocycle. Further, if $\Theta$ and $\Theta'$ correspond to two  quandle 2-cocycles, then $\Theta$ and $\Theta'$ differ by a group-theoretical 1-coboundary.
\end{theorem}

\begin{proof}
Suppose that $\Theta=\Theta_{[\alpha]}$ for $[\alpha] \in \Ho^2_{\rm Q}(X; A)$ and $\mathbf{g} \in \Ho^2_{\rm Q}(X; A) \rtimes \big( \Aut(X) \times \Aut(A)\big)$. Then for elements $[\alpha], ~^{\mathbf{g}}[\alpha] \in \Ho^2_{\rm Q}(X; A)$, there exists a unique $[\beta] \in \Ho^2_{\rm Q}(X; A)$ such that $~^{[\beta]}[\alpha]=~^{\mathbf{g}}[\alpha]$. This shows that $[\beta]^{-1}\mathbf{g} \in \mathbb{H}$, where $\mathbb{H}$ is the stabilizer subgroup of $\Ho^2_{\rm Q}(X; A) \rtimes \big( \Aut(X) \times \Aut(A) \big)$ at $[\alpha]$. This shows that
$$ \Ho^2_{\rm Q}(X; A) \rtimes \big( \Aut(X) \times \Aut(A) \big)= \Ho^2_{\rm Q}(X; A)\mathbb{H}.$$ Further, since $\Ho^2_{\rm Q}(X; A)$ acts freely on itself, it follows that $\mathbb{H}$ is the complement of $\Ho^2_{\rm Q}(X; A)$ in the group $\Ho^2_{\rm Q}(X; A) \rtimes \big(\Aut(X) \times \Aut(A) \big)$. By Lemma \ref{1-cocycle-complement}, let $f:\Aut(X) \times \Aut(A) \to \Ho^2_{\rm Q}(X; A)$ be the unique 1-cocycle corresponding to the complement $\mathbb{H}$ of $\Ho^2_{\rm Q}(X; A)$ in $\Ho^2_{\rm Q}(X; A) \rtimes \big(\Aut(X) \times \Aut(A)\big)$. Then, we have 
$$\mathbb{H}= \big\{f(\phi, \theta)(\phi, \theta) \, \mid \, (\phi, \theta) \in \Aut(X) \times \Aut(A) \big\},$$
that is, $$[\alpha]=~^{f(\phi, \theta)(\phi, \theta)}[\alpha]=~^{f(\phi, \theta)}{\big(^{(\phi, \theta)}[\alpha]\big)}.$$
It now follows from the definition of $\Theta$ (see \eqref{abelian-wells-map}) that $f(\phi, \theta)=\Theta (\phi, \theta)$, and hence $\Theta$ is a group-theoretical 1-cocycle.
\para

For the second assertion, let $\Theta=\Theta_{[\alpha]}$ and $\Theta'=\Theta'_{[\alpha']}$, where $[\alpha], [\alpha'] \in \Ho^2_{\rm Q}(X; A)$. Then for any $(\phi, \theta) \in \Aut(X) \times \Aut(A)$, we have 
$$ ^{\Theta(\phi, \theta)} \big(^{(\phi, \theta)}[\alpha] \big)= [\alpha]~\textrm{and}~ ^{\Theta'(\phi, \theta)} \big(^{(\phi, \theta)}[\alpha'] \big)= [\alpha'].$$
Since $\Ho^2_{\rm Q}(X; A)$ acts transitively on itself by left multiplication, there exists a unique $[\beta] \in \Ho^2_{\rm Q}(X; A)$ such that $ ^{[\beta]}[\alpha']=[\alpha]$. This gives 
$$ ^{[\beta]^{-1}}{\big(^{\Theta(\phi, \theta)} \big(^{^{(\phi, \theta)}{[\beta]}}{\big(^{(\phi, \theta)}[\alpha'] \big)} \big)\big)}= ^{\Theta'(\phi, \theta)} \big(^{(\phi, \theta)}[\alpha'] \big).$$
Further, since $[\beta]^{-1} \Theta(\phi, \theta) {^{(\phi, \theta)}{[\beta]}}$, $\Theta'(\phi, \theta)$ and  $^{(\phi, \theta)}[\alpha']$ all lie in $\Ho^2_{\rm Q}(X; A)$, which acts freely on itself, we must have $$[\beta]^{-1} \Theta(\phi, \theta) {^{(\phi, \theta)}{[\beta]}}= \Theta'(\phi, \theta).$$
Thus, $\Theta$ and $\Theta'$ differ by a group-theoretical 1-coboundary, which completes the proof.
$\blacksquare$    \end{proof}          

\begin{remark}
{\rm 
Subsequent to Theorem \ref{abelian extension main theorem}, many works, such as \cite{MR4889273, MR4644858, MR4688775, MR4604853}, have established Wells-type exact sequences for various other solutions to the Yang--Baxter equation.}
\end{remark}

\bigskip
\bigskip


\section{Quandle extensions arising from group extensions}\label{groups-to-quandles-extensions-section}
In this section, we relate extensions of groups to extensions of quandles by applying the functors described in Section \ref{adjoint-functor-section} of Chapter \ref{chapter-adjoint-groups-quandles}. We begin with the following result \cite[Proposition 9.1]{MR4282648}.

\begin{proposition}\label{core conj extension}
Let $1 \to A \to E \to G \to 1$ be an extension of groups and $w(x,y) = y x^{-1} y$ or $y^{-n} x y^n$ for some $n \in \mathbb{Z}$. Then $ Q_w(E) \cong Q_w(G) \times_\alpha Q_w(A)$ for some dynamical 2-cocycle $\alpha$. 
\end{proposition}

\begin{proof}
Let $\pi:E \to G$ be the quotient homomorphism and $\kappa: G \to E$ a set-theoretic section such that $\kappa(1)=1$. By functoriality, 
we have a surjective quandle homomorphism $Q_w(\pi): Q_w(E) \to Q_w(G)$ such that
$$\pi \big (\kappa(x)*\kappa(y) \big)= Q_w(\pi) \big (\kappa(x)*\kappa(y) \big)= x*y= Q_w(\pi) \big (\kappa(x*y) \big)= \pi \big (\kappa(x*y) \big)$$
for all $x, y\in G$. Thus, there exists a unique element, say, $\mu(x, y) \in A$ such that $$\kappa(x)*\kappa(y)= \kappa(x*y)~\mu(x, y).$$
Note that $\mu(x, x)=1$ for all $x \in G$. Further, every element of $E$ can be written uniquely as $\kappa(x)s$ for some $x \in G$ and $s \in A$. For each $x \in G$, there is a bijection $\gamma_x: Q_w(\pi)^{-1}(x) \to A$ given by $\gamma_x\big(\kappa(x)s\big)=s$.  Hence, by Proposition \ref{AGlemma}, we have an isomorphism
$$ Q_w(E) \cong Q_w(G) \times_\alpha A,$$ where $\alpha$ is the dynamical 2-cocycle given by $\alpha_{x, y}(s, t)= \gamma_{x*y} \big(\gamma_x^{-1}(s)* \gamma_y^{-1}(t) \big)$ for $x, y \in G$ and $s, t \in A$. 

\begin{enumerate}
\item[Case 1:] If $w(x,y) = y x^{-1} y$, then
\begin{eqnarray*}
\alpha_{x, y}(s, t) &=& \gamma_{x*y} \big(\gamma_x^{-1}(s)* \gamma_y^{-1}(t) \big)\\
 &=& \gamma_{x*y} \big( \big(\kappa(x)s \big)* \big(\kappa(y)t \big)\big)\\
&=& \gamma_{x*y} \big( \kappa(y)t ~s^{-1}\kappa(x)^{-1} ~\kappa(y)t \big)\\
&=& \gamma_{x*y} \big( \kappa(y)\kappa(x)^{-1} \kappa(y) ~^{\kappa(x)^{-1} \kappa(y)}(t s^{-1}) ~t \big)\\
&=& \gamma_{x*y} \big( \big(\kappa(x)*\kappa(y)\big) ~^{\kappa(x)^{-1} \kappa(y)}(t s^{-1}) ~t \big)\\
&=& \gamma_{x*y} \big( \kappa(x*y)~\mu(x, y) ~^{\kappa(x)^{-1} \kappa(y)}(t s^{-1}) ~t \big)\\
&=&  \mu(x, y) ~^{\kappa(x)^{-1} \kappa(y)}(t s^{-1}) ~t.
\end{eqnarray*}
Recall that the quandle operation in $Q_w(G) \times_\alpha A$ is given by $(x, s)*(y, t)=\big(x*y, \alpha_{x, y}(s, t)\big)$ for $x, y \in G$ and $s, t \in A$. For each fixed $x \in G$, we have $$(x, s)*(x, t)=\big(x, \alpha_{x, x}(s, t)\big)=(x, ts^{-1}t)=(x, s*t),$$ which agrees with the quandle operation in $Q_w(A)$.
\para

\item[Case 2:] If $w(x,y) = y^{-n} x y^n$ for some $n \in \mathbb{Z}$, then
\begin{eqnarray*}
\alpha_{x, y}(s, t) &=& \gamma_{x*y} \big( \gamma_x^{-1}(s)* \gamma_y^{-1}(t) \big)\\
 &=& \gamma_{x*y} \big( \big(\kappa(x)s \big)* \big(\kappa(y)t \big)\big)\\
&=& \gamma_{x*y} \big( (\kappa(y)t)^{n} (\kappa(x)s) (\kappa(y)t)^{-n} \big)\\
&=& \gamma_{x*y} \Big( \big(\kappa(y)^n\kappa(x)\kappa(y)^{-n}\big)~\big(^{\kappa(y)^n\kappa(x)^{-1}\kappa(y)^{-n+1}}t\big)~\big(^{\kappa(y)^n\kappa(x)^{-1}\kappa(y)^{-n+2}}t \big)\cdots\\
& & \cdots \big(^{\kappa(y)^n\kappa(x)^{-1}}t \big)~\big(^{\kappa(y)^n}{(st^{-1})} \big)~\big(^{\kappa(y)^{n-1}}{t^{-1}} \big) \cdots \big(^{\kappa(y)^2}{t^{-1}} \big)~\big(^{\kappa(y)}{t^{-1}} \big) \Big)\\
&=& \gamma_{x*y} \Big( \big(\kappa(x)* \kappa(y) \big)~\big(^{\kappa(y)^n\kappa(x)^{-1}\kappa(y)^{-n+1}}t\big)~\big(^{\kappa(y)^n\kappa(x)^{-1}\kappa(y)^{-n+2}}t \big)\cdots\\
& & \cdots \big(^{\kappa(y)^n\kappa(x)^{-1}}t \big)~\big(^{\kappa(y)^n}{(st^{-1})} \big)~\big(^{\kappa(y)^{n-1}}{t^{-1}} \big) \cdots \big(^{\kappa(y)^2}{t^{-1}} \big)~\big(^{\kappa(y)}{t^{-1}} \big) \Big)\\
&=&  \mu(x, y) ~ \Big(\big(^{\kappa(y)^n\kappa(x)^{-1}\kappa(y)^{-n+1}}t\big)~\big(^{\kappa(y)^n\kappa(x)^{-1}\kappa(y)^{-n+2}}t \big)\cdots\\
& & \cdots \big(^{\kappa(y)^n\kappa(x)^{-1}}t \big)~\big(^{\kappa(y)^n}{(st^{-1})} \big)~\big(^{\kappa(y)^{n-1}}{t^{-1}} \big) \cdots \big(^{\kappa(y)^2}{t^{-1}} \big)~\big(^{\kappa(y)}{t^{-1}} \big) \Big).
\end{eqnarray*}
Taking $x=y=1$ and using the facts that $\mu(1, 1)=1=\kappa(1)$, we obtain $$(1, s)*(1, t)=\big(1, \alpha_{1, 1}(s, t)\big)=\big(1, t^{n}st^{-n}\big)=(1, s*t),$$ which agrees with the quandle operation in $Q_w(A)$. Hence, we have $Q_w(E) \cong Q_w(G) \times_\alpha Q_w(A)$.
\end{enumerate}
This completes the proof.
$\blacksquare$    \end{proof}          

An immediate consequence of Proposition \ref{core conj extension} is the following result, which generalises \cite[Proposition 3.9]{MR3544543}.

\begin{corollary}
If $G$ and $A$ are groups, then 
$$ \Core(G \times A) \cong \Core(G) \times \Core(A) \quad \textrm{and} \quad  \Conj_n(G \times A) \cong \Conj_n(G) \times \Conj_n(A)$$
for each $n \in \mathbb{Z}$.
\end{corollary}

Next, we prove an analogous result for generalised Alexander quandles \cite[Proposition 9.3]{MR4282648}.

\begin{proposition}\label{Alex-prop}
Let $1 \to A \to E \to G \to 1$ be an extension of groups and $f \in \Aut(E)$ such that $f(A)=A$. If $f_1\in \Aut(G)$ and $f_2 \in \Aut(A)$ are the automorphisms induced by $f$,  then $\Alex(E, f) \cong \Alex(G, f_1) \times_\alpha \Alex(A, f_2)$ for some dynamical 2-cocycle $\alpha$. 
\end{proposition}

\begin{proof}
Let $\pi: E \to G$, $\kappa:G \to E$ and $\mu: G \times G \to A$ be as in Proposition \ref{core conj extension}, that is, $\kappa(x)*\kappa(y)= \kappa(x*y)~\mu(x, y)$ for all $x, y \in G$. Note that $\Alex(\pi): \Alex(E, f) \to \Alex(G, f_1)$ is a surjective quandle homomorphism. For each $x \in G$, define a bijection $\gamma_x: \Alex(\pi)^{-1}(x) \to A$ by $\gamma_x \big(\kappa(x)s\big)=s$.  Again, by Proposition \ref{AGlemma}, we obtain
$$ \Alex(E, f) \cong \Alex(G, f_1) \times_\alpha A,$$ where $\alpha$ is the dynamical 2-cocycle given by $\alpha_{x, y}(s, t)= \gamma_{x*y} \big(\gamma_x^{-1}(s)* \gamma_y^{-1}(t) \big)$ for $x, y \in G$ and $s, t \in A$. We check that
\begin{eqnarray*}
\alpha_{x, y}(s, t) &=& \gamma_{x*y} \big( \big(\kappa(x)s \big)* \big(\kappa(y)t \big)\big)\\
&=& \gamma_{x*y} \Big(f \big( \kappa(x)s t^{-1}\kappa(y)^{-1} \big)\big(\kappa(y)t \big) \Big)\\
&=& \gamma_{x*y} \Big(f \big(\kappa(x)\kappa(y)^{-1}\big)\kappa(y)~^{\kappa(y)^{-1}f(\kappa(y))}{f(st^{-1})}~t\Big)\\
&=&\gamma_{x*y} \Big(\big(\kappa(x)*\kappa(y)\big)~^{\kappa(y)^{-1}f(\kappa(y))}{f(st^{-1})}~t\Big)\\
&=& \gamma_{x*y} \Big( \kappa(x*y)~\mu(x, y)~^{\kappa(y)^{-1}f(\kappa(y))}{f(st^{-1})}~t\Big)\\
&=&  \mu(x, y)~^{\kappa(y)^{-1}f(\kappa(y))}{f_2(st^{-1})}~t.
\end{eqnarray*}
Taking $x=y=1$ and using the facts that $\mu(1, 1)=1=\kappa(1)$, we get $$(1, s)*(1, t)=\big(1, \alpha_{1, 1}(s, t)\big)=\big(1, f_2(st^{-1})t\big)=(1, s*t),$$ which is the desired quandle structure on $\Alex(A, f_2)$.
$\blacksquare$    \end{proof}          

As a consequence of Proposition \ref{Alex-prop}, we recover the following special case of \cite[Proposition 3.7]{MR3544543}.

\begin{corollary}
If $G$ and $A$ are groups and  $(f_1, f_2) \in \Aut(G) \times \Aut(A)$, then 
$$\Alex_{(f_1, f_2)}(G \times A) \cong \Alex_{f_1}(G) \times \Alex_{f_2}(A).$$
\end{corollary}
\medskip

It is natural to ask whether an extension of quandles induces a corresponding extension of their adjoint groups. The following example show that this is not the case in general.

\begin{example}{\rm 
Let $\R_4 = \{ 0, 1, 2, 3 \}$ be the dihedral quandle of order four. If $X = \{ u, v \}$ is the two element trivial quandle, then there is a surjective quandle homomorphism $\pi: \R_4 \to X$ given by $\pi(0)=\pi(2)=u$ and $\pi(1)=\pi(3)=v$. Then, by Proposition \ref{AGlemma}, we have 
$$\R_4 \cong X \times_\alpha S,$$ where $\alpha$ is a dynamical 2-cocycle  and $S$ is a set of two elements. Since both $X$ and $S$ are two element trivial quandles, we have $\Adj(X) \cong \Adj(S) \cong \mathbb{Z} \times \mathbb{Z}$.  We claim that $\Adj(\R_4)$ is not an extension of $\mathbb{Z} \times \mathbb{Z}$ by $\mathbb{Z} \times \mathbb{Z}$.
\para
We write $\R_4 = Y_0 \sqcup Y_1$ as a disjoint union of fibers, where $Y_0 = \{0, 2\}$ and $Y_1 = \{1, 3 \}$. Both $Y_1, Y_2$ are the connected components as well as trivial subquandles of $\R_4$. Note that $S_{0} = S_{2}$ with
$$
S_{0}(1) = 3  ~ \textrm{and}~  S_{0}(3) = 1.
$$
Similarly, $S_{1} = S_{3}$ with
$$
S_{1}(0) = 2  ~ \textrm{and}~ S_{1}(2) = 0.
$$
Since  $Y_0,Y_1$ and $X$ are trivial quandles, we have
\begin{eqnarray*}
\Adj(Y_0) &=& \big\langle \textswab{a}_0, \textswab{a}_2 \, \mid \, \textswab{a}_0 \textswab{a}_2 = \textswab{a}_2 \textswab{a}_0 \big\rangle, \\
\Adj(Y_1) &=& \big\langle \textswab{a}_1, \textswab{a}_3 \, \mid \, \textswab{a}_1 \textswab{a}_3 = \textswab{a}_3 \textswab{a}_3 \big\rangle,\\
\Adj(X) &=& \big\langle f_u, f_v \, \mid \, f_u f_v = f_v f_u \big\rangle.
\end{eqnarray*}
Note that the group $\Adj(\R_4)$ is generated by its subgroups $\Adj(Y_0), \Adj(Y_1)$ and has relations
$$
\textswab{a}_1^{-1} \textswab{a}_0 \textswab{a}_1 = \textswab{a}_2,~~~\textswab{a}_1^{-1} \textswab{a}_2 \textswab{a}_1 = \textswab{a}_0,~~~ \textswab{a}_0^{-1} \textswab{a}_1 \textswab{a}_0 = \textswab{a}_3,~~~\textswab{a}_2^{-1} \textswab{a}_1 \textswab{a}_2 = \textswab{a}_3,~~~ \textswab{a}_3^{-1} \textswab{a}_0 \textswab{a}_3 = \textswab{a}_2,~~~\textswab{a}_3^{-1} \textswab{a}_2 \textswab{a}_3 = \textswab{a}_0.
$$
Clearly, the elements $b_0:=\textswab{a}_0 \textswab{a}_2^{-1}$ and  $b_1:=\textswab{a}_1 \textswab{a}_3^{-1}$ lie in the kernel of the homomorphism
$$
\Adj(\pi) : \Adj(\R_4) \to \Adj(X).
$$
We remove the generators $\textswab{a}_2$ and $\textswab{a}_3$ using the equalities $\textswab{a}_2 = \textswab{a}_0 b_0^{-1}$ and $\textswab{a}_3 = \textswab{a}_1 b_1^{-1}$. It follows that the other relations have the form
$$
b_0 = \textswab{a}_1^{-1} \textswab{a}_0^{-1} \textswab{a}_1 \textswab{a}_0,~~b_1 = \textswab{a}_0^{-1} \textswab{a}_1^{-1} \textswab{a}_0 \textswab{a}_1,~~b_1 = b_0^{-1},~~b_0^2 = 1,~~\textswab{a}_0^{-1}b_0 \textswab{a}_0 = b_0,~~\textswab{a}_1^{-1}b_1 \textswab{a}_1 = b_1.
$$
Hence, $\Adj(\R_4)$ has a presentation
$$
\Adj(\R_4) = \big\langle \textswab{a}_0, \textswab{a}_1, b_0 \, \mid \, b_0^2 = 1,~~\textswab{a}_1 \textswab{a}_0 = \textswab{a}_0 \textswab{a}_1 b_0,~~\textswab{a}_1^{-1} \textswab{a}_0 = \textswab{a}_0 \textswab{a}_1^{-1} b_0,~~ [\textswab{a}_0, b_0] =   [\textswab{a}_1, b_0] = 1 \big\rangle
$$
and there is an exact sequence
$$
1 \to \mathbb{Z}_2 \to \Adj(\R_4) \to  \mathbb{Z} \times \mathbb{Z} \to 1,
$$
where $\mathbb{Z}_2 = \langle b_0 \rangle$ and $\mathbb{Z} \times \mathbb{Z}= \langle f_u, f_v \rangle$.}
\end{example}
\bigskip
\bigskip

\section{Maps from group cohomology to quandle cohomology}\label{Homomorphism group cohomology to quandle cohomology}
In view of Remark \ref{remark on trunk map}, if $X$ is a quandle and $\mathcal{A}$ is the category of abelian groups, then a trunk map $A: \mathcal{T}'(X) \to \mathcal{A}$ assigns an abelian group $A_x$ to each $x \in X$ and  group homomorphisms $\eta_{x,y} \colon A_x \rightarrow A_{x*y}$ and $\xi_{x,y} \colon A_y \rightarrow A_{x*y}$ such that 
\begin{equation}\label{homo-quandle-module-1}
\eta_{x*y,z}\, \eta_{x,y}=\eta_{x*z,y*z} \, \eta_{x,z} \quad \textrm{and} \quad 
\eta_{x*y,z}\, \xi_{x,y} = \xi_{x*z,y*z} \, \eta_{x,z}
\end{equation}
for all $x,y,z \in X$.
\para 

A \index{quandle module}{\it quandle module} over $X$ is a trunk map
$A: \mathcal{T}'(X) \to \mathcal{A}$ such that each
$\eta_{x,y}$ is an isomorphism and the identities
\begin{equation}\label{homo-quandle-module-3}
\xi_{x*y, z}(s) =  \big(\eta_{x*z,y*z}\, \xi_{x, z}(s)\big)~ \big(\xi_{x*z, y*z}\, \xi_{y, z}(s)\big)
\end{equation}
and 
\begin{equation}\label{homo-quandle-module-4}
\eta_{z, z}(s)~\xi_{z, z}(s)=s
\end{equation}
hold for all $x,y,z \in X$ and $s \in A_z$. A quandle module over $X$ is called {\it homogeneous} if $A_x=A$ for all $x \in X$. Further, a homogeneous quandle module is called {\it trivial} if $\eta_{x,y}$ is the identity map and $\xi_{y, x}$ is the trivial map for all $x, y \in X$.
\para

\begin{example}\label{homogen-module-example}
{\rm 
Any abelian group $A$ can be turned into a homogeneous quandle module over a quandle $X$ by taking $\eta_{x,y}(s)=s^{-1}$ and $\xi_{x,y}(s)=s^2$ for all $x, y \in X$ and $s \in A$.}
\end{example}

For the rest of this section, we consider only homogeneous quandle modules over a quandle. If $A$ is a homogeneous quandle module over $X$, then a {\it factor set} is a map $\mu: X \times X \to A$ satisfying
\begin{equation}\label{factor-set-def-condition}
\mu(x*y,z) ~ \eta_{x*y,z}\big(\mu(x,y)\big) = \eta_{x*z,y*z}\big(\mu(x,z)\big)  ~\mu(x*z,y*z)~ \xi_{x*z, y*z}\big(\mu(y,z)\big)
\end{equation}
and 
\begin{equation}\label{factor-set-def-condition-2}
\mu(x,x) = 1
\end{equation}
for all $x,y,z \in X$. Identities \eqref{homo-quandle-module-1}-\eqref{homo-quandle-module-4} imply that each factor set $\mu$ defines a quandle extension $E=X \times_\alpha A$ of $X$ by $A$, where
\begin{equation}\label{factor-set-quandle-operation}
(x, s)* (y, t)= \big(x*y, ~\eta_{x,y}(s)~\mu(x, y)~\xi_{x,y}(t) \big)
\end{equation}
for all $x, y \in X$ and $s, t \in A$. Conversely, it follows from \cite[Proposition 3.1]{MR2155522} that every extension of a quandle $X$ by a homogeneous quandle module $A$ over $X$ is determined by a factor set. 
\para

Two factor sets $\mu_1,\mu_2: X \times X \to A$ are {\it cohomologous} if there exists a map $\lambda:X \to A$ such that
\begin{equation}\label{factor-set-equivalence}
\mu_2(x, y) = \mu_1(x, y)~\eta_{x,y} \big(\lambda(x) \big)~\xi_{x,y} \big(\lambda(y) \big)~\lambda(x*y)^{-1}
\end{equation}
 for all $x,y \in X$. It is not difficult to see that two factor sets are cohomologous if and only if the corresponding extensions are equivalent  \cite[Theorem 3.6]{MR2155522}. Thus, the set $\mathcal{H}^2_{\rm Q}(X;A)$ of cohomology classes of factor sets, which forms an abelian group with point-wise multiplication,  is in bijective correspondence with the set of equivalence classes of quandle extensions of $X$ by $A$. When $A$ is a trivial homogeneous quandle module over $X$, then $\mathcal{H}^2_{\rm Q}(X;A)=\Ho^2_{\rm Q}(X;A)$ and the quandle operation \eqref{factor-set-quandle-operation} agrees with \eqref{abelian-extension-operation}.
\para

We recall some group cohomology required for our final two results (see Section \ref{obstruction to SLB arising from RBG} of Chapter \ref{chap RBG and YBE}). Let $G$ be a group acting on an abelian group $A$ from the left by automorphisms. The second cohomology of $G$ with coefficients in $A$ is defined as $\Ho^2_{\rm Grp}(G; A)=Z^2(G; A)/B^2(G; A)$, where
$$Z^2(G; A) =\big\{\nu:G \times G \to A  \, \mid \,   \big(x \cdot \nu(y,z) \big)\nu(x, yz)= \nu(xy, z)\nu(x, y)~ \textrm{for all}\ x, y, z\in G\big\}$$
is the group of  2-cocycles and
$$B^2(G;A) = \big\{\nu:G \times G \to A  \, \mid \,  \textrm{there exists}~\lambda:G \to A~\textrm{with}~ \nu(x, y)=~ \big(x \cdot \lambda(y) \big)\lambda(xy)^{-1} \lambda(x)~\textrm{for all}~x, y\in G \big\}$$
is the group of 2-coboundaries.
\para

It is well-known that $\Ho^2_{\rm Grp}(G;A)$ classifies extensions of groups $1 \to A \to E \to G \to 1$ that give rise to the given action of $G$ on $A$. In particular, $E$ is abelian if and only if any group-theoretical  2-cocycle $\nu$ classifying the extension is symmetric, that is, 
$$\nu(x, y)=\nu(y, x)$$ for all $x, y \in G$. Since any group-theoretical 2-cocycle cohomologous to a symmetric group-theoretical 2-cocycle is itself symmetric, we can consider the set $\Ho^2_{\rm SGrp}(G;A)$ of all symmetric group-theoretical cohomology classes of degree 2. In fact, it turns out that $\Ho^2_{\rm SGrp}(G;A)$ is a subgroup of $\Ho^2_{\rm Grp}(G;A)$. With this set-up, we prove the following result \cite[Theorem 10.2]{MR4282648}.

\begin{theorem}\label{homo from group cohom  to quandle cohom}
Let $G$ and $A$ be abelian groups with $A$ viewed as a trivial $G$-module. Then there is a natural group homomorphism $\Lambda:\Ho^2_{\rm SGrp}(G;A) \to \mathcal{H}^2_{\rm Q} \big(\Core(G);A\big )$ given by $\Lambda \big([\nu] \big) = [\check{\nu}]$.
\end{theorem}

\begin{proof}
Let $1 \to A \to E \to G \to 1$ be an extension of abelian groups and $\kappa:G \to E$ a set-theoretic section with $\kappa(1)=1$. Let  $\nu: G \times G \to A$ be a symmetric group-theoretical 2-cocycle corresponding to the  set-theoretic section $\kappa$, that is, 
\begin{equation}\label{tau-nu}
\kappa(x)\kappa(y)= \kappa(xy)~\nu(x, y)
\end{equation}
for all $x, y \in G$. Define $\check{\nu}: G \times G \to A$ by 
$$\kappa(x)*\kappa(y)=\kappa(x*y)~\check{\nu}(x, y)$$
for all $x, y \in G$. Then, we see that
\begin{eqnarray}\label{mu-nu-relation}
\nonumber \check{\nu}(x, y) &=& \big( \kappa(x)*\kappa(y) \big)\kappa(x*y)^{-1}\\
\nonumber &=& \kappa(y)\kappa(x)^{-1} \kappa(y)\kappa(yx^{-1}y)^{-1}\\
\nonumber  &=& \kappa(y)\kappa(x)^{-1} \kappa(yx^{-1})^{-1} ~\nu(yx^{-1}, y),~\textrm{using \eqref{tau-nu} and the fact that $E$ is abelian}\\
 &=& \nu(yx^{-1}, x)^{-1}~ \nu(yx^{-1}, y).
\end{eqnarray}
As in Example \ref{homogen-module-example}, we view $A$ as a homogeneous quandle module over $\Core(G)$ by taking $\eta_{x,y}(s)=s^{-1}$ and $\xi_{y,x}(s)=s^2$ for all $x, y \in G$ and $s \in A$. By Proposition \ref{core conj extension}, $\Core(E)$ is a quandle extension of $\Core(G)$ by the homogeneous quandle module $A$. Further, $\check{\nu}$ is the factor set of the extension $\Core(E)$ relative to the set-theoretic section $\kappa$ in the sense of \eqref{factor-set-def-condition}. 
\para

We claim that if $\nu_1, \nu_2: G \times G \to A$ are cohomologous  group-theoretical  2-cocycles, then the corresponding factor sets $\check{\nu_1}, \check{\nu_2}$ are cohomologous in the sense of \eqref{factor-set-equivalence}. Let $\lambda:G \to A$ be a map such that $\nu_2(x, y)=\nu_1(x,y) \lambda(x)\lambda(xy)^{-1}\lambda(y)$ for all $x, y \in G$. Then, using \eqref{mu-nu-relation}, we have 
\begin{eqnarray*}
\check{\nu_2}(x, y) &=& \nu_2(yx^{-1}, x)^{-1} \nu_2(yx^{-1}, y)\\
 &=& \big(\nu_1(yx^{-1}, x)^{-1}~  \nu_1(yx^{-1},y) \big) \lambda(y)^2 \lambda(x)^{-1} \lambda(yx^{-1}y )^{-1}\\
 &=& \check{\nu_1}(x, y) \lambda(y)^2 \lambda(x)^{-1} \lambda(x*y )^{-1}\\
 &=& \check{\nu_1}(x, y) ~\xi_{x,y}\big( \lambda(y) \big)~\eta_{x,y} \big(\lambda(x)\big) ~\lambda(x*y )^{-1}
\end{eqnarray*}
for all $x, y \in G$. Thus, by \eqref{factor-set-equivalence}, $\check{\nu_1}$ and $\check{\nu_2}$ are cohomologous and determine equivalent extensions. This defines a map $$\Lambda:\Ho^2_{\rm SGrp}(G;A) \to \mathcal{H}^2_{\rm Q} \big(\Core(G);A\big)$$ given by $\Lambda\big([\nu] \big)= [\check{\nu}]$. If $[\nu_1], [\nu_2] \in \Ho^2_{\rm SGrp}(G;A)$, then we see that
\begin{eqnarray*}
\Lambda\big([\nu_1][\nu_2] \big) &=& \Lambda\big([\nu_1\nu_2] \big) \\
&=& [\check{(\nu_1\nu_2)}]\\
&=& [\check{\nu_1}\check{\nu_2}]\\
&=& [\check{\nu_1}][\check{\nu_2}]\\
&=& \Lambda\big([\nu_1]\big) \Lambda\big( [\nu_2] \big),
\end{eqnarray*}
and hence $\Lambda$ is a group homomorphism. 
\para

For naturality of  $\Lambda$, let $G_1, G_2, A_1, A_2$ be abelian groups and $f:G_2 \to G_1$, $h:A_1 \to A_2$ be group homomorphisms. For each group-theoretical 2-cocycle $\nu:G_1 \times G_1 \to A_1$, we define
$$\nu'(x, y)=h \big(\nu\big(f(x), f(y)\big) \big)$$ for all $x, y \in G_2$. A routine check shows that $\nu \mapsto \nu'$ induces a group homomorphism $\Ho^2_{\rm SGrp}(G_1;A_1) \to \Ho^2_{\rm SGrp}(G_2;A_2)$ given by $[\nu] \mapsto [\nu']$. Similarly, for each factor set $\mu: \Core(G_1) \times \Core(G_1) \to A_1$, define 
$$\mu'(x, y)=h \big(\mu\big(f(x), f(y)\big) \big)$$ for all $x, y \in \Core(G_2)$. A direct check shows that $\mu'$ satisfies the factor set conditions \eqref{factor-set-def-condition} and \eqref{factor-set-def-condition-2}. Since cohomologous factor sets are mapped to cohomologous factor sets, we obtain a group homomorphism  $\mathcal{H}^2_{\rm Q} \big(\Core(G_1); A_1\big) \to \mathcal{H}^2_{\rm Q} \big(\Core(G_2); A_2\big)$ given by $[\mu] \mapsto [\mu']$. Finally, it follows from the construction of the preceding maps that the diagram 
$$
\xymatrix{
\Ho^2_{\rm SGrp}(G_1;A_1) \ar[r]^{[\nu] \mapsto [\nu']} \ar[d]^{\Lambda} & \Ho^2_{\rm SGrp}(G_2;A_2) \ar[d]^{\Lambda} \\
\mathcal{H}^2_{\rm Q} \big(\Core(G_1); A_1\big) \ar[r]^{[\mu] \mapsto [\mu']}  & \mathcal{H}^2_{\rm Q} \big(\Core(G_2); A_2\big)}
$$
commutes, and the map $\Lambda$ is natural.
$\blacksquare$    \end{proof}          

Let $G$ be a group and $A$ a trivial $G$-module. For a group-theoretical 2-cocycle $\nu \in Z^2(G;A)$, define $\breve{\nu}:G \times G \to A$  by 
\begin{equation}\label{nu-nu-breve-relation}
\breve{\nu}(x, y)=\nu(x, y)~ \nu(y, y^{-1}xy)^{-1}
\end{equation}
 for $x, y \in G$. It has been proved in \cite[Theorem 7.1]{MR1990571} that $\breve{\nu} \in Z^2_{\rm Q}\big( \Conj(G);A\big)$. We establish the following analogue of Theorem \ref{homo from group cohom  to quandle cohom}.

\begin{theorem}\cite[Theorem 10.3]{MR4282648}
Let $G$ be a group and $A$ an abelian group viewed as a trivial $G$-module. Then there is a natural group homomorphism $\Gamma:\Ho^2_{\rm Grp}(G;A) \to \Ho^2_{\rm Q}\big(\Conj(G);A\big)$ given by $\Gamma \big([\nu] \big) = [\breve{\nu}]$.
\end{theorem}

\begin{proof}
Let $\nu_1, \nu_2$ be cohomologous  group-theoretical  2-cocycles. Then there exists a map $\lambda:G \to A$ such that $\nu_2(x, y)=\nu_1(x,y) \lambda(x)\lambda(xy)^{-1}\lambda(y)$ for all $x, y \in G$. Using \eqref{nu-nu-breve-relation}, we obtain 
\begin{eqnarray*}
\breve{\nu_2}(x, y) &=& \nu_2(x, y) ~\nu_2(y, y^{-1}xy)^{-1}\\
 &=& \big(\nu_1(x,y) ~\nu_1(y, y^{-1}xy)^{-1} \big)\lambda(x)\lambda(y^{-1}xy)^{-1}\\
 &=& \breve{\nu_1}(x, y)\lambda(x)\lambda(x*y)^{-1}
\end{eqnarray*}
for all $x, y \in G$. Thus, $\breve{\nu_1}$ and $\breve{\nu_2}$ are cohomologous quandle 2-cocycles. This defines a map $$\Gamma:\Ho^2_{\rm Grp}(G;A) \to \Ho^2_{\rm Q} \big(\Conj(G);A\big)$$ given by $\Gamma\big([\nu] \big)= [\breve{\nu}]$. If $[\nu_1], [\nu_2] \in 
\Ho^2(G;A)$, then we see that
\begin{eqnarray*}
\Gamma\big([\nu_1][\nu_2] \big) &=& \Gamma\big([\nu_1\nu_2] \big) \\
&=& [\breve{(\nu_1\nu_2)}]\\
&=& [\breve{\nu_1}\breve{\nu_2}]\\
&=& [\breve{\nu_1}][\breve{\nu_2}]\\
&=& \Gamma\big([\nu_1]\big) \Gamma\big( [\nu_2] \big),
\end{eqnarray*}
and $\Gamma$ is a group homomorphism. The naturality of $\Gamma$ follows in a manner similar to that of $\Lambda$  in Theorem  \ref{homo from group cohom  to quandle cohom}.
$\blacksquare$    \end{proof}          

\begin{remark}
{\rm 
We note that a relation between group extensions of a group $G$ and quandle extensions of the quandle $\Alex(G,\iota_x)$ has been obtained in \cite{MR4203539}, where $\iota_x$ is the inner automorphism of $G$ induced by an element $x$.
}    
\end{remark}
\medskip


\chapter*{Epilogue}
As we conclude this monograph, we continue to encounter a steady stream of new papers on the algebraic theory of Yang--Baxter equation including, skew braces, Rota--Baxter groups, racks, quandles, their generalisations and related structures. It is worth noting that numerous captivating papers have surfaced, which could have easily been included in this monograph. However, a comprehensive discussion of all these developments would lead us astray, and attempting an encyclopaedic treatment would prove impossible for us. We selected topics aligned with our research interests and sincerely apologize to the researchers whose valuable contributions we unfortunately had to exclude. We trust that the monograph has sufficiently equipped curious readers to embark on further reading and exploration of the subject on their own. The avenues for exploration are boundless, and our hope is that we have illuminated some of these paths along the way.

\bibliographystyle{amsplain}

\printindex

\end{document}